\documentclass{article}
\usepackage[T1]{fontenc}
\usepackage{amsthm, fullpage}
\usepackage{amsmath,enumerate}
\usepackage{amssymb}
\usepackage{amsfonts,txfonts}
\usepackage[utf8]{inputenc}
\usepackage{eucal}
\usepackage[all]{xy}
\usepackage{hyperref}

\usepackage{mathrsfs}
\DeclareMathOperator{\ord}{ord}
\DeclareMathOperator{\gl}{gl}

\DeclareMathOperator{\coh}{Coh}
\DeclareMathOperator{\colim}{colim}
\DeclareMathOperator{\Tr}{Tr}
\DeclareMathOperator{\tor}{tor}

\DeclareMathOperator{\Hom}{Hom}
\DeclareMathOperator{\loc}{loc}
\DeclareMathOperator{\qc}{qc}
\DeclareMathOperator{\Supp}{Supp}
\DeclareMathOperator{\Sm}{Sm}
\DeclareMathOperator{\Sat}{Sat}

\newtheorem{prop}{Proposition}[subsubsection]
\newtheorem{theo}[prop]{Theorem}
\newtheorem{coro}[prop]{Corollary}
\newtheorem{lemm}[prop]{Lemma}
\newtheorem{lem}[prop]{Lemma}
\newtheorem*{lemm*}{Lemma}

\newtheorem{exs}[prop]{Examples} 

\theoremstyle{definition}

\newtheorem{empt}[prop]{}
\newtheorem{dfn}[prop]{Definition}
\newtheorem{rem}[prop]{Remark}
\newtheorem{ntn}[prop]{Notation}
\newtheorem{ex}[prop]{Examples}
\newtheorem*{rem*}{Remark}

\theoremstyle{thm}
\newtheorem{thm}[prop]{Theorem}
\newtheorem*{thm*}{Theorem}
\newtheorem*{lem*}{Lemma}
\newtheorem*{cor*}{coro}
\newtheorem*{prop*}{Proposition}

\theoremstyle{dfn}
\newtheorem*{dfn*}{Definition}

\numberwithin{equation}{prop}

\newcommand{\riso}{ \overset{\sim}{\longrightarrow}\, }
\newcommand{\liso}{ \overset{\sim}{\longleftarrow}\, }

\newcommand{\Spec}{\mathrm{Spec}\,}
\newcommand{\Spf}{\mathrm{Spf}\,}

\renewcommand{\sp}{\mathrm{sp}}
\renewcommand{\det}{\mathrm{det}}
\newcommand{\gr}{\mathrm{gr}}

\newcommand{\rig}{\mathrm{rig}}

\newcommand{\coker}{\mathrm{Coker}}

\renewcommand{\AA}{{\mathcal{A}}}

\newcommand{\B}{{\mathcal{B}}}
\newcommand{\E}{{\mathcal{E}}}
\newcommand{\G}{{\mathcal{G}}}
\renewcommand{\H}{{\mathcal{H}}}
\newcommand{\M}{{\mathcal{M}}}
\newcommand{\T}{{\mathfrak{T}}}

\newcommand{\A}{\mathbb{A}}
\renewcommand{\P}{\mathbb{P}}
\newcommand{\F}{\mathbb{F}}

\newcommand{\D}{\mathcal{D}}
\newcommand{\DD}{\mathbb{D}}
\renewcommand{\L}{\mathbb{L}}
\newcommand{\R}{\mathbb{R}}
\newcommand{\Q}{\mathbb{Q}}
\newcommand{\bbQ}{\mathbb{Q}}
\newcommand{\Z}{\mathbb{Z}}
\newcommand{\bbN}{\mathbb{N}}

\newcommand{\hdag}{  \phantom{}{^{\dag} }    }

\usepackage{accents}

\newcommand\underrvec[1]{{\underaccent{\rightarrow}{#1}}}

\newcommand{\stack}[1]{\cite[\href{https://stacks.math.columbia.edu/tag/#1}{Tag #1}]{stackproject}}

\def\bbD{{\mathbb D}}

\def\bbL{{\mathbb L}}

\def\bbN{{\mathbb N}}

\def\bbQ{{\mathbb Q}}
\def\bbR{{\mathbb R}}

\def\bbY{{\mathbb Y}}
\def\bbZ{{\mathbb Z}}

\newcommand{\cA}{\mathcal{A}}
\newcommand{\cB}{\mathcal{B}}
\newcommand{\cC}{\mathcal{C}}
\newcommand{\cD}{\mathcal{D}}
\newcommand{\cE}{\mathcal{E}}
\newcommand{\cF}{\mathcal{F}}

\newcommand{\cI}{\mathcal{I}}
\newcommand{\cJ}{\mathcal{J}}

\newcommand{\cM}{\mathcal{M}}
\newcommand{\cN}{\mathcal{N}}
\newcommand{\cO}{\mathcal{O}}

\newcommand{\cV}{\mathcal{V}}
\newcommand{\cW}{\mathcal{W}}

\newcommand{\fA}{\mathfrak{A}}
\newcommand{\fB}{\mathfrak{B}}
\newcommand{\fC}{\mathfrak{C}}
\newcommand{\fD}{\mathfrak{D}}

\newcommand{\fM}{\mathfrak{M}}

\newcommand{\fP}{\mathfrak{P}}
\newcommand{\fQ}{\mathfrak{Q}}

\newcommand{\fS}{\mathfrak{S}}
\newcommand{\fT}{\mathfrak{T}}
\newcommand{\fU}{\mathfrak{U}}
\newcommand{\fV}{\mathfrak{V}}

\newcommand{\fX}{\mathfrak{X}}
\newcommand{\fY}{\mathfrak{Y}}
\newcommand{\fZ}{\mathfrak{Z}}

\newcommand{\fm}{\mathfrak{m}}

\newcommand{\fp}{\mathfrak{p}}

\def\rmb{{\mathrm b}}
\def\rml{{\mathrm l}}

\def\debrom{
\makeatletter
\renewcommand{\theenumi}{(\roman{enumi})}
\renewcommand{\labelenumi}{\theenumi}
\makeatother\begin{enumerate}[(a)]}
\def\finrom{\end{enumerate}}

\begin{document}

\title{Arithmetic $\cD$-modules over Laurent series fields: the relative case}
\author{Daniel Caro}
\date{August 20, 2026}

\maketitle

\begin{abstract}
Let $k$ be a perfect field of characteristic $p>0$.
Within Berthelot's theory of arithmetic $\mathcal{D}$-modules,
we construct a $p$-adic formalism of Grothendieck's six operations for quasi-projective schemes over the basis $\mathrm{Spec}\, k[[t]]$. 
\end{abstract}

\tableofcontents

\bigskip

\section*{Introduction}
Let $\cV$ be a complete discrete valuation ring of mixed characteristic $(0,p)$,  $\pi$ be a uniformizer, $k:= \cV/\pi \cV$ be its residue field and $K$ be its fraction field. We suppose $k$ perfect.
In order to build a $p$-adic formalism of Grothendieck six operations for $k$-varieties (i.e. separated $k$-schemes of finite type), 
Berthelot introduced an arithmetic avatar of the theory of modules over the differential operator rings. 
The objects appearing in his theory are called arithmetic $\D$-modules or complexes of arithmetic $\D$-modules
(for an introduction, see \cite{Beintro2}). 

Within Berthelot's arithmetic $\D$-modules theory,
such a $p$-adic formalism was already known in different contexts. Let us describe these known cases. 
With N. Tsuzuki  (see \cite{caro-Tsuzuki}),
the author got such a formalism  for overholonomic $F$-complexes of arithmetic $\D$-modules 
(i.e. complexes together with a Frobenius structure)
over realizable $k$-varieties (i.e. 
$k$-varieties which can be embedded into a proper formal $\cV$-scheme). 
Another example was given later with holonomic $F$-complexes of arithmetic $\D$-modules over quasi-projective varieties (\cite{caro-stab-holo}).
In a wider geometrical context, T. Abe 
established a six functors formalism for admissible stacks, 
namely algebraic stacks of finite type with finite diagonal morphism
(see \cite[2.3]{Abe-Langlands}). By cohomological descent, he reduced to the case of quasi-projective $k$-varieties. 
Again, some Frobenius structures are involved in his construction. 
Finally, without Frobenius structure, in \cite{caro-unip}, 
we explained how to build 
such a $p$-adic formalism of Grothendieck's six functors, e.g. with quasi-unipotent complexes of arithmetic $\D$-modules (see \cite{caro-unip}).

Recently, Lazda and P\'al have developped in their book \cite{Lazda-Pal-Book} a theory of overconvergent isocrystals on schemes of finite type over 
$\Spec ~ k[[t]]$. Their constructions are inspired by Berthelot's one. 
One main difference is the use of adic spaces instead of Tate's analytic rigid spaces. 
In the same way as Berthelot's category of overconvergent isocrystals, these  overconvergent isocrystals are stable under
tensor products, pull-backs, duality, extraordinary pull-backs, rigid cohomology. 
But, as in the classical case, we do not have the stability under push-forwards by a closed immersion. Hence, two of Grothendieck's six operations (the push-forward
and the extraordinary push-forward) are missing. 
In order to obtain succefully a $p$-adic formalism of Grothendieck's six operations valid for schemes of finite type over $\Spec k[[t]]$, 
the purpose of this paper is to extend their work in the theory of arithmetic $\cD$-modules. 
If $X$ is a $\Spec k[[t]]$-scheme of finite type, then a theory of arithmetic $\cD$-modules on $X /\Spec k$  corresponds to the  ``absolute'' case whereas the ``relative'' case concerns the study on  $X /\Spec k [[t]]$.
The absolute theory has already been studied in \cite{caro-semistable}. 
Here, we focus on a ``relative'' theory.

\bigskip
Let us now clarify the content of the paper. The notion of weak completion was defined by Monsky-Washnitzer for commutative ring (see \cite{MonskyWashnitzer}). 
In the first chapter we study the weak completion in the non commutative case. We introduce the ring of overconvergent formal power series  
 $D [[T _1, \dots, T _m]] ^{\mathrm{oc}}$ in commuting variables $T _1,\dots, T _m$ on a filtered ring $(D, D _n)$ so that $\gr D$ is noetherian.
We prove that $D [[T _1, \dots, T _m]] ^{\mathrm{oc}}$ is noetherian  (see \ref{oc-noeth}) and is a Zariski ring (see \ref{cor-Zar-ring-oc}).
We check that $D [[T]] ^{\mathrm{oc}}$ (i.e. with one variable) is regular if so is $D$ (see \ref{cor-regular-ringpre}).
When $D$ is a $\cV$-algebra, then we denote by $(D, D _n) ^\dag$ the image of the map $D [[T]] ^{\mathrm{oc}}\to \widehat{D}$ sending $T$ to $\pi$ in the $p$-adic completion of $D$. 
Let $R$ be a commutative $\cV$-algebra which is separated for the $p$-adic topology. Next, we construct the weak completion $D ^\dag$ of an $R$-ring $D$ in the same way as 
Monsky-Washnitzer in the commutative case (see \ref{dfn-wc}). If $D$ is filtered (as above), then we have $D ^\dag \subset (D, D _n) ^\dag$ and we give explicit conditions so that the inclusion is an equality, e.g. in the following situation which implies the desired noetherian property.
Let $\mathfrak{S}: = \Spf \cV [[t]]$, let $s= \Spf \cV$ (resp. $\eta = \Spf \cV [[t]] \{ \frac{1}{t}\}$)  be the formal scheme over $\cV$ which corresponds to the closed point (resp. generic point) of $\fS$.
Let $\fS ^{\sharp}:= (\fS , M _{s})$,  where $M_s$ means the log structure induced by the divisor $V(s)$. 
Let $\fX ^\sharp$ be an affine  $p$-torsion-free log smooth formal $\fS ^\sharp$-scheme.  Then we get that $\Gamma (\fX, \cD ^{(m)} _{\fX ^\sharp/\fS ^{\sharp}}) ^\dag$ is  noetherian (see \ref{lem-dagDDm}). 

Let $\fX$ be a $\fS$-formal scheme  of finite type (see \cite[10.13.3]{EGAI})
and flat over $\fS $ whose generic fiber $\fX _\eta$ is smooth  over $\eta= \Spf \cV [[t]] \{ \frac{1}{t}\}$.
Such an object $\fX$ is called a weakly smooth $\fS$-formal scheme.
We build the sheaf $\widetilde{\cO} _{\fX} $ such that if $\fU$ is an affine open of $\fX$ we have
$\Gamma ( \fU, \widetilde{\cO} _{\fX} ) = ( \Gamma ( \fU, \cO _{\fX} ) _t ) ^\dag$, where the symbol $\dag$ means the $p$-adic weak completion as $\Gamma ( \fU, \cO _{\fX} )$-algebra (see \ref{sheaf-ThmB}).
This is a coherent sheaf whose sections over affine open subsets are noetherian (see \ref{123-coherentB}) and regular (see \ref{lem-regularbis}).
We check that Theorems of type $A$ and $B$ holds, i.e., when $\fX$ is affine, 
1) the functors $M \mapsto \widetilde{\cO} _{\fX} \otimes _{\Gamma (\fX,\widetilde{\cO} _{\fX})} M $ and $\cM \mapsto \Gamma (\fX,\cM)$  induce canonically exact quasi-inverse equivalences of categories between the category of
 finitely generated (resp. and projective) $\Gamma (\fX,\widetilde{\cO} _{\fX})$-modules and coherent (resp. and projective) $\widetilde{\cO} _{\fX}$-modules
and 2) for any coherent $\widetilde{\cO} _{\fX}$-module $\cM$, for any integer $n \geq 1$,  $H ^{n} (\fX,\cM)=0$ (see \ref{equcatBDelta}).
Suppose now $\fX ^{\sharp}/\fS ^{\sharp}$ is log smooth. Let $m\in \bbN$ and 
set $\widetilde{\cD} ^{(m)} _{\fX^\sharp/\fS^\sharp} := \widetilde{\cO} _\fX \otimes _{\cO _\fX}\cD ^{(m)} _{\fX^\sharp/\fS^\sharp}$.
If  $\fU  ^\sharp$ is an affine open of $\fX^\sharp$  such that $\fU ^\sharp/\fS ^\sharp$ has logarithmic coordinates,
then the ring  $(\Gamma (\fU, \widetilde{\cD} ^{(m)}  _{\fX^\sharp/\fS^\sharp}) ^{\dag}$ is (right and left) noetherian 
(see \ref{noeth-weak completion}). We also give a simple concrete description of the later ring 
(see \ref{prop-wkcpdscplevelm}) and we prove this  ring  endowed with the $\pi$-adic filtration is a left and right Zariski ring
 (see \ref{cor-Zar-ring}) and is regular (see \ref{cor-regular-ring}).
Set $\widetilde{\cD}  _{\fX^\sharp/\fS^\sharp} := \widetilde{\cO} _\fX \otimes _{\cO _\fX}\cD  _{\fX^\sharp/\fS^\sharp}$. 
We get a nice description of $(\Gamma (\fU, \widetilde{\cD}   _{\fX^\sharp/\fS^\sharp}) ^{\dag}$ (see \ref{loc-m2dagbis})
and we prove that the extension  $\Gamma (\fU, \widetilde{\cD} _{\fX^\sharp/\fS^\sharp} ) ^\dag \to  \Gamma (\fU  _\eta, \cD  ^\dag _{\fX _\eta/\eta} ) $ is right and left faithfully flat
(see \ref{lem-faithfullflat}).

Suppose  $\fX ^\sharp$ is affine and $\fX ^\sharp/\fS ^\sharp$ has logarithmic coordinates.  We denote by $\mathfrak{B}$ the basis of open subsets of  $\fX$ consisting of {\it principal} open subsets. 
We get a presheaf $\widetilde{\cD} ^{(m)\dag} _{\fX^\sharp/\fS^\sharp}$ on $\mathfrak{B}$  defined by putting
$\Gamma ( \mathfrak{D} (f) ,  \widetilde{\cD} ^{(m)\dag} _{\fX^\sharp/\fS^\sharp}):=  \Gamma (\mathfrak{D} (f) ,  \widetilde{\cD} ^{(m)} _{\fX^\sharp/\fS^\sharp}) ^\dag $
where $f \in \Gamma (\fX, \cO _{\fX})$ and where $\dag$ on the right term means  the $p$-adic weak completion as $\Gamma ( \mathfrak{D} (f) , \cO _{\fX})$-ring.
Put $\widetilde{D} ^{(m)\dag} _{\fX^\sharp/\fS^\sharp}:= \Gamma (\fX, \widetilde{\cD} ^{(m)\dag} _{\fX^\sharp/\fS^\sharp})$. 
Let $M$ be a $\widetilde{D} ^{(m)\dag} _{\fX^\sharp/\fS^\sharp}$-module of finite type. 
We prove  (see \ref{Ddagissheaf}) that the presheaf $M ^\Delta$ on $\mathfrak{B}$ defined by setting 
$$\Gamma ( \mathfrak{D} (f) ,  M ^\Delta):=  \Gamma (\mathfrak{D} (f) ,  \widetilde{\cD} ^{(m)\dag} _{\fX^\sharp/\fS^\sharp}) \otimes _{\widetilde{D} ^{(m)\dag} _{\fX^\sharp/\fS^\sharp}} M$$
is in fact a sheaf (in particular $\widetilde{\cD} ^{(m)\dag} _{\fX^\sharp/\fS^\sharp}$ is a sheaf). Moreover,  for every principal open subset $\fU $ of $\fX$ and integer $i>0$, we have $H ^{i} (\fU,M ^\Delta)=0$ (theorem of type $B$).
Theorem of type $A$ holds: 
The functors $M \mapsto M ^{\Delta}$ and $\cM \mapsto \Gamma (\fX,\cM)$ 
induce exact quasi-inverse equivalences  between the category of
$\widetilde{\cD} ^{(m)\dag} _{\fX ^\sharp/\fS ^{\sharp}}$-coherent sheaves
and that of finitely generated 
$\widetilde{D} ^{(m)\dag} _{\fX ^\sharp/\fS ^{\sharp}}$-modules (see \ref{cor-thADdagm}).
Inverting $p$, we get the sheaf $\widetilde{\cD} ^{(m)\dag} _{\fX^\sharp/\fS^\sharp,\bbQ}$ and the theorems of type $A$ and $B$ hold for 
$\widetilde{\cD} ^{(m)\dag} _{\fX^\sharp/\fS^\sharp,\bbQ}$ (see \ref{cor-thADdagmQ}).

Let $u \colon \fZ \hookrightarrow \fX$ be a closed immersion of codimension $d$ of $\fS$-formal schemes $\fP$ of finite type (see \cite[10.13.3]{EGAI})
and flat over $\fS $ whose generic fibers are smooth  over $\eta= \Spf \cV [[t]] \{ \frac{1}{t}\}$.
Put  $\widetilde{\omega} _{\fZ /\fX} := \mathcal{H}om _{\widetilde{\cO} _{\fZ}}  ( \wedge ^n (\widetilde{\cI}  /\widetilde{\cI}  ^2), \widetilde{\cO} _{\fZ})$, where $\cI$ is the ideal defining $u$.
We prove (see \ref{FLI}) that the {\it fundamental local isomorphism} holds by adding some tildes: 
Let $\cM$ be a coherent  $\widetilde{\cO} _{\fX}$-module. We have the isomorphism
$$ \mathcal{E} xt ^d _{\widetilde{\cO} _{\fX}} (\widetilde{\cO} _{\fZ}, \cM)  \riso  \cM \otimes _{\widetilde{\cO} _{\fX}} \widetilde{\omega} _{\fZ /\fX} . $$
Furthermore, if $\cM$ is $\widetilde{u} ^*$-acyclic then $\mathcal{E} xt ^{d-r} _{\widetilde{\cO} _{\fX}} (\widetilde{\cO} _{\fZ}, \cM)  =0$, for any $r \not = 0$.

Let $\fX ^\sharp$ be a formal log scheme   which is  log smooth and saturated (see Definitions \cite[I.3.12 and II.2.10]{Tsuji-saturated}) over $\fS ^\sharp$ and whose generic fiber
$\fX _\eta$ 
is smooth over $\eta$.
 We prove Frobenius descent holds adding singularites: 
the functors $\cE \to F ^* \cE$  (resp.  $\cM \to F ^\flat \cM$) and $\cF \to F _+\cF $ induce exact  equivalences from the category of  left (resp. right) $\widetilde{\cD} ^{\dag} _{\fX ^{\sharp (s)}/\fS ^{\sharp}}$-modules and that of  left (resp. right) 
$\widetilde{\cD} ^{\dag} _{\fX ^{\sharp}/\fS ^{\sharp}}$-modules (see \ref{eqcatFrobdescexplF+} ). 

In order to extend such equivalences, we introduce the notion of  strictly semi-stable formal schemes as follows:
Let  $X$ be an integral $S$-scheme of finite type which is ``strictly semi-stable over $S$'' in the sense of \cite[2.16]{dejong}.
Since $X _s$ is a divisor with strict normal crossing on $X$, then  
$X ^{\sharp}:= (X, M _{X _s})$ is log smooth over $S ^{\sharp}$, where $M _{X _s}$ is the log structure associated with the strict normal crossing divisor $X _s$ of $X$ (\cite[3.14]{Kato-logFontaine-Illusie}). 
A fine $p$-adic formal log scheme $\fX ^\sharp$ over $\fS ^\sharp$ is said to be ``strictly semistable over $\fS ^\sharp$'' if $\fX ^\sharp$ is log smooth over $\fS ^\sharp$ and if $X ^\sharp$ (the reduction modulo $\pi$) is strictly semistable over $S ^\sharp$. 
We prove that for any strictly semi-stable affine log formal scheme over $\fS ^\sharp$, its global section are regular (see \ref{lem-regular}).
Let $f \colon \fY ^\sharp \to \fX ^{\sharp}$ be a projective (in the sense of \ref{Tr+PNdef}) morphism of strictly semistable log formal schemes over $\fS ^\sharp$.
We suppose that the induced morphism  $f _0\colon Y \to X$ is finite, surjective and radicial. 
We establish that the functors  $f _+$ and $f ^!$ are quasi-inverse equivalences of categories between 
the category of coherent left $\widetilde{\cD} ^{\dag} _{\fX ^{\sharp}/\fS ^{\sharp},\bbQ}$-modules 
which are also $\widetilde{\cO}  _{\fX,\bbQ}$-coherent,
and that of coherent left 
$\widetilde{\cD} ^{\dag} _{\fY ^{\sharp}/\fS ^{\sharp},\bbQ}$-modules
which are also $\widetilde{\cO}  _{\fY,\bbQ}$-coherent (see \ref{univhomeo-eqcat-isoc}).

Let $k [[t]] \to l[[u]]$ be a finite morphism of complete discrete valuation rings, $S ' \to S$ be the corresponding morphism of schemes and 
$\mathfrak{S}'\to \mathfrak{S}$ be a lifting. 
Let $\fX ^{\sharp}$ be a strictly semistable log formal scheme over $\fS ^\sharp$, $\fX ^{\prime \sharp} := \fX ^{\sharp} \times _{\fS^{\sharp}} \fS ^{\prime \sharp}$ be the log smooth formal log scheme over $\fS ^{\prime \sharp}$.
Let $Z$ be a divisor of $X$ containing $X _s$ and $Z':= f ^{-1} (Z)$.
We prove the following finite descent property:  if $\cE$ is  a $\cD ^\dag _{\fX^{\sharp}/\fS^{\sharp}} (\hdag Z) _{\bbQ}$-coherent module 
then it is a coherent $\widetilde{\cD} ^\dag _{\fX^{\sharp}/\fS^{\sharp}, \bbQ}$-module if and only if 
$f _Z ^! (\E) $ is a coherent $\widetilde{\cD} ^\dag _{\fX^{\prime \sharp}/\fS ^{\prime \sharp}, \bbQ}$-module (see \ref{desc-coh-chgbase}).

We  introduce the following notion:
a morphism $u \colon \fX ^\sharp \to \fP ^{\sharp}$ of strictly semistable log formal schemes over $\fS ^\sharp$ 
is {\it weak closed immersion} if the underlying morphism of log formal schemes $\underline{u}$ is a closed immersion  (see Definition \ref{dfn-wclimm}).
Recall in the terminology of log schemes, a weak closed immersion is not an immersion. 
We check that the ideal defined by such a weak closed immersion is regular (see \ref{prop-regular}).
Let $u \colon \fX ^\sharp \to \fP ^{\sharp}$ be a weakly closed immersion of strictly semistable log formal schemes over $\fS ^\sharp$.
We extend Kashiwara theorem as follows: 
the functors $u ^{!}$ and $u _{+}$ induce quasi-inverse equivalence of categories of 
coherent $\widetilde{\cD} ^{\dag} _{\fP ^{\sharp}/\fS ^\sharp,\bbQ}$-modules with support in $X$ 
and that of coherent $\widetilde{\cD} ^{\dag} _{\fX ^{\sharp}/\fS ^\sharp,\bbQ}$-modules.
These functors $u ^{!}$ and $u  _{+}$ are acyclic over these categories (see \ref{Berthelot-Kashiwara-full})

Let $\fP ^\sharp$ be a separated strictly semi-stable log formal scheme over $\fS ^\sharp$.
Let $u _0\colon X ^\sharp \to P ^\sharp$ be a weakly closed immersion of strictly semistable log schemes over $S ^\sharp$. 

Let $(\fP ^{\sharp} _{\alpha}) _{\alpha \in \Lambda}$ be an open covering of  $\fP ^{\sharp}$.
We set $\fP ^{\sharp} _{\alpha \beta}:= \fP ^{\sharp} _\alpha \cap \fP ^{\sharp} _\beta$,
$\fP ^{\sharp} _{\alpha \beta \gamma}:= \fP ^{\sharp} _\alpha \cap \fP ^{\sharp} _\beta \cap \fP ^{\sharp} _\gamma$,
$X ^{\sharp} _\alpha := X ^{\sharp} \cap P ^\sharp _\alpha$,
$X ^{\sharp}_{\alpha \beta } := X ^{\sharp} _\alpha \cap X ^{\sharp} _\beta$ et
$X ^{\sharp}_{\alpha \beta \gamma } := X ^{\sharp} _\alpha \cap X ^{\sharp} _\beta \cap X ^{\sharp} _\gamma $.
We denote by $Y _\alpha := (X ^{\sharp} _\alpha) _\eta$,
$Y _{\alpha \beta} := Y _\alpha \cap Y _\beta$,
$Y _{\alpha \beta \gamma} := Y _\alpha \cap Y _\beta \cap Y _\gamma $,
$j _\alpha$ : $ Y _\alpha
\hookrightarrow X ^{\sharp} _\alpha$,
$j _{\alpha \beta} $ :
$Y _{\alpha \beta}  \hookrightarrow  X ^{\sharp} _{\alpha \beta}$
and
$j _{\alpha \beta \gamma} $ :
$Y _{\alpha \beta \gamma }  \hookrightarrow  X ^{\sharp} _{\alpha \beta \gamma} $
the canonical open immersions.
We suppose that for every $\alpha\in \Lambda$, $X ^{\sharp} _\alpha$ is affine, 
(for instance when the covering $(\fP ^{\sharp} _{\alpha}) _{\alpha \in \Lambda}$ is affine).
Since $P$ is separated, for any $\alpha,\beta ,\gamma \in \Lambda$,
$X ^{\sharp}_{\alpha \beta }$ and $X ^{\sharp}_{\alpha \beta \gamma }$ are also affine.

We construct  the category $\mathrm{Coh} ((\fX ^{\sharp}  _\alpha )_{\alpha \in \Lambda}/\cE ^{\dag}  _K)$ whose objects are families $(\cE _\alpha) _{\alpha \in \Lambda}$
of coherent  $\smash{\widetilde{\cD}} ^{\dag} _{\fX ^{\sharp} _{\alpha} /\fS ^\sharp,\bbQ} $-modules
endowed with glueing data. We denote by $\mathrm{Coh} (X, \fP ^{\sharp}/\cE ^{\dag}  _K)$ the category of coherent $\smash{\widetilde{\cD}} ^{\dag} _{\fP ^{\sharp} /\fS ^\sharp,\bbQ} $-modules with support in $X$. 
We prove that we have canonical quasi-inverse equivalences $u ^!  _0$ and $ u _{0+}$ of categories between 
$\mathrm{Coh} ((\fX ^{\sharp} _\alpha) _{\alpha \in \Lambda}/\cE ^{\dag}  _K)$ and $\mathrm{Coh} (X,\, \fP ^{\sharp} /\cE ^{\dag}  _K)$
(see \ref{prop1}), which corresponds to an extension of Berthelot-Kashiwara's theorem by using glueing data.

Let $T$ be a divisor of $P$ containing $P _s$, $\fU ^\sharp$ (resp. $Y$) be the open subset of $\fP ^{\sharp}$ (resp. $X$) complementary to the support of $T$.
In the third chapter, we prove the functor $\sp _*$ induces an equivalence of categories between the category of  convergent isocrystals on $U/\eta$, and the category of 
coherent  $\cD ^\dag _{\fU/\eta, \bbQ} $-modules which are $\cO _{\fU,\bbQ} $-locally projective of finite type (see \ref{thm-eqcat-cvisoc}).
Concerning the overconvergent case, let    $\mathrm{MIC} ^{\dag\dag} (\fP ^\sharp, T/\fS ^{\sharp}) $ be the category of coherent  $\cD ^\dag _{\fP ^\sharp/\fS ^{\sharp}} (\hdag T) _\bbQ$-modules 
which are $\cO _{\fP} (\hdag T) _\bbQ$-coherent (see notation \cite[11.2.1.4]{Car25}). Its objects are called {\it overcoherent isocrystals} on $(\fP ^\sharp, T/\fS ^{\sharp})$.
Let $\mathrm{Isoc} ^{\dag} ((Y, X,\fP)/\cE ^{\dag}  _K)$ be the category of overconvergent 
isocrystals on $(Y, X,\fP)/\cE ^{\dag}  _K$ (see \cite[2.54.3]{Lazda-Pal-Book}), which is equivalent to 
$\mathrm{MIC} ^{\dag} (Y, X,\fX ^\sharp/\cE ^{\dag}  _K) $,  the category of coherent  $j ^\dag \cO _{]X[}$-modules together with an overconvergent connection (see \cite[2.61]{Lazda-Pal-Book}).
When $X = P$, we establish (see \ref{letterBerthelotCaro2007}) that the following properties hold.
\begin{enumerate}[(a)]
\item The functors $\sp _*$ and $\sp ^*$ induces quasi-inverse equivalences between 
$\mathrm{MIC} ^{\dag} (Y, X,\fP ^\sharp/\cE ^{\dag}  _K) $ and 
$\mathrm{MIC} ^{\dag\dag} (\fP ^\sharp, T/\fS ^{\sharp}) $.
\item Let $\cE$ be a coherent $\cD ^\dag _{\fP ^\sharp/\fS ^{\sharp}} (\hdag T) _\bbQ$-module.
Then $\cE \in \mathrm{MIC} ^{\dag\dag} (\fP ^\sharp, T/\fS ^{\sharp})  $ if and only if 
$\cE | \fU$ is $\cO _{\fU,\bbQ}$-coherent. 
\end{enumerate}

By gluing the canonical equivalence of categories $\sp _*$  (see \ref{ntn-dfnsp+}), we get the following one: 
$$\sp _+ \colon  \mathrm{Isoc} ^{\dag} ((Y, X,\fP)/\cE ^{\dag}  _K) \cong \mathrm{MIC} ^{\dag \dag} (Y, X,\fP ^\sharp/\cE ^{\dag}  _K).$$
This equivalence commutes with Frobenius:  if  $E \in \mathrm{Isoc} ^{\dag} ((Y, X,\fP)/\cE ^{\dag}  _K)$,
and  $E ^\vee$ is its dual, then we have the functorial canonical isomorphism in $E$ : $ \sp _{+} (E ^\vee ) \riso \widetilde{\bbD}  \circ \sp _{+} (E )$. (see \ref{propspetdualsansfrob})

Let $\fX ^{\sharp}$ be a strict semistable log formal schemes over $\fS ^\sharp$. Let $Z$ be a divisor of $X$ containing $X _s$.
In the forth chapter, we establish $\cB _{\fX}^{(\bullet)} (Z)  \in  \smash{\underrightarrow{LM}}   _{\bbQ, \mathrm{coh}} (\overset{^\mathrm{l}}{} \smash{\widehat{\cD}} _{\fX ^{\sharp}/\fS ^{\sharp} } ^{(\bullet)} (s))$.
One main ingredient is the local cohomological functor with support in a strictly semistable weakly closed sub log scheme of $X ^{\sharp}$ (see \ref{ntn-GammaZ}) and its coherence (see \ref{coh-smoothsubsch}).
Using the coherence of the constant coefficient, we can define the local cohomological functor with support in any subscheme of a smooth formal schemes over $\fS $ (see \ref{ntn-HdagnX}). 
We prove that local cohomological functors commute with pushforwards, extraordinary pullbacks and tensor products (see \ref{2.2.18}).
Let $f \colon \fP ' \to \fP$  be a realizable morphism of smooth $\fS$-formal schemes (see definition in the sense of \ref{realizablefscheme}).
In the same way as \cite[13.2.3.4]{Car25}, 
we can prove that for any $\cE ^{\prime (\bullet)} \in \underrightarrow{LD} ^{\rmb} _{\bbQ,\coh} ( \widehat{\cD} _{\fP ^{\prime }/\fS} ^{(\bullet)} (P _s'))$
with proper support over $P$ (i.e., if $X'$ is the support of $\cE ^{\prime (\bullet)}$ in the sense \ref{dfn-support} then the composite
$X ' \hookrightarrow  P' \overset{f}{\to} P$ is proper),  the object
$f ^{ (\bullet)}_+(\cE ^{\prime (\bullet)} ) $
belongs to $\underrightarrow{LD} ^{\rmb} _{\bbQ,\coh} ( \widehat{\cD} _{\fP/\fS} ^{(\bullet)} (P _s))$.
Moreover, such pushforwards commutes with extraordinary inverse image by a smooth morphism (see \ref{theo-iso-chgtbase})
or duality (see the relative duality isomorphism \ref{basechange}).

\medskip 
Suppose now there exists an automorphism  $\sigma \colon \cV \riso \cV$ which is a lifting of the $r$th Frobenius power of $k$ for some integer $r\geq 1$.  
In the sixth chapter, we adapt the construction given in  \cite{caro-unip} of a formalism of Grothendieck six functors. 
For instance, we introduce  the category of {\it oriented rings} which is a replacement of 
the category $\mathrm{DVR}  (\cV)$ whose objects $\cW$ are  $\cV$-algebras which is a complete DVR of mixed characteristic $(0,p)$ with perfect residue field $l$ 
so that there exists an automorphism  $\tau  \colon \cW \riso \cW$ which is a lifting of the $s$th Frobenius power of $l$ for some integer $s\geq 1$ which commutes with $\sigma$.
Roughly speaking the basis $\cV$ in our context is replaced by $\cV [[t]]$.
Via these orientations, we also introduce a notion of {\it special descent of the base}. Such notion allows us to use de Jong's desingularization theorem.
Via Theorem \ref{dfnquprop} and the example \ref{ex-datastableevery},
we explain how to build a data of coefficients 
which contains the constant coefficient, 
which is 
local,
stable under devissages, direct summands, 
local cohomological functors, 
pushforwards, extraordinary pullbacks, 
base change, tensor products, duals,
cohomology
and
special descent of the base.
In order to obtain data that are just as stable, but which also contain isocrystals converging
on $k$-smooth varieties, we propose a second construction
that involves slightly more external tensor products (see Theorem \ref{dfnqupropbis}).

Finally, in the last chapter, we obtain a
$p$-adic formalism for the six Grothendieck functors on pairs of varieties $S$ $(Y,X)$ that
can be embedded in a frame of the form
$(Y,X,\fP)$, where $\fP$ is a realizable smooth formal scheme (i.e., one that can be embedded in a proper smooth scheme over $\fS$)
over $\fS$,
$X$ is a closed subscheme of the special fiber of $\fP$, and $Y$ is an open set in $X$.
For a sufficiently stable data of coefficients $\fC$, a coefficient
on $(Y,X,\fP)$ is a coefficient on $\fP$ whose support lies in $X$ and which has overconvergent singularities along $X \setminus Y$
(i.e., which is isomorphic to its image under $\R \underline{\Gamma} ^\dag _{X \setminus Y}$).
We prove that the coefficients of such a $\fC$
on $(Y,X,\fP)$ are independent of the choice of the
enveloping frame $(Y,X)$  (\ref{ind-CYW}).
When $\fC$ is stable by cohomology,
this independence preserves the $t$-structures.
When $X$ is proper over $S$, then
the coefficients of this $\fC$
on $(Y,X)$ are independent (up to canonical equivalence of categories)
of the choice of those proper varieties $X$ containing $Y$.
This gives rise to a formalism for the six Grothendieck operations on varieties $S$
(which can be contained within a frame).

\subsection*{Acknowledgment}
The author was supported by the Institut Universitaire de France.

\subsection*{Notation}
Let $\cV$ be a complete discrete valuation ring of mixed characteristic $(0,p)$,  $\pi$ be a uniformizer, $k:= \cV/\pi \cV$ be its residue field and $K$ be its fraction field. We suppose $k$ perfect.
A formal $\cV$-scheme $\fX$ means is a (noetherian following Grothendieck's terminology of EGA)  $p$-adic formal scheme endowed with  a structural morphism of $p$-adic formal schemes $\fX \to \Spf \cV$.
 
We put $S := \Spec k [[t]]$, $\mathfrak{S}: = \Spf \cV [[t]]$. 
We denote by $s= \Spec k$ (resp. $\eta = \Spec k ((t))$) the scheme which corresponds to the closed point (resp. generic point) of $S$.
By abuse of notation, we still denote by $s= \Spf \cV$ (resp. $\eta = \Spf \cV [[t]] \{ \frac{1}{t}\}$) 
the formal scheme over $\cV$ which corresponds to the closed point (resp. generic point) of $\fS$.
Let $S ^{\sharp}:= (S , M _{s})$, $\fS ^{\sharp}:= (\fS , M _{s})$ where $M_s$ means the log structure induced by the special fiber $V(s)$, i.e. $M _s = \{a \in \cO _{\fS}\text{, such that $a | \eta \in \cO _{\fS} ^*$.}\}$.  We have $\Gamma (\fS, M _s) = \cV [[t]] \setminus \pi  \cV [[t]]$.
For any integer $i \geq 0$, set $\cV _i \coloneqq \cV /\pi ^{i+1}$ and $S _i := \Spec \cV _i [[t]]$.
For any integer $r\geq0$ we denote by 
$S _{(r)}:= \Spec k [[t]] [T _0, \dots, T _{r}]/( T _0 \cdots T _{r} -t)$
and 
$S _{(r)} ^{\sharp}:= (S _{(r)}, M _{V (t)})$, 
$\fS _{(r)}:= \Spec \cV [[t]] \{ T _0, \dots, T _{r}\}/(T _0 \cdots T _{r} -t)$
and 
$\fS _{(r)} ^{\sharp}:= (\fS _{(r)}, M _{V (t)})$,
where $T _1, \dots , T _{r+1}$ are indeterminates.

Sheaves will be denoted with calligraphic letters and their global sections with the associated straight letter. 
By default, a module means a left module. 
We denote by a hat the $p$-adic completion and if 
$\cE$ is an abelian sheaf of groups, we set $\cE _{\bbQ}:= \cE \otimes _{\Z} \bbQ$. 
Let  $\AA$ be a sheaf of rings.
If $*$ is one of the symboles $+$, $-$, ou $\mathrm{b}$, 
$D ^* ( \AA )$ means the derived category of the complexes of 
(left) $\AA$-modules satisfying the corresponding condition of vanishing of cohomological spacesd. 
When we would like to clarify between right and left,  we will write
$D ^* ( \overset{ ^\mathrm{l}}{}\AA )$ ou $D ^* ( \AA \overset{ ^\mathrm{r}}{})$.
We denote by $D ^{\mathrm{b}} _{\mathrm{coh}} ( \AA )$
the subcategory of  $D  ( \AA )$
of bounded and coherent complexes.

Logarithmic $\fS$-formal schemes will be indicated with gothic letters and their reduction modulo $\pi$ 
with the associated roman letter. 
Logarithmic $\fS$-formal schemes or 
logarithmic $S$-schemes will be quasi-compact. 
Finally, when 
$f \colon \fX ^{\sharp} \to \fP ^{\sharp}$
is a log smooth morphism of log formal  $\fS ^\sharp$-schemes,
for any integer $i \in \bbN$,
we denote by  $ f _{i}   \colon X ^{\sharp} _i \to P ^{\sharp} _i$
the induced morphism modulo $\pi ^{i+1}$.

For $n$ be an integer and $P$ be a fine monoid.  
We denote by $A _P$ the log-scheme whose underlying scheme is $\Spec \Z [ P]$ and whose
a pre log structure is defined by the canonical map $P \to \Z [P]$.
We denote by $\mathfrak{A} _P$ the $p$-adic completion of $A _P$.

Let  $T$ be a fine log scheme over $S$
and $\T$ be a $p$-adic formal fine log scheme over $\fS$.  
We denote by $A _{P,T} $ the $T$-log-scheme $A _P \times _{\Spec \Z} T$
and by $\mathfrak{A} _{P,\T }$ the formal log-scheme 
$\mathfrak{A} _P \times _{\Spf \Z _p} \T$.  

More specifically, we put
$A ^r := A _{\bbN ^r}$,
$\mathfrak{A} ^r := \mathfrak{A} _{\bbN ^r}$,
$A ^r _T := A _{\bbN ^r, T}$,
$\mathfrak{A} ^r _\T := \mathfrak{A} _{\bbN ^r,\T}$.

\section{Weakly completion of level $m$}

\subsection{Noncommutative weak completion of an $R$-ring}

\subsubsection{Overconvergent formal power series}
\label{ntn-ord-D-deg}
In this subsection, let us fix some notation. 
Let  $(D, (D _n) _{n\in \bbN} )$ be a filtered ring
whose (increasing unless otherwise specified) filtration is exhaustive and non negative.
We suppose $\gr D$ is left (resp. right) noetherian.
We denote by $\sigma \colon D \to \gr D:= \oplus _{n\in \bbN} D _{n}/ D _{n-1}$ (we put $D _{-1}=0$) the principal symbol
and we put $\ord = \deg \circ \sigma$, where $\deg \colon \gr D \to \bbN$
is the degree.
We still denote by $\sigma \colon D ^r \to (\gr D) ^r$ the map induced by $\sigma$ coordinate by coordinate.
We denote by $\deg  \colon  (\gr D) ^r \to  \bbN$ the map defined by putting  $\deg (P _1, \dots, P _r) := \max _{i=1} ^r \deg P _i$.
We put $\ord = \deg \circ \sigma$.
Since $\gr D$ is left (resp. right) noetherian,
then $D$ is left (resp. right) notherian. 
We will need the following lemma for which we have not found a reference in the litterature:
\begin{lem}
\label{lem-noetherian}
Let $M$ be a left (resp. right) $D$-module which is a $D$-submodule of a free $D$-module of finite rank. 
Then $M$ has a finite number of generators $g _1,\dots, g _n$ so that
any $g\in M$ may be written 
$g = \sum _{i=1} ^{n} a _i g _i$ 
(resp. $g = \sum _{i=1} ^{n} g _i a _i$) 
with $a _i \in D$ and
$\ord (a _i) \leq \ord (g) - \ord (g _{i})$.
\end{lem}

\begin{proof}
Consider the left case. Let $M$ be a left sub $D$-module of $D ^r$.
Let $D ^* := (\gr D) [T]$ be the ring of polynomials with coefficient over $\gr D$.
This is a graded ring as follows
$D ^* := \oplus _{n\in\bbN} D _{n} ^*$, with 
$D _{n} ^*:= \oplus _{i+j=n} (D _{i}/ D _{i-1}) T ^j$.
We denote by $\deg ^*$ the degree on $D ^*$ corresponding to this filtration.
Consider the mapping 
$D ^r \to (D ^* ) ^r$ defined by
$f = (f _1, \dots, f _r) \mapsto f ^* := ( \sigma (f _1) T ^{\ord (f)-\ord(f_1)}, \dots, \sigma (f _r) T ^{\ord (f)-\ord(f_r)})$.
By definition, 
$\deg ^* (\sigma (f _i) T ^{\ord (f)-\ord(f_i)}) = \ord (f)$, i.e. 
the degree of any coordinate of $f ^*$ is equal to order of $f$.

Let $M ^*$ be the left sub-$D ^*$-module of $(D ^*) ^r$
generated by $\alpha ^* $ with $\alpha \in M$.
Since $\gr D$ is noetherian,
so is $D ^*$.
Hence, there exist $g _1, \dots, g _n\in M$ such that
$g _1 ^*, \dots, g _n ^*$ generate $M ^*$. 
Let $g \in M$. 
We can find $c _1,\dots, c _n \in D ^*$  
such that 
$g ^* = \sum _{i=1} ^{n} c _i  g _i ^*$.
Since the coordinates of $g ^*$ have degree $\ord (g)$,
since the coordinates of $g ^*_i $ have degree $\ord (g _i)$,
we can suppose $c _i$ homogenous of degree $\ord (g) - \ord (g _i)$.

Replacing $T$ by $1$, we get the evaluation morphisms
$\mathrm{ev} _{T=1} \colon  D ^* \to \gr D$,
$\mathrm{ev} _{T=1} \colon  (D ^*) ^r \to (\gr D) ^r$.
Put $b _i := \mathrm{ev} _{T=1} (c _i)$.
Since $\sigma (g) = \mathrm{ev} _{T=1} (g ^*) $, 
we get 
$\sigma(g) = \sum _{i=1} ^{n} b _i \sigma(g _i)$.
We have $\deg (b _i) \leq \deg ^* (c _i) = \ord (g) - \ord (g _i)$.
For any integer $n\in \bbN$, let  $[ ] _{n} \colon  \gr D \to D _{n}/D _{n-1}$ be the canonical projection, and the same over  $(\gr D ) ^r$ coordinates by coordinates. 
Let $a _i \in D _{\deg(b _i)}$ such that  $[b _i ] _{\deg(b _i)} = \sigma (a _i)$. Then $\ord (g - \sum _{i=1} ^{n} a _i g _i) < \ord g$. 
Indeed, since  $\ord (a _i)  = \deg (b _i) \leq \ord (g) - \ord (g _i) $,  then  $\ord (a _i g _i) \leq \ord (g) $ and  $[\sigma (a _i g _i )] _{\ord(g)} = [b _i \sigma(g _i)] _{\ord(g)}$.
Hence, we get $[\sigma (g - \sum _{i=1} ^{n} a _i g _i )] _{\ord(g)}  = [\sigma (g )] _{\ord(g)} -\sum _{i=1} ^{n} [\sigma (a _i g _i )] _{\ord(g)} =
[\sigma (g )] _{\ord(g)} -\sum _{i=1} ^{n} [b _i \sigma(g _i)] _{\ord(g)} = [\sigma (g ) -\sum _{i=1} ^{n} b _i \sigma(g _i)] _{\ord(g)} =  0$. We conclude the proof by induction on the order of $g$.
\end{proof}

\begin{ntn} \label{ntn-ord-D-degbis}
We denote by $D [[T _1, \dots, T _m]]$ the ring of formal power series in commuting variables $T _1,\dots, T _m$ over $D$. 
Recall that as a set, this is $D ^{\bbN ^m}$,  the product $( P '' _{\underline{n}} ) _{\underline{n}\in \bbN ^m}:=  ( P _{\underline{n}} ) _{\underline{n}\in \bbN ^m} (P'  _{\underline{n}} ) _{\underline{n}\in \bbN ^m}$
is defined by  $P '' _{\underline{n}}:= \sum _{\underline{n} _1+\underline{n} _2=\underline{n}}P _{\underline{n} _1}  P '_{\underline{n} _2} $. The variable $T _i$ corresponds to $( P _{\underline{n}} ) _{\underline{n}\in \bbN ^m}$
so that $P _{\underline{n}}=0$ except when $\underline{n}=\epsilon _i$, with $\epsilon _i= (0,\dots, 0,1,0 \dots, 0)$ with the $1$ at the $i$th place. 
Moreover, the ring $D [[T _1, \dots, T _m]]$ is canonically  endowed with a $D$-bimodule structure.
For any  $\underline{n}=(n _1,\dots, n _m) \in \bbN ^m$,  we set $\underline{T} ^{\underline{n}}:=  T _1 ^{n _1} \cdots T _m ^{n _m}$ and  $|\underline{n}| =n_1 +\dots + n _m$.  Let $P \in D [[T _1, \dots, T _m]]$. 
There exists a unique family  of elements of $D$, denoted by $(P _{\underline{n}} ) _{\underline{n}\in \bbN ^m}$, such that $P = \sum _{\underline{n}\in \bbN ^m} P _{\underline{n}} \underline{T} ^{\underline{n}}$.
Such element $P$ is said to be homogeneous of degree $d$ in $T _1,\dots, T _m$ if $P = \sum _{|\underline{n}|=d} P _{\underline{n}} \underline{T} ^{\underline{n}}$.
We denote by $D [T _1, \dots, T _m]$ the polynomial ring in commuting variables $T _1,\dots, T _m$ over $D$, i.e. the subset of $D [[T _1, \dots, T _m]]$ of elements $P$ such that  $P _{\underline{n}} =0$ except for a finite number of $\underline{n}$.
We denote by $S ^{(d)}$ the sub $D$-bimodule of $D [[T _1, \dots, T _m]]$ consisting of elements having homogeneous degree $d$ in $T _1,\dots, T _m$.  Since $S ^{(d)}$ is a free left (or right) $D$-module of finite type,  we get the function  $\ord \colon S ^{(d)} \to \bbN$. 
If $P \in D [[T _1, \dots, T _m]]$, we denote by $P ^{(d)} \in S ^{(d)}$ the elements such that 
$P = \sum _{d\in \bbN} P ^{(d)}$.

We define the ring (see \ref{ntnL(a,b)}) of overconvergent formal power series in commuting variables $T _1,\dots, T _m$ over $D$, denoted  by $D [[T _1, \dots, T _m]] ^{\mathrm{oc}}$,  as the subset of  $D [[T _1, \dots, T _m]]$ consisting of the elements $P$ satisfying the following property : 
there exist a constant $c >0$ such that $\ord P ^{(d)} \leq c (d +1)$ for any $d\in \bbN$. 
We recall that following Hilbert Basis Theorem  $D [T _1, \dots, T _m]$ and $D [[T _1, \dots, T _m]]$ are noetherian (see \cite[3.21]{Rotman-HomAlg}). We will see that so is $D [[T _1, \dots, T _m]] ^{\mathrm{oc}}$ (see \ref{oc-noeth}).
\end{ntn}

\begin{ntn} \label{ntnL(a,b)}
We keep Notation \ref{ntn-ord-D-degbis}. Let $r$ be an integer and  $L:= (D [[ T _1, \dots, T _m]] ^{\mathrm{oc}}) ^r$. For any $a,b \in \R $  we denote by  $L (a,b)$ the subset of $L$ of element  $(P _1,\dots, P _r)$ such that 
$\ord P _i ^{(d)} \leq a d +b$ for any $i= 1,\dots, r$ and any $d \in \bbN$. The following properties are straightforward and implies that  $D [[T _1, \dots, T _m]] ^{\mathrm{oc}}$ do is a ring. 

\begin{enumerate}[(a)]
\item For any $a,a',b,'\in \R$ such that $a\leq a'$ and $b\leq b'$, we have $L (a,b)\subset L (a',b')$. Hence, $D [[ T _1, \dots, T _m]] ^{\mathrm{oc}} = \cup _{a,b\in \bbR} L (a,b)$.
\item For any $a,b\in \R$, for any $f, g \in L (a,b)$ we have $f+g \in L (a,b)$.

\item \label{eq-L(a,b)Tn} Let  $\underline{n} \in \bbN ^m$, $a,b\in \R$ and $f\in L$. We have  $f\in  L (a,b) \Leftrightarrow \underline{T} ^{\underline{n}} f \in L (a,b-a|\underline{n}|)$.

\item \label{eq-L(a,b)Tnbis}
For any $P \in D _n$, for any $a,b\in \R$, for any $g \in L(a,b)$, we have $Pg \in L(a,b+n)$.

\item \label{eq-L(a,b)Tnter} 
For any  $a,b,c\in \R$, for any $f \in L (a,b)$, $g \in L (a,c)$,  $fg \in L (a,b+c)$.
\end{enumerate}
\end{ntn}

\begin{lem}
\label{oc-noeth}
With the notation \ref{ntn-ord-D-degbis}, $D [[T _1, \dots, T _m]] ^{\mathrm{oc}}$ is left (resp. right) noetherian.
\end{lem}

\begin{proof}
 The proof is an analogue of that of \cite{Fulton-noetherian}:
Let $I$ be a left ideal of $D [[T _1, \dots, T _m]] ^{\mathrm{oc}}$.
We denote by $I ^{\geq n}$ the left ideal of $D [[T _1, \dots, T _m]] ^{\mathrm{oc}}$
of elements $P\in I$ such that $P ^{(d)}=0$ for any $d <n$.
We denote by $M ^{(n)}$ the left sub-$D$-module of $S ^{(n)}$ 
of the elements $x \in S ^{(n)}$ such that there exists an element $P$ of $I ^{\geq n}$ satislying $x = P ^{(n)}$.
Since $S ^{(n)} M ^{(n')} \subset M ^{(n+n')}$ for any $n,n'\in\bbN$,
then $M:= \oplus _{n\in\bbN} M ^{(n)}$ is the homogeneous ideal of the polynomial ring $D [T _1, \dots, T _m]$
generated by all the $M ^{(n)}$.
Since $D [T _1, \dots, T _n]$ is noetherian (see \ref{ntn-ord-D-degbis}), 
then there exists a large enough integer $N$ such that for any  $d \geq N$ we have
$M ^{(d)}= S ^{(d-N)} M ^{(N)}$.
We denote by $T ^{(n)}$ the subset of $S ^{(n)}$ of the elements of the form $T _1 ^{d _1} \cdots T _m ^{d _m}$,
with $\sum d _i = n$. 
Since an element of $D$ commutes with the indeterminates $T _1,\dots, T _m$, 
we get
$M ^{(d)}= <T ^{(d-N)} M ^{(N)}>$, 
where $<T ^{(d-N)} M ^{(N)}>$ is the subgroup generated by the elements of the form 
$P x$, with $P \in T ^{(d-N)}$ and 
$x \in M ^{(N)}$.

Let $r$ be an integer large enough such that for any $0\leq j \leq N$, there exist
$Q _{j,1},\dots, Q _{j,r} \in I ^{\geq j}$ such that,
the elements $Q _{j,1} ^{(j)},\dots, Q _{j,r} ^{(j)}$ generate $M  ^{(j)}$.
Making another choice if needed, we can suppose that 
$Q _{N,1} ^{(N)},\dots, Q _{N,r} ^{(N)}$ satisfy the property of 
\ref{lem-noetherian}.
Since we have also, 
$M ^{(d)}= <T ^{(d-N)} M ^{(N)}>$, then such property is satisfied for $d \geq N$, i.e., 
 for any 
$h \in M  ^{(d)}$ we can write $h$ of the form
$h= \sum _{k=1} ^{r} a _k Q _{N,k} ^{(N)}$ 
with $a _k \in S ^{(d-N)}$ and
$\ord (a _k ) \leq \ord (h) - \ord Q _{N,k} ^{(N)}$.

Let us check that the family $Q _{j,k} $ for $0\leq j \leq N$, $1\leq k \leq r$ generate $I$.
Let $g \in I$.  There exists an element $\widetilde{g}$ of the left ideal generated by $Q _{j,k} $ for $0\leq j \leq N-1$, $1\leq k \leq r$,
such that $g -\widetilde{g}\in I ^{\geq N}$. Hence, we reduce to check that $Q _{N,k} $ for  $1\leq k \leq r$ generate $I ^{\geq N}$.
Suppose $g \in I ^{\geq N}$. Let  $c \geq 1$ be such that for any $0\leq k \leq r$ we have $\ord Q _{N,k}  ^{(N+d)} \leq c (d +1)$ for any $d\in \bbN$. 
Increasing $c$ if necessary, we can suppose that $\ord g  ^{(N+d)} \leq c (d+1)$ for any $d\in \bbN$. 

By induction on $d\in\bbN$, we construct for any $0\leq i \leq d$, $1\leq k \leq r$, some elements $a _{k} ^{(i)}\in S ^{(i)}$
such that $\ord (a _{k} ^{(i)}) \leq 2c (i+1)$ and $g - \sum _{k=1} ^{r} \sum _{i=0} ^{d} a _{k} ^{(i)}  Q _{N,k} \in I ^{\geq N+d+1}$.
First, let us check the property for $d=0$. Since $Q _{N,1} ^{(N)},\dots, Q _{N,r} ^{(N)}$ satisfy the property of 
\ref{lem-noetherian}, then we can write $g ^{(N)}$ of the form $g ^{(N)}= \sum _{k=1} ^{r} a _k  ^{(0)} Q _{N,k} ^{(N)}$ 
with $a _k ^{(0)} \in D$ and
$\ord (a _k ^{(0)}) \leq \ord (g ^{(N)}) - \ord Q _{N,k} ^{(N)} \leq c \leq 2c$.
We also have
$g - \sum _{k=1} ^{r}a _{k} ^{(0)} Q _{N,k} \in I ^{\geq N+1}$.
Now, suppose that we have constructed such elements $a _{k} ^{(i)}\in S ^{(i)}$
for any
$0\leq i \leq d-1$, $1\leq k \leq r$.
Put $h:= g - \sum _{k=1} ^{r} \sum _{i=0} ^{d-1} a _{k} ^{(i)}  Q _{N,k} \in I ^{\geq N+d}$.
We get
$h^{(N+d)}= g ^{(N+d)} - \sum _{k=1} ^{r} \sum _{i=0} ^{d-1} a _{k} ^{(i)} Q _{N,k}  ^{(N+d -i)}$.
We compute 
$\ord ( g ^{(N+d)}) \leq c ( d+1) \leq 2c (d +1)$ and
$\ord ( a _{k} ^{(i)} Q _{N,k}  ^{(N+d -i)})=
\ord ( a _{k} ^{(i)} ) +\ord (Q _{N,k}  ^{(N+d -i)})
\leq 2c (i +1) + c (d -i +1) =
c ( d +i +3)\leq 2c( d+1)$ for $i \leq d -1$.
Hence $\ord h^{(N+d)} \leq  2c( d+1)$.
Since $h^{(N+d)}\in M ^{(N+d)}$, then we can write 
$h ^{(N+d)}= \sum _{k=1} ^{r} a _k  ^{(d)} Q _{N,k} ^{(N)}$ 
with $a _k ^{(d)} \in S ^{(d)}$ and
$\ord (a _k ^{(d)}) \leq \ord (h ^{(N+d)}) - \ord Q _{N,k} ^{(N)}  \leq 2c (d +1)$.
Since 
$h \in I ^{\geq N+d}$
and 
$\sum _{k=1} ^{r}  a _{k} ^{(d)}  Q _{N,k} \in I ^{\geq N+d}$, 
we get also
$h -\sum _{k=1} ^{r}  a _{k} ^{(d)}  Q _{N,k}  \in I ^{\geq N+d+1}$.
This concludes the induction.
We get the formula
$g = \sum _{k=1} ^{r} a _k Q _{N,k}$,
where $a _k :=\sum _{i=0} ^{\infty} a _{k} ^{(i)} \in D [[T _1, \dots, T _m]] ^{\mathrm{oc}}$.

\end{proof}

From now, we focus to the case where the number of variables $m$ is be equal to $1$.

\begin{lem}
\label{lem-Jacobson-oc}
Let $R$ be a  left (resp. right) noetherian (non necessarily commutative) ring, 
let $x$ be an element of the center of $R$. 
Let $\widehat{R}:= \underleftarrow{\lim} _i R / x ^i R$ be the $x$-adic completion of $R$.
We suppose that $x$ belongs to the Jacobson radical $J (R)$ of $R$ (see \cite[13]{Isaacs}).
\begin{enumerate}[(a)]

\item Any left (resp. right) $R$-module $M$ of finite type is separated for the $x$-adic topology.
Let $M'$ be a submodule of $M$. Then the topology induced by $M$ on $M'$ is the $x$-adic topology of $M'$.
Moreover, $M'$  is a closed subset of $M$. 

\item For any left (resp. right) $R$-module $M$ of finite type, 
the extension 
$M  \to \widehat{M}$ is faithfully flat 
(in the sense of \cite[I.3.1, Definition 1]{bourbaki}).

\end{enumerate}

\end{lem}

\begin{proof}
1) Put $N:= \cap _{n\in \bbN} x^n M$.  Since $R$ is noetherian, since $x$ is in the center of $R$,
since $M$ is an $R$-module of finite type,  then by using \cite[3.2.3.(i)]{Be1} there exists an integer $n $ such that 
$N \cap x ^{n} M \subset x N$. This yields the inclusions $N \subset N \cap x ^{n} M \subset x N \subset J (R) N\subset N$ which are therefore equalities.
From Nakayama's lemma (in the form of Theorem  \cite[13.11]{Isaacs}),  this implies that $N=0$, i.e. $M$ is separated.  Let $M'$ be a submodule of $M$. 
It follows from  \cite[3.2.3.(i)]{Be1} that the induced by $M$ topology  on $M'$ is the $x$-adic topology of $M'$. As $M$, the quotient of $M/M'$ is separated for the $x$-adic topology.
This yields that  $M'$ is a closed subset of $M$.

2) Let us check that $\widehat{R}$ is right faithfully flat $R$-module. From the part 1), we already know that the extension $R \to \widehat{R}$ is injective. 

i) Let us prove the extension $R \to \widehat{R}$ is left and right flat. From \cite[3.2.2.2]{Be1}, for any integer $n\in \bbN$, the canonical morphism
$R/x ^{n}R \to \widehat{R} / x ^n \widehat{R}$ is an isomorphism.  Hence the $x$-adic topology on $R$ is equal to  the topology induced by the $x$-adic topology of $\widehat{R}$.
Moreover, since $R$ is noetherian, since $x$ is in the center of $R$, then from \cite[3.2.3.(iv)]{Be1} we are done. 

ii) \footnote{useless step}Let $I$ be a right ideal of $R$. By flatness,  the canonical morphism $I \otimes _{R}\widehat{R} \to I \widehat{R}$ is an isomorphism.
From \cite[3.2.3.(iii)]{Be1}, the canonical morphism $I \otimes _{R}\widehat{R} \to \widehat{I}$ is an isomorphism.  From \cite[3.2.2.2]{Be1}, for any integer $n\in \bbN$,
the canonical morphism $I/x ^{n}I \to \widehat{I} / x ^n \widehat{I}$ is an isomorphism.  Hence, $I$ is dense in $\widehat{I}$. From the part 1),  $I$ is closed in $R$. 
We check similarly that  $I \widehat{R}$ is closed in $\widehat{R}$ (indeed using \cite[3.2.2]{Be1} we check that  $\widehat{R}$ is noetherian and moreover $x \in J (\widehat{R})$ because 
$\widehat{R}$ is separated and complete for the $x$-adic topology).  This yields that the closure of $I$ in $\widehat{R}$ is $\widehat{I}=I \widehat{R}$.

ii) Let us prove $I= I\widehat{R} \cap R$. Since the converse inclusion is obvious, it remains to check $I\supset I\widehat{R} \cap R$.
Let $a \in I\widehat{R} \cap R$.  Since $a \in I\widehat{R} $, then there exist $i _n \in I$, $b _n \in \widehat{R}$ such that  $a = i_n + x ^n b _n$. 
Since $a -i _n \in R \cap  x ^n\widehat{R} = x ^n R$,  then there exist $c _n \in R$ such that  $a = i_n + x ^n c _n$. Hence, $a$ is in the closure of $I$ in $R$.  Since $I$ is closed in $R$, this yields $a \in I$ and we are done.

iii) Suppose $I \widehat{R} = \widehat{R}$. Then from ii) we get $I =R$.  From \cite[I.3.1, Proposition 1.(d)]{bourbaki}, we conclude using i). 
\end{proof}

\begin{lem} \label{Jacobson-oc}
With the notation of \ref{ntn-ord-D-degbis}, let $J  (D [[ T]] ^{\mathrm{oc}})$ be the Jacobson radical of $D [[ T]] ^{\mathrm{oc}}$ (see \cite[13]{Isaacs}).
\begin{enumerate}[(a)]
\item Then $T D [[ T]] ^{\mathrm{oc}} \subset J  (D [[ T]] ^{\mathrm{oc}}) $.

\item Moreover, any left (resp. right) $D [[ T]] ^{\mathrm{oc}}$-module $M$ of finite type is separated for the $T$-adic topology and any submodule  of such $M$ is a closed subset. 

\item The $T$-adic completion of $D [[ T]] ^{\mathrm{oc}}$ is $D [[ T]]$ and the 
extension  $D [[ T]] ^{\mathrm{oc}} \to D [[ T]]$ is right and left faithfully flat. 
\end{enumerate}
\end{lem}

\begin{proof}
Let us check the first part. Let $f \in D [[ T]] ^{\mathrm{oc}}$.  Let $P := \sum _{n\in \bbN} f ^n T ^n \in D [[ T]]$ (the sum converges since  $D [[ T]]$ is complete with respect to the $T$-adic topology).  We have $(1 -Tf ) P = 1$.  
Let us check that $P \in D [[ T]] ^{\mathrm{oc}}$. Choose $c \geq 1$ such that  $f \in L (c,1)$ (recall notation \ref{ntnL(a,b)}).
From \ref{ntnL(a,b)}.\ref{eq-L(a,b)Tn} and \ref{ntnL(a,b)}.\ref{eq-L(a,b)Tnter}, we compute that  $f ^n T ^n \in L ( c , n (1-c))\subset  L ( c , 1)$.  Hence, $P \in L ( c , 1) \subset D [[ T]] ^{\mathrm{oc}}$. 
This yields that  $1 -Tf $ is invertible in $D [[ T]] ^{\mathrm{oc}}$. From Theorem \cite[13.9]{Isaacs}, this implies  $T D [[ T]] ^{\mathrm{oc}} \subset J (D [[ T]] ^{\mathrm{oc}})$.
Using \ref{ntnL(a,b)}.\ref{eq-L(a,b)Tn}, we get $D [[ T]] ^{\mathrm{oc}}\cap \underline{T} ^{\underline{n}} D [[ T]] ^{\mathrm{oc}} = \underline{T} ^{\underline{n}} D [[ T]] ^{\mathrm{oc}}$ for any $\underline{n} \in \bbN ^m$. 
Since $D [[ T]]$ is $T$-adically complete, this yields that the $T$-adic completion of $D [[ T]] ^{\mathrm{oc}}$ is $D [[ T]]$.
The rest of the lemma is a consequence of \ref{lem-Jacobson-oc}.
\end{proof}

\begin{coro} \label{cor-Zar-ring-oc}
The ring $D [[ T]] ^{\mathrm{oc}}$ endowed with the $T$-filtration is a left and right Zariski ring (see definition \cite[II.2.1.1]{ZarFilt}).
\end{coro}

\begin{proof}
By taking the $T$-adic filtration on $D [[ T]] ^{\mathrm{oc}}$, i.e. $F _n D [[ T]] ^{\mathrm{oc}} := T ^{\max(0,-n)}D [[ T]] ^{\mathrm{oc}}$, we get  $\gr D [[ T]] ^{\mathrm{oc}} = \gr D [[ T]] $.
Hence,   since  $\gr D [[ T]] ^{\mathrm{oc}}  = \gr D [[ T]]$  is moreover Noetherian, since  $D [[ T]] ^{\mathrm{oc}} \to D [[ T]]$ is faithfully flat (see \ref{Jacobson-oc}),  we conclude by using \cite[II.2.1.2.(4)]{ZarFilt}.
\end{proof}

\begin{dfn} \label{regular}
We recall $D$ is left regular (resp. right regular) if every finitely generated left (resp. right) $D$-module has finite projective dimension.  
\end{dfn}

\begin{coro} \label{cor-regular-ringpre}
Suppose  $D$ left (resp. right) regular (see \ref{regular}). 
The ring $D [[ T]] ^{\mathrm{oc}}$ is left (resp. right) regular.
\end{coro}

\begin{proof}
By taking the $T$-adic filtration on $D [[ T]] ^{\mathrm{oc}}$, i.e. $F _n D [[ T]] ^{\mathrm{oc}} := T ^{\max(0,-n)}D [[ T]] ^{\mathrm{oc}}$,
we get $\gr D [[ T]] ^{\mathrm{oc}} = \gr D [[ T]]  = D [ T]$. 
Following Hilbert's Syzygy Theorem, since  $D$ is  left (resp. right) regular, then the ring $D [T]$ is left (resp. right) regular (see \cite[9.34]{Rotman-HomAlg}).  
Since $D [[ T]] ^{\mathrm{oc}}$ is a left and right Zariski ring
(see \ref{cor-Zar-ring-oc}), then we conclude using Theorem \cite[Theorem 4 of II.3.1 ]{ZarFilt}.
\end{proof}

The goal of the rest of the subsection is to check Corollary \ref{convergencepre}. 
It will be useful later, 
e.g. to check Proposition \ref{convergence} which is a key point in the proof of 
Theorem \ref{thmAcoh}.

\begin{lem} \label{lem1-ThmA}
Let $r$ be an integer, $L:= (D [[ T ]] ^{\mathrm{oc}}) ^r$. Let  $0\to K \to L \to M \to 0$ be an exact sequence of  left $D [[ T]] ^{\mathrm{oc}}$-modules such that $M$ has no  $T$-torsion.  Then there exists $a\in\R _+$ satisfying the following property : 
\begin{enumerate}[(a)]
\item  [(*)]    for any $c \geq a$, $d\in \R _+$, $i \in \bbN$, and $f \in ( T ^i L +K) \cap L (c,d)$ (see notation \ref{ntnL(a,b)}), there exists  $g \in T ^i L \cap L (c, d)$ such that $f -g \in K$. 
\end{enumerate}
\end{lem}

\begin{proof}
Since $M$ has no $T$-torsion,  then $T ^i K = K \cap T ^i L$ for any $i\in \bbN$ and  $K /TK$ is a left $D$-submodule of  $L/TL = D ^r$. 
If $g \in L$ (resp. $g \in K$), we denote by $\overline{g}$ its image in $L/TL$ (resp. $K/TK$). From Lemma \ref{lem-noetherian},  there exists $g _1,\dots, g _t \in K$ such that 
for any $x\in K/TK$ there exists  $a _1,\dots, a _t \in D$ such that  $a = \sum _{j=1} ^{t} a _j \overline{g} _j$ and  $\ord (a _j) \leq \ord (\overline{g})  -\ord (\overline{g} _j) $.
Let $a\in \R _+$ be large enough so that $g _j \in L (a,\ord (\overline{g} _j))$ for any $j=1,\dots, t$.

We proceed by induction on $i \in \bbN$, let us check that such choice of $a$ suits.
When $i = 0$, this is obvious.
Let  $c \geq a$, $d\in \R _+$, $i \geq 1$,
and $f \in ( T ^i L +K) \cap L (c,d)$.
By induction hypothesis, 
there exists 
$h \in 
T ^{i-1} L \cap L (c, d )$ such that $f -h \in K$.
Let $h ' \in L$ such that $h = T ^{i-1} h'$.
Since $M$ has no $T$-torsion, 
since $f \in T ^i L +K$ and $f - T ^{i-1} h' \in K$,
then we check that $h' \in T L +K$.
Hence, $\overline{h'} \in K /TK$.
We can write 
$\overline{h'} = 
\sum _{j=1} ^{t}
a _j \overline{g} _j$
with $a _1,\dots, a _t \in D$ such that 
$\ord (a _j) \leq \ord (\overline{h'})  - \ord (\overline{g} _j) $.
Using \ref{ntnL(a,b)}.\ref{eq-L(a,b)Tn}, since $h\in L (c, d )$, then we get
$h' \in L (c, d +(i-1)c)$. 
This yields
$\ord (\overline{h'}) \leq d +(i-1)c$.
Hence,
$ \ord (a _j) \leq d +(i-1)c -\ord (\overline{g} _j) $.
Put 
$g' := h' -
\sum _{j=1} ^{t}
a _j g _j \in T L$
and 
$g := T ^{i-1} g' \in T ^i L$.
We have $f -g \in K$.
Moreover, using \ref{ntnL(a,b)}.\ref{eq-L(a,b)Tnbis}, we get
$a _j g _j \in L ( a,d +(i-1)c)\subset L ( c, d  +(i-1)c)$.
Hence, $g' \in L ( c,d +(i-1)c)$
and then $g = T ^{i-1} g' \in L ( c,d +(i-1)c -(i-1)c)
=
L ( c,d)$.
Hence, $g$ satisfies every desired properties. 
\end{proof}

\begin{lem}
\label{lem2-thmA}
Let $r$ be an integer and 
$L:= (D [[ T ]] ^{\mathrm{oc}}) ^r$.
Let 
$0\to K \to L \overset{\phi}{\to} M \to 0$ be an exact sequence of 
left (resp. right) $D [[ T]] ^{\mathrm{oc}}$-modules such that $M$ has no 
$T$-torsion. Let $c, d \in \R _+$.
For any $j\in \bbN$,  
let $g _j \in (T ^j L +K)\cap L ( c, (j+1)d)$ (see notation \ref{ntnL(a,b)}).
Then $\sum _{j\in \bbN} \phi ( g _j)$ converges in $M$ for the $T$-adic topology.
\end{lem}

\begin{proof}
Since $L ( c, D) \subset L ( c +1, D)$ for any $c$ and $D$, we can suppose that $c \geq a$ where $a$ 
satifies the condition of Lemma \ref{lem1-ThmA}.
Hence, there exists 
$h _j \in T ^j L \cap L ( c,  (j+1)d)$ such that $g _j - h_j \in K$, i.e.
such that 
$\phi ( h _j) = \phi ( g _j) $.
Set 
$h := \sum _{j\in \bbN} h _j \in (D [[ T ]] ) ^r$.
For any integer $n\in \bbN$, 
since 
$h ^{(n)}
= 
 \sum _{j\leq n} h _j ^{(n)}$, 
 since 
$\ord  h _j ^{(n)} \leq cn + (j+1)d$, 
then 
$\ord h ^{(n)}
\leq 
(c+d)n + d$. 
Hence, 
$h \in L$.
We have 
$\phi (h) 
= 
\phi (\sum _{j\in \bbN} h _j )
= 
\sum _{j\in \bbN} \phi ( g _j) $, 
where the equality holds a priori
in $D [[ T]]  \otimes _{D [[ T]] ^{\mathrm{oc}}} M$, 
the $T$-adic completion of $M$.
Since $\phi (h ) \in M$, then we are done.
\end{proof}

\begin{prop}
\label{prop-thmA}
Let $M$ be a left $D [[ T]] ^{\mathrm{oc}}$-module of finite type. 
Let $c,d\in \R _+$,
 $\mu _1, \dots, \mu _r \in M$.
For any $j\in \bbN$ and $i=1,\dots, r$ , let $P _{i,j} \in D [[ T]] ^{\mathrm{oc}} (c,  (j+1)d)$.
If, for any $j\in \bbN$, we have  $\sum _{i=1} ^{r} P _{i,j} \mu _i \in T ^{j} M$ 
then 
 $\sum _{j=0} ^{\infty} \sum _{i=1} ^{r} P _{i,j} \mu _i \in M$, i.e. converges to an element of $M$ for the $T$-adic topology.
\end{prop}

\begin{proof}
1) First suppose $M$ has no $T$-torsion. 
Put $L:= (D [[ T]] ^{\mathrm{oc}}) ^r$. 
Let $\phi \colon L \to M$ be the morphism of left $D [[ T]] ^{\mathrm{oc}}$-module given
by $\mu _1, \dots, \mu _r \in M$. 
Let $N$ be the left $D [[ T]] ^{\mathrm{oc}}$-submodule of $M$ generated 
by $\mu _1, \dots, \mu _r \in M$.
By using \cite[3.2.3.(i)]{Be1}, 
since $D [[ T]] ^{\mathrm{oc}}$ is noetherian and since $T$ is in the center of $D [[ T]] ^{\mathrm{oc}}$, 
 there exists an integer $n _0$ such that 
$N \cap T ^{n _0 +j} M \subset T ^j N$ for any $j\in \bbN$.
Hence, 
replacing $M$ by $N$ if necessary, 
we can suppose $\phi$ surjective. 
Let $K := \ker \phi$. 
For any $j \in \bbN$, we remark that 
$g _j:= (P _{1j},\dots, P _{rj}) \in L ( c, (j+1) d)$,
$\phi (g _j) = \sum _{i=1} ^{r} P _{i,j} \mu _i\in T ^{j} M$ and then
$g _j \in (T ^j L + K)$.
Using \ref{lem2-thmA}, we get that 
$\sum _{j=0} ^{\infty} \sum _{i=1} ^{r} P _{i,j} \mu _i \in M$.

2) Since $T$ is in the center of $D [[ T]] ^{\mathrm{oc}}$, we remark that the subset $M'$ of $M$
of $T$-torsion elements is a left $D [[ T]] ^{\mathrm{oc}}$-module. 
Using \ref{oc-noeth}, 
$M'$ is a left $D [[ T]] ^{\mathrm{oc}}$-module of finite type. 
Since $T$ is in the center of $D [[ T]] ^{\mathrm{oc}}$, this yields that for $e\in \bbN$ large enough
$T ^e M' = 0$.
This yields $T ^e M$ has no $T$-torsion. 
Hence, the sequence 
 $\sum _{j=e} ^{n} \sum _{i=1} ^{r} P _{i,j} \mu _i \in T ^e M$
 converges to the element $\sum _{j=e} ^{\infty} \sum _{i=1} ^{r} P _{i,j} \mu _i \in T ^e M$.
This yields that 
$\sum _{j=0} ^{\infty} \sum _{i=1} ^{r} P _{i,j} \mu _i \in M$.
\end{proof}

\begin{ntn} \label{ntn-Dndag}
We suppose moreover that $D$ is a $\cV$-algebra.  We denote by $\widehat{D}$ the $\pi$-adic completion of $D$.  Since $\pi$ is in the center of $D$, 
we get a morphism of $\cV$-algebras of the form  $D [[T ]] ^{\mathrm{oc}} \to \widehat{D}$ which sends $\sum _{d} P _d T ^{d}$ to $\sum _{d} P _d \pi ^{d}$, where $P _d \in D$.
We denote by  $(D , (D _n) _{n\in \bbN}) ^\dag$  or simply by $(D , D _n) ^\dag$ the image of this morphism.  From Lemma \ref{oc-noeth}, $(D , (D _n) _{n\in \bbN}) ^\dag$ is (right and left) noetherian.
Moreover,  for any $a,b\in \R _+$, we denote by  $(D , D _n) ^\dag (a,b)$ the image of  $D [[T ]] ^{\mathrm{oc}} (a,b)$ (see Notation \ref{ntnL(a,b)}) via the homomorphism  $D [[T ]] ^{\mathrm{oc}} \to \widehat{D}$.
\end{ntn}

\begin{coro} \label{convergencepre}
We keep Notation \ref{ntn-Dndag}. Let $M$ be a left $(D , D _n) ^\dag$-module of finite type.  Let  $c,d\in \R _+$,  $\mu _1, \dots, \mu _r \in M$. Let    $Q _{i,j} \in (D , D _n) ^\dag (c,  (j+1)d)$ for any $j\in \bbN$ and $i=1,\dots, r$.
If, for any $j \in \bbN$,   $\sum _{i=1} ^{r} Q _{i,j} \mu _i \in p ^{j} M$,  then   $\sum _{j=0} ^{\infty} \sum _{i=1} ^{r} Q _{i,j} \mu _i \in M$, i.e. converges to an element of $M$ for the $p$-adic topology.
\end{coro}

\begin{proof}
This is a straightforward consequence of \ref{prop-thmA}. 
\end{proof}

 \subsubsection{Definition of the non commutative weak completion and first properties}
Let $R$ be a commutative $\cV$-algebra which is separated for the $p$-adic topology.
For commodity, we introduce the following non standard definition:
\begin{dfn} [Category of $R$-rings]
\begin{enumerate}[(a)]
\item 
An $R$-ring is the data of a $\cV$-algebra (not necessarily commutative) $D$ which is separated for the $p$-adic topology  and of a {\it monomorphism} of $\cV$-algebras of the form  $R\to D$. 

\item Let $D$ and $D'$ be two $R$-rings.  A morphism of $\cV$-algebras of the form $f\colon D\to D'$ is a morphism  of $R$-rings if it commutes with the structural morphisms $R \to D$ and $R \to D'$. When $f$ is injective, we say that $D$ is a sub-$R$-ring of $D'$. 
\end{enumerate}
\end{dfn}

\begin{rem}
By convention, for us if $D$ is an $R$-ring, then $R$ is a sub $R$-ring of $D$.
\end{rem}

\begin{empt}
Let $D$ be an $R$-ring. 
Let $r\in \bbN$. Let $d _1, \dots, d _r$ be $r$-elements of $D\setminus R$ (when $r=0$, this is the empty set), 
where $R$ is viewed as a subset of $D$ via the monomorphism $R \to D$. 
We denote by $\mathfrak{W} (d _1,\dots, d _r/R)$ the set of words over the alphabet $R \cup \{d _1, \dots, d _r \}$
(when $r=0$, $\mathfrak{W} (d _1,\dots, d _r/R)$ is the set of words over the alphabet $R$).

We denote by $( \mathfrak{V} (d _1,\dots, d _r/R), +,\star)$ the free (associative, unital) $\cV$-algebra on the set 
$R \cup \{d _1, \dots, d _r \}$.
We recall that $( \mathfrak{V} (d _1,\dots, d _r/R), +,\star)$ is a free $\cV$-module
which has
$\mathfrak{W} (d _1,\dots, d _r/R)$ as a basis and
that the 
$\cV$-bilinear morphism 
$\star\colon \mathfrak{V} (d _1,\dots, d _r/R)\times \mathfrak{V} (d _1,\dots, d _r/R)
\to \mathfrak{V} (d _1,\dots, d _r/R)$ is defined by concatenation.
By abuse of notation,  the image of an element $P \in \mathfrak{W} (d _1,\dots, d _r/R)$
in $\mathfrak{V} (d _1,\dots, d _r/R)$
is still denoted by $P$.

The canonical morphism of $\cV$-algebras of the form
$\mathfrak{V} (d _1,\dots, d _r/R) 
\to D$,
which is by definition the identity on the letters of $R \cup \{d _1, \dots, d _r \}$, 
will be denoted by 
$P \mapsto P (d _1,\dots, d _r) $.

\end{empt}

\begin{ntn}
Let $D$ be an $R$-ring. 
Let $r\in \bbN$ and let $d _1, \dots, d _r$ be $r$ elements of $D\setminus R$.

\begin{enumerate}[(a)]

\item For any $i\in \{ 1,\dots, r\}$,
we define 
the degree $\deg _{d _i}$ with respect to $d _i$ in $\mathfrak{W} (d _1,\dots, d _r/R)$ as follows
\begin{enumerate}[(a)]

 \item $\deg _{d _i} (0)=-\infty$ ;
 \item $\deg _{d _i} (a)=0$ if $a \in R\setminus  \{ 0\}$;
 \item $\deg _{d _i} (d _i)=1$ ;
  \item $\deg _{d _i} (P\star Q)=
 \deg _{d _i} (P)+
 \deg _{d _i} (Q)$ for any $P,Q\in  \mathfrak{W} (d _1,\dots, d _r/R)$.
\end{enumerate}

\item We define  the degree  in $\mathfrak{W} (d _1,\dots, d _r/R)$  by setting  $\deg (P) := \sum _{i =1} ^{r} \deg _{d _i} (P)$  for any $P \in \mathfrak{W} (d _1,\dots, d _r/R)$.
Remark that $\deg$ might depend on $R$  but does not depend on $d _1,\dots, d _r$.

\item We get the function
$\deg\colon
\mathfrak{V} (d _1,\dots, d _r/R)
\to \bbN \cup \{-\infty\}$ 
(resp. $\deg _{d _i}\colon
\mathfrak{V} (d _1,\dots, d _r/R)
\to \bbN \cup \{-\infty\}$) 
defined 
by sending $0$ to $-\infty$ and 
a finite sum of the form $\sum _{P} a _P P$ (with $a _P\in \cV \setminus \{0\}$ and $P\in \mathfrak{W} (d _1,\dots, d _r/R)$)
to the maximal of $\deg (P)$ (resp. $\deg _{d _i} (P)$).

\end{enumerate}

\end{ntn}

\begin{empt}
\label{ineq-star}
Let $D$ be an $R$-ring. 
Let $d _1, \dots, d _r$ be some elements of $D\setminus R$.
We check 
$\deg _{d _i} (P\star Q)
\leq
 \deg _{d _i} (P)+
 \deg _{d _i} (Q)$
 and 
$\deg (P\star Q)
\leq
 \deg (P)+
 \deg (Q)$
 for any $P,Q\in  \mathfrak{V} (d _1,\dots, d _r)$.
\end{empt}

\begin{empt}
\label{substitution}
Let $D$ be an $R$-ring. 
Let $z _1, \dots, z _r$ be some elements of $D\setminus R$,
let $P\in \mathfrak{W} (z _1,\dots, z _r)$.
Let $d _1, \dots, d _s\in D\setminus R$.
For any $\alpha =1,\dots, r$, let 
$\underline{Q} _\alpha 
=
(Q _{\alpha,1}, \dots, Q _{\alpha, n _\alpha})
\in \mathfrak{V} (d _1,\dots, d _s) ^{n _\alpha}$, 
where 
$n _\alpha:= \deg _{z _\alpha}P$. 
We denote by 
$P (\underline{Q} _1,\dots, \underline{Q} _r)$ 
the element of  $\mathfrak{V} (d _1,\dots, d _s)$ 
where the $\beta$th term  $z _\alpha$ in $P$ (reading from the left to the right)  is substitued by $Q _{\alpha,\beta}$.
Using \ref{ineq-star}, we compute 
\begin{gather}
\notag
\deg _{d_i} P (\underline{Q} _1,\dots, \underline{Q} _r)
 \leq  
 \sum _{\alpha =1} ^r
 \sum _{\beta =1} ^{n _\alpha}  
 \deg _{d_i} (Q _{\alpha,\beta}),\\
\label{deg-comp}
\deg P (\underline{Q} _1,\dots, \underline{Q} _r) 
\leq  
\sum _{\alpha =1} ^r 
 \sum _{\beta =1} ^{n _\alpha}
\deg (Q _{\alpha,\beta}).
\end{gather}
\end{empt}

\begin{lem}
Let $D$ be an $R$-ring,
 $\widehat{D}:= \underleftarrow{\lim}_ n \, D / \pi ^{n+1}D$ the associated $p$-adically separated completed $\widehat{R}$-ring. 
Let $d _1, \dots, d _s$ be some elements of $D\setminus R$.
Let $z _1, \dots, z _r$ be some elements of $D\setminus R$ 
which can be written in the form
$z _\alpha = \sum _{k _\alpha=0} ^{\infty} \pi ^{k _\alpha}Q _{\alpha, k _\alpha} (d _1,\dots, d _s)$,
where $Q _{\alpha, k _\alpha} \in  \mathfrak{V} (d _1,\dots, d _s/R)$ and $\alpha = 1,\dots, r$.
Let $P\in \mathfrak{W} (z _1,\dots, z _r/R)$. 
For any $\alpha = 1,\dots, r$, 
for any $\underline{k} _{\alpha} \in \bbN ^{n _\alpha}$
with $n _\alpha:= \deg _{z _\alpha}P$,
put 
$|\underline{k} _{\alpha}   |
=
\sum _{\beta =1} ^{n _\alpha}
k _{\alpha, \beta}$, 
put
$Q _{\alpha, \underline{k} _\alpha}
:=
(Q _{\alpha, k  _{\alpha,1}},
\dots,
Q _{\alpha, k _{\alpha,n _\alpha}})
\in \mathfrak{V} (d _1,\dots, d _s/R) ^{n _\alpha}$.
Then we have the equality in $\widehat{D}$ (a priori, the right term is an element of $\widehat{D}$,
and recall $D$ is included in $\widehat{D}$ so the left term is also an  element of $\widehat{D}$):
\begin{equation}
\label{computation-compo-pi}
P(z _1,\dots, z _r) 
=
\sum _{\underline{k} _{1} \in \bbN ^{n _1} ,
\dots 
,
\underline{k} _{r} \in \bbN ^{n _r}}
\pi ^{|\underline{k} _1| +\dots+|\underline{k} _r| } 
P ( Q _{1, \underline{k} _1}, \dots, Q _{r, \underline{k} _r})
(d _1,\dots, d _s)
\end{equation}
\end{lem}

\begin{proof}
Put 
$S _{\alpha, \underline{k} _\alpha}
:=
(\pi ^{k  _{\alpha,1}} Q _{\alpha, k  _{\alpha,1}},
\dots,
\pi ^{k  _{\alpha, n _\alpha}} Q _{\alpha, k _{\alpha,n _\alpha}})
\in \mathfrak{V} (d _1,\dots, d _s/R) ^{n _\alpha}$.
By $\cV$-multilinearity, we get the formula in 
$\mathfrak{V} (d _1,\dots, d _s/R) $ :
\begin{equation}
\notag
\label{computation-compo-pi-pre}
P ( S _{1, \underline{k} _1}, \dots, S _{r, \underline{k} _r})
=
\pi ^{|\underline{k} _1| +\dots+|\underline{k} _r| } 
P ( Q _{1, \underline{k} _1}, \dots, Q _{r, \underline{k} _r}).
\end{equation}
By additivity (and convergence of the sums in $\widehat{D}$),
we get the equality
$$P(z _1,\dots, z _r) =
\sum _{\underline{k} _{1} \in \bbN ^{n _1} ,
\dots 
,
\underline{k} _{r} \in \bbN ^{n _r}}
P ( S _{1, \underline{k} _1}, \dots, S _{r, \underline{k} _r})(d _1,\dots, d _s).$$
Hence, we conclude.
\end{proof}

\begin{dfn}\label{dfn-wc}
Let $D$ be an $R$-ring,  $\widehat{D}:= \underleftarrow{\lim}_ n \, D / \pi ^{n+1}D$ the associated $p$-adically separated completed $\widehat{R}$-ring. 
The weak completion $D ^\dag$ of $D$ as $R$-ring is by definition the subset of  the $\widehat{R}$-ring $\widehat{D}$ consisting of elements $P$ having the representations
\begin{equation} \label{dfn-wc=pre}
P= \sum _{i=0} ^{\infty} \pi ^{i} P _i (d _1,\dots, d _r)
\end{equation}
with $r\in \bbN$, $d _1,\dots, d _r$ are $r$ elements of $D\setminus R$ (when $r=0$, this is empty set),   $P _i \in  \mathfrak{V} (d _1,\dots, d _r/R)$ are such that there exists a constant $c >0$  so that   $\deg (P _i) \leq c (i+1)$ for any $i$.
We remark that the weak completion does depend on the ring $R$ so ``as $R$-ring'' is important to mention.
When we consider different rings $R$ at the same time for the same $\cV$-algebra $D$,  to be more precise we might write  $D ^{\dag _{/R}}$ to clarify that this is the weak completion as $R$-ring.
 
An $R$-ring $D$ is weakly complete (as $R$-ring to be precise) if the natural morphism  $D \to D ^\dag$ is an isomorphism. 
\end{dfn}

\begin{dfn} \label{(c,d)bounded}
Let $D$ be an $R$-ring. Let $(c,d) \in \R _+  \times \R$.  An element $P$ of $\widehat{D}$ is said to be $(c,d)$-bounded (as $R$-ring) if we can write $P$ of the form
\begin{equation} \label{dfn-wc=}
P= \sum _{i=0} ^{\infty} \pi ^{i} P _i (d _1,\dots, d _r)
\end{equation}
 where $d _1,\dots, d _r \in D \setminus R$,   $P _i \in  \mathfrak{V} (d _1,\dots, d _r/R)$ are such that  $\deg (P _i) \leq c i+d$ for any $i$.
In this case such a $P$ is said to be $(c,d)$-bounded as $R$-ring with respect to $d _1,\dots, d _r$. 

\end{dfn}

\begin{empt}\label{(c,d)bounded-ring}
Let $D$ be an $R$-ring. Let $P,P' \in D ^\dag$. Let $c,c'\in\bbR _+$, $d,d'\in \bbR$. If   $P $ is $(c , d)$-bounded
and  $P '$ is $(c', d ')$-bounded  then   $P+ P '$ is $(c '', d '')$-bounded  and $P  P '$ is $(c '', d  + d ')$-bounded,  where $c ''= \max \{ c , c '\}$  and   $d ''= \max \{ d , d '\}$. In particular, this yields that  $D ^\dag$ is a $R$-ring.
\end{empt}

\begin{ex}
\label{ex-comp-wcomp}
Let $D$ be an $R$-ring. Then $\widehat{D}:= \underleftarrow{\lim}_ n \, D / \pi ^{n+1}D$ is weakly complete since it is separated for the $p$-adic topology
and since every sum of the form \ref{dfn-wc=} satisfying the  condition of the definition converges. 
\end{ex}

\begin{dfn}
Let $S$ be a subset of a weakly complete $R$-ring $D'$. We denote by $S ^\dag$ the subset of $D'$ which consists
of all infinite sums of the type  \ref{dfn-wc=} subject to the restriction that each $d _i$ lie in $S\setminus R$. 
When $S \subset R$, $S ^\dag = R$. In general, we get  $R \subset S ^\dag$ (take $r= 0$).
We say that $S ^\dag$  is {\it the weakly generated $R$-ring by $S$ in $D'$}. This is justified by the next theorem. 
\end{dfn}

\begin{thm}
\label{thmSdag}
Let $S$ be a subset of a weakly complete $R$-ring $D'$.
Then $S ^{\dag}$ is the smallest weakly complete sub-$R$-ring of $D'$ containing $S$. 
\end{thm}

\begin{proof}
This is obvious by definition that $S ^\dag$ is a sub-$R$-ring of $ D'$ and $S ^\dag$ is separated and contains $S$.  We remark that if $T$ is a subset of $D'$ which contains $S$ then $T ^\dag $ contains $S ^\dag$.
Hence, we get the  inclusion $S ^\dag \subset (S ^\dag ) ^\dag$. Now, let us check the inclusion  $(S ^\dag ) ^\dag \subset  S ^\dag $. Let $P \in (S ^\dag ) ^\dag$. Since $R \subset S ^\dag$, we reduce to suppose that  $P \not \in R$.
We can write $P= \sum _{i=0} ^{\infty} \pi ^{i}P _i (z _1,\dots, z _r)$  where $z _1,\dots, z _r \in S ^{\dag}\setminus R$,   $P _i \in  \mathfrak{V} (z _1,\dots, z _r/R)$ and for some constant $c >0$ we have $\deg(P _i) \leq c (i+1)$ for any $i\in\bbN$.
By definition of $S ^\dag$, increasing $c$ if necessary, we can write $z _\alpha = \sum _{k _\alpha=0} ^{\infty} \pi ^{k _\alpha}Q _{\alpha, k _\alpha} (d _1,\dots, d _s)$  where $d _1,\dots, d _s \in S\setminus R$,   
$Q _{\alpha, k _\alpha}\in  \mathfrak{V} (d _1,\dots, d _s/R)$ is such that   $\deg (Q _{\alpha, k _\alpha}) \leq c (k _\alpha+1)$.
 By definition, for any $i\in \bbN$, there exists $N _i \geq 0$, such that   we can write $P _i$ as a finite sum of the form  $P _i = \sum _{j= 0} ^{N _i}a _{i,j} P  _{i,j}$ 
where  $a _{i,j} \in \cV \setminus\{0\}$ and $P_{i,j} \in  \mathfrak{W} (z _1,\dots, z _r/R)$  are pairwise distinct for any $0\leq j \leq N _i$.
 We set   $n _{ij,\alpha}:= \deg _{z _\alpha} (P  _{i,j})$ for any $\alpha=1,\dots, r$  and   $n _{ij}:=\deg ( P _{ij}) = \sum _{\alpha =1} ^{r} n _{ij,\alpha} \leq \deg ( P _{i}) \leq c (i+1)$.
 Using the formula \ref{computation-compo-pi}, we compute: 
\begin{gather} \notag
P  =  \sum _{i=0} ^{\infty} \sum _{j= 0} ^{N _i} \pi ^{i} a _{i,j} P  _{i,j}  (z _1,\dots, z _r)
\\ \label{2.1.11.1} \overset{\ref{computation-compo-pi}}{=} \sum _{i=0} ^{\infty} \sum _{j= 0} ^{N _i} \sum _{\underline{k} _{ij,1} \in \bbN ^{n _{ij,1}}}  \dots  \sum _{   \underline{k} _{ij,r} \in \bbN ^{n _{ij,r}}} 
\pi ^{i}  \pi ^{|\underline{k} _{ij,1}| +\dots+|\underline{k} _{ij,r}| }  a _{i,j} P _{i,j} ( Q _{1, \underline{k} _{ij,1}}, \dots, Q _{r, \underline{k} _{ij,r}}) (d _1,\dots, d _s).
\end{gather}
For any fixed non negative integers $i$ and $k$, we set
$$E _{ik}:= \{ (j, \underline{k} _{ij,1}, \dots, \underline{k} _{ij,r})\in \bbN ^{n _{ij}+1} \text{such that } 0\leq j \leq N _i, \underline{k} _{ij,1} \in \bbN ^{n _{ij,1}} , \dots , \underline{k} _{ij,r} \in \bbN ^{n _{ij,r}},  |\underline{k} _{ij,1}| +\dots+|\underline{k} _{ij,r}| =k\}.$$
We set  $R _{i,k}:= \sum _{(j, \underline{k} _{ij,1}, \dots, \underline{k} _{ij,r}) \in E _{ik}} a _{i,j}P _{i,j} ( Q _{1, \underline{k} _{ij,1}}, \dots, Q _{r, \underline{k} _{ij,r}})$.
Using the equality \ref{2.1.11.1}, we get $P = \sum _{i=0} ^{\infty}  \sum _{k=0} ^{\infty} \pi ^{i+k}  R _{i,k}(d _1,\dots, d _s)$.
From \ref{deg-comp}, for any  $(j, \underline{k} _{ij,1}, \dots, \underline{k} _{ij,r}) \in E _{ik}$,  we get the first inequality
\begin{gather*}
\deg  P _{i,j} ( Q _{1, \underline{k} _{ij,1}}, \dots, Q _{r, \underline{k} _{ij,r}})  \leq   \sum _{\alpha =1} ^r   \sum _{\beta =1} ^{n _{ij,\alpha}} \deg (Q _{\alpha, k _{ij,\alpha, \beta}}) \leq   \sum _{\alpha =1} ^r   \sum _{\beta =1} ^{n _{ij,\alpha}}  c ( k _{ij,\alpha, \beta} +1))
 \\ =   c ( k + n _{ij})   \leq  c ( k + c (i+1))  \leq c (1+c) ( k+i+1)
\end{gather*}
Hence, 
$\deg ( R _{i,k})
 \leq c ' ( k+i+1)$, where $c ':=c (1+c)$.
By putting $R _l:= \sum _{l=i+k} R _{i,k}$,
we get  
$P =\sum _{l=0} ^{\infty} \pi ^{l}R _l (d _1,\dots, d _s)$
with 
$\deg ( R _{l})
\leq c ' (l +1)$, which implies that $P \in S ^\dag$.

\end{proof}

\begin{thm}
\label{DdagWC}
Let $D$ be an $R$-ring. Then the weak completion $D ^\dag$  of $D$ is weakly complete.
\end{thm}

\begin{proof}
We apply the theorem \ref{thmSdag} to the case where $S$ is equal to  $D$ and $D'$ is $\widehat{D}$.
\end{proof}

\begin{coro}
\label{coro-Sdag=Ddag}
Let $D$ be an $R$-ring.
Let $S =\{d _i \ ;\ i\in I \} $ be a subset of $D\setminus R$.
We suppose that for any $d \in D$, there exist a finite subset $\{ i _1, \dots, i _r\}$ of $I$, 
$P \in \mathfrak{V} (d _{i _1},\dots, d _{i _r})$ such that 
$P (d _{i _1},\dots, d _{i _r})=d$. 
Then $S ^\dag =D ^\dag$.
\end{coro}

\begin{proof}
By hypothesis, $D \subset S ^\dag$.
By applying Theorems \ref{thmSdag} and \ref{DdagWC}, 
we conclude.
\end{proof}

\begin{dfn} \label{dfn-wcft}
Let $D$ be an $R$-ring. We say that $D$ is a weakly complete of finite type $R$-ring (is an $R$-wcft ring for short)  if $D$ is weakly complete as $R$-ring and if there exists 
a finite subset  $S$ of $D\setminus R$ such that $S ^\dag =D $.
\end{dfn}

\begin{prop}
The weak completion as $R$-ring functor  from the category of $R$-rings to the category of weakly complete $R$-rings is the left adjoint functor of the canonical forgetful functor. 
\end{prop}

\begin{proof}
Let $\phi \colon D \to E$ be a morphism of $R$-rings so that $E$ is weakly complete.  We have to check that $\phi$ has a unique factorization  $\phi ^\dag \colon D ^\dag \to E$. 
The uniqueness is obvious. For the existence, it is enough to check that the $p$-adic completion of $\phi$,  $\widehat{\phi} \colon \widehat{D}\to \widehat{E}$ send $D ^\dag$ to $E$ ($=E ^\dag$), which is easy.
\end{proof}

\begin{coro} \label{coro-indu-lim}
The weak completion as $R$-ring functor commutes with inductive limits. 
\end{coro}

\begin{dfn}
A filtered $R$-ring $(D , D _n )$ is the data of  a non negative exhausted filtered ring $(D , (D _n) _{n\in \bbN} )$ such that $D$ is separated for the $p$-adic topology  and of a structure of $R$-ring on $D _0$ (in particular $D$ is an $R$-ring).
\end{dfn}

\begin{prop}
\label{comp-Ddag-DDndag}
Let $(D , D _n )$ be a filtered $R$-ring. 
\begin{enumerate}[(a)]
\item We get the inclusion 
$D ^\dag \subset  (D , D _n) ^\dag$, 
where $D ^\dag$  is the weak completion of $D$ as $R$-ring (see \ref{dfn-wc})
and $(D , D _n) ^\dag$ is defined in \ref{ntn-Dndag}.

\item Suppose there exist
$d _1,\dots, d _r \in D\setminus R$, $c,d\in \R _+$ satisfying the following property : 
for any $n \in \bbN$, for any $x \in D _n$ 
there exists  $P \in  \mathfrak{V} (d _1,\dots, d _r/R)$ such that 
$x= P (d _1,\dots, d _r)$ and 
  $\deg (P) \leq c n+d$. Putting $S := \{ d _1,\dots, d _r \}$, 
  we have $(D , D _n) ^\dag =S^\dag =D ^\dag$.

\end{enumerate}
\end{prop}

\begin{proof}
1) Let 
$x \in D ^\dag$ having the representation
$x= \sum _{i=0} ^{\infty} \pi ^{i} P _i (d _1,\dots, d _r)$
 where $d _1,\dots, d _r \in D\setminus R$, 
 $P _i \in  \mathfrak{V} (d _1,\dots, d _r/R)$ are such that there exists a constant $c >0$
 so that 
 $\deg (P _i) \leq c (i+1)$ for any $i$.
Let $n _0\in \bbN$ large enough such that 
$d _1,\dots, d _r \in D _{n _0}$.
Since we have also $R \subset D _0$, then 
$\ord P _i (d _1,\dots, d _r) \leq n _0 \deg (P _i) 
\leq n _0 c (i+1)$.
Hence, 
$ \sum _{i=0} ^{\infty} T ^{i} P _i (d _1,\dots, d _r) 
\in D [[T ]] ^{\mathrm{oc}}$.
This yields that $x \in (D , D _n) ^\dag$.

2) Let $ \sum _{i=0} ^{\infty} T ^{i} x _i
\in D [[T ]] ^{\mathrm{oc}}$, where $x _i \in D$ and there exist  $a,b\in \R _+$ such that  $x _i \in D _{ai+b}$.
Hence, there exists  $P _i \in  \mathfrak{V} (d _1,\dots, d _r/R)$ such that 
$x _i = P _i  (d _1,\dots, d _r)$ and 
  $\deg (P _i ) \leq c (ai+b)+d=cai +(cb +d)$. 
This yields that 
$ \sum _{i=0} ^{\infty} \pi ^{i} x _i \in S^\dag $.
\end{proof}

\begin{coro}
\label{cor-comp-Ddag-DDndag}
Let $D ^\dag$ be a weakly complete of finite type $R$-ring. 
Choose a subset $S =\{d _1, \dots, d _r \} $ of $D ^\dag \setminus R$ such that $S ^\dag =D ^\dag$.
We denote by $D$ the $R$-subring of $D ^\dag$ generated by $S$. 
This $R$-ring $D$ is in fact a filtered $R$-ring : 
for any $n\in \bbN$, set 
$$D _n:= \{d \in D \; |\; 
\text{$\exists P \in \mathfrak{V} (d _{1},\dots, d _{r})$ such that 
$P (d _{1},\dots, d _{r})=d$ and  $\deg (P) \leq n$.}\} .$$
Then $(D , D _n) ^\dag =D ^\dag$.
\end{coro}

\begin{proof}
This is the corollary of \ref{comp-Ddag-DDndag}.2 in the special case where $c=1$ and $d=0$.
\end{proof}

\begin{empt} \label{ex-B-wc}
For example, if $A$ is a commutative $R$-algebra and $t \in A$, we have  the surjection  $A [T] \to A _t$. From the filtration by the degree of $A [T]$, we get an induced filtration $(F _n A _t) _{n\in \bbN}$ on  $A _t$. 
The weak completion $A _{[t]}$ of $A _t$ as $A$-algebra is equal to  $( A _t, (F _n A _t) _{n\in \bbN}) ^\dag$. Hence, we can apply \ref{convergencepre} to $A _{[t]}$.
\end{empt}

\subsubsection{Weak completion of differential operators of level $m$:local description}

Let $\fX ^\sharp$ be a $p$-torsion-free affine  log smooth formal $\fS ^\sharp$-scheme.  
Put $D ^{(m)}:= \Gamma (\fX, \cD ^{(m)} _{\fX ^\sharp/\fS ^{\sharp}})$, $D := \Gamma (\fX, \cD _{\fX ^\sharp/\fS ^{\sharp}})$ and $A:= \Gamma (\fX, \cO _{\fX})$. 
Let $D ^{\dag}$ and $D ^{(m)\dag}$ be the weak completion of   $D $ and  $D ^{(m)}$ as $A$-ring. 

\begin{empt}
Let $A [X _1,\dots, X _r] _{(m)}$ be the commutative $A$-algebra defined as follows :  $A [X _1,\dots, X _r] _{(m)}$ is free as $A$-module and has a basis  
$\underline{X} ^{<\underline{k}> _{(m)}}$ with $\underline{k}\in \bbN ^r$ such that  $\underline{X} ^{<\underline{k}> _{(m)}} \cdot\underline{X} ^{<\underline{k'}> _{(m)}} = 
\left< \begin{smallmatrix}  \underline{k}+\underline{k'}  \\   \underline{k} \end{smallmatrix}\right> _{(m)} \underline{X} ^{<\underline{k}+\underline{k'}> _{(m)}}$.
We denote by $(k )_i:=(0,\dots, 0,k,0,\dots, 0)$ where the $k$ is at the $i$th place, we put  $X _i ^{<k> _{(m)}}:= \underline{X} ^{<(k)_i> _{(m)}}$.  We compute  $\underline{X} ^{<\underline{k}> _{(m)}} = \prod _i X _i ^{<k _i> _{(m)}}$.
Since $A$ is a noetherian commutative ring and since $A [X _1,\dots, X _r] _{(m)}$ is an $A$-algebra generated by  $X _i ^{<p ^{j}> _{(m)}}$ with $0 \leq j \leq m$ and $1\leq i \leq r$ (same computation than \cite[2.2.5.1]{Be1}), then 
the $A$-algebra $A [X _1,\dots, X _r] _{(m)}$ is noetherian.

Let $f= \sum _{\underline{k}} a _{\underline{k}} \underline{X} ^{<\underline{k}> _{(m)}}  \in A [X _1,\dots, X _r] _{(m)}$.  If $f$ is not zero, we define $\deg f$, the degree of $f$, to be  the greatest of $|\underline{k}|$ such that $a _{\underline{k}}\not =0$.
\end{empt}

\begin{empt} \label{ex-grD-A(m)}
Since $\fX ^\sharp$ is separated and quasi-compact,  the functor  $\Gamma (\fX ^\sharp , -)$ commutes with inductive limits.  Hence, setting  $D ^{(m)}_n := \Gamma (\fX, \cD ^{(m)} _{\fX ^\sharp/\fS ^{\sharp},n})$, we get 
$D ^{(m)}= \cup D ^{(m)}_n$. We get the filtered ring $(D, (D _n) _{n\in \bbN} )$  given by the order of the differential operators.  Using Proposition \cite[2.3.2.a)]{these_montagnon}, we check that $ \gr D ^{(m)}:= \oplus _{n\in \bbN} D ^{(m)} _{n}/ D ^{(m)} _{n-1}$ is a commutative $A$-algebra generated by 
$\oplus _{n=1} ^{p ^m}D ^{(m)} _{n}/ D ^{(m)} _{n-1}$. Since $A$ is noetherian and  $D ^{(m)} _{n}/ D ^{(m)} _{n-1}$ are finitely generated $A$-modules, this yields that  $ \gr D ^{(m)}$ is noetherian.  Hence, we are in the situation of \ref{ntn-ord-D-deg}.

In fact, when   $ \fX^\sharp$ has logarithmic coordinates $u _1,\dots, u _r$,   we have the following explicit description of $ \gr D ^{(m)}$.  Put $\xi _{i} ^{<k> _{(m)}}:= \sigma (\partial _{i} ^{<k> _{(m)}})$,
where $\sigma \colon D ^{(m)} \to \gr D ^{(m)}:= \oplus _{n\in \bbN} D ^{(m)}_{n}/ D ^{(m)}_{n-1}$  is the canonical mapping.  When $k=1$, we simply denote it by $\xi _i$. 
By Proposition \cite[2.3.2.a)]{these_montagnon}, $\gr D  ^{(m)}= A [\xi _1,\dots, \xi _r] _{(m)}$, which is noetherian.
\end{empt}

\begin{empt} [$p$-adic norm, $v _p$, $v _\pi$]
\label{padic-norm}
We denote by $e$ the absolute ramification index of $\cV$. Let $B$ be a $p$-torsion free commutative $\cV$-algebra 
We endow $B _K$ with a norm induced by the $p$-adic topology on $B$ as follows.
For any $b \in B _K$, we set 
$\| b \| = p ^{- v _p (b)}$, where $E _b =\{n \in \bbZ \; ; \; b \in p ^n  B \}$ and
$v _p (b):= \max E _b$ if $b \not = 0$ and $v _p (0)=+\infty $. 
This norm is called the $p$-adic norm on $B _K$ given by $B$.

Remark moreover that $\| b b' \| \leq \| b \| \ \| b '\| $ and then $v _p \colon  B _K \to \bbZ \cup \{+\infty\}$ is a quasi-valuation. We call $v _p$ to be the $p$-adic quasi-valuation of  $B _K$ induced by $B$.
We also define $v _\pi (b):= \max \{n \in \bbZ \; ; \; b \in \pi ^nB  \}$ if $b\not =0$ and $v _\pi (b)=+\infty$ otherwise.
We say that $v _\pi\colon  B _K \to \bbZ \cup \{+\infty\}$
the $\pi$-adic quasi-valuation of 
$B _K$ induced by $B$.
For any $b \in  B _K$, we have 
\begin{equation}
\label{vp-vpi}
v _p (b) \leq \frac{v _{\pi} (b)}{e } < v _p (b)+1.
\end{equation}
We denote by $| -|$ the $p$-adic norm on $\bbQ$ given by $\bbZ _{(p)}$.
Hence, by convention $| p | = p ^{-1}$. 
\end{empt}

\begin{empt}\label{eq-condDdag}
We suppose that $\fX ^\sharp/\fS ^{\sharp}$ has local logarithmic coordinates $u _1,\dots, u _d \in M _\fX$ i.e. 
we have a log étale morphism over $\fS ^{\sharp}$ of the form  $\fX ^{\sharp} \to \fS ^{\sharp} \times \mathfrak{A} _{\bbN ^d}$ given by $u _1,\dots, u _d \in M _\fX$. Let $\partial _{\sharp 1},\dots, \partial _{\sharp d}$ be the corresponding logarithmic derivations. 
The sheaf $\cD ^{(m)} _{\fX ^{\sharp} /\fS ^{\sharp}}$ is $\cO _{\fX}$-free with basis  $\{ \underline{\partial} _{\sharp} ^{<\underline{i}> _{(m)}},\, \underline{i} \in \bbN ^d\}$.

Let  $z$ be an element of $\widehat{D} ^{(m)}$,  the $p$-adic completion  of $D ^{(m)}$. For any $c \geq 0$, the following conditions are equivalent:

\begin{enumerate}[(a)]
\item We can write $z$ of  the form
\begin{equation} \label{wkcpdscplevelmpre0}
z = \sum _{\underline{i}\in \bbN ^d} a _{\underline{i}} \underline{\partial} _{\sharp} ^{<\underline{i}> _{(m)}}
\end{equation}
where $a _{\underline{i}} \in A$ satisfy $v _{p} ( a _{\underline{i}}) \geq \frac{|\underline{i}|}{c} -1$.

\item We can write $z$ of  the form
\begin{equation} \label{wkcpdscplevelmpre0bis}
z= \sum _{n=0} ^{\infty} p ^{n} z _n
\end{equation}
where $z _n \in D ^{(m)} _{c(n+1)}$.
\end{enumerate}
Remark we can also replace $p$ and $v _p$ by respectively $\pi$ and $v _{\pi}$.
With notation \ref{ntn-Dndag}, $z  \in (D ^{(m)}, (D ^{(m)} _n) _{n\in \bbN}) ^\dag$  if and only if there exists a constant $c>0$ such that   $z$ satisfies one the above equivalent conditions. 
\end{empt}

\begin{lem} \label{lem-dagDDm}
We suppose that $\fX ^\sharp$ has local logarithmic coordinates  and we keep notation \ref{eq-condDdag}.
 Let $S _{(m)}:=  \{\partial _{\sharp i} ^{[p ^j]}\ |\ i=1,\dots, d\; ;\; j= 0,\dots, m\}$ and $S := \cup _{m\in \bbN} S _{(m)}$.

\begin{enumerate}[(a)]
\item We have $S ^{\dag } = \cup _{m\in \bbN}S _{(m)} ^\dag$,  $S _{(m)} ^\dag = D ^{(m)\dag}$,  $  D ^{(m)\dag} \subset \widehat{D} ^{(m)} \subset D ^{(m+1)\dag}$, $S ^\dag = D ^{\dag}$  and  $D ^\dag = \Gamma (\fX, \cD _{\fX/\T} ^\dag)$.

\item We have the equality $(D ^{(m)}, (D ^{(m)} _n) _{n\in \bbN}) ^\dag = D ^{(m)\dag}$.

\item \label{noeth-weak completion-without}
The ring $D ^{(m)\dag}$ is (right and left) noetherian. 

\end{enumerate}
\end{lem}

\begin{proof}
a) The equality  $S _{(m)} ^\dag = D ^{(m)\dag}$ is a consequence of  \ref{coro-Sdag=Ddag} and of  the fact that  $S _{(m)} $ generates  $D ^{(m)}$ as $A$-ring (see \cite[2.3.1]{these_montagnon}).
Similarly we get $S ^\dag = D ^{\dag}$. The equality $S ^{\dag } = \cup _{m\in \bbN}S _{(m)} ^\dag$ is clear (e.g. see \ref{coro-indu-lim}).
 
Fix $m \in \bbN$.  The inclusion $D ^{(m)\dag} \subset \widehat{D} ^{(m)}$ is tautological. Let us prove now the inclusion $\widehat{D} ^{(m)} \subset D ^{(m+1)\dag}$.
For any $\underline{i}= (i _1,\dots, i _d)$, put $\lambda _{\underline{i}}:= \prod _{j=1} ^{d}\frac{[i _j/p ^{m}] !}{[i _j/p ^{m +1}]!}$. From \cite[2.4.3.1]{Be1}, we compute 
$v _p (\lambda _{\underline{i}})  \geq  \frac{|\underline{i}|}{ p ^{m+1} } -\frac{pd}{p-1} -\sum _{j=1} ^{d}\log _p (i _j +1) \geq  \frac{|\underline{i}|}{ p ^{m+1} } -\frac{pd}{p-1} -d \log _p (|\underline{i}| +1) $.
Hence, for $c >0$ large enough, we have  $v _p (\lambda _{\underline{i}})  \geq   \frac{|\underline{i}|}{c} -1$ (remark since $v _p (\lambda _{\underline{i}})  \geq 0$, then this later property has only to be satisfied for  $|\underline{i}| \geq c$).
An element $z$ in  $\widehat{D} ^{(m)}$  can be written in the form  $z=\sum _{\underline{i}\in \bbN ^d} a _{\underline{i}} \underline{\partial} _{\sharp} ^{<\underline{i}> _{(m)}}$,
where $a _{\underline{i}} \in A$ converges to $0$ when $|\underline{i}|$ goes to infinity.  We get $z=\sum _{\underline{i}\in \bbN ^d} b _{\underline{i}}   \underline{\partial} _{\sharp} ^{<\underline{i}> _{(m+1)}}$
where $b _{\underline{i}} = a _{\underline{i}} \lambda _{\underline{i}}$ are such that there exists a constant $c>0$ such that  $v _{p } ( b _{\underline{i}}) \geq \frac{|\underline{i}|}{c} -1$.
From (see \cite[2.3.1]{these_montagnon}), $\underline{\partial} _{\sharp} ^{<\underline{i}> _{(m+1)}}$ is a product of elements in  $S _{(m+1)}$ (remark also that the number of factors is less than $|\underline{i}|$).  Hence, $z \in D ^{(m+1)\dag}$.

Finally, since  $\Gamma (\fX, \cD _{\fX/\T} ^\dag)= \cup _{m\in\bbN} \widehat{D} ^{(m)}$,  since $D ^\dag = \cup _{m\in\bbN} D ^{(m)\dag}$, then we get  $D ^\dag = \Gamma (\fX, \cD _{\fX/\T} ^\dag)$.

b) Let us check the second part of the Lemma.  
Following \cite[2.3.1]{these_montagnon},  $\underline{\partial} _{\sharp} ^{<\underline{i}> _{(m)}}$ is a product of elements of $S _{(m)}$ so that the number of factors is less than $|\underline{i}|$.
Using \ref{comp-Ddag-DDndag}, this yields $(D ^{(m)}, (D ^{(m)} _n) _{n\in \bbN}) ^\dag = D ^{(m)\dag}$.

c) Since $(D ^{(m)}, (D ^{(m)} _n) _{n\in \bbN}) ^\dag = D ^{(m)\dag}$, then we get the surjective morphism of $A$-rings $D ^{(m)} [[T ]] ^{\mathrm{oc}} \to D ^{(m)\dag}$ which sends $\sum _{d} P ^{(d)} T ^{d}$ to $\sum _{d} P ^{(d)} \pi ^{d}$.
Hence, using \ref{oc-noeth}, we are done. 
\end{proof}

\subsection{Overconvergent singularities along the special fiber of weakly smooth $\fS$-formal schemes} 
\label{dfn-Sm-eta}   We denote by $\Sm ^\dag _\eta$ the category of $\fS$-formal schemes $\fP$ of finite type (see \cite[10.13.3]{EGAI})
and flat over $\fS $ whose generic fiber $\fP _\eta$ are smooth  over $\eta= \Spf \cV [[t]] \{ \frac{1}{t}\}$. Such objects are called weakly smooth $\fS$-formal schemes.

\subsubsection{The sheaf of overconvergent sections along the special fiber of weakly smooth $\fS$-formal schemes}

The goal of this subsection is to build the sheaf $\widetilde{\cO} _{\fP} $,  where $\fP$ is  a weakly smooth $\fS$-formal scheme.

\begin{ntn} \label{def-[t]}
Let $C$ be a noetherian commutative $\cV [[t]]$-algebra, complete with respect to the $p$-adic topology. We denote by $C _{[t]}$ the weak completion of $C _t$ as $C$-algebra.
\end{ntn}

\begin{empt}\label{prop-flat-[t]}
Since $C _{[t]}$ is a weakly complete finitely generated $C$-algebra, 
since $C$ is noetherian, then  $C _{[t]}$ is noetherian (see \cite{Fulton-noetherian})
and the homomorphism $C \to C _{[t]}$ is flat (see \cite[Lemma 2.1]{MonskyWashnitzer}).
Since $\pi C _{[t]}$ is included in the Jacobson radical of $C$ (see \cite[1.6]{MonskyWashnitzer}), then 
using \cite[8.14]{matsumura} we get the extension $C _{[t]} \to C _{\{t\}}$ is faithfully flat, 
where $C _{\{t\}}$ is the $p$-adic completion of $C _t$ (this is the notation of \cite{EGAI}).
\end{empt}

\begin{ntn} \label{ntnMtilde}
In this rest of the section, we will keep the following notation.  Let $\fX$ be an affine formal scheme of finite type and flat over $\fS $. 
Put $A := \Gamma (\fX , \cO _{\fX})$. and $\widetilde{A}:= A _{[t]}$.
Let $f\in A$.  As in \cite{EGAI}, we denote by $\mathfrak{D} (f):= \Spf A _{\{f\}}$ the open subset of $\fX$. A {\it principal open formal subscheme}  of $\fX$ is by definition a subset of the form $\mathfrak{D} (f)$.
Put $\widetilde{A} _{<f>}:= (A _{\{f\}})_{[t]}$ and  $\widetilde{A} _{\{f\}}$ the $p$-adic completion of $\widetilde{A} _f$ (or of $\widetilde{A} _{<f>}$, or of $A _{ft}$).  With \ref{def-[t]}, since $A$ is noetherian, then the extension $\widetilde{A} _{<f>} \to \widetilde{A} _{\{f\}}$ 
is  faithfully flat.  Since $\widetilde{A} \to \widetilde{A} _{\{f\}}$ is flat, this implies that the extension  $\widetilde{A}\to \widetilde{A} _{<f>} $ is flat.

Let $f,g\in A$ such that $\mathfrak{D} (f)\subset \mathfrak{D} (g)$. We have the  morphism of complete $A$-algebras of the form $A _{\{f\}} \to A _{\{g\}}$. 
Since $\widetilde{A} _{<g>} $ is also weakly complete as $A _{\{f\}} $-algebra,  we obtain the factorization $\widetilde{A} _{<f>}  \to \widetilde{A} _{<g>} $. 
\end{ntn}

\begin{ntn}
\label{functor-sim-ff exact}
We denote by $\mathfrak{B}$ the basis of open subsets of  $\fX$ consisting of
{\it principal} open subsets. 
If $M$ is a $\widetilde{A}$-module, 
we get a presheaf $M ^{\Delta}$ on $\mathfrak{B}$ 
defined by $\Gamma (\mathfrak{D} (f) , M ^{\Delta}) = \widetilde{A} _{<f>} \otimes _{\widetilde{A}} M$.
We get the functor
$\Delta\colon M \mapsto M ^{\Delta}$ from the category of 
$\widetilde{A}$-module to the category of presheaves of $\widetilde{A} ^{\Delta}$-modules on $\mathfrak{B}$.
Since $\widetilde{A}\to \widetilde{A} _{<f>} $ is flat, the functor $\Delta$ is exact. 
Moreover, $\Delta$ is fully faithful. Indeed a morphism of $\widetilde{A} ^{\Delta}$-modules
of the form $\phi \colon M ^{\Delta} \to N ^{\Delta}$ 
is equal to $\Gamma (\fX, \phi) ^\Delta$, i.e. 
the global section functor $\Gamma (\fX, -)$ gives a quasi-inverse functor from the essential image of 
$\Delta$ to the category of $\widetilde{A}$-modules. 
\end{ntn}

\begin{lem} \label{prelem1-lemma7}
Let $P \in \widetilde{A} _{<f>}$.  Let $(c,d) \in \R _+ ^2$ such that $P$ is $(c,d)$-bounded  as $A _{\{f\}}$-algebra with respect to $\frac{1}{t}$ (see \ref{(c,d)bounded}). 
Then  we can write  $P  =  \sum _{j=0} ^{+\infty}  \pi ^{j} \frac{a _{j}}{f ^{m _{j }}t ^{n _j}} $ where $a _j \in  A$, $(m _j) _{j\in \bbN}$ is an increasing sequence of integers and $(n _j) _{j\in \bbN}$ is an increasing sequence of integers such that $n _j \leq c j+d$. 
\end{lem}

\begin{proof}
By definition,  we can write  $P= \sum _{k=0} ^{+\infty} \pi ^{k} \frac{c _k}{ t ^{n _k}} $ where $c _k \in  A _{\{f\}}$ and $(n _k) _{k\in \bbN}$ is an increasing sequence of integers such that $n _k \leq c k+d$. We can write $c _k$ of the form 
$c _k = \sum _{i=0} ^{+\infty} \pi ^{i}  \frac{a _{k,i}}{ f ^{m _{k,i}}} $,  where $a _{k,i} \in A$ and $(m _{k,i}) _{i\in \bbN}$ is an increasing sequence of integers.
We get  $P =  \sum _{j=0} ^{+\infty}  \pi ^{j}  \sum _{k+i=j} \frac{a _{k,i}}{f ^{m _{k,i }}t ^{n _k}}  =  \sum _{j=0} ^{+\infty}  \pi ^{j} \frac{a _{j}}{f ^{m _{j }}t ^{n _j}} $ where  $m _j := \max \{ m _{k,i }\ ; \ k+i=j\}$ and  $ a _j :=  \sum _{k+i=j} a _{k,i}f ^{m _j -m _{k,i }}t ^{n _j -n _k}\in A$.
\end{proof}

\begin{empt}
The ring $A _f [T]$  of polynomials with coefficients in $A _f$ has a canonical filtrations $(A _f [T] _n) _{n\in \bbN}$ with $A _f [T] _n$ equal to the $A _f$-submodule of polynomials of degree $\leq n$.
Considering the mapping of $A _f$-algebras $\alpha \colon A _f [T] \to A _{f t}$ sending $T$ to $1/t$, we get a filtration on $A _{f t}$
by setting $A _{f t,n} := \alpha (A _f [T] _n)$.  Lemma \ref{prelem1-lemma7} can be translated with the equality $\widetilde{A} _{<f>}=(A _{f t} , A _{f t,n}) ^\dag $ (see notation \ref{ntn-Dndag}).
\end{empt}

\begin{empt}
[Cech cohomology] \label{cech-cohomology}
Let $\cM$ be a presheaf on $\mathfrak{B}$. Let $\fU$ be an element of $\mathfrak{B}$, i.e. a principal open subset of $\fX$. Let $\fU _\bullet = (\fU _0, \dots, \fU _m)$ be an open cover of $\fU$ by elements of $\mathfrak{B}$.
For any $0\leq i _0 < \cdots < i _r \leq m$, we put $\fU _{i _0,\dots, i _r} := \fU _{i _0} \cap \dots \cap\fU _{i _r} $. We set $C ^{r} (\fU _{\bullet}, \cM):= 0$  for $r \geq m+1$.
We set for any $0\leq r \leq m$
\begin{equation} \label{not1-cech-cohomology} 
C ^r (\fU _{\bullet}, \cM): = \underset{0\leq i _0 < \cdots < i _r \leq m}{\oplus} \Gamma (\fU _{i _0,\dots, i _r}, \cM).
\end{equation}
For any $0\leq r \leq m$, we have the coboundary map $\partial ^{r} \colon C ^{r} (\fU _{\bullet}, \cM) \to  C ^{r+1} (\fU _{\bullet}, \cM)$  by setting 
\begin{equation} \label{not2-cech-cohomology}
(\partial ^r \alpha ) _{i _0, \dots, i _{r+1}} =  \sum _{k= 0} ^{r+1}  (-1) ^{k} \alpha _{i _0 ,\dots, \widehat{i} _k,\dots, i _{r+1}}| \fU _{i _0, \dots, i _{r+1}}.
\end{equation}
If there is no risk of confusion, we will simply write $\partial$ this map. We set $H ^0 (\fU _{\bullet}, \cM) := \ker \partial ^{0}$. For any $1\leq r \leq m$, we set $H ^r (\fU _{\bullet}, \cM) := \ker \partial ^{r} /\mathrm{im} \partial ^{r-1}$.
\end{empt}

\begin{rem} \label{rem-CrotimesM}
Let $M$ be a $\widetilde{A}$-module.  Then we have the canonical isomorphism $C ^r (\fU _{\bullet}, \widetilde{A} ^\Delta) \otimes _{\widetilde{A}} M \riso C ^r (\fU _{\bullet}, M ^\Delta)$ which fits in the commutative diagram
\begin{equation} \notag
\xymatrix{     {C ^r (\fU _{\bullet}, \widetilde{A} ^\Delta) \otimes _{\widetilde{A}} M}  \ar[d] ^-{\partial ^{r}\otimes Id _M} \ar[r] ^-{\sim} &  {C ^r (\fU _{\bullet}, M ^\Delta)}   \ar[d] ^-{\partial ^{r}}
\\   {C ^{r+1} (\fU _{\bullet}, \widetilde{A} ^\Delta) \otimes _{\widetilde{A}} M}  \ar[r] ^-{\sim} &   {C ^{r+1} (\fU _{\bullet}, M ^\Delta).}     }
\end{equation}
\end{rem}

\begin{lem}  \label{piC(n)}
Let $C \to D$ be a monomorphism of commutative $\cV$-algebras such that $D$ has no $p$-torsion, $C /\pi C \to D /\pi D$ is injective.  Fix some elements $f, g \in D$, $n ,m \in \bbN$ and 
set $E := \{ b\in D \; ;\; f ^n g ^m b \in C\}$. For any integer $i \in \bbN$,   we have  the equality $E \cap \pi ^i D = \pi ^i E$.
\end{lem}

\begin{proof}
1) First, let us check  $C \cap \pi ^i D = \pi ^i C$ (this is the particular case where $n=m=0$). The inclusion  $\pi ^i C \subset C \cap \pi ^i D$ is obvious.  Conversely, we check by induction on $i \in \bbN$ that 
$C \cap \pi ^i D \subset  \pi ^i C$. When $i= 0$, this is clear. Suppose $i \geq 1$.  Let $a \in C \cap \pi ^i D$. There exists $b \in D$ such that $a  = \pi ^i b$.  Since $C /\pi C \to D /\pi D$ is injective, 
there exists $a ' \in C$ such that  $a = \pi a'$. Since $D$ has no $p$-torsion,  $a' = \pi ^{i-1} b$. Using the induction hypothesis, we get that $a' \in \pi ^{i-1}C$. Hence, $a \in \pi ^i C$.

2) Let us check the general case. 
The inclusion 
$\pi ^i E \subset E \cap \pi ^i D $ is obvious.
Conversely, let $c \in E \cap \pi ^i D $. 
There exists $b \in D$ such that 
$c  = \pi ^i b$. Moreover, 
$f ^n g ^m c \in C$. Hence, 
$f ^n g ^m c=  \pi ^i f ^n g ^m b \in 
C \cap \pi ^i D = \pi ^i C $. 
We get an element $a \in C$ such that
$f ^n g ^m c=  \pi ^i a$. 
Since $D$ has no $p$-torsion, this yields 
$a = f ^n g ^m b$ and in particular 
$b \in E$.
Hence, $c  = \pi ^i b\in \pi ^i E$.
\end{proof}

\begin{lem}
\label{lem-vanishingCech}
Let $\fU$ be an element of $\mathfrak{B}$, i.e. a principal open subset of $\fX$.
Let $\fU _\bullet = (\fU _0, \dots, \fU _m)$ be an open cover of $\fU$ by elements of $\mathfrak{B}$.
Suppose $(A/\pi A) _t$ is integral. Then, with the notation of \ref{cech-cohomology}, 
\begin{enumerate}[(a)]
\item the natural map $\Gamma (\fU, \widetilde{A} ^{\Delta}) \to H ^0 ( \fU _\bullet , \widetilde{A} ^{\Delta})$ is bijective ; 
\item $H ^{i} (\fU _\bullet , \widetilde{A} ^{\Delta})= 0$ for all $i >0$.
\end{enumerate}
\end{lem}

\begin{proof}
The proof is analogous to that of \cite[Lemma 7]{meredith-weakformalschemes}.
Since our context is slightly different, for the convenience of the reader, let us write up a proof. 
We can suppose $\fU= \fX$.  Let us fix some notation. For any integer $0\leq i\leq m$, let $f _i$ be such that $\mathfrak{D} (f _i) = \fU _i$.
Put $C ^{-1}:= \widetilde{A}$ (beware this is not the usual notation $C ^{-1}:=0$, 
but this will be convenient in the proof to avoid making two cases), $\widetilde{A} _{i _0,\dots, i _r}:= \widetilde{A} _{< f _{i _0} \cdots f _{i _{r}}>}$ and 
$$
C ^r:= C ^r (\fU _{\bullet}, \widetilde{A} ^{\Delta})= \underset{0\leq i _0 < \cdots < i _r \leq m}{\oplus} \widetilde{A} _{i _0,\dots, i _r}
$$
for $r \geq 0$. 
We have to check the exactness of the sequence
$0 \to C ^{-1} \to C ^{0} \overset{\partial}{\to}\cdots \overset{\partial}{\to} C ^{m}\to 0$.
Fix $r \geq 0$ an integer. 

Let $x\in C ^r$ be a cocycle (i.e. such that $\partial (x) = 0$) and $x _{i _0, \dots,i _r}$ be its component on $\widetilde{A} _{i _0,\dots, i _r}$. 
Let $c \in \R _+$, $d \in \R$ such that $x$ is $(c,d)$-bounded with respect to $\frac{1}{t}$, i.e. such that
any components $x _{i _0, \dots,i _r}$ are $(c,d)$-bounded (as $A _{\{f _{i _0}\cdots f _{i _r}\}}$-algebra with respect to $\frac{1}{t}$).
We have to prove that there exists
$y\in C ^{r-1} $
which satisfies 
$x= \partial (y)$. 
In fact, we will prove more precisely that such $y$ can be choosen 
$(2c, c+d)$-bounded with respect to $\frac{1}{t}$.

1) a)  Using \ref{prelem1-lemma7}, there exist two increasing sequences of integers $(m _j) _{j\in \bbN}$ and $(n _j) _{j\in \bbN}$  such that $n _j \leq c j+d$ and for any 
$0\leq i _0 < \cdots < i _r \leq m$,  the  components $x _{i _0, \dots,i _r}$  can be written in the form  
$x _{i _0 ,\dots , i _r} =  \sum _{j=0} ^{+\infty}  \pi ^{j} \frac{a _{j; i _0, \dots,  i _r}}{ (f _{i _0} \cdots f _{i _r})  ^{m _{j }}t ^{n _j}} $ where $a _{j;i _0, \dots,  i _r} \in  A$. 
We define $x _s \in C ^{r}$ by putting
\begin{equation}
\label{def-xs}
x _{s; i _0 ,\dots , i _r}
=
 \sum _{j=0} ^{2 ^{s+1}-1}  \pi ^{j} \frac{a _{j; i _0, \dots,  i _r}}{ (f _{i _0} \cdots f _{i _r})  ^{m _{j }}t ^{n _j}} .
\end{equation}
We put $M _s := m _{2 ^{s+1}-1}$,
$N _s := n _{2 ^{s+1}-1}$. 

b ) For any integer $u\in [-1,m]$, 
we denote by $C ^{u} _{(s)}$ the sub-A-algebra of $C ^{u}$
of elements $z=(z _{i _0, \dots, i _u})$ such that 
$t ^{N _s} (f _{i _0}\cdots f _{i _u}) ^{M _s} z _{i _0, \dots, i _u} \in A$.
We remark that if $z=(z _{i _0, \dots, i _u})\in C ^{u} _{(s)}$, then 
$f _{i _l} ^{N _s} z _{i _0, \dots, i _u} \in \widetilde{A} _{i _0,\dots, \widehat{i} _{l},\dots, i _{u}}$
for any $l=0,\dots, u$.
We notice that 
$x _s \in C ^{r} _{(s)}$ an that we have the factorization 
$\partial \colon C ^{u} _{(s)} \to C ^{u+1} _{(s)}$ 
(see notation \ref{not2-cech-cohomology}).

b) Since $A /\pi A $ is flat $k[[t]]$-algebra, then $A /\pi A  \to (A /\pi A ) _t$ is injective. Since $ (A /\pi A ) _t$ is integral, then  the canonical morphism 
$A /\pi A \to  \widetilde{A} _{i _0,\dots, i _r}/\pi \widetilde{A} _{i _0,\dots, i _r}= (A /\pi A)   _{t f _{i _0} \cdots f _{i _{r}}}$ is injective and  since $\widetilde{A} _{i _0,\dots, i _r}$ has no $\pi$-torsion, then from Lemma \ref{piC(n)},  for any integers $n \in \bbN$
and $u\in [-1,m]$,  we get $\pi ^n C ^{u} \cap C ^{u} _{(s)} = \pi ^n C ^{r} _{(s)}$.

2)  a) Since $f _0,\dots, f _m$ generate $A$, for any integer $N\geq 1$ there exist $r _{0,N},\dots, r _{m,N} \in A$ so that  $\sum _{i=0} ^{m} r _{i,N} f _i ^N = 1$ (indeed, this is obvious for $N=1$, then we use the $(m+1)N$th power of  this latter relation).

b) For any fixed integers $0\leq i,u\leq m$, we define the mapping  $\kappa ^{(s)} _i \colon C ^{u} _{(s)} \to C ^{u-1} _{(s)}$ by setting for any  $z = (z _{i _0, \dots,i _{u}}) \in C ^{u} _{(s)} $,
$$\kappa ^{(s)} _i (z)  _{i _0, \dots,i _{u-1}} = (-1) ^{j}r _{i, M _s} f _i ^{M _s} z _{i _0, \dots, i_{j-1} , i, i _{j}, \dots, i _{u-1}} $$
if there exists $j$ such that $i _{j-1} <i< i _{j}$ and $\kappa ^{(s)} _i (z)  _{i _0, \dots,i _{u-1}} =0$ otherwise. 

c) A standard computation  gives $\kappa ^{(s)} _i \circ \partial + \partial \circ \kappa ^{(s)} _i = r _{i, M _s} f _i ^{M _s}$, where $r _{i, M _s} f _i ^{M _s}$ means the multiplication by 
$r _{i, M _s} f _i ^{M _s}$ in $C ^{u} _{(s)}$ component by component. Indeed, fix  $0\leq i _0 < \cdots < i _u \leq m$. Suppose there exists $j$ such that $i _{j-1} <i< i _{j}$. On one hand, we compute
\small
\begin{gather}
\notag (\kappa ^{(s)} _i \circ \partial ) (z) _{i _0,\dots, i _{u}} = (-1) ^{j}r _{i, M _s} f _i ^{M _s} \partial (z) _{i _0, \dots, i_{j-1} , i, i _{j}, \dots, i _{u}}
\\ \notag = (-1) ^{j}r _{i, M _s} f _i ^{M _s}  \left ( (-1) ^j z _{i _0, \dots, i _{u}} + \sum _{k=0} ^{j-1} ( -1) ^k z _{i _0, \dots, \widehat{i} _k,\dots, i_{j-1} , i, i _{j}, \dots, i _{u}} + \sum _{k=j} ^{u}( -1) ^{k+1} z _{i _0, \dots, i_{j-1} , i, i _{j}, \dots, \widehat{i} _k,\dots, i _{u}} \right ).
\end{gather}
\normalsize
On the other hand, we have
\small
\begin{gather}
\notag (\partial \circ \kappa ^{(s)} _i  ) (z) _{i _0,\dots, i _{u}} = \sum _{k=0} ^{u} ( -1) ^k  \kappa ^{(s)} _i (z)_{i _0, \dots, \widehat{i} _k,\dots, i _{u}}
\\ \notag  = r _{i, M _s} f _i ^{M _s}  \left ( \sum _{k=0} ^{j-1} ( -1) ^k (-1) ^{j-1} z _{i _0, \dots, \widehat{i} _k,\dots, i_{j-1} , i, i _{j}, \dots, i _{u}} + \sum _{k=j} ^{u}( -1) ^{k} (-1) ^{j} z _{i _0, \dots, i_{j-1} , i, i _{j}, \dots, \widehat{i} _k,\dots, i _{u}} \right ).
\end{gather}
Adding both equalities, we get the desired formula. The case where there exists $j$ such that $i _j= i$ is checked with similar computation.

d) Put $\kappa ^{(s)} :=  \sum _{i=0} ^{m} \kappa ^{(s)} _i\colon C ^{u} _{(s)} \to C ^{u-1} _{(s)}$. We get $\kappa ^{(s)}  \circ \partial + \partial \circ \kappa ^{(s)} = id _{C ^{u} _{(s)}}$. 
For any integer $n\in \bbN$, we have the factorization $\kappa ^{(s)} \colon \pi ^n C ^{u} _{(s)} \to \pi ^n  C ^{u-1} _{(s)}$.

\medskip

3) We construct  the sequence $(y _{s}) _{s\in \bbN}$ of cochains  of $C ^{r -1}$ satisfying the following conditions
\begin{enumerate}[(i)]
\item  $y _{0} \in C ^{r-1} _{(0)}$ and for $s \geq 1$,  $y _{s} \in \pi ^{2 ^{s}} C ^{r-1} _{(s)}$  ; \item $\partial (\sum _{k=0} ^{s} y _k) = x \mod \pi ^{2 ^{s+1}}C ^{r} $, for any $s \geq 0$.
\end{enumerate}
{\it Proof.}  We proceed by induction on $s \in \bbN$.  Since $x _0 \in C ^{r} _{(0)}$ (see Notation \ref{def-xs}). Since $ x \equiv x _0 \mod \pi ^{2} C ^{r}$,  then we get $y _0 := \kappa ^{(0)} ( x _0) \in C ^{r-1} _{(0)}$.
since $\partial (x) = 0$, we get  $\partial (x _0) \equiv 0 \mod \pi ^{2} C ^{r+1}$. Since $\pi ^{2} C ^{r+1} \cap C ^{r+1} _{(0)} = \pi ^{2} C ^{r+1} _{(0)}$ (see part 1.b), then 
$\partial (x _0) \equiv 0 \mod \pi ^{2} C ^{r+1} _{(0)}$. Hence, $\kappa ^{(0)}  \circ \partial (x _0) \equiv 0 \mod \pi ^{2} C ^{r} _{(0)}$. Since $\kappa ^{(0)}  \circ \partial + \partial \circ \kappa ^{(0)} = id _{C ^{r} _{(0)}}$ (see part 2b), 
we get  $\partial \circ \kappa ^{(0)} ( x _0)  \equiv x _0 \mod \pi ^{2} C ^{r} _{(0)}$, i.e.  $\partial ( y _0)  \equiv x _0 \mod \pi ^{2} C ^{r}_{(0)}$. Hence,  $\partial ( y _0)  \equiv x  \mod \pi ^{2} C ^{r}$.

Suppose now that $y _0,\dots, y _{s-1}$ have been constructed for $s\geq 1$. Using the properties (i) and (ii) of the induction hypothesis,  since  $\pi ^{2 ^{s}} C ^{r} \cap C ^{r} _{(s)} = \pi ^{2 ^{s}} C ^{r} _{(s)}$ (see the part 1b) of the proof),
we get a element $z _s \in \pi ^{2 ^s}C ^{r} _{(s)}$ by putting 
$$z _{s}:= x _{s} - \partial (\sum _{k=0} ^{s-1} y _k).$$
Since $x _{s} \equiv x \mod \pi ^{2 ^{s+1}} C ^{r} $, then $\partial (x _{s}) \equiv \partial (x ) \mod \pi ^{2 ^{s+1}} C ^{r+1} $. Since $\partial (x )=0$ and $\partial (z _{s})= \partial (x _{s})$, this implies  $\partial (z _{s})\equiv 0 \mod \pi ^{2 ^{s+1}} C ^{r+1}$.
Since  $\pi ^{2 ^{s+1}} C ^{r+1} \cap C ^{r+1} _{(s)} = \pi ^{2 ^{s+1}} C ^{r+1} _{(s)}$, this  yields  $\partial (z _{s})\equiv 0 \mod \pi ^{2 ^{s+1}} C ^{r+1} _{(s)}$. Hence,  $\kappa ^{(s)} \circ \partial ( z _s) \equiv 0 \mod \pi ^{2 ^{s+1}} C ^{r} _{(s)}$.
Hence, with the part 2) of the proof, $\partial \circ \kappa ^{(s)} ( z _s) \equiv z _s \mod \pi ^{2 ^{s+1}} C ^{r}_{(s)}$. By putting  $y _s := \kappa ^{(s)} ( z _s)$, this means  we get an element $y _s \in  \pi ^{2 ^s }C ^{r-1} _{(s)}$
such that $\partial ( y _s) \equiv x _{s} - \partial (\sum _{k=0} ^{s-1} y _k) \mod \pi ^{2 ^{s+1}} C ^{r} _{(s)}$.

4) Using part 3) of the proof, we get the element $y:=  \sum _{k=0} ^{\infty} y _k\in C ^{r-1} $ which satisfies  $x= \partial (y)$.  We have  $N _s \leq c (2 ^{s+1}-1) +d = 2 c 2 ^s + (d -c) \leq 2 c 2 ^s + (c+d)$.
For $s \geq 1$, this yields that any element of  $\pi ^{2 ^{s}} C ^{r-1} _{(s)}$ is $(2c, c+d)$-bounded with respect to $\frac{1}{t}$ (as $C ^{r-1} (\fU _{\bullet}, \cO _{\fX})$-algebra if $r \geq 1$ or  as $A$-algebra if $r =0$).
Moreover, $N _0= n _1 \leq c +d$. Hence, $C ^{r-1} _{(0)}$ is also  $(2c, c+d)$-bounded with respect to $\frac{1}{t}$. Hence, $y$ is $(2c, c+d)$-bounded with respect to $\frac{1}{t}$.
\end{proof}

\begin{rem}
\label{rem-flat-exseq}
The lemma \ref{lem-vanishingCech} means that we have the exact sequence
$$0 \to \Gamma (\fU, \widetilde{A} ^{\Delta}) 
\to 
C ^{0} (\fU _{\bullet}, \widetilde{A} ^{\Delta})
\overset{\partial}{\to}\cdots \overset{\partial}{\to} 
C ^{m}(\fU _{\bullet}, \widetilde{A} ^{\Delta})\to 0.$$
The remark \ref{rem-CrotimesM} yields that for any flat 
$\widetilde{A}$-module $M$, we get the exact sequence  
$$0 \to \Gamma (\fU, M ^{\Delta}) 
\to 
C ^{0} (\fU _{\bullet}, M ^{\Delta})
\overset{\partial}{\to}\cdots \overset{\partial}{\to} 
C ^{m}(\fU _{\bullet}, M ^{\Delta})\to 0.$$
\end{rem}

\begin{rem}
\label{remoflem-vanishingCech}
With the notation of the proof of Lemma \ref{lem-vanishingCech},
we have check more precisely  that the element $y$ such that 
$x= \partial (y)$ can be chosen $(2c, c+d)$-bounded if $x$ is $(c,d)$-bounded.
This detail on Lemma \ref{lem-vanishingCech} will be used later. 
\end{rem}

\begin{lem}
\label{B is regular}
With notation \ref{ntnMtilde}, we suppose that
the generic fiber $\fX _\eta$ is a formal scheme smooth over $\eta$.
Then the ring $\widetilde{A}$ is regular.
\end{lem}

\begin{proof}
Since $\fX _\eta/\eta$ is smooth, 
since 
$\Gamma (\fX _\eta, \cO _{\fX _\eta}) \riso A _{\{t\}}$, 
then the ring 
$\widetilde{A} /\pi \widetilde{A} \riso A _{\{t\}} / \pi A _{\{t\}}$ is regular. 
Since $\cV$ is regular and $\pi \widetilde{A}$ is included in the Jacobson ideal of $\widetilde{A}$
(see \cite[1.6]{MonskyWashnitzer}), 
then by using \cite[Lemma 6.1]{MonskyWashnitzer} the ring $\widetilde{A}$ is regular.
\end{proof}

\begin{thm}
\label{sheaf-ThmB}
With notation \ref{ntnMtilde}, we suppose that
the generic fiber $\fX _\eta$ is a smooth over $\eta$ formal scheme.
Let $M$ be a finitely generated $\widetilde{A}$-module. 
\begin{enumerate}[(a)]
\item The presheaf
$M ^{\Delta}$ 
is in fact a sheaf. 
\item For every principal open subset $\fU $ of $\fX$ and $i>0$, we have
$H ^{i} (\fU,M ^{\Delta})=0$.

\end{enumerate}
\end{thm}

\begin{proof}
Since $\widetilde{A}$ is regular (see Lemma \ref{B is regular}), then $M$ has finite projective dimension.  In the case where $M$ is projective, $M$ is flat and then  we check the Theorem from Remark \ref{rem-flat-exseq}.
In general, we check the Theorem by induction on the projective dimension of $M$ (we copy the proof of \cite[Theorem 14 of Chapter 2]{meredith-weakformalschemes}).
\end{proof}

\begin{lem}
\label{affinesubset}
With notation \ref{ntnMtilde}, we suppose that
the generic fiber $\fX _\eta$ is a smooth over $\eta$ formal scheme.
Let $M$ be a finitely generated $\widetilde{A}$-module. 
Let $\fU$  be an affine open set of $\fX$. 
Then $\Gamma (\fU, M ^{\Delta} )= \Gamma (\fU,\cO _\fX) _{[t]} \otimes _{\widetilde{A}}M$.
\end{lem}

\begin{proof}
There exists $f _1,\dots, f _r$ so that  $\fU= \cup _{i=1} ^{r} \mathfrak{D} ( f _i)$. Let $A'$ so that $\fU = \Spf A'$. Set  $\widetilde{A}' := A' _{[t]}$, $M':= \widetilde{A}' \otimes _{\widetilde{A}} M$ and $M ^{\prime \Delta}$  be the sheaf on 
$\Spf A'$ defined as in \ref{ntnMtilde}. Let $f ' _i$ be the image of $f _i$ in $A'$. Since $\mathfrak{D} ( f '_i)= \mathfrak{D} ( f_i)$ we have $A _{\{ f _i\}} =A ' _{\{ f '_i\}}$.
Since $ M ^{\Delta}$ is a sheaf on $\fX$ then $\Gamma (\fU, M ^{\Delta})=\ker  ( \prod _i \Gamma (\mathfrak{D} ( f _i), M ^{\Delta}) \to  \prod _{i,j} \Gamma (\mathfrak{D} ( f _i f  _j), M ^{\Delta}) )$.
But, since  $M ^{\prime \Delta}$ is a sheaf on $\fU$ then we have $M'=\Gamma (\fU, M ^{\prime \Delta}) = \ker  ( \prod _i \Gamma (\mathfrak{D} ( f '_i), M ^{\prime \Delta}) \to  \prod _{i,j} \Gamma (\mathfrak{D} ( f '_i f  '_j), M ^{\prime \Delta}) )$.
Moreover,   $\Gamma (\mathfrak{D} ( f _i), \widetilde{A} ^{\Delta}) =   (A _{\{f _i\}} ) _{[t]} = (A ' _{\{ f '_i\}} )_{[t]} = \Gamma (\mathfrak{D} ( f '_i), \widetilde{A} ^{\prime \Delta}) $.
Hence,  $\Gamma (\mathfrak{D} ( f _i), M ^{\Delta})=\Gamma (\mathfrak{D} ( f _i), \widetilde{A} ^{\Delta}) \otimes _{\widetilde{A}} M \riso \Gamma (\mathfrak{D} ( f' _i), \widetilde{A} ^{\prime \Delta}) \otimes _{\widetilde{A}'} M' = \Gamma (\mathfrak{D} ( f' _i), M ^{\prime \Delta}) $.
Similarly,  $\Gamma (\mathfrak{D} ( f _i f _j), M ^{\Delta})= \Gamma (\mathfrak{D} ( f '_i f '_j), M ^{\prime \Delta})$. Hence, we get the  desired equality  $M'= \Gamma ( \fU,M ^{\Delta})$.
\end{proof}

\begin{ntn} \label{ntnB}
Let $\fX $ be an object of $\Sm ^\dag _\eta$. 
We denote by $\mathfrak{Aff} _{\fX}$ 
the category of affine open subsets of $\fX$. We get a presheaf $\widetilde{\cO} _{\fX} $ on $\mathfrak{Aff} _{\fX}$ defined by $\Gamma (\fU , \widetilde{\cO} _{\fX}) = \Gamma (\fU , \cO _{\fX} ) _{[t]}$ for any $\fU$ in $\mathfrak{Aff} _{\fX}$.
From Lemma \ref{affinesubset},  $\widetilde{\cO} _{\fX} $ is in fact a sheaf. We still denote by $\widetilde{\cO} _{\fX} $  the induced sheaf on $\fX$. 
Let $M$ be a finitely generated $\widetilde{A}$-module. Then we have:
$$M ^{\Delta}= \widetilde{\cO} _{\fX} \otimes _{\widetilde{A} } M.$$

\end{ntn}

\begin{prop} \label{123-coherentB}
Let $\fX $ be an object of $\Sm ^\dag _\eta$.  We keep notation of \ref{ntnB}.
\begin{enumerate}[(a)]
\item For any $\fU \in \mathfrak{Aff} _{\fX}$,  $\Gamma (\fU, \widetilde{\cO} _{\fX})$ is Noetherian ; 
\item For any $\fU ,\fU'\in \mathfrak{Aff} _{\fX}$, with $\fU' \subset \fU$, the homomorphism   $\Gamma (\fU, \widetilde{\cO} _{\fX})   \to   \Gamma (\fU', \widetilde{\cO} _{\fX})$  is flat ; 
\item The sheaf $\widetilde{\cO} _{\fX}$ is coherent.
\item For any $x \in X$, the rings $\widetilde{\cO} _{\fX,x}$ and $j _* (\cO _{\fX _\eta} ) _x$ are local. 
\item For any $x \in X$, $\pi\widetilde{\cO} _{\fX,x}$ is included in the maximal ideal of $\widetilde{\cO} _{\fX,x}$
and $\pi j _* (\cO _{\fX _\eta} ) _x$ is included in the maximal ideal of $j _* (\cO _{\fX _\eta} ) _x$.
\end{enumerate}
\end{prop}

\begin{proof}
The first and second statement are already known.  The third one is a consequence of Proposition \cite[3.1.1]{Be1}. 

Let us prove the statement (d). Let $\fU = \Spf A$ be an open subset of $\fX$ and $x \in U$. Let $\fp _x$ be the open ideal of $A$ associated with $x$ and $S:= A \setminus \fp _x$. 
The homomorphism $A \to \widetilde{A} _{<f>}$ is flat for any $f \in S$ (see \ref{prop-flat-[t]}). Hence, so is $A \to \widetilde{\cO} _{\fX,x}= \varinjlim _{f\in S}\widetilde{A} _{<f>}$.
Set $\fm _x:=\fp _x \otimes _A \widetilde{\cO} _{\fX,x}$. Let $z \in \widetilde{\cO} _{\fX,x}\setminus \fm _x$.  Since $\widetilde{\cO} _{\fX,x}= \varinjlim _{f\in S}\widetilde{A} _{<f>}$, then there exists $f \in A$, 
$z _f \in \widetilde{A} _{<f>} \setminus (\fp _x \otimes _A \widetilde{A} _{<f>} )$ whose image in $ \widetilde{\cO} _{\fX,x}$ is $z$.
Let $\overline{A} :=A /\pi A $ and for any $a \in A$ let denote by $\overline{a}$ its image in $\overline{A}$. We get $\widetilde{A} _{<f>} /\pi \widetilde{A} _{<f>} = \overline{A} _{\overline{f}\overline{t}}$. 
Then  the image of $z _f$ in $\widetilde{A} _{<f>} /\pi \widetilde{A} _{<f>}$,  can be written of the form $ \overline{a}/ (\overline{f}\overline{t}) ^k$  for some $a \in A \setminus \fp _x$.
Set $g := af $ and $y := a ^{k+1} /(g t) ^k\in A _{gt}\subset \widetilde{A} _{<g>}$.
We compute $a ^{k-1} f ^{2k} t ^k/g ^k\in A _{gt}$ is an inverse of $y$.  Since $y$ is the image of $a /(ft)^k\in A _{ft}$  in $A _{gt}$, this yields that
the image of $z _f$ in $\widetilde{A} _{<g>} /\pi \widetilde{A} _{<g>}$ is equal to the image of $y$ and is therefore invertible. Hence, the image of $z _f$ in $\widetilde{A} _{<g>}$ in invertible. This implies $z$ is invertible in $ \widetilde{\cO} _{\fX,x}$.
We prove similarly that $j _* (\cO _{\fX _\eta} ) _x= \varinjlim _{f\in S}A _{\{f t \}}$ is local with maximal ideal $\fp _x \otimes _A j _* (\cO _{\fX _\eta} ) _x$.

It remains to check (e).  For any principal open formal subscheme $\fU= \Spf A$, since $\Gamma (\fU , \widetilde{\cO} _{\fX})= A _{[t]}$ then it follows from  \cite[1.6]{MonskyWashnitzer} that
$\pi \Gamma (\fU , \widetilde{\cO} _{\fX})$ is included in the Jacobson ideal of $\Gamma (\fU , \widetilde{\cO} _{\fX})$. Hence we are done. 
\end{proof}

\begin{lem} \label{faithfullyflat1}
Let $\fX $ be an object of $\Sm ^\dag _\eta$. We denote by  $j \colon \fX _\eta \to \fX$ the open immersion of $\fS$-formal schemes. 
For any affine open subset $\fU$ of $\fX$, the homomorphism $\Gamma (\fU, \widetilde{\cO} _{\fX} ) \to \Gamma (\fU, j _* (\cO _{\fX _\eta} ))$ is faithful flat. 
Moreover, the homomorphism of ringed sheaves $\widetilde{\cO} _{\fX} \to j _* (\cO _{\fX _\eta} )$ is faithfully flat, 
i.e.  for any $x \in \fX$  the homomorphism $\widetilde{\cO} _{\fX,x} \to j _* (\cO _{\fX _\eta} ) _x$ is faithfully flat (see the definition given after  \cite[Lemma 4.3.8]{Be1}).
\end{lem}

\begin{proof}
Let $\fU$ be an affine open subset of $\fX$, $A:= \Gamma ( \fU , \cO _{\fX})$. Then  $\Gamma ( \fU , \widetilde{\cO} _{\fX})= A _{[t]}$ and  $\Gamma ( \fU , j _* (\cO _{\fX _\eta} )) = A _{\{t\}}$.
Since $A _{[t]} \to  A _{\{t\}}$ is faithfully flat, then we conclude using Lemma \cite[4.3.8]{Be1}.
\end{proof}

\subsubsection{Theorem of type $A$ for coherent $\widetilde{\cO} _\fX$-modules}
Let $\fX /\fS $ be an object of $\Sm ^\dag _\eta$. We suppose $\fX$ affine. 
Put $A := \Gamma (\fX , \cO _{\fX})$. Put $\widetilde{A} := A _{[t]}$ and $A _{\{t\}}$ its $p$-adic completion.

\begin{prop} \label{convergence}
Let $M$ be a finite $\widetilde{A}$-module. Let $\mu _1, \dots , \mu _u\in M$ and $c>0$ be a constant.  For $i =1,\dots, u$ let   $(a _{i,j} ) _{j\in \bbN}$ be a sequence of $A$,  $(n _{i,j} ) _{j\in \bbN}$ be a sequence of $\bbN$ such that $n _{i,j} \leq c (j+1)$
and  $\sum _{i=1} ^{u}  \frac{a _{i,j}}{t ^{n _{i,j}}} \mu _i \in p ^j M$.  Then $\sum _{j=0} ^{\infty} \sum _{i=1} ^{u}  \frac{a _{i,j}}{t ^{n _{i,j}}} \mu _i$ converges in $M$ for the $p$-adic topology.
\end{prop}

\begin{proof}
This is a consequence of \ref{convergencepre} and \ref{ex-B-wc}.
\end{proof}

\begin{lem}
\label{lem-finite-pres-tildeM}
Let $\cM $ be a sheaf of $\widetilde{A} ^{\Delta}$-module on $\fX$ of (globally of) finite presentation.
Put $M := \Gamma (\fX , \cM)$. Then $M$ is a $\widetilde{A}$-module of finite type and
the canonical morphism 
$M ^{\Delta}\to \cM$ is an isomorphism.
\end{lem}

\begin{proof}
By hypothesis, 
$\cM $ is the cokernel of a morphism of 
$\widetilde{A} ^{\Delta}$-modules of the form 
$\phi \colon (\widetilde{A} ^{\Delta}) ^r \to (\widetilde{A} ^{\Delta}) ^s$.
Let $P$ be the $\widetilde{A}$-module of finite type equal to the cokernel
of 
$\Gamma (\fX, \phi )\colon \widetilde{A} ^r \to \widetilde{A} ^s$.
From \ref{functor-sim-ff exact}, 
by applying the functor $\Delta$ to 
$\Gamma (\fX, \phi )$ we get $\phi$. Hence, 
since $\Delta$ is exact, 
we get that the cokernel of 
$\phi \colon (\widetilde{A} ^{\Delta}) ^r \to (\widetilde{A} ^{\Delta}) ^s$
is canonically isomorphic to $P ^\Delta$, i.e. 
$\cM \riso P ^\Delta$.
Hence, $\Gamma (\fX,\cM )\riso 
\Gamma (\fX,P ^\Delta) 
=P$.
\end{proof}

\begin{thm}
\label{thmAcoh}
Let $\cM $ be an $\widetilde{\cO} _{\fX}$-coherent sheaf on $\fX$. 
Then there exists a finitely generated $\widetilde{A}$-module $M$ such that 
$\cM = M ^{\Delta}$.
\end{thm}

\begin{proof}
0) It is enough to show that the natural homomorphism $\Gamma (\fX, \cM) \to \Gamma (\fX, \M/\pi \cM)$ is surjective.
Indeed, suppose this map is surjective. 

The sheaf $\widetilde{\cO} _{X}:= \widetilde{\cO} _{\fX}/\pi \widetilde{\cO} _{\fX}$ is a quasi-coherent $\cO _X$-algebra, 
and $\widetilde{O} _{X}:= \Gamma (X, \widetilde{\cO} _{X}) =(A/\pi A) _t$  is  noetherian. 
Hence, the theorem of type $A$ for coherent $\widetilde{\cO} _{X}$-modules holds, i.e. 
$\Gamma (X, -)$ and $\widetilde{\cO} _{X} \otimes _{ \widetilde{O} _{X}}-$
are quasi-inverse equivalences between the category of coherent $\widetilde{\cO} _{X}$-modules and 
$\widetilde{O} _{X}$-modules of finite type. 
Since $\cM/\pi \cM$ is a coherent sheaf of $\widetilde{\cO} _{X}$-modules, this yields there exist elements $\overline{f} _1,\dots, \overline{f} _s \in \Gamma (X, \cM/\pi \cM)$ which generate $\cM/\pi \cM$. 

The elements $\overline{f} _1,\dots, \overline{f} _s$ generate the $\widetilde{\cO} _{X,x}$-module 
$(\cM/\pi \cM) _x$, and $\widetilde{\cO} _{X,x}= \widetilde{\cO} _{\fX,x}/\widetilde{\cO} _{\fX,x} $ and $(\cM/\pi \cM) _x = \cM _x /\pi \cM _x$.
Lift each element $\overline{f} _i$ back to an element $f _i \in \Gamma (X, \cM)$. 
Hence, since $\pi\widetilde{\cO} _{\fX,x}$ is included in the maximal ideal of $\widetilde{\cO} _{\fX,x}$ (see \ref{123-coherentB}), 
since $\cM _x$ is a finitely generated $\widetilde{\cO} _{\fX,x} $-module (because $\cM$ is locally the sheaf associated to a finitely presented $\widetilde{\cO} _{\fX}$-module),
then it follows from Nakayama's lemma that $f _1,\dots, f _s$ generate $\cM _x$ as an $\widetilde{\cO} _{\fX,x} $-module.
Thus there exists a surjective homomorphism of coherent sheaves of $\widetilde{\cO} _{\fX} $-modules of the form $\alpha \colon \widetilde{\cO} _{\fX} ^s  \to  \cM$. Repeating the above
argument for the coherent $\widetilde{\cO} _{\fX} $-module $\ker \alpha$, we see there exists an integer
$r$ and a homomorphism of coherent sheaves $\beta \colon \widetilde{\cO} _{\fX} ^s \to \widetilde{\cO} _{\fX} ^s$
such that $\cM = \coker \beta$.  Hence $\cM = \widetilde{\cO} _{\fX} \otimes _{\widetilde{A}} \coker \Gamma (\fX, \beta)$.
\bigskip 

I) {\it We fix some notation.}  Choose a covering of $\fX= \cup _{i=0} ^{m}\fU _i $ by principal open sets $\fU _i = \mathfrak{D} (f _i)$ such that 
$\cM | \fU _i$ is a finitely presented $\widetilde{\cO} _{\fU _i}$-module. Using \ref{lem-finite-pres-tildeM}, this implies that
$M _i:= \Gamma (\fU _i , \cM )$ is a finitely generated $\widetilde{A} _{<f _i>}$-module such that $\cM | \fU _i = M _i ^{\Delta}$.
Put $\widetilde{A} _i := \widetilde{A} _{<f _i>}$,  $\fU _{i,j} := \fU _i \cap \fU _j$,  $\widetilde{A} _{i,j} := \widetilde{A} _{<f _i f _j >}$, $M _{i,j} :=  \Gamma ( \fU _{i,j},\cM) $.  
As in the proof of \ref{lem-vanishingCech}, for any integer $n\geq 1$ there exist $r _{0,n},\dots, r _{m,n} \in A$ so that  $\sum _{i=0} ^{m} r _{i,n} f _i ^n = 1$.
For each $i=0,\dots, m$, select generators $\mu _{i,1},\dots, \mu _{i,r} $ for $M_i$ over $\widetilde{A} _{i}$ for $r$ large enough.
Let $\mu _i = \overset{t}{}(\mu _{i,1},\dots, \mu _{i,r})$ be the corresponding column vector with coefficients in $M _i$. 

Let $N _{i,j} = (N _{i,j,\alpha,\beta}) _{\alpha \beta}\in M _{r, r} (\widetilde{A} _{i,j})$ such that  $\mu _i | _{\fU _{ij}} = N _{i,j} \mu _j | _{\fU _{ij}}$, i.e. we have in $M _{ij}$ the equalities 
$\mu _{i,\alpha} | _{\fU _{ij}} = \sum _{\beta =1} ^{r} N _{i,j,\alpha,\beta}\mu _{j,\beta} | _{\fU _{ij}}$.
Using  \ref{prelem1-lemma7},  there exists a constant $c \geq 0$ so that we can write  $N _{i,j} =  \sum _{q=0} ^{+\infty}  \frac{\pi ^{q} }{t ^{n _q}} U _{i,j,q}$ where $U _{i,j,q} \in  M _{r, r} (A _{f _i f _j} ) $, 
and $(n _q) _{q\in \bbN}$ is an increasing sequence of integers such that $n _q \leq c (q+1)$.  Put  $N _{i,j,h} :=  \sum _{q=0} ^{2 ^{h+1}-1}  \frac{\pi ^{q} }{t ^{n _q}} U _{i,j,q}$.

Let $\tau \in \Gamma (\fX, \cM /\pi \cM)$. Since $\tau | \fU _i\in \Gamma (\fU _i, \cM /\pi \cM) = M _i /\pi M _i$, there exists $g ^{0, i} \in M _{1, r} (A _{t f _i})$ such that the image of $g ^{0, i} \mu _i \in M _i$ in $M _i /\pi M _i$ is $\tau | \fU _i$. Increasing $c$ 
if necessary,   we can suppose  that  $ t ^{4c}g ^{0,i}  \in M _{1,r} (A _{f _i}) $.

We remark that $g ^{0, j} \mu _j |\fU _{ij}$ and $g ^{0, i} \mu _i |\fU _{ij}$  induce the same element of  $\Gamma (\fU _{ij}, \cM /\pi \cM)$.
Since $g ^{0, i} \mu _i |\fU _{ij} = g ^{0, i} N _{ij} \mu _j |\fU _{ij}$,  this yields $(g ^{0,i}N _{i,j} -g ^{0,j}) \mu _j | _{\fU _{ij}} \in \pi  M _{ij}$ for any $i,j$.
\bigskip 

II) By induction on $h \geq 1$, 
for any $i =0,\dots, m$ we construct  $g ^{0,i}, \dots, g ^{h-1,i}\in M _{1,r} ( A _{t f _i})$ satisfying the conditions 
for any $i ,j=0,\dots, m$ :
\begin{enumerate}[(1)]
\item $\sum _{s=0} ^{h-1} \left ( 
g ^{s,i}N _{i,j} -g ^{s,j} 
\right ) \mu _j | _{\fU _{ij}}
\in \pi ^{2 ^h -1} M _{ij}$ ; 
\item $g ^{s,i}\mu _i \in \pi ^{2 ^s -1} M _{i}$ for any $s =0,\dots, h-1$ ;
\item  $t ^{c 2 ^{s+2}} g ^{s,i} \in M _{1,r} (A _{f _i })$ for any $s =0,\dots, h-1$. 
\end{enumerate}
For $h =1$, this is already done. 
Suppose constructed $g ^{0,i}, \dots, g ^{h-1,i}$ satisfying the above conditions.
Set    
\begin{equation}\label{thmAcoh-proof-dfn-whij}
w _{h,i,j}:=\sum _{s=0} ^{h-1} \left (  g ^{s,i} N _{i,j,h} -g ^{s,j}  \right ) \in M _{1,r} ( A _{t f _i f _j})
\end{equation}
For an integer $a _h$ large enough, we get $f _i ^{a _h} w _{h,i,j} \in M _{1,r} ( A _{t f _j})$ for any $i$ and $j$.
Since  $N _{i,j,h} = N _{i,j} \mod \pi ^{2 ^{h+1} } M _{r,r} ( \widetilde{A} _{ij})$,  by using the induction property (1), we get  $w _{h,i,j} \mu _j | _{\fU _{ij}} \in \pi ^{2 ^h -1} M _{ij}$.
Since we have also $w _{h,i,j} \mu _j \in  (M _{j}) _{f _i}$, this yields $w _{h,i,j} \mu _j \in \pi ^{2 ^h -1} (M _{j}) _{f _i}$ (indeed, the map $(M _{j}) _{f _i} \to M _{ij}$ induces the isomorphism
$(M _{j}) _{f _i} / \pi ^n (M _{j}) _{f _i}  \riso M _{ij}/ \pi ^n M _{ij}$ for any inger $n \geq 1$).
Hence, increasing $a _h$ if necessary, we get  $f _i ^{a _h} w _{h,i,j} \mu _j \in \pi ^{2 ^h -1} M _{j}$. Hence, putting 
\begin{equation}\label{thmAcoh-proof-dfn-fh}
g ^{h,j}:= \sum _{i=0} ^{m} r _{i,a _h} f _i ^{a _h} w _{h,i,j} \in M _{1,r} ( A _{t f _j}),
\end{equation}
we get $ g ^{h,j}\mu _j\in \pi ^{2 ^h -1} M _{j}$, i.e. the property (2) holds.

Since $t ^{c 2 ^{h+1}}  N _{i,j,h}  \in  M _{r,r} ( A _{f _i f _j})$,  then using the induction hypothesis (3), we compute $t ^{c 2 ^{h+2}} w _{h, i,j} \in M _{1,r} ( A _{f _j})$. This implies $t ^{c 2 ^{h+2}} g ^{h,j} \in M _{1,r} (A _{f _j })$ for any $j$, i.e. the property (3) holds.

It remains to prove the property $(1)$.
Put $\fU _{i,j,l} := \fU _i \cap \fU _j \cap \fU _l$,  $\widetilde{A} _{i,j,l} := \widetilde{A} _{<f _i f _j f _l >}$, $M _{i,j,l} :=  \Gamma ( \fU _{i,j,l},\cM) $.  
Since $N _{l,i} |  _{\fU _{ijl}} N _{i,j}  |  _{\fU _{ijl}}= N _{l,j} |  _{\fU _{ijl}}$ in $M _{r,r} (\widetilde{A} _{i,j,l})$, then 
$N _{l,i,h} |  _{\fU _{ijl}} N _{i,j,h}  |  _{\fU _{ijl}}\equiv N _{l,j,h} |  _{\fU _{ijl}} \mod  \pi ^{2 ^{h+1}} M _{r,r} (\widetilde{A} _{i,j,l})$.
Since  $\pi ^{2 ^{h+1}} A _{t f _i f _j }= A _{t f _i f _j } \cap \pi ^{2 ^{h+1}} \widetilde{A} _{i,j,l}$, this yields 
\begin{equation}\label{thmAcoh-proof-equiv}(f _l ^{a _h} N _{l,i,h} )N _{i,j,h}  \equiv (f _l ^{a _h}N _{l,j,h})  \mod  \pi ^{2 ^{h+1}} M _{r,r} (A _{t f _i f _j } ).
\end{equation}
We conclude via the computation :
\begin{gather} \notag
(g ^{h,i} N _{i,j} -g ^{h,j} ) \mu _j | _{\fU _{ij}} \equiv (g ^{h,i} N _{i,j,h} -g ^{h,j} ) \mu _j | _{\fU _{ij}} \mod  \pi ^{2 ^{h+1}} M _{ij}
\\ \notag
\overset{\ref{thmAcoh-proof-dfn-fh}}{=} \sum _{l=0} ^{m} r _{l,a _h}  \left ( f _l ^{a _h}  w _{h,l,i} N _{i,j,h} - f _l ^{a _h} w _{h,l,j}  \right ) \mu _j | _{\fU _{ij}}
\\ \notag
\overset{\ref{thmAcoh-proof-dfn-whij}}{=} \sum _{l=0} ^{m} r _{l,a _h}  \sum _{s=0} ^{h-1} \left (  \left ( g ^{s,l} f _l ^{a _h}N _{l,i,h} -g ^{s,i} f _l ^{a _h}\right) N _{i,j,h} -  \left ( g ^{s,l} f _l ^{a _h} N _{l,j,h} -g ^{s,j} f _l ^{a _h}\right)  \right ) \mu _j | _{\fU _{ij}}
\\  \notag
\overset{\ref{thmAcoh-proof-equiv}}{\equiv} \sum _{l=0} ^{m} r _{l,a _h}  f _l ^{a _h} \sum _{s=0} ^{h-1} \left ( -g ^{s,i}  N _{i,j,h} + g ^{s,j} \right ) \mu _j | _{\fU _{ij}} \mod  \pi ^{2 ^{h+1}} M _{ij}
\\ \notag
= \sum _{s=0} ^{h-1} \left ( -g ^{s,i}  N _{i,j,h} + g ^{s,j} \right ) \mu _j | _{\fU _{ij}} \equiv \sum _{s=0} ^{h-1} \left ( -g ^{s,i}  N _{i,j} + g ^{s,j} \right ) \mu _j | _{\fU _{ij}} \mod  \pi ^{2 ^{h+1}} M _{ij}.
\end{gather}

III) From the properties
$(2)$ and $(3)$ of the part II, using Proposition \ref{convergence}, the sum
$\sum _{s=0} ^{\infty} 
g ^{s,i} \mu _i$ converges in $M _{i}$.
Moreover, in $M _{ij}$ we get from the property $(1)$ of the part II the first equality in $M _{ij}$:
$\sum _{s=0} ^{\infty}
g ^{s,j}  \mu _j  | _{\fU _{ij}}
=
\sum _{s=0} ^{\infty}
g ^{s,i}N _{i,j} 
\mu _j | _{\fU _{ij}}
=
\sum _{s=0} ^{\infty}
g ^{s,i} \mu _i | _{\fU _{ij}}$.
Hence, there exists an element $x$ of 
$\Gamma (\fX ,\cM)$ which induces 
$\sum _{s=0} ^{\infty}
g ^{s,i} \mu _i$ on $M _i$ for any $i$.
Since the image of $g ^{0, i} \mu _i \in M _i$ in $M _i /\pi M _i$ is $\tau | \fU _i$, then using the Property II.$(2)$ 
the image of $\sum _{s=0} ^{\infty} g ^{s,i} \mu _i$ in  $M _i/\pi M _i$ is $\tau |\fU _i$. Hence, 
the image of $x$ in  $\Gamma (\fX ,\M/\pi \cM)$ is $\tau$.
\end{proof}

\begin{coro} \label{equcatBDelta}
The following properties hold.
\begin{enumerate}[(a)]
\item 
The functors $M \mapsto M ^{\Delta}$ and $\cM \mapsto \Gamma (\fX,\cM)$  induce canonically exact quasi-inverse equivalences of categories between the category of
 finitely generated (resp. and projective) $\widetilde{A}$-modules and coherent (resp. and projective) $\widetilde{A} ^{\Delta}$-modules.
\item For any coherent $\widetilde{A} ^{\Delta}$-module $\cM$, for any integer $n \geq 1$, 
$H ^{n} (\fX,\cM)=0$.
\end{enumerate}
\end{coro}

\begin{proof}
The second part of the corollary is  a consequence of Theorems \ref{sheaf-ThmB}.(b) and  \ref{thmAcoh}. We get the exactness of $\cM \mapsto \Gamma (\fX,\cM)$. The exactness of $M \mapsto M ^{\Delta}$ is obvious.
We conclude using Lemma \ref{lem-finite-pres-tildeM} and \ref{sheaf-ThmB}.(a).
\end{proof}

\begin{coro}\label{coh-gen-x}
Let $\cM$ be a coherent $\widetilde{\cO} _{\fX}$-module. Let $a _1,\dots, a _n \in \Gamma (\fX, \cM)$ and $x \in \fX$. 
If  the sequence $a _{1,x},\dots, a _{n,x}$ generates $\cM _x$ then there exists an affine open $\fU$ containing $x$ such that $a _1 |\fU,\dots, a _n |\fU$ generates $\cM |\fU$.
\end{coro}
\begin{proof} Following \ref{equcatBDelta}, $M:= \Gamma (\fX,\cM)$  is a finitely generated $\widetilde{A}$-modules
and the canonical morphism $M ^\Delta \to \cM$ is an isomorphism. Hence, we reduce to prove it for $M ^\Delta$ instead of $ \cM$, which is 
straightforward from \ref{thmAcoh}.
\end{proof}

\begin{coro}\label{IDeltatilderegular0}
Let $\cM$ be a coherent $\widetilde{\cO} _{\fX}$-module. Let $a _1,\dots, a _n \in \Gamma (\fX, \cM)$.  The sequence $a _{1},\dots, a _{n}$  generates $\cM$ if and only if 
the sequence $a _{1}|\fX _\eta,\dots, a _{n}|\fX _\eta$ generates $\cM|\fX _\eta$.
\end{coro}
\begin{proof}This is a consequence  of the theorems of type $A$ for coherent $\cO _{\fX _\eta}$-modules and coherent $\widetilde{\cO} _{\fX}$-modules (see \ref{sheaf-ThmB}.(a) and \ref{thmAcoh}) and of the faithful flatness of $\widetilde{A} \to A _{\{t\}}$.
\end{proof}

\begin{prop} \label{Homcohcohiscoh}
Let $\cM$ and $\cN$ be two coherent $\widetilde{\cO} _{\fX}$-modules. Put $M := \Gamma (\fX, \cM)$,  $N := \Gamma (\fX, \cN)$
The canonical morphism  $\widetilde{\cO} _{\fX} \otimes _{\widetilde{A}} \mathrm{Hom} _{\widetilde{A}} ( M, N) \to  \mathcal{H}om _{\widetilde{\cO} _{\fX}} ( \M, \cN)$ is an isomorphism.
\end{prop}

\begin{proof}
Since $M$ is an $\widetilde{A}$-module of finite presentation,  then using the five Lemma  to check that the canonical morphism
$\widetilde{\cO} _{\fX}\otimes _{\widetilde{A}}\mathrm{Hom} _{\widetilde{A}} ( M, N)\to  \mathcal{H}om _{\widetilde{\cO} _{\fX}} (\widetilde{\cO} _{\fX} \otimes _{\widetilde{A}} M, \widetilde{\cO} _{\fX} \otimes _{\widetilde{A}} N)$
is an isomorphism,  we reduce to the obvious case where  $M$ is a free $\widetilde{A}$-module. We conclude by using \ref{equcatBDelta}.
\end{proof}

\begin{coro}
Let $M$ and $N$ be  two finitely generated $\widetilde{A}$-modules. We have the isomorphism
$\widetilde{\cO} _{\fX} \otimes _{\widetilde{A}} \R \mathrm{Hom} _{\widetilde{A}} ( M, N) \riso \R\mathcal{H}om _{\widetilde{\cO}_\fX} ( M ^\Delta, N ^\Delta) .$
Hence, for any $r\in \bbN$, $ \mathrm{Ext} ^r _{\widetilde{A}} ( M, N) ^\Delta \riso \mathcal{E}xt ^r _{\widetilde{\cO}_\fX} ( M ^\Delta, N ^\Delta)$.
\end{coro}

\begin{proof}
Let $L _\bullet \riso M$ be a left resolution of $M$ by free $\widetilde{A}$-module of finite type.
This yields the left resolution  
$L _\bullet ^\Delta \riso M^\Delta$
of $M^\Delta$ by free $\widetilde{\cO} _{\fX}$-module of finite type.
Hence, 
$\widetilde{\cO} _{\fX}
\otimes _{\widetilde{A}}
\R \mathrm{Hom} _{\widetilde{A}} ( M, N)
\riso
\widetilde{\cO} _{\fX}
\otimes _{\widetilde{A}}
 \mathrm{Hom} _{\widetilde{A}} ( L _\bullet, N)
\underset{\ref{Homcohcohiscoh}}{\riso}
\mathcal{H}om _{\widetilde{\cO}_\fX} ( L _\bullet ^\Delta, N ^\Delta)
\liso
\R\mathcal{H}om _{\widetilde{\cO}_\fX} ( M ^\Delta, N ^\Delta)$.
\end{proof}

\begin{coro} \label{prop-regularbis0}
Let $u \colon \fX \to \fY$ be a morphism of $\Sm ^\dag _\eta$ which is a finite morphism.
Let $\cM$ be coherent $ \widetilde{\cO} _{\fX}$. 
Then $u _* \cM$ is a coherent $ \widetilde{\cO} _{\fY}$-module. If moreover $\fY= \Spf B$, $\fX= \Spf A$ and $\phi$ is the forgetful functor from the category of $A$-modules to the category of 
$B$-module, we have $u _* (M ^\Delta )= (\phi (M) ) ^\Delta$. 
\end{coro}

\subsubsection{Regularity}
\begin{rem} \label{rem-ff}
Let $A \to B$ be a faithfully flat morphism of commutative noetherian rings. If $B$ is regular then so is $A$. Indeed, in that case  $\Spec B \to \Spec A$ is surjective, flat and then we can apply \cite[0.17.3.3.(i) and 0.17.3.6]{EGAIV1}.
\end{rem}

\begin{lem}\label{lem-regularbis}
Let $\fX /\fS $ be an object of $\Sm ^\dag _\eta$. For any $\fU \in \mathfrak{Aff} _{\fX}$, $\Gamma (\fU, \widetilde{\cO} _{\fX})$ is regular. For any point $x$ of $\fX $, the rings $j _* (\cO _{\fX _\eta} ) _x$ and $\widetilde{\cO} _{\fX, x}$ are regular. 
\end{lem}

\begin{proof}
We can suppose $\fX$ affine and $\fU= \fX$. Put  $A:= \Gamma ( \fX , \cO _{\fX})$.  Since $\Gamma (\fU, \widetilde{\cO} _{\fX}) = A _{[t]}$, the first statement is Lemma \ref{B is regular}.
Let  $\mathfrak{p}$ bet the  open prime ideal of $A$ corresponding to $x$. Put $S := A \setminus \mathfrak{p}$ and $C := A _{\{ t\}}$. Since $\fX _\eta$ is smooth,  since $\cV$ is regular and $\pi C$ is included in the Jacobson ideal of $C$, 
from \cite[Lemma 6.1]{MonskyWashnitzer}, we get the regularity of $C$. (Remark that since $A _{[t]}\to A _{\{ t\}}=C$ is faithfully flat, from the remark \ref{rem-ff}
we get the regularity of $A _{[t]}$, which gives a second proof of Lemma \ref{B is regular}). Then $S ^{-1}C=C _\mathfrak{p}$ and its $p$-adic completion $C \{ S ^{-1}\}$  are regular. 
We compute $j _* (\cO _{\fX _\eta} ) _x = \underset{f \in S}{\underrightarrow{\lim}} C _{\{ f\}} = C _{\{S\}}$ (see notation of \cite[0.7.6]{EGAI}).
Since  $C _{\{S\}}\to C \{ S ^{-1}\}$ is a (faithfully) flat homomorphism of local noetherian ring (see \cite[0.7.6.17]{EGAI}),  since  $C \{ S ^{-1}\}$ is regular, then using  \ref{rem-ff}
we get that  $C _{\{S\}}$ is regular, i.e.  $j _* (\cO _{\fX _\eta} ) _x$ is regular.  Since $\widetilde{\cO} _{\fX, x} \to j _* (\cO _{\fX _\eta} ) _x$ is (faithfully) flat of local ring (see \ref{faithfullyflat1}) then using  \ref{rem-ff}
we get the regularity of the ring $\widetilde{\cO} _{\fX, x}$.
\end{proof}

\begin{prop} \label{coh-perf}
With \ref{123-coherentB}. Let $\fX /\fS $ be an object of $\Sm ^\dag _\eta$. 
For any $\fU \in \mathfrak{Aff} _{\fX}$,  $\gl. \dim(\Gamma (\fU,\widetilde{\cO} _{\fX}))= \dim \fX _\eta$, where $\gl. \dim$ is the global dimension (see definition \cite[2.3.4.1]{Car25}).
Then  $\tor \dim \widetilde{\cO} _{\fX} = \dim \fX _\eta$, where $\tor \dim \widetilde{\cO} _{\fX} $ is the tor-dimension of $\widetilde{\cO} _{\fX} $ (see definition \cite[1.4.3.5]{Car25}). 
For any coherent $\widetilde{\cO} _{\fU}$-module $\cE$,   there exists a resolution of length $\leq \dim \fX _\eta$ of $\cE$ by projective coherent
$\widetilde{\cO} _{\fU}$-modules. We have $D ^\mathrm{b} _{\mathrm{coh}} (\widetilde{\cO} _{\fX}) = D ^\mathrm{b} _{\mathrm{parf}} (\widetilde{\cO} _{\fX})$.
\end{prop}

\begin{proof}
The equality $\gl. \dim(\Gamma (\fU,\widetilde{\cO} _{\fX}))= \dim \fX _\eta$ is a consequence of the regularity of  $\Gamma (\fU, \widetilde{\cO} _{\fX})$  for any $\fU \in \mathfrak{Aff} _{\fX}$ (see \ref{lem-regularbis}).
Since $\Gamma (\fU,\widetilde{\cO} _{\fX})$ is noetherian, using \cite[4.1.5]{Weibel-HomologicalAlg} we get $\gl. \dim(\Gamma (\fU,\widetilde{\cO} _{\fX}))=\tor \dim(\Gamma (\fU,\widetilde{\cO} _{\fX}))$.
Via \cite[1.4.3.22]{Car25}, this yields $\tor \dim \widetilde{\cO} _{\fX} = \dim \fX _\eta$.
The last equality comes from \cite[1.4.3.29]{Car25}. We finish using \cite[1.4.3.22]{Car25} and the theorem of type $A$ (see \ref{equcatBDelta}).
\end{proof}

\begin{lem} \label{lem-eqconregseq}
Let $\fX /\fS $ be an object of $\Sm ^\dag _\eta$.
We suppose $\fX$ affine. 
Put 
$A:= \Gamma ( \fX , \cO _{\fX})$,
$\widetilde{A} := \Gamma ( \fX , \widetilde{\cO} _{\fX})$. 
Let $t _1,\dots, t _r \in \widetilde{A}$ be a sequence. 
To avoid ambiguity, we denote by $\widehat{t} _1,\dots, \widehat{t} _r$ the image of $t _1,\dots, t _r$ in $A _{\{t\}}=\Gamma ( \fX _\eta, \cO _{\fX _\eta})$, 
the $p$-adic completion of $\widetilde{A}$. 
Let $\cM$ be a coherent $\widetilde{\cO} _{\fX} $-module,
$M := \Gamma (\fX, \cM)$, 
$\widehat{M}:= \Gamma (\fX _\eta, \cM)$. 
The following assertions are equivalent
\begin{enumerate}[(a)]
\item \label{lem-eqconregseq1} the sequence $t _1,\dots, t _r $ is $M$-regular (resp. $M$-quasi-regular) in the sense of \cite[0.15.1.7]{EGAIV1} ; 
\item \label{lem-eqconregseq2} the sequence $t _1,\dots, t _r $ is $\cM$-regular (resp. $\cM$-quasi-regular) in the sense of \cite[0.15.2.2]{EGAIV1} ;
\item \label{lem-eqconregseq3} the sequence $\widehat{t} _1,\dots, \widehat{t} _r $ is $\widehat{M}$-regular (resp. $\widehat{M}$-quasi-regular) in the sense of \cite[0.15.1.7]{EGAIV1}  ;
\item \label{lem-eqconregseq4} the sequence $\widehat{t} _1,\dots, \widehat{t} _r $ is $\cM | _{\fX _\eta} $-regular (resp. $\cM | _{\fX _\eta} $-quasi-regular) in the sense of \cite[0.15.2.2]{EGAIV1}. 
\end{enumerate}

\end{lem}

\begin{proof}
Since $\widetilde{A} \to A _{\{t\}}$ is faithfully flat, then we get the equivalence between 
the conditions \ref{lem-eqconregseq1} and \ref{lem-eqconregseq3} (use \cite[0.15.1.14]{EGAIV1}).
Let 
$\cJ$ be the coherent ideal of 
$\widetilde{\cO} _{\fX}$ generated by $t _1,\dots, t _r $.
Using the Theorem of type $A$ \ref{thmAcoh}, using the fact that 
the global section functor is exact on the category of coherent $\widetilde{\cO} _{\fX}$-modules, 
we check the equivalence between 
the conditions \ref{lem-eqconregseq1} and \ref{lem-eqconregseq2}.
Similarly, we check the conditions \ref{lem-eqconregseq3} and \ref{lem-eqconregseq4} are equivalent.
\end{proof}

\begin{dfn}
\label{regularideal-EGAIV}
Let $\fX /\fS $ be an object of $\Sm ^\dag _\eta$.
Let  $\cJ$ be a ideal of  $\widetilde{\cO} _{\fX}$.
We say that $\cJ$ is {\it regular} (resp. {\it quasi-regular}) if $\cJ$ is regular with respect to the ringed space $(\fX, \widetilde{\cO} _{\fX})$ in the sense of  \cite[IV.16.9.1]{EGAIV4}, i.e.,  if for any $x \in \Supp (\widetilde{\cO} _{\fX} / \cJ)$, there exist
an open set $\fU $ of $\fX$ containing $x$, an $\widetilde{\cO} _{\fU} $-regular  (resp. $\widetilde{\cO} _{\fU} $-quasi-regular) sequence 
$t _1,\dots, t _r$ of elements of $\Gamma (\fU, \cJ)$ which generate $\cJ |\fU$.
\end{dfn}

\begin{prop}\label{prop-coh-reg-id} Let $\fX /\fS $ be an object of $\Sm ^\dag _\eta$.
Let  $\cJ$ be a coherent ideal of  $\widetilde{\cO} _{\fX}$.
The ring $\cJ$ is regular 
if and only if  for any $x \in \Supp (\widetilde{\cO} _{\fX} / \cJ)$ the ideal $\cJ _x$ of $\widetilde{\cO} _{\fX,x} $ is regular (in the sense of \cite[IV.16.9.7]{EGAIV4}). 
\end{prop}
\begin{proof}
Since $\widetilde{\cO} _{\fX,x} $ is local then $\cJ _x$ is a regular ideal of $\widetilde{\cO} _{\fX,x} $
if and only if there exists an $\widetilde{\cO} _{\fX,x} $-regular sequence of $\cJ _x$ generating $\cJ _x$. Hence, if $\cJ$ is regular 
then for any $x \in \Supp (\widetilde{\cO} _{\fX} / \cJ)$ the ideal $\cJ _x$ of $\widetilde{\cO} _{\fX,x} $ is regular. Conversely, let $x \in \Supp (\widetilde{\cO} _{\fX} / \cJ)$. 
By hypothesis, there exists an $\widetilde{\cO} _{\fX,x} $-regular sequence $a _{1,x},\dots, a _{n,x}$ of $\cJ _x$ generating $\cJ _x$.
Let $\fU$ be an open containing $x$ such that $a _{1,x},\dots, a _{n,x}$ comes from some elements  $a _{1},\dots, a _{n}$ of $\Gamma (\fU, \cJ)$.
Via \ref{coh-gen-x}, shrinking $\fU$ if necessary we can suppose that the sequence $a _{1},\dots, a _{n}$ generates $\cJ |\fU$.
Following \cite[0.15.2.4]{EGAIV1}, shrinking $\fU$ if necessary, since $\cJ$ is coherent then we get that the sequence $a _{1},\dots, a _{n}$ is $\widetilde{\cO} _{\fU}$-regular. 
Hence, we are done.
\end{proof}

\begin{prop} \label{prop-regularbis}
Let $u \colon \fX \to \fP $ be a morphism of $\Sm ^\dag _\eta$ which is a closed immersion. We denote by $\cI$ the ideal defining $u$.  Then  $u _*  \widetilde{\cO} _{\fX}\riso \widetilde{\cO} _{\fP} /\widetilde{\cI} $
and $\widetilde{\cI} $ is a regular ideal of $\widetilde{\cO} _{\fP}$.
\end{prop}

\begin{proof}
By using \ref{prop-regularbis0}, we get the isomorphism: $u _*  \widetilde{\cO} _{\fX} \riso \widetilde{\cO} _{\fP} /\widetilde{\cI} $. 
Since $\widetilde{\cI} $ is a coherent ideal of $ \widetilde{\cO} _{\fP}$, then following \ref{prop-coh-reg-id}, we have to check that for any point $x$ of $\fX$ the ideal  $\widetilde{\cI}  _x$ of $ \widetilde{\cO} _{\fP,x}$ is regular.
Following \ref{123-coherentB} and  \ref{lem-regularbis},  the rings $\widetilde{\cO} _{\fP,x}$ and $ \widetilde{\cO} _{\fP,x} / \widetilde{\cI}  _x \riso  \widetilde{\cO} _{\fX,x}$ are regular. 
This yields from \cite[IV.19.1.2]{EGAIV4} that $\widetilde{\cI}  _x$ is a regular ideal. 
\end{proof}

\subsection{Weak completion of differential operators with overconvergent singularities over weakly smooth log formal $\fS ^\sharp$-schemes}   \label{ntn-Smdagdag}
We denote by $\Sm ^\dag _{\fS ^\sharp}$ the follow category. Its objects are formal log schemes $\fX ^\sharp$  which are log smooth over $\fS ^\sharp$ and whose generic fiber $\fX _\eta$ (i.e. the fiber of $\fX ^\sharp$
over $\eta = \Spf \cV [[t]] \{ \frac{1}{t}\}$) is smooth over $\eta$ (in particular $M _{\fX} | \fX _\eta = \cO ^* _{\fX _\eta}$). 
Let $\fX ^\sharp$ and $\fY ^\sharp$ be two objects of $\Sm ^\dag _{\fS ^\sharp}$. A morphism $\fX ^\sharp\to \fY ^\sharp$ of $\Sm ^\dag _{\fS ^\sharp}$ is a morphism $\fX \to \fY$ of the underlying morphism of formal log-schemes.
Such objects are called weakly smooth log formal $\fS ^\sharp$-schemes. 
Hence, $\Sm ^\dag _{\fS ^\sharp}$ is a full subcategory of $\Sm ^\dag _{\eta}$ (see Definition \ref{dfn-Sm-eta}).
Let $m\in \bbN \cup \{ \infty\}$.

\subsubsection{Local descriptions, noetherianity, regularity}
Let $\fX ^\sharp$ be an object of $\Sm ^\dag _{\fS ^\sharp}$. Let $j \colon \fX _\eta \hookrightarrow \fX ^{\sharp}$ be the open immersion.

\begin{ntn}\label{Dtilde}
If $\cM$ is an $\cO _{\fX}$-module, we put $\widetilde{\M}:= \cM \otimes _{\cO _\fX} \widetilde{\cO} _\fX$, where $\widetilde{\cO} _\fX$ is defined in \ref{ntnB}. For instance, we set 
$\widetilde{\Omega} _{\fX ^\sharp/\fS ^\sharp}:= \widetilde{\cO} _{\fX}   \otimes _{\cO _{\fX}  }\Omega _{\fX^\sharp/\fS^\sharp}$,
$\widetilde{\cD} ^{(m)} _{\fX^\sharp/\fS^\sharp} := \widetilde{\cO} _\fX \otimes _{\cO _\fX}\cD ^{(m)} _{\fX^\sharp/\fS^\sharp}$,
$\widetilde{\cD} ^{(m)} _{\fX^\sharp/\fS^\sharp,n} := \widetilde{\cO} _\fX \otimes _{\cO _\fX}\cD ^{(m)} _{\fX^\sharp/\fS^\sharp,}$,
$\widetilde{\cD} _{\fX^\sharp/\fS^\sharp} := \widetilde{\cO} _\fX \otimes _{\cO _\fX} \cD _{\fX^\sharp/\fS^\sharp}$.
We compute that $\widetilde{\cO} _\fX$ is endowed with a canonical structure of  $\cD ^{(m)} _{\fX^\sharp/\fS^\sharp}$-module so that the inclusion
$\widetilde{\cO} _\fX \hookrightarrow  j _* \cO _{\fX _\eta}$ is $\cD ^{(m)} _{\fX^\sharp/\fS^\sharp}$-linear. Hence, 
this $\cD ^{(m)} _{\fX^\sharp/\fS^\sharp}$-module  structure of $\widetilde{\cO} _\fX$ is compatible with its structure of $\cO _\fX$-algebra. This implies 
$\widetilde{\cD} ^{(m)} _{\fX^\sharp/\fS^\sharp} $ is a sheaf of rings (see \cite[2.3.5]{Be1}). Since $\cO _\fX \to \widetilde{\cO} _\fX$ is flat, then 
the ring $\widetilde{\cD} ^{(m)} _{\fX ^\sharp /\fS ^\sharp}$ is endowed with the filtration by the order $(\widetilde{\cD} ^{(m)}_{\fX ^\sharp /\fS ^\sharp,n} ) _{n\in \bbN}$ which satisfies 
$$\widetilde{\cD} ^{(m)}_{\fX ^\sharp /\fS ^\sharp,n} \cdot \widetilde{\cD} ^{(m)}_{\fX ^\sharp /\fS ^\sharp,n'}  \subset\widetilde{\cD} ^{(m)}_{\fX ^\sharp /\fS ^\sharp,n+n'} $$  
for any integers $n,n'$.
\end{ntn}

\begin{lem} \label{Otilde*}
We have $\widetilde{\cO} ^* _{\fX } =\widetilde{\cO} _{\fX} \cap j _{*} \cO _{\fX _\eta} ^* $.
\end{lem}

\begin{proof}
Since this is local, we can suppose $\fX ^\sharp$ affine.  Put $A := \Gamma (\fX, \cO _{\fX})$. We get  $ A  _{[t]}=\Gamma (\fX, \widetilde{\cO} _{\fX})$, $ A ^* _{[t]}=\Gamma (\fX, \widetilde{\cO} ^* _{\fX})$, 
$A ^* _{\{t\}} = \Gamma (\fX, j _{*} \cO ^* _{\fX _\eta})$.  Since the canonical morphism $A  _{[t]} /p A  _{[t]} \to  A _{\{t\}} / p  A _{\{t\}}$ is an isomorphism (see \cite[Thm 1.4]{MonskyWashnitzer}),
since $x \in A  _{[t]}$ (resp. $x \in A  _{\{t\}}$) is invertible if and only if the image of $x $ in $A  _{[t]} /p A  _{[t]}$ (resp. $A _{\{t\}} / p  A _{\{t\}}$) is invertible  then $A ^* _{[t]} = A _{[t]} \cap A ^* _{\{t\}} $. Hence, we are done. 
\end{proof}

\begin{lem} \label{rem-QcapY=BX}
We have $\widetilde{\cO} _{\fX } =\widetilde{\cO} _{\fX,\bbQ} \cap j _{*} \cO _{\fX _\eta}$,
$\widetilde{\cD} ^{(m)} _{\fX ^\sharp/\fS ^{\sharp}} =\widetilde{\cD} ^{(m)} _{\fX ^\sharp/\fS ^{\sharp},\bbQ}  \cap j _{*} \cD ^{(m)} _{\fX _\eta/\eta}$,
$\widetilde{\cD} ^{(m)} _{\fX ^\sharp/\fS ^{\sharp},n} =\widetilde{\cD} ^{(m)} _{\fX ^\sharp/\fS ^{\sharp},\bbQ}  \cap j _{*} \cD ^{(m)} _{\fX _\eta/\eta,n}$,
$\widetilde{\cD}  _{\fX ^\sharp/\fS ^{\sharp}}  = \widetilde{\cD}  _{\fX ^\sharp/\fS ^{\sharp},\bbQ}  \cap j _{*} \cD  _{\fX _\eta/\eta}$
and $\widetilde{\cD}  _{\fX ^\sharp/\fS ^{\sharp},n}  = \widetilde{\cD}  _{\fX ^\sharp/\fS ^{\sharp},\bbQ}  \cap j _{*} \cD  _{\fX _\eta/\eta,n}$. 
\end{lem}

\begin{proof}
Let us check the first equality.  Since this is local, we can suppose $\fX ^\sharp$ affine. 
Put $A := \Gamma (\fX, \cO _{\fX})$, $\widetilde{A} := A _{[t]}=\Gamma (\fX, \widetilde{\cO} _{\fX})$. Then $A _{\{t\}} = \Gamma (\fX, j _{*} \cO _{\fX _\eta})$. 
Since the canonical morphism $\widetilde{A} /p \widetilde{A} \to  A _{\{t\}} / p  A _{\{t\}}$
is an isomorphism (see \cite[Thm 1.4]{MonskyWashnitzer}) then $\widetilde{A} \cap p A _{\{t\}} = p \widetilde{A}$.
Hence, we have checked $\widetilde{\cO} _{\fX }  = \widetilde{\cO} _{\fX,\bbQ} \cap j _{*} \cO _{\fX _\eta}$.

Let us check the second equality. Since this is local, we can suppose $\fX ^\sharp$ affine and there exist log coordinates 
$u _1,\dots, u _d \in \Gamma ( \fX , \cO _{\fX})$ of $\fX ^{\sharp} /\fS ^{\sharp} $.
Let $\partial _{\sharp 1},\dots, \partial _{\sharp d}$ be the corresponding logarithmic derivations. 
The sheaf $\cD ^{(m)} _{\fX^\sharp/\fS^\sharp}$ is $\cO _{\fX ^\sharp}$-free with the basis 
$\{ \underline{\partial} _{\sharp} ^{<\underline{i}> _{(m)}},\, \underline{i} \in \bbN ^d\}$
and  $\cD ^{(m)} _{\fX _\eta/\eta}$ is $\cO _{\fX _\eta}$-free with the basis 
$\{ \underline{\partial} _{\sharp} ^{<\underline{i}> _{(m)}},\, \underline{i} \in \bbN ^d\}$
(where we still denote by $\underline{\partial} _{\sharp} ^{<\underline{i}> _{(m)}}$
the image of $\underline{\partial} _{\sharp} ^{<\underline{i}> _{(m)}}$ 
via  $\widetilde{\cD} ^{(m)} _{\fX ^\sharp/\fS ^{\sharp}}  \hookrightarrow j _{*} \cD ^{(m)} _{\fX _\eta/\eta}$).
Hence, we get the second equality from the first one.  We proceed similarly for the other ones. 
\end{proof}

\begin{lem}
We have the canonical map  $d \colon \widetilde{\cO} _{\fX}  \to \widetilde{\Omega} _{\fX ^\sharp /\fS ^\sharp}  $ making commutative the diagram
\begin{equation}
\label{dtildeOmegasharp}
\xymatrix{ { \cO _{\fX}}  \ar@{^{(}->}[r] \ar[d] ^-{d} & { \widetilde{\cO} _{\fX}}  \ar@{^{(}->}[r] \ar@{.>}[d] &  {j _* \cO _{\fX _\eta}}  \ar[d] ^-{j _* d _\eta}
\\  {\Omega _{\fX ^\sharp /\fS ^\sharp} }  \ar@{^{(}->}[r] & {\widetilde{\Omega} _{\fX ^\sharp /\fS ^\sharp}}  \ar@{^{(}->}[r] &  {j _* \Omega _{\fX _\eta/\eta} ,}     }
\end{equation}
where $d$ and $d _\eta$ are the structural morphisms, and $\Omega _{\fX ^\sharp /\fS ^\sharp} \subset \widetilde{\Omega} _{\fX ^\sharp /\fS ^\sharp}$ is the canonical $\cO _\fX$-linear map.
\end{lem}

\begin{proof}
Since $j ^* \Omega _{\fX ^\sharp /\fS ^\sharp}  =\Omega _{\fX _\eta/\eta}$, we get by functoriality the commutativity of the outline of the diagram \ref{dtildeOmegasharp}. 
Since the Lemma is local, we can suppose $\fX = \Spf A$ and that  there exist log coordinates  $u _1,\dots, u _d$ of $\fX ^{\sharp} /\fS ^{\sharp} $.
Let $\partial _{\sharp 1},\dots, \partial _{\sharp d}$ be the corresponding logarithmic derivations,  which corresponds to the dual basis  $d \log u _1, \dots, d \log u _d$ of $\Omega _{\fX ^\sharp /\fS ^\sharp} $. 
Then  $d \colon  \cO _{\fX}  \to  \Omega _{\fX ^\sharp /\fS ^\sharp} $ is given by  $a \mapsto  \sum _{j=1} ^d \partial _{\sharp j} (a) d \log u _j$.
Let  $t _1,\dots, t _d \in \Gamma (\fX _\eta, \cO _{\fX ^*_\eta})$ be the image of respectively $u _1,\dots, u _d$ of $\fX ^{\sharp} /\fS ^{\sharp} $
via $\Gamma (\fX, M _{\fX}) \to \Gamma (\fX, \cO _{\fX}) \subset \Gamma (\fX _\eta, \cO _{\fX _\eta})$. It follows from \ref{Otilde*} that 
 $t _1,\dots, t _d \in \Gamma (\fX , \widetilde{\cO} ^* _\fX)$.
Let $\partial _{1},\dots, \partial _{d}$ be the corresponding derivations of $\fX _\eta /\eta$,  which corresponds to the dual basis  $d t _1, \dots, d t  _d$ of $\Omega _{\fX _\eta /\eta} $.
Then,  $d _\eta \colon \cO _{\fX _\eta} \to  \Omega _{\fX _\eta /\eta} $ is given by  $c \mapsto  \sum _{j=1} ^d \partial _{j} (c) d t _j$.

Let  $b \in \Gamma (\fX, \widetilde{\cO} _\fX)=A _{[t]}$. There exists a sequence $(a _i) _{i\in \bbN}$ of elements of $A$, a real number $c$ such that $v _p (a _i) \geq \frac{i}{c}-1$ 
and $b=\sum _{i\in \bbN} \frac{a _i}{ t ^i}$.  Since the restriction of  $\partial _{\sharp j} $ on $\fX _\eta$ is equal to  $ t _j\partial _{j} $, then we compute
$t _j \partial _{j} (b) = \sum _{i\in \bbN} \frac{ \partial _{\sharp j} (a _i)}{ t ^i} \in \Gamma (\fX, \widetilde{\cO} _\fX)$,  for any $j =1,\dots, d$.
Since $d  t _j = t _j d \log u _j \in \Gamma (\fX, \widetilde{\Omega} _{\fX ^\sharp /\fS ^\sharp})$,  then  $\partial _{j} (b) d  t _j = \left (\sum _{i\in \bbN} \frac{ \partial _{\sharp j} (a _i)}{ t ^i} \right ) d \log u _j $,  for any $j =1,\dots, d$.
Hence, we are done.
\end{proof}

\begin{dfn}
\label{dfn-ovcvloccoor}
Let $t _1,\dots, t _d \in \Gamma (\fX, \widetilde{\cO} _{\fX}  )$. We say that   $t _1,\dots, t _d$ are ``overconvergent local coordinates of $\fX ^\sharp/\fS ^\sharp$'' if   $t _1 |\fX _\eta,\dots, t _d|\fX _\eta$ are local coordinates of $\fX _\eta/\eta$. 
\end{dfn}

\begin{ex}[Local logarithmic coordinates]  \label{ntn-loc-coord-partial}
Let $x$ be a point of $X$. There exists an affine open $\fU ^\sharp $ of $\fX ^\sharp$ containing $x$, there exist log coordinates 
$u _1,\dots, u _d \in \Gamma ( \fU , \cO _{\fU})$ of $\fU ^{\sharp} /\fS ^{\sharp} $.
Then $u _1,\dots, u _d$ are also overconvergent local coordinates of $\fU ^\sharp/\fS ^\sharp$.
We finish this example with the following remark. Let $\partial _{\sharp 1},\dots, \partial _{\sharp d}$ be the corresponding logarithmic derivations. 
Since, the sheaf $\cD ^{(m)} _{\fU^\sharp/\fS^\sharp}$ is $\cO _{\fU ^\sharp}$-free with basis 
$\{ \underline{\partial} _{\sharp} ^{<\underline{i}> _{(m)}},\, \underline{i} \in \bbN ^d\}$
then  $\widetilde{\cD} ^{(m)} _{\fU ^\sharp/\fS^\sharp}$ is 
$\widetilde{\cO} _{\fU ^\sharp}$-free with basis $\{ \underline{\partial} _{\sharp} ^{<\underline{i}> _{(m)}},\, \underline{i} \in \bbN ^d\}$.
\end{ex}

\begin{rem}\label{rem-freeness}
Let $A \to B$ be a faithfully flat morphism and $M$ be an $A$-module of finite type.  
If $x _1,\dots, x _n \in M$ induce a $B$-basis of $B\otimes _A M$,
then $x _1,\dots, x _n $ is a $A$-basis of $M$. In particular, $M$ is free.
\end{rem}

We can improve the remark \ref{rem-freeness} as follows (it might not be useful here):
\begin{lem} \label{freeness}
Let $A$ be a commutative noetherian $\cV$-algebra such that
its $p$-adic completion $\widehat{A} $ has no $p$-torsion and 
$A \to \widehat{A} $ is faithfully flat. 
Let $M$ be an $A$-module of finite type such that 
its $p$-adic completion $\widehat{M}$ is a free $\widehat{A} $-module. 
Then $M$ is a free $A$-module. 
\end{lem}

\begin{proof}
Let $x _1,\dots, x _n \in M$ such that $\overline{x} _1, \dots, \overline{x} _n$ is a  $A / pA$-basis of $M / p M = \widehat{M}/p \widehat{M}$. 
Let $\phi \colon A ^n \to M$ be the $A$-linear map defined by the sequence $x _1,\dots, x _n$. Applying the functor $\widehat{A} \otimes _{A} -$, 
we get the map $\widehat{\phi} \colon \widehat{A} ^n \to \widehat{M}$. From \cite[3.2.2.(ii)]{Be1}, the map $\widehat{\phi} $ is surjective.
Since $\widehat{M}$ has no $p$-torsion and $\widehat{A} ^n $ is separated (for the $p$-adic topology), 
then $\widehat{\phi} $ is injective. Since the extension $A \to \widehat{A} $ is faithfully flat,
the map $\phi$ is also an isomorphism.  
\end{proof}

\begin{ntn}
[Local description of the sheaf of differential operators with overconvergent coordinates]
\label{localcoord-ovcv} Suppose  $\fX ^\sharp$ is affine and that there exist
overconvergent local coordinates $t _1,\dots, t _d  \in \Gamma (\fX, \widetilde{\cO} _{\fX})$  of $\fX ^\sharp/\fS ^\sharp$.

Let $\cI _\Delta$ be the ideal of $\cO _{\fX  \times _{\fS } \fX }$ defining the diagonal immersion $\Delta _{\fX / \fS  } \colon  \fX  \hookrightarrow \fX  \times _{\fS  } \fX $.
Remark that the map  $\fX ^\sharp \times _{\fS ^\sharp} \fX ^\sharp \to \fX  \times _{\fS } \fX $ induces the isomorphism
$(\fX ^\sharp \times _{\fS ^\sharp} \fX ^\sharp ) _\eta \riso ( \fX  \times _{\fS } \fX ) _\eta \riso \fX _\eta  \times _{\eta } \fX _\eta$.
Since  $\Omega _{\fX ^\sharp /\fS ^\sharp} |\fX _\eta= \Omega _{\fX  /\fS } |\fX _\eta =  \Omega _{\fX _\eta  /\eta }$, 
then the homomorphism of $\cO _\fX$-modules of finite type $\Omega _{\fX  /\fS } \to \Omega _{\fX ^\sharp /\fS ^\sharp} $ becomes an isomorphism
after the extension $\cO _{\fX} \to \widetilde{\cO} _{\fX}$,  i.e. $\widetilde{\Omega} _{\fX ^\sharp /\fS ^\sharp}= \widetilde{\Omega} _{\fX  /\fS }$.
From \ref{prop-regularbis}, $\widetilde{\cI}  _\Delta$ is regular with respect to the ringed space $(\fX, \widetilde{\cO} _{\fX})$. In fact, we get more precisely the following. 

We denote by $\tau _i : = 1 \otimes t _i - t _i \otimes 1 \in  \Gamma( \fX  \times _{\fS } \fX, \widetilde{\cO} _{\fX \times _{\fS } \fX } )  $ 
and by $d t _i$ the image of $\tau _i $ in $\widetilde{\Omega} _{\fX ^\sharp /\fS ^\sharp}= \widetilde{\cI}  _\Delta / \widetilde{\cI}  _\Delta ^2$.
Remark that the map $d \colon \widetilde{\cO} _{\fX} \to  \widetilde{\cI}  _\Delta / \widetilde{\cI}  _\Delta ^2$ defined by $a\mapsto 1 \otimes a - a \otimes 1 \mod  \widetilde{\cI}  _\Delta ^2$
 is equal to the factorisation of   the map   $d \colon \widetilde{\cO} _{\fX}  \to \widetilde{\Omega} _{\fX ^\sharp /\fS ^\sharp}  $ making commutative the diagram \ref{dtildeOmegasharp}.
Since  $t _1 |\fX _\eta,\dots, t _d|\fX _\eta$ are local coordinates of $\fX _\eta/\eta$,  then  $\tau _1 | \fX _\eta  \times _{\eta } \fX _\eta,\dots, \tau _d | \fX _\eta  \times _{\eta } \fX _\eta$ is an 
 $\cO _{\fX _\eta\times _{\eta } \fX _\eta }$-regular sequence which generates $\widetilde{\cI}  _\Delta | \fX _\eta\times _{\eta } \fX _\eta $.
Hence, from \ref{lem-eqconregseq} and  \ref{IDeltatilderegular0},  $\tau _1 ,\dots, \tau _d $ is an  $\widetilde{\cO} _{\fX \times _{S } \fX }$-regular sequence which generates  $\widetilde{\cI}  _\Delta$.
In particular $ \widetilde{\Omega} _{\fX ^\sharp /\fS ^\sharp}$  is  $\widetilde{\cO} _{\fX}$-free  with   $d t _1, \dots, d t _d$ 
 as  a basis (or this is a consequence of \ref{freeness}).   Let $\partial _1, \dots, \partial _d\in  \widetilde{\mathcal{T}} _{\fX ^\sharp /\fS ^\sharp} $ be the dual basis of  $d t _1,\dots, d t _d\in  \widetilde{\Omega} _{\fX ^\sharp /\fS ^\sharp}$.
We denote (if there is no risk of confusion) by  $\partial  _1, \dots , \partial  _d\in \mathcal{T} _{\fX _\eta/\eta,}$ the dual basis of $d t _1 |\fX _\eta,\dots, d t _d|\fX _\eta \in  \Omega _{\fX _\eta /\eta }$.

We set $\underline{\partial} ^{[\underline{i}]} := \frac{\underline{\partial} ^{\underline{i}}}{\underline{i}!}  \in \widetilde{\cD} _{\fX ^\sharp /\fS ^\sharp,\bbQ} $. 
Its image in $\cD _{\fX _\eta/\eta,\bbQ}$ lies in fact in $\cD _{\fX _\eta/\eta}$ and is still denoted by 
$\underline{\partial} ^{[\underline{i}]}$. 
Since  $t _1 |\fX _\eta,\dots, t _d|\fX _\eta$ are local coordinates of $\fX _\eta/\eta$, then the $\cO _{\fX _\eta}$-module (for both structures) $\cD _{\fX _\eta/\eta} $ is free 
with $\{\underline{\partial} ^{ [\underline{i}]} \ ; \ \underline{i} \in \bbN ^{r}\}$ as a basis.

We set $\underline{\partial} ^{<\underline{i}> _{(m)}}:=  \frac{[\underline{i}/p ^{m} ] !}{\underline{i}!} \underline{\partial} ^{\underline{i}}\in \widetilde{\cD} ^{(m)} _{\fX ^\sharp /\fS ^\sharp,\bbQ}$. 
Its image in $\cD ^{(m)} _{\fX _\eta/\eta,\bbQ}$ lies in fact in $\cD ^{(m)} _{\fX _\eta/\eta}$ and is still denoted by 
$\underline{\partial} ^{<\underline{i}> _{(m)}}$. 
Since  $t _1 |\fX _\eta,\dots, t _d|\fX _\eta$ are local coordinates of $\fX _\eta/\eta$, then the $\cO _{\fX _\eta}$-module (for both structures) $\cD ^{(m)} _{\fX _\eta/\eta} $ is free 
with $\{\underline{\partial} ^{<\underline{i}> _{(m)}} \ ; \ \underline{i} \in \bbN ^{r}\}$ as a basis.

Using the equalities of \ref{rem-QcapY=BX}, we get 
$\underline{\partial} ^{<\underline{i}> _{(m)}}\in \widetilde{\cD} ^{(m)}_{\fX ^\sharp /\fS ^\sharp,n} $,
$\underline{\partial} ^{[\underline{i}]} \in \widetilde{\cD} _{\fX ^\sharp /\fS ^\sharp,n} $ if $| \underline{i}| \leq n$.
Since $\widetilde{\cD} ^{(m)}_{\fX _\eta /\fS _\eta} $ is an $\cO _{\fX _\eta}$-module locally free of finite type, 
since $\cD _{\fX _\eta /\fS _\eta} $ (resp. $\cD ^{(m)}_{\fX _\eta /\fS _\eta} $) is $\cO _{\fX _\eta } $-free with  
$\{ \underline{\partial} ^{[\underline{i}]} \ ; \ \underline{i} \in \bbN ^{r}\}$ (resp. $\{ \underline{\partial} ^{<\underline{i}> _{(m)}} \ ; \ \underline{i} \in \bbN ^{r}\}$) as a basis, 
then using \ref{rem-freeness} we get that $ \widetilde{\cD} _{\fX ^\sharp /\fS ^\sharp} $ (resp. $ \widetilde{\cD} ^{(m)}_{\fX ^\sharp /\fS ^\sharp} $) is $\widetilde{\cO} _{\fX } $-free with  
$\{ \underline{\partial} ^{[\underline{i}]} \ ; \ \underline{i} \in \bbN ^{r}\}$ (resp. $\{ \underline{\partial} ^{<\underline{i}> _{(m)}} \ ; \ \underline{i} \in \bbN ^{r}\}$) as a basis. 
\end{ntn}

\begin{rem}
\label{bounded-twodescrp}
Suppose $\fX$ affine, and set
 $A:= \Gamma (\fX, \cO _{\fX})$.
 Let $x \in A \{ \frac{1}{t}\}$.
 Let $c,d \in \R$ with $c >0$. With notation \ref{padic-norm}, the following assertion are equivalent
 \begin{enumerate}[(a)]
\item The element $x$ is $(c, d)$-bounded as $A$-algebra with respect to $\frac{1}{t}$, i.e. following \ref{(c,d)bounded}, 
we can write (non uniquely)  $x$ in the form 
$x = \sum _{n\in \bbN} \pi ^n \frac{a _{n} }{t ^{\alpha _n}}$
with $a _{n} \in A$ and $\alpha _n \in \Z _{\leq cn+d}$,
for any  $n \in \bbN$.

 \item We can write (non uniquely)  $x$ in the form 
$x = \sum _{j\in \bbN}  \frac{a _{j}}{t ^{j}}$
with $a _{j} \in A$
such 
$v _{\pi} ( a _{j}) \geq \frac{j-d}{c}$,
for any  $j \in \bbN$.

 \end{enumerate}

\end{rem}

\begin{lem} \label{(c,c-i)bounded} 
Suppose  $\fX ^\sharp$ is affine and that there exist
overconvergent local coordinates $t _1,\dots, t _d  \in \Gamma (\fX, \widetilde{\cO} _{\fX})$  
of $\fX ^\sharp/\fS ^\sharp$.
Let  $z\in \Gamma (\fX _\eta, \widehat{\cD} ^{(m)}  _{\fX _\eta/\eta})$, 
 $A:= \Gamma (\fX, \cO _{\fX})$.
With Notation \ref{localcoord-ovcv},
we can write uniquely 
$z$ in the form
$z = \sum _{\underline{i}\in \bbN ^d}  b _{\underline{i}} \underline{\partial} ^{<\underline{i}> _{(m)}}$
(resp. $z = \sum _{\underline{i}\in \bbN ^d}  \underline{\partial} ^{<\underline{i}> _{(m)}}b _{\underline{i}} $), 
with $b _{\underline{i}}\in A \{ \frac{1}{t}\}$ converging to $0$. 
 Let $c _0,d _0 \in \R$ with $c _0 >0$. The following assertion are equivalent

\begin{enumerate}[(a)]
\item We have $b _{\underline{i}}\in A _{[t]}$ 
and $b _{\underline{i}}$  is $(c _0, d _0-|\underline{i}|)$-bounded
as $A$-algebra with respect to $\frac{1}{t}$, 
for any $\underline{i} \in \bbN ^d$.

\item It is possible to write $z$ in the form 
$z = \sum _{\underline{i}\in \bbN ^d} \sum _{j\in \bbN} a _{\underline{i},j} \frac{1}{t ^{j}}\underline{\partial} ^{<\underline{i}> _{(m)}}$
(resp. $z = \sum _{\underline{i}\in \bbN ^d} \sum _{j\in \bbN} \frac{1}{t ^{j}}\underline{\partial} ^{<\underline{i}> _{(m)}} a _{\underline{i},j}$)
with $a _{\underline{i},j} \in A$ such that 
$v _{\pi} ( a _{\underline{i},j}) \geq \frac{j+|\underline{i}| -d _0}{c _0}$, for any $\underline{i} \in \bbN ^d$ and $j \in \bbN$.
\end{enumerate}

\end{lem}

\begin{proof}
This is straightforward from the remark \ref{bounded-twodescrp}.
\end{proof}

\begin{lem}
\label{comp-nota} Suppose $m\in \bbN$.
Suppose  $\fX ^\sharp$ is affine, and $\fX ^\sharp/\fS ^\sharp$ has logarithmic coordinates. Set $A:= \Gamma (\fX, \cO _{\fX})$, $\widetilde{A}:= A _{[t]}=\Gamma (\fX, \widetilde{\cO} _{\fX})$.
Fix $t _1,\dots, t _d  \in \widetilde{A}$   some overconvergent local coordinates  of $\fX ^\sharp/\fS ^\sharp$ (see Example \ref{ntn-loc-coord-partial}).
We keep the notation \ref{localcoord-ovcv} concerning the basis of $ \widetilde{\cD} ^{(m)}_{\fX ^\sharp /\fS ^\sharp} $  as $\widetilde{\cO} _{\fX } $-module given by $\{ \underline{\partial} ^{<\underline{i}> _{(m)}} \ ; \ \underline{i} \in \bbN ^{r}\}$. 

There exists $c _0\in \R _{>0}$, $d _0 \in \R$ satisfying the following property:  for any $e \in \R$, for any
$(c _0,e)$-bounded element $b \in \widetilde{A}$ as $A$-algebra with respect to $1/t$ (see the definition \ref{(c,d)bounded}), for any $\underline{i}\in \bbN ^d$ such that $|\underline{i}|\leq p ^{m}$,
the elements $\beta _{\underline{i}, \underline{j}}  \in \widetilde{A}$ (resp. $\beta ' _{\underline{i}, \underline{j}}  \in \widetilde{A}$) defined by the formula
\begin{equation} \label{partialb-bounded}
\underline{\partial} ^{<\underline{i}> _{(m)}}b = \sum _{\underline{j}\leq \underline{i}} \beta _{\underline{i}, \underline{j}}  \underline{\partial} ^{<\underline{j}> _{(m)}}
(\mathrm{resp. } \;  b \underline{\partial} ^{<\underline{i}> _{(m)}} = \sum _{\underline{j}\leq \underline{i}} \underline{\partial} ^{<\underline{j}> _{(m)}} \beta ' _{\underline{i}, \underline{j}} ),
\end{equation}
are $( c_0, e + 2 d _0)$-bounded as $A$-algebra with respect to $1/t$.
\end{lem}

\begin{proof}
Using the transposition, the respective case is a consequence of the non respective case.  Let us check this latter case. 
Let  $u _1,\dots, u _d$ be log coordinates of $\fX ^{\sharp} /\fS ^{\sharp} $ and let  $\underline{\partial} _{\sharp} ^{<\underline{j}> _{(m)}}$ be the corresponding  $\cO _{\fX}$-basis of  $\cD ^{(m)}_{\fX ^\sharp /\fS ^\sharp}$
(beware the image of $u _1,\dots, u _d$ in  $\Gamma (\fX, \widetilde{\cO} _{\fX})$  are not necessarily equal to $t _1,\dots, t _d $).
Let $\underline{i}\in \bbN ^d$ such that $|\underline{i}|\leq p ^{m}$, i.e.  $\underline{\partial} ^{<\underline{i}> _{(m)}}  \in \widetilde{\cD} ^{(m)}_{\fX ^\sharp /\fS ^\sharp, p ^m}$. 
Since  $ \widetilde{\cD} ^{(m)}_{\fX ^\sharp /\fS ^\sharp, |\underline{i}|} $  is $\widetilde{\cO} _{\fX } $-free with $\{ \underline{\partial} ^{<\underline{j}> _{(m)}} \ ; \ |\underline{j}| \leq |\underline{i}|\}$ 
or  $\{ \underline{\partial} _{\sharp} ^{<\underline{j}> _{(m)}} \ ; \ |\underline{j}| \leq |\underline{i}|\}$ as a basis, then  there exist
$b _{\underline{i}, \underline{j}} \in \widetilde{A}$ and  $b '_{\underline{i}, \underline{j}} \in \widetilde{A}$ such  
$\underline{\partial} _\sharp ^{<\underline{i}> _{(m)}} = \sum _{\underline{j}} b _{\underline{i}, \underline{j}}  \underline{\partial} ^{<\underline{j}> _{(m)}}$
and $\underline{\partial} ^{<\underline{i}> _{(m)}} = \sum _{\underline{j}} b '_{\underline{i}, \underline{j}}  \underline{\partial} _\sharp ^{<\underline{j}> _{(m)}}$,
where both sums are taken over $\underline{j} \in \bbN ^d$ such that  $|\underline{j}| \leq |\underline{i}|$.
Let $(c _0, d _0) \in \R _+ ^2$ such that for any  $\underline{i}, \underline{j} \in \bbN ^d$ such that  $|\underline{i}|\leq p ^{m}$ and  $|\underline{j}| \leq |\underline{i}|$,
the elements  $b _{\underline{i}, \underline{j}} $ and  $b '_{\underline{i}, \underline{j}} $ are $(c _0,d _0)$-bounded.

Let $b \in A _{[t]}$  be a $(c _0,e)$-bounded  element  for some $e \in \R _+$. Since $\underline{\partial} _\sharp ^{<\underline{i}> _{(m)}} (t) = 0$ for any $\underline{i}\in \bbN ^d$, 
we get $\underline{\partial} _\sharp ^{<\underline{i}> _{(m)}} ( b) $ is $(c _0,e)$-bounded.  For any $\underline{j}\in \bbN ^d$, 
we have the formula $\underline{\partial} _\sharp ^{<\underline{j}> _{(m)}} b  = \sum _{\underline{h}\leq \underline{j}} \left \{ \begin{smallmatrix}  \underline{j}\\  \underline{h} \end{smallmatrix} \right\}
\underline{\partial} _\sharp ^{<\underline{j} -\underline{h}> _{(m)}} ( b)   \underline{\partial} _\sharp ^{<\underline{h}> _{(m)}}$.
With the formulas of \ref{(c,d)bounded-ring} concerning the product and the sum of bounded elements,  we can write
\begin{equation}
\notag \underline{\partial} ^{<\underline{i}> _{(m)}}b = \sum _{\underline{j}\leq \underline{i}} \beta _{\underline{i}, \underline{j}}  \underline{\partial} ^{<\underline{j}> _{(m)}},
\end{equation}
where $\beta _{\underline{i}, \underline{j}}  \in A _{[t]}$ are $( c_0, e + 2 d _0)$-bounded.
\end{proof}

\begin{prop}
[Weak completion of level $m$]
\label{prop-wkcpdscplevelm}
Suppose $m\in \bbN$. Suppose  $\fX ^\sharp$ is affine, and $\fX ^\sharp/\fS ^\sharp$ has logarithmic coordinates. 
Let 
$D ^{(m)} := \Gamma (\fX, \cD ^{(m)} _{\fX^\sharp/\fS^\sharp} )$, $\widetilde{D} ^{(m)} := \Gamma (\fX, \widetilde{\cD} ^{(m)}  _{\fX^\sharp/\fS^\sharp})$,  $A:= \Gamma (\fX, \cO _{\fX})$.
Fix $t _1,\dots, t _d  \in A _{[t]}$   some overconvergent local coordinates  of $\fX ^\sharp/\fS ^\sharp$. With notation \ref{localcoord-ovcv}, let 
 $T _{(m)}:= 
 \{\frac{1}{t}\} 
 \cup 
\{\partial _i ^{[p ^j]}\ |\ i=1,\dots, d\; ;\; j= 0,\dots, m\}$.
The symbole
$\dag$ 
means 
weak completion as $A$-ring.

\begin{enumerate}[(i)]
\item We have the equalities  $(\widetilde{D} ^{(m)})^{\dag } = (D ^{(m)} [\frac{1}{t}] )^{\dag } = T _{(m)} ^{\dag } $ ; 

\item Let  $z$ be an element of $\widehat{D} ^{(m)}\{ \frac{1}{t}\}$,  the $p$-adic completion  of $\widetilde{D} ^{(m)}$. We have $z \in (\widetilde{D} ^{(m)}) ^{\dag }$ if and only if we can write $z$ in the form
\begin{equation} \label{wkcpdscplevelm}
z = \sum _{k\in \bbN}\sum _{\underline{i}\in \bbN ^d}  \frac{a _{\underline{i},k}}{t ^{k}}\underline{\partial} ^{<\underline{i}> _{(m)}}
(resp.\, z = \sum _{k\in \bbN}\sum _{\underline{i}\in \bbN ^d}  \underline{\partial} ^{<\underline{i}> _{(m)}}\frac{a _{\underline{i},k}}{t ^{k}})
\end{equation}
where $a _{\underline{i},k} \in A$ are such that there exists a constant $c>0$ satisfying
$v _{\pi} ( a _{\underline{i},k}) \geq \frac{k+|\underline{i}|}{c} -1$ (resp. $v _{p} ( a _{\underline{i},k}) \geq \frac{k+|\underline{i}|}{c} -1$), for any $\underline{i} \in \bbN ^d$ and $k \in \bbN$.

\item \label{prop-wkcpdscplevelm(iii)}
Let $z =  \sum _{\underline{i}\in \bbN ^d} \widetilde{a} _{\underline{i}} \underline{\partial} ^{<\underline{i}> _{(m)}} $ (resp.  $z =  \sum _{\underline{i}\in \bbN ^d} \underline{\partial} ^{<\underline{i}> _{(m)}}  \widetilde{a} _{\underline{i}} $),
where $\widetilde{a} _{\underline{i}} \in A _{\{t\}}$, be an element of $\widehat{D}\{ \frac{1}{t}\}$. We have  $z \in \widetilde{D} ^{\dag}$ if and only if there exists a constant $c>0$ such that 
$\widetilde{a} _{\underline{i}} \in A _{[t]} $ are $(c, c -|\underline{i}|)$-bounded.
\end{enumerate}
\end{prop}

\begin{proof}
i) From the inclusions $T _{(m)}  \subset D ^{(m)} [\frac{1}{t}]  \subset  \widetilde{D} ^{(m)}$, we get  $T _{(m)} ^\dag  \subset (D ^{(m)} [\frac{1}{t}] )^\dag \subset  (\widetilde{D} ^{(m)})^\dag $.
Let $\underline{k}\in \bbN ^d$. Following \cite[2.2.5.1]{Be1},  $\underline{\partial} ^{<\underline{k}> _{(m)}}$ is a product of elements in  $\{\partial _i ^{[p ^j]}\ |\ i=1,\dots, d\; ;\; j= 0,\dots, m\}$,
we get  $\underline{\partial} ^{<\underline{k}> _{(m)}}  \in T _{(m)} ^\dag $. Since $\widetilde{D} ^{(m)}$ is $\widetilde{A}$-free with $\{ \underline{\partial} ^{<\underline{k}> _{(m)}}, \, \underline{k}\in \bbN ^d\}$ as basis, 
this yields  $\widetilde{D} ^{(m)} \subset  T _{(m)} ^\dag $. Using Theorem \ref{thmSdag}, this implies $(\widetilde{D} ^{(m)}) ^\dag \subset  T _{(m)} ^\dag $. Hence, we get the desired equalities.

ii) The equivalence between the second respective and non respective cases follows from \ref{vp-vpi}. 
Since the check of the first respective case is similar, we reduce to prove the non respective case. 

a)  Using \cite[2.2.5.1]{Be1},  we can check that for any $\underline{k}\in \bbN ^d$ the element $\underline{\partial} ^{<\underline{k}> _{(m)}}$ is a product of elements in 
$\{\partial _i ^{[p ^j]}\ |\ i=1,\dots, d\; ;\; j= 0,\dots, m\}$  and that the number of factors is less than $|\underline{k}|$.
Hence, in the same way as  \ref{bounded-twodescrp}, the hypothesis that there exists a constant $c>0$ such that $z$ can be written in the form \ref{wkcpdscplevelm} with $v _{\pi} ( a _{\underline{i},k}) \geq \frac{k+|\underline{i}|}{c} -1$
implies  $z$ is $(c,c)$-bounded as $A$-ring with respect to the elements of  $T _{(m)}$ and then  this yields  $z\in  T _{(m)} ^\dag$. 

b) Conversely, suppose
$z \in (\widetilde{D} ^{(m)}) ^\dag \subset \Gamma (\fX _\eta, \widehat{\cD} ^{(m)}  _{\fX _\eta/\eta})$.
We can write uniquely 
$z$ in the form
$z = \sum _{\underline{i}\in \bbN ^d}  b _{\underline{i}} \underline{\partial} ^{<\underline{i}> _{(m)}}$, with
$b _{\underline{i}}\in A \{ \frac{1}{t}\}$ converging to $0$.  
Following \ref{(c,c-i)bounded}, 
we have to check that there exists some $c >0$ such that 
$b _{\underline{i}}$ are $(c , c-|\underline{i}|)$-bounded
(as $A$-algebra with respect to $\frac{1}{t}$), 
for any $\underline{i} \in \bbN ^d$.

1) Let $\fU _\bullet = (\fU _0, \dots, \fU _m)$ be an open cover of $\fX$ by principal open subsets of $\fX$.
Following  Lemma \ref{(c,c-i)bounded}, if $b _{\underline{i}} |\fU _j$ are $(c _0, c _0 -|\underline{i}|)$-bounded for any $j = 0,\dots, m$,
then using Remark \ref{remoflem-vanishingCech},  we get that  $b _{\underline{i}} $ are  $(2c _0, 2c _0 -|\underline{i}|)$-bounded.
Hence, we can suppose that 
$\fX ^\sharp$ 
has logarithmic local coordinates.

2) 
Let $c _0\in \R _{>0}$, $d _0 \in \R$ satisfying the property of Lemma \ref{comp-nota}. We can suppose $d _0 \geq 1$.   Choose an order $\delta _1, \dots, \delta _n$ of the elements of 
$\{\underline{\partial} ^{[\underline{i}]}\ |\ |\underline{i}|\leq p ^{m}\}$. Put $T ' _{(m)}:= \{\frac{1}{t}\} \cup \{\underline{\partial} ^{[\underline{k}]}\ |\ |\underline{k}|\leq p ^{m}\}$.
Since $T _{(m)}\subset T ' _{(m)} \subset \widetilde{D} ^{(m)}$, then it follows from 1) the equality $ (T ' _{(m)} )^\dag=(\widetilde{D} ^{(m)}) ^\dag$. Hence $z \in (T' _{(m)} )^\dag$ and we can write
\begin{equation} \label{dfn-wc=cons1}
z= \sum _{i=0} ^{\infty} \pi ^{i}P _i (\frac{1}{t}, \delta _1, \dots, \delta _n)
\end{equation}
 where  $P _i \in  \mathfrak{V} (\frac{1}{t}, \delta _1,\dots, \delta _n)$ are such that there exists a constant $c >0$  so that 
 $\deg _{(\frac{1}{t}, \delta _1,\dots, \delta _n)} (P _i) \leq c (i+1)$ for any $i$. Since $z$ is fixed, increasing $c _0$ if necessary,  we can suppose $c$  equal to $c _0$. 
By definition, since $t$ is in the center of $\widetilde{D} ^{(m)}$, $P _i (\frac{1}{t}, \delta _1, \dots, \delta _n)$is a finite sum of elements of the form  $\frac{a _0}{t ^n}  \delta _{i _1}a _1 \delta _{i _2} a _2 \dots \delta _{i _N} a _N$,
where $a _j \in A$ and $n +N\leq c _{0} (i + 1)$. Set $\widetilde{A}:= A _{[t]}=\Gamma (\fX, \widetilde{\cO} _{\fX})$. Using the formulas \ref{(c,d)bounded-ring} and  \ref{partialb-bounded},  since $a _j $ are $(c _0, 0)$-bounded,  we check that 
$\delta _{i _1}a _1 \delta _{i _2} a _2 \dots \delta _{i _N} a _N = \sum _{\underline{j} \in \bbN ^r, |\underline{j}| \leq N} \beta _{\underline{j}} \underline{\delta} ^{\underline{j}} $,
where  $\underline{\delta} ^{\underline{j}} := \delta _1^{j _1} \cdots \delta _r^{j _r}$ and  $\beta _{\underline{j}} \in \widetilde{A}$ are $(c _0, 2 N d _0)$-bounded (and therefore
$\frac{a _0}{t ^n}  \beta _{\underline{j}} \in \widetilde{A}$ are $(c _0, 2 N d _0+n)$-bounded).  
Hence, since $n+N \leq  c _{0} (i + 1)$, then $2 N d _0 +n \leq 2d _0 ( n +N)  \leq 2 d _0 c _{0} (i + 1)$. Hence, we can write 
$$P _i (\frac{1}{t}, \delta _1, \dots, \delta _n)= \sum _{ \underline{j} \in \bbN ^r,  |\underline{j}| \leq c _0 (i +1)} b _{i, \underline{j}}  \underline{\delta} ^{\underline{j}} ,$$
where  $b _{i, \underline{j}} \in \widetilde{A}$ are $(c _0, 2  c _{0} (i + 1) d _0)$-bounded.  Hence, we can write  $b _{i, \underline{j}}=  \sum _{n=0} ^{+\infty}  \pi ^{n} \frac{a _{i, \underline{j},n} }    {t ^{\alpha _{i, \underline{j} ,n}} } $
where  $a _{i, \underline{j},n}  \in A$, and $(\alpha _{i, \underline{j},n}) _{n\in \bbN}$ is an increasing sequence of integers such that  $\alpha _{i, \underline{j} ,n} 
\leq c _0 n +2 c _{0} (i + 1) d _0 $. 
Hence,
$$z= \sum _{i=0} ^{\infty} \pi ^{i} P _i (\frac{1}{t}, \delta _1, \dots, \delta _n)=
\sum _{i=0} ^{\infty} 
\sum _{\underset{ |\underline{j}| \leq c _0 (i +1)}{ \underline{j} \in \bbN ^r}}
\sum _{n=0} ^{+\infty}  
\pi ^{i+n} \frac{a _{i, \underline{j},n} }    {t ^{\alpha _{i, \underline{j} ,n}} }
\underline{\delta} ^{\underline{j}} $$
Put 
$Q _h := \sum _{i+n=h} \sum _{ |\underline{j}| \leq c _0 (i +1)}
 \frac{a _{i, \underline{j},n} }    {t ^{\alpha _{i, \underline{j} ,n}} }
\underline{\delta} ^{\underline{j}} $.
The terms in this sum satisfy
$\alpha _{i, \underline{j} ,n} + |\underline{j}| 
\leq 
(c _0 n +2  c _{0} (i + 1) d _0) 
+
c _0 (i+1) 
=
c _0 
(n
+
(2  d _0  + 1 )(i+1)) 
\leq 
c _{0}
(1+2 d _0) (h+1)$.
Hence, 
 $\deg _{(\frac{1}{t}, \delta _1,\dots, \delta _n)} (Q _h) \leq 
c _{0}
(1+2 d _0) (h+1)$
for any $h$.
Since 
$z= \sum _{h=0} ^{\infty} \pi ^{h} Q _h (\frac{1}{t}, \delta _1, \dots, \delta _n)$,
since 
we can write 
$\underline{\delta} ^{\underline{j}}
= 
\sum _{\underline{l}\leq p ^{m} \underline{j} }
\lambda _{\underline{j}, \underline{l}} 
\underline{\partial} ^{<\underline{l}> _{(m)}}$, 
with 
$\lambda _{\underline{j}, \underline{l}}  \in \Z _{(p)}$
(see the formula \cite[2.2.5.1]{Be1}),
then
$z$ satisfies the conditions of \ref{wkcpdscplevelm}.

iii) The last statement is a consequence of  2) and of Lemma \ref{(c,c-i)bounded}.

\end{proof}

\begin{empt}\label{wkcpdscplevelmbis-empt}
With notation \ref{prop-wkcpdscplevelm}, 
we get a filtration on $D ^{(m)}$ given by the order: 
$D _n ^{(m)} := \Gamma (\fX, \cD ^{(m)} _{\fX^\sharp/\fS^\sharp,n})$.
 This yields a filtration on  $D := D ^{(m)} [U]$, the polynomial ring in the commuting variable $U$ over $D  ^{(m)}$ (see \ref{ntn-ord-D-degbis})
 given by the $\cV$-submodule: $D _n := \oplus _{i+j=n} D _{i}^{(m)}  U^j$. We get $\gr D = (\gr D ^{(m)} ) [U]$.
Following \ref{ex-grD-A(m)}, the ring $\gr  D ^{(m)}$ is commutative and noetherian. Hence, so is $\gr D $.

From the $D ^{(m)}$-linear surjective map $\theta\colon D \to  D ^{(m)} _t $ given by $U\mapsto \frac{1}{t}$, we get an exhaustive filtration $D ^{(m)} _{t,n}:= \theta ( D _n)$ on $D ^{(m)} _{t}$.
Since $\gr \theta \colon \gr D \to  \gr D ^{(m)} _t $ is an epimorphism, then $\gr D ^{(m)} _t $ is noetherian.

Following \ref{wkcpdscplevelm},  for any element $z \in (\widetilde{D} ^{(m)}) ^{\dag }$ there exists $c >0$ such that $z$ can be written of the form
\begin{equation} \label{wkcpdscplevelmbis}
z = \sum _{n\in \bbN} \pi ^n P _n 
\end{equation}
where $P _n  \in D ^{(m)} _{t,c(n+1)}$.
Hence, 
\begin{equation} \label{wkcpdscplevelmter}
\Gamma ( \fX,  \widetilde{\cD} ^{(m)\dag} _{\fX^\sharp/\fS^\sharp}) = (D ^{(m)}  _{t}, D ^{(m)}  _{t, n }) ^\dag.
\end{equation}
\end{empt}

\begin{prop}
\label{noeth-weak completion}
Suppose $m\in \bbN$. We keep the notations and hypotheses of  \ref{prop-wkcpdscplevelm}.
The $A$-ring $\widetilde{D} ^{(m)\dag}$ is (right and left) noetherian. 
\end{prop}

\begin{proof}
With its notation, following \ref{wkcpdscplevelmbis-empt} the ring $\gr D ^{(m)}  _{t} $ is commutative and noetherian and we have a surjective morphism of $D ^{(m)}  _{t}$-rings of the form $D ^{(m)}  _{t}  [[T  ]] ^{\mathrm{oc}} \to \widetilde{D} ^{(m)\dag }$.
Following Lemma \ref{oc-noeth}, $D [[T  ]] ^{\mathrm{oc}} $ is left and right noetherian. Hence, we are done. 
\end{proof}

\begin{prop} \label{Jacobson}
Suppose $m\in \bbN$. With the notation of \ref{prop-wkcpdscplevelm}, let $J  (\widetilde{D} ^{(m)\dag})$ be the Jacobson radical of $\widetilde{D} ^{(m)\dag}$ (see \cite[13]{Isaacs}).

\begin{enumerate}[(a)]
\item \label{Jacobson(a)} Then $p \widetilde{D} ^{(m)\dag} \subset J  (\widetilde{D} ^{(m)\dag}) $.

\item \label{Jacobson(b)} Moreover, any left (resp. right) $\widetilde{D} ^{(m)\dag}$-module $M$ of finite type is separated and any submodule  of such $M$ is a closed subset. 

\item \label{Jacobson(c)} The extension 
$\widetilde{D} ^{(m)\dag} \to \widehat{D} ^{(m)}\{ \frac{1}{t}\}$ is right and left faithfully flat (see  \cite[I.3.1, Definition 1]{bourbaki}).

\end{enumerate}

\end{prop}

\begin{proof}
From Theorem \cite[13.9]{Isaacs}, to check the first assertion, we have to show that for any $P \in\widetilde{D} ^{(m)\dag}$,
$1-\pi P$ is invertible in $\widetilde{D} ^{(m)\dag}$.
With notation \ref{wkcpdscplevelmbis-empt}, there exists $Q \in D ^{(m)}  _{t}  [[T  ]] ^{\mathrm{oc}} $ whose image in $(D ^{(m)}  _{t}, D ^{(m)}  _{t, n }) ^\dag$ is $P$.
Following Lemma \ref{Jacobson-oc}, $1-TQ$ is invertible in $Q \in D ^{(m)}  _{t}  [[T  ]] ^{\mathrm{oc}} $. Hence, $1 -\pi P$ is invertible.
We get the other assertions by using \ref{lem-Jacobson-oc}.
\end{proof}

\begin{coro} \label{cor-Zar-ring}
Suppose $m\in \bbN$.
With the notation \ref{prop-wkcpdscplevelm}, the ring $\widetilde{D} ^{(m)\dag}$ endowed with the $\pi$-adic filtration is a left and right Zariski ring (see definition \cite[II.2.1.1]{ZarFilt}).
\end{coro}

\begin{proof}
By taking the $\pi$-adic filtrations on $\widetilde{D} ^{(m)\dag}$ and $\widehat{D} ^{(m)}\{ \frac{1}{t}\}$,  we get the isomorphism of graded rings 
$\gr \widetilde{D} ^{(m)\dag} \riso \gr \widehat{D} ^{(m)}\{ \frac{1}{t}\} \riso \Gamma (X _\eta, \cD ^{(m)} _{X _\eta/\eta} ) [T]$
(see for instance the proof of \cite[4.4.4]{Be2} for the last isomorphism).  Since $\Gamma (X _\eta, \cD ^{(m)} _{X _\eta/\eta} ) [T]$ is Noetherian, since 
$\widetilde{D} ^{(m)\dag} \to \widehat{D} ^{(m)}\{ \frac{1}{t}\}$ is faithfully flat,  using \cite[II.2.1.2.(4)]{ZarFilt},  we conclude.
\end{proof}

\begin{coro} \label{cor-regular-ring}
Suppose $m\in \bbN$.
With the notation \ref{prop-wkcpdscplevelm}, the ring $\widetilde{D} ^{(m)\dag}$ is left and right regular. 
\end{coro}

\begin{proof}
Following the proof of \ref{cor-Zar-ring}, 
we have the isomorphism of graded rings  $\gr \widetilde{D} ^{(m)\dag}  \riso \Gamma (X _\eta, \cD ^{(m)} _{X _\eta/\eta} ) [T]$. Since $\Gamma (X _\eta, \cD ^{(m)} _{X _\eta/\eta} ) [T]$ is left and right regular 
(and then left and right $\gr$-regular)  then so is $\gr \widetilde{D} ^{(m)\dag}$ (use Proposition \cite[I.7.2.1]{ZarFilt} or \cite[A.I.2.7 or D.VII.9]{GradedRingTheory}).
Hence, since $\widetilde{D} ^{(m)\dag}$ is a left and right Zariski ring, using Theorem \cite[II.3.1.4]{ZarFilt}, 
we can conclude. 
\end{proof}

Let us finish by a further description.
\begin{prop} \label{loc-m2dagbis}
We keep the notations and hypotheses of  \ref{prop-wkcpdscplevelm}. Let $D := \Gamma (\fX, \cD _{\fX^\sharp/\fS^\sharp}  )$,  $\widetilde{D} := \Gamma (\fX, \widetilde{\cD} _{\fX^\sharp/\fS^\sharp} )$,
$T :=\{\frac{1}{t}\} \cup   \{\partial _i^{[p ^j]}\ |\ i= 1,\dots, d \; ; \; j \in \bbN\}$.

\begin{enumerate}[(i)]
\item  We have the equalities  $ \widetilde{D} ^{\dag} =(D [\frac{1}{t}] )^{\dag} = T ^{\dag } = \cup _{m\in \bbN}T _{(m)} ^{\dag} =\cup _{m\in \bbN} (\widetilde{D} ^{(m)\dag}) $.

\item 
Let $z$ be an element of $\widehat{D}\{ \frac{1}{t}\}$,
the $p$-adic completion  of $\widetilde{D}$.
We have 
$z \in \widetilde{D} ^{\dag}$
if and only if we can write $z$ in the form
\begin{equation}
\label{wkcpdscp}
z = \sum _{k\in \bbN}\sum _{\underline{i}\in \bbN ^d} a _{\underline{i},k} \frac{1}{t ^{k}}\underline{\partial} ^{[\underline{i}]}
\
(resp. \, 
z 
= 
\sum _{k\in \bbN}\sum _{\underline{i}\in \bbN ^d}  \frac{1}{t ^{k}}\underline{\partial} ^{[\underline{i}]}a _{\underline{i},k})
\end{equation}
where $a _{\underline{i},k} \in A$ are such that there exists a constant $c>0$ satisfying
$v _{\pi} ( a _{\underline{i},k}) \geq \frac{k+|\underline{i}|}{c} -1$
(resp. $v _{p} ( a _{\underline{i},k}) \geq \frac{k+|\underline{i}|}{c} -1$).

\item \label{loc-m2dagbis-(iii)}
Let $z = \sum _{\underline{i}\in \bbN ^d} 
\widetilde{a} _{\underline{i}}
\underline{\partial} ^{[\underline{i}]} $
(resp. $z = \sum _{\underline{i}\in \bbN ^d} 
\underline{\partial} ^{[\underline{i}]} \widetilde{a} _{\underline{i}} $),
where $\widetilde{a} _{\underline{i}} \in A _{\{t\}}$,
be an element of $\widehat{D}\{ \frac{1}{t}\}$.
We have 
$z \in \widetilde{D} ^{\dag}$
if and only if there exists a constant $c>0$ such that 
$\widetilde{a} _{\underline{i}} \in A _{[t]} (c, c -|\underline{i}|)$.

\end{enumerate}
\end{prop}

\begin{proof}
I) The equality $T = \cup _{m\in \bbN} T _{(m)}$ is obvious. 
Using Corollary \ref{coro-indu-lim}, 
this yields 
$T ^{\dag } = \cup _{m\in \bbN}T _{(m)} ^{\dag}$.
From \ref{prop-wkcpdscplevelm}, 
we get 
$\cup _{m\in \bbN}T _{(m)} ^{\dag} 
=\cup _{m\in \bbN} \widetilde{D} ^{(m)\dag} $.
Since
$\cup _{m\in \bbN} \widetilde{D} ^{(m)}=
\widetilde{D} $, 
using again Corollary \ref{coro-indu-lim},
then we get 
$\widetilde{D} ^{\dag}
=
\cup _{m\in \bbN} \widetilde{D} ^{(m)\dag}$.
Since  $ \widetilde{D} ^{\dag} \supset  (D [\frac{1}{t}] )^{\dag} \supset   T ^{\dag}  = \widetilde{D} ^{\dag}$,  these inclusions are in fact equalities.

II) Let us check the description of $\widetilde{D} ^{\dag}$ of \ref{wkcpdscp}. 
The equivalence between the second respective and non respective cases follows from \ref{vp-vpi}. 
Since the check of the first respective case is similar, we reduce to prove the non respective case. 
Let $z\in \widehat{D}\{ \frac{1}{t}\}$.

1) Suppose we can write 
$z = \sum _{k\in \bbN}\sum _{\underline{i}\in \bbN ^d} a _{\underline{i},k} \frac{1}{t ^{k}}\underline{\partial} ^{[\underline{i}]}$
where $a _{\underline{i},k} \in A$ are such that there exists a constant $c>0$ such that 
$v _{p } ( a _{\underline{i},k}) \geq \frac{k+|\underline{i}|}{c} -1$. Let us prove $z \in  \widetilde{D} ^{\dag}$.

a) For any such fixed $c>0$, we can find $m _0$ such that for any $m\geq m _0$,
$l \geq p ^m$ we have 
\begin{align} \label{formula1}
& \frac{1}{p ^{m}(p-1)} < \frac{1}{2c} ; 
\\ \label{formula2} & \frac{l}{p ^{m}(p-1)} < \frac{l}{c} -1.
\end{align}
Indeed, for $m$ large enough, we have \ref{formula1}.  For such $m$, the formula \ref{formula2} is satisfied for  $l \geq p ^m$ if and only if it is satisfied for $l = p ^m$. 
Increasing $m$ is necessary, we get  $\frac{p ^m}{p ^{m}(p-1)} < \frac{p ^m}{c} -1$ and we are done.

b ) For any $\underline{i}= (i _1,\dots, i _d)\in \bbN ^d$, put $q ^{(m _0)} _{\underline{i}} !:= \prod _{j=1} ^{d}[i _j/p ^{m _0}] !$,
$b _{\underline{i},k} := \frac{a _{\underline{i},k} }{q ^{(m _0)} _{\underline{i}} !}$. We check in this step that $b _{\underline{i},k}\in A$. 
First, from \cite[2.4.3.1]{Be1}, we compute  $v _p (q ^{(m _0)} _{\underline{i}} !)  \leq  \frac{\phi ( \underline{i})}{ p ^{m _0} (p-1)}$,
where $$\phi ( \underline{i}) := \sum _{j=1} ^{d} i _j \cdot \min \{[i _j/p ^{m _0}], 1  \}\leq |\underline{i}|.$$  
We have $\phi ( \underline{i})=0$ if and only if  $i _j < p ^{m _0}$ for any integer $j$.
When $\phi ( \underline{i})=0$, $q ^{({m _0})} _{\underline{i}} !=1$ and then $b _{\underline{i},k}\in A$. Suppose now $\phi ( \underline{i})\not=0$.
In that case, we have $k+|\underline{i}| \geq p ^{m _0}$ and then we can use using the formula \ref{formula2} for $l= k+|\underline{i}| $.
Hence, we get $v _{p } ( a _{\underline{i},k}) \geq \frac{k+|\underline{i}|}{c} -1 >  \frac{k+|\underline{i}|}{ p ^{{m _0}} (p-1)} \geq \frac{\phi ( \underline{i})}{ p ^{{m _0}} (p-1)} $. 
This yields  that $v _{p} (	b _{\underline{i},k}	) \geq 0$ and  hence $b _{\underline{i},k}\in A$. 

c) For any $\underline{i}$ and $k$, we compute	 $v _p (	b _{\underline{i},k}	)  \geq \frac{k+|\underline{i}|}{c} -1 - \frac{\phi ( \underline{i})}{ p ^{m _0} (p-1)}
> \frac{2(k+|\underline{i}|) -\phi ( \underline{i})}{2c} -1 \geq  \frac{k+|\underline{i}|}{2c} -1$ (use  \ref{formula1} for the  first inequality).
Since $a _{\underline{i},k} \frac{1}{t ^{k}}\underline{\partial} ^{[\underline{i}]} = b _{\underline{i},k} \frac{1}{t ^{k}}\underline{\partial} ^{<\underline{i}> _{(m _0)}} $, 
using Proposition \ref{prop-wkcpdscplevelm},  we get that  $z \in \widetilde{D} ^{(m _0)\dag}$ and we are done.

2) Conversely, suppose $z \in  \widetilde{D} ^{\dag}$.  Then there exists $m$ such that  $z \in  \widetilde{D} ^{(m)\dag}$. 
Since  $\underline{\partial} ^{<\underline{i}> _{(m)}} $ divides  $\underline{\partial} ^{[\underline{i}]}$, then using the description  \ref{wkcpdscplevelm}, we are done. 

III) This is a consequence of \ref{(c,c-i)bounded} .
\end{proof}

\begin{lem} \label{lem-faithfullflat}
We keep the notations and hypotheses of  \ref{prop-wkcpdscplevelm}. 
The extension  $\Gamma (\fX, \widetilde{\cD} _{\fX^\sharp/\fS^\sharp} ) ^\dag \to  \Gamma (\fX _\eta, \cD  ^\dag _{\fX _\eta/\eta} ) $ is right and left faithfully flat.
\end{lem}

\begin{proof}
We have 
$\Gamma (\fX, \widetilde{\cD} _{\fX^\sharp/\fS^\sharp} )  ^{\dag}
=\cup _{m\in \bbN} (\Gamma (\fX, \widetilde{\cD} ^{(m)} _{\fX^\sharp/\fS^\sharp} ) ) ^\dag $
and 
$\Gamma (\fX _\eta, \cD  ^\dag _{\fX _\eta/\eta} ) 
=
\cup _{m\in \bbN} (\Gamma (\fX, \widehat{\cD} ^{(m)} _{\fX _\eta/\eta} ) )  
$.
Following 
\ref{Jacobson}, 
the extension 
$\Gamma (\fX, \widetilde{\cD} ^{(m)} _{\fX^\sharp/\fS^\sharp} ) ^{\dag}
\to 
\Gamma (\fX, \widehat{\cD} ^{(m)} _{\fX _\eta/\eta} )$
is faithfully flat. 
Since the inductive limits of faithfully flat extensions is faithfully flat, 
we are done.

\end{proof}

\subsubsection{Sheafications}

Let $\fX ^\sharp$ be an object of $\Sm ^\dag _{\fS ^\sharp}$. Let $j \colon \fX _\eta \hookrightarrow \fX ^{\sharp}$ be the open immersion. Let $m\in \bbN$.
\begin{empt} \label{presheafDupDag}
Suppose  $\fX ^\sharp$ is affine and $\fX ^\sharp/\fS ^\sharp$ has logarithmic coordinates.  We denote by $\mathfrak{B}$ the basis of open subsets of  $\fX$ consisting of {\it principal} open subsets. 
We get a presheaf $\widetilde{\cD} ^{(m)\dag} _{\fX^\sharp/\fS^\sharp}$ on $\mathfrak{B}$  defined by putting
$\Gamma ( \mathfrak{D} (f) ,  \widetilde{\cD} ^{(m)\dag} _{\fX^\sharp/\fS^\sharp}):=  \Gamma (\mathfrak{D} (f) ,  \widetilde{\cD} ^{(m)} _{\fX^\sharp/\fS^\sharp}) ^\dag $
where $\dag$ on the right term means  the weak completion as $\Gamma ( \mathfrak{D} (f) , \cO _{\fX})$-ring.
\end{empt}

\begin{lem}
\label{lem-vanishingCechD}
With Notation \ref{presheafDupDag}, 
let $\fU$ be an element of $\mathfrak{B}$, i.e. a principal open subset of $\fX$.
Let $\fU _\bullet = (\fU _0, \dots, \fU _m)$ be an open cover of $\fU$ by elements of $\mathfrak{B}$.
Then: 
\begin{enumerate}[(a)]
\item the natural map $\Gamma (\fU, \widetilde{\cD} ^{(m)\dag} _{\fX^\sharp/\fS^\sharp}) 
\to H ^0 ( \fU _\bullet , \widetilde{\cD} ^{(m)\dag} _{\fX^\sharp/\fS^\sharp})$ (see Notation \ref{cech-cohomology}) is bijective ; 
\item $H ^{i} (\fU _\bullet , \widetilde{\cD} ^{(m)\dag} _{\fX^\sharp/\fS^\sharp})= 0$ for all $i >0$.
\end{enumerate}
\end{lem}

\begin{proof}
We can suppose $\fU= \fX$.
We have to check the exactness of the sequence
$0 \to C ^{-1} (\fU _{\bullet}, \widetilde{\cD} ^{(m)\dag} _{\fX^\sharp/\fS^\sharp}) 
\overset{\partial}{\to} C ^{0}  (\fU _{\bullet}, \widetilde{\cD} ^{(m)\dag} _{\fX^\sharp/\fS^\sharp})
\overset{\partial}{\to}\cdots \overset{\partial}{\to} C ^{m} (\fU _{\bullet}, \widetilde{\cD} ^{(m)\dag} _{\fX^\sharp/\fS^\sharp})\to 0$,
where 
$C ^{-1} (\fU _{\bullet}, \widetilde{\cD} ^{(m)\dag} _{\fX^\sharp/\fS^\sharp}) 
:= \Gamma (\fX, \widetilde{\cD} ^{(m)\dag} _{\fX^\sharp/\fS^\sharp}) $, 
$C ^r (\fU _{\bullet}, \widetilde{\cD} ^{(m)\dag} _{\fX^\sharp/\fS^\sharp})
:= \underset{0\leq i _0 < \cdots < i _r \leq m}{\oplus} \Gamma (\fU _{i _0,\dots, i _r}, \widetilde{\cD} ^{(m)\dag} _{\fX^\sharp/\fS^\sharp})$
for any $r \geq 0$.
Let $f _i$ be so that $\mathfrak{D} (f _i) = \fU _i$.
Put $ \widetilde{A} _{i _0,\dots, i _r}:= \widetilde{A} _{< f _{i _0} \cdots f _{i _{r}}>}$.
Let 
$P \in C ^r (\fU _{\bullet}, \widetilde{\cD} ^{(m)\dag} _{\fX^\sharp/\fS^\sharp})$
be a cocycle (i.e. such that $\partial (P) = 0$) and $P _{i _0, \dots,i _r}$ be its component on 
$\Gamma (\fU _{i _0,\dots, i _r}, \widetilde{\cD} ^{(m)\dag} _{\fX^\sharp/\fS^\sharp})$. 
Using \ref{prop-wkcpdscplevelm}.(\ref{prop-wkcpdscplevelm(iii)}).
there exists $c \in \R _+$ large enough such that any components $P _{i _0, \dots,i _r}$  can be uniquely written in the form  
$P _{i _0, \dots,i _r}  = \sum _{\underline{i}\in \bbN ^d}  b _{i _0, \dots,i _r, \underline{i}} \underline{\partial} ^{<\underline{i}> _{(m)}}$, where $b _{i _0, \dots,i _r, \underline{i}}\in \widetilde{A} _{i _0,\dots, i _r}$ are 
$(c, c -|\underline{i}|)$-bounded (as $A _{\{f _{i _0}\cdots f _{i _r}\}}$-algebra with respect to $\frac{1}{t}$).
Put  $b _{\underline{i}}:= (b _{i _0, \dots,i _r, \underline{i}}) \in C ^r:= \underset{0\leq i _0 < \cdots < i _r \leq m}{\oplus} \widetilde{A} _{i _0,\dots, i _r}$.  From the proof of \ref{lem-vanishingCech}, 
using the unicity of the writing in the form $P _{i _0, \dots,i _r}  = \sum _{\underline{i}\in \bbN ^d}  b _{i _0, \dots,i _r, \underline{i}} \underline{\partial} ^{<\underline{i}> _{(m)}}$,
we check there exists  $c _{\underline{i}}= (c _{i _0, \dots,i _{r-1}, \underline{i}}) \in C ^{r-1}$ such that $\partial (c _{\underline{i}}) = b _{\underline{i}}$. Moreover, following the remark \ref{remoflem-vanishingCech}, 
since $b _{\underline{i}}$ is $(c, c -|\underline{i}|)$-bounded then we can choose $c _{\underline{i}}$ so that $c _{\underline{i}}$ is $(2c, 2c -|\underline{i}|)$-bounded.
Hence, we get $Q _{i _0, \dots,i _{r-1}}  := \sum _{\underline{i}\in \bbN ^d}  c _{i _0, \dots,i _{r-1}, \underline{i}} \underline{\partial} ^{<\underline{i}> _{(m)}} \in \Gamma (\fU _{i _0,\dots, i _{r-1}}, \widetilde{\cD} ^{(m)\dag} _{\fX^\sharp/\fS^\sharp})$.
Putting $Q:= (Q _{i _0, \dots,i _{r-1}} ) \in C ^{r-1} (\fU _{\bullet}, \widetilde{\cD} ^{(m)\dag} _{\fX^\sharp/\fS^\sharp})$, we have $\partial (Q) =P$. 
\end{proof}

\begin{thm}
\label{Ddagissheaf}
We keep Notation and hypotheses \ref{presheafDupDag}. 
Put $\widetilde{D} ^{(m)\dag} _{\fX^\sharp/\fS^\sharp}:= \Gamma (\fX, \widetilde{\cD} ^{(m)\dag} _{\fX^\sharp/\fS^\sharp})$. 
Let $M$ be a $\widetilde{D} ^{(m)\dag} _{\fX^\sharp/\fS^\sharp}$-module of finite type. 
The presheaf $M ^\Delta$ on $\mathfrak{B}$ defined by setting 
$$\Gamma ( \mathfrak{D} (f) ,  M ^\Delta):=  \Gamma (\mathfrak{D} (f) ,  \widetilde{\cD} ^{(m)\dag} _{\fX^\sharp/\fS^\sharp}) \otimes _{\widetilde{D} ^{(m)\dag} _{\fX^\sharp/\fS^\sharp}} M$$
is in fact a sheaf. Moreover,  for every principal open subset $\fU $ of $\fX$ and $i>0$, we have $H ^{i} (\fU,M ^\Delta)=0$.
\end{thm}

\begin{proof}
Since $\widetilde{D} ^{(m)\dag} _{\fX^\sharp/\fS^\sharp}$ is left and right regular
(see \ref{cor-regular-ring}),
we can use the same arguments as for Theorem \ref{sheaf-ThmB}.
\end{proof}

\begin{lem}
\label{affinesubsetD}
We keep Notation and hypotheses \ref{presheafDupDag}.
Let $M$ be a $\widetilde{D} ^{(m)\dag} _{\fX^\sharp/\fS^\sharp}$-module of finite type. 
Let $\fU$  be an affine open set of $\fX$. 
Then $\Gamma (\fU, M ^{\Delta} )= 
\Gamma (\fU ,  \widetilde{\cD} ^{(m)} _{\fX^\sharp/\fS^\sharp}) ^\dag 
\otimes _{\widetilde{D} ^{(m)\dag} _{\fX^\sharp/\fS^\sharp}}M$,
where $\Gamma (\fU ,  \widetilde{\cD} ^{(m)} _{\fX^\sharp/\fS^\sharp}) ^\dag $ means 
the weak completion as $\Gamma ( \fU, \cO _{\fX})$-ring of $\Gamma (\fU ,  \widetilde{\cD} ^{(m)} _{\fX^\sharp/\fS^\sharp})$.
\end{lem}

\begin{proof}
In the same way as \ref{affinesubset},
this is a consequence of \ref{Ddagissheaf}.
\end{proof}

\begin{ntn} \label{ntnB'Ddagm}
We denote by $\mathfrak{B}'$ the basis of opens $\fU ^\sharp$ of  $\fX ^\sharp$ such that $\fU ^\sharp$ is affine and $\fU ^\sharp/\fS ^\sharp$ has logarithmic coordinates.
We get a presheaf $\widetilde{\cD} ^{(m)\dag} _{\fX^\sharp/\fS^\sharp}$ on $\mathfrak{B}'$  defined by putting
$$\Gamma ( \fU ,  \widetilde{\cD} ^{(m)\dag} _{\fX^\sharp/\fS^\sharp}):=  \Gamma (\fU ,  \widetilde{\cD} ^{(m)} _{\fX^\sharp/\fS^\sharp}) ^\dag $$ 
where $ \Gamma (\fU ,  \widetilde{\cD} ^{(m)} _{\fX^\sharp/\fS^\sharp}) ^\dag $ is the weak completion as $\Gamma ( \fU, \cO _{\fX})$-ring of $ \Gamma (\fU ,  \widetilde{\cD} ^{(m)} _{\fX^\sharp/\fS^\sharp})$.

Remark that the presheaf $\widetilde{\cD} ^{(m)\dag} _{\fX^\sharp/\fS^\sharp}$ does not depend on the log structure $M _{\fX}$ on $\fX$ such that $\fX ^\sharp = (\fX , M _{\fX})$ is log smooth over $\fS ^{\sharp}$
and might simply be denoted by $\widetilde{\cD} ^{(m)\dag} _{\fX/\fS}$.
\end{ntn}

\begin{prop} \label{Dmdagisasheaf}
The presheaf  $\widetilde{\cD} ^{(m)\dag} _{\fX^\sharp/\fS^\sharp}$ defined above in \ref{ntnB'Ddagm} is in fact a sheaf. Moreover, we have the equality:
$$M ^\Delta= \widetilde{\cD} ^{(m)\dag} _{\fX^\sharp/\fS^\sharp} \otimes _{\widetilde{D} ^{(m)\dag} _{\fX^\sharp/\fS^\sharp}} M.$$
\end{prop}

\begin{proof}
This is a consequence of \ref{affinesubsetD}.
\end{proof}

\begin{prop} \label{123-coherent}
We keep notation of \ref{ntnB'Ddagm}.
\begin{enumerate}[(a)]
\item For any $\fU \in \mathfrak{B}'$,  $\Gamma (\fU, \widetilde{\cD} ^{(m)\dag} _{\fX^\sharp/\fS^\sharp})$ is right and left Noetherian; 
\item \label{123-coherent(b)} For any $\fU ,\fU'\in \mathfrak{B}'$, with $\fU' \subset \fU$, the homomorphism  $\Gamma (\fU, \widetilde{\cD} ^{(m)\dag} _{\fX^\sharp/\fS^\sharp})   \to   \Gamma (\fU', \widetilde{\cD} ^{(m)\dag} _{\fX^\sharp/\fS^\sharp})$  is right and left flat; 
\item The sheaf $\widetilde{\cD} ^{(m)\dag} _{\fX^\sharp/\fS^\sharp}$ is right and left coherent.
\end{enumerate}
\end{prop}

\begin{proof}
The first assertion was checked in \ref{noeth-weak completion}. From \ref{Jacobson}.(\ref{Jacobson(c)}),  the homomorphisms 
$\Gamma (\fU, \widetilde{\cD} ^{(m)\dag} _{\fX^\sharp/\fS^\sharp}) \to  \Gamma (\fU _\eta, \widehat{\cD} ^{(m)} _{\fX _\eta/\eta})$ is faithfully flat. 
Since  $\Gamma (\fU _\eta, \widehat{\cD} ^{(m)} _{\fX _\eta/\eta})  \to  \Gamma (\fU '_\eta, \widehat{\cD} ^{(m)} _{\fX _\eta/\eta})$ is flat, we get the second assertion.  The last one is a consequence of Proposition \cite[3.1.1]{Be1}.
\end{proof}

\begin{dfn}
We put $\widetilde{\cD} ^{\dag} _{\fX^\sharp/\fS^\sharp}  = \underrightarrow{\lim} _m \widetilde{\cD} ^{(m)\dag} _{\fX^\sharp/\fS^\sharp}$, where the inductive limit is taken in the category of sheaves. 
\end{dfn}

\begin{prop} \label{rem-QcapY=BXdag}
We have 
$\widetilde{\cD} ^{(m)\dag} _{\fX ^\sharp/\fS ^{\sharp}} =\widetilde{\cD} ^{(m)\dag} _{\fX ^\sharp/\fS ^{\sharp},\bbQ}  \cap j _{*} \widehat{\cD} ^{(m)} _{\fX _\eta/\eta}$,
$\widetilde{\cD}  ^{\dag}_{\fX ^\sharp/\fS ^{\sharp}}  = \widetilde{\cD}  ^{\dag }_{\fX ^\sharp/\fS ^{\sharp},\bbQ}  \cap j _{*} \cD  ^\dag _{\fX _\eta/\eta}$.
\end{prop}

\begin{proof}
This is a consequence of \ref{rem-QcapY=BX}.
\end{proof}

\begin{coro}
\label{coro}
Let $\fU$ be an affine open subsets of $\fX$ such that $\fU ^\sharp/\fS ^\sharp$ has logarithmic coordinates.
Then  
$$\Gamma ( \fU ,  \widetilde{\cD} ^{\dag} _{\fX^\sharp/\fS^\sharp}) \riso \Gamma (\fU ,  \widetilde{\cD}  _{\fX^\sharp/\fS^\sharp}) ^\dag  ,$$ where the symbol $\dag$ on the right side 
means the weak completion as $\Gamma (\fU, \cO _{\fX})$-ring.
\end{coro}

\begin{proof}
Since $\fU$ is quasi-compact and quasi-separated,  the functor $\Gamma ( \fU ,-)$ commutes with filtrant inductive limits.  Hence, we get
$\Gamma ( \fU ,  \widetilde{\cD} ^{\dag} _{\fX^\sharp/\fS^\sharp}) \riso \underrightarrow{\lim} _m \Gamma ( \fU ,  \widetilde{\cD} ^{(m)\dag} _{\fX^\sharp/\fS^\sharp})
\underset{\ref{Dmdagisasheaf}}{\riso}\underrightarrow{\lim} _m\Gamma (\fU ,  \widetilde{\cD} ^{(m)} _{\fX^\sharp/\fS^\sharp}) ^\dag 
\underset{\ref{coro-indu-lim}}{\riso} \Gamma (\fU ,  \widetilde{\cD}  _{\fX^\sharp/\fS^\sharp}) ^\dag  $.
\end{proof}

\begin{coro}\label{lem-faithfullflat-cor1}
The extension $\widetilde{\cD} _{\fX^\sharp/\fS^\sharp} ^\dag \to j _* \cD  ^\dag _{\fX _\eta/\eta} $is right and left faithfully flat.
\end{coro}

\begin{proof}
This is a consequence of \ref{lem-faithfullflat} and \ref{coro}.
\end{proof}

\begin{prop}\label{DdagVSDhat}
We have the equalities $\widetilde{\cO} _{\fX,\bbQ} =\cO _{\fX}(\hdag X _s) _\bbQ$, and $\widetilde{\cD} ^{\dag} _{\fX ^{\sharp }/\fS ^{\sharp},\bbQ} =\cD ^{\dag} _{\fX ^{\sharp }/\fS ^{\sharp}} (\hdag X _s) _\bbQ$.
We have the inclusions  
$\cB ^{(m)} _{\fX} (X _s)  \widehat{\otimes} _{\cO _{\fX}} \widehat{\cD} ^{(m)} _{\fX ^{\sharp }/\fS ^{\sharp},\bbQ}\subset  \widetilde{\cD} ^{(m+1)\dag} _{\fX ^{\sharp }/\fS ^{\sharp},\bbQ}$ for any $m \in \bbN$.
\end{prop}

\begin{proof}
Since this is local, we can suppose $\fX= \Spf A$ and $\fX ^\sharp/\fS ^\sharp$ has logarithmic coordinates $u  _1, \dots, u  _d \in \Gamma (\fX, M _{\fX ^\sharp})$.
We denote by  $t  _1, \dots, t  _d $ the image of $u  _1, \dots, u  _d $ in $ \Gamma (\fX, \cO _{\fX})$.

I) Let us prove $\widetilde{\cO} _{\fX,\bbQ} =\cO _{\fX}(\hdag X _s) _\bbQ$. 
1) For any $m $, since $\Gamma ( \fX, \cB ^{(m)} _{\fX} (X _s))= A \{ \frac{p}{ t ^{p ^{m+1}}} \} \subset A [ \frac{1}{t} ] ^\dag = \Gamma ( \fX, \widetilde{\cO} _{\fX} )$,
then we get $\cB ^{(m)} _{\fX} (X _s)  \subset  \widetilde{\cO} _{\fX}$. Taking the limit on $m$ we get  $ \cO _{\fX}(\hdag X _s) _\bbQ \subset \widetilde{\cO} _{\fX}$.

2)  We prove $\widetilde{\cO} _{\fX}\subset \cO _{\fX}(\hdag X _s) _\bbQ$.

i) Let $b \in A _{[t]}$. There exists $c >0$ such that we can write $b = \sum _{n\geq 0} \frac{a _n}{ t ^{n}}$, with $a _n \in A$ are such that  $v _{p} (a _n) \geq \frac{n}{c} -1$ for any $n$.
Let $m _0 $ such that $2c< p ^{m _0+1}$. For $n _0$ large enough, we have $v _p (a _n) \geq \frac{2n}{p^{m _0+1}} +1$ for any $n \geq n _0$.

ii) We have $\sum _{n\geq n _0} \frac{a _n}{ t ^{n}}\in A [ \frac{p}{ t ^{p ^{m _0+1}}} ] ^{\dag}\subset A \{ \frac{p}{ t ^{p ^{m _0+1}}} \}=\Gamma ( \fX, \cB ^{(m _0)} _{\fX} (X _s))$.
Indeed, let $q _n:= [n / p ^{m _0+1}]$ and  $r _n := n - p^{m _0+1} q _n$. We get $ \frac{a _n}{ t ^{n}}  = \frac{a _n t ^{p ^{m _0+1}-r _n}}{ p ^{q _n+1}} (\frac{p}{ t ^{p ^{m _0+1}}}) ^{q _n+1}$,
with $\frac{a _n t ^{r _n}}{ p ^{q+1}} \in A$ and  $v _p (\frac{a _n t ^{r _n}}{ p ^{q+1}}) \geq \frac{n}{p^{m _0+1}}$.

iii)  Since $\sum _{n<n _0} \frac{a _n}{ t ^{n}} \in A _t \subset \Gamma ( \fX, \cB ^{(m _0)} _{\fX} (X _s) _\bbQ)$, then 
it follows from ii) that $b \in \Gamma ( \fX, \cB ^{(m)} _{\fX} (X _s) _\bbQ)$. Hence, $\widetilde{\cO} _{\fX}\subset  \cup _{m\in \bbN}\cB ^{(m)} _{\fX} (X _s) _\bbQ =\cO _{\fX}(\hdag X _s) _\bbQ$.

3) From 1) and 2), we get  $\widetilde{\cO} _{\fX,\bbQ} = \cO _{\fX}(\hdag X _s) _\bbQ$.

II) We now prove $\widetilde{\cD} ^{\dag} _{\fX ^{\sharp }/\fS ^{\sharp}}\subset \cD ^{\dag} _{\fX ^{\sharp }/\fS ^{\sharp}}(\hdag X _s) _\bbQ$.

i)  Let  $P \in  \Gamma (\fX, \widetilde{\cD} ^{\dag} _{\fX ^{\sharp}/\fS ^{\sharp}})$. For $m$ large enough, 
$P \in  \Gamma (\fX, \widetilde{\cD} ^{(m)\dag} _{\fX ^{\sharp}/\fS ^{\sharp}})$. Following \ref{wkcpdscplevelm}, we can write 
$P = \sum _{\underline{k}\in \bbN ^d} \sum _{l\in \bbN} a _{\underline{k},l} \frac{1}{t ^{l}}\underline{\partial} ^{<\underline{k}> _{(m)}}$,
where $a _{\underline{k},l} \in A$ are such that there exists a constant $c>0$ satisfying
$v _{p} ( a _{\underline{k},l}) \geq \frac{l+|\underline{k}|}{c} -1$,
for any $\underline{k} \in \bbN ^d$
and $l \in \bbN$.
Increasing $m \geq 0$ if necessary, there exists $N$ large enough such that for any $l+|\underline{k}| \geq N$, we have  $\frac{l+|\underline{k}|}{c} -1 \geq \frac{2(l+|\underline{k}|)}{p^{m+1}} +1$.
Hence, for any $|\underline{k}| \geq N$, we have$v _{p} ( a _{\underline{k},l}) \geq \frac{2l}{p^{m+1}} +1$.
As for the step 1.2.(ii), this yields $\sum _{l\in \bbN}a _{\underline{k},l} \frac{1}{t ^{l}} \in \Gamma ( \fX, \cB ^{(m)} _{\fX} (X _s))$ for any $|\underline{k}| \geq N$.
Hence, 
$$\sum _{|\underline{k}|\geq N} \sum _{l\in \bbN} a _{\underline{k},l} \frac{1}{t ^{l}}\underline{\partial} ^{<\underline{k}> _{(m)}}
\in  \Gamma (\fX, \cB ^{(m)} _{\fX} (X _s)  \widehat{\otimes} _{\cO _{\fX}} \widehat{\cD} ^{(m)} _{\fX ^{\sharp }/\fS ^{\sharp}})
\subset \Gamma (\fX,\cD ^{\dag} _{\fX ^{\sharp }/\fS ^{\sharp}} (\hdag X _s) _\bbQ).$$

ii)  We get from I.2) the first inclusion  $\widetilde{\cD} ^{(m)} _{\fX ^{\sharp }/\fS ^{\sharp}} \subset \cO _{\fX}(\hdag X _s) _\bbQ \otimes _{\cO _{\fX,\bbQ}} \cD  _{\fX ^{\sharp }/\fS ^{\sharp},\bbQ}
\subset \cD ^{\dag} _{\fX ^{\sharp }/\fS ^{\sharp}}(\hdag X _s) _\bbQ$.
Since $\sum _{|\underline{k}|< N} \sum _{l\in \bbN} a _{\underline{k},l} \frac{1}{t ^{l}}\underline{\partial} ^{<\underline{k}> _{(m)}} \in  \Gamma ( \fX, \widetilde{\cD} ^{(m)} _{\fX ^{\sharp }/\fS ^{\sharp}} )$, this implies 
$\sum _{|\underline{k}|< N} \sum _{l\in \bbN} a _{\underline{k},l} \frac{1}{t ^{l}}\underline{\partial} ^{<\underline{k}> _{(m)}} \in  \Gamma ( \fX, \cD ^{\dag} _{\fX ^{\sharp }/\fS ^{\sharp}}(\hdag X _s) _\bbQ)$.

iii) From i) and ii), we obtain $P \in  \Gamma ( \fX, \cD ^{\dag} _{\fX ^{\sharp }/\fS ^{\sharp}}(\hdag X _s) _\bbQ)$.

III) It remains to check $\cB ^{(m)} _{\fX} (X _s)  \widehat{\otimes} _{\cO _{\fX}}\widehat{\cD} ^{(m)} _{\fX ^{\sharp }/\fS ^{\sharp},\bbQ}\subset \widetilde{\cD} ^{(m+1)\dag} _{\fX ^{\sharp }/\fS ^{\sharp},\bbQ}$ 
(indeed, this yields the inverse inclusion of the part II). Let  $P \in \Gamma (\fX, \cB ^{(m)} _{\fX} (X _s)  \widehat{\otimes} _{\cO _{\fX}} \widehat{\cD} ^{(m)} _{\fX ^{\sharp }/\fS ^{\sharp}}$. 
We can write $P = \sum _{\underline{k}\in \bbN ^d} b _{\underline{k}} \underline{\partial} ^{<\underline{k}> _{(m)}}$, with  $b _{\underline{k}} \in \Gamma (\fX,\cB ^{(m)} _{\fX} (X _s) )$ 
converging to $0$ when $|\underline{k}|$ goes to infinity. 
Since $b _{\underline{k}} \in \Gamma (\fX,\cB ^{(m)} _{\fX} (X _s) )$, then we can write $b _{\underline{k}} =  \sum _{l\in \bbN} a _{\underline{k},l} \frac{1}{t ^{l}}$,
 with $a _{\underline{k},l} \in A$ are such that $v _p (a _{\underline{k},l}) \geq [ l/p ^{m+1}] +1 \geq \frac{l}{ p^{m+1}}$.

For any $\underline{k}= (k _1,\dots, k _d)\in \bbN ^d$, put $\lambda _{\underline{k}}:= \prod _{j=1} ^{d}\frac{[k _j/p ^{m}] !}{[k _j/p ^{m +1}]!}$.
We have checked in the proof of \ref{lem-dagDDm},
for $c>0$ large enough, we have 
$v _p (\lambda _{\underline{k}}) 
\geq 
 \frac{|\underline{k}|}{c} -1$.
 We can suppose $c \geq p^{m+1}$.
We get
$P = \sum _{\underline{k}\in \bbN ^d} \sum _{l\in \bbN} a ' _{\underline{k},l} \frac{1}{t ^{l}}\underline{\partial} ^{<\underline{k}> _{(m+1)}}$, where
$a '_{\underline{k},l} := a _{\underline{k},l} \lambda _{\underline{k}}\in A$
are such that  
$v _{p } ( a '_{\underline{k},l}) \geq \frac{l+|\underline{k}|}{c} -1$.

\end{proof}

\begin{prop}\label{Jacobsonx}
Let $x \in \fX$. Let $J  (\widetilde{\cD} ^{(m)\dag} _{\fX^\sharp/\fS^\sharp,x})$ be the Jacobson radical of $\widetilde{\cD} ^{(m)\dag} _{\fX^\sharp/\fS^\sharp,x}$.
Then $p\widetilde{\cD} ^{(m)\dag} _{\fX^\sharp/\fS^\sharp,x}\subset J  (\widetilde{\cD} ^{(m)\dag} _{\fX^\sharp/\fS^\sharp,x}) $.
\end{prop}

\begin{proof}
This a consequence of \ref{Jacobson}.(\ref{Jacobson(a)}).
\end{proof}

\begin{rem}
This is not clear if $\widetilde{\cD} ^{(m)\dag} _{\fX^\sharp/\fS^\sharp,x}$ is noetherian. Hence, the statement of \ref{Jacobsonx} in the context of $\widetilde{\cD} ^{(m)\dag} _{\fX^\sharp/\fS^\sharp,x}$-modules might be wrong.
\end{rem}

\subsubsection{Theorem of type $A$ for coherent $\widetilde{\cD} ^{(m)\dag} _{\fX ^\sharp/\fS ^{\sharp}}$-modules}

Let $\fX ^\sharp$ be an object of $\Sm ^\dag _{\fS ^\sharp}$. Let $j \colon \fX _\eta \hookrightarrow \fX ^{\sharp}$ be the open immersion.  
W suppose  $\fX ^\sharp$ is affine, and $\fX ^\sharp/\fS ^\sharp$ has logarithmic coordinates. 
Let $D ^{(m)} := \Gamma (\fX, \cD ^{(m)} _{\fX^\sharp/\fS^\sharp} )$, $\widetilde{D} ^{(m)} := \Gamma (\fX, \widetilde{\cD} ^{(m)}  _{\fX^\sharp/\fS^\sharp})$,  $A:= \Gamma (\fX, \cO _{\fX})$.
Fix $t _1,\dots, t _d  \in A _{[t]}$   some overconvergent local coordinates  of $\fX ^\sharp/\fS ^\sharp$ and we keep  notation \ref{localcoord-ovcv}.

\begin{ntn} \label{filtrationDft} 
Let $f\in A$ and $\fU := \mathfrak{D} (f) $. We compute that the left $A _{f}$-module $A _{f} \otimes _{A} D ^{(m)} $ is in fact a subring of the $A _{\{f\}}$-ring
$A _{\{f\}} \otimes _{A}D ^{(m)}  = \Gamma (\fU, \cD ^{(m)} _{\fU^\sharp/\fS^\sharp} )$.
We denote by  $D ^{(m)} _{f} := A _{f} \otimes _{A} D ^{(m)} $ this subring.
The ring $D ^{(m)}  _{f}$ is endowed with the filtration by the order $D ^{(m)}  _{fn}:=A _{f} \otimes _{A} \Gamma (\fX, \cD ^{(m)} _{\fX^\sharp/\fS^\sharp,n} )$, for any $n \in \bbN$.
We have $\gr D ^{(m)}  _{f} = A _{f} \otimes _{A} \gr D ^{(m)}  $.
Following \ref{ex-grD-A(m)}, $\gr D ^{(m)} $ is a commutative noetherian $A _f$-algebra. Hence so is $\gr D ^{(m)}  _{f} $.

Since $t$ is in the center of the $A _{\{f\}}$-ring $A _{\{f\}} \otimes _{A}D ^{(m)}$, then  $(A _{\{f\}} \otimes _{A}D ^{(m)} ) _t$ is a ring.
We compute similarly that the left $A _{ft}$-module $A _{ft} \otimes _{A} D ^{(m)} $ is in fact a subring of  $(A _{\{f\}} \otimes _{A}D ^{(m)} ) _t$. We denote by  $D ^{(m)} _{ft} := A _{ft} \otimes _{A} D ^{(m)} $ this subring.

We get a filtration on $D ^{(m)}  _{f} [T] $ by putting $(D ^{(m)}  _{f} [T]  ) _n := \oplus _{n _1 + n_2 = n}  D ^{(m)}  _{f, n _1} T ^{n _2}$. We have $\gr (D ^{(m)}  _{f} [T] ) = (\gr D ^{(m)} _f) [T]$.
Hence, $\gr D ^{(m)}  _{f} [T] $ is commutative and noetherian. Consider the canonical morphism of rings  $\alpha \colon D ^{(m)}  _{f} [T]  \to D ^{(m)}  _{t f} $ which sends 
$\sum P _k T ^k $ to $\sum P _k \frac{1}{t ^k}  $. We get the exhaustive filtration $ D ^{(m)}  _{t f, n }: = \alpha ((D ^{(m)}  _{f} [T]  ) _n)$ on $D ^{(m)}  _{t f, n }$.
For any $P \in D ^{(m)}  _{t f}\setminus \{ 0\}$,  we put  $\ord (P)  := \min \{ n \in \bbN \ ;\ P \in D ^{(m)}  _{t f, n }\} $.
Since $\gr \alpha \colon \gr D ^{(m)}  _{f} [T]  \to  \gr D ^{(m)}  _{t f} $ is an epimorphism, then $\gr D ^{(m)}  _{t f}$ is noetherian and we can therefore apply the section \ref{ntn-ord-D-deg}.
\end{ntn}

\begin{lem} \label{lem-widetildef}
With notation and hypotheses of \ref{filtrationDft}, with notation \ref{ntn-Dndag} we have the equality
$$\Gamma ( \mathfrak{D} (f) ,  \widetilde{\cD} ^{(m)\dag} _{\fX^\sharp/\fS^\sharp}) = (D ^{(m)}  _{t f}, D ^{(m)}  _{t f, n }) ^\dag.$$ 
\end{lem}

\begin{proof}
1) Let us prove the inclusion $\Gamma ( \mathfrak{D} (f) ,  \widetilde{\cD} ^{(m)\dag} _{\fX^\sharp/\fS^\sharp})\subset (D ^{(m)}  _{t f}, D ^{(m)}  _{t f, n }) ^\dag$.
Let $P \in \Gamma ( \mathfrak{D} (f) ,  \widetilde{\cD} ^{(m)\dag} _{\fX^\sharp/\fS^\sharp})$.  Following \ref{wkcpdscplevelm}, we can write 
$P= \sum _{k\in \bbN}\sum _{\underline{i}\in \bbN ^d} \alpha _{\underline{i},k} \frac{1}{t ^{k}}\underline{\partial} ^{<\underline{i}> _{(m)}}$
where $\alpha _{\underline{i},k} \in A _{\{f\}}$ are such that there exists a constant $c>0$ satisfying $v _{\pi} ( \alpha _{\underline{i},k}) \geq \frac{k+|\underline{i}|}{c} -1$, for any $\underline{i} \in \bbN ^d$ and $k \in \bbN$.
Since $\alpha _{\underline{i},k} \in A _{\{f\}}$, then we can write $\alpha _{\underline{i},k} = \sum _{n \geq  0 } \pi ^n \frac{a _{\underline{i},k,n}}{f ^{m _{\underline{i},k,n}}}$,
where $(m _{\underline{i},k,n}) _{n\in \bbN}$ is an increasing sequence,  and where $a _{\underline{i},k,n}\in A$ are such that $v _\pi (a _{\underline{i},k,n})=0$.
Since $v _{\pi} ( \alpha _{\underline{i},k}) \geq \frac{k+|\underline{i}|}{c} -1$, then  $a _{\underline{i},k,n} =0$ if  $n< \frac{k+|\underline{i}|}{c} -1$. Set 
$$P _n :=  \sum _{k\in \bbN, \underline{i}\in \bbN ^d}  \frac{a _{\underline{i},k,n}}{f ^{m _{\underline{i},k,n}}} \frac{1}{t ^{k}}\underline{\partial} ^{<\underline{i}> _{(m)}}
= \sum _{\underset{k + |\underline{i}|\leq c (n+1)}{k\in \bbN, \underline{i}\in \bbN ^d}}  \frac{a _{\underline{i},k,n}}{f ^{m _{\underline{i},k,n}}} \frac{1}{t ^{k}}\underline{\partial} ^{<\underline{i}> _{(m)}}.$$
Hence, $P _n \in D ^{(m)}  _{t f, c (n +1)}$ and $P = \sum _{n\in \bbN} \pi ^n P _n$. This means  $P \in (D ^{(m)}  _{t f}, D ^{(m)}  _{t f, n }) ^\dag$.
Hence, we are done.

2) Conversely, set $D _{\fU}^{(m)}:= \Gamma (\fU, \cD ^{(m)} _{\fU^\sharp/\fS^\sharp} )$. The ring  $D _{\fU,t}^{(m)}$ endowed with the filtration $(D _{\fU,t,n}^{(m)}) _{n\in\bbN}$ 
which is induced by the order filtration (see \ref{wkcpdscplevelmbis-empt}).
Since $(D ^{(m)}  _{t f}, D ^{(m)}  _{t f, n })\subset ( D _{\fU,t}^{(m)},D _{\fU,t,n}^{(m)})$, then we get 
$$(D ^{(m)}  _{t f}, D ^{(m)}  _{t f, n }) ^\dag \subset ( D _{\fU,t}^{(m)},D _{\fU,t,n}^{(m)})^\dag \overset{\ref{wkcpdscplevelmter}}{=}\Gamma ( \fU,  \widetilde{\cD} ^{(m)\dag} _{\fU^\sharp/\fS^\sharp}) 
=\Gamma ( \mathfrak{D} (f) ,  \widetilde{\cD} ^{(m)\dag} _{\fX^\sharp/\fS^\sharp}).$$ 
\end{proof}

\begin{rem}
The formula \ref{wkcpdscplevelmter} is a particular case (take $f =1$) of 
Lemma \ref{lem-widetildef}.
\end{rem}

\begin{lem}
\label{lem-finite-pres-tildeMbis}
Let $\cM $ be a (globally of) finite presentation $\widetilde{\cD} ^{(m)\dag} _{\fX^\sharp/\fS^\sharp}$-module.
Put $M := \Gamma (\fX , \cM)$. Then $M$ is a $\Gamma (\fX, \widetilde{\cD} ^{(m)\dag} _{\fX^\sharp/\fS^\sharp})$-module of finite type and the canonical morphism $M ^{\Delta}\to \cM$ is an isomorphism.
\end{lem}

\begin{proof}
Using  \ref{Ddagissheaf} and \ref{123-coherent}.(\ref{123-coherent(b)}), this is checked in the same way as \ref{lem-finite-pres-tildeM}.
\end{proof}

\begin{thm}\label{thmAcohDtilde}
Let $\cM $ be a $\widetilde{\cD} ^{(m)\dag} _{\fX ^\sharp/\fS ^{\sharp}}$-coherent sheaf on $\fX$. 
Then there exists a finitely generated $\Gamma (\fX, \widetilde{\cD} ^{(m)\dag} _{\fX ^\sharp/\fS ^{\sharp}})$-module $M$ such that $\cM = M ^{\Delta}$.
\end{thm}

\begin{proof}
We may follow the proof of \ref{thmAcoh}. 

0) Using Nakayama's lemma for non commutative ring  (e.g. see \cite[13.12]{Isaacs}) and \ref{Jacobsonx}, in the same way as the step 0) of the proof of \ref{thmAcoh}, we can check it is enough to show that the natural homomorphism
$\Gamma (\fX, \cM) \to \Gamma (\fX, \M/\pi \cM)$ is surjective.

I) {\it We fix some notation.} 
Choose a covering of $\fX= \cup _{i=0} ^{m}\fU _i $ by principal open sets $\fU _i = \mathfrak{D} (f _i)$ such that 
$\cM | \fU _i$ is a finitely presented $\widetilde{\cD} ^{(m)\dag} _{\fX ^\sharp/\fS ^{\sharp}}| \fU _i$-module. 
Using \ref{lem-finite-pres-tildeMbis}, this implies that
$M _i:= \Gamma (\fU _i , \cM )$ is a finitely generated $\Gamma (\fU _i, \widetilde{\cD} ^{(m)\dag} _{\fX ^\sharp/\fS ^{\sharp}})$-module
such that
$\cM | \fU _i
=
M _i ^{\Delta}$.
Set
$D :=\Gamma (\fX, \cD ^{(m)\dag} _{\fX ^\sharp/\fS ^{\sharp}})$,
$\widetilde{D} _i := \Gamma (\fU _i, \widetilde{\cD} ^{(m)\dag} _{\fX ^\sharp/\fS ^{\sharp}})$, 
$\fU _{i,j} := \fU _i \cap \fU _j$, 
$\widetilde{D} _{i,j} := \Gamma (\fU _{i,j}, \widetilde{\cD} ^{(m)\dag} _{\fX ^\sharp/\fS ^{\sharp}})$,
$M _{i,j} :=  \Gamma ( \fU _{i,j},\cM) $.  
For each $i=0,\dots, m$, select generators $\mu _{i,1},\dots, \mu _{i,r} $ for $M_i$ over $\widetilde{D} _{i}$.
Let $\mu _i = \overset{t}{}(\mu _{i,1},\dots, \mu _{i,r})$ be the corresponding column vector with coefficients in $M _i$.

Let $N _{i,j} = (N _{i,j,\alpha,\beta}) _{\alpha \beta}\in M _{r, r} (\widetilde{D} _{i,j})$ such that  $\mu _i | _{\fU _{ij}} = N _{i,j} \mu _j | _{\fU _{ij}}$,
Using \ref{lem-widetildef},  there exist a constant $c \geq 0$ so that we can write  $N _{i,j} =  \sum _{q=0} ^{+\infty} \pi ^{q}  U _{i,j,q}$ where $U _{i,j,q} \in  M _{r, r} (D _{t f _i f _j} ) $ is such that 
$\ord ( U _{i,j,q})  \leq c (q+1)$. Put  $N _{i,j,h} :=  \sum _{q=0} ^{2 ^{h+1}-1}  \pi ^{q}  U _{i,j,q}$.

Let $\tau \in \Gamma (\fX, \cM /\pi \cM)$. Since $\tau | \fU _i\in \Gamma (\fU _i, \cM /\pi \cM) = M _i /\pi M _i$, then there exists $g ^{0, i} \in M _{1, r} (D _{t f _i})$
such that the image of $g ^{0, i} \mu _i \in M _i$ in $M _i /\pi M _i$ is $\tau | \fU _i$. Increasing $c$  if necessary,   we can suppose  that 
$\ord (g ^{0,i} ) \leq 4c$.

We remark that $g ^{0, j} \mu _j |\fU _{ij}$ and $g ^{0, i} \mu _i |\fU _{ij}$ induce the same element of  $\Gamma (\fU _{ij}, \cM /\pi \cM)$. Since $g ^{0, i} \mu _i |\fU _{ij} = g ^{0, i} N _{ij} \mu _j |\fU _{ij}$, then we get $(g ^{0,i}N _{i,j} -g ^{0,j}) \mu _j  \in \pi  M _{ij}$ for any $i,j$.

\medskip

II) Copying the computations of the part (II) of \ref{thmAcoh}, 
we construct by induction on $h \geq 1$ 
for any $i =0,\dots, m$, some elements $g ^{0,i}, \dots, g ^{h-1,i}\in M _{1,r} ( D _{t f _i})$ satisfying the conditions 
for any $i ,j=0,\dots, m$ :
\begin{enumerate}[(1)]
\item $\sum _{s=0} ^{h-1} \left ( 
g ^{s,i}N _{i,j} -g ^{s,j} 
\right ) \mu _j 
\in \pi ^{2 ^h -1} M _{ij}$ ; 
\item $g ^{s,i}\mu _i \in \pi ^{2 ^s -1} M _{i}$ for any $s =0,\dots, h-1$ ;
\item  $ord ( g ^{s,i}) \leq c 2 ^{s+2}$ for any $s =0,\dots, h-1$. 
\end{enumerate}

\medskip

III) From the properties 
II.2 and II.3), 
using Propositions \ref{convergencepre} and \ref{lem-widetildef}, the sum
$\sum _{s=0} ^{\infty} 
g ^{s,i} \mu _i$ converges in $M _{i}$.
We conclude in the same way as the proof of
\ref{thmAcoh}.

\end{proof}

\begin{coro}
\label{cor-thADdagm}
The functors $M \mapsto M ^{\Delta}$ and $\cM \mapsto \Gamma (\fX,\cM)$ 
induce exact quasi-inverse equivalences  between the category of
$\widetilde{\cD} ^{(m)\dag} _{\fX ^\sharp/\fS ^{\sharp}}$-coherent sheaves 
and that of finitely generated 
$\widetilde{D} ^{(m)\dag} _{\fX ^\sharp/\fS ^{\sharp}}$-modules.
\end{coro}

\begin{proof}
The exactness follows from  \ref{Ddagissheaf} and \ref{123-coherent}.(\ref{123-coherent(b)}).
\end{proof}

\subsubsection{Theorem of type $A$ for coherent $\widetilde{\cD} ^{(m)\dag} _{\fX ^\sharp/\fS ^{\sharp},\bbQ}$-modules}

Let $\fX ^\sharp$ be an object of $\Sm ^\dag _{\fS ^\sharp}$.
We suppose $\fX ^\sharp$ is Noetherian.
Let $j \colon \fX _\eta \hookrightarrow \fX ^{\sharp}$ be the open immersion.

\begin{lem} 
[see {\cite[3.4.3]{Be1}} ]
\label{343Be1}
Let $\fX = \fX' \cup \fX ''$ be an open covering of $\fX$. 
Let $\cM '$ be a $p$-torsion-free coherent $\widetilde{\cD} ^{(m)\dag} _{\fX ^{\prime \sharp}/\fS ^{\sharp}}$-module,
$\cM ''$ be a $p$-torsion-free coherent $\widetilde{\cD} ^{(m)\dag} _{\fX ^{\prime \prime \sharp}/\fS ^{\sharp}}$-module,
and $\epsilon 
\colon 
\cM ' _\bbQ | \fX' \cap \fX''
\riso 
\cM '' _\bbQ
| \fX' \cap \fX''$ be an isomorphism of 
coherent $\widetilde{\cD} ^{(m)\dag} _{\fX ^{\prime \prime } \cap \fX ^{\prime \sharp}/\fS ^{\sharp},\bbQ}$-modules.
Then there exists a $p$-torsion-free
coherent $\widetilde{\cD} ^{(m)\dag} _{\fX ^{\sharp}/\fS ^{\sharp}}$-module $\cM$ which is extending $\M'$, and an isomorphism
$\cM _\bbQ | \fX''
\riso 
\cM '' _\bbQ$
extending $\epsilon$.
\end{lem}

\begin{proof}
We can follow the proof of \cite[3.4.3]{Be1} word for word. 
For the reader, let us recall briefly the construction.
Put $\cE ':= \cM ' | \fX' \cap \fX''$, $\E'' := \cM '' | \fX' \cap \fX''$.
Replacing $\epsilon$ by $p ^N \epsilon $ for $N$ large enough is necessary, 
we can suppose 
$\epsilon ( \E' ) \subset \cE ''$.
Let $i$ be an integer so that 
$p ^{i+1} \E'' / \epsilon ( \E' ) =0$.
Let $\cF  _i$ be the image of the morphism of coherent 
$\widetilde{\cD} ^{(m)\dag} _{\fX ^{\prime \prime} \cap \fX ^{\prime \sharp}/\fS ^{\sharp}}$-modules
$\E' \to  \E'' /\pi ^{i+1} \cE ''$.
Since 
$\cF _i$
is killed by $p ^{i+1}$, 
then $\smash{\overline{\E}} ' _i$
is a coherent
$\cD ^{(m)} _{X _i ^{\prime \prime} \cap X _i ^{\prime \sharp}/S _i ^{\sharp}}$-module.
Hence, there exists a coherent $\cD ^{(m)} _{X _i ^{\prime \prime \sharp}/S _i ^{\sharp}}$-submodule $\G _i$
of $\cM '' / p ^{i+1} \cM ''$ such that $\G _i  | \fX' \cap \fX'' =\cF  _i$.
Let $\cN$ be the kernel of the morphism of coherent $\widetilde{\cD} ^{(m)\dag} _{\fX ^{\prime \prime \sharp}/\fS ^{\sharp}}$-modules
$\cM '' \to  (\cM '' / p ^{i+1} \cM '') / \G _i$. Then $\epsilon $ induces the isomorphism
$\epsilon \colon \cM ' | \fX' \cap \fX'' \riso \cN | \fX' \cap \fX'' $. 
Glueing $\M' $ and $\cN$, we get the $p$-torsion-free
coherent $\widetilde{\cD} ^{(m)\dag} _{\fX ^{\prime \prime} \cap \fX ^{\prime \sharp}/\fS ^{\sharp}}$-module $\cM$.
Since $p ^{i+1} \M'' \subset \cN$, then 
$\cN _\bbQ = \cM '' _\bbQ$ and we are done.

\end{proof}

\begin{prop}
\label{345Be1}
We suppose $\fX$ to be noetherian.
Let $\coh (\widetilde{\cD} ^{(m)\dag} _{\fX ^\sharp/\fS ^{\sharp}}) $ 
(resp. $\coh (\widetilde{\cD} ^{(m)\dag} _{\fX ^\sharp/\fS ^{\sharp},\bbQ}) $)
be the category of coherent $\widetilde{\cD} ^{(m)\dag} _{\fX ^\sharp/\fS ^{\sharp}}$-modules
(resp. coherent $\widetilde{\cD} ^{(m)\dag} _{\fX ^\sharp/\fS ^{\sharp},\bbQ}$-modules).
Let $\coh (\widetilde{\cD} ^{(m)\dag} _{\fX ^\sharp/\fS ^{\sharp}}) _\bbQ$ be the category of 
coherent $\widetilde{\cD} ^{(m)\dag} _{\fX ^\sharp/\fS ^{\sharp}}$-modules up to isogeny. 
Then $\coh (\widetilde{\cD} ^{(m)\dag} _{\fX ^\sharp/\fS ^{\sharp}}) _\bbQ$ is equivalent to the category of 
$p$-torsion-free coherent $\widetilde{\cD} ^{(m)\dag} _{\fX ^\sharp/\fS ^{\sharp}}$-modules up to isogeny. 
Moreover, 
the functor 
$\cM \mapsto \cM _\bbQ$ induces an equivalence of categories between 
$\coh (\widetilde{\cD} ^{(m)\dag} _{\fX ^\sharp/\fS ^{\sharp}}) _\bbQ$ 
and 
$\coh (\widetilde{\cD} ^{(m)\dag} _{\fX ^\sharp/\fS ^{\sharp},\bbQ}) $.
\end{prop}

\begin{proof}
Using Lemma \ref{343Be1}, 
we can follow the proof of \cite[3.4.5]{Be1}.
\end{proof}

\begin{coro}
\label{cor-thADdagmQ}
Suppose $\fX$ is affine. Put 
$\widetilde{D} ^{(m)\dag} _{\fX ^\sharp/\fS ^{\sharp},\bbQ}:= 
\Gamma ( \fX, \widetilde{\cD} ^{(m)\dag} _{\fX ^\sharp/\fS ^{\sharp},\bbQ})$.
The functors $M \mapsto \widetilde{\cD} ^{(m)\dag} _{\fX ^\sharp/\fS ^{\sharp},\bbQ} \otimes 
_{\widetilde{D} ^{(m)\dag} _{\fX ^\sharp/\fS ^{\sharp},\bbQ}}
M$ and $\cM \mapsto \Gamma (\fX,\cM)$ 
induce quasi-inverse equivalences  between the category of
$\widetilde{\cD} ^{(m)\dag} _{\fX ^\sharp/\fS ^{\sharp},\bbQ}$-coherent sheaves
and that of finitely generated 
$\widetilde{D} ^{(m)\dag} _{\fX ^\sharp/\fS ^{\sharp},\bbQ}$-modules.
\end{coro}

\begin{proof}
This is a consequence of \ref{cor-thADdagm} and \ref{345Be1}.
\end{proof}

\subsection{Frobenius descent in the context of weakly completion of level $m$}

\subsubsection{The extraordinary pullbacks $f ^\flat$ and $f ^{\sharp}$}

This subsection will be useful later. 
More precisely, in order to prove Lemma \ref{Fflat-omega}, we will need 
\ref{prop-III.8.2}.

\begin{lem}\label{dfn-f-tilde}
Let $f \colon \fY \to \fX$ be a morphism of  $\Sm ^\dag _\eta$.
Then there exists canonically a morphism of locally ringed spaces of the form
$\widetilde{f} \colon  (\fY, \widetilde{\cO} _\fY)  \to  (\fX, \widetilde{\cO} _\fX) $ whose underlying morphism of topological spaces is that of $f$ and making commutative the diagram
\begin{equation}
\notag 
\xymatrix{
{(\fY, \widetilde{\cO} _\fY)  } 
\ar[r] ^-{\widetilde{f}}
\ar[d] ^-{}
& 
{(\fX, \widetilde{\cO} _\fX)} 
\ar[d] ^-{}
\\ 
{\fY} 
\ar[r] ^-{} 
& 
{\fX.} 
}
\end{equation}
\end{lem}

\begin{proof}
Let $\fX' =\Spf A$ (resp. $\fY'=\Spf B$) be an open affine set of $\fX$ (resp. $\fY$) such that  $f ( \fY') \subset \fX'$.  Since $f ( \fY') \subset \fX'$, we get the homomorphism of algebras $A \to B$ and then  $A _t \to B _t$.
Hence, $A _{[t]}   \to  B _{[t]} $, where $A _{[t]} $ is the weak completion of $A _t$ as $A$-algebra, and $B _{[t]} $ is the weak completion of $B _t$ as $B$-algebra.
Hence, with \cite[0.3.5.1]{EGAI}, we are done.
\end{proof}

\begin{prop} \label{D2Dmdagpre}
Let $f \colon \fY \to \fX$ be a finite surjective morphism of  $\Sm ^\dag _\eta$.
\begin{enumerate}[(a)]
\item Then the morphism $\widetilde{f} \colon  (\fY, \widetilde{\cO} _\fY)  \to  (\fX, \widetilde{\cO} _\fX) $ (see \ref{dfn-f-tilde}) is flat, 
i.e. for any $y \in \fY$,  $\widetilde{\cO} _{\fX,f(y)} \to  \widetilde{\cO} _{\fY,y}$ is flat. 
In particular, the functor  $\widetilde{f} ^* = \widetilde{\cO} _{\fY,\bbQ}  \otimes _{f ^{-1} \widetilde{\cO} _{\fX ,\bbQ}}  f ^{-1} -$ is exact.

\item Moreover, the sheaf $f _*\widetilde{\cO} _{\fY}$ is locally free of finite type  $\widetilde{\cO} _{\fX}$-module.
\end{enumerate}
\end{prop}

\begin{proof}
Let us check the first statement.  Since $f _{\eta 0} \colon Y _\eta \to X _\eta$ is a finite and surjective morphism of regular schemes,  then from \cite[4.3.11]{Liu-livre-02} (or frome the more general statement \cite[IV.15.4.2]{EGAIV3})
we get that $f _{\eta 0}$ is flat.  Since $\cO _{\fY _\eta}$ has no $p$-torsion and is (weakly) $p$-adically complete, then 
$f _\eta$ is flat (e.g. use \cite[Lemma 2.1]{MonskyWashnitzer}). Let $y$ be a point of $\fY$ and $x= f (y)$.  We have to check that 
$\widetilde{\cO} _{\fX,x} \to  \widetilde{\cO} _{\fY,y}$ is flat. Since this is local, we can suppose $\fX=\Spf A$, $\fY=\Spf B$  affine. 
Let respectively $\mathfrak{p} _x$ and $\mathfrak{q} _y$ be the prime ideal of $A$ and $B$ corresponding  to $x$ and $y$. 
Both canonical homomorphisms 
$(j _* \cO _{\fX _\eta}) _x =  \underrightarrow{\lim} _{a \not \in  \mathfrak{p} _x} \Gamma (\mathfrak{D} (a), j _* \cO _{\fX _\eta}) 
= \underrightarrow{\lim} _{a \not \in  \mathfrak{p} _x} \Gamma (\mathfrak{D} (a), \cO _{\fX _\eta}) 
\to \underrightarrow{\lim} _{a \not \in  \mathfrak{p} _x} \Gamma (\mathfrak{D} (f _\eta ^* (a), \cO _{\fY _\eta}) 
\to  \underrightarrow{\lim} _{b \not \in  \mathfrak{q} _y}  \Gamma (\mathfrak{D} (b), \cO _{\fY _\eta})   =  (j _* \cO _{\fY _\eta}) _y   $  are flat. 
Using \ref{faithfullyflat1}, this yields that  $\widetilde{\cO} _{\fX,x} \to  \widetilde{\cO} _{\fY,y}$ is flat.

Since the second part of the lemma is local, since $f _\eta$ is flat, then we can suppose $B$ is a free $A$-module of finite type. 
Since $A _{[t]} \to  A _{\{t\}}$ and $B _{[t]} \to  B _{\{t\}}$ are faithfully flat, then we conclude.
\end{proof}

\begin{ntn}
Let $f \colon \fY \to \fX$ be a finite morphism of  $\Sm ^\dag _\eta$. We denote by  $\overline{f} \colon  (\fY, \widetilde{\cO} _\fY)  \to  (\fX, f _* \widetilde{\cO} _\fY) $ the morphism whose composition with 
$(\fX, f _* \widetilde{\cO} _\fY)  \to  (\fX, \widetilde{\cO} _\fX) $ is  $\widetilde{f}$.
\end{ntn}

\begin{lem}
\label{flatnessovf}
Let $f \colon \fY \to \fX$ be a finite morphism of 
$\Sm ^\dag _\eta$.
Then,  $\overline{f}$ is flat, i.e. for any $y \in \fY$, 
$( f _* \widetilde{\cO} _\fY) _{f (y)} \to \widetilde{\cO} _{\fY, y}$ is flat.   
\end{lem}

\begin{proof}
We can suppose 
$\fX = \Spf A$, $\fY= \Spf B$, 
$y= \mathfrak{q}$,
$x:= f(y)= \mathfrak{p}= \phi ^{-1} (\mathfrak{q})$,
where $\phi \colon A \to B$ is the homomorphism corresponding to $f$.
We compute
$( f _* \widetilde{\cO} _\fY) _{f (y)}
= 
\underrightarrow{\lim}
_{a \not \in \mathfrak{p}}
B \{ 1/ \phi (a)\}  _{[t]}$,
and
$\widetilde{\cO} _{\fY, y}
=
\underrightarrow{\lim}
_{b \not \in \mathfrak{q}}
B \{ 1/ b\}  _{[t]}$.
Since 
$\phi ( A \setminus \mathfrak{p})
\subset 
B \setminus \mathfrak{q}$,
since 
$B _{[t]} \to 
B \{ 1/ b\}  _{[t]}$ are flat
(because so is the $p$-adic completion, 
and $B \{ 1/ b\}  _{[t]} 
\to 
B \{ 1/ tb\}  $ are faithfully flat),
then 
$( f _* \widetilde{\cO} _\fY) _{f (y)} \to \widetilde{\cO} _{\fY, y}$ is flat.   
\end{proof}

\begin{lem}
\label{lem-dagABfiniteotimes}
Let $\phi \colon A \to B$ be a finite morphism of 
$p$-torsion-free commutative $\cV [[t]]$-algebras. 
Then 
$(B _t ) ^{\dag/A}=(B _t ) ^{\dag/B}$
and 
$(B _t ) ^{\dag/B}=
B \otimes _{A} (A _t ) ^{\dag/A}$.
\end{lem}

\begin{proof}
We have always
$(B _t ) ^{\dag/A}\subset(B _t ) ^{\dag/B}$.
Since $\phi$ is finite, 
we get the surjectivity from an easy computation.

Since $\phi$ is finite,  $B \otimes _{A}  \widehat{A _t} \to  \widehat{B _t}$ is a morphism of separated complete algebras. Since $B \to B _t$ is flat and $B$ is flat over $\cV$ then $B _t$ has no $p$-torsion.
Since the reduction modulo $\pi$ of  $B \otimes _{A}  \widehat{A _t} \to  \widehat{B _t}$ is an isomorphism,  since $B _t$ has no $p$-torsion, then  $B \otimes _{A}  \widehat{A _t} \to  \widehat{B _t}$  is an isomorphism.

Since $(A _t ) ^{\dag/A} \to \widehat{A _t}$ is faithfully flat, the left vertical arrow of the diagram
\begin{equation} \notag
\xymatrix{
{B \otimes _{A} (A _t ) ^{\dag/A}}  \ar@{^{(}->}[d] ^-{} \ar[r] ^-{} &  {(B _t ) ^{\dag/B}}  \ar@{^{(}->}[d] ^-{}
\\  {B \otimes _{A}  \widehat{A _t}}  \ar[r] ^-{\sim} &   { \widehat{B _t}}      }
\end{equation}
is injective.  Hence, the map  $B \otimes _{A} (A _t ) ^{\dag/A} \to (B _t ) ^{\dag/B}$ is injective.  The surjectivity is an easy computation. 
\end{proof}

\begin{prop}
\label{prop-f*O-fprop}
Let $f \colon \fY \to \fX$ be a finite morphism
$\Sm ^\dag _\eta$. 
\begin{enumerate}[(a)]
\item The sheaf of rings $f _* \widetilde{\cO} _{\fY} $ is $\widetilde{\cO} _{\fX} $-coherent.
\item \label{prop-f*O-fprop(b)} Let $\cM$ be a $f _* \widetilde{\cO} _{\fY} $-module. 
The sheaf $\cM$ is $f _* \widetilde{\cO} _{\fY} $-coherent if and only if it is 
$\widetilde{\cO} _{\fX} $-coherent (via $\widetilde{\cO} _{\fX} \to f _* \widetilde{\cO} _{\fY}$). 
\item If $\fX$ is affine, then the functors 
$\Gamma (\fX, -)$ and 
$f _* \widetilde{\cO} _{\fY} \otimes _{\Gamma ( \fY,  \widetilde{\cO} _{\fY})} -$
induce quasi-inverse equivalence of categories between 
the category of coherent $f _* \widetilde{\cO} _{\fY}$-modules 
and that of coherent $\Gamma ( \fY,  \widetilde{\cO} _{\fY})$-modules.
\item For any coherent $\widetilde{\cO} _{\fY} $-module $\cN$, for any integer $n\geq1$,
$R ^n f _* (\cN) = 0$.
\item \label{prop-f*O-fprop(e)} The exact functors $\overline{f} _*= f_*$
and $\overline{f} ^*$  
induce quasi-inverse equivalence between the category of 
coherent $f _* \widetilde{\cO} _{\fY} $-modules and that of 
coherent $\widetilde{\cO} _{\fY} $-modules.
The functors $\R \overline{f} _*$
and $\overline{f} ^*$  
induce quasi-inverse equivalence between 
$D ^\mathrm{b} _{\mathrm{coh}}
(\widetilde{\cO} _{\fY})$
and
$D ^\mathrm{b} _{\mathrm{coh}}
(f _*\widetilde{\cO} _{\fY})$.

\end{enumerate}

\end{prop}

\begin{proof}
I) Suppose $\fX$ affine. Let us prove the part $3)$ of the proposition. Put  $A := \Gamma (\fX, \cO _{\fX})$, $B := \Gamma (\fY, \cO _{\fY})$, and $\phi\colon A \to B$ the homomorphism corresponding to $f$.
i) We prove in this step that  the canonical morphism  $\widetilde{\cO} _{\fX} \otimes _{A _{[t]}} B _{[t]} \to  f _* \widetilde{\cO} _{\fY}$ is an isomorphism (i.e. we check the part 1) of the proof).
For any $a \in A$,  $\Gamma ( \mathfrak{D} (a), f _* \widetilde{\cO} _{\fY}) =  \Gamma ( \mathfrak{D} (\phi (a)), \widetilde{\cO} _{\fY}) = (B _{\{\phi (a)\}}) _{[t]}$.
Since $B _{\{\phi (a)\}}$ has no $\pi$-torsion, since the homomorphism $A _{\{ a\}} \otimes _A  B \to B _{\{\phi (a)\}}$ is an isomorphism modulo $\pi$,  then this is an isomorphism.
Hence,  $A _{\{ a\}} \to B _{\{\phi (a)\}}$ is finite. This yields  $(A _{\{ a\}}) _{[t]} \otimes _{A } B \riso  (A _{\{ a\}}) _{[t]} \otimes _{A _{\{ a\}}} B _{\{\phi (a)\}}  \underset{\ref{lem-dagABfiniteotimes}}{\riso}  (B _{\{\phi (a)\}} ) _{[t]}$.
Since $\phi$ is finite,  using Lemma \ref{lem-dagABfiniteotimes},  we get  $A _{[t]} \otimes _{A } B \riso  B _{[t]}$
and then  $(A _{\{ a\}}) _{[t]} \otimes _{A } B \riso  (A _{\{ a\}}) _{[t]} \otimes _{A _{[t]}} B _{[t]}$. Hence, the canonical morphism
$(A _{\{ a\}}) _{[t]} \otimes _{A _{[t]}} B _{[t]} \to (B _{\{\phi (a)\}} ) _{[t]}$ is an isomorphism.

Moreover, since $B _{[t]}$ is an $A _{[t]}$-module of finite type, then 
$\widetilde{\cO} _{\fX} \otimes _{A _{[t]}} B _{[t]}$ is $\widetilde{\cO} _{\fX}$-coherent
and 
$\Gamma ( \mathfrak{D} (a), \widetilde{\cO} _{\fX} \otimes _{A _{[t]}} B _{[t]})
=
(A _{\{ a\}}) _{[t]} \otimes _{A _{[t]}} B _{[t]}$.
Hence, by applying the functor
$\Gamma ( \mathfrak{D} (a), -)$
to the morphism
$\widetilde{\cO} _{\fX} \otimes _{A _{[t]}} B _{[t]}
\to 
f _* \widetilde{\cO} _{\fY}$,
we get an isomorphism.

ii) Let $\cM$ be a $f _* \widetilde{\cO} _{\fY} $-module which is  $\widetilde{\cO} _{\fX} $-coherent. Then, the canonical morphism $\widetilde{\cO} _{\fX} \otimes _{A _{[t]}} \Gamma (\fX, \cM) \to  \cM$  is an isomorphism.
Since $\widetilde{\cO} _{\fX} \otimes _{A _{[t]}} B _{[t]} \riso  f _* \widetilde{\cO} _{\fY}$, then we get 
$\widetilde{\cO} _{\fX} \otimes _{A _{[t]}} \Gamma (\fX, \cM) \riso  (\widetilde{\cO} _{\fX} \otimes _{A _{[t]}} B _{[t]})  \otimes  _{B _{[t]}} \Gamma (\fX, \cM)
\riso  f _* \widetilde{\cO} _{\fY} \otimes  _{B _{[t]}} \Gamma (\fX, \cM) $. Hence, the canonical morphism
$f _* \widetilde{\cO} _{\fY} \otimes  _{B _{[t]}} \Gamma (\fX, \cM) \to \cM$ is an isomorphism.

iii) Let $M$ be a coherent $B _{[t]}$-module. Since $\widetilde{\cO} _{\fX} \otimes _{A _{[t]}} B _{[t]} \riso  f _* \widetilde{\cO} _{\fY}$,
we get  $\widetilde{\cO} _{\fX} \otimes _{A _{[t]}} M \riso  (\widetilde{\cO} _{\fX} \otimes _{A _{[t]}} B _{[t]})  \otimes  _{B _{[t]}} M \riso  f _* \widetilde{\cO} _{\fY} \otimes  _{B _{[t]}} M$.
Since $M$ is a coherent $A _{[t]}$-module,  from theorem of type $A$ this yields 
$M \riso  \Gamma (\fX, f _* \widetilde{\cO} _{\fY} \otimes  _{B _{[t]}} M)$. Hence, we are done.

II) Using the step I.i) and Proposition \cite[3.1.1]{Be1},  we get the coherence of the sheaf of ring $f _* \widetilde{\cO} _{\fY} $, i.e. the part 1) of the proposition.
Moreover, from the step I.i),  $f _* \widetilde{\cO} _{\fY} $ is a coherent $\widetilde{\cO} _{\fX}$-module.
Hence, if $\cM$ is a coherent $f _* \widetilde{\cO} _{\fY} $-module,  then $\cM$  is $\widetilde{\cO} _{\fX} $-coherent.
Conversely,  let $\cM$ be a $f _* \widetilde{\cO} _{\fY} $-module which is  $\widetilde{\cO} _{\fX} $-coherent. Then from the step I.ii), 
we check that  $\cM$ is $f _* \widetilde{\cO} _{\fY} $-coherent. Hence, we have checked the part $2)$.

III) Let $\cN$ be a coherent $\widetilde{\cO} _{\fY} $-module. 

i) The part 4) of the proposition is a consequence of \ref{equcatBDelta} and of the fact that  $R ^n f _* (\cN) $ is the sheaf associated to the  presheaf given by $\fU \mapsto H ^n ( f ^{-1} (\fU), \cN)$.

ii) It follows from the part 4) of the proposition  that the functor $f _*$ is exact over the category of coherent $\widetilde{\cO} _{\fY} $-modules.
Since $\overline{f}$ is flat  (see \ref{flatnessovf}), we get the exactness of $\overline{f} ^*$. 

iii) The fact that $f _* (\cN) $ is $f _* \widetilde{\cO} _{\fY} $-coherent is local.
Hence, we reduce to the case where $\fX$ is affine. Using theorem of type $A$ for  coherent $\widetilde{\cO} _{\fY} $-modules, we get $\widetilde{\cO} _{\fY}  \otimes _{\Gamma ( \fY,  \widetilde{\cO} _{\fY})}  \Gamma ( \fY,  \cN) \to \cN$ is an isomorphism.
Using the five Lemma, since both functors are exact on $\cN$,  we get that the canonical morphism
\begin{equation}
\label{prop-f*O-fprop-proof} f _* (\widetilde{\cO} _{\fY} ) \otimes _{\Gamma ( \fY,  \widetilde{\cO} _{\fY})}  \Gamma ( \fY,  \cN) \to  f _* (\widetilde{\cO} _{\fY}  \otimes _{\Gamma ( \fY,  \widetilde{\cO} _{\fY})}  \Gamma ( \fY,  \cN))
\end{equation}
is an isomorphism. Hence, $f _* (\cN) $ is $f _* \widetilde{\cO} _{\fY} $-coherent.

iv) Let $\cM$ be a coherent $f _* \widetilde{\cO} _{\fY} $-module. Then the canonical morphism  $\cM \to \overline{f} _* \circ \overline{f} ^* (\cM)$ is an isomorphism. Indeed, since this is local we reduce to the case where $\fX$ is affine.
In that case,  following the step I.ii), the canonical morphism $f _* \widetilde{\cO} _{\fY} \otimes  _{B _{[t]}} \Gamma (\fX, \cM) \to \cM$ is an isomorphism.
This yields $\overline{f} ^* (\cM) \riso  \widetilde{\cO} _{\fY}\otimes  _{B _{[t]}} \Gamma (\fX, \cM)$. 
This yields the first isomorphism:$\overline{f} _* \circ \overline{f} ^* (\cM) \riso  \overline{f} _* (\widetilde{\cO} _{\fY}\otimes  _{B _{[t]}} \Gamma (\fX, \cM))
\underset{\ref{prop-f*O-fprop-proof}}{\riso}  f _* (\widetilde{\cO} _{\fY} )\otimes  _{B _{[t]}} \Gamma (\fX, \cM) \riso \cM$, the last one was checked at the step I.ii)..

iv) We check similarly that  for any coherent $\widetilde{\cO} _{\fY} $-module $\cN$,  the canonical morphism  $\overline{f} ^* \circ \overline{f} _* (\cN) \to \cN$ is an isomorphism. 

v) The similar statement concerning derived categories is an obvious consequence.
\end{proof}

\begin{ntn}
Let $f \colon \fY \to \fX$ be a finite morphism of  $\Sm ^\dag _\eta$. From \ref{coh-perf}, $D ^\mathrm{b} _{\mathrm{coh}} (\widetilde{\cO} _{\fX}) =D ^\mathrm{b} _{\mathrm{parf}} (\widetilde{\cO} _{\fX})$.
Hence, using \ref{prop-f*O-fprop}.1, we get  $f _* \widetilde{\cO} _\fY \in D ^\mathrm{b} _{\mathrm{parf}} (\widetilde{\cO} _{\fX})$. For any  $\cM \in D ( \widetilde{\cO} _\fX)$, we put
\begin{equation} \label{dfnfflat}
f ^\flat (\cM)  :=  \overline{f} ^*  \R \mathcal{H} om _{\widetilde{\cO} _\fX} (f _* \widetilde{\cO} _\fY , \cM).
\end{equation}
Moreover, since $f _* \widetilde{\cO} _\fY \in D ^\mathrm{b} _{\mathrm{parf}} (\widetilde{\cO} _{\fX})$, then for any $\cM \in  D ^\mathrm{b} _{\mathrm{coh}} (\widetilde{\cO} _{\fX})$, 
we have  $\R \mathcal{H} om _{\widetilde{\cO} _\fX} (f _* \widetilde{\cO} _\fY , \cM) \in D ^\mathrm{b} _{\mathrm{coh}} (\widetilde{\cO} _{\fX})$.
Hence, by \ref{prop-f*O-fprop}.(\ref{prop-f*O-fprop(b)}) this yields $\R \mathcal{H} om _{\widetilde{\cO} _\fX} (f _* \widetilde{\cO} _\fY , \cM) \in D ^\mathrm{b} _{\mathrm{coh}} (f _* \widetilde{\cO} _\fY)$.
Via \ref{prop-f*O-fprop}.(\ref{prop-f*O-fprop(e)}), this yields $f ^\flat (\cM) \in  D ^\mathrm{b} _{\mathrm{coh}} ( \widetilde{\cO} _\fY)$.
\end{ntn}

\begin{empt}Let $f \colon \fY \to \fX$ be a finite morphism of  $\Sm ^\dag _\eta$.
Let $\cM \in  D ^\mathrm{b} _{\mathrm{coh}} (\widetilde{\cO} _{\fX})$. Denoting by $\mathrm{ev _1}$ the evalution at $1$, we get 
\begin{equation} \label{f*fflatadj}
\R f _*  f ^\flat (\cM)  =  \R f _* \overline{f} ^*  \R \mathcal{H} om _{\widetilde{\cO} _\fX} (f _* \widetilde{\cO} _\fY , \cM)
\underset{\ref{prop-f*O-fprop}.(\ref{prop-f*O-fprop(e)})}{\liso} \R \mathcal{H} om _{\widetilde{\cO} _\fX} (f _* \widetilde{\cO} _\fY , \cM) \underset{\mathrm{ev _1}}{\longrightarrow} \M.
\end{equation}
\end{empt}

\begin{lem} \label{II5.5HaRD}
Let $f \colon \fY \to \fX$ be a morphism of  $\Sm ^\dag _\eta$. Let $\cN \in  D ^\mathrm{b} _{\mathrm{coh}} (\widetilde{\cO} _{\fY})$, $\cN ' \in  D ^\mathrm{+}  (\widetilde{\cO} _{\fY})$. The canonical morphism 
\begin{equation}
\R f _* \R \mathcal{H} om _{\widetilde{\cO} _\fY} (\cN, \cN' ) \to \R \mathcal{H} om _{f _* (\widetilde{\cO} _\fY)} (\R f _* (\cN), \R f _* (\cN ')),
\end{equation}
defined as in \cite[II.5.5]{HaRD}, is an isomorphism in  $D ^\mathrm{+}  (f _*\widetilde{\cO} _{\fY})$.
\end{lem}

\begin{proof}
Since this is local and since $\cN \in D ^\mathrm{b} _{\mathrm{parf}} (\widetilde{\cO} _{\fX})$, then we reduce to the case where $\cN =\widetilde{\cO} _{\fY}$, which is clear. 
\end{proof}

\begin{lem} \label{cohf*O=parf}
Let $f \colon \fY \to \fX$ be a finite morphism of  $\Sm ^\dag _\eta$.
\begin{enumerate}[(a)]
\item Then  $\tor \dim f _* \widetilde{\cO} _{\fY} = \dim \fY _\eta$.

Let  $\fU \in \mathfrak{Aff} _{\fX}$, $\fV:= f  ^{-1}(\fU)$ and $g \colon \fV \to \fU$ be the induced morphism. For any coherent $f _* \widetilde{\cO} _{\fV}$-module $\cE$, 
there exists a resolution of length $\leq \dim \fY _\eta$ of $\cE$ by projective coherent
$f _* \widetilde{\cO} _{\fV}$-modules. 
We have $D ^\mathrm{b} _{\mathrm{coh}} (f _* \widetilde{\cO} _{\fY} ) = D ^\mathrm{b} _{\mathrm{parf}} (f _* \widetilde{\cO} _{\fY} )$.

\item Let $\cM \in  D ^+ (\widetilde{\cO} _{\fX})$, and $\cE \in  D ^\mathrm{b} _{\mathrm{coh}} (f _*\widetilde{\cO} _{\fY})$. The canonical morphism 
\begin{equation}
\label{cohf*O=parf(2)} \R \mathcal{H} om _{\widetilde{\cO} _\fX} (\E, \cM ) \to \R \mathcal{H} om _{f _*\widetilde{\cO} _\fY} (\E, \R \mathcal{H} om _{\widetilde{\cO} _\fX} (f _* \widetilde{\cO} _\fY , \cM) ) 
\end{equation}
is an isomorphism.
\end{enumerate}
\end{lem}

\begin{proof}
Using \ref{prop-f*O-fprop}, we check the first statement in the same way as \ref{coh-perf}. We  construct the morphism \ref{cohf*O=parf(2)} by using a resolution of $\cM$ by injective $\widetilde{\cO} _{\fX}$-modules. By using \cite[I.7.1.(iv)]{HaRD}, 
since  $D ^\mathrm{b} _{\mathrm{parf}} (f _*\widetilde{\cO} _{\fY}) = D ^\mathrm{b} _{\mathrm{coh}} (f _*\widetilde{\cO} _{\fY})$, to check this is an isomorphism, we reduce to the case where  $\cE$ is a free  $f _*\widetilde{\cO} _{\fY}$-module,  which is obvious.
\end{proof}

\begin{prop}
Let $f \colon \fY \to \fX$ be a finite morphism of  $\Sm ^\dag _\eta$. Let $\cM \in  D ^\mathrm{b} _{\mathrm{coh}} (\widetilde{\cO} _{\fX})$,
and $\cN \in  D ^\mathrm{b} _{\mathrm{coh}} (\widetilde{\cO} _{\fY})$. We have the canonical isomorphism in $D ^\mathrm{b} _{\mathrm{coh}} (f _*\widetilde{\cO} _{\fY}):$
\begin{equation} \label{dualityisoflat}
\R f _* \R \mathcal{H} om _{\widetilde{\cO} _\fY} (\cN, f ^\flat (\cM) ) \riso \R \mathcal{H} om _{\widetilde{\cO} _\fX} (\R f _* (\cN), \cM ).
\end{equation}
\end{prop}

\begin{proof}
We construct the isomorphism \ref{dualityisoflat} by composition as follows:
\begin{gather}
\notag
\R f _* \R \mathcal{H} om _{\widetilde{\cO} _\fY} (\cN, f ^\flat (\cM) ) \underset{\ref{II5.5HaRD}}{\riso} \R \mathcal{H} om _{f _*\widetilde{\cO} _\fY} (\R f _*\cN, \R f _*f ^\flat (\cM) )   \underset{\ref{prop-f*O-fprop}.(\ref{prop-f*O-fprop(e)})}{\riso}
\\  \notag \R \mathcal{H} om _{f _*\widetilde{\cO} _\fY} (\R f _*\cN, \R \mathcal{H} om _{\widetilde{\cO} _\fX} (f _* \widetilde{\cO} _\fY , \cM)   ) \underset{\ref{cohf*O=parf(2)}}{\liso} \R \mathcal{H} om _{\widetilde{\cO} _\fX}  (\R f _* (\cN), \cM ).
\end{gather}
\end{proof}

\begin{prop}
\label{prop-III.6.2}
Let $g \colon \fZ \to \fY$,
$f \colon \fY \to \fX$
be two finite morphisms of 
$\Sm ^\dag _\eta$.
For any 
$\cM \in 
D ^\mathrm{b} _{\mathrm{coh}}
(\widetilde{\cO} _{\fX})$,
we have the canonical isomorphism
$g ^\flat \circ f ^\flat (\cM)
\riso 
(f\circ g) ^\flat (\cM)$.

\end{prop}

\begin{proof}The isomorphism is the following composite one:
\begin{gather}
\notag
(f\circ g) ^\flat (\cM)
\riso
\overline{g} ^* \overline{f} ^* 
\R \mathcal{H} om _{\widetilde{\cO} _\fX}
(f _* g _*\widetilde{\cO} _\fZ , \cM)
\underset{\ref{dualityisoflat}}{\liso}
\overline{g} ^* \overline{f} ^* 
\R f _* \R \mathcal{H} om _{\widetilde{\cO} _\fY}
(g _*\widetilde{\cO} _\fZ, f ^\flat (\cM) )
\\
\underset{\ref{prop-f*O-fprop}.(\ref{prop-f*O-fprop(e)})}{\riso} \overline{g} ^*  \R \mathcal{H} om _{\widetilde{\cO} _\fY} (g _* \widetilde{\cO} _\fZ ,  f ^\flat (\cM) ) = g ^\flat \circ f ^\flat (\cM).
\end{gather}

\end{proof}

\begin{lem}
\label{lem-flatExtcomm}
Let $Y \to X$ be a morphism of topological spaces, $\AA \to \B$  (resp. $\AA \to \AA'$) be a morphism (resp. a flat morphism) of sheaves of rings on $X$. 
Set  $\B' := \AA '\otimes _{\AA} \B$. Let $\cM \in D ^+ (\AA)$. 
\begin{enumerate}[(a)]
\item We have the morphism of $D ^+ (\B ')$ of the form 
\begin{equation} \label{lem-flatExtcomm-iso}
\B ' \otimes _{\B} \R \mathcal{H} om _{\AA} (\B, \cM) \to  \R \mathcal{H} om _{\AA'} (\B', \AA' \otimes _{\AA}\cM).
\end{equation}
When  $\B \in D ^{\mathrm{b}} _{\mathrm{parf}} (\AA)$, this morphism is an isomorphism.
 \item We have the morphism of $D ^+ (u ^{-1}\B)$ of the form 
\begin{equation} \label{lem-flatExtcomm-iso2}
u ^{-1}\R \mathcal{H} om _{\AA} (\B, \cM) \to  \R \mathcal{H} om _{u ^{-1}\AA} (u ^{-1}\B, u ^{-1}\cM).
\end{equation}
When  $\B \in D ^{\mathrm{b}} _{\mathrm{parf}} (\AA)$, this morphism is an isomorphism.
\end{enumerate}
\end{lem}

\begin{proof}
If $\cI$ is an $\AA$-module, we compute that the canonical morphism $\mathcal{H} om _{\AA} (\B, \cI) \to  \mathcal{H} om _{\AA'} (\AA' \otimes _{\AA} \B, \AA' \otimes _{\AA}\cI)$ is $\B$-linear. 
Hence, taking a right resolution of $\cM$ by injective $\AA$-modules,  we construct in  $D ^+ (\B )$ the morphism 
$\R \mathcal{H} om _{\AA} (\B, \cM) \to  \R \mathcal{H} om _{\AA'} (\B', \AA' \otimes _{\AA}\cM)$, which  yields the desired morphism by extension.

Suppose $\B \in D ^{\mathrm{b}} _{\mathrm{parf}}
(\AA)$. It is enough to check that 
the composition of \ref{lem-flatExtcomm-iso}
with the isomorphism
$\AA ' \otimes _{\AA}
\R \mathcal{H} om _{\AA}
(\B, \cM)
\riso 
\B ' \otimes _{\B}
\R \mathcal{H} om _{\AA}
(\B, \cM)$
gives an isomorphism of
$D ^+ (\AA ')$.
To check this latter property, we reduce to the case where
$\B$ is a free $\AA$-module, which is obvious.

The second part of the lemma is checked similarly.
\end{proof}

\begin{prop}
\label{prop-III.6.4}
Let $f \colon \fY \to \fX$ be a finite morphism of 
$\Sm ^\dag _\eta$,
$u\colon \fX ' \to \fX$ be a morphism of 
$\Sm ^\dag _\eta$ such that 
$u _\eta$ is smooth. 
Let $\fY ':= \fY \times _{\fX} \fX '$, 
$v\colon \fY'\to \fY$, and 
$g\colon \fY' \to \fX'$ be the projections.
For any 
$\cM \in 
D ^\mathrm{b} _{\mathrm{coh}}
(\widetilde{\cO} _{\fX})$,
we have the canonical isomorphism
$\widetilde{v} ^* \circ f ^\flat (\cM)
\riso 
g ^\flat \circ \widetilde{u} ^*(\cM)$.

\end{prop}

\begin{proof}
By adjunction, we have the morphism
$\widetilde{u} ^* f _* \widetilde{\cO} _{\fY}
\to 
g _* \widetilde{v} ^*  \widetilde{\cO} _{\fY}
=
g _*  \widetilde{\cO} _{\fY'}$.
Since this is a morphism of coherent 
$\cO _{\fY'}$-module, since its restriction to 
$\fY ' _\eta$ is an isomorphism,
then it is an isomorphism.
This yields the morphism 
of ringed spaces
$\widetilde{u} '\colon 
(\fX' , g _* \widetilde{\cO} _{\fY'})
\to
(\fX , f _* \widetilde{\cO} _{\fY})$.
We get
$\widetilde{v} ^* \circ f ^\flat (\cM)
=
\widetilde{v} ^*
\overline{f} ^* 
\R \mathcal{H} om _{\widetilde{\cO} _\fX}
(f _* \widetilde{\cO} _\fY , \cM)
\riso
\overline{g} ^* 
\widetilde{u} ^{\prime *}
\R \mathcal{H} om _{\widetilde{\cO} _\fX}
(f _* \widetilde{\cO} _\fY , \cM)
\underset{\ref{lem-flatExtcomm}}{\riso}
\overline{g} ^* 
\R \mathcal{H} om _{\widetilde{\cO} _{\fX'}}
(g _* \widetilde{\cO} _{\fY'}, \widetilde{u} ^*\cM))
=
g ^\flat \circ \widetilde{u} ^*(\cM)$.
\end{proof}

\begin{prop} [Fundamental local isomorphism]\label{FLI}
Let $u \colon \fZ \hookrightarrow \fX$ be a closed immersion of codimension $d$ of $\Sm ^\dag _\eta$.
Put  $\widetilde{\omega} _{\fZ /\fX} := \mathcal{H}om _{\widetilde{\cO} _{\fZ}}  ( \wedge ^n (\widetilde{\cI}  /\widetilde{\cI}  ^2), \widetilde{\cO} _{\fZ})$, where $\cI$ is the ideal defining $u$.
\begin{enumerate}[(a)]
\item We have the exact sequence of locally free of finite type $\widetilde{\cO} _{\fX}$-modules: 
\begin{equation} \label{exactsequenceOvCvbis}
0 \to  \widetilde{\cI}  / \widetilde{\cI}  ^2  \to  \widetilde{u} ^* \widetilde{\Omega} _{\fX /\fS }  \to \widetilde{\Omega} _{\fZ  /\fS }   \to  0
\end{equation}

\item Let $\cM$ be a coherent  $\widetilde{\cO} _{\fX}$-module. We have the isomorphism
$$ \mathcal{E} xt ^d _{\widetilde{\cO} _{\fX}} (\widetilde{\cO} _{\fZ}, \cM)  \riso  \cM \otimes _{\widetilde{\cO} _{\fX}} \widetilde{\omega} _{\fZ /\fX} . $$
Furthermore, if $\cM$ is $\widetilde{u} ^*$-acyclic then $\mathcal{E} xt ^{d-r} _{\widetilde{\cO} _{\fX}} (\widetilde{\cO} _{\fZ}, \cM)  =0$, for any $r \not = 0$.
\end{enumerate}
\end{prop}

\begin{proof}
The first part is similar to \ref{exactsequenceOvCv}.
The proof of the second part is identical to the standard theorem of \cite[III.7.2]{HaRD}: 
$\cM \otimes _{\widetilde{\cO} _{\fX}}
\widetilde{\omega} _{\fZ /\fX}
\riso 
\mathcal{H}om _{\widetilde{\cO} _{\fZ}} 
( \wedge ^n (\widetilde{\cI}
 /\widetilde{\cI}
 ^2),
\cM / \widetilde{\cI}
 \cM)$.
Since $\widetilde{\cI}
$ is regular (see \ref{prop-regularbis}),
$\widetilde{\cI}
$ is locally generated by a regular sequence. 

When $\widetilde{\cI}
$ is generated by a regular sequence $\underline{f}$,
using the Koszul complex given by such regular sequence, 
we get the isomorphism
$\mathcal{E} xt ^{d-r} _{\widetilde{\cO} _{\fX}}
(\widetilde{\cO} _{\fZ},
\cM) 
\riso 
\mathcal{T} or _{r} ^{\widetilde{\cO} _{\fX}}
(\widetilde{\cO} _{\fZ},
\cM)$ (see \cite[III.1.1]{EGAIII1}).
Hence, when 
 $\cM$ is $\widetilde{u} ^*$-acyclic then
$\mathcal{E} xt ^{d-r} _{\widetilde{\cO} _{\fX}}
(\widetilde{\cO} _{\fZ},
\cM) 
=0$,
for any $r \not = 0$. 
For $r=0$, we get 
$\phi _{\underline{f}}
\colon 
\mathcal{E} xt ^d _{\widetilde{\cO} _{\fX}}
(\widetilde{\cO} _{\fZ},
\cM) 
\riso 
\cM / \widetilde{\cI}
 \cM$.
The isomorphism
$\widetilde{\cO} _{\fZ}\riso \wedge ^n (\widetilde{\cI}
 /\widetilde{\cI}
 ^2)$,
given by 
$a \mapsto a f _1 \wedge \dots \wedge f _d$
induces 
$\psi _{\underline{f}}
\colon 
\cM / \widetilde{\cI}
 \M
\liso
\mathcal{H}om _{\widetilde{\cO} _{\fZ}} 
( \wedge ^n (\widetilde{\cI}
 /\widetilde{\cI}
 ^2),
\cM / \widetilde{\cI}
 \cM)$.
This yields
$\mathcal{E} xt ^d _{\widetilde{\cO} _{\fX}}
(\widetilde{\cO} _{\fZ},
\cM) 
\riso 
\cM / \widetilde{\cI}
 \cM \otimes _{\widetilde{\cO} _{\fZ}}
\widetilde{\omega} _{\fZ /\fX} $.
It remains to check that this isomorphism does not depend on the choice of 
$\underline{f}$.
This computation is a consequence of the analogue of \cite[Lemma 7.1]{HaRD}.
\end{proof}

\begin{coro}
Let $u \colon \fZ \hookrightarrow \fX$ be a closed immersion of codimension $d$ of
$\Sm ^\dag _\eta$.
Let $\cM \in 
D ^{\mathrm{b}} _{\mathrm{coh}} ( \widetilde{\cO} _\fX)$.
We have the isomorphism
\begin{equation}
\label{cor-FLI}
\eta _u \colon 
u ^\flat (\cM) 
\riso 
\L \widetilde{u} ^* (\cM) 
\otimes _{\widetilde{\cO} _{\fZ}}
\widetilde{\omega} _{\fZ /\fX} [-d].
\end{equation}

\end{coro}

\begin{ntn}
Let $f\colon \fY \to \fX$ be a morphism of 
$\Sm ^\dag _\eta$ such that 
$f _\eta$ is smooth of relative dimension $d _f$. 
We set $\omega _{\fY /\fX}
:=
\wedge ^d 
\Omega _{\fY /\fX}$.
In the same way as \cite[III.2]{HaRD},  for any 
$\cM \in D ( \widetilde{\cO} _\fX)$, we set
\begin{equation}
\label{dfnfsharp}
f ^\sharp (\cM) 
:= 
\widetilde{f} ^* (\cM) \otimes 
_{\widetilde{\cO} _\fY}
\widetilde{\omega} _{\fY /\fX} [d _f].
\end{equation}
When 
$\cM \in 
D ^{\mathrm{b}} _{\mathrm{coh}} ( \widetilde{\cO} _\fX)$, 
we have 
$f ^\sharp (\cM) 
\in 
D ^{\mathrm{b}} _{\mathrm{coh}} ( \widetilde{\cO} _\fY)$.
\end{ntn}

\begin{prop}
\label{prop-III.2.2}
Let $g\colon \fZ \to \fY$, 
$f\colon \fY \to \fX$ be two morphisms of 
$\Sm ^\dag _\eta$ such that 
$f _\eta$ and $g _\eta$ are smooth.
\begin{enumerate}[(a)]
\item  We have a canonical exact sequence of coherent locally free $\widetilde{\cO} _\fY$-sheaves of the form
\begin{equation}
\label{exseq-III.2.2}
0 
\to 
\widetilde{g} ^* (\widetilde{\Omega} _{\fY /\fX})
\to 
\widetilde{\Omega} _{\fZ /\fX}
\to 
\widetilde{\Omega} _{\fZ /\fY}
\to 
0.
\end{equation}

\item There is an isomorphism
\begin{equation}
\label{iso-III.2.2}
\zeta _{f,g}
\colon 
(f \circ g ) ^\sharp  (\cM) 
\riso 
g ^\sharp \circ f ^\sharp (\cM),
\end{equation}
for all 
$\cM \in D ( \widetilde{\cO} _\fX)$.

\end{enumerate}
\end{prop}

\begin{proof}
Since this is the case on the generic fiber, 
we check the exactness of the sequence. 
This yields 
$\widetilde{g} ^*  ( \widetilde{\omega} _{\fY /\fX})
\otimes  _{\widetilde{\cO} _\fZ}
\widetilde{\omega} _{\fZ /\fY}
\riso 
\widetilde{\omega} _{\fZ /\fX} $,
and then
\begin{gather}
\notag
g ^\sharp \circ f ^\sharp (\cM)
= 
\widetilde{g} ^*  ( 
\widetilde{f} ^* (\cM) 
\otimes  _{\widetilde{\cO} _\fY}
\widetilde{\omega} _{\fY /\fX} [d _f])
\otimes  _{\widetilde{\cO} _\fZ}
\widetilde{\omega} _{\fZ /\fY}
[d _g]
\riso 
\widetilde{f\circ g} ^*  (\cM) 
\otimes  _{\widetilde{\cO} _\fZ}
\widetilde{\omega} _{\fZ /\fX} [d _{f\circ g}]
=
(f \circ g ) ^\sharp  (\cM) 
\end{gather}

\end{proof}

\begin{lem}
\label{lem-III.8.1}
Let $f\colon \fY \to \fX$ be a morphism of 
$\Sm ^\dag _\eta$ such that 
$f _\eta$ is smooth, 
and let $u \colon \fX \to \fY$ be a section of $f$. 
Then there is a functorial isomorphism
\begin{equation}
\label{iso-III.8.1}
\psi _{u,f}
\colon
\cM 
\riso 
u ^\flat \circ f ^\sharp (\cM)
\end{equation}
for all 
$\cM \in 
D ^{\mathrm{b}} _{\mathrm{coh}} ( \widetilde{\cO} _\fX)$.

\end{lem}

\begin{proof}
From the exact sequence \ref{exseq-III.2.2}, 
we get the isomorphism
$\widetilde{\omega} _{\fY /\fX}
\riso 
\widetilde{\omega} _{\fY /\fS}
\otimes _{\widetilde{\cO} _{\fY}}
\widetilde{f} ^* (\widetilde{\omega} _{\fX /\fS} ^{-1})$.
On the other hand, 
from \ref{exactsequenceOvCvbis} and its notation, we get
$\widetilde{\omega} _{\fX /\fY}
\riso 
\widetilde{\omega} _{\fX /\fS} 
\otimes _{\widetilde{\cO} _{\fX}}
\widetilde{u} ^* (\widetilde{\omega} _{\fY /\fS} ^{-1})$.
Hence, we get
$\L \widetilde{u} ^* (\widetilde{\omega} _{\fY /\fX}) 
\otimes _{\widetilde{\cO} _{\fX}}
\widetilde{\omega} _{\fX /\fY}
\riso
\widetilde{u} ^* (\widetilde{\omega} _{\fY /\fX}) 
\otimes _{\widetilde{\cO} _{\fX}}
\widetilde{\omega} _{\fX /\fY}
\riso 
(\widetilde{u} ^* (\widetilde{\omega} _{\fY /\fS} )
\otimes _{\widetilde{\cO} _{\fX}}
\widetilde{\omega} _{\fX /\fS} ^{-1} )
\otimes _{\widetilde{\cO} _{\fX}}
(\widetilde{\omega} _{\fX /\fS} 
\otimes _{\widetilde{\cO} _{\fX}}
\widetilde{u} ^* (\widetilde{\omega} _{\fY /\fS} ^{-1}))
\riso 
\widetilde{\cO} _{\fX}$.
This yields the last isomorphism:
\begin{gather}
\notag
u ^\flat \circ f ^\sharp (\cM)
=
u ^\flat \left (
\widetilde{f} ^* (\cM) \otimes 
_{\widetilde{\cO} _\fY}
\widetilde{\omega} _{\fY /\fX} [d]
\right) 
\underset{\ref{cor-FLI}}{\riso}
\L \widetilde{u} ^* (\widetilde{f} ^* (\cM) \otimes 
_{\widetilde{\cO} _\fY}
\widetilde{\omega} _{\fY /\fX}) 
\otimes _{\widetilde{\cO} _{\fX}}
\widetilde{\omega} _{\fX /\fY}
\\
\notag
\riso 
 \cM \otimes 
_{\widetilde{\cO} _\fX}
\L \widetilde{u} ^* (\widetilde{\omega} _{\fY /\fX}) 
\otimes _{\widetilde{\cO} _{\fX}}
\widetilde{\omega} _{\fX /\fY}
\riso 
\M.
\end{gather}

\end{proof}

\begin{prop}
[Residue Isomorphism]
\label{prop-III.8.2}
Let $g\colon \fZ \to \fY$, 
$f\colon \fY \to \fX$ be two morphisms of 
$\Sm ^\dag _\eta$.
We suppose that 
$f _\eta$ is smooth, $g$ is finite and $f\circ g$ is finite. 
Then there is a functorial isomorphism
\begin{equation}
\label{III.8.2}
\psi _{f,g}
\colon 
(f \circ g ) ^\flat (\cM)
\riso 
g ^\flat \circ f ^\sharp (\cM),
\end{equation}
for all 
$\cM \in 
D ^{\mathrm{b}} _{\mathrm{coh}} ( \widetilde{\cO} _\fX)$.

\end{prop}

\begin{proof}
The proof is identical to that of \cite[III.8.2]{HaRD}:
this is a consequence of 
Propositions \ref{prop-III.6.2}, \ref{prop-III.6.4},
and Lemma \ref{lem-III.8.1}.

\end{proof}

\begin{prop}
\label{prop-III.8.4}
Let $g\colon \fZ \to \fY$, 
$f\colon \fY \to \fX$ be two morphisms of 
$\Sm ^\dag _\eta$.
We suppose that 
$f _\eta$ and $(f\circ g) _\eta$ are smooth, and $g$ is finite.
Then there is a functorial isomorphism
\begin{equation}
\label{III.8.4}
\psi _{f,g}
\colon 
(f \circ g ) ^\sharp (\cM)
\riso 
g ^\flat \circ f ^\sharp (\cM),
\end{equation}
for all 
$\cM \in 
D ^{\mathrm{b}} _{\mathrm{coh}} ( \widetilde{\cO} _\fX)$.

\end{prop}

\begin{proof}
The proof is identical to that of \cite[III.8.4]{HaRD}:
this is a consequence of 
Propositions \ref{prop-III.2.2}, \ref{prop-III.8.2}
and Lemma \ref{lem-III.8.1}.

\end{proof}

\subsubsection{Adjoint differential operators}
Let $\fX ^\sharp$ be an object of $\Sm ^\dag _{\fS ^\sharp}$. Let $j \colon \fX _\eta \hookrightarrow \fX ^{\sharp}$ be the open immersion.  Let $m\in \bbN$. 
\begin{empt}
[Adjoint operator]
\label{adjoint operator-dag(m)}
Suppose  $\fX ^\sharp$ is affine, and $\fX ^\sharp/\fS ^\sharp$ has logarithmic coordinates. 
Fix $t _1,\dots, t _d  \in \Gamma (\fX, \widetilde{\cO} _{\fX})$   some overconvergent local coordinates  of $\fX ^\sharp/\fS ^\sharp$. 
Let $P \in \Gamma (\fX, \widetilde{\cD} ^{(m)\dag}  _{\fX^\sharp/\fS^\sharp})$.
Following the description 
of $\Gamma (\fX, \widetilde{\cD} ^{(m)\dag}  _{\fX^\sharp/\fS^\sharp})$
of \ref{wkcpdscplevelm}, we can write
$P = \sum _{\underline{i}\in \bbN ^d} 
\sum _{k\in \bbN}
\frac{a _{\underline{i},k} }{t ^{k}}\underline{\partial} ^{<\underline{i}> _{(m)}}$, 
where $a _{\underline{i},k} \in A$ are such that there exists a constant $c>0$ satisfying
$v _{p} ( a _{\underline{i},k}) \geq \frac{k+|\underline{i}|}{c} -1$,
for any $\underline{i} \in \bbN ^d$
and $k \in \bbN$.
Using the respective description of  \ref{wkcpdscplevelm} we get an element of $ \Gamma (\fX, \widetilde{\cD} ^{(m)\dag}  _{\fX^\sharp/\fS^\sharp}) $ the {\it adjoint operator of $P$} by setting:
$$^t P := \sum _{\underline{i}\in \bbN ^d} \sum _{k\in \bbN}  \underline{\partial} ^{<\underline{i}> _{(m)}} \frac{(-1) ^{|\underline{i}|}a _{\underline{i},k}}{t ^{k}} \in \Gamma (\fX, \widetilde{\cD} ^{(m)\dag}  _{\fX^\sharp/\fS^\sharp}).$$
This adjoint operator is compatible with 
the morphism 
$\Gamma (\fX _\eta, \widehat{\cD} ^{(m)}  _{\fX _\eta/\eta})
\to 
\Gamma (\fX _\eta, \widehat{\cD} ^{(m)}  _{\fX _\eta/\eta})$ 
given by 
$P = \sum _{\underline{i}\in \bbN ^d} 
b _{\underline{i}}\underline{\partial} ^{<\underline{i}> _{(m)}}
\mapsto 
 \sum _{\underline{i}\in \bbN ^d} 
 ( -1) ^{|\underline{i}|}
\underline{\partial} ^{<\underline{i}> _{(m)}}
b _{\underline{i}}
=
{}^t P$,
where $b _{\underline{i}} \in \Gamma (\fX _\eta, \cO   _{\fX _\eta})$
converge to $0$ when $|\underline{i}|$ goes to infinity.
This later morphism corresponds to the projective limits of the usual adjoint operator
of \cite[1.2.2]{Be2}.
\end{empt}

\begin{lem} \label{lemfullfaith}
The extension  $\widetilde{\cD} ^{(m) \dag}  _{\fX^\sharp/\fS^\sharp} \to  j _* \widehat{\cD} ^{(m)}  _{\fX_\eta/\eta}$ is faithfully flat (in the sense of the definition after Lemma \cite[4.3.8]{Be1}).
\end{lem}
\begin{proof} 
Using Lemma \cite[4.3.8]{Be1}, this is a consequence of  \ref{Jacobson}.(\ref{Jacobson(c)}).
\end{proof}

\begin{prop}\label{twist-Ddagm}
Let  $\cE$ be a coherent left $\widetilde{\cD} ^{(m) \dag}  _{\fX^\sharp/\fS^\sharp}$-module. Let  $\cM$ be a coherent right  $\widetilde{\cD} ^{(m) \dag}  _{\fX^\sharp/\fS^\sharp}$-module.

\begin{enumerate}[(a)]
\item \label{twist-Ddagm(a)}  Then  $\cE \otimes _{\widetilde{\cO} _{\fX}} \widetilde{\omega} _{\fX ^\sharp /\fS ^\sharp}$ is endowed with a canonical right $\widetilde{\cD} ^{(m) \dag}  _{\fX^\sharp/\fS^\sharp}$-module structure, 
and  $\cM \otimes _{\widetilde{\cO} _{\fX}} \widetilde{\omega} ^{-1}_{\fX ^\sharp /\fS ^\sharp} \riso  \mathcal{H} om _{\widetilde{\cO} _{\fX}} (\widetilde{\omega} _{\fX ^\sharp /\fS ^\sharp}, \cM)$ 
is endowed with a canonical left $\widetilde{\cD} ^{(m) \dag}  _{\fX^\sharp/\fS^\sharp}$-module structure.

\item \label{twist-Ddagm(b)} Suppose  $\fX ^\sharp$ is affine, and $\fX ^\sharp/\fS ^\sharp$ has logarithmic coordinates.  Fix $t _1,\dots, t _d  \in \Gamma (\fX, \widetilde{\cO} _{\fX})$   some overconvergent local coordinates  of $\fX ^\sharp/\fS ^\sharp$. 
Let $P \in \Gamma (\fX, \widetilde{\cD} ^{(m)\dag}  _{\fX^\sharp/\fS^\sharp})$,  $x \in \Gamma (\fX, \cE)$, $y \in \Gamma (\fX, \cM)$. We have more precisely the formula
\begin{gather} \notag
(x \otimes d t _1 \wedge \cdots \wedge d t _d) \cdot P = (^t P \cdot x )\otimes d t _1 \wedge \cdots \wedge d t _d,
\\ \notag P \cdot (y \otimes (d t _1 \wedge \cdots \wedge d t _d) ^\vee ) = (y \cdot ^t P )\otimes (d t _1 \wedge \cdots \wedge d t _d) ^\vee .
\end{gather}
\item \label{twist-Ddagm(c)} The canonical isomorphisms 
\begin{gather} \label{left2right}
(\cM \otimes _{\widetilde{\cO} _{\fX}} \widetilde{\omega} ^{-1}_{\fX ^\sharp /\fS ^\sharp} ) \otimes _{\widetilde{\cO} _{\fX}} \widetilde{\omega} _{\fX ^\sharp /\fS ^\sharp}
\riso  \cM  \text{, and } (\cE \otimes _{\widetilde{\cO} _{\fX}} \widetilde{\omega} _{\fX ^\sharp /\fS ^\sharp} ) \otimes _{\widetilde{\cO} _{\fX}} \widetilde{\omega} ^{-1}_{\fX ^\sharp /\fS ^\sharp} \riso  \E
\end{gather}
are $\widetilde{\cD} ^{(m) \dag}  _{\fX^\sharp/\fS^\sharp}$-linear.
\end{enumerate}
\end{prop}

\begin{proof}
By setting $\cE _{\eta}:=  \widehat{\cD} ^{(m)}  _{\fX_\eta/\eta} \otimes _{j ^{-1}\widetilde{\cD} ^{(m) \dag}  _{\fX^\sharp/\fS^\sharp}} j ^{-1} \cE$, we get: 
\begin{equation}\label{twist-Ddagmproof1} \cE \otimes _{\widetilde{\cO} _{\fX}} \widetilde{\omega} _{\fX ^\sharp /\fS ^\sharp}  \to j _* j ^{-1} (\cE \otimes _{\widetilde{\cO} _{\fX}} \widetilde{\omega} _{\fX ^\sharp /\fS ^\sharp} ) 
\riso j _* (j ^{-1} (\cE)  \otimes _{\cO _{\fX _\eta}} \omega _{\fX _{\eta} /\eta} ) \to j _* (\cE _{\eta}  \otimes _{\cO _{\fX _\eta}} \omega _{\fX _{\eta} /\eta} ).
\end{equation}
Let us prove that the composite map \ref{twist-Ddagmproof1} is injective. Since this is local, we can Suppose $m\in \bbN$. Suppose  $\fX ^\sharp$ is affine, and $\fX ^\sharp/\fS ^\sharp$ has logarithmic coordinates. 
Let  $D ^{(m)} := \Gamma (\fX, \cD ^{(m)} _{\fX^\sharp/\fS^\sharp} )$, $\widetilde{D} ^{(m)} := \Gamma (\fX, \widetilde{\cD} ^{(m)}  _{\fX^\sharp/\fS^\sharp})$,  $A:= \Gamma (\fX, \cO _{\fX})$.
Fix $t _1,\dots, t _d  \in A _{[t]}$   some overconvergent local coordinates  of $\fX ^\sharp/\fS ^\sharp$. 
Then $ \widetilde{\omega} _{\fX ^\sharp /\fS ^\sharp}$ is a free $\widetilde{\cO} _{\fX}$-module of rank one generated by $d t _1 \wedge \cdots \wedge d t _d$
and $\omega _{\fX _{\eta} /\eta} $ is a free $\cO _{\fX _\eta}$-module generated by $d t _1 \wedge \cdots \wedge d t _d$.
Using Theorem of type $A$ of  \ref{cor-thADdagm}, we check that the canonical morphism
$\Gamma ( \fX _\eta , \widehat{\cD} ^{(m)}  _{\fX_\eta/\eta})  \otimes _{\Gamma (\fX, \widetilde{\cD} ^{(m)}  _{\fX^\sharp/\fS^\sharp})} \Gamma (\fX, \cE) \to \Gamma (\fX _\eta, \cE _{\eta}) $, 
Hence, taking the global section of \ref{twist-Ddagmproof1}, we get a morphism which is isomorphism to the canonical map 
$\Gamma (\fX, \cE) \to \Gamma ( \fX _\eta , \widehat{\cD} ^{(m)}  _{\fX_\eta/\eta})  \otimes _{\Gamma (\fX, \widetilde{\cD} ^{(m)}  _{\fX^\sharp/\fS^\sharp})} \Gamma (\fX, \cE)$ which is injective because of 
\ref{Jacobson}.(\ref{Jacobson(c)}). In fact, via the injective morphism of \ref{twist-Ddagmproof1}, 
$\cE \otimes _{\widetilde{\cO} _{\fX}} \widetilde{\omega} _{\fX ^\sharp /\fS ^\sharp}$
is a right  $\widetilde{\cD} ^{(m) \dag}  _{\fX^\sharp/\fS^\sharp}$-submodule of $j _* (\cE _{\eta}  \otimes _{\cO _{\fX _\eta}} \omega _{\fX _{\eta} /\eta} )$. More precisely, we get the formula of (\ref{twist-Ddagm(b)}) coming from the similar one
satisfied by $\cE _{\eta}  \otimes _{\cO _{\fX _\eta}} \omega _{\fX _{\eta} /\eta}$.
We check similarly the assertions concerning $\cM$ in (\ref{twist-Ddagm(a)}) and (\ref{twist-Ddagm(c)}). The last assertions are clear.
\end{proof}

\begin{coro} \label{omega-Dmodule}
$\widetilde{\omega} _{\fX ^\sharp /\fS ^\sharp}$ is endowed with a canonical right $\widetilde{\cD} ^{(m) \dag}  _{\fX^\sharp/\fS^\sharp}$-module structure.
\end{coro}

\subsubsection{Frobenius descent}

\begin{empt}
[Definition and properties of $\Sat ^\dag _{\fS ^\sharp}$]
We denote by $\Sat ^\dag _{\fS ^\sharp}$ the category of 
formal log scheme $\fX ^\sharp$  which are log smooth and saturated (see Definitions \cite[I.3.12 and II.2.10]{Tsuji-saturated}) over $\fS ^\sharp$ and whose generic fiber
$\fX _\eta$ (i.e. the fiber of $\fX ^\sharp$
over $\eta = \Spf \cV [[t]] \{ \frac{1}{t}\}$) 
is smooth over $\eta$ (in particular $M _{\fX} | \fX _\eta = \cO ^* _{\fX _\eta}$). 
Such objects are called saturated log smooth log formal $\fS ^\sharp$-schemes.

The notion of saturated morphism is  remarkably stable 
(see Propositions \cite[II.2.11--12]{Tsuji-saturated}). When we deal with cartesian products of objects in 
$\Sat ^\dag _{\fS ^\sharp}$ we do not have to be careful in which category (coherent log formal schemes, fine log formal schemes
or fine saturated log formal schemes) 
the cartesian product is considered since there are the same. 

Let $\fX ^\sharp$ be an object of $\Sat ^\dag _{\fS ^\sharp}$.
The log scheme $S ^\sharp$ is (log) regular (see definition \cite[II.4.5]{Tsuji-saturated}). 
Indeed, the regularity over the generic point $\eta$ is obvious. 
Concerning the special point, we compute $I _{\overline{s}} \cO _{S ^\sharp, \overline{s}}$ is the ideal generated by $t$, and then 
the equality involving the dimensions of  \cite[II.4.5]{Tsuji-saturated} is satisfied.
Hence, following \cite[II.4.8]{Tsuji-saturated}, $X ^\sharp$ is also (log) regular. 
Following \cite[II.4.7]{Tsuji-saturated}, this yields that $X ^\sharp$ is Cohen-Macaulay and normal. 
Recall that following \cite[II.4.2]{Tsuji-saturated}, the special fiber $X ^\sharp _s$ is reduced.
\end{empt}

\begin{ex}
Let $\fP ^\sharp$ be a strictly semistable log formal schemes over $\fS ^\sharp$.
Then
$\fP ^\sharp $  is an object of 
 $\Sat ^\dag _{\fS ^\sharp}$.
 Hence, we do not have to bother with cartesian products of strictly semistable log formal schemes over $\fS ^\sharp$.
\end{ex}

\begin{empt}
Let $X ^\sharp = (X,M)$ be a logarithmic scheme over $\F _p$.  We recall that the absolute Frobenius  $F _{X ^{\sharp}} \colon X ^{\sharp} \to X ^{\sharp}$ is defined as follows (see \cite[4.7]{Kato-logFontaine-Illusie}).
The underlying morphism of schemes   of  $F _{X ^{\sharp}}$ is the usual absolute Frobenius  $F _X \colon X \to X$,  and the homomorphism  $F ^{-1} _X (M) \to M$ is the multiplication by $p$ on $M$ under the canonical 
identification of $F ^{-1} _X (M)$ with $M$.  Let $s \geq 0$ be an integer. We denote by  $ F ^s _{X ^{\sharp}}$ the $s$th power of the absolute Frobenius of  $X ^{\sharp}$.
\end{empt}

Let $\fX ^\sharp$ be a formal log scheme log smooth over $\fS ^\sharp$.
Let $\fX _\eta$ be the open of $\fX ^\sharp$ complementary to $X _s$ i.e. the fiber of $\fX ^\sharp$
over $\eta = \Spf \cV [[t]] \{ \frac{1}{t}\}$. We suppose $\fX _\eta /\eta$ smooth.
Let $j \colon \fX _\eta \hookrightarrow \fX ^{\sharp}$ be the open immersion. 
We denote by $X ^{\sharp (s)}$ the log $S^\sharp$-scheme deduced from 
$X ^\sharp$ by base change by $ F ^s _{S ^{\sharp}}$ (in the category of coherent log $S^\sharp$-schemes).
We remark that the underlying scheme of $X ^{\sharp (s)}$
is $X ^{(s)}$, the base change of $X $ by $ F ^s _{S }$
(indeed following \cite[1.6]{Kato-logFontaine-Illusie} the forgetful functor from the category of log schemes to that of 
schemes commutes with finite inverse limits).
We denote by 
$F ^s _{X ^\sharp/ S ^\sharp}\colon 
X ^{\sharp} \to X ^{\sharp (s)}$ the relative Frobenius, i.e. such that we get the commutative diagram
$$\xymatrix{
{X ^{\sharp}} 
\ar[rd] ^-{}
\ar[r] ^-{F ^s _{X ^\sharp/ S ^\sharp}}
& 
{X ^{\sharp (s)} } 
\ar@{}[rd] ^-{}|\square
\ar[r] ^-{}
\ar[d] ^-{}
& 
{X ^{\sharp} } 
\ar[d] ^-{}
\\ 
{} 
& 
{S ^{\sharp} } 
\ar[r] ^-{F ^s _{S ^\sharp}}
& 
{S ^{\sharp} .} 
}
$$

We suppose that the $s$th power of Frobenius $k \to k$ has a lifting  $\sigma \colon \cV \riso \cV$. Then the $s$th power of Frobenius of $k [[t]]$ has the canonical lifting $F ^* _{\cV [[t]]} \colon \cV [[t]] \to \cV [[t]]$ given by 
$\sum _{i\geq 0} a _i t ^i \mapsto \sum _{i\geq 0} \sigma (a _i ) t ^{pi}$. We get the morphism of log formal schemes $F _{\fS ^\sharp}\colon \fS ^\sharp \to \fS ^\sharp$ 
whose underlying morphism of schemes $F _{\fS}$ is given by  $F ^* _{\cV [[t]]}$,  and the homorphism  $F ^{-1}_{\fS ^\sharp} ( M _{\fS ^\sharp}) =M _{\fS ^\sharp} \to M _{\fS ^\sharp}$ is induced by  $F ^* _{\cV [[t]]}$.

We denote by $\fX ^{\sharp (s)}$ the log formal $\fS^\sharp$-scheme deduced from  $\fX ^\sharp$ by base change by 
$F _{\fS ^\sharp}$ (in the category of coherent log formal $\fS^\sharp$-schemes). We remark that the underlying scheme of $\fX ^{\sharp (s)}$ is $\fX ^{(s)}$, the base change of $\fX $ by $F ^s _{\fS }$.

When $\fX$ is affine, since  $\fX ^{\sharp (s)}/\fS ^\sharp$ is log smooth, using  \cite[3.11]{Kato-logFontaine-Illusie}, there exists 
a lifting  $F \colon  \fX ^{\sharp} \to \fX ^{\sharp (s)}$ of  $F ^s _{X ^\sharp/ S ^\sharp}\colon  X ^{\sharp} \to X ^{\sharp (s)}$ making commutative the diagram
$$\xymatrix{
{\fX ^{\sharp}}  \ar[rd] ^-{} \ar@{.>}[r] ^-{F } &  {\fX ^{\sharp (s)} }  \ar@{}[rd] ^-{}|\square \ar[r] ^-{} \ar[d] ^-{} &  {\fX ^{\sharp} }  \ar[d] ^-{}
\\  {}  &  {\fS ^{\sharp} }  \ar[r] ^-{F _{\fS ^\sharp}} &  {\fS ^{\sharp} .}   }   $$
The morphism $F$ induces by restriction  the morphism of formal $\eta$-schemes  $F _\eta \colon \fX _\eta \to \fX ^{ (s)} _\eta$. The functor $F _\eta ^*$ induces an equivalence 
between the category of  coherent $\widehat{\cD} ^{(m)} _{\fX ^{(s)}_\eta /\eta}$-modules  and that of  coherent $\widehat{\cD} ^{(m+s)} _{\fX _\eta /\eta}$-modules  (see \cite[4.1]{Be2}).

Suppose $\fX$ affine and that there exist log coordinates  $u _1,\dots, u _d \in \Gamma ( \fX , M _{\fX ^{\sharp}})$ of $\fX ^{\sharp} /\fS ^{\sharp} $.
We denote by  $t _1,\dots, t _d \in \Gamma ( \fX  , \cO _{\fX })$ the overconvergent local coordinates of $\fX ^{\sharp} /\fS ^{\sharp} $)
 induced by 
$u _1,\dots, u _d$ via the morphism 
$\Gamma ( \fX , M _{\fX ^{\sharp}}) \to \Gamma ( \fX , \cO _{\fX})$. 
We still denote by 
$t _1,\dots, t _d \in \Gamma ( \fX _\eta , \cO _{\fX _\eta})$ the local coordinates of 
$\fX _\eta /\eta$ 
 induced by 
$t _1,\dots, t _d$ via the inclusions 
$\Gamma ( \fX , \cO _{\fX})\subset  \Gamma ( \fX , \widetilde{\cO} _{\fX}) \subset \Gamma ( \fX _\eta , \cO _{\fX _\eta})$.

Let $\alpha\colon \fX ^\sharp 
\to 
\mathfrak{A} _{\bbN ^d, \fS ^\sharp}$
be the corresponding log étale morphism.
We denote by 
$F _{\mathfrak{A} _{\bbN ^d}}
\colon 
\mathfrak{A} _{\bbN ^d}
\to \mathfrak{A} _{\bbN ^d}$
the morphism of log formal schemes induced by the multiplication by $p$
of $\bbN ^d$.
We get the lifting
$F _{\mathfrak{A} _{\bbN ^d, \fS ^\sharp}}
=
F _{\mathfrak{A} _{\bbN ^d}}
\times 
F _{\fS ^\sharp}
\colon 
\mathfrak{A} _{\bbN ^d, \fS ^\sharp}
\to 
\mathfrak{A} _{\bbN ^d, \fS ^\sharp}
$
of the $s$th power of Frobenius of
$A _{\bbN ^d, S ^\sharp}$.
We denote by 
$F _{\mathfrak{A} _{\bbN ^d}}
\colon 
\mathfrak{A} _{\bbN ^d, \fS ^\sharp}
\to 
\mathfrak{A} _{\bbN ^d, \fS ^\sharp} ^{(s)}$
the canonical $\fS ^\sharp$-morphism induced by 
$F _{\mathfrak{A} _{\bbN ^d}}$.
Since 
$\alpha ^{(s)} \colon \fX ^{\sharp (s)}
\to 
\mathfrak{A} _{\bbN ^d, \fS ^\sharp} ^{(s)}$ is log etale, using 
\cite[3.11]{Kato-logFontaine-Illusie}, there exists 
a unique morphism
$F \colon \fX ^{\sharp} \to \fX ^{\sharp (s)} $
making commutative the following diagram
$$\xymatrix @R=0,3cm{      {\fX ^{\sharp}}  \ar@{.>}[r] ^-{F } \ar[d] ^-{\alpha}
&  {\fX ^{\sharp (s)} }  \ar@{}[rd] ^-{}|\square \ar[r] ^-{} \ar[d] ^-{\alpha ^{(s)}} &  {\fX ^{\sharp} }  \ar[d] ^-{\alpha}
\\  {\mathfrak{A} _{\bbN ^d, \fS ^\sharp}}  \ar[r] ^-{F _{\mathfrak{A} _{\bbN ^d}}} & {\mathfrak{A} _{\bbN ^d, \fS ^\sharp} ^{(s)}}  \ar@{}[rd] ^-{}|\square \ar[r] ^-{} \ar[d] ^-{} &  {\mathfrak{A} _{\bbN ^d, \fS ^\sharp}}  \ar[d] ^-{}
\\  {}  &  {\fS ^{\sharp} }  \ar[r] ^-{F _{\fS ^\sharp}} &  {\fS ^{\sharp} .}  } $$
Put $q:= p^s$. Following the notation of \cite[2.4.2]{Garnier-descente-frov-expl}, we set 
$$H _{\underline{t}} :=  \frac{1}{q ^d} \sum _{\underline{\zeta} ^q =\underline{1}} \sum _{\underline{k}\in \bbN ^d} (\underline{\zeta} - \underline{1} ) ^{\underline{k}} \underline{t} ^{\underline{k}} \underline{\partial} ^{[\underline{k}]}.$$

Following \cite[4.6.2.(1)]{Garnier-descente-frov-expl}, since $(1) ^{(s)}= H _{\underline{t}}$ (see the convention of \cite[4.3.3]{Garnier-descente-frov-expl}), 
then  $H _{\underline{t}} \in \Gamma (\fX _\eta,\widehat{\cD} ^{(s)} _{\fX _\eta /\eta})$. 
Since $H _{\underline{t}}$ in  an operator in $\underline{t} ^{\underline{k}} \underline{\partial} ^{[\underline{k}]} $ with coefficient in $\bbZ$ (see \cite[2.5.3]{Garnier-descente-frov-expl}), then we get 
$H _{\underline{t}} \in \Gamma (\fX,\cB ^{(s)} _{\fX} (X _s)  \widehat{\otimes} _{\cO _{\fX}} \widehat{\cD} ^{(s)} _{\fX ^{\sharp }/\fS ^{\sharp},\bbQ})$.
Hence, it follows from \ref{DdagVSDhat} and \ref{rem-QcapY=BXdag} that 
\begin{equation}\label{HinDs+1}
H _{\underline{t}} \in \Gamma (\fX,\widetilde{\cD} ^{(s+1)\dag} _{\fX ^\sharp/\fS ^{\sharp}}).
\end{equation}

Following \cite[2.5.1.(iii)]{Garnier-descente-frov-expl} (the last one is an easy computation), for any $\underline{r}\in \bbN ^d$ such that  $|\underline{r}| _{\infty} <q$, we have 
\begin{gather}\label{251iii2GDesc1}
\sum _{|\underline{r}| _{\infty} <q} \underline{t} ^{\underline{r}} H \underline{t} ^{-\underline{r}}=1,
\\ \label{251iii2GDesc} H _{\underline{t}}  \underline{t} ^{-\underline{r}}H _{\underline{t}}= \delta _{\underline{r}, \underline{0}} H _{\underline{t}},
\qquad H _{\underline{t}}  \underline{t} ^{\underline{r}}H _{\underline{t}}= \delta _{\underline{r}, \underline{0}} H _{\underline{t}}.
\end{gather}
where $\delta _{\underline{r}, \underline{0}}$ is the Kronecker symbol.

From \cite[2.5.3.(5)]{Garnier-descente-frov-expl}, for any $|\underline{r}|_{\infty} <q$ we compute that
\begin{equation} \label{preProp472G1}
H _{\underline{t}}  \underline{t} ^{-\underline{r}}  = \sum _{\underline{k}\geq \underline{r}}  \sum _{\underset{q| \underline{u}-\underline{r}}{\underline{u}\leq \underline{k}}} (-1) ^{\underline{k}-\underline{u}} 
\left ( \begin{smallmatrix}  \underline{k}\\  \underline{u} \end{smallmatrix} \right ) \underline{t} ^{\underline{k}-\underline{r}} \underline{\partial} ^{[\underline{k}]}
= \underline{\partial} ^{[\underline{r}]} + \sum _{\underline{k}>\underline{r}}  \alpha _{\underline{k},\underline{r}}  \underline{\partial} ^{[\underline{k}]},
\end{equation}
 with $\alpha _{\underline{k},\underline{r}}\in \Z [ t _1,\dots,t_d]$.
We can view the differential operator  $ H _{\underline{t}}  \underline{t} ^{-\underline{r}}$ as an element   of $\mathcal{H}om _{\cO _{\fX ^{(s)} _\eta } }  ( F _{\eta *} \cO _{\fX _\eta }  ,\cO _{\fX ^{(s)} _\eta } )$ given by
$a \mapsto H _{\underline{t}}  \underline{t} ^{-\underline{r}} (a)$   for any $a \in F _{\eta *} \cO _{\fX _\eta } $.  From \cite[2.5.1 and 2.5.3.(5)]{Garnier-descente-frov-expl}, 
$\{ H _{\underline{t}}  \underline{t} ^{-\underline{r}},   \text{  $\underline{r} \in \bbN ^d$ such that  $|\underline{r}|_{\infty} <q$}\}$
is a basis of the free $\cO _{\fX ^{(s)} _\eta } $-module  $\mathcal{H}om _{\cO _{\fX ^{(s)} _\eta } }  ( F _{\eta *} \cO _{\fX _\eta } ,\cO _{\fX ^{(s)} _\eta } )$.

\begin{empt}
From 
\cite[4.3]{Garnier-descente-frov-expl}, 
we have a morphism of 
$\cO _{\fX ^{(s)} _\eta} $-rings :
$\widehat{\cD} ^{(m)} _{\fX ^{(s)}_\eta /\eta}
\to 
F _* 
(
\widehat{\cD} ^{(m+s)} _{\fX _\eta /\eta}
)
$
we will denote by 
$P \mapsto P ^{(s)}$.
In our context, we have the explicit formula
given by
\cite[4.6.2.(2)]{Garnier-descente-frov-expl}
$$
(\underline{\partial} ^{<\underline{k}> _{(m)}} ) ^{(s)}
=
\sum _{\underline{u}\geq \underline{k}} 
\sum _{|\underline{r}| _{\infty} <q}
(-\underline{t}) ^{(\underline{u}-\underline{k})q +\underline{r}} 
\lambda _{\underline{k},\underline{u},\underline{r}}
\underline{\partial} ^{<\underline{u}q + \underline{r} >_{(m+s)}} ,
$$
where
$\lambda _{\underline{k},\underline{u},\underline{r}}
:=
\sum _{\underline{i}\leq \underline{u}-\underline{k}}
(-1) ^{\underline{i}} 
\left (
\begin{smallmatrix}
\underline{u}q + \underline{r}\\
q (\underline{k} + \underline{i})
\end{smallmatrix}
\right )
\left (
\begin{smallmatrix}
\underline{k} + \underline{i}\\
\underline{k} 
\end{smallmatrix}
\right )
\frac{q _{\underline{k}} ^{(m)}!}{q _{\underline{u}} ^{(m)}!}
$ (use also $q _{\underline{u}} ^{(m)} = q _{\underline{u}q + \underline{r}} ^{(m+s)}$).
Using \cite[4.6.2.(3)]{Garnier-descente-frov-expl} and \cite[1.1.0.1]{Be1}, we get the estimation
\begin{gather}\notag 
v _p (\lambda _{\underline{k},\underline{u},\underline{r}})  =v _p (\sum _{\underline{i}\leq \underline{u}-\underline{k}} (-1) ^{\underline{i}}  \left ( \begin{smallmatrix} \underline{u}q + \underline{r}\\ q (\underline{k} + \underline{i}) \end{smallmatrix} \right )
\left ( \begin{smallmatrix} \underline{k} + \underline{i}\\ \underline{k}  \end{smallmatrix} \right ) \frac{q _{\underline{k}} ^{(m)}!}{q _{\underline{u}} ^{(m)}!} )
\geq  \frac{p}{p-1} |\underline{u}-\underline{k} |  + \sum _{i= 1} ^d \frac{\sigma (k _i)-\sigma (u _i)}{p-1}  +v _p (q _{\underline{k}} ^{(m)}!) -v _p (q _{\underline{u}} ^{(m)}!)
\\  = \left ( |\underline{u}| +v _p (\underline{u}!)-   v _p (q _{\underline{u}} ^{(m)}!)\right )  - \left ( |\underline{k}|  +v _p (\underline{k}!)-   v _p (q _{\underline{k}} ^{(m)}!)\right)  \geq  |\underline{u}| -  \frac{p |\underline{k}|}{p-1} .
\end{gather}
In other word,  we can write
\begin{equation} \label{preProp472G2}
(\underline{\partial} ^{<\underline{k}> _{(m)}} ) ^{(s)} = \underline{\partial} ^{<q\underline{k} >_{(m+s)}} +  \sum _{\underline{n}>q\underline{k}}  \mu _{\underline{k},\underline{n}} \underline{\partial} ^{<\underline{n} >_{(m+s)}} ,
\end{equation}
with $\mu _{\underline{k},\underline{n}} \in \Z [ t _1, \dots, t _d]$ satisfying
\begin{equation}\label{preProp472G2bis}
v _p (\mu _{\underline{k},\underline{n}}) \geq    |q _{\underline{n}} ^{(s)}| -  \frac{p |\underline{k}|}{p-1} \geq \frac{|\underline{n}|}{(q+1)}-  \frac{p|\underline{k}|}{p-1} .
\end{equation}

\end{empt}

\begin{lem} \label{lemProp472G2}
We can write  $\underline{\partial} ^{<\underline{k} >_{(m+s)}}$ uniquely in the form
\begin{equation} \label{Prop472G2}
\underline{\partial} ^{<\underline{k} >_{(m+s)}} = \sum _{\underline{l} \geq \underline{k}} \beta _{\underline{k},\underline{l}}
(\underline{\partial} ^{<q _{\underline{l}} ^{(s)}> _{(m)}} ) ^{(s)}  H _{\underline{t}}  \underline{t} ^{- r_{\underline{l}} ^{(s)}}
\end{equation}
with $\beta _{\underline{k},\underline{l}} \in \Z [ t _1, \dots, t _d]$ satisfying $v _p (\beta _{\underline{k},\underline{l}})  \geq |q _{\underline{l}} ^{(s)} -q _{\underline{k}} ^{(s)}|$.
\end{lem}

\begin{proof}
See the proof of \cite[4.7.2]{Garnier-descente-frov-expl}. 
\end{proof}

\begin{lem}
\label{lem-factorization}
The homomorphism 
$\widehat{\cD} ^{(m)} _{\fX ^{(s)}_\eta /\eta}
\to 
F _{\eta *} 
(
\widehat{\cD} ^{(m+s)} _{\fX _\eta /\eta}
)
$
induces the homomorphism of $\cO _{\fX ^{(s)}}$-rings 
$\widetilde{\cD} ^{(m)\dag} _{\fX ^{\sharp (s)}/\fS ^{\sharp}}
\to 
F _*\widetilde{\cD} ^{(m+s)\dag} _{\fX ^{\sharp}/\fS ^{\sharp}}$.
\end{lem}

\begin{proof}
Let $P \in \widetilde{D} ^{(m)\dag} _{\fX ^{\sharp (s)}/\fS ^{\sharp}}$.  Following,  \ref{wkcpdscplevelm}, we can write  $P = \sum _{i\in \bbN}\sum _{\underline{k}\in \bbN ^d} a ' _{\underline{k},i} \frac{1}{t ^{i}}\underline{\partial} ^{<\underline{k}> _{(m)}}$
where $a ' _{\underline{k},i} \in A ^{(s)}= \Gamma (\fX ^{\sharp (s)}, \cO _{\fX ^{\sharp (s)}})$ are such that there exists a large enough constant $c>0$ satisfying 
$v _{p} ( a ' _{\underline{k},i}) \geq \frac{i+|\underline{k}|}{c} -1$, for any $\underline{k} \in \bbN ^d$ and $i \in \bbN$. Using \ref{preProp472G2}, we get: 
$$P ^{(s)} =  \sum _{i\in \bbN}\sum _{\underline{k}\in \bbN ^d}  \sum _{\underline{n}\geq q\underline{k}}  
F ^{*}(a ' _{\underline{k},i}) \frac{1}{t ^{i}} \mu _{\underline{k},\underline{n}} \underline{\partial} ^{<\underline{n} >_{(m+s)}}
= \sum _{i\in \bbN} \sum _{\underline{n}\in \bbN ^d}  \frac{1}{t ^{i}} \left ( \sum _{q\underline{k}\leq \underline{n}}  F ^{*}(a ' _{\underline{k},i})  \mu _{\underline{k},\underline{n}}\right ) \underline{\partial} ^{<\underline{n} >_{(m+s)}}  $$
Remark that increasing $c$ if necessary, the coefficient of $|\underline{n}|$ is positive.  When   $2q\underline{k}\leq \underline{n}$,  we compute:
\begin{gather*}
v _p (F ^{*}(a ' _{\underline{k},i})  \mu _{\underline{k},\underline{n}}) 
\underset{\ref{preProp472G2}}{\geq}     ( \frac{i+|\underline{k}|}{c} -1)  +  ( \frac{|\underline{n}|}{q+1} -\frac{p |\underline{k}|}{p-1}) 
 \geq   \frac{i}{c} +  \frac{|\underline{n}|}{q+1}    - (\frac{p}{p-1} - \frac{1}{c}) |\underline{k}|  -1
 \\
 \geq   \frac{i}{c} +  |\underline{n}| \left ( \frac{1}{q+1}    - \frac{1}{2q}(\frac{p}{(p-1} - \frac{1}{c}) \right )  -1
\end{gather*}
When   $2q\underline{k}\geq \underline{n}$,  we compute:
\begin{gather}\notag
v _p (F ^{*}(a ' _{\underline{k},i})  \mu _{\underline{k},\underline{n}}) 
\geq 
v _p (F ^{*}(a ' _{\underline{k},i}) ) \geq  
( \frac{i+|\underline{k}|}{c} -1)  
 \geq   \frac{i}{c} +  \frac{|\underline{n}|}{2qc}   -1.
\end{gather}
Hence, there exists $C>0$ large enough such that in both cases we get 
\begin{gather}\notag
v _p (F ^{*}(a ' _{\underline{k},i})  \mu _{\underline{k},\underline{n}})  \geq  ( \frac{i+|\underline{n}|}{C} -1) .
\end{gather}
Hence,  $P ^{(s)} \in  \widetilde{D} ^{(m+s)\dag} _{\fX ^{\sharp}/\fS ^{\sharp}}$. 
\end{proof}

\begin{empt}
We set
$F ^* _\eta
F ^\flat _\eta \widehat{\cD} ^{(m)} _{\fX ^{(s)}_\eta /\eta}
=
\cO _{\fX _\eta } 
\otimes _{\cO _{\fX ^{(s)} _\eta }}
 \left ( \widehat{\cD} ^{(m)} _{\fX ^{(s)}_\eta /\eta}
\otimes _{\cO _{\fX ^{(s)} _\eta } }
\mathcal{H}om _{\cO _{\fX ^{(s)} _\eta } } 
( F _{\eta *} \cO _{\fX _\eta } ,
\cO _{\fX ^{(s)} _\eta } ) \right )$. 
Following
\cite[4.7.2]{Garnier-descente-frov-expl}, 
the isomorphism of 
$\widehat{\cD} ^{(m+s)} _{\fX _\eta /\eta}$-bimodules
\begin{equation}
\label{Phi}
\Phi
\colon 
F ^* _\eta
F ^\flat _\eta \widehat{\cD} ^{(m)} _{\fX ^{(s)}_\eta /\eta}
\riso 
\widehat{\cD} ^{(m+s)} _{\fX _\eta /\eta}
\end{equation}

is explicitly given by 
\begin{equation}
\label{descr-Phi}
\Phi 
\left ( \sum _{|\underline{r}| _{\infty} <q}  \sum _{|\underline{u}|_{\infty} <q}  \underline{t} ^{\underline{r}}  \otimes   P _{\underline{r},\underline{u}}   \otimes   H _{\underline{t}}  \underline{t} ^{-\underline{u}}  \right )
=  \sum _{|\underline{r}| _{\infty} <q}  \sum _{|\underline{u}|_{\infty} <q}  \underline{t} ^{\underline{r}}  P _{\underline{r},\underline{u}} ^{(s)}  H _{\underline{t}}  \underline{t} ^{-\underline{u}},
\end{equation}
where in the finite sums $\underline{r},\underline{u} \in \bbN ^d$, and
$P _{\underline{r},\underline{u}}  \in \widehat{\cD} ^{(m)} _{\fX ^{(s)}_\eta /\eta}$. 

\end{empt}

We put
$F ^* F ^\flat  \widetilde{\cD} ^{(m)\dag} _{\fX ^{\sharp (s)}/\fS ^{\sharp}}
=
\widetilde{\cO} _{\fX}  
\otimes _{\widetilde{\cO} _{\fX ^{(s)}  }}  
\widetilde{\cD} ^{(m)\dag} _{\fX ^{\sharp (s)}/\fS ^{\sharp}}
\otimes _{\widetilde{\cO} _{\fX ^{(s)}}}
\mathcal{H}om _{\widetilde{\cO} _{\fX ^{(s)} } } 
(\widetilde{\cO} _{\fX} ,
\widetilde{\cO} _{\fX ^{(s)}  } ) $.

\begin{prop} \label{phidagdescrip}
The image of 
$F ^* F ^\flat  \widetilde{\cD} ^{(m)\dag} _{\fX ^{\sharp (s)}/\fS ^{\sharp}}$ via $j _*\Phi$ is
$\widetilde{\cD} ^{(m+s)\dag} _{\fX ^{\sharp}/\fS ^{\sharp}}$.
We get a 
canonical structure of 
$\widetilde{\cD} ^{(m+s)\dag} _{\fX ^{\sharp}/\fS ^{\sharp}}$-bimodule
on $F ^* F ^\flat  \widetilde{\cD} ^{(m)\dag} _{\fX ^{\sharp (s)}/\fS ^{\sharp}}$ and 
the commutative diagram of $\widetilde{\cD} ^{(m+s)\dag} _{\fX ^{\sharp}/\fS ^{\sharp}}$-bimodules 
\begin{equation}  \notag
\xymatrix{
{j _* F ^* _\eta F ^\flat _\eta \widehat{\cD} ^{(m)} _{\fX ^{(s)}_\eta /\eta}}  \ar[r] ^-{\sim} &  {j _* \widehat{\cD} ^{(m+s)} _{\fX _\eta /\eta}} 
\\  {F ^* F ^\flat  \widetilde{\cD} ^{(m)\dag} _{\fX ^{\sharp (s)}/\fS ^{\sharp}}}  \ar[r] ^-{\sim} \ar@{^{(}->}[u] ^-{} &  {\widetilde{\cD} ^{(m+s)\dag} _{\fX ^{\sharp}/\fS ^{\sharp}}.}  \ar@{^{(}->}[u] ^-{}
}
\end{equation}
In particular, any  $P \in \widetilde{\cD} ^{(m+s)\dag} _{\fX ^{\sharp}/\fS ^{\sharp}}$ can written uniquely in the form 
$P = \sum _{|\underline{r}| _{\infty} <q}  \sum _{|\underline{u}|_{\infty} <q}  \underline{t} ^{\underline{r}}  P _{\underline{r},\underline{u}} ^{(s)}  H _{\underline{t}}  \underline{t} ^{-\underline{u}}$,
where in the finite sums  $P _{\underline{r},\underline{u}}  \in \widetilde{\cD} ^{(m)\dag} _{\fX ^{\sharp (s)}/\fS ^{\sharp}}$. 
\end{prop}

\begin{proof}
1) From Lemma \ref{lem-factorization}, using the formula \ref{descr-Phi} we get the factorization. It remains to check the surjectivity.

2) Let $P = \sum _{i\in \bbN}\sum _{\underline{k}\in \bbN ^d} a  _{\underline{k},i} \frac{1}{t ^{i}}\underline{\partial} ^{<\underline{k}> _{(m+s)}}$
where $a _{\underline{k},i} \in A = \Gamma (\fX ^{\sharp}, \cO _{\fX ^{\sharp}})$ are such that there exists a constant $c>1$ satisfying
$v _{p} ( a  _{\underline{k},i}) \geq \frac{i+|\underline{k}|}{c} -1$,
for any $\underline{i} \in \bbN ^d$
and $k \in \bbN$.

From \ref{lemProp472G2}, we compute
\begin{gather}
P 
= 
\sum _{i\in \bbN}
\sum _{\underline{k}\in \bbN ^d} 
 \frac{1}{t ^{i}}  a  _{\underline{k},i}
\sum _{\underline{n}\geq \underline{q} ^{(s)} _{\underline{k}}}
 \sum _{|\underline{u}|_{\infty} <q}
\beta _{\underline{k},q\underline{n}+\underline{u}}
(\underline{\partial} ^{<\underline{n}}> _{(m)} ) ^{(s)}
 H _{\underline{t}}  \underline{t} ^{- \underline{u}}
 \\
 =
  \sum _{|\underline{u}|_{\infty} <q}
\sum _{\underline{n}\in \bbN ^d} 
  \sum _{i\in \bbN}
  \sum _{\underline{k}\leq q \underline{n}+\underline{u}}
  \frac{1}{t ^{i}}
 a  _{\underline{k},i}
\beta _{\underline{k},q\underline{n}+\underline{u}}
(\underline{\partial} ^{<\underline{n}> _{(m)}} ) ^{(s)}
 H _{\underline{t}}
 \underline{t} ^{- \underline{u}}.
\end{gather}
Since  $\beta _{\underline{k},q\underline{n}+\underline{u}} \in \Z [ t _1, \dots, t _d]$ satisfy
$v _p (\beta _{\underline{k},q\underline{n}+\underline{u}}) 
\geq |\underline{n} -q _{\underline{k}} ^{(s)}|$,
increasing $c$ is necessary, 
using similar computations than in the proof of \ref{lem-factorization}, 
we can write 
$\sum _{\underline{k}\leq q \underline{n}+\underline{u}}
a  _{\underline{k},i}
\beta _{\underline{k},q\underline{n}+\underline{u}}
=
 \sum _{|\underline{r}| _{\infty} <q}
 \underline{t} ^{\underline{r}}
F ^* ( a'  _{\underline{r},\underline{u},\underline{n},i})$
with 
$ a'  _{\underline{r},\underline{u},\underline{n},i} \in A ^{(s)}$ are such that 
$v _{p} ( a'  _{\underline{r},\underline{u},\underline{n},i}) \geq \frac{i+|\underline{n}|}{c} -1$,
for any $|\underline{u}|_{\infty} <q, \underline{n}\in \bbN ^d, i\in \bbN$.
Hence, 
$P _{\underline{r},\underline{u}}:=\sum _{\underline{n}\in \bbN ^d} 
  \sum _{i\in \bbN}
  \frac{1}{t ^{i}}
a'  _{\underline{r},\underline{u},\underline{n},i}
\underline{\partial} ^{<\underline{n}> _{(m)}} \in  \widetilde{D} ^{(m)\dag} _{\fX ^{\sharp (s)}/\fS ^{\sharp}}$.
We get
$P ^{(s)} _{\underline{r},\underline{u}}
=\sum _{\underline{n}\in \bbN ^d} 
  \sum _{i\in \bbN}
 \frac{1}{t ^{i}}
F ^* (a'  _{\underline{r},\underline{u},\underline{n},i})
(\underline{\partial} ^{<\underline{n}> _{(m)}} ) ^{(s)}$.
We compute 
\begin{gather*}
\Phi  \left ( \sum _{|\underline{r}| _{\infty} <q}  \sum _{|\underline{u}|_{\infty} <q}  \underline{t} ^{\underline{r}}  \otimes   P _{\underline{r},\underline{u}}   \otimes   H _{\underline{t}}  \underline{t} ^{-\underline{u}}     \right )
= \sum _{|\underline{r}| _{\infty} <q}  \sum _{|\underline{u}|_{\infty} <q} \sum _{\underline{n}\in \bbN ^d}    \sum _{i\in \bbN}    
 \underline{t} ^{\underline{r}}  \frac{1}{t ^{i}} F ^* (a'  _{\underline{r},\underline{u},\underline{n},i}) (\underline{\partial} ^{<\underline{n}> _{(m)}} ) ^{(s)}  H _{\underline{t}}  \underline{t} ^{-\underline{u}}
 \\  =  \sum _{|\underline{u}|_{\infty} <q} \sum _{\underline{n}\in \bbN ^d}    \sum _{i\in \bbN}   \sum _{\underline{k}\leq q \underline{n}+\underline{u}} 
  \frac{1}{t ^{i}} a  _{\underline{k},i} \beta _{\underline{k},q\underline{n}+\underline{u}} (\underline{\partial} ^{<\underline{n}> _{(m)}} ) ^{(s)}  H _{\underline{t}}  \underline{t} ^{-\underline{u}}   =P.
\end{gather*}

\end{proof}

\begin{coro} \label{lem-factorization2}
The sheaf $F ^* \widetilde{\cD} ^{(m)\dag} _{\fX ^{\sharp (s)}/\fS ^{\sharp}}$ has a canonical structure of  
$(\widetilde{\cD} ^{(m+s)\dag} _{\fX ^{\sharp}/\fS ^{\sharp}}, \widetilde{\cD} ^{(m)\dag} _{\fX ^{\sharp (s)}/\fS ^{\sharp}})$-subimodule of  $j _* F ^* _{\eta} \widehat{\cD} ^{(m)} _{\fX ^{(s)}_\eta /\eta}$.
The canonical homomorphism of left $\widehat{\cD} ^{(m+s)} _{\fX _\eta /\eta}$-modules $ \widehat{\cD} ^{(m+s)} _{\fX _\eta /\eta} \to  F ^* _{\eta} \widehat{\cD} ^{(m)} _{\fX ^{(s)}_\eta /\eta} $  (given by $P\mapsto P \cdot (1\otimes 1)$) 
induces the  homomorphism of left $\widetilde{\cD} ^{(m+s)\dag} _{\fX ^{\sharp}/\fS ^{\sharp}}$-modules
$\widetilde{\cD} ^{(m+s)\dag} _{\fX ^{\sharp}/\fS ^{\sharp}} \to  F ^* \widetilde{\cD} ^{(m)\dag} _{\fX ^{\sharp (s)}/\fS ^{\sharp}} $ making commutative the diagram
\begin{equation} \label{lem-factorization2diag}
\xymatrix{
 {F ^* F ^\flat  \widetilde{\cD} ^{(m)\dag} _{\fX ^{\sharp (s)}/\fS ^{\sharp}}}  \ar[r] _-{\ref{phidagdescrip}} ^-{\sim}  \ar@{->>}[d] ^-{}
 &  {\widetilde{\cD} ^{(m+s)\dag} _{\fX ^{\sharp}/\fS ^{\sharp}}}  \ar@{.>>}[d] ^-{}\ar@{^{(}->}[r] ^-{} &  {j _* \widehat{\cD} ^{(m+s)} _{\fX _\eta /\eta}} \ar@{->>}[d] ^-{}
\\  {F ^*  \widetilde{\cD} ^{(m)\dag} _{\fX ^{\sharp (s)}/\fS ^{\sharp}}}  \ar@{=}[r] ^-{} & {F ^* \widetilde{\cD} ^{(m)\dag} _{\fX ^{\sharp (s)}/\fS ^{\sharp}}}  \ar@{^{(}->}[r] ^-{} &  {j _* F ^* _{\eta} \widehat{\cD} ^{(m)} _{\fX ^{(s)}_\eta /\eta}.} 
}
\end{equation}
where the left vertical arrow is given by
$ \sum _{|\underline{r}| _{\infty} <q}  \sum _{|\underline{u}|_{\infty} <q}  \underline{t} ^{\underline{r}}  \otimes   P _{\underline{r},\underline{u}}   \otimes   H _{\underline{t}}  \underline{t} ^{-\underline{u}} 
\mapsto  \sum _{|\underline{r}| _{\infty} <q}   \underline{t} ^{\underline{r}}  \otimes   P _{\underline{r},\underline{0}}   $. 

The left $\widetilde{\cD} ^{(m+s)\dag} _{\fX ^{\sharp}/\fS ^{\sharp}}$-module $F ^* \widetilde{\cD} ^{(m)\dag} _{\fX ^{\sharp (s)}/\fS ^{\sharp}}$ is locally projective of finite type.
\end{coro}

\begin{proof}
The commutativity of the diagram is clear (e.g. this is a consequence of the beginning of the proof of \cite[2.5.3]{Be2}).  The right vertical arrow is given by 
$  \sum _{|\underline{r}| _{\infty} <q}  \sum _{|\underline{u}|_{\infty} <q}  \underline{t} ^{\underline{r}}  P _{\underline{r},\underline{u}} ^{(s)}  H _{\underline{t}}  \underline{t} ^{-\underline{u}}
\mapsto  \sum _{|\underline{r}| _{\infty} <q}   \underline{t} ^{\underline{r}}  \otimes   P _{\underline{r},\underline{0}}  $. 
Hence, we get the factorisation of the diagram \ref{lem-factorization2diag}.

Consider the section given by  
Let $s\colon F ^* \widetilde{\cD} ^{(m)\dag} _{\fX ^{\sharp (s)}/\fS ^{\sharp}} \to \widetilde{\cD} ^{(m+s)\dag} _{\fX ^{\sharp}/\fS ^{\sharp}} $ the map sending 
$ \sum _{|\underline{r}| _{\infty} <q}   \underline{t} ^{\underline{r}}   \otimes   P _{\underline{r}}  $ to $\sum _{|\underline{r}| _{\infty} <q}  \underline{t} ^{\underline{r}}  P _{\underline{r}} ^{(s)} $.
This a section of the projection $\widetilde{\cD} ^{(m+s)\dag} _{\fX ^{\sharp}/\fS ^{\sharp}} \to  F ^* \widetilde{\cD} ^{(m)\dag} _{\fX ^{\sharp (s)}/\fS ^{\sharp}} $.
Moreover, since the composition of the projection with the section   is  $P\mapsto P  H _{\underline{t}} $ (use the first equality of \ref{251iii2GDesc})  is $\widetilde{\cD} ^{(m+s)\dag} _{\fX ^{\sharp}/\fS ^{\sharp}}$-linear,  then so is the section $s$.
\end{proof}

\begin{empt}
By using the last equality of \ref{251iii2GDesc}, we get 
$H _{\underline{t}}  \cdot  \sum _{|\underline{r}| _{\infty} <q}  \sum _{|\underline{u}|_{\infty} <q}  \underline{t} ^{\underline{r}}  P _{\underline{r},\underline{u}} ^{(s)}  H _{\underline{t}}  \underline{t} ^{-\underline{u}}=
\sum _{|\underline{u}|_{\infty} <q}  P _{\underline{0},\underline{u}} ^{(s)}  H _{\underline{t}}  \underline{t} ^{-\underline{u}}$.
By $\widetilde{\cD} ^{(m+s)\dag} _{\fX ^{\sharp}/\fS ^{\sharp}}$-linearity of the homomorphism $\widetilde{\cD} ^{(m+s)\dag} _{\fX ^{\sharp}/\fS ^{\sharp}} \to  F ^* \widetilde{\cD} ^{(m)\dag} _{\fX ^{\sharp (s)}/\fS ^{\sharp}} $ , this yields
\begin{equation}
\label{331Garnier-descente-frov-explbisform0} H _{\underline{t}}  \cdot   \sum _{|\underline{r}| _{\infty} <q}   \underline{t} ^{\underline{r}}  \otimes   P _{\underline{r},\underline{0}}   = 
1 \otimes   P _{\underline{0},\underline{0}}  .
\end{equation}

For any $P,Q\in  \widetilde{\cD} ^{(m)\dag} _{\fX ^{\sharp (s)}/\fS ^{\sharp}} $, by linearity we get that the homomorphism $\widetilde{\cD} ^{(m+s)\dag} _{\fX ^{\sharp}/\fS ^{\sharp}} \to  F ^* \widetilde{\cD} ^{(m)\dag} _{\fX ^{\sharp (s)}/\fS ^{\sharp}} $ sends 
$ P  ^{(s)} Q ^{(s)}$  to  $ 1 \otimes PQ$ and by linearity to $ P  ^{(s)}\cdot (1 \otimes Q) $. Hence:
\begin{equation}\label{FormP(s)otimes}
P  ^{(s)}\cdot (1 \otimes Q) =1 \otimes PQ.
\end{equation}
This yields that the map  $ \widetilde{\cD} ^{(m)\dag} _{\fX ^{\sharp (s)}/\fS ^{\sharp}}  \to  F ^* \widetilde{\cD} ^{(m)\dag} _{\fX ^{\sharp (s)}/\fS ^{\sharp}} $ is $ \widetilde{\cD} ^{(m)\dag} _{\fX ^{\sharp (s)}/\fS ^{\sharp}} $-bilinear,
where the structure of left $ \widetilde{\cD} ^{(m)\dag} _{\fX ^{\sharp (s)}/\fS ^{\sharp}} $-module of $F ^* \widetilde{\cD} ^{(m)\dag} _{\fX ^{\sharp (s)}/\fS ^{\sharp}} $ is induced by the ring homomorphism
$ \widetilde{\cD} ^{(m)\dag} _{\fX ^{\sharp (s)}/\fS ^{\sharp}} \to F _* \widetilde{\cD} ^{(m+s)\dag} _{\fX ^{\sharp}/\fS ^{\sharp}}$ given by $P \mapsto P ^{(s)}$.
\end{empt}

\begin{empt} \label{dfn-FflatF*0}
Let $\cE$ be a left $\widetilde{\cD} ^{(m)\dag} _{\fX ^{\sharp (s)}/\fS ^{\sharp}}$-module.  We get the left $\widetilde{\cD} ^{(m+s)\dag} _{\fX ^{\sharp}/\fS ^{\sharp}}$-module
$$F ^*  \cE :=  F ^*  \widetilde{\cD} ^{(m)\dag} _{\fX ^{\sharp (s)}/\fS ^{\sharp}}  \otimes _{ \widetilde{\cD} ^{(m)\dag} _{\fX ^{\sharp (s)}/\fS ^{\sharp}}}  \E. $$
\end{empt}

\begin{prop}\label{331Garnier-descente-frov-expl}
Let $\cE ^{(m+s)}$ be a $\widetilde{\cD} ^{(m+s)\dag} _{\fX ^{\sharp}/\fS ^{\sharp}}$-module. Set 
\begin{equation}\label{331Garnier-descente-frov-expl1} \cE ^{(m)}\coloneqq \mathrm{Im} ( F_* \cE ^{(m+s)} \overset{ H _{\underline{t}}\cdot}{\longrightarrow} F _* \cE ^{(m+s)}).
\end{equation}
Then $\cE ^{(m)}$ is a $\widetilde{\cD} ^{(m+s)\dag} _{\fX ^{\sharp}/\fS ^{\sharp}}$-module and the map 
\begin{equation}\label{331Garnier-descente-frov-expl2} \Phi \colon F ^*  \cE ^{(m)} \to \cE ^{(m+s)}   
\end{equation}
given by $a \otimes H _{\underline{t}} (e)\mapsto a H _{\underline{t}} (e)$, is a $\widetilde{\cD} ^{(m+s)\dag} _{\fX ^{\sharp}/\fS ^{\sharp}}$-linear isomorphism. The inverse is given by the map
\begin{equation}\label{331Garnier-descente-frov-expl3} \phi ^{-1}\colon e \mapsto \sum _{|\underline{r}| _{\infty} <q} \underline{t} ^{\underline{r}} \otimes H _{\underline{t}}  \underline{t} ^{-\underline{r}} (e).
\end{equation}
\end{prop}

\begin{proof}1) Let us prove the morphism \ref{331Garnier-descente-frov-expl2} is $\widetilde{\cD} ^{(m+s)\dag} _{\fX ^{\sharp}/\fS ^{\sharp}}$-linear. 
Since it is $\cO _{\fX}$-linear, then we reduce to check $P\cdot \Phi ( 1 \otimes 1 \otimes H _{\underline{t}}  e) =\Phi ( P \cdot  (1 \otimes 1 \otimes H _{\underline{t}}  e) )$ for any $P \in \widetilde{\cD} ^{(m+s)\dag} _{\fX ^{\sharp}/\fS ^{\sharp}}$, 
$e \in \cE ^{(m+s)}$. By additivity, we can suppose  $P = \underline{t} ^{\underline{r}}  P _{\underline{r},\underline{u}} ^{(s)}  H _{\underline{t}}  \underline{t} ^{-\underline{u}}$,
for $P _{\underline{r},\underline{u}}\in  \widetilde{\cD} ^{(m)\dag} _{\fX ^{\sharp (s)}/\fS ^{\sharp}}$, $|\underline{r}| _{\infty} <q$ and $|\underline{u}|_{\infty} <q$.
We get $P\cdot \Phi ( 1 \otimes 1 \otimes H _{\underline{t}}  e) = \underline{t} ^{\underline{r}}  P _{\underline{r},\underline{u}} ^{(s)}  H _{\underline{t}}  \underline{t} ^{-\underline{u}}  H _{\underline{t}}  e) 
\overset{\ref{251iii2GDesc}}{=} \delta _{\underline{u}, \underline{0}} \underline{t} ^{\underline{r}}  P _{\underline{r},\underline{u}} ^{(s)}H _{\underline{t}}  e$. On the other hand 
$ P \cdot  (1 \otimes 1 \otimes H _{\underline{t}}  e) =  \underline{t} ^{\underline{r}}  P _{\underline{r},\underline{u}} ^{(s)}  H _{\underline{t}}  \underline{t} ^{-\underline{u}}    \cdot  (1 \otimes 1 )  \otimes H _{\underline{t}}  e
=  \delta _{\underline{u}, \underline{0}} ( \underline{t} ^{\underline{r}}  \otimes P _{\underline{r},\underline{u}} )     \otimes H _{\underline{t}}  e
=  \delta _{\underline{u}, \underline{0}} \underline{t} ^{\underline{r}}  P _{\underline{r},\underline{u}} ^{(s)} H _{\underline{t}}   e$.

2) Using \ref{251iii2GDesc1}, we prove the map given by the formula \ref{331Garnier-descente-frov-expl3} is an inverse. Hence, we are done.
\end{proof}

\begin{prop}\label{331Garnier-descente-frov-explbis}
Let $\cE ^{(m)}$ be a $ \widetilde{\cD} ^{(m)\dag} _{\fX ^{\sharp (s)}/\fS ^{\sharp}}$-module.  For any $\underline{r}\in \bbN ^d$ such that  $|\underline{r}| _{\infty} <q$, any $e \in  \cE ^{(m)}$, we have 
$H _{\underline{t}}  ( \underline{t} ^{\underline{r}} \otimes e) =\delta _{\underline{r}, \underline{0}} ( \underline{t} ^{\underline{r}} \otimes e) .$ We have the isomorphism of $ \widetilde{\cD} ^{(m)\dag} _{\fX ^{\sharp (s)}/\fS ^{\sharp}}$-modules:
\begin{equation}\label{331Garnier-descente-frov-expl1bis} 
\cE ^{(m)}\riso  \mathrm{Im} ( F _*  F^* \cE ^{(m)} \overset{ H _{\underline{t}}\cdot}{\longrightarrow} F _*  F ^* \cE ^{(m)}),
\end{equation} 
given by $e \mapsto 1 \otimes e$.
\end{prop}

\begin{proof}By construction, we reduce to the case $\cE ^{(m)}=\widetilde{\cD} ^{(m)\dag} _{\fX ^{\sharp (s)}/\fS ^{\sharp}}$.  This is therefore a consequence of the formula \ref{331Garnier-descente-frov-explbisform0}.
\end{proof}

\begin{empt} \label{eqcatFrobdescexpl2pre} 
The functor $\cE ^{(m)} \mapsto F ^* \cE ^{(m)}$  from the category of  left  $\widetilde{\cD} ^{(m)\dag} _{\fX ^{\sharp (s)}/\fS ^{\sharp}}$-modules to that of  left 
$\widetilde{\cD} ^{(m+s)\dag} _{\fX ^{\sharp}/\fS ^{\sharp}}$-modules is an exact equivalence of categories. Locally, a quasi-inverse functor is given by
$\cE ^{(m+s)} \mapsto  \mathrm{Im} ( F_* \cE ^{(m+s)} \overset{ H _{\underline{t}}\cdot}{\longrightarrow} F _* \cE ^{(m+s)})) $.
For any left $\widetilde{\cD} ^{(m)\dag} _{\fX ^{\sharp (s)}/\fS ^{\sharp}}$-modules $\cE$ and $\cF$, the canonical morphism
\begin{equation} \label{Hom-F*const1pre}
\mathcal{H} om _{\widetilde{\cD} ^{(m)\dag} _{\fX ^{\sharp (s)}/\fS ^{\sharp}}} (\cE,\cF) \to  \mathcal{H} om _{\widetilde{\cD} ^{(m+s)\dag} _{\fX ^{\sharp}/\fS ^{\sharp}}} (F ^* \cE, F ^*\cF )
\end{equation}
is an isomorphism.
\end{empt}

\begin{prop} \label{FflatDlocproj}
Suppose $H \in \widetilde{\cD} ^{(m+s)\dag} _{\fX ^{\sharp}/\fS ^{\sharp}}$ (e.g. if $m\geq 1$: see \ref{HinDs+1}).  
The sheaf $F ^\flat  \widetilde{\cD} ^{(m)\dag} _{\fX ^{\sharp (s)}/\fS ^{\sharp}}$ is  a $(\widetilde{\cD} ^{(m)\dag} _{\fX ^{\sharp (s)}/\fS ^{\sharp}}, \widetilde{\cD} ^{(m+s)\dag} _{\fX ^{\sharp}/\fS ^{\sharp}})$-sub-bimodule of  
$j _* F ^\flat _\eta \widehat{\cD} ^{(m)} _{\fX ^{(s)}_\eta /\eta}$. Moreover,  $F ^\flat  \widetilde{\cD} ^{(m)\dag} _{\fX ^{\sharp (s)}/\fS ^{\sharp}}$ is  locally projective of finite type  as $\widetilde{\cD} ^{(m+s)\dag} _{\fX ^{\sharp}/\fS ^{\sharp}}$-module.
\end{prop}

\begin{proof}
Following \ref{phidagdescrip},  $F ^* F ^\flat  \widetilde{\cD} ^{(m)\dag} _{\fX ^{\sharp (s)}/\fS ^{\sharp}}$ is  a $\widetilde{\cD} ^{(m+s)\dag} _{\fX ^{\sharp}/\fS ^{\sharp}}$-subbimodule of  $j _* F ^* _\eta F ^\flat _\eta \widehat{\cD} ^{(m)} _{\fX ^{(s)}_\eta /\eta}$.
It follows from \ref{331Garnier-descente-frov-explbis} (resp. \cite[3.3.1]{Garnier-descente-frov-expl})
that  the image by the action on the left by $H _{\underline{t}}$ on  $F ^* F ^\flat  \widetilde{\cD} ^{(m)\dag} _{\fX ^{\sharp (s)}/\fS ^{\sharp}}$  is  $F ^\flat  \widetilde{\cD} ^{(m)\dag} _{\fX ^{\sharp (s)}/\fS ^{\sharp}}$
(resp. on  $j _{*}F ^* F ^\flat  \widetilde{\cD} ^{(m)\dag} _{\fX _{\eta} ^{(s)}/\eta}$  is  $j _{*} F ^\flat  \widetilde{\cD} ^{(m)\dag} _{\fX _\eta  (s)}/\eta$), we get the first assertion.
By composition, we get  the canonical $\widetilde{\cD} ^{(m+s)\dag} _{\fX ^{\sharp}/\fS ^{\sharp}}$-linear (and even $(\widetilde{\cD} ^{(m)\dag} _{\fX ^{\sharp (s)}/\fS ^{\sharp}}, \widetilde{\cD} ^{(m+s)\dag} _{\fX ^{\sharp}/\fS ^{\sharp}})$-bilinear) homomorphism
$F ^\flat  \widetilde{\cD} ^{(m)\dag} _{\fX ^{\sharp (s)}/\fS ^{\sharp}} \to  F ^* F ^\flat  \widetilde{\cD} ^{(m)\dag} _{\fX ^{\sharp (s)}/\fS ^{\sharp}} \riso  \widetilde{\cD} ^{(m+s)\dag} _{\fX ^{\sharp}/\fS ^{\sharp}}$
sending $\sum _{|\underline{u}|_{\infty} <q}  P _{\underline{u}} \otimes  H _{\underline{t}}  \underline{t} ^{-\underline{u}}$ to  $\sum _{|\underline{u}|_{\infty} <q}  P _{\underline{u}} ^{(s)}  \underline{t} ^{-\underline{u}}$.
We have the retraction  $\widetilde{\cD} ^{(m+s)\dag} _{\fX ^{\sharp}/\fS ^{\sharp}}  \to  F ^\flat  \widetilde{\cD} ^{(m)\dag} _{\fX ^{\sharp (s)}/\fS ^{\sharp}}$
given by  $\sum _{|\underline{r}| _{\infty} <q}  \sum _{|\underline{u}|_{\infty} <q}  \underline{t} ^{\underline{r}}  P _{\underline{r},\underline{u}} ^{(s)}    \underline{t} ^{-\underline{u}}
\mapsto \sum _{|\underline{u}|_{\infty} <q}  P _{\underline{0},\underline{u}} \otimes  H _{\underline{t}}  \underline{t} ^{-\underline{u}}$. The composition of this retraction followed by the canonical inclusion is given by 
$\sum _{|\underline{r}| _{\infty} <q}  \sum _{|\underline{u}|_{\infty} <q}  \underline{t} ^{\underline{r}}  P _{\underline{r},\underline{u}} ^{(s)}    \underline{t} ^{-\underline{u}}\mapsto 
\sum _{|\underline{u}|_{\infty} <q}  P _{\underline{r},\underline{u}} ^{(s)}    \underline{t} ^{-\underline{u}}$, which is $\widetilde{\cD} ^{(m+s)\dag} _{\fX ^{\sharp}/\fS ^{\sharp}}$-linear since this is the multiplication on the left by $H _{\underline{t}}$ (use the last formula
of \ref{251iii2GDesc}).
Hence, so is the retraction and we are done. 
\end{proof}

\begin{empt} \label{dfn-FflatF*}
Let $\cM$ be a right  $\widetilde{\cD} ^{(m)\dag} _{\fX ^{\sharp (s)}/\fS ^{\sharp}}$-module.  We have the isomorphism of  $\widetilde{\cO} _{\fX}$-modules, which justifies the notation
$$\cM \otimes _{ \widetilde{\cD} ^{(m)\dag} _{\fX ^{\sharp (s)}/\fS ^{\sharp}}}  F ^\flat  \widetilde{\cD} ^{(m)\dag} _{\fX ^{\sharp (s)}/\fS ^{\sharp}}
= \cM \otimes _{ \widetilde{\cD} ^{(m)\dag} _{\fX ^{\sharp (s)}/\fS ^{\sharp}}}  \mathcal{H} om _{\cO _{\fX ^{(s)}}}  (\cO _{\fX}, \widetilde{\cD} ^{(m)\dag} _{\fX ^{\sharp (s)}/\fS ^{\sharp}}) 
\riso  \mathcal{H} om _{\cO _{\fX ^{ (s)}}} (\cO _{\fX},  \cM)  =  F ^\flat (\cM).$$ 
Hence, we can define without problem of notation the right  $\widetilde{\cD} ^{(m+s)\dag} _{\fX ^{\sharp}/\fS ^{\sharp}}$-module 
$$F ^\flat  \cM := \cM \otimes _{ \widetilde{\cD} ^{(m)\dag} _{\fX ^{\sharp (s)}/\fS ^{\sharp}}} F ^\flat  \widetilde{\cD} ^{(m)\dag} _{\fX ^{\sharp (s)}/\fS ^{\sharp}}.$$
\end{empt}

\begin{empt}
\label{dfn-FflatF*coh}
Since $\widetilde{\cD} ^{(m)\dag} _{\fX ^{\sharp (s)}/\fS ^{\sharp}}$ is right and left coherent 
(see \ref{123-coherent}), 
we get from \ref{FflatDlocproj} (resp. \ref{lem-factorization2})
the $\widetilde{\cD} ^{(m+s)\dag} _{\fX ^{\sharp}/\fS ^{\sharp}}$-coherence of 
$F ^\flat  \widetilde{\cD} ^{(m)\dag} _{\fX ^{\sharp (s)}/\fS ^{\sharp}}$
(resp. $F ^*  \widetilde{\cD} ^{(m)\dag} _{\fX ^{\sharp (s)}/\fS ^{\sharp}}$).
Since the functors $F ^* $ and $F ^\flat$ are flat, 
we get the functor
$\cE \mapsto F ^* \cE$ 
(resp. 
$\cM \mapsto F ^\flat \cM$)
from the category of coherent 
left (resp. right) $\widetilde{\cD} ^{(m)\dag} _{\fX ^{\sharp (s)}/\fS ^{\sharp}}$-modules
to that of 
coherent 
left (resp. right) 
$\widetilde{\cD} ^{(m+s)\dag} _{\fX ^{\sharp}/\fS ^{\sharp}}$-modules.
\end{empt}

\begin{lem} \label{lem254Be2}
We have the canonical  isomorphism of $(\widetilde{\cD} ^{(m+s)\dag} _{\fX ^{\sharp}/\fS ^{\sharp}}, \widetilde{\cD} ^{(m)\dag} _{\fX ^{\sharp (s)}/\fS ^{\sharp}})$-bimodules
\begin{equation} \label{lem254Be2iso}
\mathcal{H} om _{\widetilde{\cD} ^{(m)\dag} _{\fX ^{\sharp (s)}/\fS ^{\sharp}}} (F ^\flat \widetilde{\cD} ^{(m)\dag} _{\fX ^{\sharp (s)}/\fS ^{\sharp}}, \widetilde{\cD} ^{(m)\dag} _{\fX ^{\sharp (s)}/\fS ^{\sharp}} )
\riso F ^* \widetilde{\cD} ^{(m)\dag} _{\fX ^{\sharp (s)}/\fS ^{\sharp}}.
\end{equation}
\end{lem}

\begin{proof}
By functoriality, we have the homomorphisms
of $(\widetilde{\cD} ^{(m+s)\dag} _{\fX ^{\sharp}/\fS ^{\sharp}}, \widetilde{\cD} ^{(m)\dag} _{\fX ^{\sharp (s)}/\fS ^{\sharp}})$-bimodules
\begin{equation} \label{Hom-F*const1}
\mathcal{H} om _{\widetilde{\cD} ^{(m)\dag} _{\fX ^{\sharp (s)}/\fS ^{\sharp}}}
(F ^\flat \widetilde{\cD} ^{(m)\dag} _{\fX ^{\sharp (s)}/\fS ^{\sharp}},
\widetilde{\cD} ^{(m)\dag} _{\fX ^{\sharp (s)}/\fS ^{\sharp}}
)
\underset{\ref{Hom-F*const1pre}}{\riso} 
\mathcal{H} om _{\widetilde{\cD} ^{(m+s)\dag} _{\fX ^{\sharp}/\fS ^{\sharp}}}
(F ^* F ^\flat \widetilde{\cD} ^{(m)\dag} _{\fX ^{\sharp (s)}/\fS ^{\sharp}},
F ^*\widetilde{\cD} ^{(m)\dag} _{\fX ^{\sharp (s)}/\fS ^{\sharp}}
)
\underset{\ref{phidagdescrip}}{\riso} 
F ^* \widetilde{\cD} ^{(m)\dag} _{\fX ^{\sharp (s)}/\fS ^{\sharp}}.
\end{equation}
\end{proof}

\begin{empt}\label{comp-lem254Be2}
Locally, $F ^\flat \widetilde{\cD} ^{(m)\dag} _{\fX ^{\sharp (s)}/\fS ^{\sharp}}=\mathcal{H} om _{\cO _{\fX ^{(s)}}}  (\cO _{\fX}, \widetilde{\cD} ^{(m)\dag} _{\fX ^{\sharp (s)}/\fS ^{\sharp}}) 
\riso  \widetilde{\cD} ^{(m)\dag} _{\fX ^{\sharp (s)}/\fS ^{\sharp}} \otimes _{\cO _{\fX ^{(s)}}}  
\mathcal{H} om _{\cO _{\fX ^{(s)}}}  (\cO _{\fX}, \cO _{\fX ^{(s)}}) $ is a free left $\widetilde{\cD} ^{(m)\dag} _{\fX ^{\sharp (s)}/\fS ^{\sharp}}$ with the basis
$\{1 \otimes  H _{\underline{t}}  \underline{t} ^{-\underline{r}},   \text{  $\underline{r} \in \bbN ^d$ such that  $|\underline{r}|_{\infty} <q$}\}$.
We get the dual basis $\{( 1 \otimes  H _{\underline{t}}  \underline{t} ^{-\underline{r}}) ^*,   \text{  $\underline{r} \in \bbN ^d$ such that  $|\underline{r}|_{\infty} <q$}\}$ of 
$\mathcal{H} om _{\widetilde{\cD} ^{(m)\dag} _{\fX ^{\sharp (s)}/\fS ^{\sharp}}} (F ^\flat \widetilde{\cD} ^{(m)\dag} _{\fX ^{\sharp (s)}/\fS ^{\sharp}}, \widetilde{\cD} ^{(m)\dag} _{\fX ^{\sharp (s)}/\fS ^{\sharp}} )$.
Using the formula \ref{251iii2GDesc1}, we compute the isomorphism \ref{lem254Be2iso} sends $( 1 \otimes  H _{\underline{t}}  \underline{t} ^{-\underline{r}}) ^*$ to $\underline{t} ^{\underline{r}} \otimes 1$.
\end{empt}

\begin{prop}
\label{255Be2}
\begin{enumerate}[(a)]
\item We have the canonical isomorphism of $\widetilde{\cD} ^{(m)\dag} _{\fX ^{\sharp (s)}/\fS ^{\sharp}}$-bimodules
\begin{equation} \label{255Be2iso}
F ^\flat \widetilde{\cD} ^{(m)\dag} _{\fX ^{\sharp (s)}/\fS ^{\sharp}} \otimes _{\widetilde{\cD} ^{(m+s)\dag} _{\fX ^{\sharp}/\fS ^{\sharp}}}F ^* \widetilde{\cD} ^{(m)\dag} _{\fX ^{\sharp (s)}/\fS ^{\sharp}} \riso    \widetilde{\cD} ^{(m)\dag} _{\fX ^{\sharp (s)}/\fS ^{\sharp}}.
\end{equation}

\item For any integer $n \geq 1$, 
\begin{equation}
\label{255Be2Tor=0}
\mathcal{T} or ^n _{\widetilde{\cD} ^{(m+s)\dag} _{\fX ^{\sharp}/\fS ^{\sharp}}}
(F ^\flat \widetilde{\cD} ^{(m)\dag} _{\fX ^{\sharp (s)}/\fS ^{\sharp}}
,
F ^* \widetilde{\cD} ^{(m)\dag} _{\fX ^{\sharp (s)}/\fS ^{\sharp}})
=0.
\end{equation}

\end{enumerate}

\end{prop}

\begin{proof}
We have the canonical homomorphism of $\widetilde{\cD} ^{(m)\dag} _{\fX ^{\sharp (s)}/\fS ^{\sharp}}$-bimodules
\begin{equation} \label{255Be2isoproof1}
F ^\flat \widetilde{\cD} ^{(m)\dag} _{\fX ^{\sharp (s)}/\fS ^{\sharp}} \otimes _{\widetilde{\cD} ^{(m+s)\dag} _{\fX ^{\sharp}/\fS ^{\sharp}}} \mathcal{H} om _{\widetilde{\cD} ^{(m)\dag} _{\fX ^{\sharp (s)}/\fS ^{\sharp}}}
(F ^\flat \widetilde{\cD} ^{(m)\dag} _{\fX ^{\sharp (s)}/\fS ^{\sharp}},\widetilde{\cD} ^{(m)\dag} _{\fX ^{\sharp (s)}/\fS ^{\sharp}} )\to  \widetilde{\cD} ^{(m)\dag} _{\fX ^{\sharp (s)}/\fS ^{\sharp}}.
\end{equation}
Using \ref{lem254Be2}, this correspond to a homomorphism of the form \ref{255Be2iso}.  To check this is an isomorphism, we reduce to prove it after applying the functor $F ^*$. 
 We get the map 
\begin{equation}\label{255Be2isoproof2}
F ^{*} F ^\flat \widetilde{\cD} ^{(m)\dag} _{\fX ^{\sharp (s)}/\fS ^{\sharp}} \otimes _{\widetilde{\cD} ^{(m+s)\dag} _{\fX ^{\sharp}/\fS ^{\sharp}}}F ^* \widetilde{\cD} ^{(m)\dag} _{\fX ^{\sharp (s)}/\fS ^{\sharp}}
\to  F ^{*}  \widetilde{\cD} ^{(m)\dag} _{\fX ^{\sharp (s)}/\fS ^{\sharp}}.
\end{equation}
Using \ref{phidagdescrip}, this ocrresponds to a map of the form 
\begin{equation}\label{255Be2isoproof3}
\widetilde{\cD} ^{(m+s)\dag} _{\fX ^{\sharp}/\fS ^{\sharp}} \otimes _{\widetilde{\cD} ^{(m+s)\dag} _{\fX ^{\sharp}/\fS ^{\sharp}}}F ^* \widetilde{\cD} ^{(m)\dag} _{\fX ^{\sharp (s)}/\fS ^{\sharp}}
\to  F ^{*}  \widetilde{\cD} ^{(m)\dag} _{\fX ^{\sharp (s)}/\fS ^{\sharp}}.
\end{equation}
The fact this is an isomorphism is local.  With notation \ref{comp-lem254Be2}, we compute the map \ref{255Be2isoproof1} sends 
$\sum _{|\underline{u}|_{\infty} <q} \sum _{|\underline{v}| _{\infty} <q}   (P _{\underline{u}}   \otimes   H _{\underline{t}}  \underline{t} ^{-\underline{u}}) \otimes  Q _{\underline{v}} ( 1   \otimes   H _{\underline{t}}  \underline{t} ^{-\underline{v}}) ^* $
to $\sum _{|\underline{u}|_{\infty} <q} \sum _{|\underline{v}| _{\infty} <q}   \delta _{\underline{u}, \underline{v}} Q _{\underline{v}} P _{\underline{u}}   $.
The map \ref{255Be2iso} sends 
$\sum _{|\underline{u}|_{\infty} <q} \sum _{|\underline{v}| _{\infty} <q}   (P _{\underline{u}}   \otimes   H _{\underline{t}}  \underline{t} ^{-\underline{u}}) \otimes ( \underline{t} ^{\underline{v}} \otimes Q _{\underline{v}})$
to $\sum _{|\underline{u}|_{\infty} <q} \sum _{|\underline{v}| _{\infty} <q}   \delta _{\underline{u}, \underline{v}} Q _{\underline{v}} P _{\underline{u}}   $.
Hence, the map \ref{255Be2isoproof2} sends 
$\sum _{|\underline{r}|_{\infty} <q} \sum _{|\underline{u}|_{\infty} <q} \sum _{|\underline{v}| _{\infty} <q}   (\underline{t} ^{\underline{r}} \otimes 
P _{\underline{u}}   \otimes   H _{\underline{t}}  \underline{t} ^{-\underline{u}}) \otimes ( \underline{t} ^{\underline{v}} \otimes Q _{\underline{v}})$
to $\sum _{|\underline{r}|_{\infty} <q}  \sum _{|\underline{u}|_{\infty} <q} \sum _{|\underline{v}| _{\infty} <q}   \delta _{\underline{u}, \underline{v}} \, \underline{t} ^{\underline{r}} \otimes Q _{\underline{v}} P _{\underline{u}}   $.
Using \ref{251iii2GDesc1}, we get the map \ref{255Be2isoproof3} send 
$\sum _{|\underline{v}| _{\infty} <q} 1 \otimes (\underline{t} ^{\underline{v}} \otimes Q _{\underline{v}}) $ to $ \sum _{|\underline{v}| _{\infty} <q} \underline{t} ^{\underline{v}} \otimes Q _{\underline{v}}$, i.e.,
\ref{255Be2isoproof3} is the canonical morphism which is an isomorphism. 
\end{proof}

\begin{coro}
\label{eqcatFrobdescexpl}
Let $\cM$ be a right 
$\widetilde{\cD} ^{(m)\dag} _{\fX ^{\sharp (s)}/\fS ^{\sharp}}$-module. 
Let $\cE$ be a left
$\widetilde{\cD} ^{(m)\dag} _{\fX ^{\sharp (s)}/\fS ^{\sharp}}$-module. 

\begin{enumerate}[(a)]
\item We have the $\widetilde{\cD} ^{(m)\dag} _{\fX ^{\sharp (s)}/\fS ^{\sharp}}$-linear isomorphism
\begin{equation}\label{eqcatFrobdescexpl-iso1}
F ^\flat \widetilde{\cD} ^{(m)\dag} _{\fX ^{\sharp (s)}/\fS ^{\sharp}}
\otimes ^\L _{\widetilde{\cD} ^{(m+s)\dag} _{\fX ^{\sharp}/\fS ^{\sharp}}}
F ^* \E
\riso 
F ^\flat \widetilde{\cD} ^{(m)\dag} _{\fX ^{\sharp (s)}/\fS ^{\sharp}}
\otimes _{\widetilde{\cD} ^{(m+s)\dag} _{\fX ^{\sharp}/\fS ^{\sharp}}}
F ^* \E
\riso 
\E.
\end{equation}

\item We have the $\widetilde{\cD} ^{(m)\dag} _{\fX ^{\sharp (s)}/\fS ^{\sharp}}$-linear isomorphism
\begin{equation}\label{eqcatFrobdescexpl-iso2}
F ^\flat \M
\otimes ^\L _{\widetilde{\cD} ^{(m+s)\dag} _{\fX ^{\sharp}/\fS ^{\sharp}}}
F ^* \widetilde{\cD} ^{(m)\dag} _{\fX ^{\sharp (s)}/\fS ^{\sharp}}
\riso 
F ^\flat \M
\otimes  _{\widetilde{\cD} ^{(m+s)\dag} _{\fX ^{\sharp}/\fS ^{\sharp}}}
F ^* \widetilde{\cD} ^{(m)\dag} _{\fX ^{\sharp (s)}/\fS ^{\sharp}}
\riso 
\M.
\end{equation}

\end{enumerate}

\end{coro}

\begin{coro} \label{eqcatFrobdescexpl2} 
The functor $\cE \to F ^* \cE$  (resp.  $\cM \to F ^\flat \cM$) from the category of  left (resp. right) $\widetilde{\cD} ^{(m)\dag} _{\fX ^{\sharp (s)}/\fS ^{\sharp}}$-modules to that of  left (resp. right) 
$\widetilde{\cD} ^{(m+s)\dag} _{\fX ^{\sharp}/\fS ^{\sharp}}$-modules is an exact equivalence of categories. A quasi-inverse functor is given by
$\cE '\mapsto  F ^\flat \widetilde{\cD} ^{(m)\dag} _{\fX ^{\sharp (s)}/\fS ^{\sharp}} \otimes _{\widetilde{\cD} ^{(m+s)\dag} _{\fX ^{\sharp}/\fS ^{\sharp}}} \cE '$
(resp.  $\cM ' \mapsto   \cM ' \otimes  _{\widetilde{\cD} ^{(m+s)\dag} _{\fX ^{\sharp}/\fS ^{\sharp}}} F ^* \widetilde{\cD} ^{(m)\dag} _{\fX ^{\sharp (s)}/\fS ^{\sharp}}$). 
\end{coro}

\begin{proof}
This is a consequence of \ref{phidagdescrip} and of \ref{255Be2iso}.
\end{proof}

\begin{lem} \label{Fflat-omega}
There exists a canonical isomorphism of right  $\widetilde{\cD} ^{(m+s)\dag} _{\fX ^{\sharp}/\fS ^{\sharp}}$-modules of the form 
\begin{equation} \label{mu-Fflat-omega}
\mu _{\fX ^{\sharp}/\fS ^{\sharp}} \colon  \widetilde{\omega} _{\fX ^\sharp /\fS ^\sharp} \riso  F ^\flat \widetilde{\omega} _{\fX ^{\sharp (s)}/\fS ^{\sharp}}.
\end{equation}
\end{lem}

\begin{proof}
Let $f\colon \fX ^\sharp \to \fS ^\sharp$ and  $g\colon \fX ^\sharp \to \fS ^\sharp$ be the structural morphisms. Let $d$ be the dimension of $X$ and $X ^{(s)}$. By definition (see \ref{dfnfsharp}), we get
$\widetilde{\omega} _{\fX ^\sharp /\fS ^\sharp} \riso f ^\sharp (\widetilde{\cO} _\fS)   [-d ]$,  and  $\widetilde{\omega} _{\fX ^{\sharp (s)} /\fS ^\sharp} \riso g ^\sharp (\widetilde{\cO} _\fS)   [-d ]$ .
By using Proposition \ref{prop-III.8.2},  we get the isomorphism of $\widetilde{\cO} _{\fX}$-modules $\mu _{\fX ^{\sharp}/\fS ^{\sharp}} \colon  \widetilde{\omega} _{\fX ^\sharp /\fS ^\sharp}
\riso  F ^\flat \widetilde{\omega} _{\fX ^{\sharp (s)}/\fS ^{\sharp}}$. Following \ref{omega-Dmodule} and \ref{dfn-FflatF*}, the terms of  \ref{mu-Fflat-omega} are   endowed with a canonical right $\widetilde{\cD} ^{(m) \dag}  _{\fX^\sharp/\fS^\sharp}$-module structure
Since the extension  $\widetilde{\cD} ^{(m+s) \dag}  _{\fX^\sharp/\fS^\sharp} \to  j _* \widehat{\cD} ^{(m+s)}  _{\fX_\eta/\eta}$ is fully faithful (see \ref{lemfullfaith}),  the vertical arrows 
\begin{equation} \notag
\xymatrix{   {\widetilde{\omega} _{\fX ^\sharp /\fS ^\sharp}}  \ar[r] ^-{\sim}   \ar@{^{(}->}[d] ^-{} &  {F ^\flat \widetilde{\omega} _{\fX ^{\sharp (s)} /\fS ^\sharp}}  \ar@{^{(}->}[d] ^-{}
\\  {j _* \omega _{\fX _\eta /\eta}} \ar[r] ^-{\sim}  &  {j _* F ^\flat \omega _{\fX _\eta ^{(s)}/\eta} }   }
\end{equation}
are injective. Hence, we reduce to check  the  $\widehat{\cD} ^{(m+s)}  _{\fX_\eta/\eta}$-linearity of  the canonical isomorphism $\omega _{\fX _\eta /\eta} \riso  F ^\flat \omega _{\fX _\eta ^{(s)}/\eta} $, 
which is checked from \cite[2.4.2]{Be2} by taking projective limits. 
\end{proof}

\begin{prop} \label{prop-omegaFflatF*}
Let $\cM$ be a right  $\widetilde{\cD} ^{(m)\dag} _{\fX ^{\sharp (s)}/\fS ^{\sharp}}$-module.  Let $\cE$ be a left $\widetilde{\cD} ^{(m)\dag} _{\fX ^{\sharp (s)}/\fS ^{\sharp}}$-module. 

\begin{enumerate}[(a)]
\item We have the canonical isomorphism of  right  $\widetilde{\cD} ^{(m)\dag} _{\fX ^{\sharp}/\fS ^{\sharp}}$-modules
\begin{equation} \label{prop-omegaFflatF*1}
\mu _{\E} \colon  \widetilde{\omega} _{\fX ^\sharp /\fS ^\sharp} \otimes _{\widetilde{\cO} _{\fX}} F ^* (\cE)  \riso F ^\flat ( \widetilde{\omega} _{\fX ^{\sharp (s)} /\fS ^\sharp} \otimes _{\widetilde{\cO} _{\fX ^{(s)}}}  \cE).
\end{equation}

\item We have the canonical isomorphism of  right  $\widetilde{\cD} ^{(m)\dag} _{\fX ^{\sharp}/\fS ^{\sharp}}$-modules
\begin{equation} \label{prop-omegaFflatF*2}
\nu _{\M} \colon  F ^* ( \M \otimes _{\widetilde{\cO} _{\fX ^{(s)}}} \widetilde{\omega} ^{-1} _{\fX ^{\sharp (s)} /\fS ^\sharp}) \riso F ^\flat (\cM) \otimes _{\widetilde{\cO} _{\fX}} \widetilde{\omega} ^{-1} _{\fX ^\sharp /\fS ^\sharp}.
\end{equation}
\end{enumerate}
\end{prop}

\begin{proof}In the same way as \cite[4.3.4.1]{Be2}, we construct the isomorphism \ref{prop-omegaFflatF*1} as follows:
\begin{gather}\notag   \widetilde{\omega} _{\fX ^\sharp /\fS ^\sharp} \otimes _{\widetilde{\cO} _{\fX}} F ^* (\cE)   \riso \widetilde{\omega} _{\fX ^\sharp /\fS ^\sharp} \otimes _{\widetilde{\cO} _{\fX ^{(s)}}} \cE
\underset{\ref{mu-Fflat-omega}}{\riso} \mathcal{H} om _{\cO _{\fX ^{ (s)}}} (\cO _{\fX}, \widetilde{\omega} _{\fX ^{\sharp (s)} /\fS ^\sharp}) \otimes _{\widetilde{\cO} _{\fX ^{(s)}}} \cE
\\  \label{prop-omegaFflatF*proof1} \riso \mathcal{H} om _{\cO _{\fX ^{ (s)}}} (\cO _{\fX}, \widetilde{\omega} _{\fX ^{\sharp (s)} /\fS ^\sharp} \otimes _{\widetilde{\cO} _{\fX ^{(s)}}} \cE )
=   F ^\flat ( \widetilde{\omega} _{\fX ^{\sharp (s)} /\fS ^\sharp} \otimes _{\widetilde{\cO} _{\fX ^{(s)}}}  \cE).
\end{gather}
It remains to check it is $\widetilde{\cD} ^{(m)\dag} _{\fX ^{\sharp}/\fS ^{\sharp}}$-linear.
Let $F _{\eta}\colon \fX _\eta \to \fX ^{(s)} _{\eta}$ be the map induced by $F$. 
Set $\cM \coloneqq \widetilde{\omega} _{\fX ^{\sharp} /\fS ^\sharp} \otimes _{\widetilde{\cO} _{\fX }}  \cE$,
$\cE _{\eta}:=  \widehat{\cD} ^{(m)}  _{\fX ^{(s)} _\eta/\eta} \otimes _{j ^{-1}\widetilde{\cD} ^{(m) \dag}  _{\fX^{\sharp(s)}/\fS^\sharp}} j ^{-1} \cE$,
$\cM _\eta\coloneqq  \cM \otimes _{j ^{-1}\widetilde{\cD} ^{(m) \dag}  _{\fX^{\sharp(s)} /\fS^\sharp}} j ^{-1}  \widehat{\cD} ^{(m)}  _{\fX ^{(s)}_\eta/\eta}  $, 
 $(F ^{*}\cE) _{\eta}:=  \widehat{\cD} ^{(m)}  _{\fX_\eta/\eta} \otimes _{j ^{-1}\widetilde{\cD} ^{(m) \dag}  _{\fX^\sharp/\fS^\sharp}} j ^{-1} F ^{*}\cE $,
$(F ^{\flat}\cM) _\eta\coloneqq  \cM \otimes _{j ^{-1}\widetilde{\cD} ^{(m) \dag}  _{\fX^\sharp/\fS^\sharp}} j ^{-1}  \widehat{\cD} ^{(m)}  _{\fX_\eta/\eta}  $, we get the isomorphisms:
$(F ^{*}\cE) _{\eta} \riso F _{\eta}^{*}(\cE _\eta)$ and  $(F ^{\flat}\cM) _\eta \riso (F _\eta ^{\flat}\cM _\eta) $. We obtain the canonical isomorphism:
\begin{equation}\label{prop-omegaFflatF*proof2}
\xymatrix{   { \widetilde{\omega} _{\fX ^\sharp /\fS ^\sharp} \otimes _{\widetilde{\cO} _{\fX}} F ^* (\cE)   } \ar[r] ^-{\ref{prop-omegaFflatF*proof1}} \ar@{^{(}->}[d] ^-{\ref{twist-Ddagmproof1}}
& { F ^\flat ( \widetilde{\omega} _{\fX ^{\sharp (s)} /\fS ^\sharp} \otimes _{\widetilde{\cO} _{\fX ^{(s)}}}  \cE) }  \ar@{^{(}->}[d] ^-{\ref{twist-Ddagmproof1}}
\\  {j _* (\omega _{\fX _{\eta} /\eta} )   \otimes _{\cO _{\fX _\eta}} F _{\eta}^{*} \cE _{\eta}   } \ar[r] ^-{\sim}&  { F  _{\eta} ^\flat (\omega _{\fX  _{\eta} ^{(s)} /\eta} \otimes _{\cO _{\fX  _{\eta} ^{(s)}}}  \cE  _{\eta}) ,}     }
\end{equation}
where the bottom isomorphism comes from by $p$-adic completion the standard isomorphisms  of Proposition \cite[2.4.3]{Be2}. This shows the $\widetilde{\cD} ^{(m)\dag} _{\fX ^{\sharp}/\fS ^{\sharp}}$-linearity of the top isomorphism.
Using \ref{left2right}, we get the isomorphism \ref{prop-omegaFflatF*2} from \ref{prop-omegaFflatF*1}. 
\end{proof}

\begin{coro} \label{cor-F+(m)}       We have  the isomorphism of $(\widetilde{\cD} ^{(m)\dag} _{\fX ^{\sharp}/\fS ^{\sharp}}, \widetilde{\cD} ^{(m)\dag} _{\fX ^{\sharp (s)}/\fS ^{\sharp}})$-bimodules: 
\begin{equation} \label{cor-F+(m)iso} \widetilde{\cD} ^{(m)\dag} _{\fX ^{\sharp (s)} \leftarrow \fX ^{\sharp}/\fS ^{\sharp}} := F _{\mathrm{l}} ^* (\widetilde{\cD} ^{(m)\dag} _{\fX ^{\sharp (s)}/\fS ^{\sharp}}
\otimes _{\widetilde{\cO} _{\fX ^{(s)}}} \widetilde{\omega} ^{-1} _{\fX ^{\sharp (s)} /\fS ^\sharp}) \otimes _{\widetilde{\cO} _{\fX}} \widetilde{\omega} _{\fX ^\sharp /\fS ^\sharp} \riso  F ^\flat (\widetilde{\cD} ^{(m)\dag} _{\fX ^{\sharp (s)}/\fS ^{\sharp}}).
\end{equation}
\end{coro}

\begin{proof}By functoriality, this is a consequence of \ref{prop-omegaFflatF*}.
\end{proof}

\begin{empt}  \label{infinitelevel}
Taking the colimits on the level,  we get from \ref{phidagdescrip} the isomorphism of  
$\widetilde{\cD} ^{\dag} _{\fX ^{\sharp}/\fS ^{\sharp}}$-bimodules 
$$F ^* F ^\flat  \widetilde{\cD} ^{\dag} _{\fX ^{\sharp (s)}/\fS ^{\sharp}} \riso \widetilde{\cD} ^{\dag} _{\fX ^{\sharp}/\fS ^{\sharp}}.$$
Replacing respectively  $\widetilde{\cD} ^{(m+s)\dag} _{\fX ^{\sharp}/\fS ^{\sharp}}$ by  $\widetilde{\cD} ^{\dag} _{\fX ^{\sharp}/\fS ^{\sharp}}$,
and  $\widetilde{\cD} ^{(m)\dag} _{\fX ^{\sharp (s)}/\fS ^{\sharp}}$ by  $\widetilde{\cD} ^{\dag} _{\fX ^{\sharp (s)}/\fS ^{\sharp}}$, we analogous statements of  
\ref{lem-factorization2},
\ref{FflatDlocproj},
\ref{dfn-FflatF*},
\ref{lem254Be2},
\ref{eqcatFrobdescexpl},
\ref{eqcatFrobdescexpl2},
\ref{prop-omegaFflatF*},
\ref{cor-F+(m)}.
Later, we can freely use these references for the cases without ``$(m)$''. 
\end{empt}

\begin{ntn}We have the $(\widetilde{\cD} ^{(m)\dag} _{\fX ^{\sharp (s)}/\fS ^{\sharp}}, \widetilde{\cD} ^{(m)\dag} _{\fX ^{\sharp}/\fS ^{\sharp}})$-bimodule and the $(\widetilde{\cD} ^{\dag} _{\fX ^{\sharp (s)}/\fS ^{\sharp}}, \widetilde{\cD} ^{\dag} _{\fX ^{\sharp}/\fS ^{\sharp}})$-bimodule:
$$\widetilde{\cD} ^{(m)\dag} _{\fX ^{\sharp} \to  \fX ^{\sharp (s)} /\fS ^{\sharp}}\coloneqq F ^* \widetilde{\cD} ^{(m)\dag} _{\fX ^{\sharp (s)} /\fS ^{\sharp}} ,
\qquad  \widetilde{\cD} ^{\dag} _{\fX ^{\sharp} \to  \fX ^{\sharp (s)} /\fS ^{\sharp}}\coloneqq F ^* \widetilde{\cD} ^{\dag} _{\fX ^{\sharp (s)} /\fS ^{\sharp}}.$$
\end{ntn}

\begin{ntn}
Let $\cE$ be a left
$\widetilde{\cD} ^{\dag} _{\fX ^{\sharp}/\fS ^{\sharp}}$-module. 
We set 
$$F _{+} (\cE):= 
F _* ( 
\widetilde{\cD} ^{\dag} _{\fX ^{\sharp (s)} \leftarrow \fX ^{\sharp}/\fS ^{\sharp}}
\otimes ^\L _{\widetilde{\cD} ^{\dag} _{\fX ^{\sharp}/\fS ^{\sharp}}}
\cE)
\underset{\ref{cor-F+(m)iso}}{\riso} 
F _* (
F ^\flat (\widetilde{\cD} ^{\dag} _{\fX ^{\sharp (s)}/\fS ^{\sharp}})
\otimes ^{\bbL}_{\widetilde{\cD} ^{\dag} _{\fX ^{\sharp}/\fS ^{\sharp}}}
\cE)
\underset{\ref{FflatDlocproj}}{\riso} 
F _* (
F ^\flat (\widetilde{\cD} ^{\dag} _{\fX ^{\sharp (s)}/\fS ^{\sharp}})
\otimes _{\widetilde{\cD} ^{\dag} _{\fX ^{\sharp}/\fS ^{\sharp}}}
\cE).$$
Let $\cM$ be a right
$\widetilde{\cD} ^{\dag} _{\fX ^{\sharp}/\fS ^{\sharp}}$-module. 
We set 
$$F _{+} (\cM)\coloneqq
F _* ( 
\cM 
\otimes ^\L _{\widetilde{\cD} ^{\dag} _{\fX ^{\sharp}/\fS ^{\sharp}}}
\widetilde{\cD} ^{\dag} _{\fX ^{\sharp} \to  \fX ^{\sharp (s)} /\fS ^{\sharp}}
)
=
F _* ( 
\cM 
\otimes ^\L _{\widetilde{\cD} ^{\dag} _{\fX ^{\sharp}/\fS ^{\sharp}}}
F ^* \widetilde{\cD} ^{\dag} _{\fX ^{\sharp (s)} /\fS ^{\sharp}}
)
\underset{\ref{lem-factorization2}}{\riso}  
F _* ( 
\cM 
\otimes  _{\widetilde{\cD} ^{\dag} _{\fX ^{\sharp}/\fS ^{\sharp}}}
F ^* \widetilde{\cD} ^{\dag} _{\fX ^{\sharp (s)} /\fS ^{\sharp}}
).$$
\end{ntn}

\begin{prop}
\label{F+Fflatiso}
Let $\cM '$ be a right  $\widetilde{\cD} ^{\dag} _{\fX ^{\sharp (s)}/\fS ^{\sharp}}$-module. Let $\cE '$ be a left $\widetilde{\cD} ^{\dag} _{\fX ^{\sharp (s)}/\fS ^{\sharp}}$-module. 
\begin{enumerate}[(a)]
\item We have the $\widetilde{\cD} ^{\dag} _{\fX ^{\sharp (s)}/\fS ^{\sharp}}$-linear isomorphism
\begin{equation}\label{F+Fflatiso-left}
F _+ \circ F ^* (\E') \riso  \cE '.
\end{equation}
\item We have the $\widetilde{\cD} ^{\dag} _{\fX ^{\sharp (s)}/\fS ^{\sharp}}$-linear isomorphism
\begin{equation}
\label{F+Fflatiso-right}
F _+\circ  F ^\flat \M ' \riso   \M'.
\end{equation}
\end{enumerate}
\end{prop}

\begin{proof}
This is a consequence of respectively \ref{eqcatFrobdescexpl-iso1} and \ref{eqcatFrobdescexpl-iso2}.
\end{proof}

\begin{prop}
\label{F+Fflatisobis}
Let $\cM $ be a right  $\widetilde{\cD} ^{\dag} _{\fX ^{\sharp}/\fS ^{\sharp}}$-module. Let $\cE $ be a left $\widetilde{\cD} ^{\dag} _{\fX ^{\sharp}/\fS ^{\sharp}}$-module. 
\begin{enumerate}[(a)]
\item The adjoint canonical $\widetilde{\cD} ^{\dag} _{\fX ^{\sharp}/\fS ^{\sharp}}$-linear isomorphism
\begin{equation}\label{F+Fflatisobis-left}
\cE \riso F ^*\circ F _+   (\E).
\end{equation}
is an isomorphism. 
\item The adjoint canonical $\widetilde{\cD} ^{\dag} _{\fX ^{\sharp}/\fS ^{\sharp}}$-linear isomorphism
\begin{equation} \label{F+Fflatisobis-right}
 \M \riso F ^\flat \circ F _+ \M .
\end{equation}
is an isomorphism.
\end{enumerate}
\end{prop}

\begin{proof}
This is a consequence of \ref{phidagdescrip}.
\end{proof}

\begin{coro} \label{eqcatFrobdescexplF+} 
The functors $\cE \mapsto F ^* \cE$  (resp.  $\cM \mapsto F ^\flat \cM$) and $\cF \mapsto F _+ \cF$ induce exact  equivalences from the category of  left (resp. right) $\widetilde{\cD} ^{\dag} _{\fX ^{\sharp (s)}/\fS ^{\sharp}}$-modules and that of  left (resp. right) 
$\widetilde{\cD} ^{\dag} _{\fX ^{\sharp}/\fS ^{\sharp}}$-modules. 
\end{coro}

\begin{proof}
This is a consequence of \ref{F+Fflatiso} and of \ref{F+Fflatisobis}.
\end{proof}

\begin{empt}
\label{Tr+Fiso}
We have the commutative diagram
\begin{equation}
\xymatrix{
{F _* \widetilde{\omega} _{\fX ^{\sharp}/\fS ^{\sharp}}} 
\ar[r] ^-{\Tr _F}
\ar[d] ^-{T _F}
& 
{\widetilde{\omega} _{\fX ^{\sharp(s)}/\fS ^{\sharp}}} 
\ar@{=}[d] ^-{}
\\ 
{F _+ \widetilde{\omega} _{\fX ^{\sharp}/\fS ^{\sharp}}} 
\ar[r] ^-{\Tr _{+,F}} _-{\sim}
& 
{\widetilde{\omega} _{\fX ^{\sharp(s)}/\fS ^{\sharp}},} 
}
\end{equation}
where $\Tr _F$ is the trace morphism of \ref{f*fflatadj}, 
$T _F
\colon F _* \widetilde{\omega} _{\fX ^{\sharp}/\fS ^{\sharp}}
\to
F _* ( 
\widetilde{\omega} _{\fX ^{\sharp}/\fS ^{\sharp}}
\otimes ^\L _{\widetilde{\cD} ^{\dag} _{\fX ^{\sharp}/\fS ^{\sharp}}}
\widetilde{\cD} ^{\dag} _{\fX ^{\sharp}
\to 
\fX ^{\sharp (s)} /\fS ^{\sharp}}
)
= F _+ \widetilde{\omega} _{\fX ^{\sharp}/\fS ^{\sharp}}$
is the canonical morphism induced by 
$\widetilde{\cD} ^{\dag} _{\fX ^{\sharp}/\fS ^{\sharp}}
\to 
\widetilde{\cD} ^{\dag} _{\fX ^{\sharp}
\to 
\fX ^{\sharp (s)} /\fS ^{\sharp}}
$,
and
$\Tr _{+,F}$ is the composition
$F _+ \widetilde{\omega} _{\fX ^{\sharp}/\fS ^{\sharp}} \underset{\ref{mu-Fflat-omega}}{\riso} F _+ F ^\flat \widetilde{\omega} _{\fX ^{\sharp(s)}/\fS ^{\sharp}} \underset{\ref{F+Fflatiso-right}}{\riso}  \widetilde{\omega} _{\fX ^{\sharp(s)}/\fS ^{\sharp}}$.
Since the restriction of $T _F $ on $\fX ^{(s)} _\eta$ is surjective, 
since the diagram remains commutative on $\fX ^{(s)} _\eta$,
since Virrion's trace map 
$\Tr _{+,F _\eta} \colon  F _{\eta +} \omega _{\fX _\eta/\eta}  \to  \omega _{\fX ^{(s)} _\eta/\eta}$
satisfies the same commutative diagram, then 
\begin{equation}
\label{Tr+Fiso-Treta}
\Tr _{+,F _\eta} =\Tr _{+,F } | \fX ^{(s)} _\eta .
\end{equation}

\end{empt}

\subsubsection{Universal homeomorphism}

\begin{lem} \label{sumbounded}
Let $\fX$, $\fY$ be two objects of $\Sm ^\dag _{\eta}$ which are affine $\fS$-formal schemes (see Definition \ref{dfn-Sm-eta}). Let $f\colon \fY \to \fX$ be a finite morphism.
Put $A '= \Gamma (\fX, \cO _{\fX})$,  $A = \Gamma (\fY, \cO _{\fY})$. Let $e _1\dots, e _r \in A$ wich generate $A$ as $A'$-module.  Let $b \in A _{[t]}$. Let $c>0$, $d$ be two real numbers. 

\begin{enumerate}[(a)]
\item If $b$ is $(c,d)$-bounded as $A$-algebra  with respect to $1/t$, then there exist $ b' _1, \dots, b ' _r \in A ' _{[t]}$ which are $(c,d)$-bounded as $A'$-algebra  with respect to $1/t$ and satisfy the equality $b = \sum _{j=1} ^{r} b ' _j e _j$.

\item \label{sumbounded(b)}
If we can write $b$ in the form $b=  \sum _{n\geq 0} \pi ^n b _n $,  where $b _n \in A _{[t]}$ is $(c, cn +d)$-bounded as $A$-algebra with respect to $1/t$, then  $b$ is $(c,d)$-bounded as $A$-algebra with respect to $1/t$. 
\end{enumerate}
\end{lem}

\begin{proof}
We can write  $b = \sum _{n\geq 0} \pi ^n \frac{a _n}{ t ^{\alpha _n}}$,  with $a _n \in A$, $\alpha _n\leq cn+d$ is an integer. There exist  $ a' _{n1}, \dots, a ' _{nr} \in A ' $ such that $a _n = \sum _{j=1} ^{r} a ' _{nj} e _j$.
Then $b ' _j =  \sum _{n\geq 0} \pi ^n \frac{a '_{nj}}{ t ^{\alpha _n}}$ satisfy the desired properties. The second statement is also easy to check. 
\end{proof}

\begin{prop} \label{D2Dmdag}
Let $f \colon \fY ^{\sharp}\to \fX ^{\sharp}$ be a finite surjective morphism of  $\Sm ^\dag _{\fS ^\sharp}$  (see Definition \ref{ntn-Smdagdag}).
Let   $\widetilde{f}\colon  (\fY, \widetilde{\cO} _{\fY})  \to  (\fX, \widetilde{\cO} _{\fX})$ be the flat induced  morphism of locally ringed spaces  (see \ref{dfn-f-tilde} and \ref{D2Dmdagpre}).
The sheaf $\widetilde{f} ^{ *} \widetilde{\cD} ^{(m)\dag} _{\fX ^{\sharp}/\fS ^{\sharp}}$  is a $(\widetilde{\cD} ^{(m)\dag} _{\fY ^{\sharp}/\fS ^{\sharp}},  f ^{-1} \widetilde{\cD} ^{(m)\dag} _{\fX ^{\sharp}/\fS ^{\sharp}})$
sub-bimodule of $j _*f ^{ *} _\eta \widehat{\cD} ^{(m)} _{\fX _\eta/\eta} $. In particular, we get the canonical  morphism
$\widetilde{\cD} ^{(m)\dag} _{\fY ^{\sharp}/\fS ^{\sharp}} \to  \widetilde{f} ^* \widetilde{\cD} ^{(m)\dag} _{\fX ^{\sharp}/\fS ^{\sharp}}$ making commutative the diagram
\begin{equation} \label{D2Dmdag-diag}
\xymatrix{  {\widetilde{\cD} ^{(m)} _{\fY ^{\sharp}/\fS ^{\sharp}}}  \ar@{^{(}->}[r] ^-{} \ar[d] ^-{\theta} &  {\widetilde{\cD} ^{(m)\dag} _{\fY ^{\sharp}/\fS ^{\sharp}}}  \ar@{^{(}->}[r] ^-{}  \ar@{.>}[d] ^-{\theta ^\dag}
&   {j _*\widehat{\cD} ^{(m)} _{\fY _\eta/\eta}} \ar[d] ^-{\widehat{\theta}}
\\   {\widetilde{f} ^{ *} \widetilde{\cD} ^{(m)} _{\fX ^{\sharp}/\fS ^{\sharp}}}  \ar@{^{(}->}[r] ^-{} & {\widetilde{f} ^{ *} \widetilde{\cD} ^{(m)\dag} _{\fX ^{\sharp}/\fS ^{\sharp}}}  \ar@{^{(}->}[r] ^-{} &  {j _*f ^{ *} _\eta \widehat{\cD} ^{(m)} _{\fX _\eta/\eta} }  }
\end{equation}
where $\theta \colon \widetilde{\cD} ^{(m)} _{\fY ^{\sharp}/\fS ^{\sharp}} \to  \widetilde{f} ^* \widetilde{\cD} ^{(m)} _{\fX ^{\sharp}/\fS ^{\sharp}}$
and $\widehat{\theta}  \colon \widehat{\cD} ^{(m)} _{\fY _\eta/\eta} \to  f ^* _\eta \widehat{\cD} ^{(m)} _{\fX _\eta/\eta}$ are the canonical morphisms.
\end{prop}

\begin{proof}
Since $\widetilde{f} ^{ *} \widetilde{\cD} ^{(m)\dag} _{\fX ^{\sharp}/\fS ^{\sharp}}$ 
is a right $f ^{-1} \widetilde{\cD} ^{(m)\dag} _{\fX ^{\sharp}/\fS ^{\sharp}}$
sub-module
of
$j _*f ^{ *} _\eta \widehat{\cD} ^{(m)} _{\fX _\eta/\eta} $, 
and since for any 
$P \in 
\Gamma (\fY, \widetilde{\cD} ^{(m)\dag} _{\fY ^{\sharp}/\fS ^{\sharp}})$
we have
$\widehat{\theta} (P)
=
P \cdot (1 \otimes 1)$,
then 
we reduce to check 
the commutativity of \ref{D2Dmdag-diag}.

Since this is local, we can suppose $\fX= \Spf A'$ has logarithmic coordinates 
$u ' _1, \dots, u ' _d \in \Gamma (\fX, M _{\fX ^\sharp})$, and 
$\fY= \Spf A$ has logarithmic coordinates 
$u  _1, \dots, u  _d \in \Gamma (\fY, M _{\fY ^\sharp})$.
We denote by 
$t ' _1, \dots, t ' _d $
the image of 
$u ' _1, \dots, u ' _d $
in $ \Gamma (\fX, \cO _{\fX ^\sharp})$,
and by 
$t  _1, \dots, t  _d $
the image of 
$u  _1, \dots, u  _d $
in $ \Gamma (\fY, \cO _{\fY ^\sharp})$.
We have the elements
$\underline{\tau} ^{\{ \underline{k}\} _{(m)}}
\in 
\mathcal{P} ^n _{\fY/\fS, (m)}$,
and
$\underline{\tau} ^{\prime \{ \underline{i}\} _{(m)}}
\in 
\mathcal{P} ^n _{\fX/\fS, (m)}$.
Consider the commutative diagram
\begin{equation} \notag
\xymatrix{ {\mathcal{P} ^n _{\fY/\fS, (m)}}  \ar[r] ^-{} \ar[d] ^-{} &  {\mathcal{P} ^n _{\fY ^\sharp/\fS ^\sharp, (m)}}  \ar@{^{(}->}[d] ^-{}
\\  {\widetilde{\mathcal{P}} ^n _{\fY/\fS, (m)}}  \ar@{=}[r] ^-{} &  {\widetilde{\mathcal{P}} ^n _{\fY ^\sharp/\fS ^\sharp, (m)}.}     }
\end{equation}
Since we are only considering the elements in
$\widetilde{\mathcal{P}} ^n _{\fY /\fS , (m)}$,
we will still denote by 
$\underline{\tau} ^{\{ \underline{k}\} _{(m)}}$
its image in $\mathcal{P} ^n _{\fY ^\sharp/\fS ^\sharp, (m)}$ (a priori, the top arrow is not injective).
We use the similar abuse of notation with
$\underline{\tau} ^{\prime \{ \underline{i}\} _{(m)}}$.

Using the universal property of the $m$-PD-envelope,  we get the canonical morphism  $\phi _f \colon f ^* \mathcal{P} ^n _{\fX ^{\sharp}/\fS ^{\sharp}, (m)} \to \mathcal{P} ^n _{\fY ^{\sharp}/\fS ^{\sharp} , (m)}$.
This yields by extension the morphism $\widetilde{\phi} _f\colon \widetilde{f} ^* \widetilde{\mathcal{P}} ^n _{\fX ^{\sharp}/\fS ^{\sharp}, (m)} \to  \widetilde{\mathcal{P}} ^n _{\fY ^{\sharp}/\fS ^{\sharp}, (m)}$.
By applying  $\mathcal{H}om _{\widetilde{\cO} _\fY} ( - ,\widetilde{\cO} _\fY)$ to $ d _{0*}\widetilde{\phi} _f$,  where ``$d _{0*}$'' means that we consider the left structure of $\widetilde{\cO} _\fY$-modules (which comes from the left projection) we get
$\theta _n \colon \widetilde{\cD} ^{(m)} _{\fY ^{\sharp}/\fS ^{\sharp},n} \to  \widetilde{f} ^* \widetilde{\cD} ^{(m)} _{\fX ^{\sharp}/\fS ^{\sharp},n}$. Since there are compatible,  taking the inductive limits in the order $n$ this yields
$\theta  \colon \widetilde{\cD} ^{(m)} _{\fY ^{\sharp}/\fS ^{\sharp}} \to  \widetilde{f} ^* \widetilde{\cD} ^{(m)} _{\fX ^{\sharp}/\fS ^{\sharp}}$.
The canonical map  $\theta  \colon \cD ^{(m)} _{\fY _\eta/\eta} \to  f _\eta ^* \widetilde{\cD} ^{(m)} _{\fX _\eta/\eta}$ is constructed similarly. By completion, this yields
$\widehat{\theta} \colon  \widehat{\cD} ^{(m)} _{\fY _\eta/\eta} \to   f ^{ *} _\eta \widehat{\cD} ^{(m)} _{\fX _\eta/\eta} $. The outline of the diagram  \ref{D2Dmdag-diag} is commutative.  It remains to check that 
$\widehat{\theta} (\widetilde{\cD} ^{(m)\dag} _{\fY ^{\sharp}/\fS ^{\sharp}}) \subset \widetilde{f} ^* \widetilde{\cD} ^{(m)\dag} _{\fX ^{\sharp}/\fS ^{\sharp}}$.

For any $\underline{i}$ and $\underline{k}$ such that 
$|\underline{i}| \leq |\underline{k}|$,
we compute
$$a _{\underline{k},\underline{i}}  : = \theta _{|\underline{k}| } ( \underline{\partial} ^{<\underline{k}> _{(m)}}) (1 \otimes \underline{\tau} ^{\prime\{ \underline{i}\} _{(m)}})
=  \underline{\partial} ^{<\underline{k}> _{(m)}}  (\widetilde{\phi} _f  (1 \otimes \underline{\tau} ^{\prime\{ \underline{i}\} _{(m)}}))
= \underline{\partial} ^{<\underline{k}> _{(m)}}  (1 \otimes \phi _f  (\underline{\tau} ^{\prime\{ \underline{i}\} _{(m)}})),$$
where 
$\phi _f (\underline{\tau} ^{\prime\{ \underline{i}\} _{(m)}})$
means its image in 
$\mathcal{P} ^{|\underline{k}| } _{\fY ^{\sharp} /\fS ^{\sharp}, (m)}$.
We get
$\underline{t} ^{\underline{k}} \underline{\partial} ^{<\underline{k}> _{(m)}} 
(1 \otimes \phi _f 
(\underline{\tau} ^{\prime\{ \underline{i}\} _{(m)}}))
=
 \underline{\partial} ^{<\underline{k}> _{ (m)}} _{\sharp} 
(1 \otimes \phi _f 
(\underline{\tau} ^{\prime\{ \underline{i}\} _{(m)}}))
=
 \underline{\partial} ^{<\underline{k}> _{ (m)}} _{\sharp} 
( \phi _f (\underline{\tau} ^{\prime\{ \underline{i}\} _{(m)}}))
\in 
A
$.
Hence,  $a _{\underline{k},\underline{i}} \in \frac{1}{\underline{t} ^{\underline{k}}} \cO _{\fY} \subset  \widetilde{\cO} _{\fY}$ (remark that since  $\underline{t} ^{\underline{k}} \in A \{ \frac{1}{t}\} ^*= \Gamma (\fX _\eta, \cO ^* _{\fX _\eta})$
then following \ref{Otilde*} we get $\underline{t} ^{\underline{k}} \in A _{[t]} ^*= \Gamma (\fX _\eta, \widetilde{\cO} ^* _{\fX})$). Since  the elements  $1 \otimes \underline{\tau} ^{\prime\{ \underline{i}\} _{(m)}}$ such that 
$|\underline{i}| \leq |\underline{k}| $ form a basis of the free $\widetilde{\cO} _{\fY}$-module $\widetilde{f} ^* \widetilde{\mathcal{P}} ^{|\underline{k}|} _{\fX /\fS , (m)}$,  we get the equality 
\begin{equation} \label{partial-k-=sumpartial'}
\theta ( \underline{\partial} ^{<\underline{k}> _{(m)}}) = \theta _{|\underline{k}|} ( \underline{\partial} ^{<\underline{k}> _{(m)}}) =
\frac{1}{\underline{t} ^{\underline{k}}} \sum _{|\underline{i}|\leq |\underline{k}|} a _{\underline{k},\underline{i}} \otimes \underline{\partial} ^{\prime<\underline{i}> _{(m)}}
\end{equation}
with $a _{\underline{k},\underline{i}}   \in A$.

Let  $P \in  \Gamma (\fY, \widetilde{\cD} ^{(m)\dag} _{\fY ^{\sharp}/\fS ^{\sharp}})$. Following \ref{wkcpdscplevelm}, we can write 
$P =\sum _{\underline{k}\in \bbN ^d} \sum _{l\in \bbN}   \frac{a _{\underline{k},l}}{t ^{l}}\underline{\partial} ^{<\underline{k}> _{(m)}}$,
where $a _{\underline{k},l} \in A$ are such that there exists a constant $c>0$ satisfying $v _{p} ( a _{\underline{k},l}) \geq \frac{l+|\underline{k}|}{c} -1$,
for any $\underline{k} \in \bbN ^d$ and $l \in \bbN$. For any $\underline{i} \in \bbN ^d$ and $j \in \bbN$,  we set 
$$E _{\underline{i},j} = \{ (\underline{k}, l),  \textrm{ with $l\in \bbN$, $\underline{k}\in \bbN ^d$ are such that $|\underline{i}| \leq |\underline{k} |$ and  $v _{\pi} (a _{\underline{k},l} a _{\underline{k},\underline{i}})=j$ }  \}.$$
For any  $(\underline{k}, l) \in E _{\underline{i},j}$,  we have  $c (j+1) \geq c (1 +v _{\pi} (a _{\underline{k},l})) \geq l + |\underline{k} |$.  Hence, the set $E _{\underline{i},j}$ is finite.
For any  $(\underline{k}, l) \in E _{\underline{i},j}$,  we put  $a _{\underline{i}, j, \underline{k},l} :=\pi ^{-j} a _{\underline{k},l} a _{\underline{k},\underline{i}}\in A$.
Put  $ b _{\underline{i},j} :=  \sum _{(\underline{k}, l) \in E _{\underline{i},j}} \frac{1}{t ^{l} \underline{t} ^{\underline{k}}} a _{\underline{i}, j, \underline{k},l}$,
and  $ b _{\underline{i}} :=  \sum _{j\geq 0} \pi ^j b _{\underline{i},j}$. Using   the formula  \ref{partial-k-=sumpartial'}, since $\widehat{\theta}$ is $A _{[t]}$-linear,  we get 
$$ \widehat{\theta} (P) = \sum _{\underline{k}\in \bbN ^d} \sum _{l\in \bbN}  \sum _{|\underline{i}|\leq |\underline{k}|}
\frac{a _{\underline{k},l} a _{\underline{k},\underline{i}}}{t ^{l} \underline{t} ^{\underline{k}}} \otimes \underline{\partial} ^{\prime<\underline{i}> _{(m)}}
=  \sum _{\underline{i} \in \bbN ^d}\sum _{j\geq 0} \sum _{(\underline{k}, l) \in E _{\underline{i},j}} \pi ^j \frac{a _{\underline{i}, j, \underline{k},l}}{t ^{l} \underline{t} ^{\underline{k}}}
\otimes \underline{\partial} ^{\prime<\underline{i}> _{(m)}} =  \sum _{\underline{i} \in \bbN ^d} b _{\underline{i}} \otimes \underline{\partial} ^{\prime<\underline{i}> _{(m)}}.  $$
Let $c \geq 1$ be a real number such that $1/t _1,\dots, 1/t _d \in A _{[t]} ^*$ are  $(c,c)$-bounded as $A$-algebra with respect to $1/t$. Let  $(\underline{k}, l) \in E _{\underline{i},j}$. 
We compute that  $\frac{1}{t ^{l} \underline{t} ^{\underline{k}}}$ is $(c,l+c |\underline{k}|)$-bounded as $A$-algebra with respect to $1/t$ (see the estimates of \ref{(c,d)bounded-ring}).
Since $c (j+1) \geq l + |\underline{k}|$,   and $|\underline{k}| \geq |\underline{i}|$,  then we get $l+c|\underline{k}| =  (1+c)|\underline{k}| + l -|\underline{k}| 
\leq  (1+c)|\underline{k}| + l -|\underline{i}| \leq  (1+c)(|\underline{k}| + l) -|\underline{i}|  \leq  c(1+c)(j+1) -|\underline{i}|  = c(1+c)j + c (1+c) -|\underline{i}| \leq  C j + C -|\underline{i}| $, 
where $C := c (1+c)$. Hence, $ b _{\underline{i},j}$ is  $(C,C j + C -|\underline{i}|)$-bounded as $A$-algebra with respect to $1/t$. Hence, using \ref{sumbounded}.(\ref{sumbounded(b)}), 
$ b _{\underline{i}}$ is  $(C,C -|\underline{i}|)$-bounded as $A$-algebra with respect to $1/t$. Using \ref{sumbounded}.1, there exist
$ b' _{\underline{i},1}, \dots, b ' _{\underline{i},r} \in A ' _{[t]}$ which are $(C,C -|\underline{i}|)$-bounded as $A'$-algebra  with respect to $1/t$
and satisfy the equality $b _{\underline{i}}= \sum _{j=1} ^{r} b ' _{\underline{i},j} e _j$. Hence, $\widehat{\theta} (P)= \sum _{j=1} ^{r}  e _j \otimes   \sum _{\underline{i} \in \bbN ^d} b ' _{\underline{i},j} \underline{\partial} ^{\prime<\underline{i}> _{(m)}}$,
with  $\sum _{\underline{i} \in \bbN ^d} b ' _{\underline{i},j} \underline{\partial} ^{\prime<\underline{i}> _{(m)}} \in  \Gamma (\fX,  \widetilde{\cD} ^{(m)\dag} _{\fX ^{\sharp}/\fS ^{\sharp}})$.
\end{proof}

\begin{ntn}
Let $f \colon \fY ^{\sharp}\to \fX ^{\sharp}$ be a finite surjective morphism of  $\Sm ^\dag _{\fS ^\sharp}$  (see Definition \ref{ntn-Smdagdag}).
Let $\cE$ be a left $\widetilde{\cD} ^{\dag} _{\fX ^{\sharp }/\fS ^{\sharp},\bbQ}$-module.  We get the left $\widetilde{\cD} ^{\dag} _{\fY ^{\sharp}/\fS ^{\sharp},\bbQ}$-module
$\widetilde{f} ^*  \cE :=  f ^*  \widetilde{\cD} ^{\dag} _{\fX ^{\sharp }/\fS ^{\sharp},\bbQ}  \otimes _{f ^{-1} \widetilde{\cD} ^{\dag} _{\fX ^{\sharp}/\fS ^{\sharp},\bbQ}}  f ^{-1}\cE $. We remark that the canonical morphism
\begin{equation} \notag
 \widetilde{\cO} _{\fY,\bbQ}  \otimes _{f ^{-1} \widetilde{\cO} _{\fX ,\bbQ}}  f ^{-1}\cE 
\to \widetilde{f} ^*  \widetilde{\cD} ^{\dag} _{\fX ^{\sharp }/\fS ^{\sharp},\bbQ}  \otimes _{f ^{-1} \widetilde{\cD} ^{\dag} _{\fX ^{\sharp}/\fS ^{\sharp},\bbQ}}  f ^{-1}\cE 
\end{equation}
is an isomorphism, which justifies the notation $\widetilde{f} ^*  \cE$. Moreover, the functor  $\widetilde{f} ^*$ from the category of 
left $\widetilde{\cD} ^{\dag} _{\fX ^{\sharp }/\fS ^{\sharp},\bbQ}$-modules to that of  left $\widetilde{\cD} ^{\dag} _{\fY ^{\sharp }/\fS ^{\sharp},\bbQ}$-modules is transitive with respect to the composition, 
i.e. for any object $\fZ ^{\sharp}$ of $\Sm ^\dag _{\fS ^\sharp}$, and finite, surjective morphism $g \colon \fZ \to \fY$, we have the canonical isomorphism $\widetilde{g} ^*  (\widetilde{f} ^*  (\cE ) ) \riso  \widetilde{f\circ g} ^* (\cE)$.
\end{ntn}

\begin{empt} \label{tau-dagm}
Let $\fX ^{\sharp}$, $\fY ^{\sharp}$ be two objects of $\Sm ^\dag _{\fS ^\sharp}$ (see Definition \ref{ntn-Smdagdag}). Let $f,f'\colon \fY \to \fX$ be two finite  morphisms of $\fS$-formal schemes
such that $f _{\eta 0} =f' _{\eta 0}$. Using \cite[2.1.5]{Be2}, we get the isomorphism of $(\widehat{\cD} ^{(m)} _{\fY _\eta/\eta}, f _0 ^{-1} \widehat{\cD} ^{(m)} _{\fX _\eta/\eta})$-bimodules
$$\tau _{f,f'} \colon  f ^{\prime *} _\eta \widehat{\cD} ^{(m)} _{\fX _\eta/\eta}  \riso  f ^* _\eta \widehat{\cD} ^{(m)} _{\fX _\eta/\eta} .$$
Looking at the construction of this isomorphism,  we get the following explicit local description of this isomorphism:  suppose $\fX ^\sharp$ has logarithmic coordinates 
$u  _1, \dots, u  _d \in \Gamma (\fX, M _{\fX ^\sharp})$. Let $t  _1, \dots, t  _d $ be the image of  $u  _1, \dots, u  _d $ in $ \Gamma (\fX, \cO _{\fX ^\sharp})$. Then the image of  $1 \otimes 1$ is: 
$$ \sum _{\underline{i} \in \bbN ^d} (f ^{\prime *} (\underline{t} )- f ^{*} (\underline{t}) )^{\{\underline{i} \} _{(m)}} \otimes \underline{\partial} ^{<\underline{i}> _{(m)}}.$$
Hence, we get the factorization
$$\xymatrix{
{\widetilde{f} ^* \widetilde{\cD} ^{(m)\dag} _{\fX ^{\sharp}/\fS ^{\sharp}}}  \ar@{^{(}->}[r] ^-{} &  {j _*f ^{ *} _\eta \widehat{\cD} ^{(m)} _{\fX _\eta/\eta}  } 
\\ {\widetilde{f} ^{\prime *} \widetilde{\cD} ^{(m)\dag} _{\fX ^{\sharp}/\fS ^{\sharp}}} 
\ar@{^{(}->}[r] ^-{} \ar@{.>}[u] ^-{\tau _{f,f'}} &  {j _*f ^{\prime *} _\eta \widehat{\cD} ^{(m)} _{\fX _\eta/\eta} .} \ar[u] ^-{\tau _{f,f'}} } $$
\end{empt}

\begin{lemm}\label{univhomeo-eqcat-lemm}
Let $\fY ^{\sharp}\in \Sm ^\dag _{\fS ^\sharp}$ and let $\fX$ be an affine $\fS$-formal scheme such that there exists a $S_0$-morphism of the form $X_0 \to Y_0$. Then there exists a lifting $\fX \to \fY$ of $X \to Y$.
\end{lemm}

\begin{proof}
Let $M _{\fX}= \{s \in \cO _{\fX},\text{ such that $s | _{\fX _{\eta}} \in \cO  ^* _{\fX _{\eta}}$} \}$.
Set $\fX ^{\sharp} \coloneqq (\fX, M_{\fX})$ and $X _i ^{\sharp}\coloneqq \fX ^{\sharp} \times _{\fS ^\sharp} S  _i ^{\sharp}$.
The morphism $X _0 \to Y _0$ induces $X _0 ^{\sharp}\to Y_0 ^{\sharp}$. 
Using \cite[3.11]{Kato-logFontaine-Illusie}, since $Y ^{\sharp} _{i+1}/ S ^{\sharp} _{i+1}$ is log-smooth and $X _{i+1}$ is affine, 
then a morphism $X _i ^{\sharp}\to Y_i ^{\sharp}$ lifts to a morphism $X _{i+1} ^{\sharp}\to Y _{i+1} ^{\sharp}$. Hence, we are done by taking the inductive limits on $i$.
\end{proof}

\begin{thm}
\label{univhomeo-eqcat}
Let $f \colon \fY ^{\sharp}\to \fX ^{\sharp}$ be a finite surjective morphism of  $\Sm ^\dag _{\fS ^\sharp}$ 
such that the induced morphism  $f _0 ^\flat\colon Y \to X$ is a finite, surjective and radicial morphism. 
\begin{enumerate}[(a)]
\item \label{univhomeo-eqcat(a)} The functor
$\cE \mapsto \widetilde{f} ^* \cE$ 
from the category of left $\widetilde{\cD} ^{\dag} _{\fX ^{\sharp}/\fS ^{\sharp}}$-modules
to that of left
$\widetilde{\cD} ^{\dag} _{\fY ^{\sharp}/\fS ^{\sharp}}$-modules is an equivalence of categories.

\item \label{univhomeo-eqcat(b)} The functor 
$\cF \mapsto \widetilde{f} _+ \cF$ 
from the category of  coherent left $\widetilde{\cD} ^{\dag} _{\fY ^{\sharp}/\fS ^{\sharp},\bbQ}$-modules 
to that of 
coherent left $\widetilde{\cD} ^{\dag} _{\fX ^{\sharp}/\fS ^{\sharp},\bbQ}$-modules
is an equivalence of categories.
Moreover, for any coherent left $\widetilde{\cD} ^{\dag} _{\fY ^{\sharp}/\fS ^{\sharp},\bbQ}$-module $\cF$, the canonical morphism
$\cF \to \widetilde{f} ^* \widetilde{f} _+ \cF$ is an isomorphism.
\end{enumerate}

\end{thm}

\begin{proof}
Since the second assertion is checked similarly, let us only prove the first one. Since this is local, we can suppose $X$ affine.  Recall that $X$ is affine and normal if and only if $\fX$ is affine and normal
(use \cite[Lemma 6.1]{MonskyWashnitzer}).  Since $X$ and $Y$ are normal and $f _0$ is radicial, then 
$k (X) \subset k(Y)$ is radicial (see \cite[Exercice 5.3.9]{Liu-livre-02}).  Hence,  for $s$ large enough,  $k (Y) ^{p ^s} \subset k (X)$. 
Using \cite[Lemma  3.2.25]{Liu-livre-02} this yields that  increasing $s$ if necessary, there exists a morphism  $g _0 ^\flat \colon X \to Y ^{(s)}$ making commutative the diagram of $S$-schemes
\begin{equation} \label{univhomeo-eqcat-diag}
\xymatrix{
{Y}
\ar[d] _-{F ^s _{Y/S}}
\ar[r] ^-{f _0 ^\flat} 
& 
{X} 
\ar[ld] _-{g _0 ^\flat}
\ar[d] ^-{F ^s _{X/S}}
\\ 
{Y ^{(s)}}
\ar[r] _-{f _0 ^{^\flat (s)}} 
& 
{X ^{(s)}.} 
}
\end{equation}
Since $\fY ^{\sharp (s)}/\fS ^{\sharp}$ is log smooth and $\fX$ is affine,  then following \ref{univhomeo-eqcat-lemm} there exists a lifting $g \colon \fX ^\sharp \to \fY ^{\sharp (s)}$ of $g _0 $.
Hence, using \ref{eqcatFrobdescexplF+},
the transitivity with respect to the composition of morphisms of the functor $\widetilde{f} ^*$
from the category of left $\widetilde{\cD} ^{\dag} _{\fX ^{\sharp}/\fS ^{\sharp}}$-modules and the glueing isomorphisms of \ref{tau-dagm},  we get the part (\ref{univhomeo-eqcat(a)}).

Hence, using \ref{eqcatFrobdescexplF+}, the transitivity with respect to the composition of morphisms of the functor $\widetilde{f} _+$  from the category of  coherent left $\widetilde{\cD} ^{\dag} _{\fY ^{\sharp}/\fS ^{\sharp},\bbQ}$-modules, 
the glueing isomorphisms for coherent complexes by proper morphisms,  we get the first statement of the part (\ref{univhomeo-eqcat(b)}).
Following \ref{F+Fflatiso-left}, we have the canonical isomorphism  $\cF \riso F ^* F _+ \cF\riso \widetilde{f} ^*\widetilde{g} ^* \widetilde{g} _+  \widetilde{f} _+ \cF$.
Since this later is the composition of the adjoint morphisms:
$\cF \to \widetilde{f} ^* \widetilde{f} _+ \cF   \to \widetilde{f} ^*\widetilde{g} ^* \widetilde{g} _+  \widetilde{f} _+ \cF$,
this yields that $\cF \to \widetilde{f} ^* \widetilde{f} _+ \cF$ is injective. Hence (replace $f$ by $g$ and use the fact that $ \widetilde{f} _+ \cF$ is coherent since $f$ is proper), so is 
$\widetilde{f} _+ \cF   \to \widetilde{g} ^* \widetilde{g} _+  \widetilde{f} _+ \cF$. Since $\widetilde{f} ^*$ is exact, 
$\widetilde{f} ^* \widetilde{f} _+ \cF   \to \widetilde{f} ^*\widetilde{g} ^* \widetilde{g} _+  \widetilde{f} _+ \cF$ is therefore injective. 
Hence $\cF \to \widetilde{f} ^* \widetilde{f} _+ \cF$ is an isomorphism.

\end{proof}

\begin{lem}
\label{2.2.13-carocourbe}
Let $\fX _{\eta}$ be an affine smooth formal scheme over $\eta$ endowed with local coordinates. 
Let $\cE _{\eta}$ be a coherent left $\cD ^\dag _{\fX_{\eta}/\eta,\bbQ}$-module. 
The following conditions are equivalent
\begin{enumerate}[(a)]
\item The sheaf $\cE$ is coherent $\cO _{\fX_{\eta},\bbQ}$-module. 
\item $\Gamma ( \fX _{\eta}, \cE)$ is a $\Gamma (\fX_{\eta}, \cO _{\fX _{\eta},\bbQ})$-module of finite type.
\end{enumerate}
\end{lem}

\begin{proof}
We can copy the proof  of \cite[2.2.13]{caro_courbe-nouveau} word for word.
\end{proof}

\begin{lem}
\label{cvisoc-desc-finite}
Let $f \colon \fY _{\eta} \to \fX _{\eta}$ be a finite surjective morphism of smooth formal schemes over $\eta$. 
Let $\cE$ be a left $\cD ^\dag _{\fX_{\eta}/\eta,\bbQ}$-module. 
The following conditions are equivalent
\begin{enumerate}[(a)]
\item \label{cvisoc-desc-finite(a)} $\cE$ is a coherent $\cD ^\dag _{\fX_{\eta}/\eta,\bbQ}$-module which is also 
$\cO _{\fX_{\eta},\bbQ}$-coherent. 
\item \label{cvisoc-desc-finite(b)} $\cE$ is $\cO _{\fX_{\eta},\bbQ}$-coherent. 
\item \label{cvisoc-desc-finite(c)} $f ^* (\cE)$ is $\cO _{\fY,_{\eta}\bbQ}$-coherent. 
\item \label{cvisoc-desc-finite(d)} $f ^* (\cE)$ is coherent $\cD ^\dag _{\fY_{\eta}/\eta,\bbQ}$-module which is also 
$\cO _{\fY_{\eta},\bbQ}$-coherent. 
\end{enumerate}

\end{lem}

\begin{proof}
In the same way as \cite[4.1.4]{Be1}), we get that the equivalence between (\ref{cvisoc-desc-finite(a)}) and (\ref{cvisoc-desc-finite(b)}) and between (\ref{cvisoc-desc-finite(c)}) and (\ref{cvisoc-desc-finite(d)}). 
Since this is local, we can suppose $\fX_{\eta}$ affine. Using \ref{2.2.13-carocourbe}, if  $f ^* (\cE)$ is $\cO _{\fY_{\eta} ,\bbQ}$-coherent then $\Gamma (\fY_{\eta}, f ^* (\cE))$ is $\Gamma (\fY_{\eta}, \cO _{\fY _{\eta},\bbQ})$-coherent.
Since $f$ is locally free of finite type (see the proof of \ref{D2Dmdagpre}), this yields that $\Gamma (\fX_{\eta}, \cE)$ is $\Gamma (\fX_{\eta}, \cO _{\fX _{\eta},\bbQ})$-coherent. 
Using \ref{2.2.13-carocourbe} this yields that $\cE$ is coherent. The converse being obvious, this yields 
the equivalence between (\ref{cvisoc-desc-finite(b)}) and (\ref{cvisoc-desc-finite(c)}).
\end{proof}

\begin{thm}
\label{univhomeo-eqcat-isoc}
Let $f \colon \fY ^\sharp \to \fX ^{\sharp}$ be a projective (in the sense of \ref{Tr+PNdef}) morphism of strictly semistable log formal schemes over $\fS ^\sharp$.
We suppose that
the induced morphism 
$f _0\colon Y \to X$ is a finite, surjective and radicial. 

The functors 
$f _+$ and $f ^!$ are quasi-inverse equivalences of categories between 
the category of coherent left $\widetilde{\cD} ^{\dag} _{\fX ^{\sharp}/\fS ^{\sharp},\bbQ}$-modules 
which are also $\widetilde{\cO}  _{\fX,\bbQ}$-coherent,
and that of coherent left 
$\widetilde{\cD} ^{\dag} _{\fY ^{\sharp}/\fS ^{\sharp},\bbQ}$-modules
which are also $\widetilde{\cO}  _{\fY,\bbQ}$-coherent.
\end{thm}

\begin{proof}
1) Let $\cE$ be a coherent left $\widetilde{\cD} ^{\dag} _{\fX ^{\sharp}/\fS ^{\sharp},\bbQ}$-module 
which is also $\widetilde{\cO}  _{\fX,\bbQ}$-coherent. 
Then $f ^! (\cE) $ is a coherent left 
$\widetilde{\cD} ^{\dag} _{\fY ^{\sharp}/\fS ^{\sharp}}$-module
which is also $\widetilde{\cO}  _{\fY,\bbQ}$-coherent
(use the commutation with inverse image of the functor $\mathrm{real} _\fP$ of
\ref{ntn-realP}).

2)  Using Theorem \ref{univhomeo-eqcat}, we check that  the functor $f _+$  from the category of coherent left  $\widetilde{\cD} ^{\dag} _{\fY ^{\sharp}/\fS ^{\sharp}}$-modules which are also $\widetilde{\cO}  _{\fY,\bbQ}$-coherent is  fully faithful.

3) In the case where  $\fX ^\sharp = \fX _\eta$ and $\fY ^\sharp = \fY _\eta$, we prove in this step that  the functor $f ^!$ from  the category of coherent left $\widetilde{\cD} ^{\dag} _{\fX ^{\sharp}/\fS ^{\sharp},\bbQ}$-modules  which are also 
$\widetilde{\cO}  _{\fX,\bbQ}$-coherent to that of coherent left  $\widetilde{\cD} ^{\dag} _{\fY ^{\sharp}/\fS ^{\sharp}}$-modules which are also $\widetilde{\cO}  _{\fY,\bbQ}$-coherent is an equivalence of categories.
It follows from Lemma \ref{cvisoc-desc-finite} that such a functor is well defined.  Let $\cF$ be a  coherent left  $\widetilde{\cD} ^{\dag} _{\fY ^{\sharp}/\fS ^{\sharp}}$-module which is also $\widetilde{\cO}  _{\fY,\bbQ}$-coherent.
Via Theorem \ref{univhomeo-eqcat}, we get the isomorphism $\cF \to \widetilde{f} ^* \widetilde{f} _+ \cF$.
Using again Lemma \ref{cvisoc-desc-finite}, we get that $\widetilde{f} _+ \cF$ is a coherent left $\widetilde{\cD} ^{\dag} _{\fX ^{\sharp}/\fS ^{\sharp},\bbQ}$-module which is also $\widetilde{\cO}  _{\fX,\bbQ}$-coherent.
Hence, we are done. 

4) We check in this step that the trace morphism  $\Tr _{+,f}  \colon  f _+ (\widetilde{\omega} _{\fY ^{\sharp}/\fS ^{\sharp},\bbQ})  \to  \widetilde{\omega} _{\fX ^{\sharp}/\fS ^{\sharp},\bbQ}$ of \ref{Tr+PN}
is an isomorphism.   Using \ref{lem-faithfullflat-cor1}, we reduce to the case where $\fX ^\sharp = \fX _\eta$ and $\fY ^\sharp = \fY _\eta$. For $s$ large enough, there exists a diagram such as
\ref{univhomeo-eqcat-diag} and satisfying the same conditions. Following \ref{Tr+Fiso},  $\Tr _{+,F}$ is an isomorphism. Since  $g\circ f = F$ then by using \ref{Tr+transprop} (the proof of \ref{Tr+transprop} do not use \ref{univhomeo-eqcat-isoc}) we get $\Tr _{+,g} \circ g _+ (\Tr _{+,f}) = \Tr _{+,F}$.
Hence $g _+ (\Tr _{+,f})$ is injective and $\Tr _{+,g}$ is surjective. Via Theorem \ref{univhomeo-eqcat}, this impies $\Tr _{+,f}$ is injective. 
Since $g$ satisfies the same properties than $f$ then $\Tr _{+,g}$ is injective and therefore bijective.  This yields  $g _+ (\Tr _{+,f})$ and therefore $\Tr _{+,f}$ is bijective.

5) Let $\cE$ be  a coherent left $\widetilde{\cD} ^{\dag} _{\fX ^{\sharp}/\fS ^{\sharp},\bbQ}$-modules 
which is also $\widetilde{\cO}  _{\fX,\bbQ}$-coherent.
We deduce from the step 4),  the last isomorphism 
\begin{equation}
\notag
f _+ f ^! (\cE) \riso f _+ \left (\widetilde{\cO}  _{\fY,\bbQ}\otimes _{\widetilde{\cO}  _{\fY,\bbQ}} f ^! (\cE)  \right) 
\underset{{\cite[9.4.3.1]{Car25}}}{\riso} 
f _+ (\widetilde{\cO}  _{\fY,\bbQ})\otimes _{\widetilde{\cO}  _{\fX,\bbQ}} \E
\riso 
\E,
\end{equation}
where we also use \ref{DdagVSDhat} for the second isomorphism.

6) Let $\cF$ be  a coherent left $\widetilde{\cD} ^{\dag} _{\fY ^{\sharp}/\fS ^{\sharp},\bbQ}$-modules  which is also $\widetilde{\cO}  _{\fY,\bbQ}$-coherent.
Let us prove now that $f _+ (\cF) $ is a coherent left $\widetilde{\cD} ^{\dag} _{\fX ^{\sharp}/\fS ^{\sharp},\bbQ}$-module which is also $\widetilde{\cO}  _{\fX,\bbQ}$-coherent.
By using   \ref{letterBerthelotCaro2007}, we reduce to the case where $\fX ^\sharp = \fX _\eta$ and $\fY ^\sharp = \fY _\eta$. In that case, this checked in the proof of the part 3).
Finally, following Theorem \ref{univhomeo-eqcat}, we have the isomorphism $ f ^! f _+(\cF) \riso \cF$. Hence, we are done.
\end{proof}

\subsubsection{Finite descent of coherence}
\begin{dfn} \label{dfn-morph-trait}
A morphism of complete discrete valuation rings of equal characteristic $R \to  R'$ will refer to a local ring homomorphism such that a uniformizer of $R$ is not mapped to zero (which is equivalent to saying that $R \to R'$ is injective).
A scheme $S'$ is called a {\it trait of equal characteristic}, or simply a {\it trait} since we focus on the  equal characteristic case, if it is isomorphic to a scheme of the form  $\Spec R'$, where $R'$ is a complete discrete valuation ring of equal characteristic. 
A morphism of traits is a morphism $S''\to  S'$ corresponding to a morphism of complete discrete valuation rings of equal characteristic $R '\to  R''$ as above. 
Remark that such a morphism is a homeomorphism for the underlying morphism of topological spaces. 
Such a morphism is said to be a finite extension of traits if the extension $S ''\to S'$ is finite.
\end{dfn}

\begin{empt} \label{split-basechange}
Let $k [[t]] \to l[[u]]$ be a morphism of complete discrete valuation rings. The ring $l[[u]]$ is canonically endowed with a structure of $k$-algebra via the morphism
of residue fields $k \to l$. Since $k$ is perfect, using the ideas of the proof of 
\cite[II. Proposition 8.(i)]{Serre-corpslocaux}, 
we check that $k [[t]] \to l[[u]]$ is $k$-linear.
Hence we can split $k [[t]] \to l[[u]]$ in 
$k [[t]] \to l[[t]]$ and 
$l [[t]] \to l[[u]]$.

\end{empt}

\begin{empt}
\label{ll}
Let $S' \to S$ be a finite morphism of traits.
The morphism $S' \to S$ corresponds to a finite morphism of complete valuation rings of equal characteristic
of the form $k [[t]] \to l[[u]]$.
Let $\cW$ be the unramified extension of $\cV$ whose special fiber is 
$l$. 
We denote by $s'= \Spec l$ (resp. $\eta ' = \Spec l ((u))$) the scheme which corresponds to the closed point (resp. generic point) of $S '$.
By abuse of notation, we still denote by $s '= \Spf \cW$ (resp. $\eta '= \Spf \cW [[u]] \{ \frac{1}{u}\}$) 
the formal scheme over $\cW$ which corresponds to the closed point (resp. generic point) of $\S'$.
Let $S ^{\prime \sharp}:= (S ', M _{s'})$, $\fS ^{\prime \sharp}:= (\fS ', M _{s'})$ where $M_{s'}$ means the log structure induced by the special fiber,
and $\mathfrak{S} ': = \Spf \cW [[u]]$.

1) Suppose $l=k$.
We get $\mathfrak{S} '= \Spf \cV [[u]]$.
We have $t = u ^n P _0(u)$ with 
$P _0 (u) \in k [[u ]] ^*$ and $n = [k((u)) : k((t))]$ (see \cite[chap. II, corollary $1$ of proposition $3$]{Serre-corpslocaux}).
Set $A:=k [[t]]$ and $B :=k [[u ]]$. Then the image of $1,u, \dots, u ^{n-1}$ is a basis of the $k$-vector space $B/tB$.
Then $1,u, \dots, u ^{n-1}$ are linearly independent over $A$ (indeed if $\sum a _i u ^i=0$ is a non trivial relation, then we can suppose that one of the $a _i$ is not divisible by $t$, 
which gives a contraction reducing the relation modulo $t$).
Since $B$ is complete with respect to the $t$-adic topology, we get that the $A$-submodule of $B$ generated by $1,u, \dots, u ^{n-1}$ is equal to $B$.
This means $1,u, \dots, u ^{n-1}$ is a basis of the free  $k [[t]]$-algebra $k [[u]]$.

Choose a lifting  $\mathfrak{S}'\to \mathfrak{S}$ of the morphism $S ' \to S$, i.e. choose $P (u) \in \cV [[u ]] $ a lifting of $P _0 (u)$
and the map $\cV [[u ]] \to \cV [[u ]]$ sends $t$ to $u ^n P (u)$, i.e we have the relation $t = u ^n P (u)$.
Since $P _0 (u) \in k [[u ]] ^*$, then $P (u) \in \cV [[u ]] ^*$.Since $S' \to S$ is finite, then so is $\mathfrak{S} '\to \mathfrak{S}$.
As previously, since $\cV [[u ]]$ is complete for the $\pi$-adic topology, we establich that $1,u, \dots, u ^{n-1}$ is a basis of the free  $\cV [[t]]$-algebra $\cV [[u]]$.

Let $\fX ^{\sharp}$ be an object of $\Sm ^\dag _{\fS ^\sharp}$.
Let $\fX ' _{\eta'}$ be the open of $\fX ^{\prime \sharp}$ complementary to $X '_s$ i.e. the fiber of $\fX ^{\prime \sharp}$
over $\eta '= \Spf \cV [[u]] \{ \frac{1}{u}\}$. By stability of smoothness by base change, $\fX '_{\eta'} /\eta'$ is smooth. Let $j '\colon \fX '_\eta \hookrightarrow \fX ^{\prime \sharp}$ be the open immersion. 
Let  $f \colon \fX ^{\prime \sharp}\to \fX ^\sharp$ be the projection. Since $\fS ^{\prime \sharp}\to \fS ^\sharp$ is finite and surjective, then so is 
$f$. More precisely,  $f _*\cO _{\fX'}$ is a free $\cO _{\fX}$ module, with  $1,u,\dots, u ^{n-1}$ as a basis. 
The canonical homomorphism of sheaves of rings $f ^{-1} \cD^{(m)} _{\fX^{ \sharp}/\fS^{\sharp}} \to  \cD^{(m)} _{\fX^{\prime \sharp}/\fS^{\prime \sharp}}$
induces $f ^{-1}\widetilde{\cD}^{(m)} _{\fX^{ \sharp}/\fS^{\sharp}} \to  \widetilde{\cD}^{(m)} _{\fX^{\prime \sharp}/\fS^{\prime \sharp}}$.
Since $t = u ^n P (u)$ with $P (u) \in \cV [[u ]] ^*$, then $\cO _{\fX'} [\frac{1}{t}] = \cO _{\fX'} [\frac{1}{u}] $. Hence,  $\Gamma (\fX ', \widetilde{\cD}^{(m)} _{\fX^{\prime \sharp}/\fS^{\prime \sharp}} )$
is a free $\Gamma (\fX , \widetilde{\cD}^{(m)} _{\fX^{ \sharp}/\fS^{ \sharp}} )$-module with $1,u,\dots, u ^{n-1}$ as a basis.

Suppose  $\fX ^\sharp$ is affine, and $\fX ^\sharp/\fS ^\sharp$ has logarithmic coordinates. 
Fix $t _1,\dots, t _d  \in \Gamma (\fX, \widetilde{\cO} _{\fX})$   some overconvergent local coordinates  of $\fX ^\sharp/\fS ^\sharp$. 
We denote by  $t '_1,\dots, t '_d  \in \Gamma (\fX', \widetilde{\cO} _{\fX'})$   the overconvergent local coordinates 
of $\fX^{\prime \sharp}/\fS^{\prime \sharp}$ induced by base change by  $\fS^{\prime \sharp} \to \fS^{\sharp}$. Let  $A:= \Gamma (\fX, \cO _{\fX})$,   $A':= \Gamma (\fX', \cO _{\fX'})$.
 
 Let  $P'$ be an element of $\Gamma (\fX, \widehat{\cD} ^{(m)} _{\fX^{\prime} _{\eta '}/\eta ^{\prime}} ) $, 
the $p$-adic completion  of $\Gamma (\fX, \widetilde{\cD}^{(m)} _{\fX^{\prime \sharp}/\fS^{\prime \sharp}} ) $.
Following \ref{prop-wkcpdscplevelm}, we have
$z \in \Gamma (\fX ', \widetilde{\cD}^{(m)\dag} _{\fX^{\prime \sharp}/\fS^{\prime \sharp}} )$
if and only if we can write $z$ in the form
\begin{equation}
\label{wkcpdscplevelmapp1}
P' = \sum _{k\in \bbN}\sum _{\underline{i}\in \bbN ^d} a '_{\underline{i},k} \frac{1}{u ^{k}}\underline{\partial} ^{<\underline{i}> _{(m)}}
\end{equation}
where $a '_{\underline{i},k} \in A'$ are such that there exists a constant $c>0$ satisfying
$v _{p} ( a ' _{\underline{i},k}) \geq \frac{k+|\underline{i}|}{c} -1$,
for any $\underline{i} \in \bbN ^d$
and $k \in \bbN$.

We get the factorization
$\Gamma (\fX , \widetilde{\cD}^{(m)\dag} _{\fX^{ \sharp}/\fS^{ \sharp}} )
\to 
\Gamma (\fX ', \widetilde{\cD}^{(m)\dag} _{\fX^{\prime \sharp}/\fS^{\prime \sharp}} )$
given by 
$$\sum _{k\in \bbN}\sum _{\underline{i}\in \bbN ^d} a _{\underline{i},k} \frac{1}{t ^{k}}\underline{\partial} ^{<\underline{i}> _{(m)}}
 \mapsto 
 \sum _{k\in \bbN}\sum _{\underline{i}\in \bbN ^d} f ^* (a _{\underline{i},k}) P (u) ^{-k} \frac{1}{u ^{kn}}\underline{\partial} ^{<\underline{i}> _{(m)}}$$
where $a _{\underline{i},k} \in A$ are such that there exists a constant $c>0$ satisfying
$v _{p} ( a  _{\underline{i},k}) \geq \frac{k+|\underline{i}|}{c} -1$,
for any $\underline{i} \in \bbN ^d$
and $k \in \bbN$.
We easily compute that 
$\Gamma (\fX ', \widetilde{\cD}^{(m)\dag} _{\fX^{\prime \sharp}/\fS^{\prime \sharp}} )$
is a free $\Gamma (\fX , \widetilde{\cD}^{(m)\dag} _{\fX^{ \sharp}/\fS^{ \sharp}} )$-module with $1,u,\dots, u ^{n-1}$ as a basis.  Taking the limit, we get that 
$\Gamma (\fX ', \widetilde{\cD}^{\dag} _{\fX^{\prime \sharp}/\fS^{\prime \sharp}} )$ is a free $\Gamma (\fX , \widetilde{\cD}^{\dag} _{\fX^{ \sharp}/\fS^{ \sharp}} )$-module with $1,u,\dots, u ^{n-1}$ as a basis. 

2) We consider now the extension $k [[t]] \to l[[t]]$. Since $k$ is perfect, then $k [[t]] \to l[[t]]$ and $\cV [[t]] \to \cW [[t]]$ are finite and separable. 

\end{empt}

\begin{lemm}
\label{lem-desc-coh-chgbase}
We keep notation \ref{ll}. 
Let $f \colon \fX^{\prime \sharp} \to \fX^{\sharp}$ be the canonical projection. 
Let $Z$ be a divisor of $X$ containing $X _s$ and $Z':= f ^{-1} (Z)$.
\begin{enumerate}[(a)]
\item  
\label{lem-desc-coh-chgbase0}
The canonical 
homomorphism 
$\cD ^\dag _{\fX^{\prime \sharp}/\fS ^{\prime \sharp}} (\hdag Z') _{\bbQ} 
\to 
\cD ^\dag _{\fX^{\prime \sharp}\to \fX^{\sharp}/\fS ^{\prime \sharp}\to \fS ^{\sharp} } (\hdag Z') _{\bbQ} $
is an isomorphism. 
The composite morphism
$f ^{-1}\cD ^\dag _{\fX^{\sharp}/\fS^{\sharp}} (\hdag Z) _{\bbQ} 
\to 
\cD ^\dag _{\fX^{\prime \sharp}\to \fX^{\sharp}/\fS ^{\prime \sharp}\to \fS ^{\sharp} } (\hdag Z') _{\bbQ} 
\liso
\cD ^\dag _{\fX^{\prime \sharp}/\fS ^{\prime \sharp}} (\hdag Z') _{\bbQ} $
is a homomorphism of rings. 
Hence, if $\cE$ is a coherent $\cD ^\dag _{\fX^{\sharp}/\fS^{\sharp}} (\hdag Z) _{\bbQ}$-module, then
$f _Z ^! (\cE) \riso \cD ^\dag _{\fX^{\prime \sharp}/\fS ^{\prime \sharp}} (\hdag Z') _{\bbQ} \otimes _{f ^{-1}\cD ^\dag _{\fX^{\sharp}/\fS^{\sharp}} (\hdag Z) _{\bbQ}}
f ^{-1}\cE$, where $f ^! _Z$ is the extraordinary inverse image of $\fX^{\prime \sharp} \to \fX$ above $\fS^{\prime \sharp} \to \fS$
 with overconvergent singularities along $Z$, i.e. $f ^! _Z$ is the base change inverse image.

\item 
\label{lem-desc-coh-chgbase1}
Suppose $\fX$ is affine.
Let $\cE$ be a coherent $\cD ^\dag _{\fX^{\sharp}/\fS^{\sharp}} (\hdag Z) _{\bbQ}$-module.
 Then the canonical morphisms
$$\Gamma (\fS ',\cO _{\fS'})  \otimes _{\Gamma (\fS ,\cO _{\fS})} \Gamma (\fX, \cE) 
\to 
D ^\dag _{\fX^{\prime \sharp}/\fS ^{\prime \sharp}} (\hdag Z') _{\bbQ} 
\otimes _{D ^\dag _{\fX^{\sharp}/\fS^{\sharp}} (\hdag Z) _{\bbQ}}
\Gamma (\fX, \cE) 
\to 
\Gamma (\fX', f _Z ^! (\cE)) $$
 are isomorphisms.
 Moreover, 
$D ^\dag _{\fX^{\prime \sharp}/\fS^{\prime \sharp}} (\hdag Z') _{\bbQ}$ is a faithfully flat 
$D ^\dag _{\fX^{\sharp}/\fS^{\sharp}} (\hdag Z) _{\bbQ}$-module for both left or right structure.

\item 
\label{lem-desc-coh-chgbase2}  For any  $\cD ^\dag _{\fX^{\sharp}/\fS^{\sharp}} (\hdag Z) _{\bbQ}$-module $\cE$, the canonical morphisms
\begin{equation}
\notag
f ^* (\cE) :=  \widetilde{\cO} _{\fX'} \otimes _{f ^{-1} \widetilde{\cO} _{\fX}} f ^{-1}\cE  \to \cO _{\fX'} (\hdag Z') _{\bbQ} \otimes _{f ^{-1}\cO_{\fX} (\hdag Z) _{\bbQ}} f ^{-1}\cE 
\to \cD ^\dag _{\fX^{\prime \sharp}/\fS ^{\prime \sharp}} (\hdag Z') _{\bbQ} \otimes _{f ^{-1}\cD ^\dag _{\fX^{\sharp}/\fS^{\sharp}} (\hdag Z) _{\bbQ}} f ^{-1}\cE 
\end{equation}
are isomorphisms. 

 \item 
 \label{lem-desc-coh-chgbase3}
 Let $\phi \colon \E' \to \E$ be a morphism of $ \widetilde{\cO}_{\fX}$-modules.
 Then $\phi$ is an isomorphism if and only if  $f ^* (\phi)$ 
is an isomorphism.
\end{enumerate}

\end{lemm}

\begin{proof}
Via \ref{ll}, this is similar to \cite[8.2.1]{caro-6operations} or \cite[9.2.7.2]{Car25}.
\end{proof}

\begin{prop}\label{desc-coh-chgbase}
With notation \ref{lem-desc-coh-chgbase}, let $\E$ be a   $\cD ^\dag _{\fX^{\sharp}/\fS^{\sharp}} (\hdag Z) _{\bbQ}$-coherent module.
Then $\E$ is a coherent $\widetilde{\cD} ^\dag _{\fX^{\sharp}/\fS^{\sharp}, \bbQ}$-module if and only if 
$f _Z ^! (\E) $ is a coherent $\widetilde{\cD} ^\dag _{\fX^{\prime \sharp}/\fS ^{\prime \sharp}, \bbQ}$-module.  
\end{prop}

\begin{proof}
We proceed as in  \cite[8.2.2]{caro-6operations}  or \cite[9.2.7.4]{Car25}.
\end{proof}

\section{Berthelot-Kashiwara's theorem for weakly closed immersions in the context of semistable formal schemes}

\subsection{Strictly semistable formal scheme}

\subsubsection{Complements on log smoothness}

\begin{lem}
\label{lem-nil-immersion}
Let $u\colon Z \hookrightarrow X$ be a nil-immersion (i.e. a closed immersion which is defined by a nil ideal). 
Let $\theta \colon M\to N$ be a morphism of fine logarithmic structures of $X$ such that 
$u ^* (\theta) \colon u ^{*} M \to u ^{*} N $ is a isomorphism. 
Then $\theta$ is an isomorphism.
\end{lem}

\begin{proof}
Let $\underline{x}$ be a geometric point of $X$. From \cite[II.2.2.6.2]{Ogus-Logbook},
it is sufficient to check that $\theta _{\underline{x}}\colon M_{\underline{x}}\to N _{\underline{x}}$ is an isomorphism.
Since $u$ is a nil-immersion, $u ^{-1}$ is the identity. Hence, since $u ^* (\theta)$ is an isomorphism, 
using \cite[II.1.2.1]{Ogus-Logbook}
the morphism 
$M_{\underline{x}} / \cO ^* _{X, \underline{x}}
\to N _{\underline{x}} / \cO ^* _{X, \underline{x}}$ 
is an isomorphism.
We remark that $\theta _{\underline{x}}$ is a morphism of fine monoids 
(indeed, we can suppose that $M$ and $N$ have fine charts, then we use 
\cite[I.1.2.2]{Ogus-Logbook}
and the commutation of the stalks at $\underline{x}$ with
amalgamated sums). 
Since $\theta _{\underline{x}}$ is sharp, 
from 
\cite[I.4.1.2]{Ogus-Logbook}
this yields that 
$\theta _{\underline{x}}$ is an isomorphism.
\end{proof}

\begin{lem}
\label{lem-nil-imm-flat-pre}
Let $u\colon T' \hookrightarrow T$ be a nil-immersion of fine log-schemes annihilated by a power of $p$. 
Let $f\colon X \to Y$ be a morphism of fine log-schemes over $T$ such that
$X$ is affine and flat over $T$. Let $f' \colon X '\to Y'$ be the base change of $f$ under $u$.
Assume we are given a (fine) chart $Q\to M _{Y}$, 
$f '$ has a chart (with the convention  \cite[2.9]{Kato-logFontaine-Illusie} monoids are fine)  
of the form $( P  \to M  _{X'}, Q\to M _{Y'}, Q \to P)$ such that 
$Q\to M _{Y'}$ is the morphism induced by $Q\to M _{Y}$
and the induced $Y'$-morphism 
$g' \colon X '\to A _{P}\times _{A _Q}Y'$ (remark that 
since monoids are fine, the fibered product 
$A _{P}\times _{A _Q}Y'$ in the category of log schemes or 
fine log schemes are the same)
is étale
and the kernel and the torsion part of the cokernel of $Q ^\mathrm{gp}\to P^\mathrm{gp}$ are finite groups whose orders are prime to $p$.

Then there exists an étale $Y$-morphism
$g\colon X \to A _{P}\times _{A _Q}Y$ making commutative the following diagram
\begin{equation}
\label{lem-nil-imm-flat-pre-diag}
\xymatrix{
{} 
& 
{X'} 
\ar[ld] _-{f'}
\ar@{^{(}->}[r] ^-{} 
\ar[d] ^-{g'}
& 
{X} 
\ar[rd] ^-{f}
\ar@{.>}[d] ^-{g}
\\ 
{Y'} 
& 
{ A _{P}\times _{A _Q}Y'}
\ar@{^{(}->}[r] ^-{} 
\ar[l] ^-{}
& 
{A _{P}\times _{A _Q}Y} 
\ar[r] ^-{}
& 
{Y } 
}
\end{equation}
In particular, $f$ is log smooth 
and $P$ is also a chart of $X$.
\end{lem}

\begin{proof}
From   \cite[3.11]{Kato-logFontaine-Illusie} (see its proof for a more general statement), since $X$ is affine
and since $A _{P}\times _{A _Q}Y \to Y$ is log smooth
then
we can find a $Y$-morphism of the form 
$g \colon X \to A _{P}\times _{A _Q}Y$ making commutative the diagram 
\ref{lem-nil-imm-flat-pre-diag}.
Since $X \to T$ is flat, we get from \cite[IV.17.8.2]{EGAIV4}
that the underlying morphism of schemes of $g $ is étale. 
Using \ref{lem-nil-immersion}, the morphism $g$ is strict, hence étale.
Hence $f$ is log smooth. 
The fact that $P$ is also a chart of $X$ comes from the fact that a composition of strict morphisms
is a strict morphism and then 
$X \overset{g}{\to}A _{P}\times _{A _Q}Y \to A _{P}$ is strict.
\end{proof}

\begin{lem}
\label{lem-nil-imm-flat}
Let $u\colon T' \hookrightarrow T$ be a nil-immersion of fine log-schemes annihilated by a power of $p$. 
Let $f\colon X \to Y$ be a morphism of fine log-schemes over $T$ such that
$X$ is flat over $T$. Let $f' \colon X '\to Y'$ be the base change of $f$ under $u$.
If $f'$ is log smooth then so is $f$.
 
\end{lem}

\begin{proof}
Since the Lemma is étale locally on $X$ and since $f'$ is log smooth, 
then 
we can suppose $X$ affine and that we are given a (fine) chart $Q\to M _{Y}$. 
Using \cite[3.5]{Kato-logFontaine-Illusie}, we can also suppose 
that $f '$ has a chart 
of the form $( P  \to M  _{X'}, Q\to M _{Y'}, Q \to P)$ such that 
$Q\to M _{Y'}$ is the morphism induced by $Q\to M _{Y}$
and the induced morphism 
$X '\to A _{P}\times _{A _Q}Y'$ 
is étale
and the kernel and the torsion part of the cokernel of $Q ^\mathrm{gp}\to P^\mathrm{gp}$ are finite groups whose orders are prime to $p$. 
We conclude by using 
\ref{lem-nil-imm-flat-pre}.
\end{proof}

\begin{lem}
\label{lem-nil-imm-flat-pre-formal}
Let $f\colon \fX \to \fY$ be a morphism of log formal schemes over $\cV$.
We suppose that $\fX$ is affine and has no $p$-torsion.
Assume we are given a fine chart $Q\to M _{\fY}$, 
$f _0$ has a (fine) chart 
of the form $( P  \to M  _{X_0}, Q\to M _{Y _0}, Q \to P)$ such that 
$Q\to M _{Y _0}$ is the morphism induced by $Q\to M _{\fY}$
and 
the induced $Y _0$-morphism 
$g _0 \colon X _0\to A _{P}\times _{A _Q}Y _0$ 
is étale
and 
the kernel and the torsion part of the cokernel of $Q ^\mathrm{gp}\to P^\mathrm{gp}$ are finite groups whose orders are prime to $p$. 
Then there exists an étale $\fY$-morphism
$g\colon \fX \to \mathfrak{A} _{P}\times _{\mathfrak{A} _Q}\fY$ making commutative the following diagram
\begin{equation}
\label{lem-nil-imm-flat-pre-formal-diag}
\xymatrix{
{}  
& 
{X _0} 
\ar[ld] _-{f_0}
\ar@{^{(}->}[r] ^-{} 
\ar[d] ^-{g_0}
& 
{\fX} 
\ar[rd] ^-{f}
\ar@{.>}[d] ^-{g}
\\ 
{Y_0} 
& 
{ A _{P}\times _{A _Q}Y_0}
\ar@{^{(}->}[r] ^-{} 
\ar[l] ^-{}
& 
{\mathfrak{A} _{P}\times _{\mathfrak{A} _Q}\fY} 
\ar[r] ^-{}
& 
{\fY } 
}
\end{equation}
In particular, $f$ is log smooth 
and $P$ is also a chart of $\fX$.

\end{lem}

\begin{proof}
By induction on $i\geq 0$, we construct compatible liftings
$g _i\colon X _i\to A _{P}\times _{A _Q}Y _i$
of $g _0$. Indeed, if such $g _i$ exists, 
from \ref{lem-nil-imm-flat-pre}, 
there exists 
$g _{i+1}\colon X _{i+1}\to A _{P}\times _{A _Q}Y _{i+1}$
making commutative the diagram below
\begin{equation}
\label{lem-nil-imm-flat-pre-diag-ito+1}
\xymatrix{
{} 
& 
{X_i} 
\ar[ld] _-{f_i}
\ar@{^{(}->}[r] ^-{} 
\ar[d] ^-{g_i}
& 
{X _{i+1}} 
\ar[rd] ^-{f _{i+1}}
\ar@{.>}[d] ^-{g_{i+1}}
\\ 
{Y_i} 
& 
{ A _{P}\times _{A _Q}Y_i}
\ar@{^{(}->}[r] ^-{} 
\ar[l] ^-{}
& 
{A _{P}\times _{A _Q}Y_{i+1}} 
\ar[r] ^-{}
& 
{Y_{i+1} } 
}
\end{equation}

\end{proof}

\begin{empt}
Let $f\colon \fX \to \fY$ be a morphism of affine formal schemes over $\cV$.
We recall that $f$ is flat (as morphism of locally ringed spaces) if and only 
$\Gamma (\fY ,\cO _{\fY})
\to 
\Gamma (\fX ,\cO _{\fX})$
is flat. For instance, this is a consequence of 
\cite[$0 _{III}$.10.2.2]{EGAIII1}
and \cite[$0 _{I}$.7.6.18]{EGAIII1}.
\end{empt}

\begin{prop}
\label{smooth-flat-rel}
Let $f\colon \fX \to \fY$ be a morphism of log-formal schemes over $\cV$.
We suppose that $\fX$ has no $p$-torsion, i.e. $\fX$ is flat over $\cV$.
\begin{enumerate}[(a)]
\item \label{smooth-flat-rel1} $f$ is log smooth (resp. smooth) if and only if so is $f _0$. 
\item \label{smooth-flat-rel2} $f$ is flat if and only if so is $f _0$. 
\end{enumerate}
\end{prop}

\begin{proof}
The first statement comes from \ref{lem-nil-imm-flat-pre-formal} and \cite[3.5]{Kato-logFontaine-Illusie}. 
Let us check the second one.
From \cite[III.5.2, Theorem 1]{bourbaki3-4}, since $p \cO _\fS$ is included in the radical of $\cO _\fS$, 
$f$ is flat is and only if $f _i$ is flat for any integer $i$.
Moreover, using \cite[IV.11.3.10]{EGAIV4}, since $X _i$ is flat over $\cV /\pi ^{i+1}\cV$, then 
the flatness of $f _0$ is equivalent to that of $f _i$.
\end{proof}

\begin{prop}
\label{smooth-flat}
Let $f\colon \fX ^\sharp \to \fS ^{\sharp}$ be a morphism of log-formal schemes.
The following conditions are equivalent
\begin{enumerate}[(a)]
\item \label{smooth-flat1} $f$ is log smooth. 
\item \label{smooth-flat4} $f _i$ is log smooth for any integer $i$.
\item \label{smooth-flat2} $f$ is flat and $f _0 \colon X ^{\sharp} \to S ^{\sharp}$ is log-smooth. 
\item \label{smooth-flat3} $f _0$ is smooth and $f _i\colon X _i ^{\sharp} \to S _i ^{\sharp}$ is flat for any integer $i$.
\end{enumerate}

\end{prop}

\begin{proof}
The equivalence between \ref{smooth-flat1} and \ref{smooth-flat4} is by definition. 
The equivalence between 
\ref{smooth-flat2} and \ref{smooth-flat3} is a consequence of \ref{smooth-flat-rel}. 
We get the implication $\ref{smooth-flat4}\Rightarrow \ref{smooth-flat3}$ by using 
\cite[4.4 and 4.5]{Kato-logFontaine-Illusie}.
Finally the implication $\ref{smooth-flat3}\Rightarrow \ref{smooth-flat4}$ is a consequence of \ref{lem-nil-imm-flat}.
\end{proof}

\begin{empt}
Let $f\colon \fX ^\sharp \to \fS ^{\sharp}$ be a log-smooth morphism of log-formal schemes.
Using \cite[4.3.5]{Kato-logFontaine-Illusie}, etale locally on $X$, 
there exists a chart 
$(P _X \to M _{X ^\sharp}, \bbN \to M _{S ^\sharp}, h \colon \bbN \to P)$
extending the given 
$\bbN \to M _{S ^\sharp}$
satisfying the following conditions
\begin{enumerate}[(i)]
\item the kernel and the torsion part of the cokernel of 
$h ^\mathrm{gr}
\colon \Z \to P ^\mathrm{gr}$ 
are finite groups of orders invertible on $X$ ;
\item the induced morphism 
$X ^\sharp \to S ^\sharp \times _{A _\bbN} A _P$ is étale (i.e. log étale and strict).
\end{enumerate}
Since $S ^\sharp \to \Spec k \times _{\Spec \Z} A _\bbN$ is flat and strict, 
then so is 
$S ^\sharp \times _{A _\bbN} A _P
\to 
\Spec k \times _{\Spec \Z} A _P$.
Hence, 
$X ^\sharp
\to 
\Spec k \times _{\Spec \Z} A _P$
is flat and strict.
Following the claim at the end of 
 \cite[1.5]{Kato-logFontaine-Illusie},
 this implies that $M _X$ is the submonoid sheaf of $\cO _X$ generated by 
 $\cO ^* _X$ and $P$.
\end{empt}

\subsubsection{Definitions and first properties}

\begin{empt} [Strictly semi-stable log formal schemes over $S$]
\label{dfn-sss}
Let  $X$ be an integral $S$-scheme of finite type. 
Let $X _i$, $i\in I$ be the irreducible components of $X_s$. Put $X _J:= \cap _{j \in J} X _j$ (scheme-theoretic intersection), for a nonempty subset $J$ of $I$. 
We suppose that $X$ is ``strictly semi-stable over $S$''
i.e. that this means that $X/S$ satisfy the following properties  (see \cite[2.16]{dejong}, and also \stack{0CBN}): 
\begin{enumerate}[(a)]
\item $X _\eta$ is smooth over $\kappa (\eta)$, 
\item $X _s $ is a reduced scheme, i.e. $X _s = \cup _{i\in I} X _i$ scheme-theoretically,
\item for each $i\in I$, $X _i$ is a divisor on $X$, 
\item for each nonempty $J \subset I$, the scheme $X_J$  is smooth over $k (s)$ and has codimension $\sharp J$ in $X$.
\end{enumerate}
Remark that since $k= \kappa (s)$ is perfect, then conditions $b), c), d)$ are equivalent 
to say that $X _s$ is a divisor with strict normal crossing on $X$ (see the definition \cite[2.10]{dejong}).
$X ^{\sharp}:= (X, M _{X _s})$. Recall that from \cite[3.14]{Kato-logFontaine-Illusie}, 
$(X, M _{X _s})$ is log smooth over $S ^{\sharp}$, where $M _{X _s}$ is the log structure associated with the strict normal crossing divisor $X _s$ of $X$. 
\end{empt}

\begin{dfn}
\begin{enumerate}[(a)]
\item Let $X ^\sharp$ be a fine log scheme over $S ^\sharp$.
We say that $X ^\sharp=(X, M _X)$ ``strictly semistable over $S ^\sharp$'' 
if $X$ is strictly semistable over $S $ and $M _X$ is equal to $M _{X _s}$, the log structure associated with
the strict normal crossing divisor $X _s$ of $X$.  Recall that 
$M _{X _s} $ is the submonoid of $\cO _{X}$ of sections invertible outside $X _s$.
\item Let $\fX ^\sharp$ be a fine $p$-adic formal log scheme over $\fS ^\sharp$.
We say that $\fX ^\sharp$ is ``strictly semistable over $\fS ^\sharp$'' if $\fX ^\sharp$ is log smooth over $\fS ^\sharp$ and if $X ^\sharp$ is strictly semistable over $S ^\sharp$. 
\end{enumerate}
\end{dfn}

\begin{rem}
\label{rem-semistable lifting}
Let $X ^\sharp$ be a strictly semi-stable log scheme over $S ^\sharp$ 
such that $X$ is affine. 
Using \cite[3.14]{Kato-logFontaine-Illusie} (we can copy the proof of \cite[III.6.10]{sga1}), 
there exists a semi-stable formal log scheme $\fX ^{\sharp}$ over $\fS ^{\sharp}$ whose special fiber is $X ^\sharp$.
\end{rem}

\begin{rem}
\label{rem-morphi-sssvar}
Let $X ^\sharp$ and $Y ^{\sharp}$
be two strictly semi-stable log schemes over $S ^\sharp$.

\begin{enumerate}[(a)]
\item The canonical application from the set of morphisms of log schemes of the form
$X ^\sharp\to Y ^{\sharp}$
to the set of morphism of schemes 
$X \to Y $ 
is injective. 
Moreover, if $f \colon X \to Y $ is a morphism such that 
$f ^{-1} ( Y _s) \subset X _s$ (i.e. $f$ induces the morphism 
$X _\eta \to Y _\eta$), then 
$f$ comes from a morphism of log schemes 
of the form
$X ^\sharp\to Y ^{\sharp}$.

\item The data of a morphism of $S$-schemes of the form
$X \to Y$ is equivalent to that of a morphism
of log schemes over $S ^\sharp$ of the form
$X ^\sharp\to Y ^{\sharp}$.
\end{enumerate}

\end{rem}

\begin{empt}
Let $\alpha\colon \fS ^\sharp \to \mathfrak{A} _{\bbN}$ be the chart  
defined by the map $ \bbN \to \cV [[t]]$ sending $1$ to $t$.
Let $\alpha _{(r)}\colon \fS ^\sharp _{(r)} \to \mathfrak{A} _{\bbN ^{r+1}} $ 
be the morphism defined by the map $\bbN ^{r+1} \to \cV [[t]] \{ T _0, \dots, T _{r}\}/(T _0 \cdots T _{r} -t)$
sending the canonical basis $e _0, \dots, e _{r}$ to respectively the classes of $T _0, \dots, T _{r}$.
Let $\Delta \colon \bbN \to \bbN ^{r+1}$ be the diagonal morphism. 
We notice that we have the commutative diagram with cartesian squares of formal log schemes over $\fS$
\begin{equation}
\label{S(r)chart}
\xymatrix{
{\fS ^\sharp _{(r)}} 
\ar[r] ^-{\beta _{(r)}}
\ar[d] ^-{\alpha _{(r)}}
\ar@{}[rd] ^-{}|\square
& 
{\fS ^\sharp \times\mathfrak{A} _{\bbN ^{r} } }
\ar[r] ^-{p _1}
\ar@{}[rd] ^-{}|\square
\ar[d] ^-{}
&
{\fS ^\sharp} 
\ar[d] ^-{\alpha}
\\ 
{\mathfrak{A} _{\bbN ^{r+1}} } 
\ar@/_0,3cm/[rr] _-{\mathfrak{A} _\Delta}
\ar[r] ^-{\Delta _{1,r}}
& 
{\mathfrak{A} _{\bbN } \times \mathfrak{A} _{\bbN ^{r}} }
\ar[r] ^-{p _1}
&
{\mathfrak{A} _{\bbN} ,} 
}
\end{equation}
where $\Delta _{1,r}$ is the morphism induced by the morphism of monoids $\bbN \oplus \bbN ^{r} \to \bbN ^{r+1}$
sending $(n, n _1,\dots, n _r)$ to $(n,n _1+n,\dots, n _r +n)$, $p _1$ is the first projection, 
$\beta _{(r)}$ is the $\fS ^\sharp$-morphism induced by the map $\bbN ^r \to \cV [[t]] \{ T _0, \dots, T _{r}\}/(T _0 \cdots T _{r} -t)$ 
sending $e _1, \dots, e _r$ to respectively the class of $T_1,\dots,   T _r$. 
Indeed, the cartesianity of the right square is obvious.
The cartesianity of the composition of both squares is easy: first we check it 
without log structure and then we notice that 
$\alpha$ and  $\alpha _{(r)}$ are strict. 
This yields that $\bbN ^{r+1}$ is a prelog structure of  $\fS ^\sharp _{(r)}$.
Since $\Delta _{1,r}$ is log-étale, we get that $\beta _{(r)}$ is log-étale.
\end{empt}

\begin{lem}
\label{semistable-local}
Let $\fX ^\sharp$ be a formal log scheme over $\fS ^\sharp$.
We suppose that $\fX$ is affine and has no $p$-torsion.
Assume 
there exists a smooth $S ^\sharp$-morphism of the form
$g _0 \colon X ^\sharp \to S _{(r)} ^\sharp$. 
Then there exists a smooth $\fS ^\sharp$-morphism 
$g\colon \fX ^\sharp \to \fS _{(r)} ^\sharp$ making commutative the following diagram
\begin{equation}
\label{semistable-local-diag}
\xymatrix{
{}   
& 
{X ^\sharp} 
\ar[ld] _-{f_0}
\ar@{^{(}->}[r] ^-{} 
\ar[d] ^-{g_0}
& 
{\fX ^\sharp} 
\ar[rd] ^-{f}
\ar@{.>}[d] ^-{g}
\\ 
{S  ^\sharp} 
& 
{ S _{(r)} ^\sharp}
\ar@{^{(}->}[r] ^-{} 
\ar[l] ^-{}
& 
{\fS _{(r)} ^\sharp} 
\ar[r] ^-{}
& 
{\fS  ^\sharp } 
}
\end{equation}
In particular, $\fX ^\sharp$ is strictly semistable over $\fS ^\sharp$ and $\bbN ^{r+1}$ is a prelog structure on $\fX ^{\sharp}$.

\end{lem}

\begin{proof}
The check is analogous to that of \ref{lem-nil-imm-flat-pre-formal} (use also \ref{smooth-flat-rel}).
\end{proof}

\begin{prop}
\label{etalecoordsss}
Let $\fX ^\sharp$ be a strictly semistable formal log scheme over $\fS ^\sharp$.
Let $x$ be point of $X$. 
There exists an affine open $\fU ^\sharp $ formal subscheme of $\fX ^\sharp$ which contains $x$
and an etale $\fS ^\sharp$-morphism of the form 
$h\colon \fU ^{\sharp} \to (\widehat{\mathbb{G}} _{\mathrm{m} }) _\fS ^{s}  \times _\fS \fS ^\sharp _{(r)}  $.
\end{prop}

\begin{proof}
Let $\fU ^\sharp $ be an affine open  formal subscheme of $\fX ^\sharp$ which contains $x$.
Shrinking $\fU ^\sharp $ if necessary, it is known that
there exists a smooth morphism of the form 
$U\to S _{(r)}$ where  
$U$ is the underlying scheme of $U ^\sharp$ and where $r+1$ equals to the number 
of irreducible components of $X _s$ in which $x$ is lying. 
We get an induced morphism $g _0\colon U ^\sharp \to S _{(r)} ^{\sharp}$.
We check that this morphism $g _0$ is strict and then smooth.
From \ref{semistable-local}, there exists a smooth $\fS ^\sharp$-morphism 
$g\colon \fU ^\sharp \to \fS _{(r)} ^\sharp$ 
which induces $g _0$.
Hence, shrinking $\fU ^\sharp $ if necessary, 
we can suppose that there exists
an étale morphism of the form 
$\fU ^{\sharp} \to (\widehat{\mathbb{G}} _{\mathrm{m} }) _\fS ^{s}   \times _\fS \fS ^\sharp _{(r)}  $.
\end{proof}

\begin{empt}
[Log coordinates for strictly semistable formal schemes]
\label{nota-loc-coord}
Let $\fX ^\sharp$ be a strictly semistable formal log scheme over $\fS ^\sharp$.
Let $x$ be point of $X$. 
From \ref{etalecoordsss}, there exists an affine open $\fU ^\sharp $ formal subscheme of $\fX ^\sharp$ which contains $x$
and an etale $\fS ^\sharp$-morphism of the form 
$h\colon \fU ^{\sharp} \to (\widehat{\mathbb{G}} _{\mathrm{m} }) _\fS ^{s}   \times _\fS \fS ^\sharp _{(r)}  $.
Let $u _0, \dots, u _{r}$ be the image of the class of 
$T _0, \dots, T _r \in \cV [[t]] \{ T _0, \dots, T _{r}\}/(T _0 \cdots T _{r} -t) = \Gamma (\fS ^\sharp _{(r)}, \cO _{\fS ^\sharp _{(r)}})$ under $h$. 
We get the relation 
$t = u _0 \cdots u _r$
and
 we have
$u _1,\dots, u _r\in \Gamma (\fU ^\sharp, \cM _{\fX ^{\sharp}})\subset \Gamma (\fU, \cO _{\fX})$.
Let 
$u _{r+1}, \dots, u _{d}$ be the image of the coordinates of $\widehat{\A} ^{s} _\fS $ in $\Gamma (\fU, \cO _{\fX})$ via $h$.

Composing with 
the log etale $\fS ^\sharp$-morphism 
$\beta _{(r)}
\colon 
\fS ^\sharp _{(r)}
\to 
\fS ^\sharp \times\mathfrak{A} _{\bbN ^{r} } $ (see \ref{S(r)chart}),
we get the log etale $\fS ^\sharp$-morphism 
$f\colon \fU ^{\sharp} \to (\widehat{\mathbb{G}} _{\mathrm{m} }) _\fS ^{s}   \times\mathfrak{A} _{\bbN ^{r} } \times _\fS \fS ^\sharp$.
Hence, this yields that $u _1,\dots, u _d$ are log coordinates of $\fU ^{\sharp} /\fS ^{\sharp} $.
Let $\partial _{\sharp 1},\dots, \partial _{\sharp d}$ be the corresponding logarithmic derivations of $\cD ^{(m)} _{\fU ^{\sharp} /\fS ^\sharp}$. 
Let $\partial _{1},\dots, \partial _{d}$ be the corresponding derivations of $\cD ^{(m)} _{\fU _\eta /\eta}$. 
For $i= 1,\dots, r$, set $\partial _{(r),i} = \partial _{\sharp i}$ and for $i= r+1,\dots, d$, set $\partial _{(r),i} = \partial _{i}$.
For any $\underline{i}= (\underline{j},\underline{k})\in \bbN ^d$ with $\underline{j}\in \bbN ^r$ and $\underline{k}\in \bbN ^s$, 
set $\underline{\partial} _{(r)} ^{<\underline{i}> _{(m)}}\coloneqq \underline{\partial} _{\sharp} ^{<\underline{j}> _{(m)}}\underline{\partial} ^{<\underline{k}> _{(m)}} $
The sheaf $\cD ^{(m)} _{\fU^\sharp/\fS^\sharp}$ is $\cO _{\fU ^\sharp}$-free with basis 
$\{ \underline{\partial} _{(r)} ^{<\underline{i}> _{(m)}},\, \underline{i} \in \bbN ^d\}$.
\end{empt}

\subsubsection{Weakly closed immersion and regularity}

\begin{lem}
\label{lem-regular}
Let $\fX ^\sharp$ be a strictly semi-stable log formal scheme over $\fS ^\sharp$.
Let $\fU ^{\sharp} $ be an affine open subset of $\fX ^\sharp$, 
$A:= \Gamma ( \fU ^\sharp, \cO _{\fX ^\sharp})$,
$x$ a point of $\fU ^{\sharp}$ and $\mathfrak{p}$ the corresponding open prime ideal
of $A$. 
Then $A$, $A _\mathfrak{p}$, its $p$-adic completion $\mathfrak{A} _\mathfrak{p}$ and 
$\cO _{\fX ^\sharp, x}$ are regular. 
\end{lem}

\begin{proof}
With \ref{smooth-flat}, recall that $\fX ^{\sharp}$ is flat over $\fS ^\sharp$.
We already know that $X ^{\sharp}$ and $\fS ^{\sharp}$ are regular. 
Since $\cV$ is regular and $\pi A$ is included in the Jacobson ideal of $A$, 
from \cite[Lemma 6.1]{MonskyWashnitzer}, we get the regularity of $A$.
Hence $A _\mathfrak{p}$ is regular and flat over $\cV$.  
Again with \cite[Lemma 6.1]{MonskyWashnitzer}, 
so is 
$\mathfrak{A} _\mathfrak{p}$. 
For the last one, we may use moreover \cite[$0 _\mathrm{I}.7.6.18$]{EGAI}
and 
\cite[17.3.3.(i)]{EGAIV1}.
\end{proof}

\begin{dfn} \label{dfn-wclimm} {\it (Weakly closed immersion)}.
\begin{enumerate}[(a)]
\item Let $u _0 \colon X ^\sharp \to P ^{\sharp}$ be a morphism of (resp. strictly semistable) log schemes over $S ^\sharp$.
We say that $u _0$ is a ``weakly closed immersion'' of (resp. strictly semistable) log schemes over $S ^\sharp$) 
if $\underline{u} _0$ is a closed immersion.

\item Let $u \colon \fX ^\sharp \to \fP ^{\sharp}$ be a morphism of (strictly semistable) log formal schemes over $\fS ^\sharp$.
We say that $u$ is a ``weakly closed immersion'' of (strictly semistable) log formal schemes over $\fS ^\sharp$)
if $\underline{u}$ is a closed immersion.

\end{enumerate}

\end{dfn}

\begin{prop}
\label{prop-regular}
Let $u \colon \fX ^\sharp \to \fP ^{\sharp}$ be a weakly closed immersion of strictly semistable log formal schemes over $\fS ^\sharp$.
We denote by $\cI$ the ideal defining $\underline{u}$.  
Then $\cI$ is regular (in the sense of \cite[IV.16.9.1]{EGAIV4}).
\end{prop}

\begin{proof}
Let $x$ be a point of $\fX$.
From \ref{lem-regular}, 
$\cO _{\fX ^\sharp, x}$ and $\cO _{\fP ^\sharp, x}$ are regular noetherian local rings.
This yields from \cite[IV.19.1.2]{EGAIV4}
that $\cI _x$ is a regular ideal 
(this notion of regular ideal is defined at \cite[IV.16.9.7]{EGAIV4}). 
From the remark of  \cite[IV.16.9.7]{EGAIV4}, 
there exists a regular $\cO _{\fP ^\sharp, x}$-sequence $(f _{i,x}) _{i=1,\dots n}$ 
generating  $\cI _x$. Since the regularity of $\cI $ is local, we can suppose that the sections
$(f _{i,x}) _{i=1,\dots n}$ come from sections $(f _{i}) _{i=1,\dots n}$ of 
$\Gamma (\fP, \cO _{\fP})$.
Since $\cI$ is coherent, using \cite[$0_\mathrm{IV}.15.2.4$]{EGAIV1}, 
there exists an open neighborhood $\fU$ of $x$ such that
$(f _{i} |\fU) _{i=1,\dots n}$ is regular $\cO _{\fU ^\sharp}$-sequence.
Since $\cI$ is coherent, shrinking $\fU$ if necessary, we can suppose that
the sequence $(f _{i} |\fU) _{i=1,\dots n}$ generates $\cI |\fU$.
\end{proof}

\subsection{Berthelot-Kashiwara's theorem}
Let $u \colon \fX ^\sharp \to \fP ^{\sharp}$ be a weakly closed immersion of strictly semistable log formal schemes over $\fS ^\sharp$ (see \ref{dfn-wclimm}).
We denote by $\cI$ the ideal defining $\underline{u}$, and $\cI _\eta: =\cI |_{\fP _\eta}$.

\subsubsection{Push-forwards of level $m$ of a weakly closed immersion : local case}

\begin{lem} \label{ovcv-coord-comp-imm}
We denote  by $\widetilde{u}\colon  (\fX, \widetilde{\cO} _\fX) \to (\fP, \widetilde{\cO} _\fP)$ the homomorphism of ringed topological spaces. 
\begin{enumerate}[(a)]

\item With Notation \ref{Dtilde}, 
we have the exact sequence of locally free of finite type $\widetilde{\cO} _{\fX }$-modules: 
\begin{equation}
\label{exactsequenceOvCv}
0 \to 
\widetilde{\cI} / \widetilde{\cI} ^2 
\to 
\widetilde{u} ^* \widetilde{\Omega} _{\fP ^\sharp /\fS ^\sharp} 
\to
\widetilde{\Omega} _{\fX ^\sharp /\fS ^\sharp}  
\to 
0
\end{equation}

\item Let $x$ be a point of $\fX$. 
There exist an affine open set $\fP'$ of $\fP$ containing $x$, 
sections $t _1,\dots, t _d\in \Gamma (\fP ' ,\cO _{\fP})$ such
that $t _{1},\dots, t _d$ are overconvergent local coordinates   of $\fP ^{\prime \sharp}/\fS ^\sharp$, 
$t _{r+1},\dots, t _d \in \Gamma (\fP' ,\cI)$ is a regular sequence of $\Gamma (\fP ' ,\cO _{\fP})$ generating $\Gamma (\fP' ,\cI)$
and 
$\overline{t} _{1},\dots, \overline{t} _r$ are 
overconvergent local coordinates of $\fX \cap \fP'/\fS ^\sharp$, where 
$ \overline{t} _{1},\dots, \overline{t} _r$ are the image 
of $t _{1},\dots, t _r$ in $\Gamma (\fX \cap \fP',\cO _{\fX})$.
\end{enumerate}

\end{lem}

\begin{proof}
1) We prove the first assertion. 
a) We construct the first map.
Consider the $\cO _{\fX}$-linear morphism which is the composition of the canonical maps:
$\cI  / \cI ^2
\to 
\underline{u} ^* \Omega _{\fP /\fS } 
\to
u ^* \Omega _{\fP ^\sharp /\fS ^\sharp}$.
By applying the functor
$\widetilde{\cO} _{\fX }  \otimes _{\cO _{\fX }} -$ to this morphism, 
we get a morphism canonically isomorphic to a morphism of the form
$\widetilde{\cI} / \widetilde{\cI} ^2 
\to 
u ^* \widetilde{\Omega} _{\fP ^\sharp /\fS ^\sharp} $.

b) Similarly, the second map is constructed from 
the map 
$u ^* \Omega _{\fP ^\sharp /\fS ^\sharp}
\to
\Omega _{\fX ^\sharp /\fS ^\sharp}$
by applying the functor
$\widetilde{\cO} _{\fX }   \otimes _{\cO _{\fX }} -$.

c) Let 
$\eta$ be the open of $\fS$ corresponding to the generic point,
$\fX _\eta$ (resp. $\fP _\eta$) the open of $\fX ^\sharp$ (resp. $\fP ^\sharp$) which is the fiber of $\eta$.
Since $v\colon \fX _\eta \hookrightarrow \fP _\eta$, the morphism induced by $u$, 
is a closed immersion of smooth formal schemes over $\eta$,
we have the exact sequence of locally free of finite type $\cO _{\fX _\eta}$-modules of the form:
 \begin{equation}
 \label{ex-seq-eta}
0 \to 
\cI _\eta  /\cI _\eta ^2 
\to 
v ^* \Omega _{\fP _\eta /\eta} 
\to
\Omega _{\fX _\eta /\eta} 
\to 
0.
\end{equation}
Since the extension $\widetilde{\cO} _{\fX }  \to j _{*} \cO _{\fX _\eta}$ is faithfully flat, 
since the restriction of the sequence \ref{exactsequenceOvCv} 
on $\fX _\eta$ is the exact sequence \ref{ex-seq-eta},
we conclude. 

2) Concerning the second assertion, let 
$\fP'$ be an affine open set  of $\fP$ containing $x$. 
Shrinking $\fP'$ is necessary, 
we can suppose that each terms of the exact sequence
\ref{exactsequenceOvCv} are $\widetilde{\cO} _{\fX }$-free over $\fX' := \fX \cap \fP'$
and that 
$\mathcal{I}  /\mathcal{I} ^2$ is $\cO _{\fX}$-free over $\fX '$ (use \ref{prop-regular}).
Let $t _{r+1},\dots, t _d \in \Gamma (\fP' ,\cI)$
be a regular sequence of $\Gamma (\fP ' ,\cO _{\fP})$ generating $\Gamma (\fP' ,\cI)$.
Shrinking $\fP'$ if necessary, we can suppose 
the classes of $t _{r+1},\dots, t _d $ induce a basis of  $\Gamma (\fP' ,\cI) / \Gamma (\fP' ,\cI)  ^2$.
Using Nakayama's Lemma, shrinking $\fP'$ if necessary, we can suppose that 
$t _{r+1},\dots, t _d \in \Gamma (\fP' ,\cI)$ generate 
$\Gamma (\fP' ,\cI)$.
From \ref{ntn-loc-coord-partial},
shrinking $\fP'$ is necessary, there exist 
overconvergent local coordinates $\overline{t} _{1},\dots, \overline{t} _r$  of $\fX ^\sharp/\fS ^\sharp$ 
with moreover $\overline{t} _{1},\dots, \overline{t} _r\in \Gamma (\fX', \cO _{\fX})$.
Choose some lifting $t _{1},\dots, t _r\in \Gamma (\fP', \cO _{\fP})$
of $\overline{t} _{1},\dots, \overline{t} _r$.
Put $\fX _\eta ':= \fX _\eta \cap \fP'$ and $\fP _\eta ':= \fP _\eta \cap \fP'$.
Since  $\overline{t} _{1}| \fX _\eta' ,\dots, \overline{t} _r | \fX _\eta' $  are local coordinates of $\fX _\eta' /\eta$,
since  $t _{r+1} |\fP _\eta',\dots, t _d |\fP _\eta' \in \Gamma (\fP _\eta', \cI)$, 
then $t _1|\fP _\eta',\dots, t _d|\fP _\eta'$  are local coordinates of $\fP _\eta'/\eta$ (consider the exact sequence \ref{ex-seq-eta}),
i.e.   $t _1,\dots, t _d$  are overconvergent local coordinates of $\fP ^{\prime \sharp}/\fS ^\sharp$.
\end{proof}

\begin{ntn}[Local context] \label{ntn-ZandZ'}
Let $x$ be a point of $\fX$. 
Let $\fP'$ be an affine open set  of $\fP$ containing $x$, 
let  $t _1,\dots, t _d\in \Gamma (\fP ' ,\cO _{\fP})$ be some sections such
that $t _{1},\dots, t _d$ are overconvergent local coordinates   of $\fP ^{\prime \sharp}/\fS ^\sharp$, 
$t _{r+1},\dots, t _d \in \Gamma (\fP' ,\cI)$ is a regular sequence of $\Gamma (\fP ' ,\cO _{\fP})$ generating $\Gamma (\fP' ,\cI)$
and 
$\overline{t} _{1},\dots, \overline{t} _r$ are 
overconvergent local coordinates of $\fX \cap \fP'/\fS ^\sharp$, where 
$ \overline{t} _{1},\dots, \overline{t} _r$ are the image 
of $t _{1},\dots, t _r$ in $\Gamma (\fX \cap \fP',\cO _{\fX})$.
Let 
$r \leq e'  < e \leq d$ be some integers. 
We denote by 
$\fZ$ the closed formal subscheme of $\fP'$ 
whose ideal defining $\fZ \hookrightarrow \fP'$ is 
that generated by 
$t _{e+1}, \dots, t _d$.
We denote by 
$\fZ '$ the closed formal subscheme of $\fP'$ 
whose ideal defining $\fZ ' \hookrightarrow \fP'$ is 
that generated by 
$t _{e'+1}, \dots, t _d$.
Since $X \hookrightarrow P$ is a closed immersion of regular schemes, 
using \cite[19.1.1]{EGAIV1} we get 
that $Z$ and $Z'$ are regular at the points $x \in X$. 
Since $S$ is excellent, since $P$ is of finite type over $S$, 
the set of regular points of $P$ is open.
Hence, 
shrinking $\fP'$ if necessary, 
we can suppose $Z$ 
and $Z '$ are regular. 
Using \cite[Lemma 6.1]{MonskyWashnitzer}, 
this yields that 
$\fZ$ 
and $\fZ '$ are regular. 
To simplify the notation, we will suppose
$\fP'=\fP$.

Let $\partial _1, \dots, \partial _d\in  \widetilde{\mathcal{T}} _{\fP ^\sharp /\fS ^\sharp} $ 
be the dual basis of 
$d t _1,\dots, d t _d\in  \widetilde{\Omega} _{\fP ^\sharp /\fS ^\sharp}$.
Set 
$\widetilde{D} ^{(m)}_{\fP ^\sharp /\fS ^\sharp} 
:=
\Gamma 
(\fP, \widetilde{\cD} ^{(m)}_{\fP ^\sharp /\fS ^\sharp} )$,
$A := \Gamma (\fP,  \cO _{\fP})$, 
$\widetilde{A} : =
 \Gamma (\fP, \widetilde{\cO} _{\fP })$,
 $B := \Gamma (\fZ,  \cO _{\fZ})$,
$\widetilde{B} : =
 \Gamma (\fZ, \widetilde{\cO} _{\fZ })$,
  $B' := \Gamma (\fZ',  \cO _{\fZ'})$, 
$\widetilde{B} ' : =
 \Gamma (\fZ', \widetilde{\cO} _{\fZ' })$.
 With notation \ref{localcoord-ovcv},
$\widetilde{D} ^{(m)}_{\fP ^\sharp /\fS ^\sharp} $ 
is $\widetilde{A}$-free with  the basis
$\{ \underline{\partial} ^{<\underline{k}> _{(m)}} \ ; \
\underline{k} \in \bbN ^{n}\}$.
We denote by 
$K$ the ideal of $A$ generated by 
$t _{e+1},\dots, t _d$.

For any 
$\underline{i}= (i _1,\dots, i _e) \in \bbN ^{e}$, 
we set 
$\underline{\partial} ^{<\underline{i}> _{(m)}} :=
\partial _1 ^{<i _1> _{(m)}} \cdots \partial _e ^{<i _e> _{(m)}}$.
For any 
$\underline{j}= (j _{1},\dots, j _{d-e}) \in \bbN ^{d-e}$, 
we set 
$\underline{\partial} ^{<\underline{j}> ^\star _{(m)}} :=
\partial _{e+1} ^{<j _{1}> _{(m)}} \cdots \partial _d ^{<j _{d-e}> _{(m)}}$.
We denote by 
$\widetilde{D} ^{(m)} _{\fP,\fZ /\fS} $ the free $\widetilde{A}$-submodule 
(for the left or the right structures)
of 
$\widetilde{D} ^{(m)}  _{\fP^\sharp/\fS^\sharp}$ 
whose basis is given by
$\{ \underline{\partial} ^{<\underline{i}> _{(m)}} \ ; \
\underline{i} \in \bbN ^{e}\}$.
In fact, $\widetilde{D} ^{(m)} _{\fP,\fZ /\fS} $ is 
a subring of 
$\widetilde{D} ^{(m)}  _{\fP^\sharp/\fS^\sharp}$
and 
$\widetilde{D} ^{(m)}  _{\fP^\sharp/\fS^\sharp}$ is a free right (or left) 
$\widetilde{D} ^{(m)} _{\fP,\fZ /\fS} $-module 
with the basis
$\{ \underline{\partial} ^{<\underline{j}> ^{\star}_{(m)}} \ ; \
\underline{j} \in \bbN ^{d-e}\}$.
We set 
$\widetilde{D} ^{(m)} _{\fZ /\fS} 
:=
\widetilde{D} ^{(m)} _{\fP,\fZ /\fS} /
\widetilde{D} ^{(m)} _{\fP,\fZ /\fS}  K$.
We denote by
$[\hspace{0,15cm}]\colon
\widetilde{D} ^{(m)}_{\fP ,\fZ /\fS } 
\to 
\widetilde{D} ^{(m)}_{\fP ,\fZ  /\fS } 
/
\widetilde{D} ^{(m)}_{\fP  ,\fZ /\fS }  K$
the canonical projection.
Since $K$ is generated by some elements which are in the center of 
$\widetilde{D} ^{(m)} _{\fP,\fZ /\fS}$,
then 
$\widetilde{D} ^{(m)} _{\fZ /\fS} $
is a ring and
the map
$[\hspace{0,15cm}]\colon
\widetilde{D} ^{(m)}_{\fP ,\fZ /\fS } 
\to 
\widetilde{D} ^{(m)} _{\fZ /\fS} $
is an epimorphism of rings.  
Moreover, 
$\widetilde{D} ^{(m)} _{\fZ /\fS} $ has two structures of free 
$\widetilde{B}$-module
with 
$\{ [\underline{\partial} ^{<\underline{i}> _{(m)}}] \ ; \
\underline{i} \in \bbN ^{e}\}$
as basis.
To simplify notation, we denote
$[\underline{\partial} ^{<\underline{i}> _{(m)}}]$ by 
$\underline{\partial} ^{<\underline{i}> _{(m)}}$.
\end{ntn}

\begin{lemm} \label{lem-wkcpdscplevelm-PZ123}
We keep notation \ref{ntn-ZandZ'}.
\begin{enumerate}[(a)]
\item Let 
$\widetilde{D} ^{(m)\dag} _{\fP,\fZ /\fS} $
(resp. 
$\widehat{D} ^{(m)} _{\fP _\eta,\fZ _\eta/\eta} $)
be the weak completion as $A$-ring
(resp. be the $p$-adic completion)
of $\widetilde{D} ^{(m)} _{\fP,\fZ /\fS} $. 
Let $z \in \widehat{D} ^{(m)} _{\fP _\eta,\fZ _\eta/\eta} $.
Then $z\in 
\widetilde{D} ^{(m)\dag} _{\fP,\fZ /\fS} $ if and only if 
$z$ can be written in the form
\begin{equation}
\label{wkcpdscplevelm-PZ1}
z = \sum _{k\in \bbN}\sum _{\underline{i}\in \bbN ^e}  \frac{a _{\underline{i},k}}{t ^{k}}\underline{\partial} ^{<\underline{i}> _{(m)}}
(resp.\, z = \sum _{k\in \bbN}\sum _{\underline{i}\in \bbN ^e}  \underline{\partial} ^{<\underline{i}> _{(m)}}\frac{a _{\underline{i},k}}{t ^{k}})
\end{equation}
where $a _{\underline{i},k} \in A$ are such that there exists a constant $c>0$ satisfying
$v _{\pi} ( a _{\underline{i},k}) \geq \frac{k+|\underline{i}|}{c} -1$
(resp. $v _{p} ( a _{\underline{i},k}) \geq \frac{k+|\underline{i}|}{c} -1$),
for any $\underline{i} \in \bbN ^e$
and $k \in \bbN$.

\item Let 
$\widetilde{D} ^{(m)\dag} _{\fZ /\fS} $
(resp. 
$\widehat{D} ^{(m)} _{\fZ _\eta/\eta} $)
be the weak completion as $B$-ring
(resp. be the $p$-adic completion)
of $\widetilde{D} ^{(m)} _{\fZ /\fS} $. 
Let $z \in \widehat{D} ^{(m)} _{\fZ _\eta/\eta} $.
Then 
$z\in 
\widetilde{D} ^{(m)\dag} _{\fZ /\fS} $ if and only if 
$z$ can be written in the form
\begin{equation}
\label{wkcpdscplevelm-PZ2}
z = \sum _{k\in \bbN}\sum _{\underline{i}\in \bbN ^e}  \frac{b _{\underline{i},k}}{t ^{k}}\underline{\partial} ^{<\underline{i}> _{(m)}}
(resp.\, z = \sum _{k\in \bbN}\sum _{\underline{i}\in \bbN ^e}  \underline{\partial} ^{<\underline{i}> _{(m)}}\frac{b _{\underline{i},k}}{t ^{k}})
\end{equation}
where $b _{\underline{i},k} \in B$ are such that there exists a constant $c>0$ satisfying
$v _{\pi} ( b _{\underline{i},k}) \geq \frac{k+|\underline{i}|}{c} -1$,
for any $\underline{i} \in \bbN ^e$
and $k \in \bbN$.

\item We have the isomorphism of rings 
$\widetilde{D} ^{(m)\dag} _{\fP,\fZ /\fS} 
/ \widetilde{D} ^{(m)\dag} _{\fP,\fZ /\fS}  K
 \riso 
 \widetilde{D} ^{(m)\dag} _{\fZ /\fS} $.
\end{enumerate}
\end{lemm}

\begin{proof}
Since the equivalence between the respective and the non respective cases follows from \ref{vp-vpi}, then we reduce to check the non respective one. 

i) In the same way as the proof of \ref{prop-wkcpdscplevelm}, 
we check that an element of the form 
\ref{wkcpdscplevelm-PZ1} and satisfying the corresponding hypotheses
is an element of $\widetilde{D} ^{(m)\dag} _{\fP,\fZ /\fS} $.
Conversely, an element $z$ of $\widetilde{D} ^{(m)\dag} _{\fP,\fZ /\fS} $
is also an element of $\widetilde{D} ^{(m)\dag} _{\fP^\sharp/\fS^\sharp} $.
Hence, it should have the description of \ref{prop-wkcpdscplevelm}.
Since the writing in the form \ref{prop-wkcpdscplevelm} is unique, 
noticing that we have only  the terms $\underline{\partial} ^{<\underline{i}> _{(m)}}$ for 
$\underline{i}\in \bbN ^e$ (i.e. the term of the form
$\underline{\partial} ^{<\underline{i}> _{(m)}}$ with $\underline{i}\in \bbN ^d$
should be such that 
$i _{e+1}=0,\dots, i _d= 0$), then $z$ should satisfy the hypothesis of \ref{wkcpdscplevelm-PZ1}.

ii) We have the equality 
$\widetilde{D} ^{(m)\dag} _{\fP,\fZ /\fS} 
\cap 
\widehat{D} ^{(m)} _{\fP _\eta,\fZ _\eta/\eta} K
=
\widetilde{D} ^{(m)\dag} _{\fP,\fZ /\fS}  K$.
Indeed,
since the image of $t _i$ is not a divisor of zero of 
$\overline{A}:= A / \pi A$, we check that for any $a \in A$, 
$v _{\pi} (a) = v _{\pi} (a t _i)$.
Hence, from the part i) of the proof we conclude.

iii) 
We denote by
$[\hspace{0,15cm}]\colon
\widetilde{D} ^{(m)\dag}_{\fP ,\fZ /\fS } 
\to 
\widetilde{D} ^{(m)\dag}_{\fP ,\fZ  /\fS } 
/
\widetilde{D} ^{(m)\dag}_{\fP  ,\fZ /\fS }  K$
the canonical projection which is an epimorphism of rings (because $K$ is generated by elements in the center of 
$\widetilde{D} ^{(m)\dag}_{\fP ,\fZ /\fS } $).
Following the step ii), 
the map 
$\widetilde{D} ^{(m)\dag}_{\fP ,\fZ  /\fS } 
/
\widetilde{D} ^{(m)\dag}_{\fP  ,\fZ /\fS }  K
\to 
\widehat{D} ^{(m)} _{\fP _\eta,\fZ _\eta/\eta} 
/
\widehat{D} ^{(m)} _{\fP _\eta,\fZ _\eta/\eta} K
\riso
\widehat{D} ^{(m)} _{\fZ _\eta/\eta} $
is a monomorphism of rings.
Let $z\in 
\widetilde{D} ^{(m)\dag} _{\fP,\fZ /\fS} $.
Hence, we can view a priori 
$[z]$ as an element of $\widehat{D} ^{(m)} _{\fZ _\eta/\eta} $.
From the step i), 
we can write 
$z = \sum _{k\in \bbN}\sum _{\underline{i}\in \bbN ^e}  \frac{a _{\underline{i},k}}{t ^{k}}\underline{\partial} ^{<\underline{i}> _{(m)}}$, 
where $a _{\underline{i},k} \in A$ are such that there exists a constant $c>0$ satisfying
$v _{\pi} ( a _{\underline{i},k}) \geq \frac{k+|\underline{i}|}{c} -1$,
for any $\underline{i} \in \bbN ^e$
and $k \in \bbN$.
We get in $\widehat{D} ^{(m)} _{\fZ _\eta/\eta} $
the equality 
$[z] = \sum _{k\in \bbN}\sum _{\underline{i}\in \bbN ^e}  \frac{[a _{\underline{i},k}]}{t ^{k}}\underline{\partial} ^{<\underline{i}> _{(m)}}$
(recall $\underline{\partial} ^{<\underline{i}> _{(m)}}:= [\underline{\partial} ^{<\underline{i}> _{(m)}}]$).
Since $v _{\pi} ([a _{\underline{i},k}]) \geq v _{\pi} (a _{\underline{i},k})$, 
in the same way as the proof of \ref{prop-wkcpdscplevelm} this yields that 
$[z] \in  \widetilde{D} ^{(m)\dag} _{\fZ /\fS} $.

iv) Let $z \in \widehat{D} ^{(m)} _{\fZ _\eta/\eta} $.
In the same way as the proof of \ref{prop-wkcpdscplevelm}, if $z$ can be written in the form
\ref{wkcpdscplevelm-PZ2} with the required properties then
$z\in 
\widetilde{D} ^{(m)\dag} _{\fZ /\fS} $.

v) Let $y\in 
\widetilde{D} ^{(m)\dag} _{\fZ /\fS} $.
Then, we can write by definition (see \ref{dfn-wc=pre})
\begin{equation}
\label{dfn-wc=pre-ZP}
y= \sum _{i=0} ^{\infty} \pi ^{i} Q _i (y _1,\dots, y _r)
\end{equation}
 where $y _1,\dots, y _r \in \widetilde{D} ^{(m)} _{\fZ /\fS} \setminus B$, 
 $Q _i \in  \mathfrak{V} (y _1,\dots, y _r/B)$ are such that there exists a constant $c >0$
 so that 
 $\deg (Q _i) \leq c (i+1)$ for any $i$.
Let $x _1, \dots, x _r \in \widetilde{D} ^{(m)}_{\fP ,\fZ /\fS } $
such that 
$[x _1] = y _1,\dots  ,[x _r] = y _r$.
Since 
$[\hspace{0,15cm}]\colon
\widetilde{D} ^{(m)}_{\fP ,\fZ /\fS } 
\to 
\widetilde{D} ^{(m)} _{\fZ /\fS} $
and $A \to B$ 
are epimorphisms of rings, 
then there exists 
 $P _i \in  \mathfrak{V} (x _1,\dots, x _r/A)$
 such that 
 $\deg P _i = \deg Q _i$
 and 
 $[P _i (x _1,\dots, x _r)] 
 =
 Q _i (y _1,\dots, y _r)$
 for any $i$. 
 Set $x :=\sum _{i=0} ^{\infty} \pi ^{i} P _i (x _1,\dots, x _r) \in \widetilde{D} ^{(m)\dag} _{\fP,\fZ /\fS} $.
 Then we compute in 
 $ \widetilde{D} ^{(m)\dag} _{\fZ /\fS} $ the equality
 $[x]= y$.
 
 vi) With the steps iii), iv) and v) we get the assertions 2 and 3 of our lemma.
 \end{proof}

\begin{ntn}
\label{ntn-ZandZ'infty}
We keep notation \ref{ntn-ZandZ'}.
For any 
$\underline{i}= (i _1,\dots, i _e) \in \bbN ^{e}$, 
we set 
$\underline{\partial} ^{[\underline{i}]} :=
\partial _1 ^{[i _1]} \cdots \partial _e ^{[i _e]}$.
For any 
$\underline{j}= (j _{1},\dots, j _{d-e}) \in \bbN ^{d-e}$, 
we set 
$\underline{\partial} ^{[\underline{j}] ^\star } :=
\partial _{e+1} ^{[j _{1}] } \cdots \partial _d ^{[j _{d-e}] }$.
We denote by 
$\widetilde{D}  _{\fP,\fZ /\fS} $ the free $\widetilde{A}$-submodule 
(for the left or the right structures)
of 
$\widetilde{D}   _{\fP^\sharp/\fS^\sharp}$ 
whose basis is given by
$\{ \underline{\partial} ^{[\underline{i}] } \ ; \
\underline{i} \in \bbN ^{e}\}$.
In fact, $\widetilde{D}  _{\fP,\fZ /\fS} $ is 
a subring of 
$\widetilde{D}   _{\fP^\sharp/\fS^\sharp}$
and 
$\widetilde{D}   _{\fP^\sharp/\fS^\sharp}$ is a free right (or left) 
$\widetilde{D}  _{\fP,\fZ /\fS} $-module 
which has 
$\{ \underline{\partial} ^{[\underline{j}] ^{\star}} \ ; \
\underline{j} \in \bbN ^{d-e}\}$ as basis.
We set 
$\widetilde{D}  _{\fZ /\fS} 
:=
\widetilde{D}  _{\fP,\fZ /\fS} /
\widetilde{D}  _{\fP,\fZ /\fS}  K$.
We denote by
$[\hspace{0,15cm}]\colon
\widetilde{D} _{\fP ,\fZ /\fS } 
\to 
\widetilde{D} _{\fP ,\fZ  /\fS } 
/
\widetilde{D} _{\fP  ,\fZ /\fS }  K$
the canonical projection.
Since $K$ is generated by some elements which are in the center of 
$\widetilde{D}  _{\fP,\fZ /\fS}$,
then 
$\widetilde{D}  _{\fZ /\fS} $
is a ring and
the map
$[\hspace{0,15cm}]\colon
\widetilde{D} _{\fP ,\fZ /\fS } 
\to 
\widetilde{D}  _{\fZ /\fS} $
is an epimorphism of rings.  
Moreover, 
$\widetilde{D}  _{\fZ /\fS} $ has two structures of free 
$\widetilde{B}$-module
with 
$\{ [\underline{\partial} ^{[\underline{i}] }] \ ; \
\underline{i} \in \bbN ^{e}\}$
as basis.
To simplify notation, we simply denote
$[\underline{\partial} ^{[\underline{i}] }]$ by 
$\underline{\partial} ^{[\underline{i}] }$.

\end{ntn}

\begin{prop}
\label{lem-wkcpdscplevelm-PZ123infty}
We keep notation \ref{ntn-ZandZ'infty}.
 Let 
$\widetilde{D} ^{\dag} _{\fZ /\fS} $
(resp. 
$\widehat{D} _{\fZ _\eta/\eta} $)
be the weak completion as $B$-ring
(resp. be the $p$-adic completion)
of $\widetilde{D}  _{\fZ /\fS} $. 
Let $z \in \widehat{D}  _{\fZ _\eta/\eta} $.
Then 
$z\in 
\widetilde{D} ^{\dag} _{\fZ /\fS} $ if and only if 
$z$ can be written in the form
\begin{equation} \label{wkcpdscpinfty}
z = \sum _{k\in \bbN}\sum _{\underline{i}\in \bbN ^e}  \frac{b _{\underline{i},k}}{t ^{k}}\underline{\partial} ^{[\underline{i}]}
(resp.\, z = \sum _{k\in \bbN}\sum _{\underline{i}\in \bbN ^e}  \underline{\partial} ^{[\underline{i}]}\frac{b _{\underline{i},k}}{t ^{k}})
\end{equation}
where $b _{\underline{i},k} \in B$ are such that there exists a constant $c>0$ satisfying $v _{\pi} ( b _{\underline{i},k}) \geq \frac{k+|\underline{i}|}{c} -1$ (resp. $v _{p} ( b _{\underline{i},k}) \geq \frac{k+|\underline{i}|}{c} -1$), for any $\underline{i} \in \bbN ^e$ and $k \in \bbN$.
\end{prop}

\begin{proof}
Using \ref{wkcpdscplevelm-PZ2}, we can proceed as for the proof of \ref{wkcpdscp}.
\end{proof}

\begin{lem}
\label{empt-wkcpdscptheta}
We keep notation \ref{ntn-ZandZ'infty}.
Let $\theta \in B\cap \widetilde{B} ^*$ and let $z$ be an element of $\widehat{D}  _{\fZ _\eta/\eta}$.
Suppose we can write $z$ in the form
\begin{equation}
\label{wkcpdscptheta}
z = \sum _{k,l\in \bbN}\sum _{\underline{i}\in \bbN ^d} b _{\underline{i},k,l} \frac{1}{t ^{k} \theta ^l}\underline{\partial} ^{[\underline{i}]}
\end{equation}
where $b _{\underline{i},k,l} \in B$ are such that there exists a constant $c>0$ such that 
$v _{\pi} ( b _{\underline{i},k,l}) \geq \frac{k+l+|\underline{i}|}{c} -1$.
Then, 
$z \in \widetilde{D} ^{\dag} _{\fZ /\fS} $.
\end{lem}

\begin{proof}
For $m$ large enough (i.e. satisfying \ref{formula1} and \ref{formula2}), following the part II)1)  of the proof of \ref{loc-m2dagbis}, 
we can write $z$ of the form  $z = \sum _{k,l\in \bbN}\sum _{\underline{i}\in \bbN ^d} b _{\underline{i},k,l} \frac{1}{t ^{k} \theta ^l}\underline{\partial} ^{<\underline{i}> _{(m)}}$, $b _{\underline{i},k,l} \in B$ are such that
$v _{\pi} ( b _{\underline{i},k,l}) \geq \frac{k+l+|\underline{i}|}{2c} -1$. Since from \cite[2.2.5.1]{Be1},  $\underline{\partial} ^{<\underline{i}> _{(m)}}$ is a product of elements in 
$\{\partial _i ^{[p ^j]}\ |\ i=1,\dots, d\; ;\; j= 0,\dots, m\}$ (remark also that the number of factors is less than $|\underline{i}|$), then  $z \in  (\{ \frac{1}{\theta} \} \cup \{ \frac{1}{t} \}  \cup \{\partial _i ^{[p ^j]}\ |\ i=1,\dots, d\; ;\; j= 0,\dots, m\})^\dag 
\subset  \widetilde{D} ^{(m)\dag} _{\fZ /\fS}  \subset  \widetilde{D} ^{\dag} _{\fZ /\fS} $,  where the symbol $\dag$ means the weak completion as $B$-ring.
\end{proof}

\begin{ntn}
\label{ntn-DZZ'}
For any 
$\underline{i}'= (i '_1,\dots, i '_{e'}) \in \bbN ^{e'}$, 
we set 
$\underline{\partial} ^{<\underline{i}'> _{(m)}} :=
\partial _1 ^{<i '_1> _{(m)}} \cdots \partial _{e'} ^{<i '_{e'}> _{(m)}}$.
For any 
$\underline{j}'= (j '_{1},\dots, j '_{e-e'}) \in \bbN ^{e-e'}$, 
we set 
$\underline{\partial} ^{<\underline{j}'> ^\star _{(m)}} :=
\partial _{e'+1} ^{<j '_{1}> _{(m)}} \cdots \partial _e ^{<j '_{e-e'}> _{(m)}}$.
We denote by 
$\widetilde{D} ^{(m)} _{\fZ,\fZ ^{\prime}/\fS} $ the free $\widetilde{B}$-submodule 
(for the left or the right structures)
of 
$\widetilde{D} ^{(m)}  _{\fZ/\fS}$ 
whose basis is given by
$\{ \underline{\partial} ^{<\underline{i}'> _{(m)}} \ ; \
\underline{i}' \in \bbN ^{e'}\}$.
In fact, $\widetilde{D} ^{(m)} _{\fZ,\fZ ^{\prime}/\fS} $ is 
a subring of 
$\widetilde{D} ^{(m)}  _{\fZ/\fS}$
and 
$\widetilde{D} ^{(m)}  _{\fZ /\fS}$ is a free right (or left) 
$\widetilde{D} ^{(m)} _{\fZ,\fZ ^{\prime}/\fS} $-module 
with the basis
$\{ \underline{\partial} ^{<\underline{j}'> ^{\star}_{(m)}} \ ; \
\underline{j} '\in \bbN ^{e-e'}\}$.

We denote by 
$K'$ the ideal of $B$ generated by 
$t _{e'+1},\dots, t _e$.
We get the isomorphism of $B'$-rings
$\widetilde{D} ^{(m)} _{\fZ ^{\prime}/\fS} 
\riso
\widetilde{D} ^{(m)} _{\fZ,\fZ ^{\prime}/\fS } /
\widetilde{D} ^{(m)} _{\fZ,\fZ ^{\prime}/\fS }  K'$.
We set
$\widetilde{D} ^{(m)} _{\fZ\leftarrow\fZ ^{\prime} /\fS}
:=
\widetilde{D} ^{(m)} _{\fZ/\fS} /
\widetilde{D} ^{(m)} _{\fZ/\fS}  K'$,
$\widetilde{D} ^{(m)\dag} _{\fZ\leftarrow\fZ ^{\prime} /\fS}
:=
\widetilde{D} ^{(m)\dag} _{\fZ/\fS} /
\widetilde{D} ^{(m)\dag} _{\fZ/\fS}  K'$,
$\widehat{D} ^{(m)} _{\fZ _\eta\leftarrow\fZ ^{\prime} _\eta/\eta}
:=
\widehat{D} ^{(m)} _{\fZ _\eta/\eta} /
\widehat{D} ^{(m)} _{\fZ _\eta/\eta}  K'$.
We denote by
\begin{equation}\label{ntn-DZZ'[]} [\hspace{0,15cm}]\colon \widetilde{D} ^{(m)} _{\fZ /\fS} \to  \widetilde{D} ^{(m)} _{\fZ\leftarrow\fZ ^{\prime} /\fS},
\text{ and }[\hspace{0,15cm}]\colon \widetilde{D} ^{(m)} _{\fZ,\fZ ^{\prime}/\fS}  \to  \widetilde{D} ^{(m)} _{\fZ ^{\prime}/\fS } 
\end{equation}
the canonical projections (and similarly with ``$\widetilde{D} ^{(m)} $'' replaced by 
``$\widetilde{D} ^{(m)\dag}$'').

The object $\widehat{D} ^{(m)} _{\fZ _\eta\leftarrow\fZ ^{\prime} _\eta/\eta}$
is a 
$(\widehat{D} ^{(m)} _{\fZ _\eta/\eta}, \widehat{D} ^{(m)} _{\fZ ^{\prime} _\eta/\eta})$-bimodule
and 
$\widetilde{D} ^{(m)} _{\fZ\leftarrow\fZ ^{\prime} /\fS}$
is a sub
$(\widetilde{D} ^{(m)} _{\fZ/\fS},\widetilde{D} ^{(m)} _{\fZ ^{\prime} /\fS})$-module of 
$\widehat{D} ^{(m)} _{\fZ _\eta\leftarrow\fZ ^{\prime} _\eta/\eta}$.
More precisely, via the isomorphism $\widetilde{D} ^{(m)} _{\fZ ^{\prime}/\fS} 
\riso
\widetilde{D} ^{(m)} _{\fZ,\fZ ^{\prime}/\fS } /
\widetilde{D} ^{(m)} _{\fZ,\fZ ^{\prime}/\fS }  K'$,
the structure of 
right 
$\widetilde{D} ^{(m)} _{\fZ ^{\prime}/\fS}  $-module 
on 
$\widetilde{D} ^{(m)} _{\fZ\leftarrow\fZ ^{\prime} /\fS}$
is given by  the product
$[P] [Q] := [PQ]$, 
for any 
$P \in \widetilde{D} ^{(m)} _{\fZ/\fS }$, 
$Q \in \widetilde{D} ^{(m)} _{\fZ,\fZ ^{\prime}/\fS }$
(which is well defined because $K'$ is generated by 
some elements belonging to the center of 
$\widetilde{D} ^{(m)} _{\fZ,\fZ ^{\prime}/\fS }$). 
Moreover, we remark that 
$\widetilde{D} ^{(m)} _{\fZ\leftarrow\fZ ^{\prime} /\fS}$
is a free right 
$\widetilde{D} ^{(m)} _{\fZ ^{\prime}/\fS}  $-module 
which has 
$\{ [\underline{\partial} ^{<\underline{j}> ^\star _{(m)}}] \ ; \
\underline{j} \in \bbN ^{e-e'}\}$ as basis.
\end{ntn}

\begin{prop}
\label{lem-wkcpdscplevelm-PZ123bis}
We keep notation \ref{ntn-DZZ'}.
\begin{enumerate}[(a)]
\item 
Let 
$\widetilde{D} ^{(m)\dag} _{\fZ,\fZ ^{\prime}/\fS} $
(resp. $\widehat{D} ^{(m)} _{\fZ _\eta,\fZ ^{\prime} _\eta /\eta} $)
be the weak completion as $B$-ring (resp. the $p$-adic completion)
of 
$\widetilde{D} ^{(m)} _{\fZ,\fZ ^{\prime}/\fS} $. 
Let $z \in \widehat{D} ^{(m)} _{\fZ _\eta,\fZ ^{\prime} _\eta /\eta} $.
Then $z\in 
\widetilde{D} ^{(m)} _{\fZ,\fZ ^{\prime}/\fS}  $ if and only if 
$z$ can be written in the form
\begin{equation}
\label{wkcpdscplevelm-PZ1bis}
z = \sum _{k\in \bbN}\sum _{\underline{i}\in \bbN ^{e'}}  \frac{b _{\underline{i},k}}{t ^{k}}\underline{\partial} ^{<\underline{i}> _{(m)}}
(resp.\, z = \sum _{k\in \bbN}\sum _{\underline{i}\in \bbN ^{e'}}  \underline{\partial} ^{<\underline{i}> _{(m)}}\frac{b _{\underline{i},k}}{t ^{k}})
\end{equation}
where $b _{\underline{i},k} \in B$ are such that there exists a constant $c>0$ satisfying
$v _{\pi} ( b _{\underline{i},k}) \geq \frac{k+|\underline{i}|}{c} -1$ (resp. $v _{p} ( b _{\underline{i},k}) \geq \frac{k+|\underline{i}|}{c} -1$),
for any $\underline{i} \in \bbN ^{e'}$
and $k \in \bbN$.

\item \label{lem-wkcpdscplevelm-PZ123bis(b)} Moreover, we get the isomorphisms of rings
$\widetilde{D} ^{(m)\dag} _{\fZ,\fZ ^{\prime}/\fS} 
/
\widetilde{D} ^{(m)\dag} _{\fZ,\fZ ^{\prime}/\fS}  K'
\riso
\widetilde{D} ^{(m)\dag} _{\fZ ^{\prime}/\fS} $.
\end{enumerate}

\end{prop}

\begin{proof}
Using \ref{wkcpdscplevelm-PZ2}, 
we check 
\ref{wkcpdscplevelm-PZ1bis} (in the same way as the fact that 
\ref{wkcpdscplevelm-PZ1} is a consequence of \ref{prop-wkcpdscplevelm}).
We obtain the second statement from the first one, as \ref{wkcpdscplevelm-PZ2} was a consequence
of \ref{wkcpdscplevelm-PZ1}. 

\end{proof}

\begin{lem}
\label{lem-DZZ'subbimod}
The object $\widetilde{D} ^{(m)\dag} _{\fZ\leftarrow\fZ ^{\prime} /\fS}$
is a sub
$(\widetilde{D} ^{(m)\dag} _{\fZ/\fS},\widetilde{D} ^{(m)\dag} _{\fZ ^{\prime} /\fS})$-bimodule of 
$\widehat{D} ^{(m)} _{\fZ _\eta\leftarrow\fZ ^{\prime} _\eta/\eta}$.

\end{lem}

\begin{proof}
The object 
$\widetilde{D} ^{(m)\dag} _{\fZ/\fS} $
is endowed with a canonical 
$(\widetilde{D} ^{(m)\dag} _{\fZ/\fS} , \widetilde{D} ^{(m)\dag} _{\fZ,\fZ ^{\prime}/\fS} )$-bimodule structure. 
Since 
we have 
$\widetilde{D} ^{(m)\dag} _{\fZ,\fZ ^{\prime}/\fS} 
/
\widetilde{D} ^{(m)\dag} _{\fZ,\fZ ^{\prime}/\fS}  K'
\riso
\widetilde{D} ^{(m)\dag} _{\fZ ^{\prime}/\fS} $ (see \ref{lem-wkcpdscplevelm-PZ123bis}.2),
we get  a structure 
of 
$(\widetilde{D} ^{(m)\dag} _{\fZ/\fS},\widetilde{D} ^{(m)\dag} _{\fZ ^{\prime} /\fS})$-bimodule
on 
$\widetilde{D} ^{(m)\dag} _{\fZ\leftarrow\fZ ^{\prime} /\fS}:=
\widetilde{D} ^{(m)\dag} _{\fZ/\fS} /\widetilde{D} ^{(m)\dag} _{\fZ/\fS} K'$.
More precisely, the structure of 
right 
$\widetilde{D} ^{(m)\dag} _{\fZ ^{\prime}/\fS}  $-module 
on 
$\widetilde{D} ^{(m)\dag} _{\fZ\leftarrow\fZ ^{\prime} /\fS}$
is given by  the product
$[P] [Q] := [PQ]$, 
for any 
$P \in \widetilde{D} ^{(m)\dag} _{\fZ/\fS }$, 
$Q \in \widetilde{D} ^{(m)\dag} _{\fZ,\fZ ^{\prime}/\fS }$
(which is well defined because $K'$ is generated by 
some elements belonging to the center of 
$\widetilde{D} ^{(m)\dag} _{\fZ,\fZ ^{\prime}/\fS }$). 
In the same way as the step ii) of the proof of 
\ref{lem-wkcpdscplevelm-PZ123} (use also the description \ref{wkcpdscplevelm-PZ2}),
we get the equality 
$\widetilde{D} ^{(m)\dag} _{\fZ/\fS} 
\cap 
\widehat{D} ^{(m)} _{\fZ _\eta/\eta} K '
=
\widetilde{D} ^{(m)\dag} _{\fZ/\fS}  K'$.
This yields the injection 
$\widetilde{D} ^{(m)\dag} _{\fZ\leftarrow\fZ ^{\prime} /\fS}
\hookrightarrow
\widehat{D} ^{(m)} _{\fZ _\eta\leftarrow\fZ ^{\prime} _\eta/\eta}$.
Since the structure of 
$(\widehat{D} ^{(m)} _{\fZ _\eta/\eta},\widehat{D} ^{(m)} _{\fZ _\eta ^{\prime} /\eta})$-bimodule
on 
$\widehat{D} ^{(m)} _{\fZ _\eta\leftarrow\fZ _\eta ^{\prime} /\eta}=
\widehat{D} ^{(m)} _{\fZ _\eta/\eta} /\widehat{D} ^{(m)} _{\fZ _\eta/\eta} K'$
is defined similarly, we can conclude.

\end{proof}

\begin{rem}
Let $z\in 
\widetilde{D} ^{(m)\dag} _{\fZ /\fS} $. Then 
$z$ can be written in the form
\begin{equation}
\label{wkcpdscplevelm-PZ2bis}
z = \sum _{k\in \bbN}\sum _{\underline{i}\in \bbN ^e}  \underline{\partial} ^{<\underline{i}> _{(m)}}\frac{b _{\underline{i},k}}{t ^{k}}
\end{equation}
where $b _{\underline{i},k} \in B$ are such that there exists a constant $c>0$ satisfying
$v _{\pi} ( b _{\underline{i},k}) \geq \frac{k+|\underline{i}|}{c} -1$ (resp. $v _{p} ( b _{\underline{i},k}) \geq \frac{k+|\underline{i}|}{c} -1$), for any $\underline{i} \in \bbN ^e$
and $k \in \bbN$.
We get in 
$\widehat{D} ^{(m)} _{\fZ _\eta\leftarrow\fZ ^{\prime} _\eta/\eta}$ the formula: 
\begin{gather}
\notag
[z] 
=
\sum _{\underline{i} ''\in \bbN ^{e-e'}}
[\underline{\partial} ^{<\underline{i}''> ^\star_{(m)}}]
\left [ \sum _{k\in \bbN}\sum _{\underline{i} '\in \bbN ^{e'}} 
\underline{\partial} ^{<\underline{i} '> _{(m)}}
b _{(\underline{i} ', \underline{i} ''),k} \frac{1}{t ^{k}}
\right ]
\\
\label{wkcpdscplevelm-app2}
=
\sum _{\underline{i} ''\in \bbN ^{e-e'}}
[\underline{\partial} ^{<\underline{i}''> ^\star_{(m)}}]
\cdot
\left (
 \sum _{k\in \bbN}\sum _{\underline{i} '\in \bbN ^{e'}} 
\underline{\partial} ^{<\underline{i} '> _{(m)}}
[b _{(\underline{i} ', \underline{i} ''),k} ]\frac{1}{t ^{k}}\right )
=
\sum _{\underline{i} ''\in \bbN ^{e-e'}}
[\underline{\partial} ^{<\underline{i}''> ^\star_{(m)}}]
\cdot
P _{\underline{i}''},
\end{gather}
where $\cdot$ in the last equality means the action of 
$ P _{\underline{i}''}:=  \sum _{k\in \bbN}\sum _{\underline{i} '\in \bbN ^{e'}} 
\underline{\partial} ^{<\underline{i} '> _{(m)}}
[b _{(\underline{i} ', \underline{i} ''),k} ]\frac{1}{t ^{k}}
\in \widetilde{D} ^{(m)\dag} _{\fZ ^{\prime} /\fS}$
on 
the element $[\underline{\partial} ^{<\underline{i}''> ^\star_{(m)}}] $ of
$\widetilde{D} ^{(m)\dag} _{\fZ\leftarrow\fZ ^{\prime} /\fS}
\subset \widehat{D} ^{(m)} _{\fZ _\eta\leftarrow\fZ ^{\prime} _\eta/\eta}$.
The writing of $z$ in the form
\ref{wkcpdscplevelm-app2} is unique.
\end{rem}

\begin{lem}
\label{u+coh=cohetc-lem}
Let $\iota  \colon \fZ ^{\prime} \to \fZ $ be the induced morphism. 
Let $E$ be a  left $\widetilde{D} ^{(m)\dag} _{\fZ '/\fS }$-module of finite type, 
$\widehat{E}:=\widehat{D} ^{(m)} _{\fZ '_\eta/\eta} \otimes _{\widetilde{D} ^{(m)\dag} _{\fZ '/\fS }} E$.
The push forward by $\iota$ of level $m$ of $E$
(the push forward by $\iota _\eta$ of level $m$ of $\widehat{E}$) is defined by setting
\begin{gather}
\notag
\iota ^{(m)\dag} _{+} (E)
:= 
\widetilde{D} ^{(m)\dag}  _{\fZ  \leftarrow \fZ ^{\prime}/\fS} \otimes _{\widetilde{D} ^{(m)\dag} _{\fZ '/\fS }}E,
\, 
\notag
\iota ^{(m)} _{\eta +} (\widehat{E} )
:= 
\widehat{D} ^{(m)}  _{\fZ _\eta \leftarrow \fZ ^{\prime} _\eta/\eta}
 \otimes _{\widehat{D} ^{(m)} _{\fZ ^{\prime} _\eta/\eta}}
 \widehat{E}.
\end{gather}

\begin{enumerate}[(a)]

\item The $\widetilde{D} ^{(m)\dag} _{\fZ  /\fS }$-module $\iota ^{(m)\dag} _{+} (E)$ is of finite type
(resp. the $\widehat{D} ^{(m)} _{\fZ  _\eta/\eta}$-module 
$\iota ^{(m)} _{\eta +} (\widehat{E} ) $ is of finite type). 

\item  We have the isomorphism
$$\iota ^{(m)} _{\eta +} (\widehat{E} )
\riso 
\widehat{D} ^{(m)} _{\fZ _\eta/\eta} \otimes _{\widetilde{D} ^{(m)\dag} _{\fZ /\fS }}
\iota ^{(m)\dag} _{+} (E).$$

\item For any 
left $\widetilde{D} ^{(m)\dag} _{\fZ  /\fS }$-module $G$, 
the module 
$\mathcal{H} ^0 \iota ^! (G):= 
\mathrm{Hom} _{\widetilde{B} }
(\widetilde{B}  /K ', G)$
is endowed with a structure of 
$\widetilde{D} ^{(m)\dag} _{\fZ  ^{\prime}/\fS }$-module.

\item \label{noptors-iso}
Suppose $E$ has no $p$-torsion. Then,
the canonical homomorphism 
$E \to \mathcal{H} ^0 \iota ^! \iota ^{(m)\dag} _{+} (E)$
is an isomorphism.

\item 
\begin{enumerate}[(i)]
\item The rings $\widetilde{D} ^{(m)\dag} _{\fZ  /\fS }$ and $D^{(m)\dag} _{\fZ _\eta/\eta} $ are noetherian and $p\widetilde{D} ^{(m)\dag} _{\fZ  /\fS } \subset J (\widetilde{D} ^{(m)\dag} _{\fZ  /\fS })$,
$p\widetilde{D} ^{(m)\dag} _{\fZ _\eta/\eta} \subset J (\widetilde{D} ^{(m)\dag} _{\fZ _\eta/\eta})$.
\item Any left or right  $\widetilde{D} ^{(m)\dag} _{\fZ /\fS} $-module (resp. $\widetilde{D} ^{(m)\dag} _{\fZ _\eta/\eta} $-module) $M$ of finite type is separated and any submodule  of such $M$ is a closed subset. 
\item \label{Jacobson-ZZ'}
The homomorphisms of rings $\widetilde{D} ^{(m)\dag} _{\fZ /\fS}  \to  D^{(m)\dag} _{\fZ _\eta/\eta}  \to  \widehat{D} ^{(m)} _{\fZ _\eta/\eta} $ are left and right faithfully flat. 
\end{enumerate}

\end{enumerate}

\end{lem}

\begin{proof}
1) Since $\widetilde{D} ^{(m)\dag}  _{\fZ  \leftarrow \fZ ^{\prime}/\fS}$ is  a left $\widetilde{D} ^{(m)\dag}  _{\fZ /\fS}$-module of finite type,
since  $\widehat{D} ^{(m)}  _{\fZ _\eta \leftarrow \fZ ^{\prime} _\eta/\eta}$ is a left  $\widehat{D} ^{(m)} _{\fZ  _\eta/\eta}$-module of finite type,  the first statement is obvious.
Since  $\widehat{D} ^{(m)}  _{\fZ _\eta \leftarrow \fZ ^{\prime} _\eta/\eta} \riso \widehat{D} ^{(m)} _{\fZ _\eta/\eta} \otimes _{\widetilde{D} ^{(m)\dag} _{\fZ /\fS }} \widetilde{D} ^{(m)\dag}  _{\fZ  \leftarrow \fZ ^{\prime}/\fS}$,
then we check the second assertion by using elementary properties of the tensor product. Let $G$ be a left $\widetilde{D} ^{(m)\dag} _{\fZ  /\fS }$-module.
This induces by restriction a structure of   left $\widetilde{D} ^{(m)\dag} _{\fZ,\fZ ^{\prime}/\fS}  $-module on $G$.
Since $K'$ is generated by $t _{e'+1},\dots, t _e$, we get $\mathcal{H} ^0 \iota ^! (G) = \cap _{\alpha=e' +1} ^{e} \ker ( G \overset{t _\alpha}{\longrightarrow} G)$.
Since $t _{e'+1},\dots, t _e$ belong to the center of $\widetilde{D} ^{(m)\dag} _{\fZ,\fZ ^{\prime}/\fS}  $, then  $\mathcal{H} ^0 \iota ^! (G)$ is a sub-$\widetilde{D} ^{(m)\dag} _{\fZ,\fZ ^{\prime}/\fS} $-module of $G$.
Since $\mathcal{H} ^0 \iota ^! (G)$ is vanished by $\widetilde{D} ^{(m)\dag} _{\fZ,\fZ ^{\prime}/\fS} K'$, we obtain the third statement using  \ref{lem-wkcpdscplevelm-PZ123bis}.(\ref{lem-wkcpdscplevelm-PZ123bis(b)}).

2) Let us now check the  part (d) when $e= d$ and $E$ is moreover supposed to be separated for the $p$-adic topology.

i) Since  $\widehat{D} ^{(m)} _{\fZ _\eta/\eta}$ is noetherian  (because $\fZ _\eta/\eta$ is smooth, see \cite[3]{Be1}), since $\widehat{D} ^{(m)} _{\fZ _\eta/\eta} \otimes _{\widetilde{D} ^{(m)\dag} _{\fZ /\fS }} E$
is a  $\widehat{D} ^{(m)} _{\fZ _\eta/\eta}$-module of finite type,  then  $\widehat{D} ^{(m)} _{\fZ _\eta/\eta} \otimes _{\widetilde{D} ^{(m)\dag} _{\fZ /\fS }} E$ is separated and complete for the $p$-adic topology (see \cite[3.2.3.(v)]{Be1}).
Since  $\widehat{D} ^{(m)} _{\fZ _\eta/\eta}$ is the $p$-adic separated completion of  $\widetilde{D} ^{(m)\dag} _{\fZ /\fS }$, this yields that  $\widehat{D} ^{(m)} _{\fZ _\eta/\eta} \otimes _{\widetilde{D} ^{(m)\dag} _{\fZ /\fS }} E$
is the $p$-adic separated completion of $E$,
which justifies the notation $\widehat{E}:=\widehat{D} ^{(m)} _{\fZ _\eta/\eta} \otimes _{\widetilde{D} ^{(m)\dag} _{\fZ /\fS }} E$. Since $E$ is separated for the $p$-adic topology,
then the morphism  $E \to \widehat{E}$ is injective. Let $\widehat{E} _{\mathrm{ntor}}$ be the quotient of  $\widehat{E}$ by its $p$-torsion elements. 
Then $\widehat{E} _{\mathrm{ntor}}$ is a $\widehat{D} ^{(m)} _{\fZ _\eta/\eta}$-module of finite type  (use the same arguments as those in \cite[3.4.4]{Be1}). Since $E$ has no $p$-torsion, the morphism $E \to \widehat{E} _{\mathrm{ntor}}$ is injective.
Let  $\phi \colon \iota ^{(m)\dag} _{+} (E) \to  \widehat{D} ^{(m)} _{\fP _\eta/\eta} \otimes _{\widetilde{D} ^{(m)\dag}_{\fP ^\sharp /\fS ^\sharp}} \iota ^{(m)\dag} _{+} (E)
\riso \iota ^{(m)} _{\eta +} (\widehat{E}) \to   \iota ^{(m)} _{\eta +} (\widehat{E} _{\mathrm{ntor}})$ be the composite morphism.  By functoriality, we get the commutative diagram
\begin{equation} \notag
\xymatrix{ {E}  \ar@{^{(}->}[d] ^-{} \ar[r] ^-{} &  {\mathcal{H} ^0 \iota ^! \iota ^{(m)\dag} _{+} (E)}  \ar[d] ^-{\mathcal{H} ^0 \iota ^!  \phi}
\\  {\widehat{E} _{\mathrm{ntor}}}  \ar[r] ^-{\sim} &  {\mathcal{H} ^0 \iota ^! \iota ^{(m)} _{\eta +} (\widehat{E} _{\mathrm{ntor}} ),}    }
\end{equation}
where the bottom arrow is an isomorphism following \cite[2.3.1]{Caro:2011fk}. This yields that the top arrow is injective.

ii) An element $\widehat{y}$ of 
$\iota ^{(m)} _{\eta +} (\widehat{E} _{\mathrm{ntor}})$ can we written uniquely  in the form 
$\sum _{\underline{i} ''\in \bbN ^{d-e'}}
[\underline{\partial} ^{<\underline{i}''> ^\star_{(m)}}]
\otimes 
\widehat{x} _{\underline{i}''}$,
where $\widehat{x} _{\underline{i}''} \in \widehat{E} _{\mathrm{ntor}}$ converges to $0$ for the $p$-adic topology 
when $|\underline{i}''|$ goes to infinity.
Let $u _1,\dots, u _r$ be some generators of $E$ as $\widetilde{D} ^{(m)\dag} _{\fZ /\fS }$-module. 
Let 
$y \in \iota ^{(m)\dag} _{+} (E)$. 
Using the formula \ref{wkcpdscplevelm-app2} , we can write $y$ of the form
$y =
\sum _{\alpha = 1} ^{r}
\left (\sum _{\underline{i} ''\in \bbN ^{d-e'}}
[\underline{\partial} ^{<\underline{i}''> ^\star_{(m)}}]
\cdot
P _{\underline{i}'' , \alpha} \right )
\otimes u _\alpha$, 
where 
$P _{\underline{i}'',\alpha}
\in 
\widetilde{D} ^{(m)\dag} _{\fZ  /\fS}$. 
We compute
$\phi (y) 
=
\sum _{\underline{i} ''\in \bbN ^{d-e'}}
[\underline{\partial} ^{<\underline{i}''> ^\star_{(m)}}]
\otimes 
\left (
\sum _{\alpha = 1} ^{r}
P _{\underline{i}'' , \alpha} 
\cdot u _\alpha
\right ).
$
If $y \in\mathcal{H} ^0 \iota ^! \iota ^{(m)\dag} _{+} (E)$
since $\phi$ is $\widetilde{D} ^{(m)\dag}_{\fP ^\sharp /\fS ^\sharp}$-linear,
then we get 
$\phi (y) \in \mathcal{H} ^0 \iota ^! \iota ^{(m)} _{\eta +} (\widehat{E} _{\mathrm{ntor}})$.
Hence, using 
\cite[2.3.1]{Caro:2011fk},
we get 
$\sum _{\alpha = 1} ^{r}
P _{\underline{i}'' , \alpha} 
\cdot u _\alpha
=0
$
when 
$| \underline{i}''| \not =0$.
Set 
$y _n := \sum _{\alpha = 1} ^{r}
\left (\sum _{|\underline{i} ''| \leq n}
[\underline{\partial} ^{<\underline{i}''> ^\star_{(m)}}]
\cdot
P _{\underline{i}'' , \alpha} \right )
\otimes u _\alpha$ for any $n\in \bbN$.
We compute
$y _n
=
\sum _{|\underline{i} ''| \leq n}
[\underline{\partial} ^{<\underline{i}''> ^\star_{(m)}}]
\otimes 
\left (
\sum _{\alpha = 1} ^{r}
P _{\underline{i}'' , \alpha} 
\cdot u _\alpha
\right )
=
y _0
$.
This implies
$y - y _0 \in 
\cap _{a=0} ^\infty
p ^a \iota ^{(m)\dag} _{+} (E)$.
Following \ref{Jacobson}, 
$\iota ^{(m)\dag} _{+} (E)$
is separated for the $p$-adic topology. 
Hence, $y = y _0$, i.e. 
the morphism
$E \to \mathcal{H} ^0 \iota ^! \iota ^{(m)\dag} _{+} (E)$
is surjective.

3) 
We prove that the ring $\widetilde{D} ^{(m)\dag} _{\fZ  /\fS }$ is noetherian (which extends \ref{noeth-weak completion}).

We can suppose that $\fZ '$ and $\fZ$ are replaced respectively by  $\fZ$ and $\fP$, i.e. $\iota$ is the morphism $\fZ\to \fP$. 
Let $I$ be a left ideal of $\widetilde{D} ^{(m)\dag} _{\fZ  /\fS }$.
Let $(I _n) _n$ be an increasing sequence of left ideals of finite type of $\widetilde{D} ^{(m)\dag} _{\fZ  /\fS }$ contained in $I$.
Since $\widetilde{D} ^{(m)\dag} _{\fZ  /\fS }$ is separated for the $p$-adic topology and has no $p$-torsion,
then so are $I$ and $I _n$.
Let 
$\overline{I} _n 
:= 
\cap _{a \in \bbN} 
( I _n + p ^a \widehat{D} ^{(m)} _{\fZ _\eta/\eta})$
be the closure of $I _n$ (for the induced topology) 
in $\widehat{D} ^{(m)} _{\fZ _\eta/\eta}$.
Since 
$\widehat{D} ^{(m)} _{\fZ _\eta/\eta}$
is noetherian, 
then $\overline{I} _n $ is a 
left ideal of 
$\widehat{D} ^{(m)} _{\fZ _\eta/\eta}$ of finite type.
 Moreover, 
an element $\widehat{y}$ of 
$\iota ^{(m)} _{\eta +} (\overline{I} _n )$ 
can we written uniquely  in the form 
$\sum _{\underline{i} ''\in \bbN ^{d-e'}}
[\underline{\partial} ^{<\underline{i}''> ^\star_{(m)}}]
\otimes 
\widehat{x} _{\underline{i}''}$,
where $\widehat{x} _{\underline{i}''} \in \overline{I} _n$ converges to $0$ for the $p$-adic topology 
when $|\underline{i}''|$ goes to infinity.

i) We check that the canonical morphism
$\psi \colon
\iota ^{(m)\dag} _{+} (I _n)
\to 
\iota ^{(m)} _{\eta +} (\overline{I} _n )$
is injective for any integer $n$.

Let $u _1,\dots, u _r$ be some generators of $I _n$ as $\widetilde{D} ^{(m)\dag} _{\fZ /\fS }$-module. 
Let $y \in \iota ^{(m)\dag} _{+} (I _n)$.
Following \ref{wkcpdscplevelm-app2} , we can write $y$ of the form 
$y =
\sum _{\alpha = 1} ^{r}
\left (\sum _{\underline{i} ''\in \bbN ^{d-e'}}
[\underline{\partial} ^{<\underline{i}''> ^\star_{(m)}}]
\cdot
P _{\underline{i}'' , \alpha} \right )
\otimes u _\alpha$, 
where 
$P _{\underline{i}'',\alpha}
\in 
\widetilde{D} ^{(m)\dag} _{\fZ  /\fS}$.
Then 
$\psi (y )
= 
\sum _{\alpha = 1} ^{r}
\left (\sum _{\underline{i} ''\in \bbN ^{d-e'}}
[\underline{\partial} ^{<\underline{i}''> ^\star_{(m)}}]
\cdot
P _{\underline{i}'' , \alpha} \right )
\otimes u _\alpha
= 
\sum _{\underline{i} ''\in \bbN ^{d-e'}}
[\underline{\partial} ^{<\underline{i}''> ^\star_{(m)}}]
\otimes 
\left (
\sum _{\alpha = 1} ^{r}
P _{\underline{i}'' , \alpha} 
\cdot u _\alpha
\right ).
$
Suppose 
$\psi (y )=0$. 
Then, for any 
$\underline{i}''$, 
we have
$\sum _{\alpha = 1} ^{r}
P _{\underline{i}'' , \alpha} 
\cdot u _\alpha
 =0$.
 This yields that
 $\sum _{\alpha = 1} ^{r}
\left (\sum _{|\underline{i} ''| \leq n}
[\underline{\partial} ^{<\underline{i}''> ^\star_{(m)}}]
\cdot
P _{\underline{i}'' , \alpha} \right )
\otimes u _\alpha
=
0$. 
Following \ref{Jacobson} (which can be applied because
$\fP ^\sharp  /\fS ^\sharp$ is log smooth), 
$\iota ^{(m)\dag} _{+} (I _n)$
is separated for the $p$-adic topology. 
Hence, $y= 0$.

ii) Since the functor $\iota ^{(m)} _{\eta +}$ is exact
(use for instance \cite[3.2.4]{Be1}
to check that 
$\widehat{D} ^{(m)}  _{\fZ _\eta \leftarrow \fZ ^{\prime} _\eta/\eta}$
is a flat right
$\widehat{D} ^{(m)} _{\fZ ^{\prime} _\eta/\eta}$-module), then 
using the part 3)i) we get that
both maps 
$\iota ^{(m)\dag} _{+} (I _n)
\to 
\iota ^{(m)\dag} _{+} (I _{n+1})
\to 
\iota ^{(m)\dag} _{+} (\widetilde{D} ^{(m)\dag} _{\fZ /\fS })$
are injective. 
Following \ref{noeth-weak completion}, 
$\widetilde{D} ^{(m)\dag} _{\fP ^\sharp  /\fS ^\sharp}$
is noetherian. Hence, for $n$ large enough, 
$\iota ^{(m)\dag} _{+} (I _n)
\riso
\iota ^{(m)\dag} _{+} (I _{n+1})$.
Since 
$I _n$
are separated for the $p$-adic topology and have no torsion,  then, using part 2) of the proof, we get that for $n$ large enough $I _n= I _{n+1}$.

4) Using the  special case of $\fP =\fP _\eta$, we get from 3) that  $D^{(m)\dag} _{\fZ _\eta/\eta} $ is noetherian. 
To check the remaining properties of the part (e) of the Lemma, we proceed in the same way as \ref{Jacobson}.
For instance in the same way as \ref{Jacobson}.(\ref{Jacobson(c)}), we can check that 
the homomorphisms of rings $\widetilde{D} ^{(m)\dag} _{\fZ /\fS}  \to   \widehat{D} ^{(m)} _{\fZ _\eta/\eta} $ 
and $D^{(m)\dag} _{\fZ _\eta/\eta}  \to  \widehat{D} ^{(m)} _{\fZ _\eta/\eta} $ are left and right faithfully flat, and then we are done.

5) To check the assertion (d),  we can follow the part 2) of the proof: indeed, thanks to the assertion e)ii) of the lemma that we have proved at the step 4), 
we can remove the hypothesis that $E$ is separated for the $p$-adic topology, and this extend \ref{Jacobson} that we used in the part 2) of the proof).

\end{proof}

\begin{lem}
\label{lem-flat-bimodule}
Let $M$ be a $(\widetilde{D} ^{(m)\dag} _{\fZ  }  , \widetilde{D} ^{(m)\dag} _{\fZ ^{\prime}}  )$-bimodule
such that $M $ is a  left $\widetilde{D} ^{(m)\dag} _{\fZ  }  $-module of finite type.
We suppose that the structure of $(\widehat{D} ^{(m)} _{\fZ  _\eta}  , \widetilde{D} ^{(m)\dag} _{\fZ ^{\prime}} )$-bimodule
of
$\widehat{D} ^{(m)} _{\fZ  _\eta}
\otimes 
_{\widetilde{D} ^{(m)\dag} _{\fZ  }  } M$
extends to a structure of 
$(\widehat{D} ^{(m)} _{\fZ  _\eta}, \widehat{D} ^{(m)} _{\fZ ^{\prime} _\eta} )$-bimodule
such that $\widehat{D} ^{(m)} _{\fZ  _\eta} \otimes _{\widetilde{D} ^{(m)\dag} _{\fZ  }  } M$
is a flat right $\widehat{D} ^{(m)} _{\fZ ^{\prime} _\eta}$-module. 
Then 
$M$ is a flat right $\widetilde{D} ^{(m)\dag} _{\fZ ^{\prime}} $-module.
\end{lem}

\begin{proof}
Let $N ' \to N$ be a monomorphism of left $\widetilde{D} ^{(m)\dag} _{\fZ ^{\prime}} $-modules.
Since the extension
$\widetilde{D} ^{(m)\dag} _{\fZ  }   \to \widehat{D} ^{(m)} _{\fZ  _\eta}$ is faithfully flat (see \ref{u+coh=cohetc-lem}.\ref{Jacobson-ZZ'}), 
we reduce to check that the top horizontal morphism
\begin{equation}
\xymatrix{
{\widehat{D} ^{(m)} _{\fZ  _\eta} 
\otimes _{\widetilde{D} ^{(m)\dag} _{\fZ  } }
M 
\otimes _{\widetilde{D} ^{(m)\dag} _{\fZ ^{\prime}} } 
N'}
\ar[r] ^-{} 
\ar[d] ^-{\sim}
& 
{\widehat{D} ^{(m)} _{\fZ  _\eta} 
\otimes _{\widetilde{D} ^{(m)\dag} _{\fZ  } }
M 
\otimes _{\widetilde{D} ^{(m)\dag} _{\fZ ^{\prime}} } 
N} 
\ar[d] ^-{\sim}
\\ 
{(\widehat{D} ^{(m)} _{\fZ  _\eta} 
\otimes _{\widetilde{D} ^{(m)\dag} _{\fZ  } }
M ) 
\otimes _{\widehat{D} ^{(m)} _{\fZ ^{\prime} _\eta}}
\widehat{D} ^{(m)} _{\fZ ^{\prime} _\eta}
\otimes _{\widetilde{D} ^{(m)\dag} _{\fZ ^{\prime}} } 
N'} 
\ar[r] ^-{}
& 
{(\widehat{D} ^{(m)} _{\fZ  _\eta} 
\otimes _{\widetilde{D} ^{(m)\dag} _{\fZ  } }
M ) 
\otimes _{\widehat{D} ^{(m)} _{\fZ ^{\prime} _\eta}}
\widehat{D} ^{(m)} _{\fZ ^{\prime} _\eta}
\otimes _{\widetilde{D} ^{(m)\dag} _{\fZ ^{\prime}} } 
N}  
}
\end{equation}
is injective. 
Since $\widehat{D} ^{(m)} _{\fZ  _\eta} 
\otimes _{\widetilde{D} ^{(m)\dag} _{\fZ  } }
M $ is $\widehat{D} ^{(m)} _{\fZ ^{\prime} _\eta}$-flat and 
since the extension 
$\widetilde{D} ^{(m)\dag} _{\fZ ^{\prime}}  \to 
\widehat{D} ^{(m)} _{\fZ ^{\prime} _\eta}$ is flat, the bottom horizontal arrow is injective. 
Hence, we are done.

\end{proof}

\begin{prop}
\label{flat-u+preZZ'}
The 
$(\widetilde{D} ^{(m)\dag} _{\fZ/\fS},\widetilde{D} ^{(m)\dag} _{\fZ ^{\prime} /\fS})$-bimodule 
$\widetilde{D} ^{(m)\dag} _{\fZ\leftarrow\fZ ^{\prime} /\fS}$ 
is flat for its underlying $\widetilde{D} ^{(m)\dag} _{\fZ ^{\prime} /\fS}$-module structure. 
\end{prop}

\begin{proof}
Since 
$\widehat{D} ^{(m)} _{\fZ  _\eta} \otimes _{\widetilde{D} ^{(m)\dag} _{\fZ  }  } 
\widetilde{D} ^{(m)\dag} _{\fZ\leftarrow\fZ ^{\prime} /\fS}
\riso 
\widehat{D} ^{(m)}  _{\fZ _\eta \leftarrow \fZ ^{\prime} _\eta/\eta}$, 
we check the proposition by applying Lemma \ref{lem-flat-bimodule}
in the case 
$M := \widetilde{D} ^{(m)\dag} _{\fZ\leftarrow\fZ ^{\prime} /\fS}$.
\end{proof}

\begin{coro}
\label{cor-flat-u+preZZ'}
With notation \ref{u+coh=cohetc-lem}, 
for any   left $\widetilde{D} ^{(m)\dag} _{\fZ '/\fS }$-module $E$, we put
$\iota ^{(m)\dag} _{+} (E)
:= 
\widetilde{D} ^{(m)\dag}  _{\fZ  \leftarrow \fZ ^{\prime}/\fS} \otimes _{\widetilde{D} ^{(m)\dag} _{\fZ '/\fS }}E$.
The functor $\iota ^{(m)\dag} _{+} $ is exact. 
\end{coro}

\begin{coro}
\label{u+coh=cohetc-cor}
Let $\iota  \colon \fZ ^{\prime} \to \fZ $ be the induced morphism. 
Let $E$ be a  left $\widetilde{D} ^{(m)\dag} _{\fZ '/\fS ,\bbQ}$-module of finite type, 
$\widehat{E}:=\widehat{D} ^{(m)} _{\fZ '_\eta/\eta,\bbQ} \otimes _{\widetilde{D} ^{(m)\dag} _{\fZ '/\fS ,\bbQ}} E$.
The push forward by $\iota$ of level $m$ of $E$
(the push forward by $\iota _\eta$ of level $m$ of $\widehat{E}$) is defined by setting
\begin{gather}
\notag
\iota ^{(m)\dag} _{+} (E)
:= 
\widetilde{D} ^{(m)\dag}  _{\fZ  \leftarrow \fZ ^{\prime}/\fS,\bbQ} \otimes _{\widetilde{D} ^{(m)\dag} _{\fZ '/\fS ,\bbQ}}E,
\, 
\notag
\iota ^{(m)} _{\eta +} (\widehat{E} )
:= 
\widehat{D} ^{(m)}  _{\fZ _\eta \leftarrow \fZ ^{\prime} _\eta/\eta,\bbQ}
 \otimes _{\widehat{D} ^{(m)} _{\fZ ^{\prime} _\eta/\eta,\bbQ}}
 \widehat{E}.
\end{gather}

\begin{enumerate}[(a)]

\item The $\widetilde{D} ^{(m)\dag} _{\fZ  /\fS ,\bbQ}$-module $\iota ^{(m)\dag} _{+} (E)$ is of finite type
(resp. the $\widehat{D} ^{(m)} _{\fZ  _\eta/\eta,\bbQ}$-module 
$\iota ^{(m)} _{\eta +} (\widehat{E} ) $ is of finite type). 

\item  \label{u+coh=cohetc-cor2)} We have the isomorphism
$$\iota ^{(m)} _{\eta +} (\widehat{E} )
\riso 
\widehat{D} ^{(m)} _{\fZ _\eta/\eta,\bbQ} \otimes _{\widetilde{D} ^{(m)\dag} _{\fZ /\fS ,\bbQ}}
\iota ^{(m)\dag} _{+} (E).$$

\item For any 
left $\widetilde{D} ^{(m)\dag} _{\fZ  /\fS ,\bbQ}$-module $G$, 
the module 
$\mathcal{H} ^0 \iota ^! (G):= 
\mathrm{Hom} _{\widetilde{B} }
(\widetilde{B}  /K ', G)$
is endowed with a structure of 
$\widetilde{D} ^{(m)\dag} _{\fZ  ^{\prime}/\fS ,\bbQ}$-module.

\item \label{noptors-isoQ}
The canonical homomorphism 
$E \to \mathcal{H} ^0 \iota ^! \iota ^{(m)\dag} _{+} (E)$
is an isomorphism.

\item 
\begin{enumerate}[(a)]
\item The ring $\widetilde{D} ^{(m)\dag} _{\fZ  /\fS ,\bbQ}$ is noetherian.
\item \label{Jacobson-ZZ'Q}
The homomorphisms of rings
$\widetilde{D} ^{(m)\dag} _{\fZ /\fS,\bbQ} 
\to 
D^{(m)\dag} _{\fZ _\eta/\eta,\bbQ} 
\to 
\widehat{D} ^{(m)} _{\fZ _\eta/\eta,\bbQ} $
are left and right faithfully flat. 

\end{enumerate}

\item The functor $\iota ^{(m)\dag} _{+}:= 
\widetilde{D} ^{(m)\dag}  _{\fZ  \leftarrow \fZ ^{\prime}/\fS,\bbQ} \otimes _{\widetilde{D} ^{(m)\dag} _{\fZ '/\fS ,\bbQ}}-$ is exact.

\item For any $m' \geq m$, 
the homomorphism of rings
$\widetilde{D} ^{(m)\dag} _{\fZ /\fS,\bbQ} 
\to 
\widetilde{D} ^{(m')\dag} _{\fZ /\fS,\bbQ} $
is left and right  flat. 

\item We have the canonical isomorphism
\begin{equation}
\label{u+coh=cohetc-cor-iso8}
\widetilde{D} ^{(m')\dag} _{\fZ /\fS,\bbQ}  \otimes _{\widetilde{D} ^{(m)\dag} _{\fZ /\fS,\bbQ} }
\iota ^{(m)\dag} _{+} (E)
\riso 
\iota ^{(m')\dag} _{+} (\widetilde{D} ^{(m')\dag} _{\fZ ' /\fS,\bbQ}  \otimes _{\widetilde{D} ^{(m)\dag} _{\fZ '/\fS,\bbQ} } E).
\end{equation}

\end{enumerate}

\end{coro}

\begin{proof}
Except for the last two statements, this is a straightforward consequence of \ref{u+coh=cohetc-lem} or of \ref{cor-flat-u+preZZ'}.
The seventh statement is a consequence of the fifth one and 
of the flatness of 
$\widehat{D} ^{(m)} _{\fZ _\eta/\eta,\bbQ} 
\to 
\widehat{D} ^{(m')} _{\fZ _\eta/\eta,\bbQ} $ (see \cite[3]{Be1}).
Since the functor involved in the last statement are exact, 
we reduce to the case where 
$E  = \widetilde{D} ^{(m)\dag} _{\fZ '/\fS,\bbQ}$, which is obvious.
\end{proof}

\subsubsection{Push-forwards of level $m$ of a weakly closed immersion: sheafification}
Let $u \colon \fX ^\sharp \to \fP ^{\sharp}$ be a weakly closed immersion of strictly semistable log formal schemes over $\fS ^\sharp$.
We denote by $\cI$ the ideal defining $\underline{u}$, and
$\cI _\eta: =\cI |_{\fP _\eta}$.
We denote by 
$j \colon \fX _\eta \to \fX ^{\sharp}$,
$j \colon \fP _\eta \to \fP ^{\sharp}$ the canonical open immersions.
Let $v \colon \fX _\eta \to \fP _\eta$ be the closed immersion of smooth  formal schemes over $\eta$ induced by $u$.

We set
$\widetilde{\cD} ^{(m)}  _{\fP^\sharp \leftarrow \fX ^\sharp/\fS^\sharp}:=
u ^{-1}(\widetilde{\cD} ^{(m)}  _{\fP^\sharp/\fS^\sharp}/ \widetilde{\cD} ^{(m)}  _{\fP^\sharp/\fS^\sharp} \cI)$,
$\cD ^{(m)}  _{\fP _\eta \leftarrow \fX _\eta/\eta}
:=
v ^{-1} (\cD ^{(m)}  _{\fP _\eta/\eta}/ \cD ^{(m)}  _{\fP _\eta/\eta} \cI _\eta )$,
$\widetilde{\cD} ^{(m)\dag}  _{\fP^\sharp \leftarrow \fX ^\sharp/\fS^\sharp}:=
u ^{-1} (\widetilde{\cD} ^{(m)\dag}  _{\fP^\sharp/\fS^\sharp}/ \widetilde{\cD} ^{(m)\dag}  _{\fP^\sharp/\fS^\sharp} \cI)$,
$\cD ^{(m)\dag}  _{\fP _\eta \leftarrow \fX _\eta/\eta}
:=
v ^{-1} (\cD ^{(m)\dag}  _{\fP _\eta/\eta}/ \cD ^{(m)\dag}  _{\fP _\eta/\eta} \cI _\eta )$,
$\widehat{\cD} ^{(m)}  _{\fP _\eta \leftarrow \fX _\eta/\eta}
:=
v ^{-1} (\widehat{\cD} ^{(m)}  _{\fP _\eta/\eta}/ \widehat{\cD} ^{(m)}  _{\fP _\eta/\eta} \cI _\eta )$.

\begin{empt}
\label{flat-basis-sheaf}
Let $\cM$ be a right $\widetilde{\cD} ^{(m)\dag}  _{\fP ^\sharp/\fS^\sharp}$-module. 
Suppose there exists a basis of open sets 
$\mathfrak{B}$ of $\fP$ such that, for any $\fU \in \mathfrak{B}$,
$\cM ( \fU)$ is a flat 
$\widetilde{\cD} ^{(m)\dag}  _{\fP ^\sharp/\fS^\sharp}(\fU)$-module.
Then, for any $x \in \fP$, 
$\cM _x$ is a flat $(\widetilde{\cD} ^{(m)\dag}  _{\fP ^\sharp/\fS^\sharp} )_x$-module
(use \cite[Proposition 9 in  section 2.7]{bourbaki}). 
Hence $\cM$ is a flat $\widetilde{\cD} ^{(m)\dag}  _{\fP ^\sharp/\fS^\sharp}$-module.

\end{empt}
 
\begin{lem}
\label{flat-Ddag/O}
The sheaf $\widetilde{\cD} ^{(m)\dag}  _{\fP^\sharp/\fS^\sharp}$
is a flat $\cO _{\fP}$-module.
\end{lem}

\begin{proof}
When $\fP= \Spf A$ is affine with local coordinates, since the composition
$\widetilde{D} ^{(m)}  _{\fP^\sharp/\fS^\sharp}
\to 
\widetilde{D} ^{(m)\dag}  _{\fP^\sharp/\fS^\sharp}
\to 
\widehat{D} ^{(m)}  _{\fP _\eta/\eta}$ is flat,
since the latter morphism is faithfully flat, 
then 
$\widetilde{D} ^{(m)}  _{\fP^\sharp/\fS^\sharp}
\to 
\widetilde{D} ^{(m)\dag}  _{\fP^\sharp/\fS^\sharp}$ is flat. 
Since both maps 
$A \to \widetilde{A}\to \widetilde{D} ^{(m)}  _{\fP^\sharp/\fS^\sharp}$ are flat, 
then so is 
$A \to \widetilde{D} ^{(m)\dag}  _{\fP^\sharp/\fS^\sharp}$.
We conclude using the similar arguments of \ref{flat-basis-sheaf}.
\end{proof}

\begin{lem}
\label{GammaDXP=sDXP}
With keep Notation \ref{ntn-DZZ'} in the case where 
$\fZ ' = \fX$ and $\fZ = \fP$.
We have the isomorphism 
$\Gamma ( \fX ^\sharp, \widetilde{\cD} ^{(m)\dag}  _{\fP^\sharp \leftarrow \fX ^\sharp/\fS^\sharp})
\riso
\widetilde{D} ^{(m)\dag}  _{\fP \leftarrow \fX /\fS}$.
\end{lem}

\begin{proof}
 Because it is a finite type submodule of a coherent module, $\widetilde{\cD} ^{(m)\dag}  _{\fP^\sharp/\fS^\sharp} \cI$ 
is a coherent  left $\widetilde{\cD} ^{(m)\dag}  _{\fP^\sharp/\fS^\sharp}$-module.
Since 
$\widetilde{\cD} ^{(m)\dag}  _{\fP^\sharp/\fS^\sharp}$
is a flat $\cO _{\fP}$-module (see \ref{flat-Ddag/O})
then 
$\widetilde{\cD} ^{(m)\dag}  _{\fP^\sharp/\fS^\sharp} \otimes _{\cO _{\fP}} \cI
\riso 
\widetilde{\cD} ^{(m)\dag}  _{\fP^\sharp/\fS^\sharp} \cI$.
Moreover,  using theorem of type $A$  for coherent $\cO _\fP$-modules, we get the first isomorphism
$$\widetilde{\cD} ^{(m)\dag}  _{\fP^\sharp/\fS^\sharp} \otimes _{\cO _{\fP}} \cI
\riso 
\widetilde{\cD} ^{(m)\dag}  _{\fP^\sharp/\fS^\sharp} \otimes _{\widetilde{D} ^{(m)\dag}  _{\fP^\sharp/\fS^\sharp}}
( \widetilde{D} ^{(m)\dag}  _{\fP^\sharp/\fS^\sharp} \otimes _A I)
\riso
\widetilde{\cD} ^{(m)\dag}  _{\fP^\sharp/\fS^\sharp} \otimes _{\widetilde{D} ^{(m)\dag}  _{\fP^\sharp/\fS^\sharp}}
( \widetilde{D} ^{(m)\dag}  _{\fP^\sharp/\fS^\sharp} I),$$
the second one comes from the flatness of 
$ \widetilde{D} ^{(m)\dag}  _{\fP^\sharp/\fS^\sharp}$ as $A$-module. 
Hence, 
using theorem of type $A$ for coherent $\widetilde{\cD} ^{(m)\dag}  _{\fP^\sharp /\fS^\sharp}$-modules (see \ref{thmAcohDtilde}),
we get
$\Gamma (\fP,  \widetilde{\cD} ^{(m)\dag}  _{\fP^\sharp/\fS^\sharp} \cI)
= 
\widetilde{D} ^{(m)\dag}  _{\fP^\sharp/\fS^\sharp} I$.

Since 
$\widetilde{\cD} ^{(m)\dag}  _{\fP^\sharp/\fS^\sharp}$ and
$\widetilde{\cD} ^{(m)\dag}  _{\fP^\sharp/\fS^\sharp}\cI$ are
coherent 
left $\widetilde{\cD} ^{(m)\dag}  _{\fP^\sharp/\fS^\sharp}$-modules, 
then so is 
$\widetilde{\cD} ^{(m)\dag}  _{\fP^\sharp \leftarrow \fX ^\sharp/\fS^\sharp}$.
Moreover, 
we get 
$\Gamma ( \fX ^\sharp, \widetilde{\cD} ^{(m)\dag}  _{\fP^\sharp \leftarrow \fX ^\sharp/\fS^\sharp})
\riso
\Gamma ( \fP ^\sharp, \widetilde{\cD} ^{(m)\dag}  _{\fP^\sharp/\fS^\sharp})
/ 
\Gamma ( \fP ^\sharp, \widetilde{\cD} ^{(m)\dag}  _{\fP^\sharp/\fS^\sharp} \cI)
\riso
\widetilde{D} ^{(m)\dag}  _{\fP^\sharp/\fS^\sharp} 
/
\widetilde{D} ^{(m)\dag}  _{\fP^\sharp/\fS^\sharp} I
=:\widetilde{D} ^{(m)\dag}  _{\fP^\sharp \leftarrow \fX ^\sharp/\fS^\sharp}$.
\end{proof}

\begin{prop} \label{ntn-Dleftarrowetc}
The following properties hold.
\begin{enumerate}[(a)]
\item The sheaf 
$\widetilde{\cD} ^{(m)}  _{\fP^\sharp \leftarrow \fX ^\sharp/\fS^\sharp}$
is a $(\widetilde{\cD} ^{(m)} _{\fP ^\sharp/\fS^\sharp}, \widetilde{\cD} ^{(m)} _{\fX ^\sharp/\fS^\sharp} )$-sub-bimodule
of
$ j_* \cD ^{(m)}  _{\fP _\eta \leftarrow \fX _\eta/\eta}$.

\item The sheaf
$\widetilde{\cD} ^{(m)\dag}  _{\fP^\sharp \leftarrow \fX ^\sharp/\fS^\sharp}$
is a 
$(\widetilde{\cD} ^{(m)\dag} _{\fP ^\sharp/\fS^\sharp}, \widetilde{\cD} ^{(m)\dag} _{\fX ^\sharp/\fS^\sharp} )$-sub-bimodule
of 
$ j_* \widehat{\cD} ^{(m)}  _{\fP _\eta \leftarrow \fX _\eta/\eta}$.

\end{enumerate}
\end{prop}

\begin{proof}
Since the first statement is checked similarly, let us focus on the second one. Since   this is local,  we can suppose we are in the local siutation of  Notation \ref{ntn-ZandZ'}.
Using \ref{lem-DZZ'subbimod} and \ref{GammaDXP=sDXP}, we can conclude.
\end{proof}

\begin{lem}
\label{faithflatX2mor}
The homomorphisms of rings
$\widetilde{\cD} ^{(m)\dag} _{\fX ^\sharp/\fS^\sharp} 
\to 
j _* \cD^{(m)\dag} _{\fX _\eta/\eta} 
\to 
j _*\widehat{\cD} ^{(m)} _{\fX _\eta/\eta} $
are left and right faithfully flat (see Definition after
\cite[4.3.8]{Be1}). 
\end{lem}

\begin{proof}
This is a consequence of \cite[4.3.8]{Be1} and \ref{u+coh=cohetc-lem}.\ref{Jacobson-ZZ'}.
\end{proof}

\begin{lem}
\label{lemu+coh=cohetc}
Let $\cE$ be a coherent left $\widetilde{\cD} ^{(m)\dag} _{\fX ^{\sharp}/\fS ^\sharp}$-module, 
$\cE _\eta:=\cE | \fX _\eta$, and
$\widehat{\E} _\eta$ the $p$-adic completion of $\cE  _\eta$.
\begin{enumerate}[(a)]
\item  The canonical morphism
$j _* \cD ^{(m)\dag} _{\fX _\eta/\eta} \otimes _{\widetilde{\cD} ^{(m)\dag} _{\fX ^{\sharp}/\fS ^\sharp}} \E
\to j _* \cE _\eta$ is an isomorphism.

\item  \label{lemu+coh=cohetc(b)} The canonical morphism
$j _* \widehat{\cD} ^{(m)} _{\fX _\eta/\eta} \otimes _{\widetilde{\cD} ^{(m)\dag} _{\fX ^{\sharp}/\fS ^\sharp}} \E
\to j _* \widehat{\E} _\eta$ is an isomorphism.
\end{enumerate}

\end{lem}

\begin{proof}
The first canonical morphism is the composition
$j _* \cD ^{(m)\dag} _{\fX _\eta/\eta} \otimes _{\widetilde{\cD} ^{(m)\dag} _{\fX ^{\sharp}/\fS ^\sharp}} \E
\to 
j _* j ^*  (j _* \cD ^{(m)\dag} _{\fX _\eta/\eta} \otimes _{\widetilde{\cD} ^{(m)\dag} _{\fX ^{\sharp}/\fS ^\sharp}} \cE)
\riso
j _* j ^* \cE$, 
which is functorial in $\cE$. 
The fact that this is an isomorphism is local.
Hence,  
we can suppose $\cE$  is the cokernel of a morphism of the form 
$(\widetilde{\cD} ^{(m)\dag} _{\fX ^{\sharp}/\fS ^\sharp}) ^r 
\to 
(\widetilde{\cD} ^{(m)\dag} _{\fX ^{\sharp}/\fS ^\sharp}) ^s$.
Since $j$ is an open immersion defined by the complementary of a Cartier  divisor, 
then $j_*$ is exact. 
This yields that the functor
$\cE \mapsto j _* j ^* \cE$
is exact.
Moreover, following Lemma \ref{faithflatX2mor}
we get the exactness of the functor
$\cE \mapsto j _* \cD ^{(m)\dag} _{\fX _\eta/\eta} \otimes _{\widetilde{\cD} ^{(m)\dag} _{\fX ^{\sharp}/\fS ^\sharp}} \cE$.
Hence, 
we reduce to the case where 
$\cE$ is a free 
$\widetilde{\cD} ^{(m)\dag} _{\fX ^{\sharp}/\fS ^\sharp}$-module of finite type. 
In that case, since 
$ \cD ^{(m)\dag} _{\fX _\eta/\eta} = j ^* \widetilde{\cD} ^{(m)\dag} _{\fX ^{\sharp}/\fS ^\sharp}$, 
this is clear that the homomorphism
$j _* \cD ^{(m)\dag} _{\fX _\eta/\eta} \otimes _{\widetilde{\cD} ^{(m)\dag} _{\fX ^{\sharp}/\fS ^\sharp}} \E
\to j _* \cE _\eta$
is an isomorphism.
We proceed similarly for the second statement.
\end{proof}

\begin{prop}
\label{flat-u+pre}
The $(\widetilde{\cD} ^{(m)\dag} _{\fP ^\sharp/\fS^\sharp}, \widetilde{\cD} ^{(m)\dag} _{\fX ^\sharp/\fS^\sharp} )$-bimodule
$\widetilde{\cD} ^{(m)\dag}  _{\fP^\sharp \leftarrow \fX ^\sharp/\fS^\sharp}$
is flat for its underlying $\widetilde{\cD} ^{(m)\dag} _{\fX ^\sharp/\fS^\sharp} $-module structure. 
\end{prop}

\begin{proof}
This is a consequence of \ref{flat-basis-sheaf}, \ref{GammaDXP=sDXP} and 
Proposition \ref{flat-u+preZZ'}. 
\end{proof}

\begin{prop}
\label{u+coh=cohetc}
Let $\cE$ be a coherent left $\widetilde{\cD} ^{(m)\dag} _{\fX ^{\sharp}/\fS ^\sharp}$-module, 
$\cE _\eta:=\cE | \fX _\eta$, and
$\widehat{\E} _\eta$ the $p$-adic completion of $\cE  _\eta$.
The push forward by $u$ of level $m$ of $\cE$
(the push forward by $v$ of level $m$ of $\cE _\eta$, resp. of $\widehat{\E} _\eta$) is defined by setting
\begin{gather}
\notag
u ^{(m)\dag} _{+} (\cE)
:= u _* (\widetilde{\cD} ^{(m)\dag}  _{\fP^\sharp \leftarrow \fX ^\sharp/\fS^\sharp} \otimes _{\widetilde{\cD} ^{(m)\dag} _{\fX ^{\sharp}/\fS ^\sharp}}\cE),
v ^{(m)\dag} _{+} (\cE _\eta)
:= v _* (\cD ^{(m)\dag}  _{\fP_\eta \leftarrow \fX _\eta/\eta} \otimes _{\cD ^{(m)\dag} _{\fX _\eta/\eta}}\cE _\eta),
\\ 
\notag
v ^{(m)} _{+} (\widehat{\E} _\eta)
:= v _* (\widehat{\cD} ^{(m)}  _{\fP_\eta \leftarrow \fX _\eta/\eta} \otimes _{\widehat{\cD} ^{(m)} _{\fX _\eta/\eta}}\widehat{\E} _\eta).
\end{gather}

\begin{enumerate}[(a)]
\item The $\widetilde{\cD} ^{(m)\dag} _{\fP ^{\sharp}/\fS ^\sharp}$-module $u ^{(m)\dag} _{+} (\cE)$ is coherent
(resp. the $\cD ^{(m)\dag} _{\fP _\eta/\eta}$-module 
$v ^{(m)\dag} _{+} (\cE _\eta) $ is coherent),
(resp. the $\widehat{\cD} ^{(m)} _{\fP _\eta/\eta}$-module 
$v ^{(m)} _{+} (\widehat{\E} _\eta) $ is coherent). 

\item  
We have the isomorphism
$$j _* v ^{(m)\dag} _{+} (\cE _\eta)
\riso 
j _* \cD ^{(m)\dag} _{\fP _\eta/\eta} \otimes _{\widetilde{\cD} ^{(m)\dag} _{\fP ^{\sharp}/\fS ^\sharp}}
u ^{(m)\dag} _{+} (\cE).$$

\item  We have the isomorphism
$$j _* v ^{(m)} _{+} (\widehat{\E} _\eta)
\riso 
j _* \widehat{\cD} ^{(m)} _{\fP _\eta/\eta} \otimes _{\widetilde{\cD} ^{(m)\dag} _{\fP ^{\sharp}/\fS ^\sharp}}
u ^{(m)\dag} _{+} (\cE).$$
\end{enumerate}

\end{prop}

\begin{proof}
Let us check the first statement. Since the respective case is checked similarly, let us only consider the non respective case. Since this is local, we can suppose that $\cE$ is the cokernel of a morphism of the form 
$(\widetilde{\cD} ^{(m)\dag} _{\fX ^{\sharp}/\fS ^\sharp}) ^r  \to  (\widetilde{\cD} ^{(m)\dag} _{\fX ^{\sharp}/\fS ^\sharp}) ^s$. From \ref{flat-u+pre}, we get that the functor $u ^{(m)\dag} _{+}$ is exact. 
Hence $u ^{(m)\dag} _{+} (\cE)$ is the cokernel of  $u _* (\widetilde{\cD} ^{(m)\dag}  _{\fP^\sharp \leftarrow \fX ^\sharp/\fS^\sharp}) ^r \to u _* (\widetilde{\cD} ^{(m)\dag}  _{\fP^\sharp \leftarrow \fX ^\sharp/\fS^\sharp}) ^s$.
Since $u _* (\widetilde{\cD} ^{(m)\dag}  _{\fP^\sharp \leftarrow \fX ^\sharp/\fS^\sharp}) ^s$ is $\widetilde{\cD} ^{(m)\dag} _{\fP ^{\sharp}/\fS ^\sharp}$-coherent, then we are done.

Let us check the second statement.  Using Lemma \ref{lemu+coh=cohetc} (for $\fP ^\sharp$ instead of $\fX ^\sharp$), we reduce to check 
the isomorphism $v ^{(m)\dag} _{+} (\cE _\eta) \riso   j ^*u ^{(m)\dag} _{+} (\cE)$. Since the functors  $\cE \mapsto v ^{(m)\dag} _{+} (\cE _\eta)$ and $\cE \mapsto j ^*u ^{(m)\dag} _{+} (\cE)$ are exact, 
we reduce to the case where  $\cE = \widetilde{\cD} ^{(m)\dag} _{\fX ^{\sharp}/\fS ^\sharp}$, i.e. we reduce to check the straightforward isomorphism
$
\cD ^{(m)\dag}  _{\fP _\eta/\eta}/ \cD ^{(m)\dag}  _{\fP _\eta/\eta} \cI _\eta  \riso  j ^* (\widetilde{\cD} ^{(m)\dag}  _{\fP^\sharp/\fS^\sharp}/ \widetilde{\cD} ^{(m)\dag}  _{\fP^\sharp/\fS^\sharp} \cI)$.
We proceed similarly for the last statement.
\end{proof}

\begin{ntn}
\begin{enumerate}[(a)]
\item For any coherent left $\widetilde{\cD} ^{(m)\dag} _{\fX ^{\sharp}/\fS ^\sharp,\bbQ}$-module $\cE$, we put  $u ^{(m)\dag} _{+} (\cE)
:= u _* (\widetilde{\cD} ^{(m)\dag}  _{\fP^\sharp \leftarrow \fX ^\sharp/\fS^\sharp,\bbQ} \otimes _{\widetilde{\cD} ^{(m)\dag} _{\fX ^{\sharp}/\fS ^\sharp,\bbQ}}\cE)$.

\item For any  coherent $\widetilde{\cD} ^{(m)\dag} _{\fP ^{\sharp}/\fS ^\sharp,\bbQ}$-module $\G$,  put 
$\mathcal{H} ^0 u ^! (\G):=  \mathcal{H} om _{\widetilde{\cO} _{\fP}} (\widetilde{\cO} _{\fP} /\widetilde{\cI} , \G).$
\end{enumerate}
\end{ntn}

\begin{prop}
\label{u_+faithful}
Let $\cE$ be a coherent left $\widetilde{\cD} ^{(m)\dag} _{\fX ^{\sharp}/\fS ^\sharp,\bbQ}$-module.
The canonical homomorphism 
$\cE \to \mathcal{H} ^0 u ^! u ^{(m)\dag} _{+} (\cE)$
is an isomorphism.
The functor $u _+ \colon \mathrm{Coh} (\widetilde{\cD} ^{(m)\dag} _{\fX ^{\sharp}/\fS ^\sharp, \bbQ}) 
\to 
\mathrm{Coh} (\widetilde{\cD} ^{(m)\dag} _{\fP ^{\sharp}/\fS ^\sharp,\bbQ})$
is exact and fully faithful.
\end{prop}

\begin{proof}
From the canonical map
$\cE \to \widetilde{\cD} ^{(m)\dag}  _{\fP^\sharp \leftarrow \fX ^\sharp/\fS^\sharp,\bbQ} \otimes _{\widetilde{\cD} ^{(m)\dag} _{\fX ^{\sharp}/\fS ^\sharp,\bbQ}}\cE$ given by $x\mapsto [1] \otimes x$,  we get the homomorphism
$\cE \to \mathcal{H} ^0 u ^! u ^{(m)\dag} _{+} (\cE) \subset  \widetilde{\cD} ^{(m)\dag}  _{\fP^\sharp \leftarrow \fX ^\sharp/\fS ^\sharp ,\bbQ} \otimes _{\widetilde{\cD} ^{(m)\dag} _{\fX ^{\sharp}/\fS ^\sharp ,\bbQ}}\cE$.
It remains to prove that this is an isomorphism. Since this is local, using Lemma \ref{ovcv-coord-comp-imm},  we can suppose that we are in the local context of \ref{ntn-ZandZ'}. 
Set $E := \Gamma (\fX, \cE)$.
Using theorem of type $A$ for coherent $\widetilde{\cD} ^{(m)\dag}  _{\fX^\sharp/\fS^\sharp,\bbQ}$-modules (see \ref{cor-thADdagmQ}),
the canonical morphism
$\widetilde{\cD} ^{(m)\dag}  _{\fX^\sharp /\fS^\sharp,\bbQ} 
\otimes _{\widetilde{D} ^{(m)\dag}  _{\fX^\sharp /\fS^\sharp,\bbQ}} E
\to 
\cE$
is an isomorphism.
Since 
$\widetilde{D} ^{(m)\dag}  _{\fP^\sharp \leftarrow \fX ^\sharp/\fS^\sharp,\bbQ}$
is a coherent left
$\widetilde{D} ^{(m)\dag}  _{\fP^\sharp/\fS ^\sharp ,\bbQ}$-module,
then using \ref{GammaDXP=sDXP} we check that the canonical morphism
$\widetilde{\cD} ^{(m)\dag}  _{\fP^\sharp /\fS ^\sharp ,\bbQ} 
\otimes _{\widetilde{D} ^{(m)\dag}  _{\fP^\sharp /\fS ^\sharp ,\bbQ}} 
\widetilde{D} ^{(m)\dag}  _{\fP^\sharp \leftarrow \fX ^\sharp/\fS ^\sharp ,\bbQ}
\to
\widetilde{\cD} ^{(m)\dag}  _{\fP^\sharp \leftarrow \fX ^\sharp/\fS ^\sharp ,\bbQ}$
is an isomorphism.
Hence, we get the isomorphisms
$$\widetilde{\cD} ^{(m)\dag}  _{\fP^\sharp \leftarrow \fX ^\sharp/\fS ^\sharp ,\bbQ} \otimes _{\widetilde{\cD} ^{(m)\dag} _{\fX ^{\sharp}/\fS ^\sharp ,\bbQ}}\E
\riso 
\widetilde{\cD} ^{(m)\dag}  _{\fP^\sharp \leftarrow \fX ^\sharp/\fS ^\sharp ,\bbQ} \otimes _{\widetilde{D} ^{(m)\dag} _{\fX ^{\sharp}/\fS ^\sharp ,\bbQ}}E
\riso
\widetilde{\cD} ^{(m)\dag}  _{\fP^\sharp /\fS ^\sharp ,\bbQ} 
\otimes _{\widetilde{D} ^{(m)\dag}  _{\fP^\sharp /\fS ^\sharp ,\bbQ}} 
(\widetilde{D} ^{(m)\dag}  _{\fP^\sharp \leftarrow \fX ^\sharp/\fS ^\sharp ,\bbQ}\otimes _{\widetilde{D} ^{(m)\dag} _{\fX ^{\sharp}/\fS ^\sharp ,\bbQ}}E).$$
Hence, 
$\Gamma (\fP, u ^{(m)\dag} _{+} (\cE)) 
= 
\widetilde{D} ^{(m)\dag}  _{\fP^\sharp \leftarrow \fX ^\sharp/\fS ^\sharp ,\bbQ}\otimes _{\widetilde{D} ^{(m)\dag} _{\fX ^{\sharp}/\fS ^\sharp ,\bbQ}}E
=:
u ^{(m)\dag} _{+} (E)$.
This yields 
$\Gamma (\fP, \mathcal{H} ^0 u ^! u ^{(m)\dag} _{+} (\cE)) 
=
\mathcal{H} ^0 u ^! u ^{(m)\dag} _{+} (E)$ (see notation \ref{u+coh=cohetc-cor} for the latter term).
We conclude using 
\ref{u+coh=cohetc-lem}.\ref{noptors-iso}.

\end{proof}

\subsubsection{Berthelot-Kashiwara's theorem for weakly closed immersions of strictly semistable log formal schemes over $\fS ^\sharp$}
\begin{ntn}
\label{tildePnpre}
Let $\Delta _{\fP^\sharp/\fS^\sharp}\colon \fP ^\sharp\hookrightarrow \fP ^\sharp\times _{\fS^\sharp} \fP^\sharp$ be the diagonal immersion.
From \cite[3.3.5]{Caro-Vauclair},
we have the exact closed immersion of fine $\fS ^\sharp$-log formal schemes of the form 
$\fP ^\sharp \hookrightarrow \Delta _{\fP^\sharp /\fS^\sharp} ^n$ (which is the $n$-th infinitesimal neighborhood of 
$\Delta _{\fP^\sharp/\fS^\sharp}$)
and 
the strict and finite morphism of fine log formal $\cV$-schemes 
$p ^n _0\colon \Delta _{\fP ^\sharp/\fS^\sharp} ^n \to \fP^\sharp$
(resp. $p ^n _1\colon \Delta _{\fP^\sharp /\fS^\sharp} ^n \to \fP^\sharp$).
We denote by 
$\mathcal{P} ^{n} _{\fP^\sharp/\fS^\sharp}$ the coherent $\cO _\fP$-algebra such that 
$\Spf \mathcal{P} ^{n} _{\fP^\sharp/\fS^\sharp} =  \underline{\Delta} ^n _{\fP^\sharp/\fS^\sharp}  $.

Via ring homomorphism $\cD _{\fP ^\sharp/\fS ^\sharp} \to  j _* \cD _{\fP _\eta / \eta}$ we get a structure of  left $\cD _{\fP ^\sharp/\fS ^\sharp}$-module on $j _* \cO _{\fP _\eta}$.  We compute that 
$\widetilde{\cO} _{\fP}$ is a $\cD _{\fP ^\sharp/\fS ^\sharp}$-submodule of $j _* \cO _{\fP _\eta}$. Since this action satisfies Leibniz's formula, we get the isomorphism 
$\widetilde{\cO} _{\fP}
\otimes _{\cO _{\fP}} 
\mathcal{P}  ^n  _{\fP^\sharp/\fS^\sharp}
\riso 
\mathcal{P}  ^n  _{\fP^\sharp/\fS^\sharp}
\otimes _{\cO _{\fP}} 
\widetilde{\cO} _{\fP}$
given by the structure of 
left $\cD _{\fP ^\sharp/\fS ^\sharp}$-module on 
$\widetilde{\cO} _{\fP}$
is in fact an isomorphism of 
$\mathcal{P}  ^n  _{\fP^\sharp/\fS^\sharp}$-algebras. 
We put 
$\widetilde{\mathcal{P}}  ^n  _{\fP^\sharp/\fS^\sharp}
:= 
\widetilde{\cO} _{\fP}
\otimes _{\cO _{\fP}} 
\mathcal{P}  ^n  _{\fP^\sharp/\fS^\sharp}$.

Using \cite[16.3.9]{EGAIV4}, we check that 
$\mathcal{P}  ^n  _{\fP/\fS}
:=
\cO _{\fP \times _\fS \fP}/ \cI _\Delta ^{n+1}$ 
is a coherent $\cO _{\fP}$-module for the right or left structures. 
From \ref{prop-regularbis}, 
$\widetilde{\mathcal{P}} ^n  _{\fP/\fS}:=
\widetilde{\cO} _{\fP}
\otimes _{\cO _{\fP}} 
\mathcal{P}  ^n  _{\fP/\fS}
\riso 
\widetilde{\cO} _{\fP \times _\fS \fP}/ \widetilde{\cI}
 _\Delta ^{n+1}$.
The canonical morphism 
$\mathcal{P}  ^n  _{\fP/\fS}
\to 
\mathcal{P}  ^n  _{\fP^\sharp/\fS^\sharp}$
induces the morphism of coherent $\widetilde{\cO} _{\fP}$-modules
$\phi \colon \widetilde{\mathcal{P}}  ^n  _{\fP/\fS}
\to 
\widetilde{\mathcal{P}}   ^n  _{\fP^\sharp/\fS^\sharp}$.
Since $\phi |\fP _\eta$ is an isomorphism, then so is $\phi$.
Since $\widetilde{\cI}
 _\Delta$ is regular, we get
$\widetilde{\mathcal{P}}  ^n  _{\fP^\sharp/\fS^\sharp}=
\widetilde{\cO} _{\fP \times _{\fS} \fP} / \widetilde{\cI}
 _\Delta ^{n+1}
=
\oplus _{i=0} ^{n}
\widetilde{\cI}
 _\Delta ^i / \widetilde{\cI}
 _\Delta ^{i+1}$.

\end{ntn}

\begin{ntn}
\label{tildePn}
With notation \ref{tildePnpre}, suppose $\fP ^{\sharp} /\fS ^{\sharp} $
has log coordinates $u _1,\dots, u _d$.

1) Let $n\in \bbN$ be an integer. If $a \in M _{\fP ^\sharp}$, we denote by
$\mu  ^n (a)$ the unique section of
$\ker ( \cO _{\Delta ^n _{\fP^\sharp/\fS^\sharp}} ^{*} \to  \cO _{\fP ^\sharp} ^{*} )$
such that we get in $M ^n _{\fP^\sharp/\fS^\sharp}$ the equality
$p _1 ^{n*} (a)= p _0 ^{n*} (a) \mu ^n  (a)$.
We get 
$\mu  ^{n} \colon M  _{\fP ^\sharp} \to \ker ( \cO _{\Delta ^{n}  _{\fP^\sharp/\fS^\sharp}} ^{*} \to  \cO _{\fP} ^{*} )$
given by 
$a \mapsto \mu ^n  (a)$.
Put 
$\eta _{\lambda,n} := \mu ^n  ( u _{\lambda}) -1$.
Following \cite[3.3.6]{Caro-Vauclair},
we have  the isomorphism of $\cO _{\fP}$-algebras: 
   \begin{align}
   \notag
   \cO _{\fP} [T _1 ,\dots ,T _r] _n
   &
   \riso 
   \mathcal{P}  ^n  _{\fP^\sharp/\fS^\sharp}  \\
   \label{loc-desc-Pnfformal}
   T _\lambda &\mapsto
   \eta _{\lambda,n},
   \end{align}
where    $\cO _{\fP} [T _1 ,\dots ,T _r] _n :=    \cO _{\fP} [T _1 ,\dots ,T _r] / (T _1, \dots, T _r ) ^{n+1}$
and where
the structure of $\cO _\fP$-module of 
$\mathcal{P}  ^n  _{\fP^\sharp/\fS^\sharp}$
is given by $p ^{n} _1$ or $p ^{n} _0$.
We get the dual $\cO _{\fP}$-basis 
$\{ \underline{\partial} _{\sharp}  ^{[\underline{i}]} \ ; \
\underline{i} \in \bbN ^{d}\}$ 
of 
$\cD _{\fP ^\sharp /\fS ^\sharp,n}:= 
\mathcal{H}om _{\cO _\fP}
(p _{0*} ^n \mathcal{P}  ^n  _{\fP^\sharp/\fS^\sharp}, \cO _{\fP})$.
When $n =1$, we set
$d \log u _\lambda : = \eta _{\lambda,1}$.
Then
$d \log u _1, \dots, d \log u _d$ 
is a basis of 
$\Omega _{\fP ^\sharp /\fS ^\sharp}$
($\subset    \mathcal{P}  ^1  _{\fP^\sharp/\fS^\sharp} $) as $\cO _{\fP}$-module.

2) Since $M _{\fP ^{\sharp}} \subset (\widetilde{\cO} _{\fP}) ^*$, 
we can define $\mu (a)  := a ^{-1} \otimes a
\in 
\widetilde{\cO} _{\fP \times _{\fS} \fP} ^*$ for any 
$a \in M _{\fP ^{\sharp}}$.
Put 
$\eta _{\lambda} := \mu   ( u _{\lambda}) -1 \in \widetilde{\cO} _{\fP \times _{\fS} \fP}$.

3) We compute that the image of $\eta _{\lambda,n}$ via
$\mathcal{P}  ^n  _{\fP^\sharp/\fS^\sharp}
\hookrightarrow 
\widetilde{\mathcal{P}}  ^n  _{\fP^\sharp/\fS^\sharp}=
\widetilde{\cO} _{\fP \times _{\fS} \fP} / \widetilde{\cI}
 _\Delta ^{n+1}
=
\oplus _{i=0} ^{n}
\widetilde{\cI}
 _\Delta ^i / \widetilde{\cI}
 _\Delta ^{i+1}$
is the image of 
$\eta _{\lambda}$ modulo $\widetilde{\cI}
 _\Delta ^{n+1}$ 
via 
$\widetilde{\cO} _{\fP \times _{\fS} \fP} 
\twoheadrightarrow
\widetilde{\cO} _{\fP \times _{\fS} \fP} / \widetilde{\cI}
 _\Delta ^{n+1}
=
\widetilde{\mathcal{P}}  ^n  _{\fP^\sharp/\fS^\sharp}$.
In the same way as \cite[16.9.3]{EGAIV4}, 
we check that 
$\underline{\eta} ^{\underline{i}}$
with 
$|\underline{i}| =n$ induce modulo 
$\widetilde{\cI}
 _\Delta ^{n+1}$ a basis of 
$\widetilde{\cI}
 _\Delta ^n / \widetilde{\cI}
 _\Delta ^{n+1}$.
\end{ntn}

\begin{ntn} \label{otimes-matrix}
Let $R$ be a commutative ring. Let $M,N$ be two $R$-modules.  Let $\underline{x}= \overset{t}{}( x _1,\dots, x _n) \in M _{n,1} (M)$ be a column vector with $x _i \in M$,
$\underline{y}= \overset{t}{}( y _1,\dots, y _m)\in M _{m,1} (N)$ be a column vector with $y _i \in N$. 
We denote by $\underline{x} \otimes \underline{y}\in M _{nm,1} (M\otimes _R N)$ 
the column vector of $nm$ elements of $M \otimes _{R} N$ whose $(i, j)$th element is $x _i \otimes y _j$ for $1\leq i \leq n$ and $1\leq j\leq m$, 
where we order $\{1,\dots, n \} \times \{1,\dots, m\}$ with the lexicographic order.
Hence, $\underline{x} \otimes \underline{y}= ( ( x _1 \otimes \underline{y}) ,\dots, (x _n \otimes \underline{y}))$.

Let $\underline{x}'= \overset{t}{}( x '_1,\dots, x '_n) \in M _{n,1} (M)$, $\underline{y}'= \overset{t}{}( y '_1,\dots, y '_m)\in M _{m,1} (N)$, 
$A = (a _{ij}) _{1\leq i,j \leq n} \in M _{n} (R)$,$B = (b _{ij}) _{1\leq i,j \leq m} \in M _{m} (R)$
such that $\underline{x}' = A \underline{x}$ and $\underline{y}' = B \underline{y}$. Then, we have 
$\underline{x}' \otimes \underline{y} '= (A \otimes B) ( \underline{x}\otimes \underline{y} )$, 
where $A \otimes B \in M _{nm} (R)$ is the usual tensor product of the matrices $A$ and $B$.
\end{ntn}

\begin{lem}
\label{compa-infinitelevel}
We suppose  $\fP ^\sharp$ is affine 
and has logarithmic local coordinates.
Set 
$A := \Gamma (\fP,\cO _{\fP})$ 
and
$\widetilde{A}:=  \Gamma (\fP,\widetilde{\cO} _{\fP})$.
Let  $t _1,\dots, t _d  \in A  \subset \widetilde{A}$ 
be overconvergent local coordinates of $\fP ^\sharp/\fS ^\sharp$.
With Notation \ref{localcoord-ovcv}, 
then there exists $\vartheta \in A \cap \widetilde{A} ^*$ such that for any $\underline{k}\in \bbN ^{d}$, we have
\begin{equation}
\label{varthetapartialDlog}
\vartheta ^{|\underline{k}|} \underline{\partial} ^{[\underline{k}]} \in 
\Gamma (\fP, \cD _{\fP ^\sharp/\fS ^\sharp, |\underline{k}|}).
\end{equation}

\end{lem}

\begin{proof}
I) Notations. 1) Fix an integer $n \in \bbN$. 
Consider the map
$\rho \colon \{1,\dots, d \} ^{n}\to \bbN ^d$ 
defined by 
$(\lambda _1, \dots, \lambda _n) \mapsto ( \# \{i \; ; \; \lambda _i = 1 \}, \dots, \# \{i \; ; \; \lambda _i = d \})$.
The image of $\rho$ consists on elements $\underline{l}\in \bbN ^d$ such that $|\underline{l}|= n$.
Choose $\iota \colon  \mathrm{Im} (\rho) \to  \{1,\dots, d \} ^{n}$ a section of 
$\rho$.

2) We denote by $\tau _i : = 1 \otimes t _i - t _i \otimes 1 \in
 \Gamma( \fP  \times _{\fS } \fP, \cO _{\fP \times _{\fS } \fP } )  $ 
and by $d t _i$ the image of $\tau _i $ in $\widetilde{\Omega} _{\fP ^\sharp /\fS ^\sharp}=
\widetilde{\cI}
 / \widetilde{\cI}
 ^2$.
From the commutativity of the diagram \ref{dtildeOmegasharp}, 
we notice that $d t _1, \dots, d t _d \in \Gamma (\fP, \Omega _{\fP ^{\sharp}/\fS ^{\sharp}})$.
Put 
$\underline{x}
:=
\overset{t}{}
(d t _1, \dots, d t _d)$.

Let 
$u _1,\dots, u _d$ be log coordinates of $\fP ^{\sharp} /\fS ^{\sharp} $.
Put
$\underline{y}
:=
\overset{t}{}
(d \log u _1, \dots, d \log u _d)$.
Following notation \ref{tildePn},
put 
$\eta _{\lambda,n} := \mu ^n  ( u _{\lambda}) -1$, for any integer $n\in \bbN$ 
and $\lambda=1,\dots, d$.

3) Since 
$d \log u _1, \dots, d \log u _d$ 
is a basis of 
$\Omega _{\fP ^\sharp /\fS ^\sharp}$ as $\cO _{\fP}$-module, 
there exists a matrix 
$N:= ( a _{ij} ) \in M _d (A)$ 
such that 
$\underline{x}
=
N
\underline{y}
$.
Since
$d t _1, \dots, d t _d$ and
$d \log u _1, \dots, d \log u _d$ 
are two bases of 
$\widetilde{\Omega} _{\fP /\fS }$
as $\widetilde{\cO} _{\fP}$-module (see \ref{localcoord-ovcv}), 
then 
$N\in \mathrm{GL} _d (\widetilde{A})$.
Let $\vartheta := \det N \in A \cap \widetilde{A} ^*$.
Put $M:= \vartheta N ^{-1} \in  M _d (A)$. 
We get 
$
\vartheta \underline{y} 
=
M \underline{x}$.
Hence, 
$
\vartheta ^n \underline{y} ^{\otimes n}
=
M ^{\otimes n} \underline{x} ^{\otimes n}$ (see notation and formula in \ref{otimes-matrix}).
The elements  of  $\{1,\dots, d \} ^{n}$ 
are written in the form 
$\underline{j} = (j _1, \dots, j _n) $
and  are ordered following the lexicographic order. 
For any 
$\underline{j} \in \{1,\dots, d \} ^{n}$, we denote by 
$(\underline{x} ^{\otimes n}) _{\underline{j}}$
the $\underline{j}$th component of 
$\underline{x} ^{\otimes n}$, i.e. 
$(\underline{x} ^{\otimes n}) _{\underline{j}} 
=
x _{j _1} \otimes x _{j _2} \otimes \cdots \otimes x _{j _d}$.
Set $\widetilde{I}:= \Gamma (\fP, \widetilde{\cI}
)$. 
For any $\underline{i} \in \{1,\dots, d \} ^{n}$, we get the equality in 
$(\widetilde{I}/\widetilde{I} ^2) ^{\otimes n} $: 
$$\vartheta ^n (\underline{y} ^{\otimes n}) _{\underline{i}}
=
\sum _{\underline{j}\in \{1,\dots, d \} ^{n}}
a _{\underline{i}, \underline{j}} 
(\underline{x} ^{\otimes n}
) _{\underline{j}}$$
where $a _{\underline{i}, \underline{j}}$ are the coefficients of 
$M ^{\otimes n}\in M _{d ^n} (A)$.

II) For any $\underline{k}\in \mathrm{Im} (\rho)$, 
we remark $\underline{\tau} ^{\underline{k}} _n= (d t _1 )^{k _1} \cdots (d t  _d ) ^{k _d} \mod \widetilde{I} ^{n+1}$,
and 
$\underline{\eta} ^{\underline{k}} _n
= 
(d \log u _1) ^{k _1} \cdots (d \log u _d) ^{k _d} \mod \widetilde{I} ^{n+1}$.
Hence, the canonical homomorphism
$(\widetilde{I}/\widetilde{I} ^2) ^{\otimes n} \to \widetilde{I} ^n / \widetilde{I} ^{n+1}$
sends 
$(\underline{x} ^{\otimes n}) _{\underline{j}} $
(resp. 
$(\underline{y} ^{\otimes n}) _{\underline{i}}$)
to 
$\underline{\tau} ^{\rho (\underline{j})} _n $
(resp. $\underline{\eta} ^{\rho (\underline{i})} _n$).
This yields the equality in $ \widetilde{I} ^n / \widetilde{I} ^{n+1}$ :
$$\vartheta ^n \underline{\eta} ^{\rho (\underline{i})} _n
=
\sum _{\underline{j}\in \{1,\dots, d \} ^{n}}
a _{\underline{i}, \underline{j}} \underline{\tau} ^{\rho (\underline{j})} _n ,$$
where the structure of $\widetilde{A} $-module of 
$ \widetilde{I} ^n / \widetilde{I} ^{n+1}$ is the left one (if not clarified).

From now, fix $\underline{k}\in \mathrm{Im} (\rho)$.
We have $\underline{\partial} ^{[\underline{k}]} \in \mathrm{Hom} _{\widetilde{A}}
( \widetilde{I} ^n / \widetilde{I} ^{n+1}, \widetilde{A} )$.
For any 
$\underline{l}\in \mathrm{Im} (\rho)$, we get
$\underline{\partial} ^{[\underline{k}]}
(\underline{\tau} ^{\underline{l}} _n ) =
1$ if $\underline{l}= \underline{k}$ 
and $\underline{\partial} ^{[\underline{k}]}
(\underline{\tau} ^{\underline{l}} _n ) =
0$ if $\underline{l}\not = \underline{k}$.
For any 
$\underline{l} \in \mathrm{Im} (\rho)$, we have
$\underline{\eta} ^{\underline{l} } _n \in  \widetilde{I} ^n / \widetilde{I} ^{n+1}$ and we compute
$$\vartheta ^n \underline{\partial} ^{[\underline{k}]}
(\underline{\eta} ^{\underline{l} } _n) 
=
\sum _{\underline{j}\in \{1,\dots, d \} ^{n}}
a _{\iota (\underline{l}), \underline{j}} 
\underline{\partial} ^{[\underline{k}]}
(\underline{\tau} ^{\rho (\underline{j})} _n )
=
\sum _{\underline{j}\in \rho ^{-1} (  \underline{k})}
a _{\iota (\underline{l}), \underline{j}} \in A.
$$
Hence, since 
$\{ \underline{\eta} ^{\underline{l} } _n, 
\ ; \
\underline{l}\in \mathrm{Im} (\rho)
\}$ is a basis of 
$ \widetilde{I} ^n / \widetilde{I} ^{n+1}$ as 
$\widetilde{A} $-module for the left structure 
whose dual basis of 
$\mathrm{Hom} _{\widetilde{A}}
( \widetilde{I} ^n / \widetilde{I} ^{n+1}, \widetilde{A} )$
is 
$\{ \underline{\partial} _{\sharp}  ^{[\underline{l}]} \ ; \
\ ; \
\underline{l}\in \mathrm{Im} (\rho)
\}$
(see \ref{tildePn}),
then 
we get the equality in 
$\mathrm{Hom} _{\widetilde{A}}
( \widetilde{I} ^n / \widetilde{I} ^{n+1}, \widetilde{A} )$
of the form
\begin{equation}
\label{thetapartialDlogproof}
\vartheta ^n \underline{\partial} ^{[\underline{k}]}=
\sum _{\underline{l}\in \mathrm{Im} (\rho)}
\left ( \sum _{\underline{j}\in \rho ^{-1} (  \underline{k})}
a _{\iota (\underline{l}), \underline{j}} \right )
\underline{\partial} ^{[\underline{l}]} _\sharp \in 
\Gamma (\fP, \cD _{\fP ^\sharp/\fS ^\sharp,n}).
\end{equation}

\end{proof}

\begin{empt}
\label{ntn-theta}
Let $n \in \bbN$.
With Notation \ref{ntn-ZandZ'infty}, 
we denote by $\tau _i \in  \Gamma (\fP _\eta , \mathcal{P}  ^n  _{\fP _\eta/\eta} )$ 
the image of 
$1 \otimes t _i - t _i \otimes 1$ via the morphism
$\Gamma( \fP _\eta \times _{\eta } \fP _\eta, \cO _{\fP _\eta \times _{\eta } \fP _\eta} )  
 \to
 \Gamma (\fP _\eta , \mathcal{P}  ^n  _{\fP _\eta/\eta} )$, 
 for any $i = 1,\dots, d$.
 Set
$A _\eta := \Gamma (\fP _\eta , \cO _{\fP _\eta})$,
$B _\eta := \Gamma (\fZ _\eta , \cO _{\fZ _\eta})$.
Then 
$\{\underline{\tau} ^{\underline{i}}\; ;\; \underline{i} \in \bbN ^d \text{ such that } |\underline{i}|\leq n\} $ is an $A _\eta$-basis of 
$ \Gamma (\fP _\eta , \mathcal{P}  ^n  _{\fP _\eta/\eta} )$. 
We denote by 
$\{\underline{\partial} ^{[\underline{i}]}\; ;\; \underline{i} \in \bbN ^d \text{ such that } |\underline{i}|\leq n\} $ the associated dual $A _\eta$-basis of 
$D _{\fP _\eta,n}$, which 
yields the $A _\eta$-basis 
$\{\underline{\partial} ^{[\underline{i}]}\; ;\; \underline{i} \in \bbN ^d\} $ of
$D _{\fP _\eta}$.

Let $\overline{t} _{1},\dots, \overline{t} _e \in \Gamma (\fZ _\eta , \cO _{\fZ _\eta})$ be the images of 
$t _{1},\dots, t _e$ via the surjection 
$[\hspace{0,15cm}] 
\colon 
A _\eta \twoheadrightarrow  B _\eta$. 
We denote by 
$\overline{\tau} _i \in  \Gamma (\fZ _\eta , \mathcal{P}  ^n  _{\fZ _\eta/\eta} )$ 
the image of 
$1 \otimes \overline{t} _i - \overline{t} _i \otimes 1$ via the morphism
$\Gamma( \fZ _\eta \times _{\eta } \fZ _\eta, \cO _{\fZ _\eta \times _{\eta } \fZ _\eta} )  
 \to
 \Gamma (\fZ _\eta , \mathcal{P}  ^n  _{\fZ _\eta/\eta} )
$, for any $i = 1,\dots, e$.
We get the $B _\eta$-basis
$\{\underline{\overline{\tau}} ^{\underline{i}}\; ;\; \underline{i} \in \bbN ^e \text{ such that } |\underline{i}|\leq n\} $  of 
$ \Gamma (\fZ _\eta , \mathcal{P}  ^n  _{\fZ _\eta/\eta} )$.
To complete Notation \ref{ntn-ZandZ'},
for any 
$\underline{i}= (i _1,\dots, i _e) \in \bbN ^{e}$, 
we denote by 
$[\hspace{0,15cm}]\colon
D _{\fP _\eta, \fZ _\eta/\eta } 
\to 
D_{\fP _\eta ,\fZ _\eta /\eta } 
/
D_{\fP _\eta ,\fZ _\eta /\eta } K
\riso 
D_{\fZ _\eta /\eta } $ 
the epimorphism of rings.
For any $\underline{i} \in \bbN ^e$ such that $|\underline{i}|\leq n$, we get the commutative diagram
\begin{equation}
\label{ntn-theta-diag}
\xymatrix{
{A _\eta} 
\ar@{->>}[d] ^-{}
\ar[r] ^-{d _1}
&
{ \Gamma (\fP _\eta , \mathcal{P}  ^n  _{\fP _\eta/\eta} )} 
\ar[r] ^-{\underline{\partial} ^{[\underline{i}]}}
\ar@{->>}[d] ^-{}
& 
{A _\eta} 
\ar@{->>}[d] ^-{[\hspace{0,15cm}]}
\\ 
{B _\eta} 
\ar[r] ^-{d _1}
&
{ \Gamma (\fZ _\eta , \mathcal{P}  ^n  _{\fZ _\eta/\eta} )} 
\ar[r] ^-{[\underline{\partial} ^{[\underline{i}]}]}
& 
{B _\eta.} 
}
\end{equation}
Hence, 
$\{[\underline{\partial} ^{[\underline{i}]}]\; ;\; \underline{i} \in \bbN ^e \text{ such that } |\underline{i}|\leq n\} $ 
is the dual $B _\eta$-basis of 
$\{\underline{\overline{\tau}} ^{\underline{i}}\; ;\; \underline{i} \in \bbN ^e \text{ such that } |\underline{i}|\leq n\} $.
This justifies the notation 
$[\underline{\partial} ^{[\underline{i}]}] =\underline{\partial} ^{[\underline{i}]}$.
Moreover, using the commutativity of the diagram \ref{ntn-theta-diag}, 
we compute that for any $c \in A _\eta$, 
for any $ \underline{i} \in \bbN ^e$ we have
$[\underline{\partial} ^{[\underline{i}]} (c) ]
=
\underline{\partial} ^{[\underline{i}]} ( [c]).$
Following \ref{compa-infinitelevel} and its notation, 
denoting by $\theta$ the image of $\vartheta$ via $A \twoheadrightarrow B$, 
we get, for any $b\in B$, 
for any $\underline{i} \in \bbN ^e$, 	
\begin{equation}
\label{thetapartialDlog}
\theta ^{|\underline{i}|} \underline{\partial} ^{[\underline{i}]} (b) \in B.
\end{equation}

\end{empt}

\begin{ntn}
\label{ntncompa-infinitelevel}
We keep notation \ref{ntn-theta}.
Let $E$ be the left $B _t$-submodule of 
$\widetilde{D} _{\fZ /\fS }$ generated by 
$\{\underline{\partial} ^{[\underline{i}]} \;;\; \underline{i}\in \bbN ^e\}$.
Let $P\in E$.
We define $\ord (P)$ to be the smallest integer such that
we can write $P$ of the form
$P = \sum _{\underline{i}\in \bbN ^e, k\in \bbN} \frac{b _{\underline{i},k}}{t ^k} \underline{\partial} ^{[\underline{i}]}$
with $b _{\underline{i},k} \in B$ such that 
$k + |\underline{i}| \leq \ord (P)$.

Let $E _\theta$ be the left $B _{t\theta}$-submodule of 
$\widetilde{D} _{\fZ /\fS }$ generated by 
$\{\underline{\partial} ^{[\underline{i}]} \;;\; \underline{i}\in \bbN ^e\}$.
Let $P _\theta\in E _\theta$.
We define $\ord (P _\theta)$ to be the smallest integer such that
we can write $P _\theta$ of the form
$P _\theta= \sum _{\underline{i}\in \bbN ^e, k,l\in \bbN} \frac{b _{\underline{i},k,l}}{t ^k \theta ^l} \underline{\partial} ^{[\underline{i}]}$
with $b _{\underline{i},k,l} \in B$ such that 
$k + l+|\underline{i}| \leq \ord (P _\theta)$.

\end{ntn}

\begin{lem}
\label{lem-stabEthetaE}
With Notation \ref{ntncompa-infinitelevel}, let $P\in E$, and
$P _\theta\in E _\theta$.
Then 
$P _\theta P \in E _\theta $
and $\ord (P _\theta P) \leq \ord (P _\theta) + \ord ( P)$.
\end{lem}

\begin{proof}
By additivity, we can suppose 
$P _{\theta}= \frac{b _{1}}{t ^{k _1} \theta ^{l _1}} \underline{\partial} ^{[\underline{i} _1]}$,
$P = \frac{b _{2}}{t ^{k _2}} \underline{\partial} ^{[\underline{i} _2]}$
and also
$k _1 + l _1+|\underline{i} _1| =\ord (P _\theta)$,
$k _2 + |\underline{i} _2| =\ord (P)$.
We get
$P _{\theta}P = \frac{b _{1}}{t ^{k _1+k _2} \theta ^{l _1}} 
\underline{\partial} ^{[\underline{i} _1]} b _{2} 
\underline{\partial} ^{[\underline{i} _2]}$.
Since
$\underline{\partial} ^{[\underline{i} _1]} b _{2} 
=
\sum _{\underline{i} '_1+\underline{i} ''_1=\underline{i} _1}
\underline{\partial} ^{[\underline{i} '_1]} (b _2) 
\underline{\partial} ^{[\underline{i} ''_1]} 
=
\sum _{\underline{i} '_1+\underline{i} ''_1=\underline{i} _1}
\frac{\theta ^{|\underline{i}' _1 |}\underline{\partial} ^{[\underline{i} '_1]} (b _2)}{\theta ^{|\underline{i}' _1 |}} 
\underline{\partial} ^{[\underline{i} ''_1]} $, then
$P _{\theta}P =
\sum _{\underline{i} '_1+\underline{i} ''_1=\underline{i} _1}
\frac{b _{1} \theta ^{|\underline{i}' _1 |}\underline{\partial} ^{[\underline{i} '_1]} (b _2)}
{t ^{k _1+k _2} \theta ^{l _1+|\underline{i}' _1 |}} 
\underline{\partial} ^{[\underline{i} ''_1]}
\underline{\partial} ^{[\underline{i} _2]}$.
From \ref{thetapartialDlog}, 
$\theta ^{|\underline{i}' _1 |}\underline{\partial} ^{[\underline{i} '_1]} (b _2)
\in B$.
Hence, 
$\ord (P _{\theta}P) 
\leq \ord (P _{\theta})  + \ord (P) $.

\end{proof}

\begin{empt}
Let $\cE $ be a coherent $\widetilde{\cD} ^{(m)\dag} _{\fP ^{\sharp}/\fS ^\sharp,\bbQ}$-module. Let $\fU ^\sharp$ be the open of $\fP$ complementary to $\fX$. We say that $\cE$ has its support in $\fX$ if  $\cE | \fU= 0$.  The following properties are equivalent:
\begin{enumerate}[(a)]
\item $\cE$ has its support in $\fX$ 
\item $\cE |\fP _\eta$ has its support in $\fX _\eta$.
\end{enumerate}
Indeed, suppose $\cE |\fP _\eta$ has its support in $\fX _\eta$. Then $\G := \cE |\fU$ is coherent $\widetilde{\cD} ^{(m)\dag} _{\fU ^{\sharp}/\fS ^\sharp,\bbQ}$-module such that  $\G | \fU _\eta = 0$.
Since  $\widetilde{\cD} ^{(m)\dag} _{\fP ^\sharp/\fS^\sharp}  \to  j _*\widehat{\cD} ^{(m)} _{\fP _\eta/\eta} $ 
is faithfully flat (see \ref{faithflatX2mor}), using \ref{lemu+coh=cohetc}.(\ref{lemu+coh=cohetc(b)}) this implies that $\G =0$, i.e. $\cE$ has its support in $\fX$.  The converse is obvious. 
\end{empt}

\begin{lem} [Berthelot's key lemma]
\label{Key-Lemma}
We keep notation \ref{ntn-DZZ'}.  Let $s \geq 1$ and $1\leq j\leq e$ be two integers.
If $R \in M _s (\widetilde{D} ^{(m)\dag} _{\fZ/\fS})$ is any matrix, there exist an integer $m'\geq m$ and a matrix $P \in M _s (\widetilde{D} ^{(m')\dag} _{\fZ/\fS}) $ such that
\begin{enumerate}[(a)]
\item $P \equiv I _s \mod \pi M _s (\widetilde{D} ^{(m')\dag} _{\fZ/\fS}) $\\
\item $t _j ^{p ^{m}} P = P ( t _j^{p ^{m}} I _s - \pi R)$, where $R $ is viewed as an element of $M _s (\widetilde{D} ^{(m')\dag} _{\fZ/\fS}) $ via the inclusion $M _s (\widetilde{D} ^{(m)\dag} _{\fZ/\fS}) \subset M _s (\widetilde{D} ^{(m')\dag} _{\fZ/\fS})$.
\end{enumerate}

\end{lem}

\begin{proof}
Since the proof is identical in the general case, we can suppose $s=1$.  By symmetry, we can suppose $j =1$.

1) Notations and definitions. 
Following the equation \ref{wkcpdscpinfty}, we can write 
$R$ in the form
\begin{equation}
\label{wkcpdscpbis}
R =  \sum _{k\in \bbN}\sum _{\underline{i}\in \bbN ^d} 
 \frac{b _{\underline{i},k}}{t ^{k}} \underline{\partial} ^{[\underline{i}]}
\end{equation}
where 
$b _{\underline{i},k}\in B$ 
are such that there exists a constant $\alpha>0$ such that 
$v _{\pi} ( b _{\underline{i},k}) \geq \frac{|\underline{i}|+k}{\alpha} -1$.
We fix such a sum of the form \ref{wkcpdscpbis}.
Hence, we can put  
$[R] _l  := \pi ^{-l}  \sum _{\underline{i} \in T _l} 
  \frac{b _{\underline{i},k}}{t ^{k}}
 \underline{\partial} ^{[ \underline{i}]}$, 
 where $T _l$ is the (finite) subset $\bbN ^{d}$ of  elements
$\underline{i}$ such that $ v _{\pi } (b _{\underline{i},k})= l$.  
We set
$\sigma _l (R) := \sum _{n=0} ^{l} \pi ^{n}[R]_{n}$.
We have $\mathrm{ord} ( [R] _l) \leq \alpha (l+1)$ for any integer $l$.

2) Using the formula \cite[1.0.1.2]{Be0}, 
we compute that for any integer $N$, 
$$t _1 ^{p ^{m}} \partial _1 ^{[N + p ^m]}
-
\partial _1 ^{[N + p ^m]}t _1^{p ^{m}}
\equiv
- \partial _1 ^{[N]}
\mod \pi 
E .$$
This yields that for any operator 
$U\in E _\theta$ 
there exists 
an operator
$Q\in E _{\theta}$
such that 
$$t _1^{p ^{m}} Q
-
Qt _1^{p ^{m}}
\equiv
U \mod \pi E _{\theta}$$
and 
$\mathrm{ord} (Q) \leq \mathrm{ord} (U) + p^m$.
Indeed, we can write
$U =  \sum '_{\underline{i}\in \bbN ^d, k,l\in \bbN} \frac{b _{\underline{i},k,l}}{t ^k \theta ^l} \underline{\partial} ^{[\underline{i}]}$
with $b _{\underline{i},k,l} \in B$ such that 
$k + l+|\underline{i}| \leq \ord (U)$ and where the prime means that the sum is finite.
For any $\underline{i}\in \bbN ^d$, we set 
$\underline{\partial} ^{\prime [\underline{i}]}
:= \partial _ 2 ^{[i _2]} \cdots \partial _ d ^{[i _d]}$.
Using Lemma \ref{lem-stabEthetaE}, 
we get:
$\frac{b _{\underline{i},k,l}}{t ^k \theta ^l} \underline{\partial} ^{\prime [\underline{i}]} 
\left ( t _1 ^{p ^{m}} \partial _1 ^{[i _1 + p ^m]}
-
 \partial _1 ^{[i _1 + p ^m]}t _1^{p ^{m}} \right) 
\equiv
- \frac{b _{\underline{i},k,l}}{t ^k \theta ^l} \underline{\partial} ^{[\underline{i}]}
\mod \pi 
E _\theta$.
Since
$\underline{\partial} ^{\prime [\underline{i}]} $ 
commutes with 
$t _1 ^{p ^{m}} $ we get the equality
$\frac{b _{\underline{i},k,l}}{t ^k \theta ^l} \underline{\partial} ^{\prime [\underline{i}]} 
\left ( t _1 ^{p ^{m}} \partial _1 ^{[i _1 + p ^m]}
-
 \partial _1 ^{[i _1 + p ^m]}t _1^{p ^{m}} \right) 
=
t _1 ^{p ^{m}} ( \frac{b _{\underline{i},k,l}}{t ^k \theta ^l} \underline{\partial} ^{[\underline{i} + ( p ^m , 0,\dots, 0)]} )
-
( \frac{b _{\underline{i},k,l}}{t ^k \theta ^l} \underline{\partial} ^{[\underline{i} + ( p ^m , 0,\dots, 0)]} ) 
t _1^{p ^{m}}$.
Hence, 
we can take 
$Q= 
 \sum ' _{\underline{i}\in \bbN ^d, k,l\in \bbN} 
 \frac{b _{\underline{i},k,l}}{t ^k \theta ^l} \underline{\partial} ^{[\underline{i} + ( p ^m , 0,\dots, 0)]}$.

3) Put $\beta := \alpha + p ^m$. By induction on $l\geq 0$, we construct sections 
$Q _l \in E _\theta$ such that for any $l \geq 0$
\begin{enumerate}[(i)]
\item $Q _0 = 1$ ;
\item $\ord (Q _l) \leq \beta (l+1)$ ;
\item setting 
$P _l := \sum _{n=0} ^{l} \pi ^{l} Q _l$, we get 
$ t _1^{ p ^m} P _l \equiv P _l ( t _1^{ p ^m} - \pi  \sigma _{l-1} (R )) \mod \pi ^{l+1}  E _\theta$.
\end{enumerate}

By induction, we suppose $P _0, \dots, P _l$ constructed. 
Using Lemma \ref{lem-stabEthetaE}, for any integers $l _1$ and $l _2$, 
since $[R] _{l _2}\in E$ and  
$Q  _{l_1} \in E _{\theta}$ then 
$Q  _{l_1} [R] _{l _2} \in E _{\theta}$.
Hence, we get the following congruence of elements of
$E _\theta$ :
\begin{equation}
\notag
t _1^{ p ^m} P _l - P _l ( t _1^{ p ^m} - \pi \sigma _{l} (R ))
\equiv
t _1^{p ^m} P _l - P _l  t _1^{ p ^m}  + \pi 
\sum _{l _1 + l _2 \leq l} 
\pi ^{l _1 + l _2} Q  _{l_1} [R] _{l _2}
 \mod \pi ^{l+2}  E _\theta.
\end{equation}
By induction hypothesis we have
$t _1^{ p ^m} P _l - P _l ( t _1^{ p ^m} - \pi \sigma _{l} (R )) 
\equiv
t _1^{ p ^m} P _l - P _l ( t _1^{ p ^m} - \pi \sigma _{l-1} (R ))
\equiv 0 \mod \pi ^{l+1}  E _\theta$.
Hence, 
there exists $U \in E _\theta$ such that 
\begin{equation}
\label{3}
t _1^{p ^m} P _l - P _l  t _1^{ p ^m}  + \pi
\sum _{l _1 + l _2 \leq l} 
\pi ^{l _1 + l _2} Q  _{l_1} [R] _{l _2}
=
-\pi ^{l+1}  U .
\end{equation}
Since $\mathrm{ord} (Q  _{l_1}) \leq \beta (l _1 +1)$ (by induction hypothesis)
and $\mathrm{ord} ([R]  _{l_2}) \leq \alpha (l _2 +1)$, 
since 
$l _1 + l _2 \leq l$ 
(and $\beta \geq \alpha$),
then using Lemma \ref{lem-stabEthetaE} we check
$\mathrm{ord}
(Q  _{l_1} [R] _{l _2}) \leq \beta (l+1) + \alpha$.
Hence, 
$\mathrm{ord} ( U)  \leq \beta (l+1) + \alpha$.
Using the part 2) of the proof, 
there exists $Q _{l+1} \in E _\theta$ such that 
$U \equiv t _1^{p ^m} Q _{l+1}- Q _{l+1}t _1^{ p ^m} 
\mod \pi E _{\theta}
$ 
and 
$\mathrm{ord} ( Q) \leq \mathrm{ord} ( U) + p ^{m}$.
This yields that 
$\mathrm{ord} (Q _{l+1}) \leq \beta ( l +1)+ \alpha + p ^{m} = \beta (l +2)$.
Moreover, from \ref{3}, 
we get the equality of elements 
of $E _\theta$ of the form
\begin{align}
\notag
t _1^{p ^m} (P _l + \pi ^{l+1}  Q _{l+1}) 
& 
\overset{\ref{3}}{\equiv}
 (P _l + \pi ^{l+1}  Q _{l+1})  t _1^{ p ^m}  - \pi 
\sum _{l _1 + l _2 \leq l} 
\pi ^{l _1 + l _2} Q  _{l_1} [R] _{l _2} 
\mod \pi ^{l+2}  E _\theta
\\
\notag
&
\overset{\ref{lem-stabEthetaE}}{\equiv}
(P _l + \pi ^{l+1}  Q _{l+1}) (  t _1^{ p ^m} - \pi \sigma _l ( R) )
\mod \pi ^{l+2}  E _\theta.
\end{align}
Hence, putting $P _{l+1}:=P _l + \pi ^{l+1}  Q _{l+1} $,
we obtain
$ t _1^{ p ^m} P _{l+1} \equiv P _{l+1} ( t _1^{ p ^m} - \pi  \sigma _{l} (R )) \mod \pi ^{l+2}  E _\theta$.

4) Finally, using the part 3) and the remark of \ref{empt-wkcpdscptheta}, we get the element $P := \lim _{l\to \infty} \, P _l \in \widetilde{D} ^{\dag} _{\fZ/\fS}$.
Since 
$\widetilde{D} ^{\dag} _{\fZ/\fS} = \cup _{m'\in \bbN} \widetilde{D} ^{(m')\dag} _{\fZ/\fS}$, then for $m'$ large enough, 
we get 
$P \in \widetilde{D} ^{(m')\dag} _{\fZ/\fS}$.
This operator $P$ satisfies the required properties.
\end{proof}

\begin{lem}
\label{lem-Berthelot-Kashiwara}
We keep notation \ref{ntn-DZZ'} and \ref{u+coh=cohetc-lem}.  We suppose $e =e' +1$.
Let $E$ be a $\widetilde{D} ^{(m)\dag} _{\fZ/\fS ,\bbQ}$-module of finite type such that 
$\widehat{\cD} ^{(m)} _{\fZ _\eta/\eta,\bbQ}  \otimes _{\widetilde{D} ^{(m)\dag} _{\fZ/\fS ,\bbQ}} E$
is a left $\widehat{\cD} ^{(m)} _{\fZ _\eta/\eta,\bbQ} $-module of finite type 
with support in $\fX _\eta$.
Then there exist an integer $m' \geq m$, a $\widetilde{D} ^{(m)\dag} _{\fZ'/\fS ,\bbQ}$-module of finite type $E'$ 
such that 
$\widehat{\cD} ^{(m')} _{\fZ '_\eta/\eta,\bbQ}  \otimes _{\widetilde{D} ^{(m')\dag} _{\fZ '/\fS ,\bbQ}} E'$
is a left $\widehat{\cD} ^{(m')} _{\fZ '_\eta/\eta} $-module of finite type
with support in $\fX _\eta$,
and 
an isomorphism of $\widetilde{D} ^{(m')\dag} _{\fZ/\fS ,\bbQ}$-modules of the form
\begin{equation}
\notag
\iota _{+} ^{(m')\dag} (E') \riso 
\widetilde{D} ^{(m')\dag} _{\fZ/\fS ,\bbQ}\otimes _{\widetilde{D} ^{(m)\dag} _{\fZ/\fS ,\bbQ}} E.
\end{equation}

\end{lem}

\begin{proof}
We use from now notation \ref{ntncompa-infinitelevel}, 
for instance that concerning the definition of 
$\theta$ and $E _\theta$.

a) For any integer 
$m' \geq m$, 
we set
$E ^{(m')}:= \widetilde{D} ^{(m')\dag} _{\fZ/\fS,\bbQ}\otimes _{\widetilde{D} ^{(m)\dag} _{\fZ/\fS,\bbQ}} E$. 
We denote by 
$\mathcal{H} ^{0} \iota ^{ !} (\smash{{E}} ^{(m')})$ 
the set of elements 
$\smash{{E}} ^{(m')}$
which are annihilated by $t _e$. 
We check in this part that for $m'$ large enough, 
$\smash{{E}} ^{(m')}$
is 
generated 
as $\widetilde{D} ^{(m')\dag} _{\fZ/\fS ,\bbQ} $-module
by $\mathcal{H} ^{0} \iota ^{ !} (\smash{{E}} ^{(m')})$.

i)  Let 
$\smash{{K}} ^{(m')}$ be
the $\widetilde{D} ^{(m')\dag} _{\fZ/\fS,\bbQ } $-submodule
of $\smash{{E}} ^{(m')}$
generated by $\mathcal{H} ^{0} \iota ^{ !} (\smash{{E}} ^{(m')})$.
Let $i \geq 1$ be an integer, and $x\in \smash{{E}} ^{(m')}$. 
By induction on $i$, 
we check that 
if $t  _e ^{i} x= 0$ then $x \in \smash{{K}} ^{(m')}$.

When $i=1$, this is obvious. Suppose $i \geq 2$ and the property true for $i-1$. 
Suppose $t  _e ^{i} x= 0$.
Set $y:=( \partial _e t _e+ (i-1)) \cdot x$.
From $\partial _e t _e^{i-1} = t _e^{i-1} \partial _e + (i-1)t _e^{i-2}$
we obtain 
$\partial _e t _e^{i} = t _e^{i-1} \partial _e t _e+ (i-1)t _e^{i-1}
=t _e^{i-1} ( \partial _e t _e+ (i-1))$.
Since $\partial _e t _e^{i} \cdot x =0$, 
we obtain $t _e^{i-1} \cdot y = 0$. 
Hence, by using the induction hypothesis
we get that 
$y = ( \partial _e t _e+ (i-1)) \cdot x \in \smash{{K}} ^{(m')}$.
Again, by using the induction hypothesis 
$t _e\cdot x \in \smash{{K}} ^{(m')}$ and then $\partial _e t _e\cdot x \in \smash{{K}} ^{(m')}$.
Hence,  $(i-1) x \in \smash{{K}} ^{(m')}$ and then $x\in \smash{{K}} ^{(m')}$ and the induction is proved.

ii) Choose  some generators $e _1, \dots, e _s$ of $E$. We denote by 
$\smash{\overset{\circ}{E}} $ the sub-$\widetilde{D} ^{(m)\dag} _{\fZ/\fS}$-module of $E$ 
generated by $e _1, \dots, e _s$.
We denote by $\smash{\overset{\circ}{E}} ^{(m')}$ the quotient of 
$\widetilde{D} ^{(m')\dag} _{\fZ/\fS } \otimes _{\widetilde{D} ^{(m)\dag} _{\fZ/\fS }} \smash{\overset{\circ}{E}} $
by its $\pi$-torsion elements. 
Let 
$e ' _1, \dots, 1e '_s$ be the image of $e _1, \dots, e _s$ 
via 
$\smash{\overset{\circ}{E}} \to \smash{\overset{\circ}{E}} ^{(m')}$.
The $\widetilde{D} ^{(m')\dag} _{\fZ/\fS}$-module $\smash{\overset{\circ}{E}} ^{(m')}$ is $p$-torsion-free of finite type (generated by $e '_1, \dots,e' _s$) and
such that $\smash{\overset{\circ}{E}} ^{(m')} _\bbQ \riso E^{(m')}$.

By flatness we get the injective morphism
$\widehat{\cD} ^{(m)} _{\fZ _\eta/\eta}  \otimes _{\widetilde{D} ^{(m)\dag} _{\fZ/\fS }} \smash{\overset{\circ}{E}} 
\hookrightarrow
\widehat{\cD} ^{(m)} _{\fZ _\eta/\eta}  \otimes _{\widetilde{D} ^{(m)\dag} _{\fZ/\fS}} E  
\riso 
\widehat{\cD} ^{(m)} _{\fZ _\eta/\eta,\bbQ}  \otimes _{\widetilde{D} ^{(m)\dag} _{\fZ/\fS ,\bbQ}} E  $.
Since 
$\widehat{\cD} ^{(m)} _{\fZ _\eta/\eta,\bbQ}  \otimes _{\widetilde{D} ^{(m)\dag} _{\fZ/\fS ,\bbQ}} E  $
is a left $\widehat{\cD} ^{(m)} _{\fZ _\eta/\eta} $-module 
with support in $\fX _\eta$, 
then 
$\widehat{\cD} ^{(m)} _{\fZ _\eta/\eta}  \otimes _{\widetilde{D} ^{(m)\dag} _{\fZ/\fS }} \smash{\overset{\circ}{E}} $
has in support in $\fX _\eta$.
Since the morphisms
$\widetilde{D} ^{(m)\dag} _{\fZ/\fS} / \pi \widetilde{D} ^{(m)\dag} _{\fZ/\fS}
\to
\widehat{D} ^{(m)} _{\fZ _\eta /\eta} / \pi \widehat{D} ^{(m)} _{\fZ/\fS}
\riso 
D ^{(m)} _{Z _\eta /\eta}:= \Gamma ( Z _\eta, \cD ^{(m)} _{Z _\eta /\eta})$
are isomorphisms, this yields that 
$\widehat{\cD} ^{(m)} _{\fZ _\eta/\eta}  \otimes _{\widetilde{D} ^{(m)\dag} _{\fZ/\fS }} \smash{\overset{\circ}{E}} 
/
\pi ( \widehat{\cD} ^{(m)} _{\fZ _\eta/\eta}  \otimes _{\widetilde{D} ^{(m)\dag} _{\fZ/\fS }} \smash{\overset{\circ}{E}} )
\riso 
\cD ^{(m)} _{Z _\eta /\eta} \otimes _{D ^{(m)} _{Z _\eta /\eta}}(\smash{\overset{\circ}{E}}  / \pi \smash{\overset{\circ}{E}} )$. 
has his support in $X _\eta$. 
This implies that $ t _e^{p ^{m'} } e _i \equiv 0 \mod \pi \smash{\overset{\circ}{E}} $ (increase $m'$ if necessary). 
This yields 
$ t _e ^{p ^{m'} } e'  _i \equiv 0 \mod \pi \smash{\overset{\circ}{E}} ^{(m')}$.
Put $\underline{e}':= ^{t} (e' _1,\dots, e '_s)$.
Hence, there exists a matrix $R \in M _n (\widetilde{D} ^{(m')\dag} _{\fZ/\fS})$
such that 
$t _e^{p ^{m'}} \underline{e} '=
\pi 
R \underline{e}'$.

For $m''\geq m'$ be large enough, let $P\in M _s (\widetilde{D} ^{(m'')\dag} _{\fZ/\fS})$ be an element satisfying the condition of  Lemma \ref{Key-Lemma}
(in the case where $m$ is replaced by $m'$ and $m'$ by $m''$).
We put $\underline{e}'': = P (1 \otimes \underline{e}')\in M _{d1} (\smash{\overset{\circ}{E}} ^{(m'')})$, where $1 \otimes \underline{e}'$ means the column vector whose elements are the image
of $e ' _i$ via 
$\smash{\overset{\circ}{E}} ^{(m')} \to \smash{\overset{\circ}{E}} ^{(m'')}$.
Then 
$t _e^{p ^{m'}} \underline{e}'' = 
t _e^{p ^{m'}} P (1 \otimes \underline{e}')
=
P ( t _e^{p ^{m'}} I _s - \pi R)
(1 \otimes \underline{e}')=
P ( 1 \otimes (t _e^{p ^{m'}} \underline{e}' -
\pi 
R \underline{e}'))=0 $. 
Hence, following part i), 
we get
$e ''_1, \dots, e ''_s \in K ^{(m'')}$.

iii) We prove in this step that $e ''_1, \dots, e ''_s$ generate $E ^{(m'')}$
as $\widetilde{D} ^{(m')\dag} _{\fZ/\fS,\bbQ}$-module.
It is sufficient to check that 
$e ''_1, \dots, e ''_s$ generate $\smash{\overset{\circ}{E}} ^{(m'')}$
as $\widetilde{D} ^{(m')\dag} _{\fZ/\fS}$-module.
Using \ref{u+coh=cohetc-lem}.5, 
this latter property is equivalent to the fact 
that $e ''_1, \dots, e ''_s$ generate $\widehat{D} ^{(m)} _{\fZ _\eta/\eta} \otimes _{\widetilde{D} ^{(m)\dag} _{\fZ /\fS}} \smash{\overset{\circ}{E}} ^{(m'')}$
as $\widehat{D} ^{(m)} _{\fZ _\eta/\eta} $-module.
Since $\widehat{D} ^{(m)} _{\fZ _\eta/\eta} \otimes _{\widetilde{D} ^{(m)\dag} _{\fZ /\fS}} \smash{\overset{\circ}{E}} ^{(m'')}$
is the $p$-adic completion of 
$\smash{\overset{\circ}{E}} ^{(m'')}$, 
we conclude using the condition $1$ satisfied by $P$ in Lemma \ref{Key-Lemma}
and \cite[3.2.2.(ii)]{Be1}.

b) According to notation of \ref{u+coh=cohetc-cor},  we get $\iota ^{(m')\dag} _{+}  :=  \widetilde{D} ^{(m')\dag}  _{\fZ  \leftarrow \fZ ^{\prime}/\fS ,\bbQ}  \otimes  _{\widetilde{D} ^{(m)\dag} _{\fZ '/\fS,\bbQ }}-$,
a functor defined over the category of left $\widetilde{D} ^{(m')\dag} _{\fZ'/\fS,\bbQ}$-module. 
This yields the morphism
\begin{equation} \label{u_+u^!toif-limproj}
\mathrm{adj}\,:\, \iota ^{(m')\dag} _{+}\circ \mathcal{H} ^{0} \iota ^{ !} (\smash{{E}} ^{(m')}) \to  \smash{{E}} ^{(m')}
\end{equation}
given by $[P] \otimes x\mapsto P\cdot x$ (see notation of \ref{ntn-DZZ'[]}), for any $P \in \widetilde{D} ^{(m')\dag} _{\fZ/\fS,\bbQ}$, and $x \in \mathcal{H} ^{0} \iota ^{ !} (\smash{{E}} ^{(m')})$ (we notice that this is well defined).

From part a) and by noetherianity, we fix $m'$ so that $\smash{{E}} ^{(m')}$ is generated 
as $\widetilde{D} ^{(m')\dag} _{\fZ/\fS ,\bbQ} $-module
by a finite set of elements of  $\mathcal{H} ^{0} \iota ^{ !} (\smash{{E}} ^{(m')})$.
This finite set of elements corresponds to a morphism
$L ^{(m')} \to \mathcal{H} ^{0} \iota ^{ !} (\smash{{E}} ^{(m')})$
of 
$\widetilde{D} ^{(m')\dag} _{\fZ'/\fS,\bbQ}$-module, 
where $L ^{(m')}$ is a finite free 
$\widetilde{D} ^{(m')\dag} _{\fZ'/\fS,\bbQ}$-module. 
The composition of the homomorphisms of  $\widetilde{D} ^{(m')\dag} _{\fZ/\fS,\bbQ}$-modules $\iota ^{(m')\dag} _{+} ( L ^{(m')}) \to  \iota ^{(m')\dag} _{+} \circ \mathcal{H} ^{0} \iota ^{ !} (\smash{{E}} ^{(m')}) \to  \smash{{E}} ^{(m')}$
is surjective. Let $N ^{(m')}$ be its kernel. For $m'' \geq m'$, set 
$ L ^{(m'')}  :=  \widetilde{D} ^{(m'')\dag} _{\fZ '/\fS,\bbQ }  \otimes _{\widetilde{D} ^{(m')\dag} _{\fZ '/\fS,\bbQ }}  L ^{(m')}$,
and $N ^{(m'')}  :=  \widetilde{D} ^{(m'')\dag} _{\fZ /\fS,\bbQ }  \otimes _{\widetilde{D} ^{(m')\dag} _{\fZ /\fS,\bbQ }} N ^{(m')}$.
Using  \ref{u+coh=cohetc-cor-iso8}, we get the exactness of the top sequence
\begin{equation}
\label{lem-Berthelot-Kashiwara-diag1}
\xymatrix{
{0} 
\ar[r] ^-{}
& 
{N ^{(m'')}  } 
\ar[r] ^-{}
& 
{\iota ^{(m'')\dag} _{+} ( L ^{(m'')}) } 
\ar[r] ^-{}
& 
{\smash{{E}} ^{(m'')}} 
\ar[r] ^-{}
& 
{0} 
\\ 
& 
{\iota ^{(m'')\dag} _{+}\circ \mathcal{H} ^{0} \iota ^{ !} (N ^{(m'')}  )} 
\ar[r] ^-{}
\ar[u] ^-{\mathrm{adj}}
& 
{\iota ^{(m'')\dag} _{+}\circ \mathcal{H} ^{0} \iota ^{ !} (\iota ^{(m')\dag} _{+} ( L ^{(m'')}) )} 
\ar[u] ^-{\mathrm{adj}}
&
{ } 
& 
}
\end{equation}
Since the composition 
$\iota ^{(m'')\dag} _{+} ( L ^{(m'')}) 
\overset{\ref{u+coh=cohetc-cor}.\ref{noptors-isoQ}}{\riso}
 \iota ^{(m'')\dag} _{+}\circ \mathcal{H} ^{0} \iota ^{ !} (\iota ^{(m'')\dag} _{+} ( L ^{(m'')})  )
 \overset{\mathrm{adj}}{\longrightarrow}
\iota ^{(m'')\dag} _{+} ( L ^{(m'')}) $
is the identity,  we get the latter morphism is an isomorphism. For $m''$ large enough, the kernel $N ^{(m'')}$ of this map is generated
as $\widetilde{D} ^{(m'')\dag} _{\fZ/\fS,\bbQ} $-module by a finite set of elements of  $\mathcal{H} ^{0} \iota ^{ !} (N ^{(m'')})$.
In  particular,  $\iota ^{(m'')\dag} _{+}\mathcal{H} ^{0} \iota ^{ !} (N ^{(m'')}) \to N ^{(m'')}$ is surjective. It follows from \ref{u+coh=cohetc-cor}.\ref{noptors-isoQ} that $\mathcal{H} ^{0} \iota ^{ !} (N ^{(m'')})$ is a 
submodule of $L ^{(m'')}\riso \mathcal{H} ^{0} \iota ^{ !} \iota ^{(m'')\dag} _{+} ( L ^{(m'')}) $.
Hence, 
$\mathcal{H} ^{0} \iota ^{ !} (N ^{(m'')})$ is a $\widetilde{D} ^{(m')\dag} _{\fZ'/\fS,\bbQ}$-module of finite type.
Since the functor $\iota ^{(m')\dag} _{+}$ is exact (see \ref{cor-flat-u+preZZ'}) and 
$\mathcal{H} ^{0} \iota ^{ !}$ is left exact, then the bottom arrow of  the diagram  \ref{lem-Berthelot-Kashiwara-diag1}
is injective. Hence, so is 
$\iota ^{(m'')\dag} _{+}\mathcal{H} ^{0} \iota ^{ !} (N ^{(m'')}) \to N ^{(m'')}$ (use the commutativity of the diagram  \ref{lem-Berthelot-Kashiwara-diag1}).
Since it is also surjective, it is an isomorphism.
Hence, by exactness of $\iota ^{(m')\dag} _{+}$, we obtain from the diagram  \ref{lem-Berthelot-Kashiwara-diag1} the isomorphism
$E ^{(m')} \riso \iota ^{(m')\dag} _{+} ( L ^{(m'')}/\mathcal{H} ^{0} \iota ^{ !} (N ^{(m'')}) )$.
Put 
$\G: =\widehat{\cD} ^{(m)} _{\fZ _\eta/\eta,\bbQ}  \otimes _{\widetilde{D} ^{(m)\dag} _{\fZ/\fS ,\bbQ}} E ^{(m')}$, 
$\G' := \widehat{\cD} ^{(m)} _{\fZ '_\eta/\eta,\bbQ}  \otimes _{\widetilde{D} ^{(m)\dag} _{\fZ'/\fS ,\bbQ}} 
 L ^{(m'')}/\mathcal{H} ^{0} \iota ^{ !} (N ^{(m'')}) $.
 It remains to check that $\G'$ has its support  in $\fX _\eta$.
Using \ref{u+coh=cohetc-cor}.\ref{u+coh=cohetc-cor2)}, we get
$\iota ^{(m)} _{\eta +} (\G ') 
=
\widehat{\cD} ^{(m)}  _{\fZ _\eta \leftarrow \fZ ^{\prime} _\eta/\eta,\bbQ}
 \otimes _{\widehat{\cD} ^{(m)} _{\fZ ^{\prime} _\eta/\eta,\bbQ}}
\G '\riso \G$.
Since $\G$ has its support in $\fX _\eta$, then so is 
$\mathcal{H} ^{0} \iota ^{ !} (\G)$.
Since $\mathcal{H} ^{0} \iota ^{ !}  \iota ^{(m)} _{\eta +} (\G ') \riso \G '$, we conclude.
\end{proof}

The following theorem is the analogue of Berthelot's one in the standard case (for instance, see the version
\cite[A.6]{caro-stab-sys-ind-surcoh}).
\begin{thm}
\label{Berthelot-Kashiwara}
Let $\cE $ be a coherent $\widetilde{\cD} ^{(m)\dag} _{\fP ^{\sharp}/\fS ^\sharp,\bbQ}$-module
with support in $\fX$.
Then there exist an integer $m' \geq m$, a coherent $\widetilde{\cD} ^{(m')\dag} _{\fX^\sharp/\fS ^\sharp,\bbQ}$-module $\cF$ and 
an isomorphism of $\widetilde{\cD} ^{(m')\dag} _{\fP ^{\sharp}/\fS ^\sharp,\bbQ}$-modules of the form
\begin{equation}
\notag
u _{+} ^{(m')\dag} (\cF) \riso 
\widetilde{\cD} ^{(m')\dag} _{\fP ^{\sharp}/\fS ^\sharp,\bbQ} \otimes _{\widetilde{\cD} ^{(m)\dag} _{\fP ^{\sharp}/\fS ^\sharp,\bbQ}} \E.
\end{equation}

\end{thm}

\begin{proof}
Since the functor $u _+$ is fully faithful (see \ref{u_+faithful}), 
then the theorem is local.
Hence, following Lemma \ref{ovcv-coord-comp-imm}, 
we can suppose that $\fP$ is affine and 
there exists
sections $t _1,\dots, t _d\in \Gamma (\fP ,\cO _{\fP})$ such
that $t _{1},\dots, t _d$ are overconvergent local coordinates   of $\fP ^{\sharp}/\fS ^\sharp$, 
$t _{1},\dots, t _{e} \in \Gamma (\fP ,\cI)$ is a regular sequence of $\Gamma (\fP ,\cO _{\fP})$ generating $\Gamma (\fP,\cI)$
and 
$\overline{t} _{e+1},\dots, \overline{t} _d$ are 
overconvergent local coordinates of $\fX \cap \fP/\fS ^\sharp$, where 
$ \overline{t} _{e+1},\dots, \overline{t} _d$ are respectively the image 
of $t _{e+1},\dots, t _d$ in $\Gamma (\fX,\cO _{\fX})$.
We conclude using iteratively Lemma \ref{lem-Berthelot-Kashiwara}.
\end{proof}

\begin{coro}
[Berthelot-Kashiwara]
\label{Berthelot-Kashiwara-full}
We keep notation
\ref{Berthelot-Kashiwara}. 
The functors $u ^{!}$ and $u _{+}$ induce quasi-inverse equivalences of categories of 
coherent $\widetilde{\cD} ^{\dag} _{\fP ^{\sharp}/\fS ^\sharp,\bbQ}$-modules with support in $X$ 
and that of coherent $\widetilde{\cD} ^{\dag} _{\fX ^{\sharp}/\fS ^\sharp,\bbQ}$-modules.
These functors $u ^{!}$ and $u  _{+}$ are acyclic over these categories. 
\end{coro}

\begin{proof}
Replacing \cite[A.6]{caro-stab-sys-ind-surcoh} by Theorem \ref{Berthelot-Kashiwara},
we can follow the proof of \cite[A.8]{caro-stab-sys-ind-surcoh}. 
Recall $\widetilde{\cD} ^{\dag} _{\fX ^{\sharp}/\fS ^\sharp,\bbQ}=\underset{\underset{m}{\longrightarrow}}{\lim}\,
\widetilde{\cD} ^{(m)\dag} _{\fX ^{\sharp}/\fS ^\sharp,\bbQ}$
and
$\widetilde{\cD} ^{\dag} _{\fP ^{\sharp}/\fS ^\sharp,\bbQ}=\underset{\underset{m}{\longrightarrow}}{\lim}\,
\widetilde{\cD} ^{(m)\dag} _{\fP ^{\sharp}/\fS ^\sharp,\bbQ}$.
Hence, 
using \ref{Berthelot-Kashiwara} and \ref{u_+faithful},
we check that 
the functors  $\mathcal{H} ^{0} u ^{!} $ and $u _+  $ induce quasi-inverse equivalences between the category of  coherent $\widetilde{\cD} ^{\dag} _{\fP ^{\sharp}/\fS ^\sharp,\bbQ}$-modules with support in $X$ 
and that of coherent $\widetilde{\cD} ^{\dag} _{\fX ^{\sharp}/\fS ^\sharp,\bbQ}$-modules. The acyclicity of $u _+$ is already know (see \ref{u_+faithful}).
It remains to check the acyclicity of $u ^!$.  Since this is local, we reduce to the geometrical situation of the beginning of the proof
of \ref{Berthelot-Kashiwara} and we keep its notation.  Let $\cF ^{(m)}$ be a coherent $\widetilde{\cD} ^{(m)\dag} _{\fX ^{\sharp}/\fS ^\sharp,\bbQ}$.
Set  $\cF ^{(m+1)}:=  \widetilde{\cD} ^{(m)\dag} _{\fX ^{\sharp}/\fS ^\sharp} \otimes _{\widetilde{\cD} ^{(m)\dag} _{\fX ^{\sharp}/\fS ^\sharp}}\cF ^{(m)} $.
We have the morphism:
\begin{gather} \notag u ^{(m)\dag} _{+} (\cF ^{(m)})
:= u _* (\widetilde{\cD} ^{(m)\dag}  _{\fP^\sharp \leftarrow \fX ^\sharp/\fS^\sharp} \otimes _{\widetilde{\cD} ^{(m)\dag} _{\fX ^{\sharp}/\fS ^\sharp}}\cF ^{(m)})
= 
u _* (\widetilde{\cD} ^{(m)\dag} _{\fP ^{\sharp}/\fS ^\sharp,\bbQ} /\widetilde{\cD} ^{(m)\dag} _{\fP ^{\sharp}/\fS ^\sharp,\bbQ} \cI
\otimes _{\widetilde{\cD} ^{(m)\dag} _{\fX ^{\sharp}/\fS ^\sharp}}\cF ^{(m)})
\\ \notag  \to 
u _* ( \widetilde{\cD} ^{(m+1)\dag} _{\fP ^{\sharp}/\fS ^\sharp,\bbQ} /\widetilde{\cD} ^{(m+1)\dag} _{\fP ^{\sharp}/\fS ^\sharp,\bbQ} \cI
\otimes _{\widetilde{\cD} ^{(m)\dag} _{\fX ^{\sharp}/\fS ^\sharp}}\cF ^{(m)}) 
=
u _* ( \widetilde{\cD} ^{(m+1)\dag} _{\fP ^{\sharp}/\fS ^\sharp,\bbQ} /\widetilde{\cD} ^{(m+1)\dag} _{\fP ^{\sharp}/\fS ^\sharp,\bbQ} \cI
\otimes _{\widetilde{\cD} ^{(m+1)\dag} _{\fX ^{\sharp}/\fS ^\sharp}}\cF ^{(m+1)})
=:
u ^{(m+1)\dag} _{+} (\cF ^{(m+1)})
\end{gather}

It is sufficient to check that this induced map 
$$u _{+} ^{ (m)} (\cF ^{(m)}) / t _i u _{+} ^{ (m)} (\cF ^{(m)})
\to 
u _{+} ^{ (m+1)} (\cF ^{(m+1)})/t _i u _{+} ^{ (m+1)} (\cF ^{(m+1)})$$
is the zero morphism for any $i =1,\dots, e$.
Since $\partial _i ^{< k+1> _{(m+1)}} t _i
=
t _i \partial _i ^{< k+1> _{(m+1)}} 
+
\left \{
\begin{smallmatrix}
 k+1\\
 1
\end{smallmatrix}
\right\}
_{(m+1)}
\partial _i ^{<k> _{(m+1)}}$,
then 
$$\left \{
\begin{smallmatrix}
 k+1\\
 1
\end{smallmatrix}
\right\}
_{(m+1)}
\partial _i ^{<k> _{(m+1)}}
\equiv
-t _i \partial _i ^{< k+1> _{(m+1)}} 
\mod 
\widetilde{\cD} ^{(m+1)\dag} _{\fP ^{\sharp}/\fS ^\sharp,\bbQ} /\widetilde{\cD} ^{(m+1)\dag} _{\fP ^{\sharp}/\fS ^\sharp,\bbQ} \cI.$$
Since 
$\left \{
\begin{smallmatrix}
 k+1\\
 1
\end{smallmatrix}
\right\}
_{(m+1)}
=
q _{k+1} ^{(m+1)} !/q _k ^{(m+1)} ! $, we get
$\partial _i ^{<k> _{(m)}} 
=
\frac{q _k ^{(m)} !}{q _{k} ^{(m+1)} !}  \partial _i ^{<k> _{(m+1)}}
=
\frac{q _k ^{(m)} !}{q _{k+1} ^{(m+1)} !} 
\left \{
\begin{smallmatrix}
 k+1\\
 1
\end{smallmatrix}
\right\}
_{(m+1)}  \partial _i ^{<k> _{(m+1)}}
$.
Since 
$\frac{q _k ^{(m)} !}{q _{k+1} ^{(m+1)} !}  \to 0$ when $k \to \infty$,
then 
setting 
$\cF ^{(m+1)}:= \smash{\widetilde{\cD}} ^{(m+1)}  _{\fZ/\T,\bbQ} \otimes _{\smash{\widetilde{\cD}} ^{(m)}  _{\fZ/\T,\bbQ}}\cF ^{(m)}$,
we compute that the morphism do is the zero morphism.
\end{proof}

\subsection{Coherent arithmetic $\cD$-modules over a realizable strictly semistable log scheme}

\subsubsection{Glueing isomorphisms}

\begin{prop}
\label{2.1.5Be2}
Let $f,g,h\colon \fX ^\sharp \to \fY ^\sharp$ be three log smooth morphisms of strictly semi-stable log formal schemes over $\fS ^\sharp$ 
such that $f _0=  g _0=h _0$. 
Let $\cF$ be a coherent left $\smash{\cD} ^{\dag} _{\fY ^{\sharp} } (\hdag Y _s) _{ \bbQ}$-module.
There exists a canonical isomorphism of 
coherent $\smash{\cD} ^{\dag} _{\fX ^{\sharp} }(\hdag X _{s}) _{ \bbQ}$-modules of the form
\begin{equation}
\notag
\tau _{f,g}\colon 
g  ^! \cF \riso f ^{!} \cF,
\end{equation}
functorial in $\cF$ and 
such that $\tau _{f,f}=\mathrm{Id}$
$\tau _{f,h}= \tau _{f,g} \circ \tau _{g,h}$.
\end{prop}

\begin{proof}
See \cite[9.2.2.3]{Car25}.
\end{proof}

\begin{prop} \label{comp-comp-adj-immf}
  Consider the following diagram of strictly semi-stable log formal schemes over $\fS ^\sharp$:
\begin{equation}\label{deuxcarresadj}
  \xymatrix  @R=0,3cm {
  {\fP ^{\prime \prime \sharp }} \ar[r] ^g   &    {\fP ^{\prime\sharp }} \ar[r] ^f    &    {\fP ^{\sharp}}
  \\   {\fX ^{\prime \prime \sharp}} \ar[u]  ^{u''}\ar[r]^{\rmb}   &   {\fX ^{\prime \sharp}} \ar[u] ^-{u'} \ar[r]^a   &   {\fX ^{\sharp} ,} \ar[u] ^u    }
\end{equation}
where $f$, $g$, $a$ and $b$ are  log-smooth, where $u$, $u'$ and $u''$ are  weakly closed immersions (see Definition \ref{dfn-wclimm}). We suppose that the diagram \ref{deuxcarresadj} is commutative modulo $\pi$. 

\begin{enumerate}  [(a)]
 \item We have  the canonical adjunction morphisms
\begin{equation} \label{comp-comp-adj-immf-morp1}
u ^\prime _+ \circ a ^!\to f ^{!}\circ u _+ , \qquad  u ^{\prime (\bullet)} _+ \circ a ^{(\bullet)!} \to f ^{(\bullet)!}\circ u ^{(\bullet)} _+ 
\end{equation}
of functors  of respectively of the form  $D ^{\rmb} _{\coh} ( \widetilde{\cD}  ^{\dag}  _{\fX ^{\sharp}/ \fS ^{\sharp},\bbQ})   \to   D ^{\rmb} _{\coh} ( \widetilde{\cD}  ^{\dag}  _{\fP ^{\prime\sharp}/ \fS ^{\sharp},\bbQ})$
and $\underrightarrow{LD}  ^\rmb _{\bbQ, \coh} ({} ^{\rml} \widetilde{\cD} _{\fX ^{\sharp}/ \fS ^{\sharp}} ^{(\bullet)} (X _s)) \to \underrightarrow{LD}  ^\rmb _{\bbQ, \coh} ({} ^{\rml} \widetilde{\cD} _{\fP ^{\prime \sharp}/ \fS ^{\sharp}} ^{(\bullet)} (P' _s))$).
If the right square of \ref{deuxcarresadj} is cartesian modulo $\pi$ then both morphisms \ref{comp-comp-adj-immf-morp1}  are isomorphisms.

\item Denoting by $\phi \colon   u ^\prime _+ \circ a ^! \to f ^{!}\circ u _+$, (resp. $\phi '\colon  u ^{''} _+ \circ b ^! \to g ^! \circ u ^\prime _+$, 
resp. $\phi''\colon  u ^{''} _+ \circ (a \circ b ) ^! \to (f \circ g) ^!\circ u  _+$) the morphism of adjunction of the right square  \ref{deuxcarresadj} (resp. the left square, resp. the outer of \ref{deuxcarresadj}), 
then the following diagram
$$\xymatrix  @R=0,3cm { { u''  _+ \circ (a\circ b) ^!} \ar[r]_\sim \ar[d] ^-{\phi''} & { u''  _+ \circ b ^!\circ  a ^! } \ar[d]^{( g ^! \circ \phi)\circ (\phi ' \circ a ^!)} \\ {(f\circ g) ^! \circ  u _+  } \ar[r]_-{\sim} & { g ^!  \circ f ^!\circ  u _+,} }  $$
is commutative ; and similarly in the context of categories of the form $\underrightarrow{LD}  ^\rmb _{\bbQ, \coh}$ instead of $D ^{\rmb} _{\coh}$. By abuse of notation, we get the transitivity equality $\phi''=( g ^! \circ \phi)\circ (\phi ' \circ a ^!)  $.

\item Let $a ' \colon   \fX ^{\prime } \to \fX $ (resp. $f' \colon  \fX ^{\prime } \to \fX $)
be a morphism whose reduction $X ^{\prime } \to \fX $ (resp. $P ^{\prime } \to \fX $) is equal to that of $a$ (resp. $f$). Then the following diagram
$$\xymatrix  @R=0,3cm { { u ^\prime_+ a ^{ !}} \ar[r] ^-{\phi} & {f ^{!}\circ u _+} \\ { u' _+ a ^{\prime !} } \ar[r]^{\psi} \ar[u] ^-{ u ^\prime_+ (\tau _{a,a'})} _-{\sim} & {f ^{\prime !}\circ u  _+,} \ar[u] ^-{\tau _{f,f'} u _+} _-{\sim} }$$
where $\psi$ means the morphism of adjunction of the right square of  \ref{deuxcarresadj} whose $a$ and  $f$ have been replaced respectively by $a'$ and  $f'$, is commutative; 
and similarly for categories of the form  $\underrightarrow{LD}  ^\rmb _{\bbQ, \coh}$.
\end{enumerate}
\end{prop}

\begin{proof}
Thank to \ref{Berthelot-Kashiwara-full}, the check is identical to that of \cite[2.2.2]{caro-construction} or \cite[9.3.6.1]{Car25}.
\end{proof}

\subsubsection{Berthelot-Kashiwara's theorem for weakly closed immersions of realizable semistable log schemes over $S ^\sharp$}

Let $\fP ^\sharp$ be a separated strictly semi-stable log formal scheme over $\fS ^\sharp$.
Let $u _0\colon X ^\sharp \to P ^\sharp$ be a weakly closed immersion of strictly semistable log schemes over $S ^\sharp$ (see Definition \ref{dfn-wclimm}). 
To simplify notation, put $Y:= X _\eta$. 

\begin{empt}
\label{ntnPPalpha}
Let $(\fP ^{\sharp} _{\alpha}) _{\alpha \in \Lambda}$ be an open covering of  $\fP ^{\sharp}$.
We set $\fP ^{\sharp} _{\alpha \beta}:= \fP ^{\sharp} _\alpha \cap \fP ^{\sharp} _\beta$,
$\fP ^{\sharp} _{\alpha \beta \gamma}:= \fP ^{\sharp} _\alpha \cap \fP ^{\sharp} _\beta \cap \fP ^{\sharp} _\gamma$,
$X ^{\sharp} _\alpha := X ^{\sharp} \cap P ^\sharp _\alpha$,
$X ^{\sharp}_{\alpha \beta } := X ^{\sharp} _\alpha \cap X ^{\sharp} _\beta$ et
$X ^{\sharp}_{\alpha \beta \gamma } := X ^{\sharp} _\alpha \cap X ^{\sharp} _\beta \cap X ^{\sharp} _\gamma $.
We denote by $Y _\alpha := (X ^{\sharp} _\alpha) _\eta$,
$Y _{\alpha \beta} := Y _\alpha \cap Y _\beta$,
$Y _{\alpha \beta \gamma} := Y _\alpha \cap Y _\beta \cap Y _\gamma $,
$j _\alpha$ : $ Y _\alpha
\hookrightarrow X ^{\sharp} _\alpha$,
$j _{\alpha \beta} $ :
$Y _{\alpha \beta}  \hookrightarrow  X ^{\sharp} _{\alpha \beta}$
and
$j _{\alpha \beta \gamma} $ :
$Y _{\alpha \beta \gamma }  \hookrightarrow  X ^{\sharp} _{\alpha \beta \gamma} $
the canonical open immersions.
We suppose that for every $\alpha\in \Lambda$, $X ^{\sharp} _\alpha$ is affine, 
(for instance when the covering $(\fP ^{\sharp} _{\alpha}) _{\alpha \in \Lambda}$ is affine).
Since $P$ is separated, for any $\alpha,\beta ,\gamma \in \Lambda$,
$X ^{\sharp}_{\alpha \beta }$ and $X ^{\sharp}_{\alpha \beta \gamma }$ are also affine.

For any 3uple $(\alpha, \, \beta,\, \gamma)\in \Lambda ^3$, fix 
$\fX ^{\sharp} _\alpha$ (resp. $\fX ^{\sharp} _{\alpha \beta}$, $\fX ^{\sharp} _{\alpha \beta \gamma}$)
some semistable log formal scheme over $\fS ^{\sharp}$ (see \ref{rem-semistable lifting}) which is a lifting
of $X ^{\sharp} _\alpha$
(resp. $X ^{\sharp} _{\alpha \beta}$, $X ^{\sharp} _{\alpha \beta \gamma}$),
choose $p _1 ^{\alpha \beta}$ :
$\fX ^{\sharp}  _{\alpha \beta} \rightarrow \fX ^{\sharp} _{\alpha}$
(resp. $p _2 ^{\alpha \beta}$ :
$\fX ^{\sharp}  _{\alpha \beta} \rightarrow \fX ^{\sharp} _{\beta}$)
some lifting  of 
$X ^{\sharp}  _{\alpha \beta} \rightarrow X ^{\sharp} _{\alpha}$
(resp. $X ^{\sharp}  _{\alpha \beta} \rightarrow X ^{\sharp} _{\beta}$).

Similarly, for any $(\alpha,\,\beta,\,\gamma )\in \Lambda ^3$, fix some lifting (use \cite[3.11]{Kato-logFontaine-Illusie} to check the existence
of such liftings) 
$p _{12} ^{\alpha \beta \gamma}$ : $\fX ^{\sharp}  _{\alpha \beta \gamma} \rightarrow \fX ^{\sharp}  _{\alpha \beta} $,
$p _{23} ^{\alpha \beta \gamma}$ : $\fX ^{\sharp}  _{\alpha \beta \gamma} \rightarrow \fX ^{\sharp}  _{\beta \gamma} $,
$p _{13} ^{\alpha \beta \gamma}$ : $\fX ^{\sharp}  _{\alpha \beta \gamma} \rightarrow \fX ^{\sharp}  _{\alpha \gamma} $,
$p _1 ^{\alpha \beta \gamma}$ : $\fX ^{\sharp}  _{\alpha \beta \gamma} \rightarrow \fX ^{\sharp}  _{\alpha} $,
$p _2 ^{\alpha \beta \gamma}$ : $\fX ^{\sharp}  _{\alpha \beta \gamma} \rightarrow \fX ^{\sharp}  _{\beta} $,
$p _3 ^{\alpha \beta \gamma}$ : $\fX ^{\sharp}  _{\alpha \beta \gamma} \rightarrow \fX ^{\sharp}  _{\gamma} $,
$u _{\alpha}$ : $\fX ^{\sharp} _{\alpha } \hookrightarrow \fP ^{\sharp} _{\alpha }$,
$u _{\alpha \beta}$ : $\fX ^{\sharp} _{\alpha \beta} \hookrightarrow \fP ^{\sharp} _{\alpha \beta}$
and
$u _{\alpha \beta \gamma}$ : $\fX ^{\sharp} _{\alpha \beta \gamma } \hookrightarrow \fP ^{\sharp} _{\alpha \beta \gamma}$.

Via the isomorphisms of the form $\tau $ (\ref{2.1.5Be2}), 
we have the following commutative diagrams of functors: 
\begin{equation}
  \label{diag-transp12p1}
  \xymatrix  @R=0,3cm@C=1,5cm {
  {p _1 ^{\alpha \beta \gamma!}}
  \ar[r] ^-{\tau } _-{\sim}
  \ar@{=}[d]
  &
  {p _{12} ^{\alpha \beta \gamma!} \circ p _1 ^{\alpha \beta !}}
  \ar[d]^\tau _-{\sim}
  \\
  {p _1 ^{\alpha \beta \gamma!}}
  \ar[r] ^-{\tau } _-{\sim}
  &
  {p _{13} ^{\alpha \beta \gamma!} \circ p _1 ^{\alpha \gamma !},}
  }
\xymatrix  @R=0,3cm {
  {p _2 ^{\alpha \beta \gamma!}}
  \ar[r] ^-{\tau } _-{\sim}
  \ar@{=}[d]
  &
  {p _{12} ^{\alpha \beta \gamma!} \circ p _2 ^{\alpha \beta !}}
  \ar[d] ^-{\tau } _-{\sim}
  \\
  {p _2 ^{\alpha \beta \gamma!}}
  \ar[r] ^-{\tau } _-{\sim}
  &
  {p _{23} ^{\alpha \beta \gamma!} \circ p _1 ^{\beta \gamma !},}
  }
  \xymatrix  @R=0,3cm {
  {p _3 ^{\alpha \beta \gamma!}}
  \ar[r] ^-{\tau } _-{\sim}
  \ar@{=}[d]
  &
  {p _{13} ^{\alpha \beta \gamma!} \circ p _2 ^{\alpha \gamma !}}
  \ar[d] ^-{\tau } _-{\sim}
  \\
  {p _3 ^{\alpha \beta \gamma!}}
  \ar[r] ^-{\tau } _-{\sim}
  &
  {p _{23} ^{\alpha \beta \gamma!} \circ p _2 ^{\beta \gamma !}.}
  }
\end{equation}
\end{empt}

\begin{dfn}\label{defindonnederecol}
For any $\alpha \in \Lambda$, let $\cE _\alpha$ be a coherent
$\smash{\widetilde{\cD}} ^{\dag} _{\fX ^{\sharp} _{\alpha} /\fS ^\sharp,\bbQ} $-module.
A \textit{glueing data } on $(\cE _{\alpha})_{\alpha \in \Lambda}$
is the data for any $\alpha,\,\beta \in \Lambda$ of a
$\smash{\widetilde{\cD}} ^{\dag} _{\fX ^{\sharp} _{\alpha \beta} /\fS ^\sharp,\bbQ} $-linear
isomorphism
$ \theta _{  \alpha \beta} \ : \  p _2  ^{\alpha \beta !} (\cE _{\beta}) \riso p  _1 ^{\alpha \beta !} (\cE _{\alpha}),$
satisfying the cocycle condition:
$\theta _{13} ^{\alpha \beta \gamma }=
\theta _{12} ^{\alpha \beta \gamma }
\circ
\theta _{23} ^{\alpha \beta \gamma }$,
where $\theta _{12} ^{\alpha \beta \gamma }$, $\theta _{23} ^{\alpha \beta \gamma }$
and $\theta _{13} ^{\alpha \beta \gamma }$ are defined so as to make the following diagram commutative:
\small
\begin{equation}
  \label{diag1-defindonnederecol}
\xymatrix  @R=0,3cm {
{  p _{12} ^{\alpha \beta \gamma !} p  _2 ^{\alpha \beta !}  (\cE _\beta )}
\ar[r] ^-{\tau} _-{\sim}
\ar[d] ^-{p _{12} ^{\alpha \beta \gamma !} (\theta _{\alpha \beta})} _-{\sim}
&
{p _2 ^{\alpha \beta \gamma!}  (\cE _\beta )}
\ar@{.>}[d] ^-{\theta _{12} ^{\alpha \beta \gamma }}
\\
{ p _{12} ^{\alpha \beta \gamma !}  p  _1 ^{\alpha \beta !}  (\cE _\alpha)}
\ar[r]^{\tau} _-{\sim}
&
{p _1 ^{\alpha \beta \gamma!}(\cE _\alpha),}
}
\xymatrix  @R=0,3cm {
{  p _{23} ^{\alpha \beta \gamma !} p  _2 ^{\beta \gamma!}  (\cE _\gamma )}
\ar[r] ^-{\tau} _-{\sim}
\ar[d] ^-{p _{23} ^{\alpha \beta \gamma !} (\theta _{ \beta \gamma})} _-{\sim}
&
{p _3 ^{\alpha \beta \gamma!}  (\cE _\gamma )}
\ar@{.>}[d] ^-{\theta _{23} ^{\alpha \beta \gamma }}
\\
{ p _{23} ^{\alpha \beta \gamma !}  p  _1 ^{ \beta \gamma !}  (\cE _\beta)}
\ar[r]^{\tau} _-{\sim}
&
{p _2 ^{\alpha \beta \gamma!}(\cE _\beta),}
}
\xymatrix  @R=0,3cm {
{  p _{13} ^{\alpha \beta \gamma !} p  _2 ^{\alpha \gamma !}  (\cE _\gamma )}
\ar[r] ^-{\tau} _-{\sim}
\ar[d] ^-{p _{13} ^{\alpha \beta \gamma !} (\theta _{\alpha \gamma})} _-{\sim}
&
{p _3 ^{\alpha \beta \gamma!}  (\cE _\gamma )}
\ar@{.>}[d]^{\theta _{13} ^{\alpha \beta \gamma }}
\\
{ p _{13} ^{\alpha \beta \gamma !}  p  _1 ^{\alpha \gamma !}  (\cE _\alpha)}
\ar[r]^{\tau} _-{\sim}
&
{p _1 ^{\alpha \beta \gamma!}(\cE _\alpha).}
}
\end{equation}
\normalsize

We define the category $\mathrm{Coh} ((\fX ^{\sharp}  _\alpha )_{\alpha \in \Lambda}/\cE ^{\dag}  _K)$ as follows:

- an object is a family $(\cE _\alpha) _{\alpha \in \Lambda}$    of coherent  $\smash{\widetilde{\cD}} ^{\dag} _{\fX ^{\sharp} _{\alpha} /\fS ^\sharp,\bbQ} $-modules
endowed with a glueing data $ (\theta _{\alpha\beta}) _{\alpha ,\beta \in \Lambda}$,

- a morphism
$((\cE _{\alpha})_{\alpha \in \Lambda},\, (\theta _{\alpha\beta}) _{\alpha ,\beta \in \Lambda})
\rightarrow
((\cE ' _{\alpha})_{\alpha \in \Lambda},\, (\theta '_{\alpha\beta}) _{\alpha ,\beta \in \Lambda})$
is a family of $\smash{\widetilde{\cD}} ^{\dag} _{\fX ^{\sharp} _{\alpha} /\fS ^\sharp,\bbQ} $-linear morphisms $f _\alpha$ : $\cE _\alpha \rightarrow \cE '_\alpha$
commuting with glueing datas, i.e., such that the following diagrams are commutative : 
\begin{equation}
  \label{diag2-defindonnederecol}
\xymatrix  @R=0,3cm {
{ p _2  ^{\alpha \beta !} (\cE _{\beta}) }
\ar[d] ^-{p _2  ^{\alpha \beta !} (f _{\beta}) }
\ar[r] ^-{\theta _{\alpha\beta}} _-{\sim}
&
{  p  _1 ^{\alpha \beta !} (\cE _{\alpha}) }
\ar[d] ^-{p  _1 ^{\alpha \beta !} (f _{\alpha})}
\\
{p _2  ^{\alpha \beta !} (\cE '_{\beta})  }
\ar[r]^{\theta '_{\alpha\beta}} _-{\sim}
&
{ p  _1 ^{\alpha \beta !} (\cE '_{\alpha})  .}
}
\end{equation}
\end{dfn}

\begin{dfn} \label{dfnCohXPP}
We denote by  $\mathrm{Coh} (X, \fP ^{\sharp}/\cE ^{\dag}  _K)$ the category  of  coherent $\smash{\widetilde{\cD}} ^{\dag} _{\fP ^{\sharp} /\fS ^\sharp,\bbQ} $-modules with support in $X$. 
\end{dfn}

\begin{lem}
[Construction of $u ^! _0$]
\label{const-u0!}
With Notation \ref{dfnCohXPP}, there exists a canonical functor 
$$u _0 ^! \colon 
\mathrm{Coh} (X, \fP ^{\sharp}/\cE ^{\dag}  _K) \to  \mathrm{Coh} ((\fX ^{\sharp}  _\alpha )_{\alpha \in \Lambda}/\cE ^{\dag}  _K)$$
extending the usual functor $u _0 ^!$ when $X ^{\sharp}$ has a log smooth $\fS ^{\sharp}$-formal scheme lifting.
\end{lem}

\begin{proof}
Using \ref{Berthelot-Kashiwara-full} insteaf of the usual Berthelot-Kashiwara theorem, the proof of \cite[9.3.7.5]{Car25} can be adapted.
\end{proof}

\begin{lem}
\label{const-u0+}
There exists a canonical functor  
$$u _{0+} \colon \mathrm{Coh} ((\fX ^{\sharp} _\alpha) _{\alpha \in \Lambda}/\cE ^{\dag}  _K)
\rightarrow
\mathrm{Coh} (X,\, \fP ^{\sharp} /\cE ^{\dag}  _K)$$
extending the usual functor $u _{0+}$ when $X ^{\sharp}$ has a log smooth formal $\fS ^{\sharp}$-scheme lifting.
\end{lem}

\begin{proof}
Using \ref{Berthelot-Kashiwara-full} insteaf of the usual Berthelot-Kashiwara theorem, the proof of \cite[9.3.7.6]{Car25} can be adapted.
\end{proof}

\begin{thm} \label{prop1}
The functors $u ^!  _0$ and $ u _{0+}$ constructed in respectively \ref{const-u0!} and \ref{const-u0+}  are quasi-inverse equivalences of categories between 
$\mathrm{Coh} ((\fX ^{\sharp} _\alpha) _{\alpha \in \Lambda}/\cE ^{\dag}  _K)$ and $\mathrm{Coh} (X,\, \fP ^{\sharp} /\cE ^{\dag}  _K)$.
\end{thm}

\begin{proof}
Using \ref{Berthelot-Kashiwara-full} instead of the usual Berthelot-Kashiwara theorem, the proof of \cite[9.3.7.7]{Car25} can be adapted.

\end{proof}

\section{Arithmetic $\cD$-modules associated with overconvergent isocrystals}

\subsection{The case of a strictly semistable formal log scheme with divisorial boundary}
Let $\fX ^{\sharp}$ be a strictly semistable log formal scheme over $\fS ^\sharp$.

\subsubsection{Convergent isocrystals on the generic fiber}
Let $\fY$ be an open formal subscheme of $\fX _\eta$. Since $\fY/\eta$ is smooth and  since $\eta = \Spf \cV [[t]] \{ \frac{1}{t}\}$
is the formal spectrum of a complete discrete valuation ring of mixed characteristic 
$(0,p)$, then we can use \cite[11.1]{Car25}, or \cite[4]{Be1} and \cite[3]{Be0}.
In particular, from \cite[4.1.4]{Be1} and \cite[3.1.2 and 3.1.4]{Be0}, we get the following theorem.
\begin{thm}
\label{thm-eqcat-cvisoc}
The following properties hold:
\begin{enumerate}[(a)]
\item The functor $\sp _*$ induces an equivalence of categories between the category of 
convergent isocrystals on $Y/\eta$, and the category of 
coherent 
$\cD ^\dag _{\fY/\eta, \bbQ} $-modules
which are $\cO _{\fY,\bbQ} $-locally projective of finite type.

\item Let $\cE$ be a coherent 
$\cD ^\dag _{\fY/\eta, \bbQ} $-module
which is $\cO _{\fY,\bbQ} $-locally projective of finite type. We have the following properties.
\begin{enumerate}[(a)]
\item For any $m \in \bbN$, there exists 
a (coherent) $\widehat{\cD} ^{(m)} _{\fY/\eta} $-module
$\overset{\circ}{\E}$, coherent over $\cO _\fY$ together with 
an isomorphism of $\widehat{\cD} ^{(m)} _{\fY/\eta, \bbQ} $-modules
$\overset{\circ}{\E} _\bbQ \riso \cE$.
\item The module $\cE$ is $\cD  _{\fY/\eta, \bbQ} $-coherent 
 and for any $m \in \bbN$
the canonical homomorphisms 
\begin{gather}
\cE 
\to 
\widehat{\cD} ^{(m)} _{\fY/\eta, \bbQ}  \otimes _{\cD  _{\fY/\eta, \bbQ} }
\E
,
\
\
\E
\to 
\cD ^\dag _{\fY/\eta, \bbQ} 
 \otimes _{\widehat{\cD} ^{(m)} _{\fY/\eta, \bbQ} }
\E
\end{gather}
are isomorphisms.

\end{enumerate}

\end{enumerate}

\end{thm}

\subsubsection{Equivalence between overconvergent isocrystal and overcoherent isocryslals}
Let $Z$ be a divisor of $X$ containing $X _s$ and $\fY ^\sharp$ the open subset of $\fX ^{\sharp}$ complementary to the support of $Z$.

\begin{ntn}
\label{ntn-tube}
Let $\fX ^{\sharp}$ be a strictly semistable log formal scheme over $\fS ^\sharp$.
Since we will exclusively use adic (rather than rigid analytic or Berkovich spaces), 
we will henceforth (unless otherwise mentioned) use the expression of {\it rigid space} to mean an adic space locally
of finite type over the bounded open unit disc
$\mathbb{D} ^\mathrm{b} _K := \mathrm{Spa} ( \cV [[t]] \otimes _\cV K, \cV [[t]] )$.
We denote by $\fX _K$ the adic space associated to $\fX$
(following theorem \cite[A.5.4]{FujiwaraKatoBookI} and its notation that we will not use, 
$\fX _K := \mathbf{ZR} ( \fX ^{\mathrm{rig}}$)).
We denote by $\sp \colon \fX  _K \to \fX $ 
the specialization morphism.

Let $j\colon Y\hookrightarrow X$ be an open immersion of log schemes. 
Recall the tube of $Y$ in $\fX$ is 
$]Y [ _{\fX} := \overline{\sp ^{-1}(Y)}$, the closure of the open subset
$\sp ^{-1}(Y)$.
If $\fY$ is an open formal subscheme of $\fX$, 
then following \cite[II.3.1.3.(2)]{FujiwaraKatoBookI}
we have
$\fY _K= \sp ^{-1}(Y)$, where $Y$ is the reduction modulo $\pi$ of $\fY$.
\end{ntn}

\begin{empt}
With the notation \ref{ntn-tube}, suppose $\fX = \Spf A$ is affine. 
Using \cite[II.6.5.1]{FujiwaraKatoBookI}, since $(p)$ is a principal ideal of definition of $A$, 
we check that $\fX _K$ is Stein, i.e. for any coherent $\cO _{\fX _K}$-module $\cF$, 
we have $H ^1 (\fX _K, \cF)=0$
(beware that the notation of \cite[II.6.5.1]{FujiwaraKatoBookI} are somehow misleading since their 
$\cO _{\fX ^{\rig}}$ means exactly $\cO _{\mathbf{ZR} ( \fX ^{\mathrm{rig}})}$, i.e. with our notation $\cO _{\fX _K}$).
Since $A$ is regular, it is integrally closed in $A _K:= A \otimes _\cV K$. 
Hence, using \cite[A.4.3 and A.4.7]{FujiwaraKatoBookI}, we get 
 $\fX _K = \mathrm{Spa} ( A _K, A)$.
Following Theorem \cite[II.6.5.7]{FujiwaraKatoBookI}, we have Theorem A and B over $\fX _K$ for coherent
$\cO _{\fX _K}$-modules.

\end{empt}

\begin{empt}
\label{ntn-jm}
Let $\fX ^{\sharp}$ be a strictly semistable log formal scheme over $\fS ^\sharp$.
Let $Z$ be a divisor of $X$ and $Y ^\sharp$ the open log-subscheme of $X ^\sharp$ complementary to the support of $Z$.

Let $\fU = \Spf A$ be an open formal subscheme of $\fX$, let $f \in A$ such that $Y \cap U = \Spec \overline{A} _{\overline{f}}$.
For any integer $m \in \bbN$, we have the affinoid subset
$(Y \cap U) _{m} := \{x \in \fU _K;\;  v_x ( \pi ^{-1}f ^{p ^{m+1}})   \geq 1\}$,
where $v _x$ is a valuation (in the sense of Hubert, i.e. valuations are written multiplicatively) associated with $x$.
We check by an easy computation that this only depends on $Y$, $\fU$ and not on the choice of $f$ such that $\overline{f}$ is an equation of $D \cap U$ in $U$.

Let $(\fU _i= \Spf A _i) _i$ be an affine covering of $\fX$ such that there exists $f _i \in A _i$ satisfying
$Y \cap U _i = \Spec \overline{A_i} _{\overline{f_i}}$. 
Then glueing $(Y \cap U _i) _{m}$, we get the open rigid (in the sense of Hubert) subspace $Y _m$ of $\fX _K =] X [ _\fX$. 
Following \cite[2.27]{Lazda-Pal-Book},
$(Y _m) _{m\in\bbN}$ forms a cofinal system of neighbourhoods of 
$]Y[ _{\fX} $ in $ \fX _K$.
We denote by $j _Y \colon ]Y[  _{\fX} \hookrightarrow \fX _K$  the canonical open immersion 
induced by $j$.
We denote by $j _m \colon Y _m \hookrightarrow \fX _K$ be the canonical immersion. 
We set  $j ^\dag \G = j _{Y*} j ^{-1} _Y \G$ for any $\cO _{\fX _K}$-module $\G$.  
From \cite[2.19]{Lazda-Pal-Book}, we get 
$j ^\dag \G \cong \underrightarrow{\lim} _m j _{m*} j ^{-1} _m \G$.
\end{empt}

\begin{ntn}
Since we will not use the notation \cite[8.7.3.8]{Car25}, in this paper we will denote by $\cB ^{(m)} _{\fX} (Z) _{\bbQ}$ what is denoted by $\cB ^{\prime (m)} _{\fX} (Z) _{\bbQ}$ in \cite[8.7.3.5.1]{Car25} (we remove the prime). 
\end{ntn}

\begin{prop}
\label{432-442Be1}
We keep notations \ref{ntn-jm}.
\begin{enumerate}[(a)]
\item There exist canonical isomorphisms of $\cO _{\fX}$-algebras
\begin{gather}
\cB ^{(m)} _{\fX} (Z) _{\bbQ}
\riso
\sp _* j _{m *} j ^{*} _m \cO _{\fX _K}
\\
\cO _{\fX} (\hdag Z) _\bbQ \riso \sp _* j ^\dag  \cO _{\fX _K}.
\end{gather}

\item For any affine open formal subscheme $\fU \subset \fX$, 
$\Gamma ( \fU, \cB ^{(m)} _{\fX} (Z))$,
and 
$\Gamma ( \fU, \cB ^{(m)} _{\fX} (Z) _\bbQ)$
are noetherian. The extensions 
$\cO _{\fX,\bbQ}
\to 
\cB ^{(m)} _{\fX} (Z) _\bbQ$
and 
$\cB ^{(m)} _{\fX} (Z) _\bbQ
\to 
\cB ^{(m+1)} _{\fX} (Z) _\bbQ $
are flat. 

The sheaves $\cB ^{(m)} _{\fX} (Z) $, 
$\cB ^{(m)} _{\fX} (Z) _\bbQ$,
and
$\cO _{\fX} (\hdag Z) _\bbQ$ 
are coherent, coherent modules over these sheaves 
satisfy the theorems of type $A$ and $B$.

\item The functor $\sp _* \circ j _{m*}$ (resp. $\sp _*$)
induces an exact equivalence between the category of coherent $\cO _{Y _m}$-modules 
(resp. coherent $j ^\dag \cO _{\fX _K}$-modules) 
and that of 
coherent $\cB ^{(m)} _{\fX} (Z) _\bbQ$-modules
(resp. coherent $\cO _{\fX} (\hdag Z) _\bbQ$-modules).
Via this equivalence, locally free modules of finite type (resp. integrable connexions) correspond to 
locally projective modules of finite type (resp. integrable connexions).
\end{enumerate}

\end{prop}

\begin{proof}
The exactness of the functor $\sp _* \circ j _{m*}$ (resp. $\sp _*$)
is a consequence of the second part of the proposition. 
For the rest of the check,
we can copy the proof of \cite[4.3.2 and 4.4.2]{Be1}.
\end{proof}

\begin{empt}
\label{localdesc-ovcvisoc}
Suppose here $\fX = \Spf A$ is affine, there exists 
$f \in A$ such that $Z= \Spec \overline{A} / (\overline{f})$
and overconvergent local coordinates (see definition \ref{dfn-ovcvloccoor}) 
$t _1,\dots, t _d\in A$ of 
$\fX ^\sharp/\fS ^\sharp$ (remark from \ref{ntn-loc-coord-partial} that such a situation appears locally on $\fX ^\sharp$ for the Zariski topology).
Let $E \in \mathrm{Isoc} ^\dag ((Y,X,\fX)/\mathcal{E} ^{\dag} _K)$ be an overconvergent isocrystal on $(Y,X,\fX)/\mathcal{E} ^{\dag} _K$.
From \cite[2.60]{Lazda-Pal-Book}, via a canonical fully faithful functor, we can associate to 
$E$ 
a coherent $j ^\dag \cO _{\fX _K}$-module together with an integrable connection (relative to $S _K$).
Using \cite[2.68]{Lazda-Pal-Book},
for $m$ large enough, there exist a coherent $\cO _{Y _m}$-module $E _m$ together with an integrable connection relative to $S _K$
such that $E \riso j _* ( E _m | ]Y[ _\fX)$.
Following \cite[2.68]{Lazda-Pal-Book} (this corresponds to the analogue of \cite[4.3.9]{LeStum-livreRigCoh}),
the fact that $E \in \mathrm{Isoc} ^\dag ((Y,X,\fX)/\mathcal{E} ^{\dag} _K)$ 
is equivalent to the fact that the connection satisfies the following overconvergent property: 
for any $\eta <1$, there exists $m _\eta \geq 0$ large enough so that for any $e \in \Gamma ( U _{m _\eta} , E _{m _\eta})$, we have  
$||\underline{\partial} ^{[\underline{k} ]} e || \eta ^{| \underline{k}| } \to 0$
when $|\underline{k}| \to \infty$, where $||\cdot ||$ is some choice of Banach norm
on $\Gamma ( U _{m _\eta} , E _{m _\eta})$.

\end{empt}

\begin{prop}
\label{445Be1}
The functor $\sp _*$ induces an exact  fully faithful functor from the category of 
overconvergent isocrystals on $(Y,X,\fX)/\mathcal{E} ^{\dag} _K$, 
to the category of 
$\cD ^\dag _{\fX ^\sharp/\fS ^{\sharp}} (\hdag Z) _\bbQ$-modules.
Its essential image consists in 
$\cD ^\dag _{\fX ^\sharp/\fS ^{\sharp}} (\hdag Z) _\bbQ$-modules
$\cE$ such that there exists $n _0$, and 
a coherent $\cB ^{(n_0)} _{\fX} (Z) _\bbQ$-module $\cE _0$
together with an isomorphism
\begin{equation}
\label{4451Be1}
\underrightarrow{\lim} _{n \geq n _0}
\cB ^{(n)} _{\fX} (Z) _\bbQ \otimes _{\cB ^{(n _0)} _{\fX} (Z) _\bbQ}
\cE _0
\riso 
\E,
\end{equation}
satisfying the following condition:
For any $m \in \bbN$, there exists an integer $n _m \geq \max ( n _{m-1}, m)$
and a structure of 
$(\cB ^{(n _m)} _{\fX} (Z) 
\widehat{\otimes}
\widehat{\cD} ^{(m)} _{\fX ^\sharp/\fS ^\sharp}) _\bbQ$-module
on 
$\cB ^{(n _m)} _{\fX} (Z) _\bbQ \otimes _{\cB ^{(n _0)} _{\fX} (Z) _\bbQ}
\cE _0$
such that 
the homomorphisms
$$
\cB ^{(n _m)} _{\fX} (Z) _\bbQ \otimes _{\cB ^{(n _0)} _{\fX} (Z) _\bbQ}
\cE _0
\to 
\cB ^{(n _{m+1})} _{\fX} (Z) _\bbQ \otimes _{\cB ^{(n _0)} _{\fX} (Z) _\bbQ}
\cE _0
$$
are
$(\cB ^{(n _m)} _{\fX} (Z) 
\widehat{\otimes}
\widehat{\cD} ^{(m)} _{\fX ^\sharp/\fS ^\sharp}) _\bbQ$-linear
and the isomorphism
\ref{4451Be1} is 
$\cD ^\dag _{\fX ^\sharp/\fS ^{\sharp}} (\hdag Z) _\bbQ$-linear.
\end{prop}

\begin{proof}
Using \ref{432-442Be1} and \ref{localdesc-ovcvisoc}, we can follow the proof of Theorems 
\cite[4.4.5 and 4.4.12]{Be1}.
\end{proof}

\begin{empt}
\label{cor-Be44}
In fact, the results of \cite[4.4]{Be1} are still valid. 
As a corollary, we obtain for instance the fact that if $\cE$ is in the essential image, then 
for some (large enough) function $\lambda _0 \colon \bbN \to \bbN$ satisfying $\lambda _0 (m ) \geq m$ (depending on $\cE$), 
there exists $\cE ^{(\bullet) }
\in 
\smash{\underrightarrow{LM}}  _{\bbQ, \mathrm{coh}} ( \smash{\widetilde{\B}} _{\fX} ^{(\bullet)}(Z))
\cap 
\smash{\underrightarrow{LM}}  _{\bbQ, \mathrm{coh}} ( \smash{\widetilde{\cD}} _{\fX ^{\sharp}} ^{(\bullet)}(Z))$, 
together with an isomorphism 
$\underrvec{l}_{\bbQ} ^{*} 
\cE ^{(\bullet) }
\riso 
\cE$
of coherent 
$ \smash{\cD} ^\dag _{\fX ^{\sharp}} (\hdag Z) _{\bbQ} $-modules.

\end{empt}

\begin{lem}
\label{LetterLem2}
We suppose $\fX$ affine.
Let $n \geq m \geq 0$ be two integers. 
Let $\cE$ be a coherent 
$\cB ^{(n)} _{\fX} (Z) 
\widehat{\otimes}
\widehat{\cD} ^{(m)} _{\fX ^\sharp/\fS ^\sharp}$-module. 
Then $\cE$ is $\cB ^{(n)} _{\fX} (Z) $-coherent if and only if
$\Gamma (\fX, \cE) $ is a 
$\Gamma (\fX, \cB ^{(n)} _{\fX} (Z) )$-module of finite type.
\end{lem}

\begin{proof}
See  \cite[7.5.2.3]{Car25}.
\end{proof}

\begin{rem} \label{projec-mdag}
Let $\cE$ be a coherent $\cB ^{(m _0)} _{\fX} (Z) _\bbQ$-module. 
If $\cE | \fY$ is a locally projective $\cO _{\fY,\bbQ}$-module of finite type, then 
for $m \geq m_0$ large enough, 
$\cB ^{(m)} _{\fX} (Z) _\bbQ \otimes _{\cB ^{(m _0)} _{\fX} (Z) _\bbQ} \cE$
is a locally projective $\cB ^{(m)} _{\fX} (Z) _\bbQ$-module of finite type.
Indeed, since $\cO _{\fX} (\hdag Z) _\bbQ\to j _* \cO _{\fY,\bbQ}$ is faithfully flat (we can copy the proof of \cite[4.3.10]{Be1}), 
then $\cO _{\fX} (\hdag Z) _\bbQ \otimes _{\cB ^{(m _0)} _{\fX} (Z) _\bbQ} \cE$
is a projective $\cO _{\fX} (\hdag Z) _\bbQ$-module of finite type. 
We conclude using \cite[3.6.2]{Be1}.

\end{rem}

\begin{ntn}
\label{ntnMICdag2fs}
Let $\mathrm{MIC} ^{\dag\dag} (\fX ^\sharp, Z/\fS ^{\sharp}) $
be the category of coherent 
$\cD ^\dag _{\fX ^\sharp/\fS ^{\sharp}} (\hdag Z) _\bbQ$-modules 
which are $\cO _{\fX} (\hdag Z) _\bbQ$-coherent (see notation \cite[11.2.1.4]{Car25}). Its objects are called overcoherent isocrystals on $(\fX ^\sharp, Z/\fS ^{\sharp})$.
To stick with the notation \ref{ntn-realP}, we can also write $\mathrm{MIC} ^{\dag\dag} (Y, X,\fX ^\sharp/\cE ^{\dag}  _K) $ instead of 
$\mathrm{MIC} ^{\dag\dag} (\fX ^\sharp, Z/\fS ^{\sharp}) $ (this is independant on the choice of $Z $ such that $Y=X \setminus Z$ because of theorem \ref{letterBerthelotCaro2007}).
We have the morphism of ringed spaces 
$\sp \colon 
(\fX _K , j ^\dag \cO _{\fX _K})
\to 
(\fX, \cO _{\fX} (\hdag Z) _\bbQ)$
induced by the specialization morphism.
We get the inverse image functor $\sp ^*$ by setting 
$\sp ^* (\cE ) := j ^\dag \cO _{\fX _K} \otimes _{\sp ^{-1} \widetilde{\cO} _{\fX,Q}}  \sp ^{-1} (\cE)$,
for any $\cE \in \mathrm{MIC} ^{\dag\dag} (\fX ^\sharp, Z/\fS ^{\sharp})  $.
\end{ntn}

\begin{empt}
\label{empt-projff}
Let $\cE \in \mathrm{MIC} ^{\dag\dag} (\fX ^\sharp, Z/\fS ^{\sharp})  $.
It follows from \cite[10.2.2.7]{Car25} and \cite[11.2.1.3]{Car25}, that $\cE | \fY$ is locally projective of finite type over $\cO_{\fY,\bbQ}$.
Since $\cO _{\fX} (\hdag Z) _\bbQ\to j _* \cO _{\fY,\bbQ}$ is faithfully flat (we can copy the proof of \cite[4.3.10]{Be1}), 
then  $\cE$ is a locally projective $\cO _{\fX} (\hdag Z) _\bbQ$-module of finite type. 
This implies that $\cE =0$ if and only if there exists an open dense subset $\fU ^{\sharp}$ of $\fX ^{\sharp}$
such that $\cE | \fU ^{\sharp}=0$.
\end{empt}

The second part of the following theorem extends \cite[11.2.1.12]{Car25}. 
\begin{thm}
\label{letterBerthelotCaro2007}
We have the following properties.
\begin{enumerate}[(a)]
\item The functors $\sp _*$ and $\sp ^*$ induces quasi-inverse equivalences of categories between 
$\mathrm{MIC} ^{\dag} (Y, X,\fX ^\sharp/\cE ^{\dag}  _K) $ (see notation \ref{ntn-realP}) and 
$\mathrm{MIC} ^{\dag\dag} (\fX ^\sharp, Z/\fS ^{\sharp}) $.
\item Let $\cE$ be a coherent $\cD ^\dag _{\fX ^\sharp/\fS ^{\sharp}} (\hdag Z) _\bbQ$-module.
Then $\cE \in \mathrm{MIC} ^{\dag\dag} (\fX ^\sharp, Z/\fS ^{\sharp})  $ if and only if 
$\cE | \fY$ is $\cO _{\fY,\bbQ}$-coherent. 
\end{enumerate}
\end{thm}

\begin{proof}
I) The proof of the first part is similar to \cite[11.2.1.3]{Car25}.

II. We proceed in the same way as \cite[11.2.1.12]{Car25}. 1) Let $\cE$ be a coherent $\cD ^\dag _{\fX ^\sharp/\fS ^{\sharp}} (\hdag Z) _\bbQ$-module such that 
$\cE | \fY$ is $\cO _{\fY,\bbQ}$-coherent. 
It remains to check that $\cE$ is in the essential image of $\sp _*$. The proof of Berthelot of its analogue in the case of a smooth formal scheme over $\cV$ (see \cite[11.2.1.12]{Car25}) can be adapted without too much effort: 
For $n \geq m$, we set  $\widehat{\cD} ^{(m,n)} _{\fX ^\sharp/\fS ^\sharp} (Z):=  \cB ^{(n )} _{\fX} (Z)  \widehat{\otimes} \widehat{\cD} ^{(m)} _{\fX ^\sharp/\fS ^\sharp}$, 
and  $\widehat{\cD} ^{(m)} _{\fX ^\sharp/\fS ^\sharp} (Z):=  \widehat{\cD} ^{(m,m)} _{\fX ^\sharp/\fS ^\sharp} (Z)$.  For $m _0$ large enough, 
there exists a coherent $\widehat{\cD} ^{(m _0)} _{\fX ^\sharp/\fS ^\sharp} (Z) _\bbQ$-module $\cE ^{(m_0)}$ together with a
$\cD ^\dag _{\fX ^\sharp/\fS ^{\sharp}} (\hdag Z) _\bbQ$-linear isomorphism $\cD ^\dag _{\fX ^\sharp/\fS ^{\sharp}} (\hdag Z) _\bbQ \otimes _{\widehat{\cD} ^{(m_0)} _{\fX ^\sharp/\fS ^\sharp} (Z) } \cE ^{(m_0)}
\riso \cE$ (use \cite[3.6.2]{Be1} and the isomorphism $\underrightarrow{\lim} _m \widehat{\cD} ^{(m)} _{\fX ^\sharp/\fS ^\sharp} (Z)  \riso \cD ^\dag _{\fX ^\sharp/\fS ^{\sharp}} (\hdag Z) _\bbQ$).
There exists  a $p$-torsion-free coherent $\widehat{\cD} ^{(m _0)} _{\fX ^\sharp/\fS ^\sharp} (Z)$-module $\G ^{(m _0)}$ together with  a $\widehat{\cD} ^{(m _0)} _{\fX ^\sharp/\fS ^\sharp} (Z) _\bbQ$-linear isomorphism 
$\G ^{(m _0)} _\bbQ\riso \cE ^{(m_0)}$ (see \cite[3.4.5]{Be1}). For any $n\geq m \geq m _0$, we put
$\G ^{(m,n)} := \widehat{\cD} ^{(m,n)} _{\fX ^\sharp/\fS ^\sharp} (Z)  \otimes _{\widehat{\cD} ^{(m_0)} _{\fX ^\sharp/\fS ^\sharp} (Z)}\G ^{(m_0)}/ \text{$p$-torsion}$, $\G ^{(m)} := \G ^{(m,m)}$.
Following \cite[3.4.4]{Be1},  $\G ^{(m,n)}$ is $\widehat{\cD} ^{(m,n)} _{\fX ^\sharp/\fS ^\sharp} (Z) $-coherent. 
From \ref{thm-eqcat-cvisoc},  we get  $\G ^{(m,n)} _\bbQ |\fY  \riso \cE |\fY$. From \cite[11.1.1.11]{Car25}
we get that  $\G ^{(m,n)} |\fY $ is $\cO _{\fY}$-coherent.

II.2) We will now prove that for $n$ large enough $\G ^{(m,n)}$ is  $\cB ^{(n)} _{\fX} (Z)$-coherent. 
Since this is local, using \ref{etalecoordsss} we can suppose
$\fX = \Spf A$ is affine, and there exists an etale $\fS ^\sharp$-morphism of the form 
$h\colon \fX ^{\sharp} \to (\widehat{\mathbb{G}} _{\mathrm{m} }) _\fS ^{s}  \times _\fS \fS ^\sharp _{(r)}  $
and  there exist  $f \in A$ such that $Z= \Spec \overline{A} / (\overline{f})$.
With notation \ref{nota-loc-coord}, we get   log coordinates $u _1,\dots, u _d$ of $\fU ^{\sharp} /\fS ^{\sharp} $
such that $t = u _0 \cdots u _r$ and $u _{r+1},\dots u _{d} \in \cO _{\fU} ^*$, 
the sheaf $\cD ^{(m)} _{\fU^\sharp/\fS^\sharp}$ is $\cO _{\fU ^\sharp}$-free with basis 
$\{ \underline{\partial} _{(r)} ^{<\underline{i}> _{(m)}},\, \underline{i} \in \bbN ^d\}$.
Following \ref{LetterLem2}, we reduce to check that for $n$ large enough,
$\Gamma (\fX, \G ^{(m,n)})$ is a 
$\Gamma (\fX, \cB ^{(n)} _{\fX} (Z) )$-module of finite type.

Put 
$\cD ^{(m)} _{X ^\sharp/S ^\sharp} (Z)
:=
\cD ^{(m)} _{\fX ^\sharp/\fS ^\sharp} (Z)/
\pi \cD ^{(m)} _{\fX ^\sharp/\fS ^\sharp} (Z)$,
$\smash{\overline{\G}} ^{(m)}:= 
\G ^{(m)} / \pi \G ^{(m)}$,
and 
$\smash{\overline{\G}} ^{(m,n)}:= 
\G ^{(m,n)} / \pi \G ^{(m,n)}$.
Let $\overline{x} _1, \dots, \overline{x} _r 
\in 
\Gamma (X, \smash{\overline{\G}} ^{(m _0)})$
which generate 
$\smash{\overline{\G}} ^{(m _0)}$ as 
$\cD ^{(m_0)} _{X ^\sharp/S ^\sharp} (Z)$-module.

Fix $m \geq m _0$. 
From Lemma \cite[11.1.1.11]{Car25},
$\smash{\overline{\G}} ^{(m)}|Y$ is a 
nilpotent 
$\cD ^{(m)} _{Y/\eta}$-module.
Hence, there exists $h \in \bbN$ large enough so that
we get in 
$\Gamma ( Y, \smash{\overline{\G}} ^{(m)})$ the relation
$$\label{letterBerthelotCaro2007-proof1} \forall i=1,\dots, r,
\forall j=1,\dots, d,
\forall l=1,\dots, m,
\
(\partial _{j} ^{[p ^l]}) ^h 
\cdot
\overline{x} _i 
=
0, 
$$
where by abuse of notation we still denote by 
$\overline{x} _i  $ 
the image of $\overline{x} _i $ 
via the canonical map 
$\Gamma (X, \smash{\overline{\G}} ^{(m _0)})\to 
\Gamma ( Y, \smash{\overline{\G}} ^{(m)})$.

Following \cite[2.2.6]{Be1}, for any $\overline{g}\in \overline{A}$,
$\overline{g} ^{p ^{m+1}}$ is in the center of  $\Gamma ( Y , \cD ^{(m_0)} _{Y/\eta})$.
Since  $\Gamma (X, \cD ^{(m)} _{X ^\sharp/S ^\sharp} ) \subset \Gamma ( Y , \cD ^{(m_0)} _{Y/\eta})$, then 
$\overline{g} ^{p ^{m+1}}$ is in the center of  $\Gamma (X, \cD ^{(m)} _{X ^\sharp/S ^\sharp} )$.
Since $t$ divides $\overline{f}$, since $t = u _0 \cdots u _r$ and $u _{r+1},\dots u _{d} \in \cO _{\fX} ^*$
then for $n _m$ such that $p ^{n _m}\geq  h p ^{m+1}$, 
$\forall i=1,\dots, r, \forall j=1,\dots, d, \forall l=1,\dots, m,\ \overline{f} ^{p ^{n _m}} (\partial _{j} ^{[p ^l]}) ^h \in \Gamma (X, \cD ^{(m)} _{X ^\sharp/S ^\sharp} )$.
For such $n _m$, we get in $\Gamma (X, \smash{\overline{\G}} ^{(m)})$ the relation
\begin{equation}
\label{letterBerthelotCaro2007-proof2}\forall i=1,\dots, r,
\forall j=1,\dots, d,
\forall l=1,\dots, m,
\
\overline{f} ^{p ^{n _m}}
(\partial _{j} ^{[p ^l]}) ^h 
\cdot
\overline{x} _i 
=
0. 
\end{equation}
Set $\partial _{(r) j} ^{[p ^l]}= \partial _{\sharp j} ^{[p ^l]}$ if $j\leq r$ and otherwise $\partial _{j(r)} ^{[p ^l]}= \partial _{j} ^{[p ^l]}$.
Since $\overline{f} ^{p ^{m'}}$ is in the center of  $\Gamma (X, \cD ^{(m)} _{X ^\sharp/S ^\sharp} )$ for any $m' >m$ (see \cite[2.2.6]{Be1} ), 
since $t$ divides $\overline{f}$, since $t = u _0 \cdots u _r$ and $\partial _{\sharp j} ^{[p ^l]}= u _j ^{p ^{l}} \partial _{j} ^{[p ^l]}$ for $j\leq r$
then we get from \ref{letterBerthelotCaro2007-proof2}: 
\begin{equation}
\label{letterBerthelotCaro2007-proof3} \forall i=1,\dots, r,
\forall j=1,\dots, d,
\forall l=1,\dots, m,
\
\overline{f} ^{p ^{n _m}}
(\partial _{(r) j} ^{[p ^l]}) ^h 
\cdot
\overline{x} _i 
=
0. 
\end{equation}

Let 
$P =
\prod _{j=1} ^{d}
\prod _{l=1} ^{m}
(\partial _{(r)j} ^{[p ^l]}) ^{h _{jl}} \in \Gamma (X, \cD ^{(m)} _{X ^\sharp/S ^\sharp} )$
where $h _{jl} \in \bbN$.
Using \cite[2.2.6]{Be1}, if there exist $j _0$ and $l _0$ such that $h _{j_0l_0} \geq h$, then  for any $i $, we get from \ref{letterBerthelotCaro2007-proof3}
the relation
$\overline{f} ^{p ^{n _m}} P\cdot
\overline{x} _i 
=
0$ in $\Gamma (X, \smash{\overline{\G}} ^{(m)})$.

Let $x _1, \dots, x _n \in \Gamma (\fX,\G ^{(m _0)})$ be some sections lifting
respectively 
$\overline{x} _1, \dots, \overline{x} _n $.
Let 
$P =
\prod _{j=1} ^{d}
\prod _{l=1} ^{m}
(\partial _{(r)j} ^{[p ^l]}) ^{h _{jl}} \in 
\Gamma (\fX, \cD ^{(m)} _{\fX ^\sharp/\fS ^\sharp} )$
where $h _{jl} \in \bbN$ are such that
there exist
$j _0$ and $l _0$ satisfying
$h _{jl} \geq h$.
We get in 
$\Gamma (\fX,\G ^{(m)})$ the relation
$$\forall i=1,\dots, r,
\
f ^{p ^{n _m}} P\cdot
x _i 
\in 
p 
\Gamma (\fX,\G ^{(m)}),$$
where by abuse of notation we still denote by 
$x _i $ the image of $x _i$ via the canonical map 
$\Gamma (\fX,\G ^{(m _0)}) \to \Gamma (\fX,\G ^{(m)})$.
Let $T _{n _m -1} \in \cB ^{(n _m -1)} _{\fX} (Z) $
be the element such that 
$f ^{p ^{n _m}}T _{n _m -1} =p$.
Since the $\cB ^{(n _m -1)} _{\fX} (Z) $-module
$\G ^{(m,n _m -1)}$ has no $\pi$-torsion then 
it has no $f $-torsion. 
Hence, for such $P$, we get in 
$\Gamma (\fX, \G ^{(m,n _m -1)})$: 
$$\forall i=1,\dots, r,
\
 P\cdot
x _i 
\in 
T _{n _m -1}
\Gamma (\fX, \G ^{(m,n _m -1)}).$$

Let $y _1, \dots, y _s$ be the elements of the form
$\left (\prod _{j=1} ^{d}
\prod _{l=1} ^{m}
(\partial _{(r)j} ^{[p ^l]}) ^{h _{jl}} \right)
\cdot
x _i $
where $h _{jl} \in \{ 0,\dots, h-1\}$ for any $j$ and 
$l$ (beware that these elements and their number depend on $m$).
Using \cite[Proposition 2.3.1]{these_montagnon} for the log part and \cite[2.2.5]{Be1} for the non log part,
we can check $\Gamma (\fX, \cD ^{(m)} _{\fX ^\sharp/\fS ^\sharp} )$ is generated as $\Gamma (\fX, \cO _{\fX})$-module (for its left or right structure) 
by the elements of the form
$\prod _{j=1} ^{d}
\prod _{l=1} ^{m}
(\partial _{(r)j} ^{[p ^l]}) ^{h _{jl}}$,
where
$h _{jl} \in \bbN$.
Since $\overline{x} _1, \dots, \overline{x} _n $
generate 
$\smash{\overline{\G}} ^{(m _0)}$ as 
$\cD ^{(m_0)} _{X ^\sharp/S ^\sharp} (Z)$-module,
then for any $n\geq n _m -1$,
$\Gamma (\fX, \G ^{(m,n)})$
is generated as 
$\Gamma (\fX,\cB ^{(n)} _{\fX} (Z) )$-module by 
$p\Gamma (\fX, \G ^{(m,n)})$ and by
the elements of the forms
$\left (\prod _{j=1} ^{d}
\prod _{l=1} ^{m}
(\partial _{(r)j} ^{[p ^l]}) ^{h _{jl}} \right)
\cdot
x _i $
where $h _{jl}\in \bbN$.
Since $T _{n _m -1}$ divides $p$, then we get 
$$\forall n \geq n _{m} -1,
\
\Gamma (\fX, \G ^{(m,n)})
=
\sum _{i= 1} ^s
\Gamma (\fX,\cB ^{(n)} _{\fX} (Z) )
\cdot y _i
+
T _{n _m -1}\Gamma (\fX, \G ^{(m,n)}).$$
By iteration, this yields
$$\forall n \geq n _{m} -1,
\
\Gamma (\fX, \G ^{(m,n)})
=
\sum _{i= 1} ^s
\Gamma (\fX,\cB ^{(n)} _{\fX} (Z) )
\cdot y _i
+
T _{n _m -1} ^p\Gamma (\fX, \G ^{(m,n)}).$$
For any $n \geq n _m$, we have
$T _{n _m -1} ^p
=
p ^{p-1}T _{n _m} $.
We get, 
$$\forall n \geq n _{m} ,
\
\Gamma (\fX, \G ^{(m,n)})
=
\sum _{i= 1} ^s
\Gamma (\fX,\cB ^{(n)} _{\fX} (Z) )
\cdot y _i
+
p\Gamma (\fX, \G ^{(m,n)}).$$
Since 
$\Gamma (\fX, \G ^{(m,n)})$ is $p$-adically separated and complete, this yields that
$\Gamma (\fX, \G ^{(m,n)})$ is generated
as $\Gamma (\fX,\cB ^{(n)} _{\fX} (Z) )$-module by 
$y _1,\dots, y _s$.

II.3) We can suppose that the sequence $(n _m) _m$ is increasing. 
Set $\cE _m := \G ^{(m,n _m) }_\bbQ$. Following II.2), the $\widehat{\cD} ^{(m,n _m)} _{\fX ^\sharp/\fS ^\sharp} (Z) _\bbQ $-module
$\cE _m$ is $\cB ^{(n _m)} _{\fX} (Z) _\bbQ$-coherent.
Since $\cE _{m _0} |\fY \riso\cE _{m} |\fY$ is a 
coherent $\cD ^\dag _{\fY/\eta,\bbQ}$-module which is $\cO _{\fY,\bbQ}$-coherent, then 
it is locally projective of finite type over $\cO _{\fY,\bbQ}$.
Hence, with the remark \ref{projec-mdag},
increasing $m_0$ is necessary, we can suppose 
$\cE _{m _0}$ is a projective 
$\cB ^{(n _{m _0})} _{\fX} (Z) _\bbQ$-module of finite type.

Following 
\ref{445Be1},
it is sufficient to check that 
the canonical homomorphism
$$\cB ^{(n _m)} _{\fX} (Z) _\bbQ \otimes _{\cB ^{(n _0)} _{\fX} (Z) _\bbQ} \cE _{m _0}  \to \cE _m $$
is an isomorphism for any $m \geq m _0$. 
Fix $m\geq m _0$ and set $B _m := 
\Gamma (\fX, \cB ^{(n _m)} _{\fX} (Z) _\bbQ)$,
$E _m := \Gamma (\fX, \cE _m)$,
$E _{m _0,m} := \Gamma (\fX, \cB ^{(n _m)} _{\fX} (Z) _\bbQ \otimes _{\cB ^{(n _0)} _{\fX} (Z) _\bbQ} \cE _{m _0} )$.
We get the morphism
$E _{m _0,m}
\to 
E _{m}$ 
of $B _m$-modules of finite type.

Following the part II.2) and its notations,
$\Gamma (\fX, \G ^{(m,n _m)})$ is generated as 
$\Gamma (\fX,\cB ^{(n _m)} _{\fX} (Z) )$-module by $y _1,\dots, y _s$.
We remark that $y _1,\dots, y _s \in E _{m _0}$. 
Hence, the morphism $E _{m _0,m}
\to 
E _{m}$
is 
surjective. 
After applying $B _m \to 
\Gamma (\fY,\cO _{\fY,\bbQ})$
to the morphism 
$E _{m _0,m}
\to 
E _{m}$, we get an isomorphism.
Since 
$E _{m _0,m} $ is a projective $B _m$-module of finite type
and since
$B _m 
\to \Gamma (\fY,\cO _{\fY,\bbQ})$
is injective, we get the injectivity of
$E _{m _0,m}
\to 
E _{m}$. We are done.

\end{proof}

\begin{prop}
\label{cohisosurcv}
We set $\cD _{\fX ^\sharp/\fS ^{\sharp}} (\hdag Z) _\bbQ : = \cO _{\fX , \bbQ} (\hdag Z )
\otimes _{\cO _{\fX , \bbQ }} \cD _{\fX ^\sharp/\fS ^{\sharp},\bbQ}  $.
Let 
$\cE \in \mathrm{MIC} ^{\dag\dag} (\fX ^\sharp, Z/\fS ^{\sharp})  $.
Then $\cE $ is $\cD _{\fX ^\sharp/\fS ^{\sharp}} (\hdag Z) _\bbQ$-coherent and
the canonical morphism
\begin{equation}\label{cohisosurcv2}
  \cE \rightarrow 
  \cD ^\dag _{\fX ^\sharp/\fS ^{\sharp}} (\hdag Z) _\bbQ
   \otimes _{\cD _{\fX ^\sharp/\fS ^{\sharp}} (\hdag Z) _\bbQ} \E
\end{equation}
is an isomorphism.
\end{prop}

\begin{proof}
We can copy the proof of \cite[2.2.7]{caro_comparaison}.
\end{proof}

\subsubsection{Commutation with duality}
Let $\fX ^{\sharp}$ be a strictly semistable log formal scheme over $\fS ^\sharp$. Let $Z$ be a divisor of $X$ containing $X _s$ and $\fY ^\sharp$ the open subset of $\fX ^{\sharp}$ complementary to the support of $Z$.
This subsection is a straightforward extension of \cite[11.2.6]{Car25} and \cite[11.2.7]{Car25}. We will need Theorem \ref{theoprincipal} later in \ref{BK-dual}.

\begin{empt}
Let $E \in  \mathrm{Isoc} ^{\dag\dag} (Y, X,\fX ^\sharp/\cE ^{\dag}  _K) $. Let $E _\fX \in \mathrm{MIC} ^{\dag} (Y, X,\fX ^\sharp/\cE ^{\dag}  _K) $ its realization on 
the frame $(Y, X,\fX ^\sharp)$ and $\cE := \mathrm{sp} _* ( E _{\fX})$.
It follows from \cite[11.2.2.2]{Car25}, that we have $\cE \in D _{\mathrm{parf}} (\cO _{\fX } ( \hdag Z ) _{\bbQ})$,
$\cE \in D _{\mathrm{parf}} (\cD  _{\fX ^\sharp/\fS ^{\sharp}} (\hdag Z) _\bbQ)$ and $\cE \in D _{\mathrm{parf}} (\cD ^\dag _{\fX ^\sharp/\fS ^{\sharp}} (\hdag Z) _\bbQ)$.
\end{empt}

\begin{ntn}
For any $\cF \in D ( \cD  _{\fX ^\sharp/\fS ^{\sharp}} (\hdag Z) _\bbQ ) $, we set
$\bbD ^{\mathrm{alg}} _Z (\cF)=
\R \mathcal{H}om  _{\cD  _{\fX ^\sharp/\fS ^{\sharp}} (\hdag Z) _\bbQ }
( \cF , \,\cD  _{\fX ^\sharp/\fS ^{\sharp}} (\hdag Z) _\bbQ \otimes _{\cO _{\fX}}\omega _{\fX ^\sharp/\fS ^{\sharp} } ^{-1})[d _X]$ 
and
$\cF ^{\vee} =
\R \mathcal{H}om  _{\cO _{\fX } ( \hdag Z ) _{\bbQ} }
( \cF , \,\cO_{\fX , \bbQ} ( \hdag Z ))$.
\end{ntn}

\begin{prop}\label{Doevee=De}
There exists a canonical isomorphism 
$$ \theta ^{\mathrm{alg}} \colon\bbD ^{\mathrm{alg}} _Z ( \cO _{\fX } ( \hdag Z ) _{\bbQ} ) \otimes ^\L _{\cO _{\fX } ( \hdag Z ) _{\bbQ}} \cE ^\vee
\riso \bbD ^{\mathrm{alg}} _Z (\cE).$$
\end{prop}
\begin{proof}
This is similar to \cite[2.2.1]{caro_comparaison}.  
\end{proof}

\begin{lem}The following properties hold.
\begin{enumerate}[(a)]
\item  $\cO _{\fX , \bbQ } (\hdag Z) \in D _{\mathrm{parf}}( \cD  _{\fX ^\sharp/\fS ^{\sharp}} (\hdag Z) _\bbQ)$.
\item We have the canonical isomorphism:
\begin{equation}\label{dualB=B}
  \bbD ^{\mathrm{alg}} _Z ( \cO _{\fX } ( \hdag Z ) _{\bbQ}) \riso \cO _{\fX } ( \hdag Z ) _{\bbQ}.
\end{equation}
\end{enumerate}
\end{lem}
\begin{proof}
This is similar to \cite[5.20]{caro_log-iso-hol}.
\end{proof}

\begin{rem}\label{evee=De0}
From \ref{dualB=B} and \ref{Doevee=De}, we get the isomorphism $\cE ^\vee \riso \bbD ^{\mathrm{alg}} _Z (\cE)$.
\end{rem}

\begin{prop} \label{sp*dual} We have the isomorphisms
\begin{gather}
\label{sp*dualiso1} \sp ^* \mathcal{H}om _ {\cO _{\fX , \bbQ } (\hdag Z)} ( \cE , \cO _{\fX , \bbQ } (\hdag Z) )
\riso  \mathcal{H} om  _{ j ^{\dag} \cO _{\fX _K} }(\sp ^*  \cE ,  j ^{\dag} \cO _{\fX _K} ) \\
\label{sp*dualiso2}
 \mathcal{H}om _ {\cO _{\fX , \bbQ } (\hdag Z)} ( \sp _* E , \cO _{\fX , \bbQ } (\hdag Z) )
\riso \sp _* ( \mathcal{H} om  _{ j ^{\dag} \cO _{\fX _K} }( E ,  j ^{\dag} \cO _{\fX _K} )).
\end{gather}
\end{prop}

\begin{proof}
This is similar to \cite[2.2.7]{caro_comparaison}.
\end{proof}

\begin{empt}\label{rhoisom}
Consider the following morphism:
$$\rho _Z\ :\ \bbD ^{\mathrm{alg}} _Z (\cE) \rightarrow
\cD _{\fX } ^{\dag} ( \hdag Z ) _{\bbQ} \otimes _{\cD _{\fX } ( \hdag Z ) _{\bbQ}} \bbD ^{\mathrm{alg}} _Z (\cE)
\riso
\bbD  _Z ( \cD _{\fX } ^{\dag} ( \hdag Z ) _{\bbQ} \otimes _{\cD _{\fX } ( \hdag Z ) _{\bbQ}} \cE)
\rightarrow \bbD  _Z (\cE),$$
where $\bbD  _Z$ is $\cD _{\fX } ^{\dag} ( \hdag Z ) _{\bbQ}$-linear dual functor (when $Z=X _s$, we remove the divisor in the notation).
Since $\cE$ is locally projective of finite type over $\cO _{\fX , \bbQ} (\hdag Z )$,
the morphism
$ \mathcal{H}om _ {\cO _{\fX , \bbQ } (\hdag Z)} ( \cE , \,\cO _{\fX , \bbQ } (\hdag Z) )
\rightarrow
\R \mathcal{H}om _ {\cO _{\fX , \bbQ } (\hdag Z)} ( \cE , \,\cO _{\fX , \bbQ } (\hdag Z) )= \cE ^\vee$ is 
an isomorphism.
From \ref{sp*dualiso2}, we have $\cE ^\vee \in \mathrm{MIC} ^{\dag \dag} (Y, X,\fX ^\sharp/\cE ^{\dag}  _K) $.
Since $\cE ^\vee \riso \bbD ^{\mathrm{alg}} _Z (\cE)$ (see \ref{evee=De0}), 
via \ref{cohisosurcv2} we check that 
$\rho ^\dag _Z$ is an isomorphism.
\end{empt}

\begin{empt}\label{constrhodag}
Let  $\theta$ :
$\bbD   _Z ( \cO _{\fX } ( \hdag Z ) _{\bbQ} ) \otimes _{\cO _{\fX } ( \hdag Z ) _{\bbQ}} \cE ^\vee
\riso \bbD _Z (\cE)$
be the isomorphism making commutative the following diagram:
$$\xymatrix {
{\bbD ^{\mathrm{alg}} _Z ( \cO _{\fX } ( \hdag Z ) _{\bbQ} ) \otimes _{\cO _{\fX } ( \hdag Z ) _{\bbQ}} \cE ^\vee}
\ar[r] ^(0.68){\theta ^{\mathrm{alg}}} _(0.68){\sim}
\ar[d] ^{\rho _Z\otimes id} _{\sim}
&
{\bbD ^{\mathrm{alg}} _Z (\cE)}
\ar[d] ^{\rho  _Z} _{\sim}
\\
{\bbD  _Z ( \cO _{\fX } ( \hdag Z ) _{\bbQ} ) \otimes _{\cO _{\fX } ( \hdag Z ) _{\bbQ}} \cE ^\vee }
\ar@{.>}[r] ^(0.68){\theta} _(0.68){\sim}
&
{\bbD _Z (\cE)}
}$$
\end{empt}

\begin{empt}
  \label{dualisoscvdag} 
  From  \ref{dualB=B} and \ref{rhoisom}, we get the isomorphism
  $\bbD  _Z ( \cO _{\fX } ( \hdag Z ) _{\bbQ} ) \riso \cO _{\fX } ( \hdag Z ) _{\bbQ}$.
Hence, the isomorphism $\theta $ induces the following one
  $\cE ^\vee \riso \bbD _Z (\cE)$. 
\end{empt}

\begin{thm}
\label{theoprincipal}
We have the canonical isomorphism in $\mathrm{MIC} ^{\dag \dag} (Y, X,\fP ^\sharp/\cE ^{\dag}  _K)$
\begin{equation}
\label{comspdual} 
\mathrm{sp} _* ( E ^{\vee}) \riso \bbD _{Z } (\cE ).
\end{equation}
\end{thm}
\begin{proof}
This is straightforward from \ref{sp*dual}, \ref{constrhodag} and \ref{dualisoscvdag}.
\end{proof}

\subsection{The case of a realizable strictly semistable log scheme}

\subsubsection{Construction of $\sp _+$}
\label{ntnsp+}
Let $\fP ^\sharp$ be a strictly semi-stable log formal scheme over $\fS ^\sharp$.
Let $u _0\colon X ^\sharp \to P ^\sharp$ be a weakly closed immersion of strictly semistable log schemes over $S ^\sharp$ (see Definition \ref{dfn-wclimm}). 
To simplify notation, put $Y:= X _\eta$ (remark we only consider in this section overconvergent singularities along the fiber at the closed point). 
In that case, we remark that we get a smooth proper frame $(Y,X , \fP) /\cE ^{\dag} _K$ in the terminology 
of \cite[2.22]{Lazda-Pal-Book}. We will call such a frame a ``semistable'' frame. For simplicity, we will only restrict ourselves to semistable frames.

\begin{ntn}
\label{ntn-realP}
We denote by 
$\mathrm{Isoc} ^{\dag} ((Y, X,\fP)/\cE ^{\dag}  _K)$ the category of overconvergent 
isocrystals on $(Y, X,\fP)/\cE ^{\dag}  _K$ (see \cite[2.54.3]{Lazda-Pal-Book}).
We denote by 
$\mathrm{MIC}  ((Y, X,\fP)/\cE ^{\dag}  _K)$ the category of 
coherent  $j ^\dag \cO _{]X[}$-modules together with an integrable connection,
and  by  $\mathrm{MIC} ^{\dag} ((Y, X,\fP)/\cE ^{\dag}  _K)$ the full subcategory of  $\mathrm{MIC}  ((Y, X,\fP)/\cE ^{\dag}  _K)$ of
coherent  $j ^\dag \cO _{]X[}$-modules together with an overconvergent connection (see \cite[2.61]{Lazda-Pal-Book}).
Recall following \cite[Theorem 2.62]{Lazda-Pal-Book} that the realisation functor
$E \mapsto E _\fP$ induces an equivalence of categories 
\begin{equation}
\label{ntn-realPdfn} \mathrm{real} _\fP
\colon 
\mathrm{Isoc} ^{\dag} ((Y, X,\fP)/\cE ^{\dag}  _K) \cong \mathrm{MIC} ^{\dag} ((Y, X,\fP)/\cE ^{\dag}  _K).
\end{equation}
\end{ntn}

\begin{empt}
\label{glueingisocntn}
Consider the following diagrams 
\begin{equation}
  \label{notation-rig-form-diag0}
  \xymatrix  @R=0,3cm {
{Y ' }
\ar[r] ^-{j '}
\ar[d]^b
&
{X ^{\prime \sharp}}
\ar[r] ^-{i '}
\ar[d]^a
&
{\fP ^{\prime \sharp}}
\ar[d]^u
\\
{Y }
\ar[r]^{j }
&
{X ^{\sharp}}
\ar[r]^{i }
&
{\fP ^{\sharp},}
}
\hfill
\xymatrix  @R=0,3cm {
{Y '' }
\ar[r] ^-{j ''}
\ar[d]^{b'}
&
{X ^{\prime \prime\sharp}}
\ar[r] ^-{i ''}
\ar[d]^{a'}
&
{\fP ^{\prime \prime \sharp}}
\ar[d]^{u'}
\\
{Y '}
\ar[r]^{j '}
&
{X ^{\prime \sharp}}
\ar[r]^{i '}
&
{\fP ^{\prime \sharp},}
}
\end{equation}
where $u$ (resp. $u'$) is a morphism of strictly semistable log formal schemes over $\fS ^\sharp$, 
the reduction modulo $\pi$ of $i $, $i'$ and $i''$ are weakly closed immersion of strictly semistable log schemes over $S ^\sharp$  (see Definition \ref{dfn-wclimm}),
where $Y:= X _\eta $, $Y':= X ^{\prime} _\eta $ and $Y'':= X ^{\prime \prime} _\eta $, and 
$j$, $j'$, $j''$ are the canonical open immersions.
In particular, we get the  morphisms of strictly semistable frames
$f :=(b,\,a,\,u)
\colon 
(Y', X ^{\prime },\fP')
\to 
(Y, X,\fP)$,
$f ':=(b,'\,a',\,u')
\colon 
(Y'', X ^{\prime \prime },\fP'')
\to 
(Y', X',\fP')$.
We define the functor 
$f _K ^*\colon   \mathrm{MIC} ^{\dag} ((Y, X,\fP)/\cE ^{\dag}  _K)  \to   \mathrm{MIC} ^{\dag} ((Y', X',\fP')/\cE ^{\dag}  _K) $, 
{\it the pullback by $f _K$},  by setting  of  for any object  $E _{\fP}\in   \mathrm{MIC} ^{\dag} ((Y, X,\fP)/\cE ^{\dag}  _K) $:
$$f _K ^* (E _{\fP} ) 
:= 
j ^{\prime \dag} \smash{\cO} _{]X'[ _{\fP'}} \otimes  _{ u ^{-1} _K j ^\dag  \smash{\cO} _{]X[ _{\fP}}} u ^{-1} _K E _{\fP} 
\riso 
j ^{\prime \dag} \smash{\cD} _{]X'[ _{\fP'}} \otimes  _{ u ^{-1} _K j ^\dag  \smash{\cD} _{]X[ _{\fP}}} u ^{-1} _K E _{\fP} .$$
The functor $\mathrm{real} _\fP$ of \ref{ntn-realPdfn} commutes with inverse images: 
For any $E \in \mathrm{Isoc} ^{\dag} ((Y, X,\fP)/\cE ^{\dag}  _K) $ we have 
the natural isomorphism
$f ^* _K \circ \mathrm{real} _\fP (E ) \riso \mathrm{real} _{\fP '} \circ  f ^* (E)$
which can also be written
$f ^* _K (E _\fP) \riso f ^* (E)  _{\fP '} $. 

Let $E \in \mathrm{Isoc} ^{\dag} ((Y, X,\fP)/\cE ^{\dag}  _K) $. 
We have the left and the right projection 
$q _0 , q _1
\colon 
 \fP ^{ \sharp} \times _{\fS ^\sharp} \fP ^{ \sharp}
 \to 
  \fP ^{ \sharp}$. 
This yields the morphisms of frames 
$p _0 := (id,id, q _0), 
p _1:=(id,id, q _1)
\colon 
(Y, X,\fP \times _\fS \fP) \to (Y, X,\fP)$.
We denote by $E _{\fP \times _\fS \fP}$ the realization of $E$ on the frame
$(Y, X,\fP \times _\fS \fP)$. 
With the notation of \cite[2.54]{Lazda-Pal-Book},
we get the isomorphisms
$\psi _{p_0} 
\colon 
p _0 ^* E _{\fP} \riso E _{\fP \times _\fS \fP} $
and
$\psi _{p_1} 
\colon 
p _1 ^* E _{\fP} \riso E _{\fP \times _\fS \fP} $.
This yields the isomorphism 
$$\epsilon = \psi _{p_0} ^{-1} \circ \psi _{p_1} 
\colon 
p _1 ^* E _{\fP} 
\riso 
p _0 ^* E _{\fP} ,$$ which corresponds to the overconvergent connection on 
$E _{\fP}$.

Let $v\colon \fP ^{\prime \sharp} \to \fP ^\sharp$ be a second morphism of strictly semistable log formal schemes over $\fS ^\sharp$
such that 
we get the morphism of frames  $g :=(b,\,a,\,v)
\colon 
(Y', X ',\fP')
\to 
(Y, X,\fP)$.
The morphism 
$(u,v)\colon  \fP ^{\prime \sharp} \to  \fP ^{ \sharp} \times _{\fS ^\sharp} \fP ^{ \sharp}$
induces the morphism of frame
$\delta  _{u,v}= (b , a, (u,v)) 
\colon 
(Y', X',\fP')
\to 
(Y, X,\fP \times _\fS \fP) $.
We have the morphisms of frames
$f =  p _0 \circ \delta  _{u,v}$ and $g = p _1  \circ \delta  _{u,v}$.
This yields the glueing isomorphism
$$\epsilon _{u,v}:= \delta  _{u,v} ^* (\epsilon) 
\colon 
g ^{*} E _{\fP} 
\riso 
f ^{ *} E _{\fP} .$$
Using the property of the isomorphism $\psi $ of 
\cite[2.54]{Lazda-Pal-Book}, we have the equality
$\epsilon _{u,v}:= \delta  _{u,v}^* (\epsilon) = \psi _{f} ^{-1} \circ \psi _{g} $.
In particular, $\epsilon _{u,u} =id$.
Let $w\colon \fP ^{\prime \sharp} \to \fP ^\sharp$ be a third morphism of strictly semistable log formal schemes over $\fS ^\sharp$
such that 
we get the morphism of frames  $h :=(b,\,a,\,w)
\colon 
(Y', X ',\fP')
\to 
(Y, X,\fP)$.
We have the transitive formula
$\epsilon_{u,\, w} = \epsilon_{u,\,v}\circ \epsilon_{v,\,w}$.
Let $v'\colon \fP ^{\prime\prime \sharp} \to \fP ^{\prime \sharp}$ be a second morphism of strictly semistable log formal schemes over $\fS ^\sharp$
such that 
we get the morphism of frames $g ':=(b,'\,a',\,v')
\colon 
(Y'', X ^{\prime \prime },\fP'')
\to 
(Y', X',\fP')$.
We check the formulas
$ \epsilon _{u \circ u', v \circ u'} = f ^{\prime *} _K \circ \epsilon _{u,v}$
and
$\epsilon _{u' ,v'} \circ f ^* _K =  \epsilon _{u \circ u', u \circ v'}$.

\end{empt}

\begin{empt}
\label{sp*f*com-empt}
We keep notation \ref{glueingisocntn}. 
We suppose $X ^\sharp = P ^\sharp$,
$X ^{\prime \sharp} = P ^{\prime \sharp}$
and $i,i'$ are the canonical closed immersion. 
In this case, $\fP$ will be denoted by $\fX$ and
if this do not cause too much confusion, we will simply write
 $u _K ^*$ instead of $f _K ^*$.
Let $E _\fX \in \mathrm{MIC} ^{\dag} (Y, X,\fX/\cE ^{\dag}  _K) $. 
Let $\cE \in \mathrm{MIC} ^{\dag\dag} (Y, X,\fX ^\sharp/\cE ^{\dag}  _K) $ (see Notation \ref{ntnMICdag2fs}).
Following Theorem \ref{letterBerthelotCaro2007}, 
the functors $\sp _*$ and $\sp ^*$ induce quasi-inverse equivalences of categories 
between 
$ \mathrm{MIC} ^{\dag} (Y, X,\fX/\cE ^{\dag}  _K) $
and
$\mathrm{MIC} ^{\dag\dag} (Y, X,\fX ^\sharp/\cE ^{\dag}  _K) $.
We have the functor 
$$u ^! [-d _{X'/X}] \colon \mathrm{MIC} ^{\dag\dag} (Y, X,\fX ^\sharp/\cE ^{\dag}  _K)  \to  \mathrm{MIC} ^{\dag\dag} (Y', X',\fX ^{\prime \sharp}/\cE ^{\dag}  _K) $$ 
which is compatible with $u _K ^*$,  i.e. there exist canonical isomorphisms respectively of $ \mathrm{MIC} ^{\dag} (Y', X',\fX'/\cE ^{\dag}  _K) $
and $\mathrm{MIC} ^{\dag\dag} (Y', X',\fX ^{\prime \sharp}/\cE ^{\dag}  _K) $ of the form
\begin{equation}
\label{sp*f*com}
\sp _* u _K ^* (E _\fX)  \riso u ^!  \sp _* (E _\fX)[-d _{X'/X}],
\ 
u _K ^* \sp ^* (\cE) \riso \sp ^* u ^! (\cE) [-d _{X'/X}].
\end{equation}
Moreover, these isomorphisms are transitive with respect to the composition of morphisms. 
\end{empt}

\begin{prop}
\label{sp-eps-tau}
With the notation \ref{sp*f*com-empt}, 
the following diagrams
$$\xymatrix  @R=0,3cm@C=1,4cm {
{ \sp _* u ^{\prime *} _K (E )[d _{X'/X}]}
\ar[r] ^-{\sp _*(\epsilon _{u,\,u'})} _-{\sim}
\ar[d] _-{\sim} ^-{\ref{sp*f*com}}
&
{ \sp _* u ^{*} _K (E )[d _{X'/X}]}
\ar[d] _-{\sim}  ^-{\ref{sp*f*com}}
\\
{ u ^{\prime !} \sp _*  (E )}
\ar[r] ^-{\tau _{u ,u' }} _-{\sim}
&
{ u ^{!}  \sp _*  (E ),}
}
\xymatrix  @R=0,3cm@C=1,4cm {
{ u ^{\prime *} _K \sp ^* (\cE )[d _{X'/X}]}
\ar[r] ^-{\epsilon _{u,\,u'}} _-{\sim}
\ar[d] _-{\sim}   ^-{\ref{sp*f*com}}
&
{ u^{*} _K   \sp ^*(\cE )[d _{X'/X}]}
\ar[d] _-{\sim}   ^-{\ref{sp*f*com}}
\\
{\sp ^*  u ^{\prime !}  (\cE )}
\ar[r] ^-{\sp ^* (\tau _{u,u'})} _-{\sim}
&
{ u ^{!} \sp _*  (\cE ),}
}
$$
where the glueing isomorphisms $\epsilon _{u,u'}$ and $\tau _{u,u'}$ are that of \ref{glueingisocntn} and \ref{2.1.5Be2},
are commutative.
\end{prop}

\begin{proof}
This follows from the fact that both glueing isomorphisms $\epsilon _{u,u'}$ and $\tau _{u,u'}$ are built in the same way using some factorization
via the closed imbedding 
$(u,u') \colon \fX ^{\prime \sharp} \hookrightarrow (\fX ^{ \sharp} )^{(n)}$
where 
$(\fX ^{ \sharp} )^{(n)}$
is the $n$th infinitesimal neighborhood of the diagonal immersion
$\fX ^{ \sharp}\hookrightarrow
\fX ^{ \sharp} \times _{\fS ^\sharp} \fX ^{ \sharp}$
\end{proof}

\begin{dfn}
\label{dfnMICdagalphabeta}
With notation \ref{ntnPPalpha},
we define the category
  $\mathrm{MIC} ^\dag ( (\fX   _\alpha )_{\alpha \in \Lambda}/\cE ^{\dag}  _K)$ as follows:
the objects are families $(E _\alpha) _{\alpha \in \Lambda}$
  of objects $E _\alpha $
  of $\mathrm{MIC} ^{\dag} (Y _\alpha, X _\alpha,\fX _\alpha/\cE ^{\dag}  _K)$
together with a {\it glueing data}, i.e., a collection of isomorphisms in 
$\mathrm{MIC} ^{\dag} (Y _{\alpha\beta}, X _{\alpha\beta},\fX _{\alpha\beta}/\cE ^{\dag}  _K)$
 of the form
$ \eta _{  \alpha \beta} \ : \
p _{2 K}  ^{\alpha \beta *} (E _{\beta})
\riso
p  _{1 K} ^{\alpha \beta *} (E _{\alpha})$
satisfying the cocycle condition:
$\eta _{13} ^{\alpha \beta \gamma }=
\eta _{12} ^{\alpha \beta \gamma }
\circ
\eta _{23} ^{\alpha \beta \gamma }$,
where $\eta _{12} ^{\alpha \beta \gamma }$, $\eta _{23} ^{\alpha \beta \gamma }$
and $\eta _{13} ^{\alpha \beta \gamma }$ are defined so that the following diagrams 
\begin{equation}
  \label{diag1-defindonnederecolK}
\xymatrix  @R=0,3cm @C=0,45cm {
{  p _{12K} ^{\alpha \beta \gamma *} p  _{2K} ^{\alpha \beta *}  (E _\beta )}
\ar[r] ^-{\epsilon} _-{\sim}
\ar[d] ^-{p _{12K} ^{\alpha \beta \gamma *} (\eta _{\alpha \beta})} _-{\sim}
&
{p _{2K} ^{\alpha \beta \gamma*}  (E _\beta )}
\ar@{.>}[d] ^-{\eta _{12} ^{\alpha \beta \gamma }}
\\
{ p _{12K} ^{\alpha \beta \gamma *}  p  _{1K} ^{\alpha \beta *}  (E _\alpha)}
\ar[r]^-{\epsilon} _-{\sim}
&
{p _{1K} ^{\alpha \beta \gamma*}(E _\alpha),}
}
\xymatrix  @R=0,3cm @C=0,45cm {
{  p _{23K} ^{\alpha \beta \gamma *} p  _{2K} ^{\beta \gamma*}  (E _\gamma )}
\ar[r] ^-{\epsilon} _-{\sim}
\ar[d] ^-{p _{23K} ^{\alpha \beta \gamma *} (\eta _{ \beta \gamma})} _-{\sim}
&
{p _{3K} ^{\alpha \beta \gamma*}  (E _\gamma )}
\ar@{.>}[d] ^-{\eta _{23} ^{\alpha \beta \gamma }}
\\
{ p _{23K} ^{\alpha \beta \gamma *}  p  _{1K} ^{ \beta \gamma *}  (E _\beta)}
\ar[r]^-{\epsilon} _-{\sim}
&
{p _{2K} ^{\alpha \beta \gamma*}(E _\beta),}
}
\xymatrix  @R=0,3cm @C=0,45cm {
{  p _{13K} ^{\alpha \beta \gamma *} p  _{2K} ^{\alpha \gamma *}  (E _\gamma )}
\ar[r] ^-{\epsilon} _-{\sim}
\ar[d] ^-{p _{13K} ^{\alpha \beta \gamma *} (\eta _{\alpha \gamma})} _-{\sim}
&
{p _{3K} ^{\alpha \beta \gamma*}  (E _\gamma )}
\ar@{.>}[d]^-{\eta _{13} ^{\alpha \beta \gamma }}
\\
{ p _{13K} ^{\alpha \beta \gamma *}  p  _{1K} ^{\alpha \gamma *}  (E _\alpha)}
\ar[r]^-{\epsilon} _-{\sim}
&
{p _{1K} ^{\alpha \beta \gamma*}(E _\alpha)}
}
\end{equation}
are commutative.

A morphism
$f= (f _\alpha)_{\alpha \in \Lambda}\colon ((E _{\alpha})_{\alpha \in \Lambda},\, (\eta _{\alpha\beta}) _{\alpha ,\beta \in \Lambda})
\rightarrow
((E '_{\alpha})_{\alpha \in \Lambda},\, (\eta '_{\alpha\beta}) _{\alpha ,\beta \in \Lambda})$
of $\mathrm{MIC} ^\dag ( (\fX   _\alpha )_{\alpha \in \Lambda}/\cE ^{\dag}  _K)$ 
is by definition a family of 
morphisms $f _\alpha$ : $ E _\alpha \rightarrow E ' _\alpha$
commuting with the glueing data.
\end{dfn}

\begin{prop}
\label{eqcat-iso-reco}
There exists a canonical equivalence of categories
  $$u ^*  _{0K}   \colon   \mathrm{MIC} ^{\dag} ((Y, X,\fP)/\cE ^{\dag}  _K)      \cong  \mathrm{MIC} ^\dag ( (\fX   _\alpha )_{\alpha \in \Lambda}/\cE ^{\dag}  _K).$$
\end{prop}

\begin{proof}
1) Let $\phi _\alpha : =(id,id, u _\alpha) 
\colon 
(Y _\alpha, X _\alpha, \fX _\alpha)
\to 
(Y _\alpha, X _\alpha, \fP _\alpha)$ be the proper morphism of frames. 
We remark that  $\phi _\alpha$ is the composition of the morphism
induced by the graph of $u _\alpha$
$\gamma _{u _\alpha}
\colon 
(Y _\alpha, X _\alpha, \fX _\alpha)
\to 
(Y _\alpha, X _\alpha, \fX _\alpha \times _\fS \fP _\alpha)$
and of the projection 
$p _2\colon (Y _\alpha, X _\alpha, \fX _\alpha \times _\fS \fP _\alpha)
\to 
(Y _\alpha, X _\alpha, \fP _\alpha)$.
Since $p _2$ is proper and smooth (see Definition \cite[2.22]{Lazda-Pal-Book}),
then following \cite[2.56]{Lazda-Pal-Book} we get the equivalence of categories
$p _{2K} ^* \colon 
\mathrm{MIC} ^{\dag} ((Y _\alpha, X _\alpha, \fP _\alpha)/\cE ^{\dag}  _K)
\cong
\mathrm{MIC} ^{\dag} ((Y _\alpha, X _\alpha, \fX _\alpha \times _\fS \fP _\alpha)/\cE ^{\dag}  _K)$.
Since $\gamma _{u _\alpha}$ has a proper, smooth retraction (the first projection), 
then using \cite[2.56]{Lazda-Pal-Book} we get that
$\gamma _{u _\alpha K} ^*$ is also an equivalence of categories. 
Hence, by composition, so is 
$\phi _{\alpha K} ^*
\colon 
\mathrm{MIC} ^{\dag} ((Y _\alpha, X _\alpha, \fP _\alpha)/\cE ^{\dag}  _K)
\cong
\mathrm{MIC} ^{\dag} ((Y _\alpha, X _\alpha, \fX _\alpha)/\cE ^{\dag}  _K)$.

1') Let $\phi _{\alpha \beta} : =(id,id, u _{\alpha \beta}) \colon (Y _{\alpha \beta}, X _{\alpha \beta}, \fX _{\alpha \beta})
\to (Y _{\alpha \beta}, X _{\alpha \beta}, \fP _{\alpha \beta})$ be the morphism of frames.  For the same reason, we get the equivalence of categories
$\phi _{\alpha \beta K} ^* \colon  \mathrm{MIC} ^{\dag} ((Y _{\alpha \beta}, X _{\alpha \beta}, \fP _{\alpha \beta})/\cE ^{\dag}  _K)
\cong \mathrm{MIC} ^{\dag} ((Y _{\alpha \beta}, X _{\alpha \beta}, \fX _{\alpha \beta})/\cE ^{\dag}  _K)$.

2) Let $E _\fP  \in \mathrm{MIC} ^{\dag} ((Y, X,\fP)/\cE ^{\dag}  _K)$.  
Using the properties of the glueing isomorphisms of \ref{glueingisocntn}, we get canonically the object
  $ (\phi ^* _{\alpha K} ( E _\fP |_{]X _\alpha[ _{\fP _\alpha}})) _{\alpha \in \Lambda}$
  of $ \mathrm{MIC} ^\dag ( (\fX   _\alpha )_{\alpha \in \Lambda}/\cE ^{\dag}  _K)$.
The functoriality is obvious and 
this yields the canonical functor 
$u ^*  _{0K}
  \colon
  \mathrm{MIC} ^{\dag} ((Y, X,\fP)/\cE ^{\dag}  _K) 
    \cong
  \mathrm{MIC} ^\dag ( (\fX   _\alpha )_{\alpha \in \Lambda}/\cE ^{\dag}  _K)$.
Since 
$\phi _{\alpha  K} ^*$ is fully faithful 
and 
$\phi _{\alpha \beta K} ^*$ is faithful, we check easily that the functor
$u ^*  _{0K}$ is fully faithful. 

3) Let us now prove the essential surjectivity of $u ^*  _{0K}$. 
Let  $((E _{\alpha})_{\alpha \in \Lambda},\, (\eta _{\alpha\beta}) _{\alpha ,\beta \in \Lambda})$
be an object of  $\mathrm{MIC} ^\dag ( (\fX   _\alpha )_{\alpha \in \Lambda}/\cE ^{\dag}  _K)$.
Let 
$E '_\alpha \in \mathrm{MIC} ^{\dag} ((Y _\alpha, X _\alpha, \fP _\alpha)/\cE ^{\dag}  _K)$ together with 
  an isomorphism $\iota _\alpha \colon  \phi ^* _{\alpha K} (E '_\alpha) \riso E _\alpha$ of 
  $\mathrm{MIC} ^{\dag} ((Y _\alpha, X _\alpha, \fP _\alpha)/\cE ^{\dag}  _K)$.
There exists a unique isomorphism $\tau ' _{\alpha \beta}$ :
  ${E '_\beta} |_{]X _{\alpha \beta}[ _{\fP _{\alpha \beta}}}
  \riso {E '_\alpha} |_{]X _{\alpha \beta}[ _{\fP _{\alpha \beta}}}$ 
  making commutative the diagram :
  \begin{equation}
    \label{diag2-eqcat-iso-reco}
    \xymatrix  @R=0,3cm {
  {\phi ^* _{\alpha \beta K} ({E '_\beta} |_{]X _{\alpha \beta}[ _{\fP _{\alpha \beta}}})}
  \ar[r] ^-{\epsilon} _-{\sim}
  \ar@{.>}[d] ^-{\phi ^* _{\alpha \beta K}(\tau ' _{\alpha \beta} )}
  &
  { p _{2 K}  ^{\alpha \beta *} \phi ^* _{\beta K} (E '_{\beta}) }
  \ar[rr] ^-{ p _{2 K}  ^{\alpha \beta *}(\iota _\beta)} _-{\sim}
  &
  &
  { p _{2 K}  ^{\alpha \beta *} (E _{\beta}) }
    \ar[d] ^-{\eta  _{\alpha \beta} } _-{\sim}
  \\
   {\phi ^* _{\alpha \beta K}({E '_\alpha} |_{]X _{\alpha \beta}[ _{\fP _{\alpha \beta}}})}
   \ar[r]^-{\epsilon} _-{\sim}
   &
  {p  _{1 K} ^{\alpha \beta *}\phi ^* _{\alpha K}  (E ' _{\alpha})}
  \ar[rr] ^-{p  _{1 K} ^{\alpha \beta *}(\iota _\alpha)} _-{\sim}
  &
  &
  {p  _{1 K} ^{\alpha \beta *} (E _{\alpha}).}
    }
  \end{equation}
Let us check that the isomorphisms $\tau ' _{\alpha \beta}$ can be glued.
Consider the following commutative diagram
      \footnotesize
      \begin{equation}
      \label{diag3-eqcat-iso-reco}\xymatrix  @R=0,3cm @C=0,5cm {
    {\phi ^* _{\alpha \beta \gamma K} ({E '_\beta} |_{]X _{\alpha \beta\gamma}[ _{\fP _{\alpha \beta\gamma}}})}
  \ar[r] ^-{\epsilon} _-{\sim}
  \ar[d]
  ^-{\phi ^* _{\alpha \beta \gamma K}({\tau ' _{\alpha \beta}} |_{]X _{\alpha \beta \gamma}[ _{\fP _{\alpha \beta\gamma}}})}
  _-{\sim}
  &
  {p _{12K} ^{\alpha \beta \gamma *}\phi ^* _{\alpha \beta K} ({E '_\beta} |_{]X _{\alpha \beta}[ _{\fP _{\alpha \beta}}})}
  \ar[r] ^-{\epsilon} _-{\sim}
  \ar[d] ^-{p _{12K} ^{\alpha \beta \gamma *} \phi ^* _{\alpha \beta K}(\tau ' _{\alpha \beta} )}
  _-{\sim}
  &
  {p _{12K} ^{\alpha \beta \gamma *} p _{2 K}  ^{\alpha \beta *} \phi ^* _{\beta K} (E '_{\beta}) }
  \ar[r] ^-{\iota _\beta} _-{\sim}
  &
  {p _{12K} ^{\alpha \beta \gamma *} p _{2 K}  ^{\alpha \beta *} (E _{\beta}) }
  \ar[d] ^-{p _{12K} ^{\alpha \beta \gamma *} (\eta _{\alpha \beta})} _-{\sim}
  \ar[r] ^-{\epsilon} _-{\sim}
&
{p _{2K} ^{\alpha \beta \gamma*}  (E _\beta ),}
\ar[d] ^-{\eta _{12} ^{\alpha \beta \gamma }} _-{\sim}
\\
  {\phi ^* _{\alpha \beta \gamma K} ({E '_\alpha} |_{]X _{\alpha \beta \gamma}[ _{\fP _{\alpha \beta\gamma}}})}
  \ar[r] ^-{\epsilon} _-{\sim}
  &
   {p _{12K} ^{\alpha \beta \gamma *} \phi ^* _{\alpha \beta K}({E '_\alpha} |_{]X _{\alpha \beta}[ _{\fP _{\alpha \beta}}})}
   \ar[r]^-{\epsilon} _-{\sim}
   &
  {p _{12K} ^{\alpha \beta \gamma *} p  _{1 K} ^{\alpha \beta *}\phi ^* _{\alpha K}  (E ' _{\alpha})}
  \ar[r] ^-{\iota _\alpha} _-{\sim}
  &
  {p _{12K} ^{\alpha \beta \gamma *} p  _{1 K} ^{\alpha \beta *} (E _{\alpha})}
  \ar[r]^-{\epsilon} _-{\sim}
&
{p _{1K} ^{\alpha \beta \gamma*}(E _\alpha),}
   }
    \end{equation}
\normalsize
where $\phi _{\alpha \beta \gamma} : =(id,id, u _{\alpha \beta \gamma}) 
\colon 
(Y _{\alpha \beta \gamma}, X _{\alpha \beta \gamma}, \fX _{\alpha \beta \gamma})
\to 
(Y _{\alpha \beta \gamma}, X _{\alpha \beta \gamma}, \fP _{\alpha \beta \gamma})$ is the morphism of frames,
and where the middle rectangle is induced from \ref{diag2-eqcat-iso-reco} by 
applying the functor $ p _{12K} ^{\alpha \beta \gamma *}$.
We remark that the composition of the arrows of the bottom of  \ref{diag3-eqcat-iso-reco}
does not depends on the way to build it from $\iota _\alpha$ modulo $\epsilon$-isomorphisms, i.e.,
the following diagram is commutative:
\footnotesize
$$\xymatrix  @R=0,3cm@C=0,4cm {
  {\phi ^* _{\alpha \beta \gamma K} ({E '_\alpha} |_{]X _{\alpha \beta \gamma}[ _{\fP _{\alpha \beta\gamma}}})}
  \ar[r] ^(0.45){\epsilon} _(0.45){\sim}
  &
   {p _{12K} ^{\alpha \beta \gamma *} \phi ^* _{\alpha \beta K}({E '_\alpha} |_{]X _{\alpha \beta}[ _{\fP _{\alpha \beta}}})}
   \ar[r]^-{\epsilon} _-{\sim}
   &
  {p _{12K} ^{\alpha \beta \gamma *} p  _{1 K} ^{\alpha \beta *}\phi ^* _{\alpha K}  (E ' _{\alpha})}
  \ar[r] ^-{\iota _\alpha} _-{\sim}
  &
  {p _{12K} ^{\alpha \beta \gamma *} p  _{1 K} ^{\alpha \beta *} (E _{\alpha})}
  \ar[r]^-{\epsilon} _-{\sim}
&
{p _{1K} ^{\alpha \beta \gamma*}(E _\alpha)}
\\
  {\phi ^* _{\alpha \beta \gamma K} ({E '_\alpha} |_{]X _{\alpha \beta \gamma}[ _{\fP _{\alpha \beta \gamma}}})}
  \ar[r] ^(0.45){\epsilon} _(0.45){\sim}
  \ar@{=}[u]
  &
   {p _{13K} ^{\alpha \beta \gamma *} \phi ^* _{\alpha \gamma K}({E '_\alpha} |_{]X _{\alpha \gamma}[ _{\fP _{\alpha \gamma}}})}
   \ar[r]^-{\epsilon} _-{\sim}
   \ar[u] ^\epsilon _-{\sim}
   &
  {p _{13K} ^{\alpha \beta \gamma *} p  _{1 K} ^{\alpha \gamma *}\phi ^* _{\alpha K}  (E ' _{\alpha})}
  \ar[r] ^-{\iota _\alpha} _-{\sim}
  \ar[u] ^\epsilon _-{\sim}
  &
  {p _{13K} ^{\alpha \beta \gamma *} p  _{1 K} ^{\alpha \gamma *} (E _{\alpha})}
  \ar[r]^-{\epsilon} _-{\sim}
  \ar[u] ^\epsilon _-{\sim}
&
{p _{1K} ^{\alpha \beta \gamma*}(E _\alpha).}
\ar@{=}[u]
}$$
\normalsize
Similarly, we remark that so is the composition of the top arrows of  \ref{diag3-eqcat-iso-reco}.
Hence, by writing the two others analogous diagrams of  \ref{diag3-eqcat-iso-reco}, we establish that the family 
$(E ' _{\alpha}) _{\alpha \in \Lambda}$ glues to a
coherent    $j  ^\dag \cO _{]X _{\alpha }[ _{\fP _{\alpha}}}$-module $E '$
endowed with an integrable overconvergent isocrystal, i.e. 
to an object $E' \in \mathrm{MIC} ^{\dag} ((Y, X,\fP)/\cE ^{\dag}  _K) $.
Moreover, the isomorphisms $\iota _\alpha$ induise
  $u ^*  _{0K} (E ') \riso ((E _{\alpha})_{\alpha \in \Lambda},\, (\eta _{\alpha\beta}) _{\alpha ,\beta \in \Lambda})$.
\end{proof}

\begin{ntn}
We denote by  $\mathrm{MIC} ^{\dag \dag}  ((\fX ^{\sharp}  _\alpha )_{\alpha \in \Lambda}/\cE ^{\dag}  _K)$ the full subcategory of 
$\mathrm{Coh} ((\fX ^{\sharp} _\alpha) _{\alpha \in \Lambda}/\cE ^{\dag}  _K)$ whose objects $((\cE _{\alpha})_{\alpha \in \Lambda},\, (\theta _{\alpha\beta}) _{\alpha ,\beta \in \Lambda})$
are such that,  for all $\alpha \in \Lambda$, $\cE _{\alpha}$  is  $\widetilde{\cO} _{\fX _\alpha,\,\bbQ}$-coherent.
\end{ntn}

\begin{lem} \label{lem1pre-sp+plfid}
We have the canonical functor  $\sp _*\colon\mathrm{MIC} ^\dag ( (\fX   _\alpha )_{\alpha \in \Lambda}/\cE ^{\dag}  _K)
  \rightarrow   \mathrm{MIC} ^{\dag \dag}  (\fX ^{\sharp}  _\alpha )_{\alpha \in \Lambda}/\cE ^{\dag}  _K)$.
\end{lem}

\begin{proof}
Let  $((E _{\alpha})_{\alpha \in \Lambda},\, (\eta _{\alpha\beta}) _{\alpha ,\beta \in \Lambda})\in
  \mathrm{MIC} ^\dag ( (\fX   _\alpha )_{\alpha \in \Lambda}/\cE ^{\dag}  _K)$.
Let $\theta _{\alpha \beta}$ be the isomorphism making commutative the diagram
  \begin{equation}
    \label{defdonneesp*}
    \xymatrix  @R=0,3cm
{
    {\sp _*  p _{1 K}  ^{\alpha \beta !} (E _{\alpha})}
    \ar[r] _-{\sim} ^-{\ref{sp*f*com}}
    &
    {p _{1 }  ^{\alpha \beta !} \sp _*    (E _{\alpha})}
    \\
    {\sp _*  p _{2 K}  ^{\alpha \beta !} (E _{\beta})}
    \ar[r] _-{\sim} ^-{\ref{sp*f*com}}
    \ar[u] ^-{\sp _* \eta _{\alpha \beta}} _-{\sim}
    &
    {p _{2 }  ^{\alpha \beta !} \sp _*    (E _{\beta}).}
    \ar@{.>}[u]^{\theta _{\alpha \beta}}
}
  \end{equation}
We check that
  $\sp _* ((E _{\alpha})_{\alpha \in \Lambda},\, (\eta _{\alpha\beta}) _{\alpha ,\beta \in \Lambda})
  := ( (\sp _* E _{\alpha})_{\alpha \in \Lambda},\, (\theta _{\alpha\beta}) _{\alpha ,\beta \in \Lambda})$
  is an object of 
  $\mathrm{MIC} ^{\dag \dag}  (\fX ^{\sharp}  _\alpha )_{\alpha \in \Lambda}/\cE ^{\dag}  _K)$, i.e.,
the isomorphisms $\theta _{\alpha \beta}$
satisfy the cocycle condition. 
Consider the following commutative diagram:
\scriptsize
\begin{equation}
  \label{diag1-pre-sp+plfid}
  \xymatrix  @R=0,3cm {
  {\sp _*  p _{1K} ^{\alpha \beta \gamma!}(E _\alpha)}
  \ar[r] ^(0.45){\sp _* (\epsilon)} _(0.45){\sim}
  &
  {\sp _* p _{12 K} ^{\alpha \beta \gamma *}  p _{1 K}  ^{\alpha \beta *} (E _{\alpha})}
  \ar[r] _-{\sim} ^-{\ref{sp*f*com}}
  &
  {p _{12 } ^{\alpha \beta \gamma !}  \sp _* p _{1 K}  ^{\alpha \beta *} (E _{\alpha})}
  \ar[r] _-{\sim} ^-{\ref{sp*f*com}}
  &
  {p _{12 } ^{\alpha \beta \gamma !}  p _{1 }  ^{\alpha \beta !}  \sp _* (E _{\alpha})}
  \ar[r] ^-{\tau} _-{\sim}
  &
  {p _1 ^{\alpha \beta \gamma!}(\sp _* (E _{\alpha}))}
  \\
  {\sp _* p _{2K} ^{\alpha \beta \gamma!}  (E _\beta )}
  \ar[r] ^(0.45){\sp _* (\epsilon)} _(0.45){\sim}
  \ar[u] ^-{\sp _* (\eta _{12} ^{\alpha \beta \gamma })} _-{\sim}
  &
  {\sp _* p _{12 K} ^{\alpha \beta \gamma *}  p _{2 K}  ^{\alpha \beta *} (E _{\beta})}
  \ar[r] _-{\sim} ^-{\ref{sp*f*com}}
  \ar[u] ^-{\sp _* p _{12 K} ^{\alpha \beta \gamma *} (\eta _{\alpha \beta})} _-{\sim}
  &
  {p _{12 } ^{\alpha \beta \gamma !}  \sp _* p _{2 K}  ^{\alpha \beta *} (E _{\beta})}
  \ar[r] _-{\sim} ^-{\ref{sp*f*com}}
  \ar[u]^{p _{12 } ^{\alpha \beta \gamma !} \sp _* (\eta _{\alpha \beta})} _-{\sim}
  &
  {p _{12 } ^{\alpha \beta \gamma !}  p _{2 }  ^{\alpha \beta !}  \sp _* (E _{\beta})}
  \ar[u] ^-{p _{12 } ^{\alpha \beta \gamma !} (\theta _{\alpha \beta}) } _-{\sim}
  \ar[r] ^-{\tau} _-{\sim}
  &
  {p _2 ^{\alpha \beta \gamma!}  ( \sp _* (E _{\beta})).}
  \ar[u] ^-{\theta _{12} ^{\alpha \beta \gamma }} _-{\sim}
  }
\end{equation}
\normalsize
Since the isomorphisms of the form \ref{sp*f*com} are transitive with respect to the composition, 
then using 
\ref{sp-eps-tau}
we check that the composition of the arrows of the top
(resp. of the bottom) of \ref{diag1-pre-sp+plfid} is the canonical isomorphism
 $\sp _*  (p _{1 K}  ^{\alpha \beta \gamma } ) ^* (E _{\alpha})
\riso
  p _{1 }  ^{\alpha \beta \gamma!} \sp _* (E _{\alpha})$
  (resp.  $\sp _*  (p _{2 K}  ^{\alpha \beta \gamma } ) ^* (E _{\beta})
 \riso
  p _{2 }  ^{\alpha \beta \gamma!}  \sp _* (E _{\beta})$) of \ref{sp*f*com}.
By writing the two other analogous diagrams, we check the cocycle condition is satisfied.

Moreover, if $f = (f _\alpha)_{\alpha \in \Lambda}$ :
$((E _{\alpha})_{\alpha \in \Lambda},\, (\eta _{\alpha\beta}) _{\alpha ,\beta \in \Lambda})
\rightarrow ((E '_{\alpha})_{\alpha \in \Lambda},\, (\eta '_{\alpha\beta}) _{\alpha ,\beta \in \Lambda})$
is any morphism of 
$\mathrm{MIC} ^\dag ( (\fX   _\alpha )_{\alpha \in \Lambda}/\cE ^{\dag}  _K)$,
then, by functoriality of  \ref{defdonneesp*}, the family of morphisms
$\sp _* (f): = (\sp _* (f _\alpha))_{\alpha \in \Lambda}$ commute with the glueing data. 

\end{proof}

\begin{lem}
\label{lem2pre-sp+plfid}
We have  the canonical functor 
$\sp ^*\colon
\mathrm{MIC} ^{\dag \dag}  (\fX ^{\sharp}  _\alpha )_{\alpha \in \Lambda}/\cE ^{\dag}  _K)
\rightarrow 
\mathrm{MIC} ^{\dag }  (\fX   _\alpha )_{\alpha \in \Lambda}/\cE ^{\dag}  _K)$.
\end{lem}

\begin{proof}
Let   
$((\cE _{\alpha})_{\alpha \in \Lambda},\, (\theta _{\alpha\beta}) _{\alpha ,\beta \in \Lambda})\in 
\mathrm{MIC} ^{\dag \dag}  (\fX ^{\sharp}  _\alpha )_{\alpha \in \Lambda}/\cE ^{\dag}  _K)$.
Let  $\eta _{\alpha \beta}$ be the unique morphism making commutative the following diagram
  \begin{equation}
    \label{defdonneesp*2}
    \xymatrix  @R=0,3cm
{
    {\sp ^*  p _{1 }  ^{\alpha \beta !} (\cE _{\alpha})}
    &
    {p _{1 K}  ^{\alpha \beta *} \sp ^*    (\cE _{\alpha})}
    \ar[l] ^-{\sim}
    \\
    {\sp ^*  p _{2 }  ^{\alpha \beta !} (\cE _{\beta})}
    \ar[u] ^-{\sp ^* \theta _{\alpha \beta}} _-{\sim}
    &
    {p _{2 K}  ^{\alpha \beta *} \sp ^*    (\cE _{\beta}).}
    \ar@{.>}[u]^{\eta _{\alpha \beta}}
    \ar[l] ^-{\sim}
}
  \end{equation}
We set
$\sp ^* ((\cE _{\alpha})_{\alpha \in \Lambda},\, (\theta _{\alpha\beta}) _{\alpha ,\beta \in \Lambda})
:=
((\sp ^* \cE _{\alpha})_{\alpha \in \Lambda},\, (\eta _{\alpha\beta}) _{\alpha ,\beta \in \Lambda})$.
In the same way as the proof of \ref{lem1pre-sp+plfid},
(i.e. for instance we use \ref{sp*f*com} and \ref{sp-eps-tau}),
we check 
$((\sp ^* \cE _{\alpha})_{\alpha \in \Lambda},\, (\eta _{\alpha\beta}) _{\alpha ,\beta \in \Lambda}) \in
\mathrm{MIC} ^\dag ( (\fX   _\alpha )_{\alpha \in \Lambda}/\cE ^{\dag}  _K)$.
Finally, if 
  \small
  $f = (f _\alpha)_{\alpha \in \Lambda}$ :
$((\cE _{\alpha})_{\alpha \in \Lambda},\, (\theta _{\alpha\beta}) _{\alpha ,\beta \in \Lambda})
\rightarrow ((\cE '_{\alpha})_{\alpha \in \Lambda},\, (\theta '_{\alpha\beta}) _{\alpha ,\beta \in \Lambda})$
\normalsize
is a morphism, then by functoriality of  \ref{defdonneesp*2}, we prove that 
$\sp ^* (f) := (\sp ^* (f _\alpha))_{\alpha \in \Lambda}$ commute with glueing data.
\end{proof}

\begin{prop}
\label{pre-sp+plfid}
The functors  $\sp _*$ and $\sp ^*$ are quasi-inverse equivalences between the categories
$\mathrm{MIC} ^\dag ( (\fX   _\alpha )_{\alpha \in \Lambda}/\cE ^{\dag}  _K)$ and 
$\mathrm{MIC} ^{\dag \dag}  (\fX ^{\sharp}  _\alpha )_{\alpha \in \Lambda}/\cE ^{\dag}  _K)$.
\end{prop}

\begin{proof}
i) Let $((E _{\alpha})_{\alpha \in \Lambda},\, (\eta _{\alpha\beta}) _{\alpha ,\beta \in \Lambda})$
be an object of 
  $\mathrm{MIC} ^\dag ( (\fX   _\alpha )_{\alpha \in \Lambda}/\cE ^{\dag}  _K)$.
Let  $\theta _{\alpha \beta}$ (resp. $\eta ' _{\alpha \beta}$) be the canonical glueing  isomorphisms
of $\sp _* (E _{\alpha})$ (resp. $\sp ^* \sp _* (E _{\alpha})$).
It is enough to check that the adjunction isomorphisms
$\sp ^* \sp _* (E _\alpha) \riso E _\alpha$ commute with glueing data, i.e., that
the right square of the following diagram
\small
\begin{equation}\label{diag2-pre-sp+plfid}
  \xymatrix  @R=0,3cm {
    {p _{1 K}  ^{\alpha \beta !} (E _{\alpha})}
    &
    {\sp ^* \sp _*  p _{1 K}  ^{\alpha \beta !} (E _{\alpha})}
    \ar[r] _-{\sim}
    \ar[l] ^-{\sim}
    &
    {\sp ^* p _{1 }  ^{\alpha \beta !} \sp _*    (E _{\alpha})}
    &
    {p _{1 K}  ^{\alpha \beta *} \sp ^*  \sp _*   (E _{\alpha})}
    \ar[r] _-{\sim}
    \ar[l] ^-{\sim}
    &
    {p _{1 K}  ^{\alpha \beta *}  (E _{\alpha})}
    \\
    { p _{2 K}  ^{\alpha \beta !} (E _{\beta})}
    \ar[u] ^-{\eta _{\alpha \beta}} _-{\sim}
    &
    {\sp ^* \sp _*  p _{2 K}  ^{\alpha \beta !} (E _{\beta})}
    \ar[r] _-{\sim}
    \ar[l] ^-{\sim}
    \ar[u] ^-{\sp ^* \sp _* \eta _{\alpha \beta}} _-{\sim}
    &
    {\sp ^*  p _{2 }  ^{\alpha \beta !} \sp _*    (E _{\beta})}
    \ar[u]^{\sp ^* (\theta _{\alpha \beta})} _-{\sim}
    &
    {p _{2 K}  ^{\alpha \beta *} \sp ^*  \sp _*   (E _{\beta})}
    \ar[u]^{\eta ' _{\alpha \beta}} _-{\sim}
    \ar[r] _-{\sim}
    \ar[l] ^-{\sim}
    &
    {p _{2 K}  ^{\alpha \beta *} (E _{\beta}),}
    \ar[u] ^-{\eta _{\alpha \beta}} _-{\sim}
    }
\end{equation}
\normalsize
is commutative.
The commutativity of both middle square comes from respectively
 \ref{defdonneesp*} and \ref{defdonneesp*2},
whereas that of the left square is functorial.
By construction of both isomorphisms \ref{sp*f*com}, we check that the composition of the bottom (resp. top) arrowx of the diagram
\ref{diag2-pre-sp+plfid} are the identity map.
Since their arrows are isomorphisms, this yields that the diagram
\ref{diag2-pre-sp+plfid} is commutative.

ii) Let $((\cE _{\alpha})_{\alpha \in \Lambda},\, (\theta _{\alpha\beta}) _{\alpha ,\beta \in \Lambda})$ be an objet de
$\mathrm{MIC} ^{\dag \dag}  (\fX ^{\sharp}  _\alpha )_{\alpha \in \Lambda}/\cE ^{\dag}  _K)$.
In the same way as the step i), 
we check the adjunction isomorphisms
$ \cE _\alpha \riso \sp _*  \sp ^*(\cE _\alpha)$ commute with glueing data.
\end{proof}

\begin{ntn}\label{ntnMICdag2fs2}
We denote by $\mathrm{MIC} ^{\dag \dag} (Y, X,\fP ^\sharp/\cE ^{\dag}  _K) $ the full subcategory of 
$\mathrm{Coh} (X, \fP ^\sharp/\cE ^{\dag}  _K)$ whose objects $\cE$ satisfy the following condition:
for any affine open log formal subscheme  
 $\fP ^{\prime \sharp}$ of $\fP ^{\sharp}$, for any morphism of smooth $\fS ^{\sharp}$-formal log schemes
  $v$ : $ \fX ^{\prime \sharp} \hookrightarrow \fP ^{\prime \sharp} $ which reduces modulo $\pi$ to the closed imbedding 
  $X ^{ \sharp} \cap P ^{\prime \sharp} \hookrightarrow P ^{\prime \sharp}$,
the sheaf $v ^! (\cE |_{\fP^{\prime \sharp}}) $ is $\widetilde{\cO} _{\fX',\,\bbQ}$-coherent.
\end{ntn}

\begin{ntn} \label{ntnMICdag2fs3} 
We denote by  $\mathrm{MIC} ^{(\bullet)} (Y, X,\fP ^\sharp/\cE ^{\dag}  _K)$ the full subcategory of  $\smash{\underrightarrow{LM}}  _{\Q, \coh} (\widehat{\cD} _{\fP ^{\sharp}/\fS ^{\sharp}} ^{(\bullet)}(P _s))$
consisting of objects  $\cE ^{(\bullet)}$ with support in $X$ and such that  $\underrvec{l}_{\bbQ} ^{*}  (\cE ^{(\bullet)}) \in \mathrm{MIC} ^{\dag\dag} (Y, X,\fP ^\sharp/\cE ^{\dag}  _K)$ 
where  $\mathrm{MIC} ^{\dag \dag } (Y, X,\fP ^\sharp/\cE ^{\dag}  _K)$ where $\underrvec{l}_{\bbQ} ^{*} $ is defined in \ref{ntnMICdag2fs2}.  By definition, we get the equivalence of categories 
\begin{equation}\label{ntnMICdag2fs3-eq}
\underrvec{l}_{\bbQ} ^{*}  \colon \mathrm{MIC} ^{(\bullet)} (Y, X,\fP ^\sharp/\cE ^{\dag}  _K) \cong  \mathrm{MIC} ^{\dag \dag}(Y, X,\fP ^\sharp/\cE ^{\dag}  _K).
\end{equation}
\end{ntn}

\begin{empt}
The functors $u ^!  _0$ and $ u _{0+}$ constructed in respectively \ref{const-u0!} and \ref{const-u0+} 
induce quasi-inverse equivalence of categories between 
$\mathrm{MIC} ^{\dag \dag} (Y, X,\fP ^\sharp/\cE ^{\dag}  _K)$
and
$\mathrm{MIC} ^{\dag\dag} ( (\fX  ^{\sharp} _\alpha )_{\alpha \in \Lambda}/\cE ^{\dag}  _K)$, i.e., 
we have the commutative diagram
\begin{equation}
\label{eqcat-u0+!}
\xymatrix{
{\mathrm{MIC} ^{\dag \dag} (Y, X,\fP ^\sharp/\cE ^{\dag}  _K) } 
\ar[r] ^-{}
\ar@{.>}@<4ex>[d] ^-{\cong} _-{u _{0} ^!} 
& 
{\mathrm{Coh} (X, \fP ^\sharp/\cE ^{\dag}  _K)} 
\ar@{.>}@<4ex>[d] ^-{\cong} _-{u _{0} ^!} 
\\ 
{\mathrm{MIC} ^{\dag\dag} ( (\fX ^{\sharp}  _\alpha )_{\alpha \in \Lambda}/\cE ^{\dag}  _K)}
\ar[r] ^-{}
\ar@{.>}@<4ex>[u] ^-{\cong} _-{u _{0+}} 
& {\mathrm{Coh} ((\fX ^{\sharp} _\alpha) _\alpha/\cE ^{\dag}  _K).} 
\ar@{.>}@<4ex>[u] ^-{\cong} _-{u _{0+}} 
}
\end{equation}

\end{empt}

\begin{ntn}
\label{ntn-dfnsp+}
We get the canonical equivalence of categories 
$$\sp _+ \colon 
\mathrm{Isoc} ^{\dag} ((Y, X,\fP)/\cE ^{\dag}  _K)
\cong
\mathrm{MIC} ^{\dag \dag} (Y, X,\fP ^\sharp/\cE ^{\dag}  _K)$$
by composition of the equivalences
$$\mathrm{Isoc} ^{\dag} ((Y, X,\fP)/\cE ^{\dag}  _K)
\cong
\mathrm{MIC} ^{\dag} ((Y, X,\fP)/\cE ^{\dag}  _K)
\underset{u ^*  _{0K}}{\riso}
\mathrm{MIC} ^\dag ( (\fX   _\alpha )_{\alpha \in \Lambda}/\cE ^{\dag}  _K)
\underset{\sp _*}{\riso}
\mathrm{MIC} ^{\dag\dag} ( (\fX ^{\sharp}  _\alpha )_{\alpha \in \Lambda}/\cE ^{\dag}  _K)
\underset{u _{0+}}{\riso}
\mathrm{MIC} ^{\dag\dag} ((Y, X,\fP)/\cE ^{\dag}  _K).$$

\end{ntn}

\subsubsection{Commutation with duality}
\label{BK-dual}
\begin{empt}
\label{propspetdualsansfrob721}
With notation \ref{dfnMICdagalphabeta},
let $((E _{\alpha}) _{\alpha \in \Lambda}, (\eta _{\alpha \beta }) _{\alpha, \beta \in \Lambda})
\in 
\mathrm{MIC} ^{\dag}  (\fX   _\alpha )_{\alpha \in \Lambda}/\cE ^{\dag}  _K)$.
The $j ^{\dag} \cO _{\fX _{\alpha K}}$-linear dual 
of $E _{\alpha} $ is denoted by
$E _{\alpha} ^\vee : =
\mathcal{H} om  _{ j ^{\dag} \cO _{\fX _{\alpha K} }}( E ,  j ^{\dag} \cO _{\fX _{\alpha K}} )$. 
Since the $j ^\dag \cO$-linear dual commutes with pullbacks, 
the inverse of the isomorphism
$(\eta _{\alpha \beta}) ^\vee$ is canonically isomorphic a morphism of the form 
$p _{2K}  ^{\alpha \beta *} ((E _{\beta}) ^\vee) 
\riso 
p  _{1K} ^{\alpha \beta *} ( (E _{\alpha}) ^\vee)$ that we denote by $ \eta ^* _{  \alpha \beta}$.
To check the cocycle condition, by faithfulness of the restriction functor $j ^*$, 
we reduce to the known case \cite[4.3.1]{caro-construction}.
Hence, we get the dual functor 
$(-)^\vee \colon \mathrm{MIC} ^{\dag }  (\fX   _\alpha )_{\alpha \in \Lambda}/\cE ^{\dag}  _K)
\to \mathrm{MIC} ^{ \dag}  (\fX   _\alpha )_{\alpha \in \Lambda}/\cE ^{\dag}  _K)$
defined by setting 
$$((E _{\alpha}) _{\alpha \in \Lambda}, (\theta _{\alpha \beta }) _{\alpha, \beta \in \Lambda}) ^\vee
:= (( (E _{\alpha}) ^\vee) _{\alpha \in \Lambda}, (\theta ^*_{\alpha \beta }) _{\alpha, \beta \in \Lambda}).$$

Let $E \in \mathrm{Isoc} ^{\dag} ((Y, X,\fP)/\cE ^{\dag}  _K)$, and  $E ^\vee$ be its dual. We check the isomorphism
$$u _{0K} ^* (E ^\vee) \riso (u _{0K} ^* (E )  ) ^\vee,$$ 
 i.e. that the isomorphisms coming from the commutation of the dual with the pullbacks  are compatible with the glueing data, which is easy.
\end{empt}

\begin{empt}
Let $f \colon \fX ^{\prime \sharp} \to \fX ^{\sharp}$ be an open immersion of 
strictly semi-stable log formal schemes over $\fS ^\sharp$.
Let $\cE \in D ^{\mathrm{b}} _{\mathrm{coh}} ( \widetilde{\cD} ^\dag _{\fX ^\sharp/\fS ^\sharp})$.
As in \cite[3.2.8]{caro-construction}, we define the following isomorphism
\begin{gather}\notag
\xi \ :\   f ^!  \bbD   (\cE)
\riso
\R \mathcal{H} om _{\widetilde{\cD} ^\dag _{\fX ^{\prime \sharp}/\fS ^\sharp}}
(f ^{!} ( \cE),
f ^! _{\mathrm{d}}   (\widetilde{\cD} ^\dag _{\fX ^\sharp/\fS ^\sharp}  \otimes _{\widetilde{\cO} _{\fX}} \widetilde{\omega} _{\fX ^\sharp/\fS ^\sharp} ^{-1}))[d _X]
\\
\label{defDf!=f!D1bis}
\riso
\R \mathcal{H} om _{\widetilde{\cD} ^\dag _{\fX ^{\prime \sharp}/\fS ^\sharp}}
(f ^{!}  ( \cE),
(\widetilde{\cD} ^\dag _{\fX ^{\prime \sharp}/\fS ^\sharp} \otimes  _{\widetilde{\cO} _{\fX'}} \widetilde{\omega} _{\fX ^{\prime \sharp}/\fS ^\sharp} ^{-1} ) _\mathrm{t})[d _X]
\underset{\beta}{\riso}
\bbD   f ^!  (\cE),
\end{gather}
where $\beta$ is the transposition isomorphism exchanging both structures of left 
$\widetilde{\cD} ^\dag _{\fX ^{\prime \sharp}/\fS ^\sharp} \otimes  _{\widetilde{\cO} _{\fX'}} \widetilde{\omega} _{\fX ^{\prime \sharp}/\fS ^\sharp} ^{-1} $-modules of 
$\widetilde{\cD} ^\dag _{\fX ^{\prime \sharp}/\fS ^\sharp} \otimes  _{\widetilde{\cO} _{\fX'}} \widetilde{\omega} _{\fX ^{\prime \sharp}/\fS ^\sharp} ^{-1} $.
\end{empt}

\begin{empt}
With notation \ref{ntnPPalpha},
let 
$((\cE _{\alpha}) _{\alpha \in \Lambda}, (\theta _{\alpha \beta }) _{\alpha, \beta \in \Lambda})\in 
\mathrm{MIC} ^{\dag \dag} ( (\fX   _\alpha )_{\alpha \in \Lambda}/\cE ^{\dag}  _K)$.
Via the isomorphisms 
\ref{defDf!=f!D1bis}, the inverse of the isomorphism
$\bbD (\theta _{\alpha \beta})$ is canonically isomorphic to a morphism of the form 
$ \theta ^* _{  \alpha \beta} \ : \  p _2  ^{\alpha \beta !} (\bbD (\cE _{\beta})) \riso p  _1 ^{\alpha \beta !} (\bbD (\cE _{\alpha}))$ that we denote by  $ \theta ^* _{  \alpha \beta}$.
To check the cocycle condition, by faithfulness of the restriction functor $j ^*$, 
we reduce to the known case \cite[4.3.1]{caro-construction}.
Hence, we get the dual functor 
$\bbD \colon \mathrm{MIC} ^{\dag \dag} ( (\fX   _\alpha )_{\alpha \in \Lambda}/\cE ^{\dag}  _K)
\to \mathrm{MIC} ^{\dag \dag} ( (\fX   _\alpha )_{\alpha \in \Lambda}/\cE ^{\dag}  _K)$
defined by 
$\bbD ((\cE _{\alpha}) _{\alpha \in \Lambda}, (\theta _{\alpha \beta }) _{\alpha, \beta \in \Lambda})
:=
((\bbD (\cE _{\alpha})) _{\alpha \in \Lambda}, (\theta ^*_{\alpha \beta }) _{\alpha, \beta \in \Lambda})$.

\end{empt}

\begin{empt} \label{propspetdualsansfrob724}
With notation \ref{dfnMICdagalphabeta}, let $((E _{\alpha}) _{\alpha \in \Lambda}, (\eta _{\alpha \beta }) _{\alpha, \beta \in \Lambda}) \in  \mathrm{MIC} ^{\dag}  (\fX   _\alpha )_{\alpha \in \Lambda}/\cE ^{\dag}  _K)$. 
Following \ref{theoprincipal}, 
we have the canonical isomorphism
$\mathrm{sp} _* ( E _\alpha ^{\vee}) \riso \bbD (\mathrm{sp} _* (E _\alpha))$
of 
$\mathrm{MIC} ^{\dag \dag}  (\fX ^{\sharp}  _\alpha )_{\alpha \in \Lambda}/\cE ^{\dag}  _K)$.
To check that these isomorphisms are compatible with the glueing data, 
since $j ^*$ is faithful, then we reduce to the known case \cite[4.3.1]{caro-construction}.
Hence, we get the isomorphism 
$\mathrm{sp} _* 
(((E _{\alpha}) _{\alpha \in \Lambda}, (\eta _{\alpha \beta }) _{\alpha, \beta \in \Lambda}) ^\vee )
\riso
\bbD \circ \mathrm{sp} _*
((E _{\alpha}) _{\alpha \in \Lambda}, (\eta _{\alpha \beta }) _{\alpha, \beta \in \Lambda})$.
\end{empt}

\begin{empt}
\label{propspetdualsansfrob725}
With notation \ref{ntnPPalpha},
let 
$((\cE _{\alpha}) _{\alpha \in \Lambda}, (\theta _{\alpha \beta }) _{\alpha, \beta \in \Lambda})\in 
\mathrm{MIC} ^{\dag \dag} ( (\fX   _\alpha )_{\alpha \in \Lambda}/\cE ^{\dag}  _K)$.
From the duality isomorphism (see \ref{dualrelative}), 
we have the isomorphism
$u _{\alpha +} \circ \bbD (\cE _\alpha) 
\riso 
\bbD \circ u _{\alpha +}  (\cE _\alpha) 
\riso  
\bbD (u _{0+} ((\cE _{\alpha}) _{\alpha \in \Lambda}, (\theta _{\alpha \beta }) _{\alpha, \beta \in \Lambda})) |\fP _\alpha$.
By faithfulness of the restriction functor outside $j ^*$, to check theses isomorphisms glue
we reduce to the known case \cite[4.3.1]{caro-construction}.
Hence, we get 
the commutation isomorphism :
$$u _{0+} \circ \bbD ((\cE _{\alpha}) _{\alpha \in \Lambda}, (\theta _{\alpha \beta }) _{\alpha, \beta \in \Lambda}) 
\riso \bbD \circ u _{0+} ((\cE _{\alpha}) _{\alpha \in \Lambda}, (\theta _{\alpha \beta }) _{\alpha, \beta \in \Lambda})).$$
\end{empt}

\begin{prop}
\label{propspetdualsansfrob}
We keep notation \ref{ntnsp+}.
Let $E \in \mathrm{Isoc} ^{\dag} ((Y, X,\fP)/\cE ^{\dag}  _K)$,
and  $E ^\vee$ be its dual. We have the functorial canonical isomorphism
in $E$ : $ \sp _{+} (E ^\vee ) \riso
\widetilde{\bbD}  \circ \sp _{+} (E )$.
\end{prop}
\begin{proof}
Since $\sp _{+} 
=
u _{0+} \circ \sp _{*} \circ (Y, X,\fP ^\sharp/\cE ^{\dag}  _K)$, 
the proposition is a consequence of 
\ref{propspetdualsansfrob721},
\ref{propspetdualsansfrob724}, 
and
\ref{propspetdualsansfrob725}.
\end{proof}

\section{Differential coherence of the constant coefficient $\cO _{\fX} ({} ^\dag Z) _\bbQ$ when $Z$ is a divisor}
The purpose of this chapter is to prove Theorem \ref{coh-ss-div}. 
\subsection{Local cohomological functor with support in a strictly semistable weakly closed subscheme}

\subsubsection{Direct image by the specialization morphism of an overconvergent isocrystal when the boundary is not a divisor}
\label{directimagespovcv}
Let $\fX ^{\sharp}$ be a strictly semistable log formal scheme over $\fS ^\sharp$.
Let $\sp \colon \fX  _K \to \fX$ be the specialization morphism.
Let $Z$ be a closed subscheme of $X$ containing $X _s$ and let $Y:= X \setminus Z$ and $j\colon Y \subset X$ 
be the open immersion. 
Let $E$ be an abelian sheaf on $\fX _{K} $.
When $E$ is an overconvergent isocrystal on $(Y,X,\fX)/\mathcal{E} ^{\dag} _K$, i.e. 
an object of  $\mathrm{Isoc} ^\dag ((Y,X,\fX)/\mathcal{E} ^{\dag} _K)$, 
we  will check that  $\R \sp _* (E) \in D ^\mathrm{b}  ( \widetilde{\cD} ^\dag _{\fX ^{\sharp}/\fS ^{\sharp}})$.

We follow Berthelot's construction of Cech bicomplexes of \cite[4.1]{Be0}:
Let $\mathscr{X}:= (\fX _i) _{i\in I}$ be a finite affine open covering of $\fX$.  For any $i\in I$, let $\mathscr{Y} _{i}:= (Y _{i\,j _i}) _{j_i\in J _i}$ be a finite open covering of 
$Y _i:= Y \cap \fX _i$ such that  there exists $f _{i\, j _i}\in \Gamma (\fX _i, \cO_{\fX})$ satisfying $Y _{i\,j _i} = D ( f _{i\,j _i})\cap X _i$. We get the divisor  $T _{i\,j _i}:= V (f _{i\,j _i})$ 
of $X _i$ such that  $Y _{i\,j _i} = X _i \setminus T _{i\,j _i}$.  
\medskip 

Fix $h,l\in \bbN$. Let $\underline{i} =(i_0,\dots, i _h) \in I ^{1+h}$.  We set $\fX _{\underline{i}}:= \fX _{i _0}\cap \dots \cap \fX _{i _h}$, $Y _{\underline{i}}:= Y \cap \fX _{\underline{i}}$,
$u _{\underline{i}}  \colon    \fX _{\underline{i}} \to  \fX$, $u _{\underline{i} K}  \colon    \fX _{\underline{i} K}  \to  \fX _{K} $, and $J _{\underline{i}}:= J _{i _0}\times \dots \times J _{i _h}$.
For any $\underline{j} = (j _{i _0},\dots, j _{i _h})\in J _{\underline{i}}$, we set $Y _{\underline{i}\,\underline{j}}:= Y _{i _0\, j _{i _0}} \cap \dots \cap Y _{i _h\, j _{i _h}}$,
$f _{\underline{i}\,\underline{j}}:= f _{i _0\, j _{i _0}} | _{\fX _{\underline{i}}}  \cdots f _{i _h\, j _{i _h}} | _{\fX _{\underline{i}}}$. Denoting by 
$T _{\underline{i}\,\underline{j}} :=  V (f _{\underline{i}\,\underline{j}})$  the divisor of $X _{\underline{i}}$,  we have  $Y _{\underline{i}\,\underline{j}} = X _{\underline{i}} \setminus T _{\underline{i}\,\underline{j}}$.

We get the covering  $\mathscr{Y} _{\underline{i}}:= (Y _{\underline{i}\,\underline{j}}) _{\underline{j} \in J _{\underline{i}}}$ of $Y _{\underline{i}}$.
For any $\underline{\underline{j}} = (\underline{j} _0,\dots, \underline{j} _l)  \in (J _{\underline{i}} ) ^{1+l}$,  we set 
$Y _{\underline{i},\underline{\underline{j}}}:= Y _{\underline{i}\,\underline{j} _0}\cap \dots \cap Y _{\underline{i}\,\underline{j} _{l}}$,
$f _{\underline{i},\underline{\underline{j}}}: = f _{\underline{i}\,\underline{j} _0}\cdots f _{\underline{i}\,\underline{j} _{l}}$.
We denote the corresponding open immersions by 
$$j \colon Y \hookrightarrow \fX, \qquad j_{\underline{i}} \colon Y_{\underline{i}} \hookrightarrow \fX_{\underline{i}}, \qquad
j_{\underline{i} \underline{\underline{j}}} \colon Y_{\underline{i} \underline{\underline{j}}} \hookrightarrow \fX_{\underline{i}}.$$
Denoting by  $T _{\underline{i}\,\underline{\underline{j}}} :=  V (f _{\underline{i}\,\underline{\underline{j}}})$ 
the divisor of $X _{\underline{i}}$,  we have  $Y _{\underline{i}\,\underline{\underline{j}}} = X _{\underline{i}} \setminus T _{\underline{i}\,\underline{\underline{j}}}$.
As in \cite[4.1]{Be0} or \cite[12.1.1.3]{Car25}, we get the Cech bicomplexes 
$\check{C} ^{\dag \bullet \bullet} 
(\mathscr{X}, \mathscr{Y} _{\underline{i}},E)$ 
associated with the coverings 
$\mathscr{X}, \mathscr{Y} _{\underline{i}}$ of $E$ by setting 
$$\check{C} ^{\dag hl}  (\mathscr{X}, (\mathscr{Y} _{i} )_{i\in I},E):= \prod _{\underline{i} \in I ^{1+h}} u _{\underline{i} K *}  
\check{C} ^{\dag l}  ( \fX _{\underline{i}}, \mathscr{Y} _{\underline{i}}, u _{\underline{i}K} ^{*} (E) ) = \prod _{\underline{i} \in I ^{1+h}} u _{\underline{i} K*}   
\left ( \prod _{\underline{\underline{j}} \in J _{\underline{i}} ^{1+l}} j _{\underline{i},\underline{\underline{j}}}  ^{\dag}  u _{\underline{i}K } ^{*} (E)  \right), $$
We denote by  $\check{C} ^{\dag \bullet}  (\mathscr{X}, (\mathscr{Y} _{i} ) _{i\in I},E)$  the total complex of  $\check{C} ^{\dag \bullet \bullet}  (\mathscr{X}, (\mathscr{Y} _{i} ) _{i\in I},E)$.

\begin{prop}
Suppose $E$ is a coherent $j^\dag \cO_{\fX_K}$-module which is equipped with an integrable connection overconvergent along $X \setminus Y$. 
\begin{enumerate}[(a)]
\item The complex 
$\check{C} ^{\dag \bullet}  (\mathscr{X}, (\mathscr{Y} _{i} ) _{i\in I},E)$ gives a resolution of $E$.
\item For any $h,l\in \bbN$, the module $\check{C} ^{\dag hl}  (\mathscr{X}, (\mathscr{Y} _{i} )_{i\in I},  E )$  is acyclic for the functor $\sp _*$.

\item $\sp _* \check{C} ^{\dag \bullet}  (\mathscr{X}, (\mathscr{Y} _{i} )_{i\in I},E) $ is a complex of left $\widetilde{\cD} ^\dag _{\fX ^{\sharp}/\fS ^{\sharp}}$-modules.

\item We have in $D ^{\rmb} (\widetilde{\cD} ^\dag _{\fX ^{\sharp}/\fS ^{\sharp}})$ the canonical isomorphism
\begin{equation} \label{RspCheck-reso}
\R \sp _* (E)  \riso  \sp _* \check{C} ^{\dag \bullet}  (\mathscr{X}, (\mathscr{Y} _{i} )_{i\in I},E).
\end{equation}

\item Let $T$ be a closed subscheme of $X$ and $j ' \colon Y \setminus T \hookrightarrow X$ be the open immersion. We have the exact triangle in  $D ^{\rmb} (\widetilde{\cD} ^\dag _{\fX ^{\sharp}/\fS ^{\sharp}})$:
\begin{equation} \label{RspCheck-reso-functinY}
\bbR \sp _* (\underline{\Gamma}^\dag_{]T[} (E) ) \to \bbR \sp _* (E)\to \bbR \sp _* ( j ^{\prime \dag} E )\to +1.
\end{equation}
\end{enumerate}
\end{prop}

\begin{proof}
This is checked in the same way as \cite[4.1.4]{Be0} or  \cite[12.1.1.4]{Car25}.
For instance, when $I$ has only one element, we retrieve the exact sequence of \cite[2.50]{Lazda-Pal-Book}
and since $Y _{\underline{i},\underline{j}} = D (f_{\underline{i},\underline{j}})$, then 
we get  that $\check{C} ^{\dag hl}  (\mathscr{X}, \mathscr{Y} _{\underline{i}},E)$ is acyclic for the functor
$\sp _*$ and
$\sp _* ( \check{C} ^{\dag hl} 
(\mathscr{X}, \mathscr{Y} _{\underline{i}},E))$ is endowed with a structure of 
$\widetilde{\cD} ^\dag _{\fX ^{\sharp}/\fS ^{\sharp}}$-module.
\end{proof}

\subsubsection{Coherence via explicit computations}

\begin{lem}
\label{lem-internonzerodiv}
Let $C$ be an integral commutative ring, 
$x, y\in C$. 
If the image of $y$ in $C/xC$ is not a zero divisor then 
we have the equality
$C _{x} \cap C _{y} = C$.
\end{lem}

\begin{lem}
Let $X ^\sharp$ be a strictly semistable log schemes over $S ^\sharp$. Let $D _\eta$ be a divisor of $X _\eta$ and  $D$ be the closure of $D _\eta$ in $X$. 
Suppose $X= \Spec C$ affine  and suppose there exists $a\in C$ such that $D= V ( a) $.
Then $\overline{a} $ is not a zero divisor of $C / t C $, where $\overline{a}$ is the image of $a$ in $C/tC$.
\end{lem}

\begin{proof}
Let $X _1,\dots, X _r$ be the irreducible components of $X _s$.
Denote by $\mathfrak{p} _1, \dots, \mathfrak{p} _r$ the elements of
$\Spec C$ corresponding to 
the generic points $\eta _1, \dots, \eta _r$ of $X _1,\dots, X _r$.
Since $X _i$ is a divisor of $X$, since 
the irreducible components of $D$
corresponds bijectively to that of $D _\eta$, we have 
$X _i \not \subset D$. 
This yields $\eta _i \not \in D _s$.
Hence, 
$\overline{a} \not \in \mathfrak{p} _i/t\mathfrak{p} _i$. 
Since $\mathfrak{p} _1/t\mathfrak{p} _1,\dots, \mathfrak{p} _r/t\mathfrak{p} _r$
are the minimal prime ideals of $C/tC$, 
since  $C /t C$ is reduced (because $X _s$ is reduced),
then we get the primary decomposition $0= \mathfrak{p} _1/t\mathfrak{p} _1\cap \dots \cap \mathfrak{p} _r/t\mathfrak{p} _r$. 
This implies that the set of zero divisors of $C /t C$ is 
$\mathfrak{p} _1/t\mathfrak{p} _1 \cup \dots \cup \mathfrak{p} _r/t\mathfrak{p} _r$.
We are done.
\end{proof}

\begin{lem}
\label{prelemBe0421}
Let $\fX ^{\sharp}$ be a strictly semistable log formal scheme over $\fS ^\sharp$.
Suppose $\fX$ affine and $X$ irreducible. Put $A := \Gamma (\fX  ,\cO _{\fX})$,
$\overline{A}:= A / \pi A$. 
Let $f \in A$ whose image in $\overline{A} /t \overline{A}$ is not a zero divisor.
For any $c,d \in \R ^+$, with notation \ref{(c,d)bounded}, we have the equality
$A _{[t]} \cap A _{[tf]} (c,d)
=A _{[t]} (c,d)$, where $A _{[t]} (c,d)$ (resp. $A _{[tf]} (c,d)$) is the set of 
$(c,d)$-bounded elements of $A _{[t]}$ (resp. $A _{[tf]}$) as $A$-ring.
\end{lem}

\begin{proof}
The inclusion 
$A _{[t]} (c,d) 
\subset 
A _{[t]} \cap A _{[tf]} (c,d)$ is clear.
Let us check the converse.

1) We prove that  for any $x \in A _{[t]} \cap A _{[tf]} (c,d)$,
there exists $x _0 \in A _{t}$, $y _0 \in A _{[t]}$ such that 
$t ^d x _0 \in A$ and $x = x _0 + \pi y _0$.

Since $X$ is irreducible and regular,
then $\overline{A}$ is in particular integral.
We denote by 
$\overline{f}$ the image of $f$ in $\overline{A}$ and by
$\overline{x}$ the image of $x$ 
in 
$A _{[t]} /\pi A _{[t]} = \overline{A} _t  \subset \overline{A} _{t\overline{f}} 
=
A _{[tf]} /\pi A _{[tf]}$.
Since $x \in   A _{[tf]} (c,d)$, 
we get
$ (t\overline{f}) ^d \overline{x}\in \overline{A}$.
Since the image of $\overline{f}$ in 
$\overline{A} /t \overline{A}$ is not a zero divisor, we get
$t ^d \overline{x}\in \overline{A} _t \cap 
\overline{A} _{\overline{f}}
=\overline{A}$
(see \ref{lem-internonzerodiv}).
Hence there exists $a \in A$ such that 
$t ^d x -a \in \pi A _{[t]}$.
We set 
$x _0 := a / t ^d$.
Hence,  $x -x _0 \in \pi A _{[t]}$.

2) Let $x \in A _{[t]} \cap A _{[tf]} (c,d)$. 
To check that $x \in A _{[t]} (c,d) $, it is sufficient to 
 prove by induction on $n \geq 0$ that 
there exist 
$x _n \in A _t \cap A _{[t]} (c,d)$, 
$y _n \in A _{[t]} $ such that 
$x = x _n + \pi ^{n+1} y _n$.
The case $n= 0$ corresponds to the part 1). 
Suppose the induction checked for $n\geq0$.
Since $A _{[t]} (c,d) \subset A _{[tf]} (c,d)$,
we get,
$\pi ^{n+1} y _n = x -x_n \in A _{[tf]} (c,d)$.
Hence, 
$y _n \in A _{[t]} \cap A _{[tf]} (c,d + (n+1)c)$.
From the part 1),
there exists 
$z _{n+1} \in A _t$,
$y _{n+1} \in A _{[t]}$
such that 
$t ^{d + (n+1)c} z _{n+1} \in A$
and 
$y _n = z _{n+1} + \pi y _{n+1}$.
Hence $\pi ^{n+1}z _{n+1} \in  A _t \cap A _{[t]} (c,d)$ and 
then $x _{n+1}:=  x _n + \pi ^{n+1} z _{n+1} \in  A _t \cap A _{[t]} (c,d)$.
Since $x = x _{n+1} + \pi ^{n+2} y _{n+1}$, we conclude.

\end{proof}

The following Lemma will replace \cite[4.2.1]{Be0}
\begin{lem}
\label{lemBe0421}
We keep notation \ref{prelemBe0421}.
Let $(\widetilde{a} _i) _{i\in \bbN}$ be a sequence of $A _{[t]}$.
We suppose there exists two constants $c,d\geq 0$ such that 
$\widetilde{a} _i  \in A _{[t]} (c,d-i)$ for any $i\in \bbN$
and such that in $\widetilde{A} _{<f>} \coloneqq A _{[tf]}$ the following equality holds:
$$\sum _{i=0} ^{\infty} \widetilde{a} _i  f ^{-i-1}=0.$$
Then for any $j \geq 1$ we have
$\sum _{i =0} ^{j-1} \widetilde{a} _i  f ^{j-1-i} \in A _{[t]} (c,d+1-j)$.
\end{lem}

\begin{proof}
For $i \geq j$, we have
$\widetilde{a} _i  f ^{j-1-i} 
\in 
A _{[tf]} (c,d-i+ (i+1-j))
=
A _{[tf]} (c,d+1-j)$.
Using Lemma \ref{prelemBe0421}, 
this yields
$\sum _{i =0} ^{j-1} \widetilde{a} _i  f ^{j-1-i}
=
-\sum _{i \geq j}  
\widetilde{a} _i  f ^{j-1-i}
\in 
A _{[t]} \cap A _{[tf]} (c,d+1-j)
=A _{[t]} (c,d+1-j)$.

\end{proof}

\begin{lem}
\label{3.2.1Be0}
Let $\fX ^\sharp$ be an object of $\Sm ^\dag _{\fS ^\sharp}$.
Suppose  $\fX ^\sharp$ is affine, and $\fX ^\sharp/\fS ^\sharp$ has logarithmic coordinates. 
Let $t _1,\dots, t _d  \in \Gamma (\fX, \widetilde{\cO} _{\fX})$  be some overconvergent local coordinates of $\fX ^\sharp/\fS ^\sharp$. Set  $A := \Gamma (\fX, \cO _{\fX})$.

\begin{enumerate}[(a)]
\item Let  $P =  \sum _{\underline{i}\in \bbN ^d} \widetilde{a} _{\underline{i}} \underline{\partial} ^{<\underline{i}> _{(m)}} \in \Gamma (\fX, \widetilde{\cD} ^{(m)\dag}  _{\fX^\sharp/\fS^\sharp})$, 
where $\widetilde{a} _{\underline{i}} \in A _{[t]}$.
If  $\widetilde{a} _{\underline{i}} = 0$ for any $\underline{i}$ such that $i _1 =0$ then  there exists  $Q  \in  \Gamma (\fX, \widetilde{\cD} ^{(m)\dag}  _{\fX^\sharp/\fS^\sharp})$
such that  $p ^m P = Q \partial _1$.

\item 
Let $P = \sum _{\underline{i}\in \bbN ^d}  \widetilde{a} _{\underline{i}} \underline{\partial} ^{[\underline{i}]}  \in \Gamma (\fX, \widetilde{\cD} ^{\dag}  _{\fX^\sharp/\fS^\sharp})$
where $\widetilde{a} _{\underline{i}} \in A _{[t]}$.
If  $\widetilde{a} _{\underline{i}} = 0$ for any $\underline{i}$ such that $i _1 =0$ then  there exists  $Q   \in \Gamma (\fX, \widetilde{\cD} ^{\dag}  _{\fX^\sharp/\fS^\sharp,\bbQ})$
such that  $P = Q \partial _1$.
\end{enumerate}

\end{lem}

\begin{proof}
Following \ref{prop-wkcpdscplevelm},  are such that there exists a constant $c>0$ 
such that  $\widetilde{a} _{\underline{i}} \in A _{[t]} (c, c -|\underline{i}|)$.
We prove the first statement by using the computations of  \cite[3.2.1.1]{Be0}.
Taking the inductive limit on the level, the second statement is a consequence of the first one.
\end{proof}

\begin{prop}
\label{4.2.2Be0}
Let $\fX ^{\sharp}$ be a strictly semistable log formal scheme over $\fS ^\sharp$.
Let $u \colon Z ^\sharp \to X ^{\sharp}$ be a weakly closed immersion  of strictly semistable log schemes over $S ^\sharp$ such that $Z$ is a divisor of $X$.
Let $D \coloneqq  Z \cup X _s$. 
\begin{enumerate}[(a)]
\item The $\widetilde{\cD} ^\dag _{\fX ^{\sharp}/\fS ^{\sharp}}$-module 
$\cO _{\fX} (\hdag D ) _\bbQ$ is coherent ; 
\item If $\fU ^\sharp$ is an open affine formal subscheme of $\fX ^{\sharp}$ such that $\fU ^\sharp/\fS ^\sharp$ has logarithmic coordinates, 
for any overconvergent coordinates $t _1, \dots, t _d \in \Gamma (\fU ,\cO _{\fU})$ 
satisfying  $Z \cap U = V ( \overline{t} _1)$ where $\overline{t} _1$ is the  image of $t _1$ in $\Gamma (U ,\cO _{U})$, then 
the sequence
\begin{equation}
\label{4.2.2.1Be0}
(\widetilde{\cD} ^\dag _{\fU ^{\sharp}/\fS ^{\sharp}}) ^{d}
\overset{\psi}{\longrightarrow} 
\widetilde{\cD} ^\dag _{\fU ^{\sharp}/\fS ^{\sharp}}
\overset{\phi}{\longrightarrow}
\widetilde{\cO} _{\fU} (\hdag D _U) _\bbQ
\to 0,
\end{equation}
where 
$\phi (P)= P \cdot (1/t _1)$, and 
$\psi$ is defined by
\begin{equation}
\psi ( P _1,\dots, P _d) = P _1 \partial _1 t _1 + \sum _{i=2} ^d P _i \partial _i,
\end{equation}
is exact.
\end{enumerate}

\end{prop}

\begin{proof}
We can follow the proof of \cite[4.2.2]{Be0}. For the reader convenience, let us give few details.
We can suppose $\fU$ irreducible.
Since $\widetilde{\cD} ^\dag _{\fX ^{\sharp}/\fS ^{\sharp}}$ is coherent,  since such $\fU ^{\sharp}$ exists following \ref{ovcv-coord-comp-imm}), then  this is sufficient to check the second assertion.
Let $A:= \Gamma (\fU , \cO _{\fX})$. 
Then 
$A _{[t t _1],K}= \Gamma (\fU , \cO _{\fU} (\hdag D _U ) _\bbQ)$.
Let $f = \sum _{i=0} ^{\infty} a _i  (t t _1) ^{-i-1} \in A _{[t t _1]}$, such that 
$a _i \in A$ are such that there exists $c>0$ satisfying
$v _{p} ( a _{i}) \geq \frac{i}{c} -1$.
Using the description
\ref{wkcpdscp}, 
we check 
$P := \sum _{i\in \bbN}
\frac{(-1) ^i a _{i}}{t ^{i+1}}\partial _1 ^{[i]}
\in \Gamma (\fU, \widetilde{\cD} ^\dag _{\fX ^{\sharp}/\fS ^{\sharp}})$.
We compute 
$P \cdot \frac{1}{t _1} =f$. Hence, we get the surjectivity of $\phi$ of \ref{4.2.2.1Be0}.

Let $P \in \Gamma (\fU, \widetilde{\cD} ^\dag _{\fX ^{\sharp}/\fS ^{\sharp}})$
such that $P \cdot \frac{1}{t _1} =0$.
Using the description of Proposition \ref{loc-m2dagbis}.(\ref{loc-m2dagbis-(iii)}), we have 
$P = \sum _{\underline{i}\in \bbN ^d} \widetilde{a} _{\underline{i}}\underline{\partial} ^{[\underline{i}]}$,
with $\widetilde{a} _{\underline{i}} \in  A _{[t]} (c, c -| \underline{i} |)$.
Using Lemma \ref{3.2.1Be0}, in order to prove that $P$ is in the image of $\psi$ we can suppose 
$P = \sum _{i\in \bbN} \widetilde{a} _{i}\partial _1 ^{[i]}$, 
with $\widetilde{a} _{i} \in  A _{[t]} (c, c -i )$.
The equality $P \cdot \frac{1}{t _1} =0$ yields
$\sum _{i=0} ^{\infty} (-1) ^{i}\widetilde{a} _i  t _1 ^{-i-1}=0$.
Hence, using Lemma \ref{lemBe0421}, 
we get  $\widetilde{b} _j := (-1) ^{j+1} \sum _{i =0} ^{j-1} (-1) ^{i}\widetilde{a} _i  t _1 ^{-i-1} \in A _{[t]} (c,c+1-j)$.
Hence, 
$Q:= \sum _{j\geq 1} \widetilde{b} _j \partial _1 ^{[j]} \in \Gamma (\fU, \widetilde{\cD} ^\dag _{\fX ^{\sharp}/\fS ^{\sharp}})$.
Using the computations at the end of the proof of \cite[4.2.2]{Be0}, we get 
$Q t _1 = P$. 
\end{proof}

\begin{prop} \label{NCDgencoh}
Let $\fX ^{\sharp}$ be a strictly semistable log formal scheme over $\fS ^\sharp$.  Let $D $ be a divisor of $X $.
Suppose there exists $\fU ^\sharp$ an open affine formal subscheme of $\fX ^{\sharp}$ such that $\fU ^\sharp/\fS ^\sharp$ has logarithmic coordinates, there exist overconvergent coordinates  $t _1, \dots, t _d $  of $\fU ^\sharp/\fS ^\sharp$ satisfying  $t _1,\dots, t _r \in \Gamma (\fU ,\cO _{\fU})$ with $r \leq d$, and
$D \cap U  = U _s \cup V ( \overline{t} _1\cdots \overline{t} _r)$ where $\overline{t} _1, \dots, \overline{t} _r$ is the image of $t _1,\dots, t _r$ in $\Gamma (U ,\cO _{U})$. Then, we have  the exact sequence
\begin{equation}\label{NCDgencoh-exseq}
(\widetilde{\cD} ^\dag _{\fU ^{\sharp}/\fS ^{\sharp}}) ^{d} \overset{\psi}{\longrightarrow}  \widetilde{\cD} ^\dag _{\fU ^{\sharp}/\fS ^{\sharp}}\overset{\phi}{\longrightarrow} \cO _{\fU} (\hdag D ) _\bbQ \to 0,
\end{equation}
where  $\phi (P)= P \cdot (1/t _1\cdots t _r)$, and  $\psi$ is defined by 
\begin{equation} \psi ( P _1,\dots, P _d) = \sum _{i=1} ^{r} P _i \partial _i t _i + \sum _{i=r+1} ^d P _i \partial _i.
\end{equation}
\end{prop}

\begin{proof}
Using Lemmas \ref{lemBe0421} and \ref{3.2.1Be0},  we can follow the computations of \cite[4.3.2]{Be0} to conclude.
\end{proof}

\begin{rem}
In the special case where $(X,D \cup X _s)$ is a strict semi-stable pair (see definition \ref{dfn-ssp}), we will get a similar and more precise  statement at \ref{coh-ssp} (see remark \ref{ssp:exact-closed immersion}).
\end{rem}

\subsubsection{Local cohomological functor with support in a strictly semistable weakly closed subscheme}

\begin{empt}\label{ntn-GammaZ}
Let $\fX ^{\sharp}$ be a strictly semistable log formal scheme over $\fS ^\sharp$,
$j \colon X _\eta \hookrightarrow X$ be the open immersion.
Let $u \colon Z ^\sharp \to X ^{\sharp}$ be a weakly closed immersion  of strictly semistable log schemes over $S ^\sharp$.
Let $j _Z \colon  X  \setminus Z\to X$ be the open immersion. 
We set 
$$(\hdag Z) ( \widetilde{\cO} _{\fX,\bbQ}) :=\R \sp _* j _Z ^\dag (j ^\dag  \cO _{\fX _K})
\qquad \text{ and } \qquad 
\R \underline{\Gamma} ^\dag _Z \widetilde{\cO} _{\fX,\bbQ}  :=\R \sp _* \underline{\Gamma} ^\dag _Z (j ^\dag \cO _{\fX _K})$$
(see the definition of the functor $\underline{\Gamma} ^\dag _Z$ just before \cite[2.50]{Lazda-Pal-Book}).
The exact sequence 
$0 \to \underline{\Gamma} ^\dag _Z (j ^\dag \cO _{\fX _K}) \to 
j ^\dag  \cO _{\fX _K}
\to 
j _Z ^\dag (j ^\dag  \cO _{\fX _K})
\to 
0$
induces the exact triangle 
\begin{equation}
\label{extriangleloc}
\R \underline{\Gamma} ^\dag _Z \widetilde{\cO} _{\fX,\bbQ} 
\to 
\widetilde{\cO} _{\fX,\bbQ}
\to 
(\hdag Z) ( \widetilde{\cO} _{\fX,\bbQ}) 
\to 
\R \underline{\Gamma} ^\dag _Z \widetilde{\cO} _{\fX,\bbQ}  [1].
\end{equation}
For any integer $i\in \bbZ$, we set $\mathcal{H} ^{\dag i } _Z (\widetilde{\cO}_{\fX,\bbQ}  ):= H ^{i} \bbR \underline{\Gamma} ^\dag _Z \widetilde{\cO} _{\fX,\bbQ} $.
\end{empt}

\begin{prop}
\label{coh-smoothsubsch}
With notation \ref{ntn-GammaZ}, let $r$ be the codimension of the smooth closed subscheme $Z _\eta$ of $X _\eta$.
The following properties hold.
\begin{enumerate}[(a)]
\item $(\hdag Z) ( \widetilde{\cO} _{\fX,\bbQ}) ,
\R \underline{\Gamma} ^\dag _Z \widetilde{\cO} _{\fX,\bbQ}  
\in 
D ^{\mathrm{b}} _{\mathrm{coh}}
(\widetilde{\cD} ^\dag _{\fX ^{\sharp}/\fS ^{\sharp}})$, and
$\mathcal{H} ^{\dag i } _Z (\widetilde{\cO} _{\fX,\bbQ}  )=0$
for any $i\not = r$.

\item Let $x \in X$. 
Choose  an open affine formal subscheme $\fU ^\sharp$ of $\fX ^{\sharp}$ containing $x$
such that 
there exist overconvergent coordinates
$t _1, \dots, t _d \in \Gamma (\fU ,\cO _{\fU})$ 
such that 
$Z = V ( \overline{t} _1, \dots, \overline{t} _r)$ where
 $r \leq d$ and 
$\overline{t} _1, \dots, \overline{t} _r$ are the image of $t _1,\dots, t _r$ in $\Gamma (U ,\cO _{U})$
(recall such $\fU ^{\sharp}$ exists following \ref{ovcv-coord-comp-imm}).
We have the exact sequence
\begin{equation}
\label{exseq-HrZ}
(\widetilde{\cD} ^\dag _{\fU ^{\sharp}/\fS ^{\sharp}}) ^{d}
\overset{\psi}{\longrightarrow} 
\widetilde{\cD} ^\dag _{\fU ^{\sharp}/\fS ^{\sharp}}
\overset{\phi}{\longrightarrow}
\mathcal{H} ^{\dag r} _{Z \cap U} (\widetilde{\cO} _{\fU,\bbQ}  )
\to 0,
\end{equation}
where 
$\phi (P)= P \cdot (1/t _1\cdots t _r)$, and 
$\psi$ is defined by
\begin{equation}
\psi ( P _1,\dots, P _d) = \sum _{i=1} ^{r} P _i t _i + \sum _{i=r+1} ^d P _i \partial _i.
\end{equation}
\end{enumerate}
\end{prop}

\begin{proof}
It suffices to check the second statement.  We may assume, $\fU ^\sharp = \fX ^\sharp$. Let $Y = X \setminus Z$ and $j \colon Y \hookrightarrow X$ be the open immersion. 
For $i = 1,\dots, r$, put  $Y _i:  = D (\overline{t} _i)$.
For any $i _0,\dots, i _{k} \in \{ 1,\dots ,r\}$, let  $j _{i _0,\dots, i _{k}}\colon D ( \overline{t} _{i _0} \cdots \overline{t} _{i _k})=Y _{i _0}\cap \dots \cap Y _{i _k} \hookrightarrow  X$ be the open immersion.
Set $\widetilde{\cO}  _{\fX _K} \coloneqq j ^\dag  \cO _{\fX _K}$.
According to  \ref{RspCheck-reso} (applied in the case where $I$ is a singleton), $(\hdag Z) ( \widetilde{\cO} _{\fX,\bbQ}) =\R \sp _* j _Z ^\dag (\widetilde{\cO}  _{\fX _K} )$ is represented by the complex whose first term is at degree $0$
\begin{equation} \label{coh-smoothsubsch-es1}
\prod _{i=1} ^{d}
\sp _* ( j _i ^\dag \widetilde{\cO}  _{\fX _K})
\to 
\prod _{i _0 < i _1}
\sp _* ( j _{i _0 i _1} ^\dag \widetilde{\cO}  _{\fX _K})
\to
\dots
\to 
\sp _* ( j _{1\dots r} ^\dag \widetilde{\cO}  _{\fX _K})
\to 0,
\end{equation}
and then 
$\R \underline{\Gamma} ^\dag _Z \widetilde{\cO} _{\fX,\bbQ}$ is represented by the complex whose first term is at degree $0$
\begin{equation}
\label{coh-smoothsubsch-es2}
\widetilde{\cO}  _{\fX,\bbQ}
\to 
\prod _{i=1} ^{d}
\sp _* ( j _i ^\dag \widetilde{\cO}  _{\fX _K})
\to 
\prod _{i _0 < i _1}
\sp _* ( j _{i _0 i _1} ^\dag \widetilde{\cO}  _{\fX _K})
\to
\dots
\to 
\sp _* ( j _{1\dots r} ^\dag \widetilde{\cO}  _{\fX _K})
\to 0.
\end{equation}
According to \ref{NCDgencoh}, the terms of \ref{coh-smoothsubsch-es2} are $\widetilde{\cD} ^\dag _{\fX ^{\sharp}/\fS ^{\sharp}}$-coherent.
Hence, we get  $\R \underline{\Gamma} ^\dag _Z \widetilde{\cO} _{\fX,\bbQ}  \in  D ^{\mathrm{b}} _{\mathrm{coh}}  (\widetilde{\cD} ^\dag _{\fX ^{\sharp}/\fS ^{\sharp}})$.
We get the exact sequence \ref{exseq-HrZ} from the exact sequence \ref{NCDgencoh-exseq} and \ref{coh-smoothsubsch-es2} (and in the same way as the computations of \cite[4.3.4]{Be0}).
\end{proof}

\begin{coro}
\label{coro-trace-upre}
Let $u \colon \fZ ^\sharp \to \fX ^{\sharp}$ be a weakly closed immersion of strictly semistable log formal schemes over $\fS ^\sharp$.

\begin{enumerate}[(a)]
\item We have $u ^! (\hdag Z) ( \widetilde{\cO} _{\fX,\bbQ})=0$, i.e. 
by applying the functor $u ^!$ to the canonical morphism
$\R \underline{\Gamma} ^\dag _Z \widetilde{\cO} _{\fX,\bbQ}  
\to 
\widetilde{\cO} _{\fX,\bbQ}$, we get an isomorphism.

\item We have the isomorphism
$\R \underline{\Gamma} ^\dag _Z \widetilde{\cO} _{\fX,\bbQ}  
\riso 
u _+ u ^! (\widetilde{\cO} _{\fX,\bbQ})$.
\end{enumerate}

\end{coro}

\begin{proof}
 0) Since this is local in $\fX$ and $\fZ$, we can suppose $u$ is purely of codimension $r$.
 Let $\fY ^{\sharp}$ be the open formal subscheme 
of $\fX ^{\sharp}$  complementary to $Z$. 
First we check that 
$\R \underline{\Gamma} ^\dag _Z \widetilde{\cO} _{\fX,\bbQ} |\fY ^{\sharp} =0$.
Since this  is local, 
we can suppose with the notation of Proposition \ref{coh-smoothsubsch} that
$\fU ^\sharp=\fX ^{\sharp}$.
In that case, we compute that the restriction on $\fY ^{\sharp}$ of the map $\psi$ of the exact sequence \ref{exseq-HrZ}
 is surjective. Hence, 
 $\mathcal{H} ^{\dag r } _Z (\widetilde{\cO} _{\fX,\bbQ}  )  |\fY ^{\sharp} =0$
 (or we can also compute the $r$th cohomological space of the complex \ref{coh-smoothsubsch-es2}).

a) Now, let us check that 
$u ^! (\hdag Z) ( \widetilde{\cO} _{\fX,\bbQ})=0$. 
Since this is local, 
we can suppose with the notation of Proposition \ref{coh-smoothsubsch} that
$\fU ^\sharp=\fX ^{\sharp}$.
Denote by $D _{i _0 \dots i _k} $ the divisor of $X$ defined by $t _{i _0} \cdots t _{i _k}$,
using the exact sequence \ref{coh-smoothsubsch-es1}, it suffices to check 
$ \L u ^* (\hdag D _{i _0 \dots i _k}) ( \widetilde{\cO} _{\fX,\bbQ})=0$.
Let $u _i\colon \mathfrak{D} _i\hookrightarrow \fX$ be the closed immersion of formal schemes whose corresponding ideal is generated by
$t _i$. If we choose $i \in \{ i _0 ,\dots, i _k\}$, then the multiplication by
$t _i \colon (\hdag D _{i _0 \dots i _k}) ( \widetilde{\cO} _{\fX,\bbQ}) \to (\hdag D _{i _0 \dots i _k}) ( \widetilde{\cO} _{\fX,\bbQ})$ is an isomorphism.
Hence, $ \L u _i ^* (\hdag D _{i _0 \dots i _k}) ( \widetilde{\cO} _{\fX,\bbQ})=0$.
This yields $ \L u ^* (\hdag D _{i _0 \dots i _k}) ( \widetilde{\cO} _{\fX,\bbQ})=0$.

b) Hence, by applying the functor 
$u ^!$ to the exact triangle
\ref{extriangleloc}, 
we get the isomorphism
$u ^ !\R \underline{\Gamma} ^\dag _Z \widetilde{\cO} _{\fX,\bbQ} 
\riso 
u ^ ! \widetilde{\cO} _{\fX,\bbQ} $.
Since $\R \underline{\Gamma} ^\dag _Z \widetilde{\cO} _{\fX,\bbQ} $ has his support in $Z$, 
we get from the corollary \ref{Berthelot-Kashiwara-full} the first isomorphism
$\R \underline{\Gamma} ^\dag _Z \widetilde{\cO} _{\fX,\bbQ}
\riso 
u _+ u ^ !\R \underline{\Gamma} ^\dag _Z \widetilde{\cO} _{\fX,\bbQ} 
\riso 
u _+ u ^ !  \widetilde{\cO} _{\fX,\bbQ} $.

\end{proof}

\subsubsection{$\sp _+$ of the constant coefficient}

\begin{empt}
\label{com-invMIC-ntn} 
We keep notation \ref{ntnPPalpha}. Let $r$ be the codimension of the smooth closed subscheme $X _\eta$ of $P _\eta$.
Let $v \colon Y \to P _\eta$ be the morphism induced by $u$.
We get  the morphism of frames  $f :=(v,\,u,\,id) \colon  (Y, X,\fP) \to  (P _\eta, P,\fP)$.
\begin{enumerate}[(a)]
\item Denoting by
$f:=(v _0,\,u _0,\,id) \colon   (Y, X,\fP) \to (P _\eta, P,\fP)$ the morphism of frames, we get the inverse image 
$$f _K ^*   \colon   \mathrm{MIC} ^{\dag} ((P _\eta , P,\fP)/\cE ^{\dag}  _K)  \to      \mathrm{MIC} ^{\dag} ((Y, X,\fP)/\cE ^{\dag}  _K).$$ 
    Hence, we get the functor $u _{0K}  ^*\circ  | _{]X[ _\fP}   \colon     \mathrm{MIC} ^{\dag} ((P _\eta , P,\fP)/\cE ^{\dag}  _K)   \to   \mathrm{MIC} ^\dag ( (\fX   _\alpha )_{\alpha \in \Lambda}/\cE ^{\dag}  _K)$ 
(where $u _{0K}  ^*$ is the functor \ref{eqcat-iso-reco}).

\item In the same way as the construction of $u _0 ^! \colon  \mathrm{Coh} (X, \fP ^{\sharp}/\cE ^{\dag}  _K) \to  \mathrm{Coh} ((\fX ^{\sharp}  _\alpha )_{\alpha \in \Lambda}/\cE ^{\dag}  _K)$ of   \ref{const-u0!}, we define the functor  
$$u _0  ^*   \colon   \mathrm{MIC} ^{\dag \dag} ((P _\eta , P,\fP ^\sharp)/\cE ^{\dag}  _K)  \to  \mathrm{MIC} ^{\dag\dag} ( (\fX ^{\sharp}  _\alpha )_{\alpha \in \Lambda}/\cE ^{\dag}  _K)$$
as follows. Let $\cE \in  \mathrm{MIC} ^{\dag \dag} ((P _\eta , P,\fP ^\sharp)/\cE ^{\dag}  _K) $, i.e.      a coherent $\widetilde{\cD} ^\dag _{\fP ^\sharp/\fS ^\sharp ,\bbQ}$ which is also     $\widetilde{\cO}  _{\fP  ,\bbQ}$-coherent.    
We set $\cE _\alpha: = u _{\alpha } ^{*} ( \cE | \fP _\alpha):= H ^{r}u _{\alpha } ^{!} ( \cE | \fP _\alpha)  \riso u _{\alpha } ^{!} ( \cE | \fP _\alpha) [r]$.
Then  $\cE _\alpha$ is a coherent $\smash{\widetilde{\cD}} ^{\dag} _{\fX ^{\sharp} _{\alpha} /\fS ^\sharp,\bbQ} $-module,  which is also  $\widetilde{\cO}  _{\fX _{\alpha}  ,\bbQ}$-coherent.
Via the isomorphisms of the form $\tau $ (\ref{2.1.5Be2}),  we obtain the glueing  $\smash{\widetilde{\cD}} ^{\dag} _{\fX ^{\sharp} _{\alpha \beta} /\fS ^\sharp,\bbQ} $-linear
isomorphism $ \theta _{  \alpha \beta} \ : \  p _2  ^{\alpha \beta !} (\cE _{\beta}) \riso p  _1 ^{\alpha \beta !} (\cE _{\alpha}),$ satisfying the cocycle condition: $\theta _{13} ^{\alpha \beta \gamma }= \theta _{12} ^{\alpha \beta \gamma } \circ \theta _{23} ^{\alpha \beta \gamma }$
(this is checked as in the construction of $u _0 ^! \colon  \mathrm{Coh} (X, \fP ^{\sharp}/\cE ^{\dag}  _K) \to  \mathrm{Coh} ((\fX ^{\sharp}  _\alpha )_{\alpha \in \Lambda}/\cE ^{\dag}  _K)$ of   \ref{const-u0!}).

Beware we have three different functors $u _0 ^*$, $u _0 ^!$ and $u _{0K} ^*$ (defined at \ref{eqcat-iso-reco}) whose notation are very similar.

\end{enumerate}

\end{empt}

\begin{prop}
\label{spinvim*+}
With the notation \ref{com-invMIC-ntn},
we have the canonical isomorphism
$$\sp _* \circ u _{0 K} ^*
\circ  | _{]X[ _\fP} \riso
u _0 ^* \circ \sp _*$$
of functors
$\mathrm{MIC} ^{\dag} ((P _\eta, P,\fP)/\cE ^{\dag}  _K) \to  \mathrm{MIC} ^{\dag\dag} ( (\fX ^{\sharp}  _\alpha )_{\alpha \in \Lambda}/\cE ^{\dag}  _K)$.
\end{prop}

\begin{proof}
Using
\ref{sp-eps-tau},
we check that glueing data are compatible.
\end{proof}

\begin{coro}
\label{coro-sp+jdagO}
With notations \ref{ntn-GammaZ} and \ref{com-invMIC-ntn}, we have the isomorphism of $ \mathrm{Coh} (X, \fP ^{\sharp}/\cE ^{\dag}  _K) $ of the form
$$ \sp _+ (j ^\dag \cO  _{]X[ _{\fP}}) \riso \mathcal{H} ^{\dag ,r} _X \widetilde{\cO} _{\fP,\bbQ} \riso \R \underline{\Gamma} ^\dag _X \widetilde{\cO} _{\fP,\bbQ}[r].$$
\end{coro}

\begin{proof}We can copy the proof of \cite[12.2.3.3]{Car25}.
\end{proof}

\subsubsection{Application : trace morphism and relative biduality isomorphism for projective morphism of semistable formal log schemes}

\begin{empt}
\label{Tr-uclosedimmersion}
Let $u \colon \fZ ^\sharp \to \fX ^{\sharp}$ be a weakly closed immersion of strictly semistable log formal schemes over $\fS ^\sharp$.
Using \ref{coro-trace-upre}, 
from the canonical commutative morphism
$\R \underline{\Gamma} ^\dag _Z \widetilde{\cO} _{\fX,\bbQ}  
\to 
\widetilde{\cO} _{\fX,\bbQ}  $
we get the canonical one
$u _+ u ^! (\widetilde{\cO} _{\fX,\bbQ}) \to  \widetilde{\cO} _{\fX,\bbQ}  $ making commutative the diagram
\begin{equation}  \label{Tr-uclosedimmersion-diag} \xymatrix{ {u _+ u ^! (\widetilde{\cO} _{\fX,\bbQ})} \ar@{.>}[r] ^-{} & {\widetilde{\cO} _{\fX,\bbQ} } 
\\ {u _+ u ^! (\R \underline{\Gamma} ^\dag _Z \widetilde{\cO} _{\fX,\bbQ}  )} \ar[u] ^-{\sim} _-{\ref{coro-trace-upre}}\ar[r] ^-{\sim} _-{\ref{Berthelot-Kashiwara-full}}
& {\R \underline{\Gamma} ^\dag _Z \widetilde{\cO} _{\fX,\bbQ}  }   \ar[u] ^-{}.} 
\end{equation}
Since 
$u ^! (\widetilde{\cO} _{\fX,\bbQ}  ) [d _{X }]
\riso 
\widetilde{\cO} _{\fZ,\bbQ}   [d _{Z }]$,  
this yields the map:
\begin{equation}\label{Tr-uclosedimmersion-ntn} \Tr _{+,u} 
\colon u _{ +}  (\widetilde{\cO} _{\fZ,\bbQ} [d _{Z }])
\to 
\widetilde{\cO} _{\fX,\bbQ} [d _{X }].
\end{equation}
Since $\eta = \Spf \cV [[t]] \{ \frac{1}{t}\}$
is the formal spectrum of a complete discrete valuation ring of mixed characteristic 
$(0,p)$,
then we can use the results of \cite{Vir04} in the context of the smooth formal scheme $\fX _\eta/\eta$. 
Since Virrion' trace map  $\Tr _{+,u _{\eta}} \colon  u _{\eta +}  (\cO  _{\fZ _\eta,\bbQ} [d _{Z _\eta}]) \to  \cO  _{\fX _\eta,\bbQ} [d _{X _\eta}]$ (see \cite{Vir04})
satisfies the commutative diagram \ref{Tr-uclosedimmersion-diag} over $\fX _\eta$, then  
the restriction  over $\fX _\eta$ of the trace map $\Tr _{+,u} $ of \ref{Tr-uclosedimmersion-ntn} 
is equal to $\Tr _{+,u _\eta}$.
\end{empt}

\begin{lem}
\label{Tr+trans1lem}
Let $u \colon \fY ^\sharp  \to \fX ^{\sharp} $ 
and
$v \colon \fZ ^\sharp  \to \fY ^{\sharp} $ 
be two weakly closed immersions of strictly semistable log formal schemes over $\fS ^\sharp$. 
Then the diagram
\begin{equation}
\xymatrix @C=2cm {
{u _+ \circ v _{ +}  (\widetilde{\cO} _{\fZ ^\sharp ,\bbQ} [d _Z])} 
\ar[r] ^-{u _+ (\Tr _{+,v})}
\ar[d] ^-{\sim}
& 
{u _+  (\fY ^{\sharp} [d _Z])}
\ar[d] ^-{\Tr _{+,u}} 
\\ 
{(u \circ v)  _+    (\widetilde{\cO} _{\fZ ^\sharp ,\bbQ} [d _Z])} 
\ar[r] ^-{\Tr _{+,u \circ v}} 
& 
{\fX ^{\sharp} [d _X]} 
}
\end{equation}

is commutative. 
\end{lem}

\begin{proof}
Trace maps over the generic fibers are equal to 
Virrion' trace map
\cite{Vir04}.
We conclude by using the transitivity of Virrion's trace map.
\end{proof}

\begin{empt}
\label{Tr+closed-immersion}
Let $\fX$ be an object of 
$\Sm ^\dag _{\fS ^\sharp}$.
We denote by 
$ \widehat{\P} _\cV ^N:= 
 \Spf \cV \{ T _0,\dots, T _N\}$
the $N$-th projective $p$-adic formal space over $\cV$.
Let $\widehat{\P} ^N _{\fX ^\sharp} 
:= \widehat{\P} _\cV ^N \times _{\Spf \cV} \fX ^\sharp$,
and 
$f \colon 
 \widehat{\P} _\cV ^N
\to 
\Spf \cV$ 
and 
$g \colon 
\widehat{\P} ^N _{\fX ^\sharp} 
\to 
\fX  ^\sharp$ be the projections.
Following 
Virrion' trace map
\cite{Vir04}, 
we have the trace map 
$f _+ (\cO _{ \widehat{\P} _\cV ^N ,\bbQ} [N])
\to 
\cO _{\Spf \cV , \bbQ}$. 
Using base change Theorems \cite[5.3.3.1 and 9.4.4.3]{Car25},
this yields the trace map
$g _+ (\cO _{\widehat{\P} ^N _{\fX ^\sharp} ,\bbQ} [d _P])
\to 
\cO _{\fX ^\sharp} [d _X]$. 
Using \cite[9.4.3.3]{Car25}, 
we get  $\Tr _{+,g}
\colon 
g _+ (\widetilde{\cO} _{\widehat{\P} ^N _{\fX ^\sharp}  ,\bbQ} [d _P])
\to 
\widetilde{\cO} _{\fX ^\sharp} [d _X]$.
The restriction  over $\fX _\eta$ of  this morphism, $g _{\eta +} (\cO _{\widehat{\P} ^N _{\fX _\eta}    ,\bbQ} [d _P]) \to  \cO _{\fX _\eta} [d _X]$  corresponds to Virrion' trace map \cite{Vir04}.
\end{empt}

\begin{lem}
\label{Tr+trans2lem}
Let $\fX ^{\sharp}$ be a strictly semistable log formal schemes over $\fS ^\sharp$, 
$u \colon \widehat{\P} ^N _{\fX ^\sharp}  \to \widehat{\P} ^{M} _{\fX ^\sharp} $ 
be a weakly closed immersion such that $g \circ u =f$, where 
$f \colon  \widehat{\P} ^N _{\fX ^\sharp}  \to \fX ^\sharp $ and  $g \colon  \widehat{\P} ^{M} _{\fX ^\sharp}  \to  \fX ^\sharp$ are the projections.
Then the diagram
\begin{equation}
\xymatrix @C=2cm {
{g _+ \circ u _{ +}  (\widetilde{\cO} _{\widehat{\P} ^N _{\fX ^\sharp} ,\bbQ} [N])} 
\ar[r] ^-{g _+ (\Tr _{+,u})}
\ar[d] ^-{\sim}
& 
{g _+  (\widetilde{\cO} _{\widehat{\P} ^M _{\fX ^\sharp} ,\bbQ} [M])}
\ar[d] ^-{\Tr _{+,g}} 
\\ 
{f _+  (\widetilde{\cO} _{\widehat{\P} ^N _{\fX ^\sharp} ,\bbQ} [N])} 
\ar[r] ^-{\Tr _{+,f}} 
& 
{\widetilde{\cO} _{\fX,\bbQ}} 
}
\end{equation}

is commutative. 

\end{lem}

\begin{proof}
Trace maps over the generic fibers are equal to 
Virrion' trace map
\cite{Vir04}.
We conclude by using the transitivity of Virrion's trace map.
\end{proof}

Because of the relative duality isomorphism of the form \ref{dualrelative}, which is not known in a more general case, 
we will need to focus on pushforwards by projective morphisms: 
\begin{dfn} \label{Tr+PNdef}
A morphism $f \colon \fY ^\sharp \to \fX ^{\sharp}$ of strictly semistable log formal schemes over $\fS ^\sharp$ is said to be  ``projective'' if there exists a closed immersion $\underline{u} \colon \fY  \to \widehat{\P} _\cV ^N \times _{\Spf \cV} \fX $  
such that $\underline{f}$ is the composition of the projection $\underline{g} \colon  \widehat{\P} _\cV ^N \times _{\Spf \cV} \fX \to \fX $ with $\underline{u}$.
If $f$ is projective then there exists a weakly closed immersion of strictly semistable log formal schemes over $\fS ^\sharp$ (see definition \ref{dfn-wclimm}) of the form
$u \colon \fY ^\sharp \to \widehat{\P} _\cV ^N \times _{\Spf \cV} \fX ^\sharp$ such that $f$ is the composition of the projection  $g  \colon  \widehat{\P} _\cV ^N \times _{\Spf \cV} \fX ^\sharp  \to  \fX ^\sharp$ with $u$.
\end{dfn}

\begin{lem} \label{Tr+PN}
Let $f \colon \fY ^\sharp \to \fX ^{\sharp}$ be a projective morphism of strictly semistable log formal schemes over $\fS ^\sharp$.
There exists a canonical trace morphism  $\Tr _{+,f} \colon  f _+ (\widetilde{\cO} _{\fY ^\sharp ,\bbQ} [d _Y]) \to  \widetilde{\cO} _{\fX ^\sharp} [d _X]$.
The restriction over $\fX _\eta$ of  this morphism is equal to Virrion's trace map \cite{Vir04}.
\end{lem}

\begin{proof}
Let $u \colon \fY ^\sharp \to \widehat{\P} _\cV ^N \times _{\Spf \cV} \fX ^\sharp$
be a weakly closed immersion of strictly semistable log formal schemes over $\fS ^\sharp$
such that
$f$ is the composition of the projection 
$g 
\colon 
\widehat{\P} _\cV ^N \times _{\Spf \cV} \fX ^\sharp 
\to 
\fX ^\sharp$
with $u$.
Since 
$g _+ \circ u _+ 
(\widetilde{\cO} _{\fY ^\sharp ,\bbQ})
\riso
f _+ (\widetilde{\cO} _{\fY ^\sharp ,\bbQ})$, we get from 
\ref{Tr-uclosedimmersion} 
and 
\ref{Tr+closed-immersion}
the construction of $\Tr _{+,f}$ by composition.

By using Lemmas \ref{Tr+trans1lem} and \ref{Tr+trans2lem},
we check that  this morphism
$\Tr _{+,f}$ does not depend on the choice of the splitting
$f = g \circ u$.

\end{proof}

\begin{prop}
\label{Tr+transprop}
Let $f \colon \fY ^\sharp  \to \fX ^{\sharp} $ 
and
$g \colon \fZ ^\sharp  \to \fY ^{\sharp} $ 
be two projective morphisms of strictly semistable log formal schemes over $\fS ^\sharp$.
Then the diagram
\begin{equation}
\label{Tr+transpropdiag}
\xymatrix @C=2cm {
{f _+ \circ g _{ +}  (\widetilde{\cO} _{\fZ ^\sharp ,\bbQ} [d _Z])} 
\ar[r] ^-{f _+ (\Tr _{+,g})}
\ar[d] ^-{\sim}
& 
{f _+  (\fY ^{\sharp} [d _Z])}
\ar[d] ^-{\Tr _{+,f}} 
\\ 
{(f \circ g)  _+    (\widetilde{\cO} _{\fZ ^\sharp ,\bbQ} [d _Z])} 
\ar[r] ^-{\Tr _{+,f \circ v}} 
& 
{\fX ^{\sharp} [d _X]} 
}
\end{equation}

is commutative. 
\end{prop}

\begin{proof}
Trace maps over the generic fibers are equal to 
Virrion' trace map
\cite{Vir04}.
We conclude by using the transitivity of Virrion's trace map.
\end{proof}

\begin{thm}
\label{dualrelative}
Let $f \colon \fY ^\sharp \to \fX ^{\sharp}$ be a projective morphism of strictly semistable log formal schemes over $\fS ^\sharp$.
Let $\cF \in D ^\mathrm{b} _{\mathrm{coh}}
(\widetilde{\cD} ^{\dag} _{\fY ^{\sharp}/\fS ^{\sharp},\bbQ}).$
We have in 
$D ^\mathrm{b} _{\mathrm{coh}}
(\widetilde{\cD} ^{\dag} _{\fX ^{\sharp}/\fS ^{\sharp},\bbQ})$
the isomorphism
\begin{equation}
\label{dualrelative-morp}
f _+ \circ \bbD ( \cF) 
\riso 
\bbD \circ f _+ (\cF).
\end{equation}
\end{thm}

\begin{proof}
Using the trace map \ref{Tr+PN}, 
we construct the morphism \ref{dualrelative-morp} in the same way as 
that of \cite[IV.1.3]{Vir04} (for more details, see also 
\cite[1.2.7]{caro_courbe-nouveau}).
The restriction of this morphism \ref{dualrelative-morp} to $\fX _\eta$ is the same as that 
of Virrion in \cite[IV]{Vir04}. Hence, the restriction to  $\fX _\eta$ of this morphism is an isomorphism and
then by faithfulness this is an isomorphism.
 \end{proof}

\begin{coro}
\label{cor-adj-formul}
Let $f \colon \fY ^\sharp \to \fX ^{\sharp}$ be a projective morphism of strictly semistable log formal schemes over $\fS ^\sharp$.
Let $\cF \in D ^\mathrm{b} _{\mathrm{coh}}
(\widetilde{\cD} ^{\dag} _{\fY ^{\sharp}/\fS ^{\sharp},\bbQ})$,
and
$\cE \in D ^\mathrm{b} _{\mathrm{coh}}
(\widetilde{\cD} ^{\dag} _{\fX ^{\sharp}/\fS ^{\sharp},\bbQ})$.
We have  the isomorphisms
\begin{gather}
\label{cor-adj-formul-bij1}
\R \mathcal{H} om _{\widetilde{\cD} ^{\dag} _{\fY ^{\sharp}/\fS ^{\sharp},\bbQ}}
( f _+ ( \cF) , \cE) 
\riso 
\R f _* 
\R \mathcal{H} om _{\widetilde{\cD} ^{\dag} _{\fX ^{\sharp}/\fS ^{\sharp},\bbQ}}
( \cF ,f ^!  ( \cE)). 
\\
\label{cor-adj-formul-bij2}
\R \mathrm{Hom}  _{\widetilde{\cD} ^{\dag} _{\fY ^{\sharp}/\fS ^{\sharp},\bbQ}}
( f _+ ( \cF) , \cE) 
\riso 
\R \mathrm{Hom}  _{\widetilde{\cD} ^{\dag} _{\fX ^{\sharp}/\fS ^{\sharp},\bbQ}}
( \cF ,f ^!  ( \cE)). 
\end{gather}
\end{coro}

\begin{proof}
The proof is identical to that of 
\cite[IV.4.1 and IV.4.2]{Vir04}: this is a formal consequence of 
the relative duality isomorphism \ref{dualrelative}.
\end{proof}

\begin{coro}
\label{adj-morph}
Let $f \colon \fY ^\sharp \to \fX ^{\sharp}$ be a projective morphism of strictly semistable log formal schemes over $\fS ^\sharp$.
\begin{enumerate}[(a)]
\item 
Let 
$\cF 
\in 
D ^\mathrm{b} _{\mathrm{coh}}
(\widetilde{\cD} ^{\dag} _{\fY ^{\sharp}/\fS ^{\sharp},\bbQ})$. 
We have the adjunction morphism
$\cF \to f ^! f _+ (\cF)$. 
\item Let $\cE \in D ^\mathrm{b} _{\mathrm{coh}}
(\widetilde{\cD} ^{\dag} _{\fX ^{\sharp}/\fS ^{\sharp},\bbQ})$
such that 
$f ^! (\cE) \in 
D ^\mathrm{b} _{\mathrm{coh}}
(\widetilde{\cD} ^{\dag} _{\fY ^{\sharp}/\fS ^{\sharp},\bbQ})$.
We have the adjunction morphism 
$f _+ f ^! (\cE) \to \cE$. 

\item Suppose $f$ is moreover smooth. 
Then $f _+ \colon 
D ^\mathrm{b} _{\mathrm{coh}}
(\widetilde{\cD} ^{\dag} _{\fY ^{\sharp}/\fS ^{\sharp},\bbQ})
\to 
D ^\mathrm{b} _{\mathrm{coh}}
(\widetilde{\cD} ^{\dag} _{\fX ^{\sharp}/\fS ^{\sharp},\bbQ})$
is a right adjoint functor of 
$f ^!
\colon 
D ^\mathrm{b} _{\mathrm{coh}}
(\widetilde{\cD} ^{\dag} _{\fX ^{\sharp}/\fS ^{\sharp},\bbQ})
\to 
D ^\mathrm{b} _{\mathrm{coh}}
(\widetilde{\cD} ^{\dag} _{\fY ^{\sharp}/\fS ^{\sharp},\bbQ})$.
\end{enumerate}
\end{coro}

\subsection{Differential coherence of the constant coefficient}

\subsubsection{The case of a strict semi-stable pair}

\begin{prop} \label{coh-log-smoothcase}
Let $u \colon \fZ ^\sharp \to \fX ^{\sharp}$ be an exact closed immersion of strictly semistable log formal schemes over $\fS ^\sharp$.
We suppose furthermore that $Z$ is a divisor of $X$. 
Then 
$\cB _{\fX}^{(\bullet)} (Z) 
\in \smash{\underrightarrow{LM}}  ^\mathrm{b} _{\bbQ, \mathrm{coh}}
(\overset{^\mathrm{l}}{} \smash{\widehat{\cD}} _{\fX ^{\sharp}/\fS ^{\sharp} } ^{(\bullet)})$
and we have the exact triangle in 
$\smash{\underrightarrow{LD}}  ^\mathrm{b} _{\bbQ, \mathrm{coh}}
(\overset{^\mathrm{l}}{} \smash{\widehat{\cD}} _{\fX ^{\sharp}/\fS ^{\sharp}} ^{(\bullet)})$:
\begin{equation}
u _{+} ^{ (\bullet)} \circ  u ^{ (\bullet)!}  (\cO _{\fX})
\to 
\cO _{\fX}
\to 
\cB _{\fX}^{(\bullet)} (Z)
\to 
u _{+} ^{ (\bullet)} \circ  u ^{ (\bullet)!}  (\cO _{\fX})
[1].
\end{equation}

\end{prop}

\begin{proof}
The check that  $\cB _{\fX}^{(\bullet)} (Z)  \in \smash{\underrightarrow{LM}}  ^\mathrm{b} _{\bbQ, \mathrm{coh}} (\overset{^\mathrm{l}}{} \smash{\widehat{\cD}} _{\fX ^{\sharp}/\fS ^{\sharp}} ^{(\bullet)})$
is identical to the proof of the remark  \cite[B.3]{caro-stab-sys-ind-surcoh}. Hence, the mapping cone $\cC$ of  $\cO _{\fX} \to  \cB _{\fX}^{(\bullet)} (Z)$ is a complex in 
$ \smash{\underrightarrow{LD}}  ^\mathrm{b} _{\bbQ, \mathrm{coh}} (\overset{^\mathrm{l}}{} \smash{\widehat{\cD}} _{\fX ^{\sharp}/\fS ^{\sharp} } ^{(\bullet)})$.
Since this morphism is an isomorphism outside $Z$,
then using Berthelot-Kashiwara's theorem (see \cite[9.3.5.13]{Car25}),
we get $\cC \riso u _{+} ^{ (\bullet)} \circ  u ^{ (\bullet)!}  (\cC)$. Since we have in $\smash{\underrightarrow{LM}}  ^\mathrm{b} _{\bbQ, \mathrm{coh}}
(\overset{^\mathrm{l}}{} \smash{\widehat{\cD}} _{\fX ^{\sharp}/\fS ^{\sharp}} ^{(\bullet)})$ the isomorphism  $u ^{ (\bullet)!} (\cB _{\fX}^{(\bullet)} (Z)) \riso0$, then we conclude.
\end{proof}

\begin{coro} \label{coro-coh-log-smoothcase}
We keep notation \ref{coh-log-smoothcase}. For any $\cE ^{(\bullet)}  \in \smash{\underrightarrow{LD}}  ^\mathrm{b} _{\bbQ, \mathrm{qc}} (\overset{^\mathrm{l}}{} \smash{\widehat{\cD}} _{\fX ^{\sharp}/\fS ^{\sharp}} ^{(\bullet)})$, 
and we have the exact triangle of  $\smash{\underrightarrow{LD}}  ^\mathrm{b} _{\bbQ, \mathrm{qc}} (\overset{^\mathrm{l}}{} \smash{\widehat{\cD}} _{\fX ^{\sharp}/\fS ^{\sharp}} ^{(\bullet)})$:
\begin{equation}
u _{+} ^{ (\bullet)} \circ  u ^{ (\bullet)!}  (\cE ^{(\bullet)} ) \to  \cE ^{(\bullet)}  \to  \cB _{\fX}^{(\bullet)} (Z)  \widehat{\otimes} ^{\bbL} _{\cO _{\fX ^{\sharp}}} \cE ^{(\bullet)} \to  u _{+} ^{ (\bullet)} \circ  u ^{ (\bullet)!} (\cE ^{(\bullet)} ) [1].
\end{equation}
\end{coro}

\begin{proof}
Since we have $u _{+} ^{ (\bullet)} \circ  u ^{ (\bullet)!}  (\cE ^{(\bullet)} ) \riso u _{+} ^{ (\bullet)} \circ  u ^{ (\bullet)!}  (\cO _{\fX})  \widehat{\otimes} ^{\bbL} _{\cO _{\fX ^{\sharp}}} 
\cE ^{(\bullet)} $ (use \cite[9.4.3.1]{Car25}),  this is a consequence of  \ref{coh-log-smoothcase}.
\end{proof}

\begin{empt}
[Strictly semi-stable pairs over $S$]
\label{dfn-ssp}
We recall from \cite[6.3]{dejong} that $(X,Z)$ is a {\it strict semi-stable pair} if 
\begin{enumerate}[(a)]
\item $X $ is strictly semi-stable over $S$ (see \ref{dfn-sss}), 
\item Let $Z _f := \cup _{i\in I} Z _i$ be the union of the irreducible components $Z _i$ of $Z$ which dominates $S$.
For each $J \subset I$, the scheme $Z _J:= \cap _{j\in J}Z _j$ is a disjoint union of strict semi-stable varieties over $S$.
\item $Z$ is a divisor with strict normal crossing on $X$ and $Z= Z _f \cup X _s$. 
\end{enumerate}

\end{empt}

\begin{rem} \label{ssp:exact-closed immersion}
Let $(X,Z)$ be a strict semi-stable pair.  For any  irreducible component $Z _i$ of $Z _f$, the canonical closed immersion $Z _i \hookrightarrow X$ induces  the exact closed immersion of log smooth schemes
$Z _i ^\sharp \hookrightarrow X ^\sharp $, where $Z _i ^\sharp$ and $X ^\sharp$ are endowed with their canonical log structure coming from their respective special fibers. Indeed, this is a consequence of the following local description:
Let $x \in X _s$. Let $X _1, \dots, X _n $ be the irreducible components of $X _s$ containing $x$ and $Z _1, \dots, Z _m$ be the irreducible components of $Z _f$ containing $x$.
Then there exists an open neighborhood $U$ of $x$ and a smooth morphism $U \to \Spec R [ t _1, \dots, t _n, s _1, \dots s _m] / (t - t _1 \cdots t _n)$ such that $X _i \cap U$ is given by $t _i = 0$ and 
$Z _j \cap U$ is given by $s _j =0$ (this is a consequence of the local description of  \cite[6.4]{dejong} and of \cite[17.5.3]{EGAIV4}. 
Hence, we are done. 
Moreover, if $u _1,\dots, u _o$ are coordinates of $U \to \Spec R [ t _1, \dots, t _n, s _1, \dots s _m] / (t - t _1 \cdots t _n)$
then $(u _1,\dots, u _o, t _1, \dots, t _n, s _1, \dots s _m$ are overconvergent coordinates of $U/S$. Hence,  $Z _f$ satisfies the condition of the proposition \ref{NCDgencoh}.
\end{rem}

\begin{coro} \label{coh-ssp}
Let $\fX ^{\sharp}$ be a strict semistable log formal schemes over $\fS ^\sharp$. Let $Z$ be a divisor of $X$ such that $(X,Z)$ is a strict semi-stable pair.  
Then  $\cB _{\fX}^{(\bullet)} (Z _f)  \in  \smash{\underrightarrow{LM}}   _{\bbQ, \mathrm{coh}} (\overset{^\mathrm{l}}{} \smash{\widehat{\cD}} _{\fX ^{\sharp}/\fS ^{\sharp} } ^{(\bullet)})$,
and  $\cB _{\fX}^{(\bullet)} (Z)  \in  \smash{\underrightarrow{LM}}  _{\bbQ, \mathrm{coh}} (\overset{^\mathrm{l}}{} \smash{\widehat{\cD}} _{\fX ^{\sharp}/\fS ^{\sharp} } ^{(\bullet)} ( X _s))$.
\end{coro}

\begin{proof}
This is identical to the proof of first step of Theorem \cite[B.2]{caro-stab-sys-ind-surcoh}. Since the second statement is a consequence of the first one, let us chcek the first one.
Since this is local,  using Remark \ref{ssp:exact-closed immersion} we can suppose there exists an affine smooth morphism $X \to \Spec R [ t _1, \dots, t _n, s _1, \dots s _m] / (t - t _1 \cdots t _n)$
such that $X _i $ is given by $t _i = 0$ and  $Z _j $ is given by $s _j =0$, where  $X _1, \dots, X _n $ are  the irreducible components of $X _s$,  
and $Z _1, \dots, Z _m$ be the irreducible components of $Z _f$. Let $\fZ _1 ^\sharp$ be an affine strict semistable log formal schemes over $\fS ^\sharp$ 
which lifts $Z _1 ^\sharp$.  Let $u \colon \fZ _1 ^\sharp \hookrightarrow \fX ^\sharp$ be the corresponding exact closed immersion. To check the coherence of 
$\cB _{\fX}^{(\bullet)} (Z )$, we proceed by induction on $m$. When $m\leq 1$, this is already know (see \ref{coh-log-smoothcase}).
Suppose $m \geq 2$.  Let $Z ' _f := Z _2 \cup \dots \cup Z _m$, $Z ' _{1 f} := Z _1 \cap  Z '_1$, and  $Z'= Z' _f \cup X _s$, $Z' _1= Z' _{1f} \cup Z _s$. 
Then $(Z _1, Z ' _1)$ is a strict semi-stable pair and  we compute  $u ^{ (\bullet)!}  (\cB _{\fX}^{(\bullet)} (Z ') ) \riso \cB _{\fZ _1}^{(\bullet)} (Z '_1 ) [-1]$.
Using the induction hypothesis, we get  $\cB _{\fZ _1}^{(\bullet)} (Z '_1 ) \in \smash{\underrightarrow{LM}}  ^\mathrm{b} _{\bbQ, \mathrm{coh}}
(\overset{^\mathrm{l}}{} \smash{\widehat{\cD}} _{\fZ _1 ^{\sharp}/\fS ^{\sharp} } ^{(\bullet)})$. Hence,  $u ^{ (\bullet)} _+   u ^{ (\bullet)!}  (\cB _{\fX}^{(\bullet)} (Z ') )
\in  \smash{\underrightarrow{LD}}   ^\mathrm{b}  _{\bbQ, \mathrm{coh}} (\overset{^\mathrm{l}}{} \smash{\widehat{\cD}} _{\fX ^{\sharp}/\fS ^{\sharp} } ^{(\bullet)})$.
By using again the induction hypothesis, we get $\cB _{\fX}^{(\bullet)} (Z ')  \in  \smash{\underrightarrow{LM}}  
_{\bbQ, \mathrm{coh}} (\overset{^\mathrm{l}}{} \smash{\widehat{\cD}} _{\fX ^{\sharp}/\fS ^{\sharp} } ^{(\bullet)})$. Hence, from \ref{coro-coh-log-smoothcase}, we get the exact triangle in 
$\smash{\underrightarrow{LD}}   ^\mathrm{b}  _{\bbQ, \mathrm{coh}} (\overset{^\mathrm{l}}{} \smash{\widehat{\cD}} _{\fX ^{\sharp}/\fS ^{\sharp} } ^{(\bullet)})$ of the form: 
\begin{equation}
u _{+} ^{ (\bullet)} \circ  u ^{ (\bullet)!}   (\cB _{\fX}^{(\bullet)} (Z ') ) \to   (\cB _{\fX}^{(\bullet)} (Z ') )
\to  \cB _{\fX}^{(\bullet)} (Z _1 )  \widehat{\otimes} ^{\bbL} _{\cO _{\fX ^{\sharp}}}  (\cB _{\fX}^{(\bullet)} (Z ') ) 
\to  u _{+} ^{ (\bullet)} \circ  u ^{ (\bullet)!} (\cB _{\fX}^{(\bullet)} (Z ') )  [1].
\end{equation}
Since  $\cB _{\fX}^{(\bullet)} (Z )   \riso  \cB _{\fX}^{(\bullet)} (Z _1 )  \widehat{\otimes} ^{\bbL} _{\cO _{\fX ^{\sharp}}}  (\cB _{\fX}^{(\bullet)} (Z ') )$,   we get the coherence of  $\cB _{\fX}^{(\bullet)} (Z )  $.
\end{proof}

\begin{rem}
The above corollary is an analogue of \ref{coh-smoothsubsch}.
We will extend the second part of the above corollary in the general case later in \ref{coh-ss-div} when $Z$ is any divisor.
\end{rem}

\subsubsection{The general case}

\begin{ntn}
\label{lem-RGamma->id}
Let $\fX ^{\sharp}$ be a strictly semistable log formal scheme over $\fS ^\sharp$.
Let $u \colon Z ^\sharp \to X ^{\sharp}$ be a weakly closed immersion of strictly semistable log schemes over $S ^\sharp$.
Following \ref{coh-smoothsubsch}, 
$\R \underline{\Gamma} ^\dag _Z \widetilde{\cO} _{\fX,\bbQ}  
\in 
D ^{\mathrm{b}} _{\mathrm{coh}}
(\widetilde{\cD} ^\dag _{\fX ^{\sharp}/\fS ^{\sharp}})$.
Hence, using the equivalence of categories 
\cite[8.4.1.15]{Car25}, 
there exists an object of $\smash{\underrightarrow{LD}} ^{\mathrm{b}} _{\bbQ, \mathrm{coh}} ( \smash{\widehat{\cD}} _{\fX ^{\sharp}} ^{(\bullet)}(X _s))$, 
denoted by 
$\R \underline{\Gamma} ^{\dag} _Z \cB ^{(\bullet)} _{\fX} (X _s) $ such that 
$\underrvec{l}_{\bbQ} ^{*} 
\R \underline{\Gamma} ^{\dag} _Z \cB ^{(\bullet)} _{\fX} (X _s)
\riso 
\R \underline{\Gamma} ^\dag _Z \widetilde{\cO} _{\fX,\bbQ}  $.
\end{ntn}

\begin{empt}
Let $u \colon \fZ ^\sharp \to \fX ^{\sharp}$ be a weakly closed immersion of strictly semistable log formal schemes over $\fS ^\sharp$.
Since 
$u ^{ (\bullet)!} ( \cB ^{(\bullet)} _{\fX} (X _s))
\riso 
\cB ^{(\bullet)} _{\fZ} (Z _s) [d _X - d _P]
\in 
\smash{\underrightarrow{LD}} ^{\mathrm{b}} _{\bbQ, \mathrm{coh}} ( \smash{\widehat{\cD}} _{\fZ ^{\sharp}} ^{(\bullet)}(Z _s))$, 
since $u$ is proper then using \cite[9.4.2.4]{Car25} we get 
$u _{+} ^{ (\bullet)} \circ  u ^{ (\bullet)!} (  \cB ^{(\bullet)} _{\fX} (X _s) )
\in 
\smash{\underrightarrow{LD}} ^{\mathrm{b}} _{\bbQ, \mathrm{coh}} ( \smash{\widehat{\cD}} _{\fX ^{\sharp}} ^{(\bullet)}(X _s))$
and 
$\underrvec{l}_{\bbQ} ^{*}  (u _{+} ^{ (\bullet)} \circ  u ^{ (\bullet)!}  \cB ^{(\bullet)} _{\fX} (X _s))
\riso 
u _+ ( \underrvec{l}_{\bbQ} ^{*}  (u ^{ (\bullet)!}  \cB ^{(\bullet)} _{\fX} (X _s)) )
\riso 
u _+ u ^! (\widetilde{\cO} _{\fX,\bbQ})$.
Since 
$\R \underline{\Gamma} ^{\dag} _Z \cB ^{(\bullet)} _{\fX} (X _s)
 \in 
 \smash{\underrightarrow{LD}} ^{\mathrm{b}} _{\bbQ, \mathrm{coh}} ( \smash{\widehat{\cD}} _{\fX ^{\sharp}} ^{(\bullet)}(X _s))$
and
$\underrvec{l}_{\bbQ} ^{*} 
\R \underline{\Gamma} ^{\dag} _Z \cB ^{(\bullet)} _{\fX} (X _s)
\riso 
\R \underline{\Gamma} ^\dag _Z \widetilde{\cO} _{\fX,\bbQ}  $, then 
using Corollary \ref{coro-trace-upre} and the full faithfulness of the functor
$\underrvec{l}_{\bbQ} ^{*} $ of \cite[8.4.1.15]{Car25}, 
we get the isomorphism in 
$\smash{\underrightarrow{LD}} ^{\mathrm{b}} _{\bbQ, \mathrm{coh}} ( \smash{\widehat{\cD}} _{\fX ^{\sharp}} ^{(\bullet)}(X _s))$ :
\begin{equation}
\label{coro-trace-upre-limcat}
\R \underline{\Gamma} ^{\dag} _Z \cB ^{(\bullet)} _{\fX} (X _s)
\riso 
u _{+} ^{ (\bullet)} \circ  u ^{ (\bullet)!} ( \cB ^{(\bullet)} _{\fX} (X _s))
\riso 
u _{+} ^{ (\bullet)} (\cB ^{(\bullet)} _{\fZ} (Z _s)) [d _{X/P}].
\end{equation}

\end{empt}

The following theorem extends \ref{coh-ssp}. 
\begin{thm} \label{coh-ss-div}
Let $\fX ^{\sharp}$ be a strict semistable log formal schemes over $\fS ^\sharp$. Let $Z$ be a divisor of $X$ containing $X _s$.
Then  $\cB _{\fX}^{(\bullet)} (Z)  \in  \smash{\underrightarrow{LM}}   _{\bbQ, \mathrm{coh}} (\overset{^\mathrm{l}}{} \smash{\widehat{\cD}} _{\fX ^{\sharp}/\fS ^{\sharp} } ^{(\bullet)} (s))$.
\end{thm}

\begin{proof}
Following Theorem \cite[6.5]{dejong}, 
there exist a trait $S _1$ finite over $S$, an $S _1$-variety $X _1$, an alteration of schemes over $S$ (in the sense of \cite[2.20]{dejong})
$g \colon X _1 \to X$ 
and an open immersion $j _1 \colon X _1 \to \overline{X} _1$ of $S _1$-varieties, with the following properties: 
\begin{enumerate}
[(i)]
\item $\overline{X} _1$ is a projective $S _1$-variety with geometrically irreducible generic fibre, and 
\item  the pair $(\overline{X} _1, g  ^{-1} (Z ) _\mathrm{red} \cup (\overline{X} _1 \setminus j _1 (X _1))$ is strict semi-stable.
\end{enumerate}
In particular, we get that 
$(X _1, g  ^{-1} (Z ) _\mathrm{red})$ is a strict semi-stable pair
and there exists a closed immersion of the form
$u\colon X _1 \hookrightarrow \P ^n _{X}$
whose composition with the projection 
$\P ^n _{X} \to X$ is $g$.

0) Following \ref{ll}, the morphism $S _1 \to S$ has a lifting $\fS _1 ^\sharp \to \fS ^\sharp$. 
Using \ref{desc-coh-chgbase}, we can reduce to the case where $S _1 =S$.

1) i) Set  $\fP := \widehat{\P} ^n _{\fX}$,  $f \colon \fP \to \fX$ be the projection. 
From  \ref{adj-morph},  we have the adjoint morphism
$ f _{+} ^{(\bullet)} \circ f ^{!} (\widetilde{\cO} _{\fX,\bbQ}) \to  \widetilde{\cO} _{\fX,\bbQ}$
in  $D ^{\mathrm{b}} _{\mathrm{coh}}( \smash{\cD} ^\dag _{\fX ^{\sharp}} (\hdag X _s) _{\bbQ} )$.
Following \ref{coh-smoothsubsch}, we get in $D ^{\mathrm{b}} _{\mathrm{coh}}( \smash{\cD} ^\dag _{\fP ^\sharp/\fS ^{\sharp}} (\hdag P _s) _{\bbQ} )$
the morphism $\R \underline{\Gamma} ^{\dag} _{X _1}  (\widetilde{\cO} _{\fP,\bbQ}) \to \widetilde{\cO} _{\fP,\bbQ}$.
Since  $f ^{!} (\widetilde{\cO} _{\fX,\bbQ}) \riso \widetilde{\cO} _{\fP,\bbQ} [n]$, then we get 
the morphism in  $D ^{\mathrm{b}} _{\mathrm{coh}}( \smash{\cD} ^\dag _{\fX ^{\sharp}} (\hdag X _s) _{\bbQ} )$
\begin{equation} \label{cstrcf+GammaO->O}
f _{+} ( \R \underline{\Gamma} ^{\dag} _{X _1} \widetilde{\cO} _{\fP,\bbQ} [n] )  \to  \widetilde{\cO} _{\fX,\bbQ}.
\end{equation}

ii) In this step, we construct by duality from the step i) 
the morphism
$\widetilde{\cO} _{\fX,\bbQ}
\to 
f _{+} ( \R \underline{\Gamma} ^{\dag} _{X _1} \widetilde{\cO} _{\fP,\bbQ} [n])$.
 We have 
\begin{equation}
\label{RGammadual}
 \bbD  (  \R \underline{\Gamma} ^{\dag} _{X _1}  \widetilde{\cO} _{\fP,\bbQ} [n])
\underset{\ref{coro-sp+jdagO}}{\riso}
\bbD  (  \sp _+ (j ^\dag  \smash{\cO} _{]X _1[ _{\fP}}))
\underset{\ref{propspetdualsansfrob}}{\riso}
\sp _+ ((j ^\dag  \smash{\cO} _{]X _1[ _{\fP}}) ^\vee)
\riso 
\sp _+ (j ^\dag  \smash{\cO} _{]X _1[ _{\fP}})
\underset{\ref{coro-sp+jdagO}}{\riso}
 \R \underline{\Gamma} ^{\dag} _{X _1}  \widetilde{\cO} _{\fP,\bbQ} [n].
\end{equation}
This yields
$$\widetilde{\cO} _{\fX,\bbQ} 
\underset{\ref{dualisoscvdag}}{\riso}
\bbD (\widetilde{\cO} _{\fX,\bbQ})
\underset{\ref{cstrcf+GammaO->O}}{\longrightarrow} 
\bbD f _{+}  ( \R \underline{\Gamma} ^{\dag} _{X _1} \widetilde{\cO} _{\fP,\bbQ} [n])
\underset{\ref{dualrelative}}{\riso} 
f _{+} \bbD  (  \R \underline{\Gamma} ^{\dag} _{X _1}  \widetilde{\cO} _{\fP,\bbQ} [n])
\underset{\ref{RGammadual}}{\riso} 
f _{+} ( \R \underline{\Gamma} ^{\dag} _{X _1} \widetilde{\cO} _{\fP,\bbQ} [n]).$$

iii) We check  that the composite morphism
$\widetilde{\cO} _{\fX,\bbQ}
\to 
f _{+} ( \R \underline{\Gamma} ^{\dag} _{X _1} \widetilde{\cO} _{\fP,\bbQ} [n])
 \to 
\widetilde{\cO} _{\fX,\bbQ}$
in 
$D ^{\mathrm{b}} _{\mathrm{coh}}( \smash{\cD} ^\dag _{\fX ^{\sharp}} (\hdag X _s) _{\bbQ} )$
is an isomorphism.
Indeed, using the remark of \ref{empt-projff}, since this a morphism of 
the abelian category $\mathrm{MIC} ^{\dag\dag} ((X _\eta, X,\fX ^\sharp)/\cE ^{\dag}  _K) $,
we reduce to check its restriction to an open dense subset is an isomorphism.
Hence, we can suppose that $X =X _\eta$, and that 
$X _1 \to X$ is the composition of a finite surjective radicial morphism
$X _1 \to X  _2$ with a finite and etale morphism 
$X _2 \to X  $. We get some lifting $g _1\colon \fX _1 \to \fX _2$ and $g _2\colon \fX _2 \to \fX$. 
Set $g \coloneqq g _2 \circ g _1$. We reduce to prove that  $\cO _{\fX,\bbQ} \to g _+ \cO _{\fX_1,\bbQ} \to \cO _{\fX,\bbQ}$ is an isomorphism. 
Using \ref{univhomeo-eqcat-isoc}, we reduce to the case where 
$g  \colon X _1 \to X$ is finite and étale, which is well known. 

2) From \ref{lem-RGamma->id}, we get by definition the object  $ \R \underline{\Gamma} ^{\dag} _{X _1} \cB ^{(\bullet)} _{\fP} (P _s)$ of  
$\smash{\underrightarrow{LD}} ^{\mathrm{b}} _{\bbQ, \mathrm{coh}} ( \smash{\widehat{\cD}} _{\fP ^{\sharp}/\fS ^{\sharp}} ^{(\bullet)}(P _s))$
corresponding to  $\R \underline{\Gamma} ^{\dag} _{X _1} \widetilde{\cO} _{\fP,\bbQ}$  via the equivalence of categories of \cite[8.4.1.15]{Car25}. Since $f$ is proper, 
we obtain $f _{+} ^{ (\bullet)} ( \R \underline{\Gamma} ^{\dag} _{X _1} \cB ^{(\bullet)} _{\fP} (P _s) [n])
\in \smash{\underrightarrow{LD}} ^{\mathrm{b}} _{\bbQ, \mathrm{coh}} ( \smash{\widehat{\cD}} _{\fX ^{\sharp}} ^{(\bullet)}(X _s))$
and the corresponding object of $D ^{\mathrm{b}} _{\mathrm{coh}}( \smash{\cD} ^\dag _{\fX ^{\sharp}} (\hdag X _s) _{\bbQ} )$
via the equivalence of categories of \cite[8.4.1.15]{Car25} is  $f _{+} ( \R \underline{\Gamma} ^{\dag} _{X _1} \widetilde{\cO} _{\fP,\bbQ} [n])$.
Hence, we get from 1) the morphisms 
\begin{equation}\label{surcoh2.1.4-cor-isoappl1pre} \cB _{\fX}^{(\bullet)} (X _s)
\to f _{+} ^{ (\bullet)} (\R \underline{\Gamma} ^{\dag} _{X _1} \cB ^{(\bullet)} _{\fP} (P _s) [n]) \to  \cB _{\fX}^{(\bullet)} (X _s)
\end{equation}
in $\smash{\underrightarrow{LD}} ^{\mathrm{b}} _{\bbQ, \mathrm{coh}} ( \smash{\widehat{\cD}} _{\fP ^{\sharp}/\fS ^{\sharp}} ^{(\bullet)}(P _s))$ whose composite is an isomorphism.

3) Using \cite[9.1.3.2]{Car25}, by applying the functor $-\smash{\widehat{\otimes}}^\L _{\cO ^{(\bullet)}  _{\fX} }   \cB ^{ (\bullet)} _{\fX} (Z ) $ to \ref{surcoh2.1.4-cor-isoappl1pre} , 
we obtain that $\cB _{\fX}^{(\bullet)} (Z) $ is a direct factor of the left term of the isomorphism:
\begin{gather}
\label{surcoh2.1.4-cor-isoappl1}
f _{+} ^{ (\bullet)} \left (\R \underline{\Gamma} ^{\dag} _{X _1} \cB ^{(\bullet)} _{\fP} (P _s) [n] \right)
\smash{\widehat{\otimes}}^\L _{\cO ^{(\bullet)}  _{\fX} }   \cB ^{ (\bullet)} _{\fX} (Z ) 
\underset{\cite[9.4.3.3]{Car25}}{\riso} 
f _{+} ^{ (\bullet)} \left (\R \underline{\Gamma} ^{\dag} _{X _1} \cB ^{(\bullet)} _{\fP} (P _s) [n]
\smash{\widehat{\otimes}}^\L _{\cO ^{(\bullet)}  _{\fP} }  
\cB ^{ (\bullet)} _{\fP} (f ^{-1} (Z)) \right).
\end{gather}
Since $f$ is proper, it is sufficient to
check that 
$ \R \underline{\Gamma} ^{\dag} _{X _1} \cB ^{(\bullet)} _{\fP} (P _s) [n]
\smash{\widehat{\otimes}}^\L _{\cO ^{(\bullet)}  _{\fP} }  
\cB ^{ (\bullet)} _{\fP} (f ^{-1} (Z))$
is an object of 
$\smash{\underrightarrow{LM}}  _{\bbQ, \mathrm{coh}} ( \smash{\widehat{\cD}} _{\fP ^{\sharp}/\fS ^{\sharp}} ^{(\bullet)}(P _s))$.
Since this is local in $\fP ^{\sharp}$, we can suppose $\fP ^{\sharp}$ affine. 
Hence, there exists  a morphism 
$u \colon \fX _1 ^\sharp \to \fP ^{\sharp} $
of strictly semistable log formal schemes over $\fS ^\sharp$
 which is $X _1 ^\sharp \to P ^{\sharp}$ modulo $\pi$.
 Hence, using  \cite[9.4.3.3]{Car25}, we get 
\begin{gather}
\notag
\R \underline{\Gamma} ^{\dag} _{X _1} \cB ^{(\bullet)} _{\fP} (P _s) [n]
\smash{\widehat{\otimes}}^\L _{\cO ^{(\bullet)}  _{\fP} }  
\cB ^{ (\bullet)} _{\fP} (f ^{-1} (Z))
\underset{\ref{coro-trace-upre-limcat}}{\riso} 
 u _{+} ^{ (\bullet)} ( \cB _{\fX _1} ^{(\bullet)} (X _{1s} ) )
\smash{\widehat{\otimes}}^\L _{\cO ^{(\bullet)}  _{\fP} }    \cB ^{ (\bullet)} _{\fP} (f ^{-1} (Z))
\underset{\cite[9.4.3.1]{Car25}}{\riso} 
\\
\notag 
 u _{+} ^{ (\bullet)} \left (  \cB _{\fX _1} ^{(\bullet)} (X _{1s} ) 
\smash{\widehat{\otimes}}^\L _{\cO ^{(\bullet)}  _{\fX _1} }  
u ^{ (\bullet)!} (\cB ^{ (\bullet)} _{\fP} (f ^{-1} (Z))[n])
\right)
\riso 
 u _{+} ^{ (\bullet)} \left ( \cB _{\fX _1}^{(\bullet)} (X _{1s} )  
\smash{\widehat{\otimes}}^\L _{\cO ^{(\bullet)}  _{\fX _1} } 
 \cB ^{ (\bullet)} _{\fX _1} (g  ^{-1} (Z))
\right)
\riso 
 u _{+} ^{ (\bullet)} ( \cB ^{ (\bullet)} _{\fX _1} (g ^{-1} (Z))).
\end{gather}
Since $(X _1, g  ^{-1} (Z ) _\mathrm{red})$ is a strict semi-stable pair, then following  \ref{coh-ssp} we obtain
$ \cB ^{ (\bullet)} _{\fX _1} (g  ^{-1} (Z))$ is an object of
$\smash{\underrightarrow{LM}}  
_{\bbQ, \mathrm{coh}}
(\overset{^\mathrm{l}}{} \smash{\widehat{\cD}} _{\fX _1 ^{\sharp}/\fS ^{\sharp} } ^{(\bullet)}( X _{1s}))$.
 Hence, 
$  u _{+} ^{ (\bullet)} ( \cB ^{ (\bullet)} _{\fX _1} (g  ^{-1} (Z)))\in 
\smash{\underrightarrow{LM}}  _{\bbQ, \mathrm{coh}} ( \smash{\widehat{\cD}} _{\fP ^{\sharp}/\fS ^{\sharp}} ^{(\bullet)}(P _s))$.
\end{proof}

\section{Local cohomological operations and applications} 

\subsection{Local cohomological functors}
\label{LCF}

\subsubsection{Local cohomological functor with strict support over a divisor}
Let  $\fP $ be a separated, quasi-compact smooth formal schemes over $\fS $. 
\begin{empt}
\label{GammaT}
Let  $T \subset T'$ be a two divisors of $P$.
Suppose we have the commutative diagram in 
$\smash{\underrightarrow{LD}} ^\mathrm{b} _{\bbQ,\mathrm{qc}} ( \widehat{\cD} _{\fP /\fS  } ^{(\bullet)} (P _s))$
of the form
\begin{equation}
\label{prefonct-GammaT}
\xymatrix @=0,4cm{
{\cF ^{(\bullet)}  } 
\ar[r] ^-{}
& 
{\cE ^{(\bullet)}  } 
\ar[r] ^-{}
\ar[d] ^-{\phi}
& 
{(\hdag T) (\cE ^{(\bullet)} )} 
\ar[d] ^-{(\hdag T)(\phi)}
\ar[r] ^-{}
&
{\cF ^{(\bullet)}  [1]} 
\\ 
{\cF ^{\prime (\bullet)}  } 
\ar[r] ^-{}
& 
{\cE ^{\prime (\bullet)}  } 
\ar[r] ^-{}
& 
{(\hdag T) (\cE ^{\prime (\bullet)} )} 
\ar[r] ^-{}
&
{\cF ^{\prime (\bullet)}   [1]} 
}
\end{equation}
where horizontal morphisms are the canonical ones and
where both horizontal triangles are distinguished. 
In the same way as \cite[13.1.1.4]{Car25}, we prove there exists a unique morphism $\cF ^{(\bullet)}\to \cF ^{\prime (\bullet)} $ making commutative in 
$\smash{\underrightarrow{LD}} ^\mathrm{b} _{\bbQ,\mathrm{qc}} ( \widehat{\cD} _{\fP /\fS  } ^{(\bullet)} (P _s) )$ the diagram:
\begin{equation}
\label{fonct-hdagT}
\xymatrix @=0,4cm{
{\cF ^{(\bullet)}  } 
\ar[r] ^-{}
\ar@{.>}[d] ^-{\exists !}
& 
{\cE ^{(\bullet)}  } 
\ar[r] ^-{}
\ar[d] ^-{\phi}
& 
{(\hdag T) (\cE ^{(\bullet)} )} 
\ar[d] ^-{(\hdag T)(\phi)}
\ar[r] ^-{}
& 
{\cF ^{(\bullet)}  [1]} 
\ar@{.>}[d] ^-{\exists !}
\\ 
{\cF ^{\prime (\bullet)}  } 
\ar[r] ^-{}
& 
{\cE ^{\prime (\bullet)}  } 
\ar[r] ^-{}
& 
{(\hdag T) (\cE ^{\prime (\bullet)} )} 
\ar[r] ^-{}
& 
{\cF ^{\prime (\bullet)} [1].} 
}
\end{equation}
In the same way as \cite[1.1.10]{BBD}, 
this implies that the cone of 
$\cE ^{(\bullet)}
\to  
(\hdag T) (\cE ^{(\bullet)} ) $
is unique up to canonical isomorphism. 
Hence, such a complex $\cF ^{ (\bullet)}$ is unique up to canonical isomorphism.
We denote it by 
$ \R \underline{\Gamma} ^\dag _{T} (\cE ^{(\bullet)})$.
Moreover, 
the complex 
$\R \underline{\Gamma} ^\dag _{T} (\cE ^{(\bullet)}) $ is functorial in
$\cE ^{(\bullet)}$.
\end{empt}

\subsubsection{Local cohomological functor with strict support over a closed subvariety}

\begin{dfn}
\label{dfn-4.3.4}
Let $X$ be a (reduced) closed subscheme of $P$.
As well as for \cite[13.1.3.2]{Car25} or \cite[2.2]{caro_surcoherent}, we define
the local cohomological functor  
$\R \underline{\Gamma} ^\dag _{X}
\colon 
\smash{\underrightarrow{LD}} ^{\mathrm{b}} _{\bbQ,\mathrm{qc}} 
(\overset{^\mathrm{l}}{} \smash{\widehat{\cD}} _{\fP /\fS  } ^{(\bullet)} (P_s))
\to 
\smash{\underrightarrow{LD}} ^{\mathrm{b}} _{\bbQ,\mathrm{qc}} 
(\overset{^\mathrm{l}}{} \widehat{\cD} _{\fP /\fS  } ^{(\bullet)} (P _s))$
with strict support in $X$ as follows.
Beware that in our new context, since we need the coherence property \ref{coh-ss-div}, we have to add singularities along $P _s$. 
Recall  
$P$ is integral by hypothesis.
\begin{enumerate}[(a)]
\item When $X= P$, the functor $\R \underline{\Gamma} ^\dag _{X}$ 
is by definition the identity.

\item Suppose now $X \not = P$. 
Following \cite[13.1.3.6]{Car25},
$X$ is the finite intersection of divisors of $P$.
Choose some divisors
$T _1, \dots, T _{r}$ of $P$ such that
$X = \cap _{i =1} ^{r} T _i$.
For $\cE ^{(\bullet)}
\in \smash{\underrightarrow{LD}} ^\mathrm{b} _{\bbQ,\mathrm{qc}} ( \smash{\widetilde{\cD}} _{\fP /\fS } ^{(\bullet)} (P_s))$,
in the same way as \cite[13.1.3.8]{Car25}, the complex
$$\R \underline{\Gamma} ^\dag _{X} (\cE ^{(\bullet)} ):=
\R \underline{\Gamma} ^\dag _{T _r} \circ \cdots \circ 
\R \underline{\Gamma} ^\dag _{T _1} (\cE ^{(\bullet)} )$$
does not depend canonically on the choice of the divisors  $T _1,\dots, T _r$
satisfying $X = \cap _{i =1} ^{r} T _i$.
\end{enumerate}

\end{dfn}

\subsubsection{Localisation outside a closed subscheme functor}

\begin{dfn}
Let $\E ^{(\bullet)} \in \smash{\underrightarrow{LD}} ^\mathrm{b} _{\Q,\mathrm{qc}} ( \smash{\widehat{\D}} _{\fP /\fS } ^{(\bullet)} (P _s))$.
Let $X$ be a closed subscheme of  $P$.
Using \cite[1.1.10]{BBD}, 
we check that the cone of the morphism
$\R \underline{\Gamma} ^\dag _{X} (\cE ^{(\bullet)} )
\to  
\cE ^{(\bullet)} $
is unique up to canonical isomorphism (for more details, 
see \cite[4.4.3]{caro-stab-sys-ind-surcoh}).
We will denote it by 
$(\hdag X) (\cE ^{(\bullet)} )$.
We check that 
$(\hdag X) (\cE ^{(\bullet)} )$
is functorial in $X$, and $\cE ^{(\bullet)} $.
We get canonical the distinguished triangle
$$\R \underline{\Gamma} ^\dag _{X} (\cE ^{(\bullet)} )
\to 
\cE ^{(\bullet)} 
\to 
(\hdag X) (\cE ^{(\bullet)} )
 \to 
\R \underline{\Gamma} ^\dag _{X} (\cE ^{(\bullet)} )[1] .$$
\end{dfn}

\begin{thm}\label{coh-OGammaY}
Let  $X, X'$ be two closed subschemes of $P$.
We have
$(\hdag X) 
\R \underline{\Gamma} ^\dag _{X '} (\cB _{\fP}^{(\bullet)} (X _s) )
\in 
\smash{\underrightarrow{LD}} ^{\mathrm{b}}  
_{\bbQ, \mathrm{coh}}
(\overset{^\mathrm{l}}{} \widehat{\cD} _{\fP /\fS  } ^{(\bullet)} (P _s))$.
\end{thm}

\begin{proof}
This is similar to \cite[13.1.4.6]{Car25}.
\end{proof}

\begin{empt}
In the same way as \cite[2.2.16]{caro_surcoherent} or \cite[13.1.4.5]{Car25},
we get the Mayer-Vietoris distinguished triangles :
\begin{gather}\label{eq1mayer-vietoris}
  \R \underline{\Gamma} ^\dag _{X \cap X'}(\cE ^{(\bullet)}) \rightarrow
  \R \underline{\Gamma} ^\dag _{X }(\cE ^{(\bullet)}) \oplus
\R \underline{\Gamma} ^\dag _{X ' }(\cE ^{(\bullet)})  \rightarrow
\R \underline{\Gamma} ^\dag _{X \cup X '}(\E^{(\bullet)} ) \rightarrow
\R \underline{\Gamma} ^\dag _{X \cap X'}(\E^{(\bullet)} )[1],\\
 (\hdag X \cap X')(\cE ) \rightarrow  (\hdag X )(\E^{(\bullet)} ) \oplus
 (\hdag X ') (\E^{(\bullet)} )  \rightarrow   (\hdag X \cup X ')(\E^{(\bullet)} ) \rightarrow (\hdag X \cap X')(\cE ^{(\bullet)})[1].
\end{gather}

\end{empt}

 \begin{prop} \label{prop2.2.9}
Let  $D$ be a divisor of  $P$ containing $P _s$,  $X$ be a closed subscheme of $P$, $\fU$ be the open subset of  $\fP$ complementary to the support of $X$.
Let $\E ^{(\bullet)} \in  \smash{\underrightarrow{LD}} ^{\mathrm{b}} _{\Q,\mathrm{coh}} ( \smash{\widetilde{\D}} _{\fP /\fS } ^{(\bullet)} (D))$.
The following assertions are equivalent :
\begin{enumerate}[(a)]
\item We have in  $\smash{\underrightarrow{LD}} ^{\mathrm{b}} _{\bbQ,\mathrm{coh}} ( \smash{\widetilde{\cD}} _{\fU /\fS } ^{(\bullet)} (D\cap U))$
the isomorphism $\cE ^{(\bullet)}|\fU \riso 0$.
\item The canonical morphism 
$\R \underline{\Gamma} ^\dag _{X} (\cE ^{(\bullet)})
\to 
\cE ^{(\bullet)}$ is an isomorphism in 
$\smash{\underrightarrow{LD}} ^{\mathrm{b}} _{\bbQ} ( \smash{\widetilde{\cD}} _{\fP /\fS } ^{(\bullet)} (D))$.
\item We have in $\smash{\underrightarrow{LD}} ^{\mathrm{b}} _{\bbQ,\mathrm{coh}} ( \smash{\widetilde{\cD}} _{\fP /\fS } ^{(\bullet)} (D))$
the isomorphism $(\hdag X) (\cE ^{(\bullet)} )\riso 0$.
\end{enumerate}

\end{prop}

\begin{proof}This is similar to \cite[13.1.4.8]{Car25}.
\end{proof}

\begin{empt} [Support]  \label{dfn-support} Let  $D$ be a divisor of  $P$ containing $P _s$,
$\cE ^{(\bullet)} \in \underrightarrow{LD} ^{\rmb} _{\bbQ,\coh} ( \widetilde{\D} _{\fP ^\sharp/\fS} ^{(\bullet)} (D))$.
The support of $\cE ^{(\bullet)} $ is by definition the biggest closed subscheme $X$ of $P$ such that
$(\hdag X) (\cE ^{(\bullet)} )\riso 0$ (one of the equivalent conditions of \ref{prop2.2.9}).  We denote it by $\mathrm{Supp}~\cE ^{(\bullet)}$.

Remark if $\cE ^{(\bullet)} \in \smash{\underrightarrow{LM}} _{\bbQ,\coh} ( \widetilde{\D} _{\fP ^\sharp/\fS} ^{(\bullet)} (D))$,
then this is equal to the  support (for the usual definition) of the coherent $\cD  ^\dag _{\fP ^{\sharp}/\fS} (\hdag D) _{\bbQ} $-module
$\underrvec{l}_{\bbQ} ^{*}\, \cE ^{(\bullet)} $, which justifies the terminology. 
\end{empt}

\subsubsection{Local cohomological functor with strict support over a subvariety}

\begin{empt} \label{3.2.1caro-2006-surcoh-surcv}
Let  $X$, $X'$, $T$, $T'$ be closed subschemes of  $P$ such that   $X \setminus T = X' \setminus T'$.
For any $\cE  ^{(\bullet)} \in \underrightarrow{LD} ^{\rmb} _{\bbQ ,\qc} (\widehat{\cD} _{\fP /\fS  } ^{(\bullet)} (P _s))$,
in the same way as \cite[13.1.5.1]{Car25} we have the canonical isomorphism:
\begin{equation}   \label{xtx't'}
\bbR \underline{\Gamma} ^\dag _{X} (\hdag T ) (\cE  ^{(\bullet)}) \riso \bbR \underline{\Gamma} ^\dag _{X'} (\hdag T ') (\cE  ^{(\bullet)}).
\end{equation}
Setting $ Y:= X \setminus T$, we denote by $\bbR \underline{\Gamma} ^\dag _{Y} (\cE  ^{(\bullet)}) $ one of both complexes of \ref{xtx't'}.
\end{empt}

\begin{ntn}\label{ntn-HdagnX}
Let $Y$ be a subvariety of $P$ and $T$ be a divisor of $P$ containing $P _s$.
Using notation \cite[9.1.6.6]{Car25}, we get the  functor  
$$\bbR \underline{\Gamma} ^{\dag} _{Y}:=
\mathrm{Coh} _{T}  (\bbR \underline{\Gamma} ^{\dag} _{Y} ) 
\colon     
D ^{\rmb} _{\coh} (\cD  ^\dag _{\fP /\fS} (\hdag T) _{\bbQ})
\to 
D ^{\rmb}  (\cD  ^\dag _{\fP /\fS} (\hdag T) _{\bbQ}).$$
 If we want to distinguish both  functors  
$\bbR \underline{\Gamma} ^{\dag} _{Y}$,  then we will denote  in this case by 
$\bbR \underline{\Gamma} ^{(\bullet)} _{Y}
\colon     
\underrightarrow{LD} ^{\rmb} _{\bbQ, \qc}
(\widehat{\cD} _{\fP /\fS  } ^{(\bullet)} (P _s))
\to
\underrightarrow{LD} ^{\rmb} _{\bbQ, \qc}
(\widehat{\cD} _{\fP /\fS  } ^{(\bullet)} (P _s))
$.
 Finally,   for any $n\in \bbN$,  we denote by  
$$\mathcal{H} ^{\dag n} _{Y}:= H ^{n} \circ \bbR \underline{\Gamma} ^{\dag} _{Y}
\colon     D ^{\rmb} _{\coh} (\cD  ^\dag _{\fP/\fS} (\hdag T) _{\bbQ})
\to 
M (\cD  ^\dag _{\fP/\fS} (\hdag T) _{\bbQ}),$$
 where 
$M (\cD  ^\dag _{\fP/\fS} (\hdag T) _{\bbQ})$  is   the category of 
left  $\cD  ^\dag _{\fP/\fS} (\hdag T) _{\bbQ} $-modules.
\end{ntn}

\begin{prop}
\label{prop-induction-div-coh}
Let $Y$ be a subvariety of $P$. We have
$\bbR \underline{\Gamma} ^\dag _{Y} ((\cB _{\fP}^{(\bullet)} (X _s) ))
\in
\underrightarrow{LD} ^{\rmb} _{\bbQ,\coh} ( \widehat{\cD} _{\fP /\fS  } ^{(\bullet)} (P _s) )$.
\end{prop}
\begin{proof}
This is a translation of \ref{coh-OGammaY}.
\end{proof}

\begin{empt}
Let  $Y $ and $Y'$ be two subschemes of  $P$.
Let $\cE  ^{(\bullet)} ,\cF ^{(\bullet)}\in \underrightarrow{LD} ^{\rmb} _{\bbQ ,\qc} (\widehat{\cD} _{\fP /\fS  } ^{(\bullet)} (P _s))$.
In the same way as \cite[13.1.5.6]{Car25}, we get the following properties.
\begin{enumerate}[(a)]
\item 
We have the canonical isomorphism functorial in $\cE ^{(\bullet)},~ Y$, and $Y'$:
\begin{equation}
  \label{gammayY'}
  \bbR \underline{\Gamma} ^\dag _{Y} \circ \bbR \underline{\Gamma} ^\dag _{Y'} (\cE  ^{(\bullet)})
  \riso
  \bbR \underline{\Gamma} ^\dag _{Y \cap Y'} (\cE  ^{(\bullet)}).
\end{equation}

\item  
We have the canonical isomorphism functorial in $\cE ^{(\bullet)},~\cF ^{(\bullet)},~ Y$, and $Y'$:
\begin{equation}
\label{fonctYY'Gamma-iso}
\bbR \underline{\Gamma} ^\dag _{Y \cap Y'} (\cE ^{(\bullet)}
\widehat{\otimes} ^{\bbL}  _{\cO^{(\bullet)}  _{\fP} }  \cF ^{(\bullet)} )
\riso
\bbR \underline{\Gamma} ^\dag _{Y}
(\cE ^{(\bullet)})
\widehat{\otimes} ^{\bbL}  _{\cO^{(\bullet)}  _{\fP} }
\bbR \underline{\Gamma} ^\dag _{Y'}
(\cF ^{(\bullet)}).
\end{equation}
\item  If $Y '$ is an open (resp. a closed) subscheme of  $Y$, then we have the canonical homomorphism
$\bbR \underline{\Gamma} ^\dag _{Y} (\cE  ^{(\bullet)}) \to \bbR \underline{\Gamma} ^\dag _{Y'} (\cE  ^{(\bullet)})$
(resp. $\bbR \underline{\Gamma} ^\dag _{Y'} (\cE  ^{(\bullet)}) \to \bbR \underline{\Gamma} ^\dag _{Y} (\cE  ^{(\bullet)})$).
If $Y '$ is a closed subscheme of  $Y$, then we have the localization distinguished triangle
\begin{equation}
\label{tri-local}
\bbR \underline{\Gamma} ^\dag _{Y'} (\cE  ^{(\bullet)}) \to \bbR \underline{\Gamma} ^\dag _{Y} (\cE  ^{(\bullet)})
\to \bbR \underline{\Gamma} ^\dag _{Y \setminus Y'} (\cE  ^{(\bullet)}) \to +1.
\end{equation}

\end{enumerate}

\end{empt}

\subsection{Fundamental properties}

\subsubsection{Commutation with local functors}
\begin{empt} [Base change and its commutations with cohomological operations] \label{comm-chg-base}
Let $\alpha\colon  \fT '\to \fT$ be a morphism of $\mathrm{DVR}  ( \fS)$  (see notation \ref{DVRorientedpre}), $\fX$ be a separated, quasi-compact, smooth formal schemes over $\fT$, 
$\cE ^{(\bullet)} \in \smash{\underrightarrow{LD}} ^{\mathrm{b}} _{\bbQ,\mathrm{qc}} ( \smash{\widehat{\cD}} _{\fX/\fT } ^{(\bullet)}(X _s))$,
$\fX' := \fX \times _{\fT} \fT'$, and $\varpi  \colon \fX' \to \fX$ 
be the projection.
The base change of $\cE ^{(\bullet)} $ by $\alpha$ is 
the object $\varpi  ^{!} (\cE ^{(\bullet)} )= 
\varpi  ^{! (\bullet)} (\cE ^{(\bullet)})$ of  $\smash{\underrightarrow{LD}} ^{\mathrm{b}} _{\bbQ,\mathrm{qc}} (  \smash{\widehat{\cD}} _{\fX'/\fT' } ^{(\bullet)}(X '_s))$,
where $\varpi  \colon \fX' /\fT'\to \fX/\fT$  (see \cite[9.2.6.1]{Car25}).
It can also simply be denoted by 
$$ \cO _{\fT'}  \smash{\widehat{\otimes}}^\L _{\cO _{\fT}}  \cE^{(\bullet) }.$$
\end{empt}

\begin{prop} \label{com-dual-bc}
With notation \ref{comm-chg-base}, 
we have the canonical isomorphism:
\begin{equation} \label{com-dual-bc-iso}
\cO _{\fT'}  \smash{\widehat{\otimes}}^\L _{\cO _{\fT}}  (\bbD _{\fX/\fT}  (\cE^{(\bullet)})) \riso  \bbD _{\fX'/\fT'} (\cO _{\fT'}  \smash{\widehat{\otimes}}^\L _{\cO _{\fT}} \cE^{(\bullet)}).
\end{equation}
\end{prop}
\begin{proof}
This is a consequence of \cite[9.2.6.11]{Car25}.
\end{proof}

\begin{prop}[Commutation with base change]\label{bc-com-hdag-gen}
With notation \ref{comm-chg-base}, let $Y$ be a subvariety of  $X$. Let $Y':= \varpi ^{-1} (Y)$.
For any $\cE ^{ (\bullet)}  \in  \underrightarrow{LD}  ^{-} _{\Q, \qc} (\overset{^\rml}{} \widehat{\cD} _{\fX /\fS} ^{(\bullet)} (X _s) )$, we have the isomorphism 
\begin{equation}\label{bc-com-hdag-gen-iso}
\bbR \underline{\Gamma} ^\dag _{Y'} \circ  \varpi  ^{ (\bullet)! }  (\cE ^{ (\bullet)} ) \riso \varpi    ^{(\bullet)! }  \circ\bbR \underline{\Gamma} ^\dag _{Y} (\cE ^{ (\bullet)} ) .
\end{equation}
\end{prop}

\begin{proof}This follows by devissage and by Mayer-Vietoris exact triangles from \cite[9.2.6.8]{Car25}. 
\end{proof}

\begin{theo}
\label{2.2.18}
Let  $f \colon \fX ^{\prime } \to \fX $ be a morphism of separated, quasi-compact, smooth $\fS$-formal schemes.
Let $Y$ be a subscheme of $X$,
$Y':= f ^{-1} (Y)$,
$\cE ^{(\bullet)} \in  \smash{\underrightarrow{LD}} ^{\mathrm{b}} _{\bbQ,\mathrm{qc}}  (\overset{^\mathrm{l}}{} \smash{\widehat{\cD}} _{\fX /\fS  } ^{(\bullet)} (X _s))$ 
and $\cE ^{\prime (\bullet)} \in\smash{\underrightarrow{LD}} ^{\mathrm{b}} _{\bbQ,\mathrm{qc}}  (\overset{^\mathrm{l}}{} \smash{\widehat{\cD}} _{\fX ^{\prime }/\fS  } ^{(\bullet)} (X ' _s))$. 
We have the functorial in $Y$ isomorphisms 
\begin{equation} \label{commutfonctcohlocal1}
  f ^{ (\bullet)!}  \circ\R \underline{\Gamma} ^\dag _{Y}(\cE ^{(\bullet)}) 
    \riso
   \R \underline{\Gamma} ^\dag _{Y' }\circ f ^{ (\bullet)!}  (\cE ^{(\bullet)}), 
\qquad \R \underline{\Gamma} ^\dag _{Y}\circ f ^{ (\bullet)} _{+} (\cE ^{\prime (\bullet)})
\riso
f ^{ (\bullet)}  _{+} \circ \R \underline{\Gamma} ^\dag _{Y'}(\cE ^{\prime (\bullet)}).
\end{equation}

\end{theo}

\begin{proof}
We prove the theorem  in a similar way to \cite[9.2.6.8]{Car25}. 
\end{proof}

We can precise Corollary \ref{coro-trace-upre}: 

\begin{coro}
\label{pre-loc-tri-B-t1T}
Let $u \colon \fZ  \to \fX $ be a closed immersion of smooth $\fS$-formal schemes.
\begin{enumerate}[(a)]
\item For any  $\cE ^{(\bullet)} \in \smash{\underrightarrow{LD}} ^{\mathrm{b}} _{\bbQ,\mathrm{qc}}(\overset{^\mathrm{l}}{} \smash{\widehat{\cD}} _{\fX /\fS  } ^{(\bullet)} (X _s))$, we have the isomorphism
\begin{equation} \label{pre-loc-tri-B-t1T-iso}
\R \underline{\Gamma} ^\dag _{Z} (\cE ^{(\bullet)})  \riso  u _{+} ^{ (\bullet)} \circ  u ^{ (\bullet)!} (\cE ^{(\bullet)}).
\end{equation}

\item The functor $u _{+} ^{ (\bullet)} \colon  \smash{\underrightarrow{LD}} ^{\mathrm{b}} _{\bbQ,\mathrm{coh}}(\overset{^\mathrm{l}}{} \smash{\widehat{\cD}} _{\fZ /\fS  } ^{(\bullet)} (Z _s))
\to  \smash{\underrightarrow{LD}} ^{\mathrm{b}} _{\bbQ,\mathrm{coh}}(\overset{^\mathrm{l}}{} \smash{\widehat{\cD}} _{\fX /\fS  } ^{(\bullet)} (X _s))$  is fully faithful.
\end{enumerate}
\end{coro}

\begin{proof}
This is checked in the same way as \cite[13.2.1.5]{Car25}. 
\end{proof}

\subsubsection{Base change isomorphism for coherent complexes and realizable morphisms}
\label{subsec-bc-coh-real}

Let $f \colon \fP ' \to \fP$  be a morphism of smooth $\fS$-formal schemes.
\begin{dfn} \label{realizablefscheme}
We say that $f $ is {\it realizable} if there exist a proper morphism $\pi \colon   \fP'' \to \fP$ of smooth $\cV$-formal schemes,  
an immersion $u \colon   \fP ^{\prime} \hookrightarrow \fP ^{\prime\prime}$ of formal schemes
such that $f = \pi \circ u$

When $\fP= \fS$, we say that $\fP ^{\prime}$ is a realizable smooth $\fS$-formal scheme.
\end{dfn}

\begin{prop}\label{stab-propersuppcoh}
Suppose $f $ is realizable in the sense of \ref{realizablefscheme}.
For any $\cE ^{\prime (\bullet)} \in \underrightarrow{LD} ^{\rmb} _{\bbQ,\coh} ( \widehat{\cD} _{\fP ^{\prime }/\fS} ^{(\bullet)} (P _s'))$
with proper support over $P$ (i.e., if $X'$ is the support of $\cE ^{\prime (\bullet)}$ in the sense \ref{dfn-support} then the composite
$X ' \hookrightarrow  P' \overset{f}{\to} P$ is proper),
 the object
$f ^{ (\bullet)}_+(\cE ^{\prime (\bullet)} ) $
belongs to
$\underrightarrow{LD} ^{\rmb} _{\bbQ,\coh} ( \widehat{\cD} _{\fP/\fS} ^{(\bullet)} (P _s))$.
\end{prop}

\begin{proof}
This is checked in the same way as \cite[13.2.3.4]{Car25}.
\end{proof}

\begin{theo} \label{theo-iso-chgtbase}
Suppose $f$ is realizable. Let $g \colon \fQ \to \fP$ be a smooth morphism.
We denote by $\fQ ':= \fP ' \times _{\fP} \fQ$,  by $f ' \colon \fQ ^{\prime } \to \fQ$ and $g ' \colon \fQ ^{\prime}\to \fP ^{\prime}$  the two canonical projections.
Let $\E ^{\prime (\bullet)} \in \underrightarrow{LD} ^{\rmb} _{\bbQ,\coh} ( \widehat{\cD} _{\fP ^{\prime }/\fS} ^{(\bullet)}(P _S'))$ with proper support over $P$. 
There then exists a canonical isomorphism in 
$\underrightarrow{LD} ^{\rmb} _{\bbQ,\coh} ( \widehat{\cD} _{\fQ/\fS} ^{(\bullet)}(Q _s))$:
\begin{equation}
\label{basechange}
g ^{ (\bullet)!} \circ f ^{(\bullet)}_{+} (\E ^{\prime (\bullet)}) \riso f ^{\prime (\bullet)}_{+} \circ g ^{\prime (\bullet) !} (\E ^{\prime (\bullet)}). 
\end{equation}
\end{theo}

\begin{proof}
This is checked in the same way as \cite[13.2.3.7]{Car25}.
\end{proof}

\subsubsection{Relative duality isomorphism and adjunction for realisable morphisms}
\begin{theo}
[Relative duality isomorphism]
\label{rel-dual-isom}
Let $f \colon \fP ' \to \fP$  be a realizable morphism of smooth $\fS$-formal schemes. For any
$\cE ^{\prime (\bullet)}
\in
\underrightarrow{LD} ^{\rmb} _{\bbQ,\coh} ( \widehat{\cD} _{\fP ^{\prime}/\fS} ^{(\bullet)}(P' _s))$
with proper support over $P$,
we have the isomorphism of $\underrightarrow{LD} ^{\rmb} _{\bbQ,\coh} ( \widehat{\cD} _{\fP /\fS}  ^{(\bullet)}(P _s))$ of the form
$$f ^{ (\bullet)}_+ \circ \bbD (\cE ^{\prime (\bullet)} ) \riso \bbD \circ f ^{ (\bullet)}_+ (\cE ^{\prime (\bullet)} ) .  $$
\end{theo}

\begin{proof}
This is similar to \cite[13.2.4.1]{Car25}.
\end{proof}

\begin{coro}
\label{cor-adj-formulbis}
Let $\cE ' \in D ^{\rmb} _{\coh} (\cD ^{\dag} _{\fP ^{\prime}/\fS} (\hdag P _s') _{\bbQ})$ with proper support over $P$, and 
$\E \in D ^{\rmb} _{\coh} (\cD ^{\dag} _{\fP /\fS} (\hdag P _s) _{\bbQ})$.
We have
the isomorphisms
\begin{gather}   \label{cor-adj-formulbis-bij1}
\bbR \mathcal{H} om _{\cD ^{\dag} _{\fP /\fS} (\hdag P _s) _{\bbQ}} ( f _{+} ( \cE ') , \cE) \riso \bbR f _* \bbR 
\Hom _{\cD ^{\dag} _{\fP ^{\prime}/\fS} (\hdag P _s') _{\bbQ}}  ( \cE ' ,f  ^! ( \cE)).
\\ \label{cor-adj-formulbis-bij2}
\bbR \Hom  _{\cD ^{\dag} _{\fP /\fS} (\hdag P _s) _{\bbQ}} ( f _{+} ( \cE ') , \cE) \riso \bbR \mathrm{Hom}  _{\cD ^{\dag} _{\fP ^{\prime}/\fS} (\hdag P _s') _{\bbQ}} ( \cE ' ,f ^!   ( \cE)).
\end{gather}
\end{coro}

\begin{proof}
Using \ref{rel-dual-isom},
we can copy the proof of \ref{cor-adj-formul} word for word.
\end{proof}

\begin{coro} \label{adj-morphbis}
We have the following properties.
\begin{enumerate}[(a)]
\item Let $\cE ' \in D ^{\rmb} _{\coh} (\cD ^{\dag} _{\fP ^{\prime}/\fS} (\hdag P _s') _{\bbQ})$ with proper support over $P$. 
We have the adjunction morphism $\cE ' \to f ^!f _{+} (\cE ')$.
\item Let $\cE \in D ^\rmb _{\coh} (\cD ^{\dag} _{\fP  } (\hdag P _s ) _{\bbQ})$ such that $f ^!  (\cE) \in D ^\rmb _{\coh} (\cD ^{\dag} _{\fP ^{\prime }} (\hdag P _s ' ) _{\bbQ})$.
We have the adjunction morphism $f _{+} f ^!  (\cE) \to \E$.

\item Suppose $f$ is proper and smooth. Then $f _{+} \colon   D ^\rmb _{\coh} (\cD ^{\dag} _{\fP ^{\prime }} (\hdag P _s ' ) _{\bbQ}) \to D ^\rmb _{\coh} (\cD ^{\dag} _{\fP  } (\hdag P _s ) _{\bbQ})$ 
is a right adjoint functor of $f ^! \colon   D ^\rmb _{\coh} (\cD ^{\dag} _{\fP  } (\hdag P _s ) _{\bbQ}) \to D ^\rmb _{\coh} (\cD ^{\dag} _{\fP ^{\prime }} (\hdag P _s ' ) _{\bbQ})$.
\end{enumerate}
\end{coro}

\section{Stable data of coefficient}
Suppose the residue field $k$ of $\cV$ is a perfect field of characteristic $p>0$
and there exists an automorphism  $\sigma \colon \cV \riso \cV$ which is a lifting of the $r$th Frobenius power of $k$ for some integer $r\geq 1$. 

\subsection{Stability under Grothendieck's six operations}

\subsubsection{Data of coefficients}

\begin{dfn} \label{DVRorientedpre}
We denote by $\mathrm{DVR}  (\cV)$ the category whose objects $\cW$ are  $\cV$-algebras which is a complete DVR of mixed characteristic $(0,p)$ with perfect residue field $l$ 
so that there exists an automorphism  $\tau  \colon \cW \riso \cW$ which is a lifting of the $s$th Frobenius power of $l$ for some integer $s\geq 1$ which commutes with $\sigma$, i.e. 
$f \circ \sigma ^{r}= \tau ^{s} \circ f$, where $f\colon \cV \to \cW$ is the structural morphism (such object can be written by $(\cW,\tau, s)$) ; whose morphisms  $(\cW,\tau, s)\to(\cW',\tau', s')$ are 
homomorphisms of $\cV$-algebras $\alpha \colon \cW \to \cW '$ which commute with Frobenius structures, i.e.,  
$\alpha  \circ \tau ^{s}= \tau ^{s'} \circ \alpha$.
\bigskip 

We define the subcategory $\mathrm{DVR}  (\mathfrak{S} )$ of that of formal schemes over $\mathfrak{S} $  as follows:

1) An object consists of a formal scheme $\mathfrak{T} $ over $\mathfrak{S} $ of the form 
$\Spf \cW [[u]]$, where $\cW$ is an object of $\mathrm{DVR}  (\cV)$ with residue field $l$ such that 
the induced homomorphism $k [[t]] \to l [[u]]$ is a morphism of  complete discrete valuation rings (see the convention \ref{dfn-morph-trait}) and $\cV [[t]] \to \cW [[u]]$ is a homomorphism of $\cV$-algebras. 
Remark that such  $\fT \to \fS$ is a homeomorphism.

2) Let $\mathfrak{T} ' =\Spf \cW '[[u']] $ and $\fT=\Spf \cW [[u]]$ be two objects of $\mathrm{DVR}  (\mathfrak{S} )$. 
A morphism of $\mathrm{DVR}  (\mathfrak{S} )$ of the form $\mathfrak{T} ' \to \fT$  is a morphism of formal schemes over $\fS$ so that the induced morphism
$\cW \to \cW'$ is a morphism of $\mathrm{DVR}  (\cV)$ and $\cW [[u]] \to \cW' [[u']]$ is a homomorphism of $\cW$-algebras.
Remark that such  $\fT' \to \fT$ is a homeomorphism.
\end{dfn}

\begin{dfn}[Orientations]\label{def-perfectification}
Let $\cW,\cW '$ be two objects of $\mathrm{DVR}  (\cV)$.
\begin{enumerate}[(a)]
\item We set $\cW [[u ^{p ^{-\infty}}]]\coloneqq \colim _n \cW [[u ^{p ^{-n}}]]$ and $\cW \{\{ u ^{p ^{-\infty}}\}\} \coloneqq \left ( \colim _n \cW [[u ^{p ^{-n}}]]\right ) ^{\widehat{ }}$ the $p$-adic completion of $\cW [[u ^{p ^{-\infty}}]]$.
\item An {\it orientation} of $\cW [[u]]$ is an homomorphism of commutative rings $\cW [[u]] \to \cA$ such that there exists a commutative diagram of commutative rings of the form : 
\begin{equation}
\xymatrix{ 
{\cW [[u]]} \ar[rd] ^-{} \ar[r] ^-{}  & {\cA }  
\\ {}& { \cW [[ u ^{p ^{-\infty}}]]} \ar[u] ^-{\sim} }
\end{equation}
where the oblique map is the canonical one and the vertical one is an isomorphism. 
\item Let $(\cW [[u]],\cA)$ and $(\cW '[[u']],\cA')$ be respectively an orientation of  $\cW [[u]]$ and $ \cW '[[u']]$.
A {\it morphism of orientations} $(\cW [[u]],\cA)\to (\cW '[[u']],\cA')$ is a commutative diagram of commutative rings of the form: 
\begin{equation}\label{def-perfectification-diag}
\xymatrix{ {\cA} \ar[r] ^-{} & { \cA '} \\ {\cW [[u]]} \ar[r] ^-{}  \ar[u] ^-{} & { \cW '[[u']],  }   \ar[u] ^-{}  }
\end{equation}
where the vertical maps are the structural ones.
\end{enumerate}
\end{dfn}

\begin{empt}\label{fin-sepDVR}
Let $f\colon \Spf \cW ' [[u']] \to \Spf \cW [[u]]$ be a morphism of $\mathrm{DVR}  (\fS)$ and $f _0\colon \Spec l' [[u']] \to \Spec l [[u]]$ be the induced map.

\begin{enumerate}[(a)]
\item The morphism $f$ is finite if and only if $f _0$ is finite. 

\item Suppose $f$ is finite.  The fact that $f$ is surjective means that the induced morphism $\Spec l ' [[u']] \to \Spec l [[u]]$ is flat (\cite[4.3.10]{Liu-livre-02}), i.e. 
$l [[u]] \to l' [[u']]$ is flat. Since $\cW [[u]]$ has no $\pi$-torsion, via \cite[Lemma 2.1]{MonskyWashnitzer} and Krull intersection theorem (see \cite[Theorem 8.9]{matsumura}), this implies that 
$\cW [[u]] \to \cW' [[u']]$ is flat.

\item  Suppose $f _0$ is finite. Since $l$ is perfect, then $\Spec l' [[u']] \to \Spec l [[u]]$ is  étale if and only if $l ((u)) \to l '((u'))$ is separable.

\item Let $d \geq 1$ be an integer.  Suppose $l ((u)) \to l '((u'))$ is finite and separable.
Hence, $l ((u)) \to l '((u'))$ and $l [[u]] \to l '[[u']]$ are relatively perfect, i.e. the commutative squares
\begin{equation}\label{l((u))-ext-square00}
\xymatrix{ {l ((u ))} \ar[r] ^-{} & { l '((u' )) } \\ {l ((u))} \ar[r] ^-{}  \ar[u] ^-{F ^d} & { l '((u')),  }   \ar[u] ^-{F ^d}  }
\xymatrix{ {l [[u ]]} \ar[r] ^-{} & { l '[[u' ]] } \\ {l [[u]]} \ar[r] ^-{}  \ar[u] ^-{F^d} & { l '[[u']],  }   \ar[u] ^-{F ^d}  }
\end{equation}
where the vertical arrows are the $d$th power of the absolute Frobenius $F$, are cocartesian. 
Since $l$ and $l'$ are perfect, then $ F ^d ( l ((u )) )=l ((u  ^{p ^d}))$, $F ^d (l '((u')))= l '((u ^{\prime p ^d}))$ and therefore
$F ^d \colon l ((u )) \to l ((u ))$ and $F ^d \colon l '((u' )) \to l' ((u'))$   are respectively  isomorphic to a morphism of the form 
$l ((u )) \to  l ((u ^{p^{-d}}))$ and $l '((u'))\to l '(((u') ^{p^{-d}}))$ ; and similarly with $[[-]]$ instead of $(( -))$.
This implies the cocartesianity of the diagrams
\begin{equation}\label{l((u))-ext-square0}
\xymatrix{ {l ((u ^{p^{-d}}))} \ar[r] ^-{} & { l '(((u') ^{p^{-d}})) } \\ {l ((u))} \ar[r] ^-{}  \ar[u] ^-{} & { l '((u'))  ,}   \ar[u] ^-{}  }
\xymatrix{ {l [[u ^{p^{-d}}]]} \ar[r] ^-{} & { l '[[(u') ^{p^{-d}}]] } \\ {l [[u]]} \ar[r] ^-{}  \ar[u] ^-{} & { l '[[u']]  }   \ar[u] ^-{}  }
\end{equation}
where vertical arrows are the canonical embeddings. 
Taking the colimit, we get  the cocartesian map:
\begin{equation}\label{l((u))-ext-squarepre}
\xymatrix{ {l [[u ^{p^{-\infty}}]]} \ar[r] ^-{} & { l '[[(u') ^{p^{-\infty}}]] } \\ {l [[u]]} \ar[r] ^-{}  \ar[u] ^-{} & { l '[[u']]  }   \ar[u] ^-{}  }
\end{equation}
where $l [[u ^{p ^{-\infty}}]]\coloneqq \colim _d l [[u ^{p ^{-d}}]]$.
\item \label{fin-sepDVR(e)}
By convention, we get the underlying morphism $(\cW,\tau, s)\to(\cW',\tau', s')$ of $\mathrm{DVR}  (\cV)$. 
Suppose the morphism $\cW \to \cW '$ makes $\cW'$ a Cohen $\cW$-algebra  in the terminology of \cite[0.19.8.1]{EGAIV1}.
Then in the same way as \cite[Theorem 2.5]{MonskyWashnitzer}, the first points  imply that $f$ is finite and étale if and only if so is $f _0\colon \Spec l' [[u']] \to \Spec l [[u]]$.
Suppose $f$ is finite and étale. We have the commutative diagram
\begin{equation}\label{l((u))-ext-square00W}
\xymatrix{ {\cW [[u ]]} \ar[r] ^-{} & { \cW '[[u' ]] } \\ {\cW [[u]]} \ar[r] ^-{}  \ar[u] ^-{\phi} & { \cW '[[u']],  }   \ar[u] ^-{\phi'}  }
\end{equation}
where $\phi (\sum _i a _i u ^i) = \sum _i \tau ^{s'} (a _i) u ^{p ^{s s'}i}$ and $\phi '(\sum _i a ' _i (u') ^i) = \sum _i \tau ^{\prime s} (a '_i) (u') ^{p ^{s s'}i}$, which is a lifting of the right square of \ref{l((u))-ext-square00} with $d= ss'$.
Since \ref{l((u))-ext-square00W} is cocartesian modulo the uniformizer of $\cW$, since our modules are $p$-torsion free and $p$-adically complete, then \ref{l((u))-ext-square00W} is cocartesian. 
Since $\phi ( \cW [[u ]] )=\cW [[u ^{p ^{s s'}} ]] $ and $\phi ' (\cW '[[u']])=\cW '[[(u') ^{p ^{ss'}}]] $,  then $\phi$ and $\phi' $ are respectively isomorphic to  
$\cW [[u]] \to \cW [[u ^{p^{-s s'}}]] $ and $\cW '[[u']]\to \cW '[[(u') ^{p^{-s s'}}]]$. Hence, we get the cocartesian left square:
\begin{equation}\label{l((u))-ext-square}
\xymatrix{ {\cW [[u ^{p^{-s s'}}]]} \ar[r] ^-{} & { \cW '[[(u') ^{p^{-s s'}}]] } \\ {\cW [[u]]} \ar[r] ^-{}  \ar[u] ^-{} & { \cW '[[u']],  }   \ar[u] ^-{}  }
\xymatrix{ {\cW [[u ^{p^{-\infty}}]]} \ar[r] ^-{} & { \cW '[[(u') ^{p^{-\infty}}]] } \\ {\cW [[u]]} \ar[r] ^-{}  \ar[u] ^-{} & { \cW '[[u']].  }   \ar[u] ^-{}  }
\end{equation}
Taking the colim (we can replace $s$ and $s'$ by multiples), we get  the cocartesian map of the right square of \ref{l((u))-ext-square}.
\end{enumerate}
\end{empt}

\begin{dfn} We define the category of oriented rings over ``$(\cV [[  t ]], \cV [[t ^{p ^{-\infty}}]])$'' as follows.  
\begin{enumerate}[(a)]
\item An {\it oriented ring} $(\cW [[u]],\cA)$ over $(\cV [[  t ]], \cV [[t ^{p ^{-\infty}}]])$ is the data of a ring of the form $\cW [[u]] $ 
such that $\Spf \cW [[u]]$ is an object of $\mathrm{DVR}  (\mathfrak{S} )$
together with an orientation $\cW [[u]]  \to \cA$ and of a morphism of orientations 
$(\cW [[u]],\cA)\to (\cV [[  t ]], \cV [[t ^{p ^{-\infty}}]])$. 

\item A morphism of oriented rings $(\cW [[u]],\cA) \to (\cW' [[u']],\cA')$ over $(\cV [[  t ]], \cV [[t ^{p ^{-\infty}}]])$  consists of a 
morphism of orientations $(\cW [[u]],\cA) \to (\cW' [[u']],\cA')$ 
such that the induced map $\Spf \cW ' [[u']] \to \Spf \cW [[u]]$ is a morphism of $\mathrm{DVR}  (\fS)$. 
\end{enumerate}
\end{dfn}

\begin{dfn}
\label{DVR}
We set $\cV [[  t ]] ^\flat \coloneqq \cV \{\{t ^{p ^{-\infty}}\}\}$ and  $\fS^\flat \coloneqq \Spf  \cV [[t ]]^\flat$. 
We define the category $\mathrm{DVR}  ( \fS,\fS^\flat)$  as follows: 
\begin{enumerate}[(a)]
\item An object is the data of a morphism of formal schemes $\fT^\flat= \cW [[u]] ^\flat\to \fT =\Spf \cW [[u]] $ whose corresponding homomorphism
$\cW [[u]] \to \cW [[u]] ^\flat$ is induced by the $p$-adic completion of an oriented ring $(\cW [[u]],\cA)$ over $(\cV [[  t ]], \cV [[t ^{p ^{-\infty}}]])$. Remark that we have the structural commutative diagram
\begin{equation}
\notag
\xymatrix{
{\fT }  \ar[r] ^-{} & {\fS} 
\\  {\fT^\flat }  \ar[r] ^-{} \ar[u] ^-{} &  {\fS^\flat}  \ar[u] ^-{}   }
\end{equation}
of  formal schemes whose maps are homeomorphisms.  Such an object is simply denoted by  $(\fT, \fT ^\flat)$. 

\item A morphism  $(\fT ', \fT ^{\prime \flat})\to (\fT, \fT ^\flat) $ 
is the data of  a morphism $\fT '\to \fT$  of $\mathrm{DVR}  (\fS)$  
so that there exists a morphism of oriented rings $(\cW [[u]],\cA) \to (\cW' [[u']],\cA')$ over $(\cV [[  t ]], \cV [[t ^{p ^{-\infty}}]])$
making by $p$-adic completion the following commutative diagram:
\begin{equation}  \notag
\xymatrix{
{\fT'}   \ar[r] ^-{} & {\fT}
\\  {\fT ^{\prime \flat} }  \ar[r] ^-{} \ar[u] ^-{} &  {\fT ^\flat}  \ar[u] ^-{}   }
\end{equation}
with $\fT ^{\flat} \coloneqq \Spf \widehat{\cA}$ and $\fT ^{\prime \flat} \coloneqq \Spf \widehat{\cA'}$.
\end{enumerate}
\end{dfn}

\begin{dfn}\label{sp-rad-finet}
Let $(f , f ^\flat) \colon (\fT ', \fT ^{\prime \flat})\to (\fT, \fT ^\flat) $ be a morphism of $\mathrm{DVR}  (\fS,\fS ^\flat)$. We have $\fT = \Spf \cW [[u]]$, 
where $\cW$ is an object of $\mathrm{DVR}  (\cV)$ with residue field $l$ and $\fT '= \Spf \cW '[[u']]$,
where $\cW'$ is an object of $\mathrm{DVR}  (\cV)$ with residue field $l'$.
\begin{enumerate}[(a)]
\item We say that $(f , f ^\flat)$ is finite if  $f$ is finite.  
\item We say that $(f , f ^\flat)$ is special and radicial if $f\colon \fT' \to \fT$ is finite and $f ^\flat \colon \fT ^{\prime \flat} \to \fT ^{\flat} $ is the identity and
the map $\cW [[u]]\to \cW '[[u']]$ corresponding to $f$ is isomorphic to
the canonical map $\cW [[u]]\to \cW [[u ^{p ^{-d}}]]$ for some integer $d\in\bbN$.
In particular, the map $\cW \to \cW'$ is an isomorphism and $l((u)) \to l'((u'))$ is radicial which explains the terminology. 

\item We say that $(f , f ^\flat)$ is special and étale if $f$ and $f ^\flat$ are finite and étale 
and the induced morphism $\cW \to \cW '$ makes $\cW'$ a Cohen $\cW$-algebra  in the terminology of \cite[0.19.8.1]{EGAIV1} (probably this is automatic).

\item We say that $(f , f ^\flat)$ is special if  
the map $(f , f ^\flat)$ splits into a special and radicial map $ (\fT ', \fT ^{\prime \flat})\to (\fT '', \fT ^{\prime \prime \flat})$ followed by an special and étale map $ (\fT '', \fT ^{\prime \prime \flat})\to (\fT, \fT ^\flat)$.

\end{enumerate}
\end{dfn}

\begin{lemm}\label{splitspecial}
Let $(\fT, \fT ^\flat)$ be an object of $\mathrm{DVR}  ( \fS,\fS^\flat)$. 
Let $T  ' =\Spec l '[[u']]  \to T=\Spec l [[u]]$ be a surjective finite morphism. 
Then there exists a morphism $(f , f ^\flat)\colon (\fT ', \fT ^{\prime \flat})\to (\fT, \fT ^\flat) $ of $\mathrm{DVR}  ( \fS,\fS^\flat)$ such that $f \colon \fT ' \to \fT$ extends $T' \to T$. 
\end{lemm}

\begin{proof}
Let  $l [[u]] \to l' [[u']]$  be the homomorphisms corresponding to $T' \to T$.
Let  $\cW [[u]]\to \cA$ be an oriented ring whose $p$-adic completion gives $\fT^\flat\to \fT$ and $T' \to T$.
Then we can suppose  $\cW [[u]] \to \cA$ is equal to $\cW [[u]] \to\cW [[u ^{p ^{-\infty}}]]$.

Let $L$ be the separable closure of $l((u))$ in $l'((u'))$, i.e. $L/l((u))$ is separable and $l'((u'))/L$ is radicial.
Let $B$ be the integral closure of $l[[u]]$ in $L$. 
Since $L/l((u))$ is finite then  $B$ is a complete DVR with perfect residue field (see \cite[II. Proposition 3]{Serre-corpslocaux}).
Hence, $L= l''((v))$ and $B= l''[[v]]$ where $l''$ is a subfield of $l'$ containing $l$. 
Since $l'((u'))/l''((v))$ is radicial and $l''$ is perfect, then $l''((v)) \to l'((u'))$ is isomorphic to the canonical map $l''((v)) \to l''(( v ^{p^{-d}}))$ for some integer $d$ (see \cite[II. \S 4 Exercice]{Serre-corpslocaux}). 
In particular $l''=l'$. Let $\cW '$ be the Cohen $\cW$-algebra lifting $l\to l'$ (see \cite[0.19.8.2.(ii)]{EGAIV1}).
Using \cite[0.19.8.2.(i)]{EGAIV1}, we can extend $\cW \to \cW'$ to a morphism of $\mathrm{DVR}  (\cV)$ (see the convention  (see \ref{DVRorientedpre})) of the form $(\cW,\tau, s)\to(\cW',\tau', s)$. 
Set $\fT'' \coloneqq \Spf  \cW' [[v]]$ which is a lifting of $T''=\Spec l'[[v]]$. 

The map $\alpha \colon l[[u]]\to l'[[v]]$ is given by $\alpha (u) = v ^n P _0(v)$ with  $P _0 (u) \in l' [[v]] ^*$ and $n \geq 1$ an integer.
Choose  $P (u) \in \cW' [[v]] $ a lifting of $P _0 (u)$. We get a map $f \colon \cW [[u ]] \to \cW' [[v]]$ which by definition sends $u$ to $v ^n P (u)$.
The induced map $\Spf \cW ' [[v]] \to \Spf \cW [[u]]$ is a morphism of $\mathrm{DVR}  (\fS)$ (see \ref{DVRorientedpre}).
Following the third line of \ref{fin-sepDVR}.(\ref{fin-sepDVR(e)}), $\cW [[u ]] \to \cW' [[v]]$ is finite and étale. Using the right square of \ref{l((u))-ext-square}, we get 
a  special and étale morphism  $ (\fT '', \fT ^{\prime \prime \flat})\to (\fT, \fT ^\flat)$ 
such that $\fT ''\to  \fT $ is associated to the map  $f \colon \cW [[u ]] \to \cW' [[v]]$ and
such that $\fT ''\to  \fT ^{\prime \prime \flat}$ is associated to the canonical map  $\cW '[[v]] \to\cW '[[ v  ^{p ^{-\infty}}]]$.
Moreover, since 
$l''[[v]] \to l'[[u']]$ is isomorphic to the canonical map $l''[[v]] \to l''[[ v ^{p^{-d}}]]$, then we get a special and radicial map 
$ (\fT ', \fT ^{\prime \flat})\to (\fT '', \fT ^{\prime \prime \flat})$ such that $\fT ' \to \fT ''$ is a lifting of $T' \to T''$.
Hence, we are done. 
\end{proof}

\begin{empt}
\label{t-structure-coh}
Let $(\fT, \fT ^\flat)$ be an object of $\mathrm{DVR}  ( \fS,\fS^\flat)$, and 
$\fX$ a smooth $\fT$-formal scheme.
If there is no possible confusion (some confusion might arise if for example we know that $\fT \to \fS$ is finite and etale), 
for any integer $m \in \bbN$, 
we denote 
$\smash{\widehat{\cD}} _{\fX/\T } ^{(m)} (X _s)$
(resp. $\widetilde{\cD} ^\dag _{\fX/\T , \bbQ}$)
simply by 
$\smash{\widehat{\cD}} _{\fX} ^{(m)}(X _s)$
(resp. $\widetilde{\cD} ^\dag _{\fX,\bbQ}$).
Recall the equivalence of categories of \cite[8.4.1.15]{Car25}:
\begin{equation}
\label{limeqcat}
\underrvec{l}_{\bbQ} ^{*}
\colon 
\smash{\underrightarrow{LD}} ^{\mathrm{b}} _{\bbQ,\mathrm{coh}} ( \smash{\widehat{\cD}} _{\fX/\T } ^{(\bullet)} (X _s))
\cong 
D ^{\mathrm{b}}  _{\mathrm{coh}} (\widetilde{\cD} ^\dag _{\fX/\T ,\bbQ}).
\end{equation}

The category 
$D ^{\mathrm{b}}  _{\mathrm{coh}} (\widetilde{\cD} ^\dag _{\fX/\T ,\bbQ})$ 
is endowed with its usual t-structure.
Via \ref{limeqcat}, we get a t-structure on 
$\smash{\underrightarrow{LD}} ^{\mathrm{b}} _{\bbQ,\mathrm{coh}} ( \smash{\widehat{\cD}} _{\fX/\T } ^{(\bullet)}(X _s))$
whose heart is 
$\smash{\underrightarrow{LM}}  _{\bbQ,\mathrm{coh}} ( \smash{\widehat{\cD}} _{\fX/\T } ^{(\bullet)}(X _s))$
(see Notation \cite[2.2.4]{caro-stab-sys-ind-surcoh}).
In fact, via \cite[8.4.5.6]{Car25}, 
we have the canonical explicit cohomological functors
$H ^n 
\colon 
\smash{\underrightarrow{LD}} ^{\mathrm{b}} _{\bbQ,\mathrm{coh}} ( \smash{\widehat{\cD}} _{\fX/\T } ^{(\bullet)}(X _s))
\to 
\smash{\underrightarrow{LM}} _{\bbQ,\mathrm{coh}} ( \smash{\widehat{\cD}} _{\fX/\T } ^{(\bullet)}(X _s))$.
The equivalence of categories 
\ref{limeqcat} commutes with the 
cohomogical functors $H ^n $
(where the cohomogical functors $H ^n $ on 
$D ^{\mathrm{b}}  _{\mathrm{coh}} (\widetilde{\cD} ^\dag _{\fX/\T ,\bbQ})$
are the obvious ones),
 i.e. 
$\underrvec{l}_{\bbQ} ^{*}
H ^n (\cE ^{(\bullet)})$
is  canonically isomorphic
to 
$H ^n (\underrvec{l}_{\bbQ} ^{*} \,\cE ^{(\bullet)})$.

Last but not least, 
via Theorem \cite[8.4.5.6]{Car25}
we have the equivalence of categories 
$\smash{\underrightarrow{LD}} ^{\mathrm{b}} _{\bbQ,\mathrm{coh}} ( \smash{\widehat{\cD}} _{\fX/\T } ^{(\bullet)}(X _s))
\cong 
D ^{\mathrm{b}} _{\mathrm{coh}}
(
\smash{\underrightarrow{LM}} _{\bbQ} ( \smash{\widehat{\cD}} _{\fX/\T } ^{(\bullet)}(X _s))
)
$
which is also compatible with t-structures
(the t-structure on 
$D ^{\mathrm{b}} _{\mathrm{coh}}
(\smash{\underrightarrow{LM}} _{\bbQ} ( \smash{\widehat{\cD}} _{\fX/\T } ^{(\bullet)}(X _s)))$
is the canonical one as the derived category of an abelian category).
\end{empt}

\begin{dfn}
\label{dfn-datacoef}
A {\it data of coefficients $\mathfrak{C}$ over $(\fS, \fS ^\flat)$} will be the data for any object $(\fT, \fT ^\flat)$ of $\mathrm{DVR}  ( \fS,\fS^\flat)$, 
for any   smooth $\fT$-formal scheme $\fX$ 
of a full subcategory of  $\smash{\underrightarrow{LD}} ^{\mathrm{b}} _{\bbQ,\mathrm{coh}} ( \smash{\widehat{\cD}} _{\fX/\T } ^{(\bullet)}(X _s))$,
which will be denoted by $\mathfrak{C} (\fX/(\fT, \fT ^\flat))$, or simply $\mathfrak{C} (\fX)$ if there is no ambiguity with the base $(\fT, \fT ^\flat)$.
If there is no ambiguity with $(\fS, \fS ^\flat)$, we simply say a {\it data of coefficients}.
\end{dfn}

\begin{rem}\label{rem-datacoeff-2dfn}
Remark that $\cV [[t]]\{ \frac{1}{t}\}$ is a complete local ring whose residue field is $k((t))$ and whose maximal ideal is generated by a uniformizer of $\cV$ (which is also not nilpotent in $\cV [[t]]\{ \frac{1}{t}\}$). 
Hence, $\cV [[t]]\{ \frac{1}{t}\}$ is a complete  discrete valuation whose residue field is $k((t))$. Set $\eta \coloneqq \cV [[t]]\{ \frac{1}{t}\}$. Let $\fY$ be a formal scheme over $\eta$. We get 
$\smash{\widehat{\cD}} _{\fY/\fS } ^{(\bullet)}(Y _s)=\smash{\widehat{\cD}} _{\fY/\eta } ^{(\bullet)}$.
Similarly $\cV [[t ^{p ^{-\infty}}]]\{ \frac{1}{t}\}$ is a complete  discrete valuation whose residue field is $k((t ^{p ^{-\infty}}))$.
We get a perfected orientation $(\Spf \cV [[t]]\{ \frac{1}{t}\}, \Spf \cV [[t ^{p ^{-\infty}}]]\{ \frac{1}{t}\})$
of $\Spf \cV [[t]]\{ \frac{1}{t}\}$ in the sense of \cite[11.1.1]{caro-6operations}.
Hence, the restriction of a data of coefficients $\mathfrak{C}$ over $(\fS, \fS ^\flat)$ to formal schemes over $\eta$ induces 
a restricted to perfected orientation of this form of a  data of coefficients over $(\Spf \cV [[t]]\{ \frac{1}{t}\}, \Spf \cV [[t ^{p ^{-\infty}}]]\{ \frac{1}{t}\})$ in the sense of \cite[11.1.4]{caro-6operations}.
\end{rem}

\begin{exs}
\label{ex-Dcst}
We have the following data of coefficients.
\begin{enumerate}[(a)]
\item We define the data of coefficients $\mathfrak{B} _\emptyset$ as follows: 
for any object $(\fT, \fT ^\flat)$ of $\mathrm{DVR}  ( \fS,\fS^\flat)$, 
for any  smooth $\fT$-formal scheme,  
the category $\mathfrak{B} _\emptyset (\fX)$ is the full subcategory of 
$\smash{\underrightarrow{LD}} ^{\mathrm{b}} _{\bbQ,\mathrm{coh}} ( \smash{\widehat{\cD}} _{\fX/\T } ^{(\bullet)}(X _s))$
whose unique object is  $\cB ^{(\bullet)} _{\fX} (X _s)$.

\item We will need the larger data of coefficients $\mathfrak{B} _\mathrm{div}$ 
defined as follows: 
for any object $(\fT, \fT ^\flat)$ of $\mathrm{DVR}  ( \fS,\fS^\flat)$, 
for any  smooth $\fT$-formal scheme,  
the category $\mathfrak{B} _\mathrm{div}(\fX)$ is the full subcategory of 
$\smash{\underrightarrow{LD}} ^{\mathrm{b}} _{\bbQ,\mathrm{coh}} ( \smash{\widehat{\cD}} _{\fX/\T } ^{(\bullet)}(X _s))$
whose objects are of the form  
$\cB ^{(\bullet)} _{\fX} (T)$, 
where $T$ is any divisor of $X$ containing $X _s$. 
From Corollary \ref{coh-ss-div}
we have
$\cB ^{(\bullet)} _{\fX} (T)\in 
\smash{\underrightarrow{LD}} ^{\mathrm{b}} _{\bbQ,\mathrm{coh}} ( \smash{\widehat{\cD}} _{\fX/\T } ^{(\bullet)}(X _s))$. 

\item We define $\fB _\mathrm{cst}$ 
as follows: 
for any object $(\fT, \fT ^\flat)$ of $\mathrm{DVR}  ( \fS,\fS^\flat)$, 
for any  smooth formal scheme $\fX$ over $\fT$,  
the category $\fB _\mathrm{cst}(\fX)$ is the full subcategory of 
$\smash{\underrightarrow{LD}} ^{\mathrm{b}} _{\Q,\mathrm{coh}} ( \widehat{\cD} _{\fX/\fT } ^{(\bullet)}(X _s))$
whose objects are of the form  
$\R \underline{\Gamma} ^\dag _{Y} \cO _\fX ^{(\bullet)} $, 
where $Y$ is a subvariety of the special fiber of $\fX$
and the functor 
$\R \underline{\Gamma} ^\dag _{Y}$ is defined in 
\ref{3.2.1caro-2006-surcoh-surcv}.
Recall following \ref{prop-induction-div-coh}, 
theses objects are coherent.

\item We define 
$\fM _\emptyset$
(resp. $\fM _\mathrm{div}$) as follows: 
for any object $(\fT, \fT ^\flat)$ of $\mathrm{DVR}  ( \fS,\fS^\flat)$, 
for any  smooth formal scheme $\fX$ over $\fT$,  
the category 
$\fM _\emptyset (\fX)$
(resp. $\fM _\mathrm{div}(\fX)$)
is the full subcategory of 
$\smash{\underrightarrow{LD}} ^{\mathrm{b}} _{\Q,\mathrm{coh}} ( \widehat{\cD} _{\fX/\fT } ^{(\bullet)}(X _s))$
consisting of objects 
of the form  
$(\hdag T)  (\E ^{(\bullet)} )$, 
where  $\E ^{(\bullet)} \in  \mathrm{MIC} ^{(\bullet)} (Z _\eta, \fX/K) $ (see notation \ref{ntnMICdag2fs3})
with $Z$ is a $S$-smooth scheme of  $X$ and $Z _\eta \coloneqq Z \setminus Z _s$,
and where $T$ is $X _s$ (resp. a divisor of $X$ containing $X _s$.
\end{enumerate}

\end{exs}

\begin{dfn}
\label{dfn-stable-data}
In order to be precise, let us fix some terminology.
Let $\mathfrak{C}$ and $\mathfrak{D}$ be two data of coefficients over $(\fS, \fS ^\flat)$. 
\begin{enumerate}[(a)]
\item We will say that the data of coefficients $\mathfrak{C}$ is
stable under pushforwards  (resp. projective pushforwards , resp. finite and étale pushforwards)
if for any object $(\fT, \fT ^\flat)$ of $\mathrm{DVR}  ( \fS,\fS^\flat)$, 
for any   {\it realizable} in the sense of \ref{realizablefscheme}  (resp. projective, resp. finite and étale) morphism 
$g \colon \fX ' \to \fX$ of smooth $\fT$-formal schemes, 
for any objet $\cE ^{\prime (\bullet)}$ of $\mathfrak{C} (\fX')$, 
the complex $g _{+} (\cE ^{\prime (\bullet)})$ is an object of  $\mathfrak{C} (\fX)$.

\item We will say that the data of coefficients $\mathfrak{C}$  is stable under extraordinary  pullbacks 
(resp. smooth extraordinary pullbacks, resp. finite étale extraordinary pullbacks) 
if for any object $(\fT, \fT ^\flat)$ of $\mathrm{DVR}  ( \fS,\fS^\flat)$,  for any morphism    (resp. smooth morphism, resp. finite étale morphism)
$f \colon \fY \to \fX$ of  smooth $\fT$-formal schemes,  for any object $\E ^{(\bullet)}$ of $\mathfrak{C} (\fX)$,  we have $f ^{!} (\E ^{(\bullet)})\in \mathfrak{C} (\fY)$.

\item We still say that the data of coefficients $\mathfrak{C}$  satisfies 
$BK ^!$ (resp. $BK _+$) if the following property holds:
for any object $(\fT, \fT ^\flat)$ of $\mathrm{DVR}  ( \fS,\fS^\flat)$, 
for any exact closed immersion  $u \colon \fZ \hookrightarrow \fX$ of smooth $\fT$-formal schemes, 
for any objet $\cE ^{(\bullet)}$ of $\mathfrak{C} (\fX)$ with support in $\fZ$, 
we have $u ^{!} (\cE ^{(\bullet)})\in \mathfrak{C} (\fZ)$
(resp. 
for any objet $\G ^{(\bullet)}$ of $\mathfrak{C} (\fZ)$, 
we have $u _{+} (\G ^{(\bullet)})\in \mathfrak{C} (\fX)$).
Remark that $BK ^!$ and $BK _+$ hold if and only if the data of coefficients $\mathfrak{C}$ satisfies 
(an analogue of) Berthelot-Kashiwara theorem.

\item We will say that the data of coefficients $\mathfrak{C}$ is stable under   (special and radicial) base change if for any  (special and radicial) morphism   
$(\fT ', \fT ^{\prime \flat})\to (\fT, \fT ^\flat) $ of $\mathrm{DVR}  ( \fS,\fS^\flat)$  of $\mathrm{DVR}  (\fS)$,
for any  smooth $\fT$-formal scheme,  for any objet $\cE ^{(\bullet)}$ of $\mathfrak{C} (\fX)$, 
we have $ \cO _{\fT'}  \smash{\widehat{\otimes}}^\L _{\cO _{\fT}}  \E^{(\bullet) } \in \mathfrak{C} (\fX \times _{\fT } \fT ^{\prime})$
(see Notation \ref{comm-chg-base}).

\item We will say that the data of coefficients $\mathfrak{C}$ is stable under special (resp. special and radicial, resp. special and étale)  descent of the base if, for any object $(\fT, \fT ^\flat)$ of $\mathrm{DVR}  ( \fS,\fS^\flat)$,  for any  smooth formal scheme $\fX$ over $\fT$, 
for any object $\E ^{(\bullet)}\in \smash{\underrightarrow{LD}} ^{\mathrm{b}} _{\Q,\mathrm{coh}} ( \widehat{\cD} _{\fX/\fT } ^{(\bullet)}(X _s))$ such that there exists a  special (resp. special and radicial, resp. special and étale) morphism  $(\fT ', \fT ^{\prime \flat})\to (\fT, \fT ^\flat) $ satisfying  $ \cO _{\fT'}  \smash{\widehat{\otimes}}^\L _{\cO _{\fT}}  \E^{(\bullet) } \in \mathfrak{C} (\fX \times _{\fT } \fT ^{\prime})$, we have therefore  $\E^{(\bullet) }\in \mathfrak{C} (\fX )$.

\item We will say that the data of coefficients $\mathfrak{C}$ is stable under tensor products (resp. duals) if
for any object $(\fT, \fT ^\flat)$ of $\mathrm{DVR}  ( \fS,\fS^\flat)$, 
for any  smooth formal scheme $\fX$ over $\fT$, 
for any objects $\E ^{(\bullet)}$ and $\cF ^{(\bullet)}$ of $\mathfrak{C} (\fX)$
we have 
$\cF ^{(\bullet)}
\smash{\widehat{\otimes}}^\L
_{\cO  _{\fX}} \E^{(\bullet) }
\in \mathfrak{C} (\fX)$
(resp. 
$\bbD  _{\fX}(\E^{(\bullet) }) \in  \mathfrak{C} (\fX)$).

\item We will say that the data of coefficients $\mathfrak{C}$ is stable under local cohomological functors
(resp. under localizations outside a divisor), if 
for any object $(\fT, \fT ^\flat)$ of $\mathrm{DVR}  ( \fS,\fS^\flat)$, 
for any  smooth formal scheme $\fX$ over $\fT$, 
for any object $\E ^{(\bullet)}$ of $\mathfrak{C} (\fX)$, 
for any subvariety $Y$ (resp. for any divisor $T$) of the special fiber of $\fX$, 
we have
$\R \underline{\Gamma} ^\dag _{Y} \E ^{(\bullet)} \in \mathfrak{C} (\fX)$ 
(resp. $(\hdag T) (\E ^{(\bullet)} )\in \mathfrak{C} (\fX)$).

\item We will say that the data of coefficients $\mathfrak{C}$ is stable under cohomology if,
for any object $(\fT, \fT ^\flat)$ of $\mathrm{DVR}  ( \fS,\fS^\flat)$, 
for any  smooth formal scheme $\fX$ over $\fT$,
for any object  $\E ^{(\bullet)}$ of $\smash{\underrightarrow{LD}} ^{\mathrm{b}} _{\Q,\mathrm{coh}} ( \widehat{\cD} _{\fX/\fT } ^{(\bullet)}(X _s) )$, the property 
``$\E ^{(\bullet)}$ is an object of  $\mathfrak{C} (\fX)$'' is equivalent to the fact that, for any integer $n$, 
$H ^n (\E ^{(\bullet)})$ is an object of $\mathfrak{C} (\fX)$.

\item We will say that the data of coefficients $\mathfrak{C}$ is stable under shifts if,
for any object $(\fT, \fT ^\flat)$ of $\mathrm{DVR}  ( \fS,\fS^\flat)$, 
for any  smooth formal scheme $\fX$ over $\fT$,
for any object  $\E ^{(\bullet)}$ of $\mathfrak{C} (\fX)$, for any integer $n$, 
$\E ^{(\bullet)} [n]$ is an object of $\mathfrak{C} (\fX)$.

\item We will say that the data of coefficients $\mathfrak{C}$ is stable by devissages if
 $\mathfrak{C}$ is stable by shifts and if
for any object $(\fT, \fT ^\flat)$ of $\mathrm{DVR}  ( \fS,\fS^\flat)$, 
for any  smooth $\fT$-formal scheme, 
for any exact triangle 
$\cE ^{(\bullet)} _1
\to 
\cE ^{(\bullet)} _2
\to 
\cE ^{(\bullet)} _3
\to 
\cE ^{(\bullet)} _1 [1]$
of $\smash{\underrightarrow{LD}} ^{\mathrm{b}} _{\bbQ,\mathrm{coh}} ( \smash{\widehat{\cD}} _{\fX/\T } ^{(\bullet)}(X _s))$, 
if two objects are in $\mathfrak{C} (\fX)$, then so is the third one.

\item We will say that the data of coefficients $\mathfrak{C}$ is stable under direct summands if,
for any object $(\fT, \fT ^\flat)$ of $\mathrm{DVR}  ( \fS,\fS^\flat)$, 
for any  smooth $\fT$-formal scheme we have the following property: 
any direct summand in $\smash{\underrightarrow{LD}} ^{\mathrm{b}} _{\bbQ,\mathrm{coh}} ( \smash{\widehat{\cD}} _{\fX/\T } ^{(\bullet)}(X _s))$
of an object of $\mathfrak{C} (\fX)$ is an object of
$\mathfrak{C} (\fX)$.

\item We say that $\mathfrak{C}$ contains $\mathfrak{D}$ 
(or $\mathfrak{D}$ is contained in $\mathfrak{C}$) 
if for any object $(\fT, \fT ^\flat)$ of $\mathrm{DVR}  ( \fS,\fS^\flat)$, 
for any  smooth $\fT$-formal scheme
the category $\mathfrak{D} (\fX)$ is a full subcategory of $\mathfrak{C} (\fX)$.

\item We say that the data of coefficients $\mathfrak{C}$ is local 
if for any object $(\fT, \fT ^\flat)$ of $\mathrm{DVR}  ( \fS,\fS^\flat)$, 
for any  smooth $\fT$-formal scheme, 
for any open covering $(\fX _i) _{i\in I}$ of $\fX$, 
for any object $\cE ^{(\bullet)}$ of
$\smash{\underrightarrow{LD}} ^{\mathrm{b}} _{\bbQ,\mathrm{qc}} ( \smash{\widehat{\cD}} _{\fX/\T } ^{(\bullet)}(X _s))$, 
we have 
$\cE ^{(\bullet)}\in \mathrm{Ob} \mathfrak{C} (\fX)$ if and only if 
$\cE ^{(\bullet)}| \fX _i \in \mathrm{Ob} \mathfrak{C} (\fX _i)$ for any $i \in I$. 
For instance, 
the data of coefficients 
$\smash{\underrightarrow{LD}} ^{\mathrm{b}} _{\bbQ,\mathrm{coh}}$ is local.

 \item We say that the data of coefficients $\mathfrak{C}$ is quasi-local 
if for any object $(\fT, \fT ^\flat)$ of $\mathrm{DVR}  ( \fS,\fS^\flat)$, 
for any  smooth formal scheme $\fX$ over $\fT$, 
for any open immersion $j \colon \fY \hookrightarrow \fX$ 
for any object 
$\E ^{(\bullet)}\in  \mathfrak{C} (\fX)$, 
we have 
$ j ^{ !(\bullet)}\E ^{(\bullet)} \in  \mathfrak{C} (\fY)$.

\end{enumerate}

\end{dfn}

We finish the subsection with some notation.

\begin{empt}
[Duality]
\label{ntn-dual}
Let $\mathfrak{C}$ be a data of coefficients. We define its dual data of coefficients 
$\mathfrak{C}  ^{\vee}$ as follows: 
for any object $(\fT, \fT ^\flat)$ of $\mathrm{DVR}  ( \fS,\fS^\flat)$, 
for any   smooth $\fT$-formal scheme, 
the category $\mathfrak{C}  ^{\vee} (\fX)$ is the subcategory of 
$\smash{\underrightarrow{LD}} ^{\mathrm{b}} _{\bbQ,\mathrm{coh}} ( \smash{\widehat{\cD}} _{\fX/\T } ^{(\bullet)}(X _s))$
of objects $\E^{(\bullet) }$
such that $\bbD _{\fX} (\E^{(\bullet) }) \in \mathfrak{C} (\fX)$.
\end{empt}

\begin{ntn}
\label{dfn-Delta(C)}
Let $\mathfrak{C}$ be a  data of coefficients.
We denote by $\mathfrak{C} ^+$ 
the smallest data of coefficients containing $\fC$ and stable under shifts.
We define by induction on $n\in \bbN$ the data of coefficients
$\Delta _n (\mathfrak{C} )$  as follows: 
for $n=0$, we put $\Delta _0 (\mathfrak{C} )= \mathfrak{C} ^+$.
Suppose $\Delta _n (\mathfrak{C} )$ constructed for $n\in\bbN$. 
For any object $(\fT, \fT ^\flat)$ of $\mathrm{DVR}  ( \fS,\fS^\flat)$, 
for any  smooth formal scheme $\fX$ over $\fT$,  
the category 
$\Delta _{n+1} (\mathfrak{C} ) (\fX)$
is the full subcategory of 
$\smash{\underrightarrow{LD}} ^{\mathrm{b}} _{\Q,\mathrm{coh}} ( \widehat{\cD} _{\fX/\fT } ^{(\bullet)}(X _s))$
of objects $\E^{(\bullet) }$
such that 
there exists an exact triangle of the form 
$\E^{(\bullet) } \to \cF^{(\bullet) } \to \G ^{(\bullet) } \to \E^{(\bullet) } [1]$ 
where
$\cF^{(\bullet) } $ 
and
$\G ^{(\bullet) } $ 
are objects of 
$\Delta _{n} (\mathfrak{C} ) (\fX)$.
Finally, we put 
$\Delta  (\mathfrak{C} ): = \cup _{n\in\bbN} \Delta _n (\mathfrak{C} )$. 
The data of coefficients 
$\Delta  (\mathfrak{C} )$ is the smallest data of coefficients
containing $\mathfrak{C}$ and stable under devissage. 
\end{ntn}

The following lemma is obvious.
\begin{lemm}
\label{Delta-lemm-stab}
Let $\mathfrak{D}$ be a data of coefficients over $\fS$. 
If $\mathfrak{D}$ is stable under pushforwards
(resp. extraordinary pullbacks, resp. smooth extraordinary  pullbacks, 
resp. tensor products, resp. base change, 
resp. local cohomological functors, resp. localization outside a divisor)
then so is 
$\Delta(\mathfrak{D})$.
If $\mathfrak{D}$ satisfies $BK _+$  (resp. is quasi-local) then so is 
$\Delta(\mathfrak{D})$.
If $\mathfrak{D}$ satisfies $BK ^!$ and is stable under
local cohomological functors then so is 
$\Delta(\mathfrak{D})$.
\end{lemm}

\begin{ex}
\label{stab-cst}
Using the isomorphisms \ref{fonctYY'Gamma-iso}, and Theorem \ref{2.2.18}, we check 
that $\fB _\mathrm{cst} ^+$ satisfies $BK _+$, and is stable under local cohomological functors, 
extraordinary pull-backs 
and tensor products.
\end{ex}

\begin{empt}
Beware also that if $\mathfrak{D}$ is local (resp. stable under cohomology, 
resp. stable under special descent of the base, resp. satisfies $BK ^!$), then it is not clear that 
so is $\Delta(\mathfrak{D})$.

\end{empt}

Since the converse of \ref{Delta-lemm-stab} is not true, let us introduce the following definition.
\begin{dfn}
\label{dfn-DeltaP}
Let $\mathfrak{D}$ be a data of coefficients over $(\fS)$. 
Let $P$ be one of the stability property of 
\ref{dfn-stable-data}.
We say that $\mathfrak{D}$ is 
$\Delta$-stable under $P$ if there exists
a data of coefficients $\fD'$ over $\fS$ 
such that 
$\Delta(\fD ') =\Delta(\mathfrak{D})$  and $\fD '$ is stable under $P$.

\end{dfn}

\begin{lem}
\label{Delta-lemm-stab-bis}
The data of coefficients $\mathfrak{D}$ is $\Delta$-stable under pushforwards
(resp. extraordinary pullbacks, resp. smooth extraordinary  pullbacks, 
resp. tensor products, resp. base change, 
resp. local cohomological functors, resp. localization outside a divisor)
if and only if 
$\Delta(\mathfrak{D})$ is stable under pushforwards
(resp. extraordinary pullbacks, resp. smooth extraordinary  pullbacks, 
resp. tensor products, resp. base change, 
resp. local cohomological functors, resp. localization outside a divisor).
The data of coefficients $\mathfrak{D}$ satisfies $\Delta$-$BK _+$  
(resp. is $\Delta$-quasi local) if and only if 
$\Delta(\mathfrak{D})$ satisfies $BK _+$
(resp. is quasi-local). 

\end{lem}

\begin{proof}
This is a translation of Lemma \ref{Delta-lemm-stab}. 
\end{proof}

Beware, it is not clear that if $\mathfrak{D}$ satisfies $\Delta$-$BK ^!$ and is $\Delta$-stable under
local cohomological functors then 
$\Delta(\mathfrak{D})$ satisfies $BK ^!$.

\subsubsection{Overcoherence, overholonomicity (after any base change) and complements}

\begin{dfn}
\label{dfnS(D,C)}
Let $\mathfrak{C}$ and $\mathfrak{D}$ be two data of coefficients.

\begin{enumerate}[(a)]
\item We denote by 
$S _0 (\mathfrak{D}, \mathfrak{C})$
the data of coefficients defined as follows: 
for any object $(\fT, \fT ^\flat)$ of $\mathrm{DVR}  ( \fS,\fS^\flat)$, 
for any  smooth formal scheme $\fX$ over $\fT$,  
the category 
$S _0 (\mathfrak{D}, \mathfrak{C}) (\fX)$
is the full subcategory of 
$\smash{\underrightarrow{LD}} ^{\mathrm{b}} _{\Q,\mathrm{coh}} ( \smash{\widehat{\cD}} _{\fX/\T } ^{(\bullet)}(X _s))$
of objects  $\E ^{(\bullet)}$
satisfying the following properties :
\begin{enumerate}
\item [($\star$)] for any smooth morphism $f\colon \fY \to \fX$ of  smooth $\fT$-formal schemes, 
for any object 
$\cF ^{(\bullet)}
\in
\mathfrak{D} (\fY)$,
we have 
$\cF ^{(\bullet)} \smash{\widehat{\otimes}}^\L _{\cO ^{(\bullet)}  _{\fY}} f ^{!} (\E^{(\bullet) }) \in \mathfrak{C} (\fY)$.
\end{enumerate}

\item We denote by 
$S(\mathfrak{D}, \mathfrak{C})$ 
the data of coefficients defined as follows: 
for any object $(\fT, \fT ^\flat)$ of $\mathrm{DVR}  ( \fS,\fS^\flat)$, 
for any  smooth formal scheme $\fX$ over $\fT$,  
the category 
$S  (\mathfrak{D}, \mathfrak{C}) (\fX)$ 
is the full subcategory of 
$\smash{\underrightarrow{LD}} ^{\mathrm{b}} _{\Q,\mathrm{coh}} (\smash{\widehat{\cD}} _{\fX/\T } ^{(\bullet)}(X _s))$
of objects  $\E ^{(\bullet)}$
satisfying the following property :
\begin{enumerate}
\item [($\star \star$)] for any morphism  $(\fT ', \fT ^{\prime \flat})\to (\fT, \fT ^\flat) $ of $\mathrm{DVR}  ( \fS,\fS^\flat)$,
we have  $ \cO _{\fT'}  \smash{\widehat{\otimes}}^\L _{\cO _{\fT}}  \E^{(\bullet) }\in  S _0 (\mathfrak{D}, \mathfrak{C}) (\fX \times _{\fT}  \fT ')$.
\end{enumerate}
 \item Let $\sharp$ be a symbol so that  either $S _\sharp = S _0$ or $S _\sharp = S $.
\end{enumerate}
\end{dfn}

\begin{exs}
\label{ex-cst-surcoh}
We get the notions as follows:
\begin{enumerate}[(a)]

\item
We denote by 
$ \smash{\underrightarrow{LD}} ^{\mathrm{b}} _{\Q,\mathrm{ovcoh}}= 
S _0 (\fB _\mathrm{div}, \smash{\underrightarrow{LD}} ^{\mathrm{b}} _{\Q,\mathrm{coh}})$
(see the second example  of \ref{ex-Dcst}).
This notion corresponds in the perfect residue fields case to that of overcoherence as defined in
\cite[5.4]{caro-stab-sys-ind-surcoh}.
We denote by 
$ \smash{\underrightarrow{LD}} ^{\mathrm{b}} _{\Q,\mathrm{oc}}= 
S (\fB _\mathrm{div}, \smash{\underrightarrow{LD}} ^{\mathrm{b}} _{\Q,\mathrm{coh}})$.

\item \label{hstab} 
We put 
$\mathfrak{H} _0 :=S (\fB _\mathrm{div}, \smash{\underrightarrow{LD}} ^{\mathrm{b}} _{\Q,\mathrm{coh}})$ 
and 
by induction on $i \in \bbN$, 
we put $\mathfrak{H} _{i+1} :=
\mathfrak{H} _{i} \cap S(\fB _\mathrm{div}, \mathfrak{H} _{i} ^{\vee})$
(see Notation \ref{ntn-dual}).
The coefficients of 
$\mathfrak{H} _{i}$ are called 
{\it $i$-overholonomic after any base change}.
We get the data of coefficients 
$\smash{\underrightarrow{LD}} ^{\mathrm{b}} _{\Q,\mathrm{h}}:= \mathfrak{H} _{\infty}
:=  \cap _{i\in \bbN} \mathfrak{H} _{i}$
whose objects are called {\it overholonomic after any base change}.

\item \label{ovholstab}
Replacing $S $ by $S _0$ in the definition of $\smash{\underrightarrow{LD}} ^{\mathrm{b}} _{\Q,\mathrm{h}}$, 
we get a data of coefficients, denoted by 
$\smash{\underrightarrow{LD}} ^{\mathrm{b}} _{\Q,\mathrm{ovhol}}$, whose objects are called {\it overholonomic}.

\item Finally, 
we set
$ \smash{\underrightarrow{LM}}  _{\Q,\star}:= 
 \smash{\underrightarrow{LD}} ^{\mathrm{b}} _{\Q,\star} \cap 
  \smash{\underrightarrow{LM}} _{\Q,\mathrm{coh}}$,
  for $\star \in \{\mathrm{ovcoh}, \mathrm{oc}, \mathrm{h}, \mathrm{ovhol} \}$.
\end{enumerate}
\end{exs}

\begin{rem}
\label{rem-overhol}
\begin{enumerate}[(a)]

\item Let $\mathfrak{C}$ be a data of coefficients.
The data of coefficients $\mathfrak{C}$ is stable under smooth extraordinary inverse images, localizations outside a divisor
(resp. under smooth extraordinary inverse images, localizations outside a divisor, and base change) if and only if 
$S _0 (\fB _\mathrm{div}, \mathfrak{C})=\mathfrak{C}$
(resp. $S  (\fB _\mathrm{div}, \mathfrak{C})=\mathfrak{C}$).

\item By construction, 
we remark that 
$\smash{\underrightarrow{LD}} ^{\mathrm{b}} _{\Q,\mathrm{ovhol}}$
is the biggest data of coefficients 
which contains 
$\fB _\mathrm{div}$, 
is stable by devissage, 
dual functors
and the operation
$S _{0} (\fB _\mathrm{div}, -)$. 
Moreover, 
$\smash{\underrightarrow{LD}} ^{\mathrm{b}} _{\Q,\mathrm{h}}$
is the biggest data of coefficients 
which contains 
$\fB _\mathrm{div}$, 
is stable by devissage, 
dual functors
and the operation
$S  (\fB _\mathrm{div}, -)$. 
\end{enumerate}
\end{rem}

\begin{prop} \label{ovcoh-invim-prop}
Let $\mathfrak{C}$ and $\mathfrak{D}$ be two data of coefficients. We suppose $\fD$ contains $\fB _\mathrm{div}$,
satisfies $\Delta$-$BK _+$, and is $\Delta$-stable under extraordinary pullbacks and tensor products
(resp. $\fD$ contains $\fB _\emptyset$, satisfies $\Delta$-$BK _+$, and is $\Delta$-stable under extraordinary pullbacks and local cohomological functors).
 We suppose  $\fC$ is local, satisfies $BK ^!$, is stable under devissages, pushforwards, and direct summands.

In that case, the data of coefficients $S _\sharp (\mathfrak{D}, \mathfrak{C})$ is local, stable under devissages, direct summands,
local cohomological functors, extraordinary pullbacks, pushforwards ($S(\mathfrak{D}, \mathfrak{C})$ is moreover stable under base change).
\end{prop}
\begin{proof}
Copy \cite[19.1.2.12]{Car25} or \cite[11.2.12]{caro-6operations}.
\end{proof}

\begin{coro}
\label{ovcoh-invim}
Let $i \in \bbN \cup \{ \infty\}$.
The data of coefficients
$ \smash{\underrightarrow{LD}} ^{\mathrm{b}} _{\Q,\mathrm{ovcoh}}$
(resp. $ \smash{\underrightarrow{LD}} ^{\mathrm{b}} _{\Q,\mathrm{oc}}$,
resp. 
$\mathfrak{H} _{i}$)
contains
$\fB _\mathrm{cst}$,
is local, stable under devissages, direct summands, 
local cohomological functors,
extraordinary pullbacks, pushforwards (resp. and base change).
Moreover, 
$\smash{\underrightarrow{LD}} ^{\mathrm{b}} _{\Q,\mathrm{h}}$
is stable under duality.
\end{coro}

 \subsubsection{On the stability under cohomology}

\begin{ntn} \label{ntnC0}
Let $\mathfrak{C}$ be a data of coefficients.  We denote by $\mathfrak{C} ^0 $ the data of coefficients defined as follows. 
Let $(\fT, \fT ^\flat)$ be an object of $\mathrm{DVR}  ( \fS,\fS^\flat)$, $\fX$ be a smooth formal scheme  over $\fT$. 
Then $\mathfrak{C} ^0 (\fX) := \mathfrak{C}(\fX) \cap \smash{\underrightarrow{LM}}  _{\bbQ,\coh} ( \widehat{\cD} _{\fX/\fT } ^{(\bullet)}(X _s)) $.
\end{ntn}

\begin{prop}\label{wSerresubcatgen}
Let $\fX$ be a  smooth $\fS$-formal scheme.  Let $\fA$ be an abelian strictly full subcategory of $\smash{\underrightarrow{LM}}  _{\bbQ,\coh} ( \widehat{\cD} _{\fX/\fS } ^{(\bullet)}(X _s))$ containing the zero object.
We denote by $\overline{\fA}$ the (strictly) full subcategory of $\smash{\underrightarrow{LM}}  _{\bbQ,\coh} ( \widehat{\cD} _{\fX/\fS } ^{(\bullet)}(X _s))$ consisting of object 
$\cE ^{(\bullet)}$ such that there exists a filtration $0= \cE _0^{(\bullet)}\subset \cE _1^{(\bullet)}\subset \dots \cE _r ^{(\bullet)}=\cE ^{(\bullet)}$ of length $r \geq 1$ 
of $\cE ^{(\bullet)}$ by objects $\cE _i^{(\bullet)}$ of $\smash{\underrightarrow{LM}}  _{\bbQ,\coh} ( \widehat{\cD} _{\fX/\fS } ^{(\bullet)}(X _s))$  such that 
$\cE _{i+1}^{(\bullet)}/\cE _i^{(\bullet)} $ belong to $\fA$ for any $0\leq i \leq r-1$.  
Then $\overline{\fA}$ is the smallest weak Serre subcategory of $\smash{\underrightarrow{LM}}  _{\bbQ,\coh} ( \widehat{\cD} _{\fX/\fS } ^{(\bullet)}(X _s))$ 
containing $\fA$ (in the sense of \stack{02MP}). 
\end{prop}

\begin{proof}Copy \cite[19.1.2.16]{Car25}.
\end{proof}

\begin{prop}\label{wSerresubcatgen2}
Let $\mathfrak{A}$ be a data of coefficients such that  for any object $(\fT, \fT ^\flat)$ of $\mathrm{DVR}  ( \fS,\fS^\flat)$ (see notation \ref{DVR}), for any 
smooth $\fT$-formal scheme $\fX$, $\fA (\fX)$ is an abelian strictly full subcategory of $\smash{\underrightarrow{LM}}  _{\bbQ,\coh} ( \widehat{\cD} _{\fX/\fT } ^{(\bullet)}(X _s))$.
Then, $\Delta ( \fA)$ is the smallest data of coefficients containing $\fA$ and stable under devissage and cohomology. 
\end{prop}
\begin{proof}
We can copy the proof of \cite[19.1.2.17]{Car25}.
\end{proof}

\begin{prop}\label{wSerresubcatgen-loc}
Let $\mathfrak{A}$ be a data of coefficients such that  for any object $(\fT, \fT ^\flat)$ of $\mathrm{DVR}  ( \fS,\fS^\flat)$ (see notation \ref{DVR}), for any 
smooth $\fT$-formal scheme $\fX$, $\fA (\fX)$ is an abelian strictly full subcategory of $\smash{\underrightarrow{LM}}  _{\bbQ,\coh} ( \widehat{\cD} _{\fX/\fT } ^{(\bullet)}(X _s))$.
We denote by $\overline{\fA} _{\loc}$ the data of coefficients such that  for any object $(\fT, \fT ^\flat)$ of $\mathrm{DVR}  ( \fS,\fS^\flat)$, for any 
smooth $\fT$-formal scheme $\fX$, 
$\overline{\fA} _{\loc}(\fX)$ is  the full subcategory of $\smash{\underrightarrow{LM}}  _{\bbQ,\coh} ( \widehat{\cD} _{\fX/\fT } ^{(\bullet)}(X _s))$  consisting of 
objects which are locally  in $\fX$ in  the smallest weak Serre subcategory of $\smash{\underrightarrow{LM}}  _{\bbQ,\coh} ( \widehat{\cD} _{\fX/\fT } ^{(\bullet)}(X _s))$  containing $\fA (\fX)$ 
(see \ref{wSerresubcatgen}).  
Then, $\Delta (\overline{\fA} _{\loc})$ is the smallest data of coefficients containing $\fA$ which is local and stable under devissage and cohomology. 
\end{prop}
\begin{proof}
We can copy the proof of \cite[19.1.2.18]{Car25}.
\end{proof}

\begin{prop} \label{Bcst-st-cohom}
Let $\mathfrak{C}$  be a data of coefficients. Suppose that  $\fC$ is stable under cohomology, and devissage.
Then $S _\sharp (\fB _\mathrm{cst}^+ , \mathfrak{C}) =S _\sharp (\fB _\mathrm{div}, \mathfrak{C})$ is stable under devissages  and cohomology.
\end{prop}

\begin{proof}Similar to \cite[19.1.2.21]{Car25}.
\end{proof}

\begin{coro}
\label{coro-ovcoh-oc-tstr}
The data of coefficients 
$\smash{\underrightarrow{LD}} ^{\mathrm{b}} _{\Q,\mathrm{ovcoh}}$,
and 
$\smash{\underrightarrow{LD}} ^{\mathrm{b}} _{\Q,\mathrm{oc}}$ 
are stable under cohomology.
\end{coro}

\begin{prop}
\label{prop-st-cohom}
Let $\mathfrak{C}$ and $\mathfrak{D}$ be two data of coefficients.
Suppose that $\fD$ is stable under cohomology, smooth extraordinary pullbacks,
and that $\fC$ is local and stable under cohomology, devissage, extraordinary pullbacks.
Then $S _\sharp (\mathfrak{D}, \mathfrak{C})$ is stable under devissages  and cohomology.
\end{prop}

\begin{proof}Similar to \cite[19.1.2.23]{Car25} (or \cite[11.3.8]{caro-6operations}).
\end{proof}

\begin{empt}[Holonomicity]
Let $(\fT, \fT ^\flat)$ be an object of $\mathrm{DVR}  ( \fS,\fS^\flat)$, $\fX$ be a smooth formal scheme  over $\fT$,
$\E ^{(\bullet)}\in 
\smash{\underrightarrow{LM}}  _{\Q,\mathrm{coh}} 
( \widehat{\cD} _{\fX/\fT } ^{(\bullet)}(X _s))$.
We have the dual functor $\bbD ^{(\bullet)} \colon  \smash{\underrightarrow{LD}} ^{\mathrm{b}} _{\Q,\mathrm{coh}} (  \widehat{\cD} _{\fX/\fT } ^{(\bullet)}(X _s))
\to  \smash{\underrightarrow{LD}} ^{\mathrm{b}} _{\Q,\mathrm{coh}} (  \widehat{\cD} _{\fX/\fT } ^{(\bullet)}(X _s))$.
In the same way as \cite[2.8]{caro-holo-sansFrob},
we say that 
$\E ^{(\bullet)}$ is holonomic if 
for any $i \not = 0$, 
$\H ^i (\bbD ^{(\bullet)}(\E ^{(\bullet)})) =0$.
We denote by 
$\smash{\underrightarrow{LM}}  _{\Q,\mathrm{hol}} (  \widehat{\cD} _{\fX/\fT } ^{(\bullet)}(X _s))$ the strictly subcategory 
of 
$\smash{\underrightarrow{LM}}  _{\Q,\mathrm{coh}} (  \widehat{\cD} _{\fX/\fT } ^{(\bullet)}(X _s))$
of holonomic 
$\widehat{\cD} _{\fX/\fT } ^{(\bullet)}(X _s)$-modules.
By copying \cite[2.14]{caro-holo-sansFrob}, we check 
$\smash{\underrightarrow{LM}}  _{\Q,\mathrm{hol}} (  \widehat{\cD} _{\fX/\fT } ^{(\bullet)}(X _s))$ 
is in fact a Serre subcategory 
of 
$\smash{\underrightarrow{LM}}  _{\Q,\mathrm{coh}} (  \widehat{\cD} _{\fX/\fT } ^{(\bullet)}(X _s))$.

We denote by 
$\smash{\underrightarrow{LD}} ^{\mathrm{b}} _{\Q,\mathrm{hol}}
(  \widehat{\cD} _{\fX/\fT } ^{(\bullet)}(X _s))$ the strictly full subcategory of
$\smash{\underrightarrow{LD}} ^{\mathrm{b}} _{\Q,\mathrm{coh}}
(  \widehat{\cD} _{\fX/\fT } ^{(\bullet)}(X _s))  $ 
consisting of complexes 
$\E ^{(\bullet)}$ such that 
$\H ^n \E ^{(\bullet)} \in 
\smash{\underrightarrow{LM}}  _{\Q,\mathrm{hol}} (  \widehat{\cD} _{\fX/\fT } ^{(\bullet)}(X _s))$
for any 
$n\in\Z$.
This yields the t-exact equivalence of categories
$\bbD ^{(\bullet)}
\colon 
\smash{\underrightarrow{LD}} ^{\mathrm{b}} _{\Q,\mathrm{hol}}
(  \widehat{\cD} _{\fX/\fT } ^{(\bullet)}(X _s))
\cong 
\smash{\underrightarrow{LD}} ^{\mathrm{b}} _{\Q,\mathrm{hol}}
(  \widehat{\cD} _{\fX/\fT } ^{(\bullet)}(X _s))$.
Copying the proof of  \cite[3.3.5]{surcoh-hol} word for word, we check that $\smash{\underrightarrow{LM}}  _{\Q,\mathrm{oc}}  \subset \smash{\underrightarrow{LM}}  _{\Q,\mathrm{hol}}$. This yields
\begin{equation} \label{ocinchol}
\smash{\underrightarrow{LD}} ^{\mathrm{b}} _{\Q,\mathrm{oc}}  \subset \smash{\underrightarrow{LD}} ^{\mathrm{b}}  _{\Q,\mathrm{hol}}.
\end{equation}
\end{empt}

\begin{coro} \label{coro-h-tstr}
Let $i \in \bbN \cup \{ \infty\}$. The data of coefficients $\mathfrak{H} _{i}$ is stable under cohomology.
\end{coro}
\begin{proof}Similar to \cite[19.1.2.26]{Car25}.
\end{proof}

\begin{empt}
\label{t-structure-ovcoh-oc-hol-h}
Let $\mathfrak{C}$ be a data of coefficients stable under devissages  and cohomology. 
Let $(\fT, \fT ^\flat)$ be an object of $\mathrm{DVR}  (\fS,\fS ^\flat)$,
$\fX$ be a smooth formal scheme  over $\fT$.
Recall that following 
\ref{t-structure-coh}
we have a canonical t-structure on 
$ \smash{\underrightarrow{LD}} ^{\mathrm{b}} _{\Q,\mathrm{coh}} 
( \widehat{\cD} _{\fX/\fT } ^{(\bullet)}(X _s))$.
We get a canonical t-structure on 
$\fC (\fX/\fT)$ whose heart is $\fC  ^0  (\fX/\fT)$ 
and so that the t-structure of 
$\fC (\fX/\fT)$ is induced 
by that of 
$ \smash{\underrightarrow{LD}} ^{\mathrm{b}} _{\Q,\mathrm{coh}} 
( \widehat{\cD} _{\fX/\fT } ^{(\bullet)}(X _s))$,
i.e. the truncation functors are the same
and 
$\fC  ^{\geq n}  (\fX/\fT):=
\smash{\underrightarrow{LD}} ^{\geq n} _{\Q,\mathrm{coh}} 
( \widehat{\cD} _{\fX/\fT } ^{(\bullet)}(X _s))
\cap 
\fC  (\fX/\fT)$,
$\fC  ^{\leq n}  (\fX/\fT):=
\smash{\underrightarrow{LD}} ^{\leq n} _{\Q,\mathrm{coh}} 
( \widehat{\cD} _{\fX/\fT } ^{(\bullet)}(X _s))
\cap 
\fC  (\fX/\fT)$.

For instance, using \ref{coro-ovcoh-oc-tstr}
and 
\ref{coro-h-tstr}, 
we get 
  for $\star \in \{\mathrm{ovcoh}, \mathrm{oc}, \mathrm{h},\mathrm{hol} \}$
  a canonical t-structure on 
$\smash{\underrightarrow{LD}} ^{\mathrm{b}} _{\Q,\star} $.
The heart of 
$\smash{\underrightarrow{LD}} ^{\mathrm{b}} _{\Q,\star} $
is  $\smash{\underrightarrow{LM}} ^{\mathrm{b}} _{\Q,\star}$.
\end{empt}

 \subsubsection{On the stability under special descent of base}

\begin{prop} \label{oc-st-sp-desc}
Let $\mathfrak{C}$  be a data of coefficients. Suppose that  $\fC$ is stable under  special and radicial descent of the base. Then so is $S _\sharp (\fB _\mathrm{div}, \mathfrak{C})$.
\end{prop}

\begin{proof}
Since the case where $\sharp= 0$ is similar and easier, let us treat the other case. 
Let $(\fT, \fT ^\flat)$ be an object of $\mathrm{DVR}  ( \fS,\fS^\flat)$, $\fX$ be a smooth formal scheme  over $\fT$,
$\E ^{(\bullet)}\in  \smash{\underrightarrow{LD}} ^{\mathrm{b}} _{\Q,\mathrm{coh}}  ( \widehat{\cD} _{\fX/\fT } ^{(\bullet)}(X _s))$.
We suppose that there exists  a special and radical morphism $ (\fT ', \fT ^{\flat})\to (\fT, \fT ^\flat) $ such that 
$ \cO _{\fT'}  \smash{\widehat{\otimes}}^\L _{\cO _{\fT}}   \E^{(\bullet) }\in  S (\fB _\mathrm{div}, \mathfrak{C})  (\fX \times _{ \fT}  \fT ')$.
Let us check that $\E^{(\bullet) }\in S (\fB _\mathrm{div}, \mathfrak{C})  (\fX )$.

Let $(\widetilde{\fT}, \widetilde{\fT} ^{\flat})\to (\fT, \fT ^\flat)$  be a morphism of  $\mathrm{DVR}  ( \fS,\fS^\flat)$. Set $\widetilde{\fX} := \fX \times _{ \fT}  \widetilde{\fT}$.
We have to check  $\widetilde{\E} ^{(\bullet)}:= \cO _{\widetilde{\fT}}  \smash{\widehat{\otimes}}^\L_{\cO _{\fT}}  \E^{(\bullet) }\in  S _0 (\fB _\mathrm{div}, \mathfrak{C}) (\widetilde{\fX})$.
Let  $f\colon \widetilde{\fY} \to \widetilde{\fX}$  be a smooth morphism of smooth formal $\widetilde{\fT}$-schemes,  $\widetilde{Z} $ be a divisor of  $\widetilde{Y}$.
Let $(\cW [[u]], \cA )\to (\widetilde{\cW} [[\widetilde{u}]]  ,\widetilde{\cA})$ be the homomorphism of oriented rings over $(\cV [[  t ]], \cV [[t ^{p ^{-\infty}}]])$ which 
is associated with  $(\widetilde{\fT}, \widetilde{\fT} ^{\flat})\to (\fT, \fT ^\flat)$.
Then we can suppose  $\cW [[u]] \to \cA$ is equal to $\cW [[u]] \to\cW [[u ^{p ^{-\infty}}]]$
and $\widetilde{\cW} [[\widetilde{u}]] \to \widetilde{\cA} $ is equal to $\widetilde{\cW} [[\widetilde{u}]]\to\widetilde{\cW} [[\widetilde{u} ^{p ^{-\infty}}]]$.
We can suppose the homomorphism corresponding to  $\fT ' \to \fT$ is equal to the canonical map $\cW [[u]]\to \cW [[u ^{p ^{-d}}]]$ for some integer $d\in\bbN$.
For some large enough integer $\widetilde{d}$, the map $\cW [[u ^{p ^{-\infty}}]] \to \widetilde{\cW} [[\widetilde{u} ^{p ^{-\infty}}]]$ induces 
$\cW [[u ^{p ^{-d}}]] \to \widetilde{\cW} [[\widetilde{u} ^{p ^{-\widetilde{d}}}]]$. Write $\widetilde{\fT} '\coloneqq \Spf \widetilde{\cW} [[\widetilde{u} ^{p ^{-\widetilde{d}}}]]$ and we get the commutative diagram of $\mathrm{DVR}  ( \fS,\fS^\flat)$
\begin{equation} \label{oc-st-sp-desc-diag}
\xymatrix{  {(\widetilde{\fT}, \widetilde{\fT} ^{\flat})} \ar[r] ^-{} & { (\fT, \fT ^\flat)}     \\ {(\widetilde{\fT}', \widetilde{\fT} ^{\flat})} \ar[r] ^-{} \ar[u] ^-{} & { (\fT', \fT ^{\flat}),}   \ar[u] ^-{}}
\end{equation}
where vertical maps are special and radicial. 
This yields $\cO _{\widetilde{\fT}'}  \smash{\widehat{\otimes}}^\L_{\cO _{\widetilde{\fT}}}     \widetilde{\E} ^{(\bullet)} \riso  
\cO _{\widetilde{\fT}'}  \smash{\widehat{\otimes}}^\L_{\cO _{\widetilde{\fT}}}      (\cO _{\widetilde{\fT}}  \smash{\widehat{\otimes}}^\L_{\cO _{\fT}} \cE ^{(\bullet)})
\in  S  _0 (\fB _\mathrm{div}, \mathfrak{C}) (\widetilde{\fX} ')$, where $\widetilde{\fX} ' := \fX \times _{ \fT}  \widetilde{\fT} '$. Since base changes commute with extraordinary pullbacks and local cohomological functors 
(see \ref{comm-chg-base}), this implies:
$$\cO _{\widetilde{\fT}'}  \smash{\widehat{\otimes}}^\L_{\cO _{\widetilde{\fT}}}    \left ( \smash{\widetilde{\B}} _{\widetilde{\fY}} ^{(\bullet)} (\widetilde{Z} ) \smash{\widehat{\otimes}}^\L
_{\cO _{\widetilde{\fY}}} f ^{!(\bullet)} (\widetilde{\E}^{(\bullet) }) \right)  \in \fC ( \widetilde{\fY} \times _{ \widetilde{\fT}}  \widetilde{\fT} ').$$ 
Since 
$f ^{!(\bullet)} (\widetilde{\E}^{(\bullet) })
\in
\smash{\underrightarrow{LD}} ^{\mathrm{b}} _{\Q,\mathrm{coh}} (\widetilde{\fY})$,
using \ref{desc-coh-chgbase} and \cite[9.1.6.3]{Car25}, this yields
$\smash{\widetilde{\B}} _{\widetilde{\fY}} ^{(\bullet)} (\widetilde{Z} ) \smash{\widehat{\otimes}}^\L _{\cO _{\widetilde{\fY}}} f ^{!(\bullet)} (\widetilde{\E}^{(\bullet) })
\in \smash{\underrightarrow{LD}} ^{\mathrm{b}} _{\Q,\mathrm{coh}} (\widetilde{\fY})$. Since $\fC$ is stable under special descent of the base, we are done.
\end{proof}

\begin{lemm}
\label{oc-st-sp-desc2-lemm}
Let $\mathfrak{C}$  be a data of coefficients. If $\fC$ is stable under direct summand, finite and étale pushforwards, special and radicial descent of the base
then $\fC$ is stable under special descent of the base.
\end{lemm}
\begin{proof}
By definition of a special morphism (see \ref{sp-rad-finet}), we reduce to the case where $\fC$ is stable under special and étale descent of the base, which is easy.
\end{proof}

\begin{coro}
\label{oc-st-sp-desc2}
The data of coefficients $\smash{\underrightarrow{LD}} ^{\mathrm{b}} _{\Q,\mathrm{coh}}$,
$\smash{\underrightarrow{LD}} ^{\mathrm{b}} _{\Q,\mathrm{ovcoh}}$,
and 
$\smash{\underrightarrow{LD}} ^{\mathrm{b}} _{\Q,\mathrm{oc}}$  
are stable under special descent of the base.
\end{coro}

\begin{proof}This is a consequence of \ref{oc-st-sp-desc} and \ref{oc-st-sp-desc2-lemm}.
\end{proof}

\begin{prop} \label{S(D,C)stability8}
Let $\mathfrak{C}$ and $\mathfrak{D}$ be two data of coefficients. Suppose that $\fD$ is stable under smooth extraordinary pullbacks and special and radicial base change,  and that $\fC$  is included in 
$\smash{\underrightarrow{LD}} ^{\mathrm{b}} _{\Q,\mathrm{ovcoh}}$  and  is  stable under  shifts,  and special and radicial  descent of the base. 
Then  $S(\mathfrak{D}, \mathfrak{C})$ is  stable under special and radicial descent of the base. 
\end{prop}

\begin{proof}
Let $(\fT, \fT ^\flat)$ be an object of $\mathrm{DVR}  ( \fS,\fS^\flat)$, $\fX$ be a smooth formal scheme  over $\fT$, $\E ^{(\bullet)}\in  \smash{\underrightarrow{LD}} ^{\mathrm{b}} _{\Q,\mathrm{coh}}  ( \widehat{\cD} _{\fX/\fT } ^{(\bullet)}(X _s))$.
We suppose that there exists  a special and radical morphism $ (\fT ', \fT ^{\flat})\to (\fT, \fT ^\flat) $ such that 
$ \cO _{\fT'}  \smash{\widehat{\otimes}}^\L _{\cO _{\fT}}   \E^{(\bullet) }\in  S (\mathfrak{D}, \mathfrak{C})   (\fX \times _{ \fT}  \fT ')$.
Let us check that $\E^{(\bullet) }\in S (\mathfrak{D}, \mathfrak{C})  (\fX )$.

Let $(\widetilde{\fT}, \widetilde{\fT} ^{\flat})\to (\fT, \fT ^\flat)$  be a morphism of  $\mathrm{DVR}  ( \fS,\fS^\flat)$. Set $\widetilde{\fX} := \fX \times _{ \fT}  \widetilde{\fT}$.
We have to check  $\widetilde{\E} ^{(\bullet)}:= \cO _{\widetilde{\fT}}  \smash{\widehat{\otimes}}^\L_{\cO _{\fT}}  \E^{(\bullet) }\in  S  _0 (\mathfrak{D}, \mathfrak{C}) (\widetilde{\fX})$.
In the same way as the proof of \ref{oc-st-sp-desc}, we get the commutative diagram
\begin{equation} \notag
\xymatrix{  {(\widetilde{\fT}, \widetilde{\fT} ^{\flat})} \ar[r] ^-{} & { (\fT, \fT ^\flat)}     \\ {(\widetilde{\fT}', \widetilde{\fT} ^{\flat})} \ar[r] ^-{} \ar[u] ^-{} & { (\fT', \fT ^{\flat}),}   \ar[u] ^-{}}
\end{equation}
where vertical maps are special and radicial. 
This yields
$\cO _{\widetilde{\fT}'}  \smash{\widehat{\otimes}}^\L_{\cO _{\widetilde{\fT}}}     \widetilde{\E} ^{(\bullet)} \riso  
\cO _{\widetilde{\fT}'}  \smash{\widehat{\otimes}}^\L _{\cO _{\widetilde{\fT}} }      (\cO _{\widetilde{\fT}}  \smash{\widehat{\otimes}}^\L_{\cO _{\fT}} \cE ^{(\bullet)}) \in  S  _0 (\mathfrak{D}, \mathfrak{C}) (\widetilde{\fX} ')$,
where $\widetilde{\fX} ' := \fX \times _{ \fT}  \widetilde{\fT} '$. Let  $f\colon \widetilde{\fY} \to \widetilde{\fX}$  be a smooth morphism of smooth $\widetilde{\fT}$-formal schemes, 
$\widetilde{\cF} ^{(\bullet)} \in \mathfrak{D} (\widetilde{\fY})$.
Then 
$\widetilde{\cF} ^{(\bullet)}
\smash{\widehat{\boxtimes}}^\L
_{\cO _{\widetilde{\fT}} } f ^{!(\bullet)} (\widetilde{\E}^{(\bullet) })
\in 
\smash{\underrightarrow{LD}} ^{\mathrm{b}} _{\Q,\mathrm{coh}} 
( \smash{\widehat{\D}} _{\widetilde{\fY} \times _{ \widetilde{\fT}} \widetilde{\fY}} ^{(\bullet)})$.
Since $\fD$ is stable under smooth extraordinary pullbacks and special and radicial base change, since base changes commute with (exterior) tensor products and extraordinary pullbacks,
since  $\cO _{\widetilde{\fT}'}  \smash{\widehat{\otimes}}^\L_{\cO _{\widetilde{\fT}}}    \widetilde{\E} ^{(\bullet)} \in  S  _0 (\mathfrak{D}, \mathfrak{C}) (\widetilde{\fX} ')$,
and since  $\fC$ is stable under shifts,  then  $\cO _{\widetilde{\fT}'}  \smash{\widehat{\otimes}}^\L_{\cO _{\widetilde{\fT}}}     \left ( \widetilde{\cF} ^{(\bullet)}
\smash{\widehat{\boxtimes}}^\L _{\cO _{\widetilde{\fT}}}  f ^{!(\bullet)} (\widetilde{\E}^{(\bullet) }) \right ) \in  \fC(  \widetilde{\fT} '\times _{ \widetilde{\fT}} \widetilde{\fY} \times _{ \widetilde{\fT}} \widetilde{\fY}) $.
Since  $\fC$ is stable under  special and radicial descent of the base,  this yields 
$\widetilde{\cF} ^{(\bullet)} \smash{\widehat{\boxtimes}}^\L _{\cO _{\widetilde{\fT}}}  f ^{!(\bullet)} (\widetilde{\E}^{(\bullet) }) \in  \fC( \widetilde{\fY} \times _{ \widetilde{\fT}} \widetilde{\fY}) $.
Since $\fC$
is included in 
$\smash{\underrightarrow{LD}} ^{\mathrm{b}} _{\Q,\mathrm{ovcoh}}$,
since 
$\widetilde{\cF} ^{(\bullet)}
\smash{\widehat{\otimes}}^\L
_{\cO _{\widetilde{\fT}}}  f ^{!(\bullet)} (\widetilde{\E}^{(\bullet) })$
and
$\delta ^{!(\bullet)}
\left (
\widetilde{\cF} ^{(\bullet)}
\smash{\widehat{\boxtimes}}^\L
_{\cO _{\widetilde{\fT}}}  f ^{!(\bullet)} (\widetilde{\E}^{(\bullet) })
\right )$
are isomorphic up to a shift, 
then this implies 
$\widetilde{\cF} ^{(\bullet)}
\smash{\widehat{\otimes}}^\L
_{\cO _{\widetilde{\fY}}} f ^{!(\bullet)} (\widetilde{\E}^{(\bullet) })
\in
\smash{\underrightarrow{LD}} ^{\mathrm{b}} _{\Q,\mathrm{coh}} 
(\widetilde{\fY})$.
Since 
$\cO _{\widetilde{\fT}'}  \smash{\widehat{\otimes}}^\L_{\cO _{\widetilde{\fT}}}    \widetilde{\E} ^{(\bullet)}
\in 
S  _0 (\mathfrak{D}, \mathfrak{C}) (\widetilde{\fX} ')$,
then 
$\cO _{\widetilde{\fT}'}  \smash{\widehat{\otimes}}^\L_{\cO _{\widetilde{\fT}}}   
\left (
\widetilde{\cF} ^{(\bullet)}
\smash{\widehat{\otimes}}^\L
_{\cO _{\widetilde{\fY}}} f ^{!(\bullet)} (\widetilde{\E}^{(\bullet) })
\right )
\in
\mathfrak{C} 
(\widetilde{\fY}\times _{\cO _{\widetilde{\fT}}} \cO _{\widetilde{\fT}'} )$.
Since $\fC$
is closed under special and radicial descent of the base, 
we get
$\widetilde{\cF} ^{(\bullet)}
\smash{\widehat{\otimes}}^\L
_{\cO _{\fY}} f ^{!(\bullet)} (\widetilde{\E}^{(\bullet) })
\in 
\mathfrak{C} 
(\widetilde{\fY})$
i.e.,  $\widetilde{\E} ^{(\bullet)}\in  S  _0 (\mathfrak{D}, \mathfrak{C}) (\widetilde{\fX})$.
\end{proof}

\subsubsection{Constructions of stable data of coefficients}

\begin{dfn} \label{dfn-almostdual}
Let $\mathfrak{D}$ be a data of coefficients over $(\fS, \fS ^\flat)$.
We say that $\mathfrak{D}$ is almost stable under dual functors if the following property holds:
for any data of coefficients $\mathfrak{C}$ over $(\fS, \fS ^\flat)$ 
which is local, stable under special descent of the base, devissages,  direct summands and  projective pushforwards, 
if $\mathfrak{D} \subset \mathfrak{C}$
then 
$\mathfrak{D} ^{\vee} \subset \mathfrak{C}$.
Remark from the biduality isomorphism that 
the inclusion 
$\mathfrak{D} ^{\vee} \subset \mathfrak{C}$
is equivalent to 
the following one 
$\mathfrak{D} \subset  \mathfrak{C} ^\vee$.
\end{dfn}

\begin{ntn} \label{dfnqupre}
Let $\mathfrak{C}, \mathfrak{D}$ be two data of coefficients. We put  $T _0 (\mathfrak{D} ,\mathfrak{C}) :=  S(\mathfrak{D}  ,\mathfrak{C})$.
By induction on $i \in \bbN$, we set  $U _i (\mathfrak{D}  ,\mathfrak{C}):= T _i (\mathfrak{D}  ,\mathfrak{C}) \cap T _i (\mathfrak{D}  ,\mathfrak{C} ) ^{\vee}$, 
$\widetilde{T} _{i} (\mathfrak{D}  ,\mathfrak{C}) :=  S(\mathfrak{D} , U _i (\mathfrak{D}  ,\mathfrak{C}))$
and $T _{i+1} (\mathfrak{D}  ,\mathfrak{C}):=  S(\widetilde{T} _{i} (\mathfrak{D}  ,\mathfrak{C}), \widetilde{T} _{i} (\mathfrak{D}  ,\mathfrak{C}) )$.
We put $T (\mathfrak{D}  ,\mathfrak{C}) := \cap _{i\in \bbN} T _{i} (\mathfrak{D}  ,\mathfrak{C}) $.
\end{ntn}

\begin{thm}
\label{dfnquprop}
Let $\fB _\mathrm{div} \subset \mathfrak{D}\subset \fC$ be two data of coefficients.
We suppose 
\begin{enumerate}[a)]
\item The data $\fD$ is stable under cohomology ;
\item The data $\Delta (\fD)$ satisfies $BK _+$, is stable under extraordinary pullbacks, base change, 
tensor products 
and  is almost stable under dual functors ;
\item The data $\mathfrak{C}$ 
satisfies $BK ^!$, 
is local and stable under devissages, direct summands, 
pushforwards,
cohomology
and
special descent of the base. 
\end{enumerate}
Then, 
the data of coefficients $T(\mathfrak{D}, \mathfrak{C})$ 
(see Definition \ref{dfnqupre})
 is included in 
$\mathfrak{C} $,
contains 
$\mathfrak{D} $,
is 
local,
stable by devissages, direct summands, 
local cohomological functors, 
pushforwards, extraordinary pullbacks, base change, tensor products, duals,
cohomology
and
special descent of the base. 

\end{thm}

\begin{proof}
We can copy \cite[19.1.3.6]{Car25} and \cite[11.6.6]{caro-6operations}
\end{proof}

\begin{ex}
\label{ex-datastableevery}
We can choose 
$\fD = \fB _\mathrm{div} ^+$
and 
$\fC = \smash{\underrightarrow{LD}} ^{\mathrm{b}} _{\Q,\mathrm{coh}}$.
\end{ex}

\begin{ntn} \label{dfnSboxtimes(D,C)}
Let $\mathfrak{C}$ and $\mathfrak{D}$ be two data of coefficients.
\begin{enumerate}[(a)]
\item We denote by  $\boxtimes _0 (\mathfrak{D}, \mathfrak{C})$ the data of coefficients defined as follows:  for any object $(\fT, \fT ^\flat)$ of $\mathrm{DVR}  ( \fS,\fS^\flat)$, 
for any  smooth formal scheme $\fX$ over $\fT$,  the category  $\boxtimes _0 (\mathfrak{D}, \mathfrak{C}) (\fX)$ is the full subcategory of 
$\smash{\underrightarrow{LD}} ^{\mathrm{b}} _{\Q,\mathrm{coh}} ( \widehat{\cD} _{\fX/\fT } ^{(\bullet)}(X _s))$ consisting of objects  $\E ^{(\bullet)}$ satisfying the following property :
\begin{enumerate}
\item [($\star$)] for any smooth $\fT$-formal scheme $\fY$,  for any object $\cF ^{ (\bullet)} \in \mathfrak{D} (\fY)$, we have  $\E ^{(\bullet)} \smash{\widehat{\boxtimes}}^\L _{\cO _{ \fT }}  \cF^{ (\bullet) } \in \mathfrak{C} (\fX\times _{ \fT} \fY)$.
\end{enumerate}

\item We denote by  $\boxtimes (\mathfrak{D}, \mathfrak{C})$ the data of coefficients defined as follows:  for any object $(\fT, \fT ^\flat)$ of $\mathrm{DVR}  ( \fS,\fS^\flat)$, 
for any  smooth formal scheme $\fX$ over $\fT$,   the category  $\boxtimes (\mathfrak{D}, \mathfrak{C}) (\fX)$ is the full subcategory of 
$\smash{\underrightarrow{LD}} ^{\mathrm{b}} _{\Q,\mathrm{coh}} ( \widehat{\cD} _{\fX/\fT } ^{(\bullet)}(X _s))$ consisting of objects  $\E ^{(\bullet)}$ satisfying the following property :
\begin{enumerate}
\item  [($\star \star$)] for any morphism  $(\fT ', \fT ^{\prime \flat})\to (\fT, \fT ^\flat) $ of $\mathrm{DVR}  ( \fS,\fS^\flat)$,
$ \cO _{\fT'}  \smash{\widehat{\otimes}}^\L _{\cO _{\fT}}  \E^{(\bullet) } \in  \boxtimes _0 (\mathfrak{D}, \mathfrak{C}) (\fX \times _{ \cW}  \cW ')$. 
\end{enumerate}
 \item Let $\sharp$ be a symbol so that  either $\boxtimes _\sharp = \boxtimes _0$ or $\boxtimes _\sharp = \boxtimes $.
\end{enumerate}
\end{ntn}

\begin{ntn}
\label{dfnquprebis}
We put 
$T _0 := 
S  (\mathfrak{B} _{\mathrm{div}}  , \smash{\underrightarrow{LD}} ^{\mathrm{b}} _{\Q,\mathrm{coh}} )$.
By induction on $i \in \bbN$, we set 
$U _i:= T _i \cap T _i ^{\vee}$, 
$\widetilde{U}  _i := 
\boxtimes (
\boxtimes ( \mathfrak{M} _{\mathrm{div}} , U _i ),
\boxtimes ( \mathfrak{M} _{\mathrm{div}} , U _i ))$,
and 
$T _{i+1}:= S (\mathfrak{B} _{\mathrm{div}}  ,\widetilde{U}  _i)$.
We put $T  := \cap _{i\in \bbN} T _{i}  $.
\end{ntn}

\begin{thm} \label{dfnqupropbis}
The data of coefficients $T$  contains  $\mathfrak{M} _{\mathrm{div}}$, is  local, stable by devissages, direct summands,  local cohomological functors, 
pushforwards, extraordinary pullbacks,  base change, tensor products, duals, cohomology and special descent of the base. 
\end{thm}

\begin{proof}
We can copy \cite[19.1.3.10]{Car25} or  \cite[11.6.10]{caro-6operations}.
\end{proof}

 \subsection{Grothendieck six operations for arithmetic $\D$-modules over realizable pairs}

Let $(\fT, \fT ^\flat)$ be an object of $\mathrm{DVR}  (\fS,\fS ^\flat)$. We have in particular $\fT = \Spf \cW [[u]]$, where $\cW$ is an object of $\mathrm{DVR}  (\cV)$ with residue field $l$.
We denote by $T= \Spec l [[ u]]$ the special fiber of $\fT$, i.e. its reduction modulo the uniformiser of $\cW$. We denote by $T _s=\Spec l $ and $T _\eta=\Spec l[[u]]$ respectively the closed fiber and the generic fiber of $T$.
We give a  $p$-adic formalism of Grothendieck six operations over pairs of $l$-varieties $(Y,X)$ which  
can be enclosed into a frame of the form $(Y,X,\fP)$ where $\fP$ is a realizable 
smooth formal scheme over $\fT$, $P$ be the special fiber of $\fP$ (i.e. its reduction modulo the uniformiser of $\cW$), 
$X$ is a closed subscheme of $P$, $Y$ is an open of $X _{\eta}\coloneqq X \times _{T} T _{\eta}$. For a data of coefficients $\fC$, a coefficient 
over $(Y,X,\fP)$ is a coefficient over $\fP$ with support in $X$ and having overconvergent singularities along $X \setminus Y$
(i.e. which is isomorphic to  its image via $\R \underline{\Gamma} ^\dag _{X \setminus Y}$). 
For an enough stable data of coefficients $\fC$, we prove the independence with respect to the choice of the frame enclosing $(Y,X)$ of the coefficients of such $\fC$
over $(Y,X,\fP)$  (\ref{ind-CYW}).  When $\fC$ is also  stable under cohomology,  this independence preserves t-structures.  When $X$ is proper over $T$, then 
the coefficients of such $\fC$ over $(Y,X)$ is independent (up to canonical equivalence of categories)  of the choice of a proper variety $X$ enclosing $Y$. This yields a formalism of Grothendieck's six operations over 
$T _{\eta}$-varieties  (which can be enclosed into a frame).

 \subsubsection{Data of coefficients over frames}

\begin{dfn}
We define the category of frames over $\fT$ as follows. 
\begin{enumerate}[(a)]
\item 
A {\it frame} $(Y,X,\fP)$ 
over $\fT$ means that $\fP$ is a realizable (in the sense of \ref{realizablefscheme}) smooth formal scheme over $\fT$, 
$X$ is a closed subscheme of the special fiber $P$ of $\fP$ and $Y$ is an open subscheme of $X _{\eta}$, the generic fiber of $X$. 
Let 
$(Y', X', \fP')$ 
and $(Y,X,\fP)$ be two frames over $\fT$. 
A morphism $\theta= (b,a,f) \colon (Y', X', \fP')\to (Y,X,\fP)$ of frames over $\fT$ 
is the data of a  
morphism $f\colon \fP' \to \fP$ of realizable smooth formal schemes over $\fT$,
a morphism $a\colon X' \to X$ of $T$-schemes, 
and a morphism $b \colon Y' \to Y$ of $T _\eta$-schemes 
inducing the commutative diagram
$$\xymatrix{
{Y'} 
\ar[d] ^-{b}
\ar@{^{(}->}[r] ^-{}
& 
{X'}
\ar[d] ^-{a}
\ar@{^{(}->}[r] ^-{}
& 
{\fP'} 
\ar[d] ^-{f}
\\
{Y} 
\ar@{^{(}->}[r] ^-{}
& 
{X}
\ar@{^{(}->}[r] ^-{}
& 
{\fP.} 
}$$
If there is no ambiguity with $\fT$, we simply say frame or morphism of frames. 

\item A morphism 
$\theta= (b,a,f) \colon (Y', X', \fP')\to (Y,X,\fP)$ of frames over $\fT$ 
is said to be {\it complete} 
(resp. {\it strictly complete})
if $a$ is proper (resp. 
$f$ and $a$ are proper). 

\end{enumerate}

\end{dfn}

\begin{dfn}
\begin{enumerate}[(a)]
\item We define the category of {\it realizable pairs} over $\fT$ as follow. 
A realizable pair $(Y, X)$ over $\fT$ means the two first data of a frame over $\fT$ of the form $(Y,X, \fP)$.
A frame of the form $(Y,X, \fP)$ is said to be enclosing $(Y,X)$. 
A morphism of realizable pairs $u=(b,a)\colon (Y', X') \to (Y, X)$ over $\fT$
is the data of a morphism of $T$-schemes of the form  $a\colon X' \to X$ such that $a (Y' ) \subset Y$ and $b\colon Y'  \to Y$ is the induced morphism.

\item A morphism of realizable pairs $u=(b,a)\colon (Y', X') \to (Y, X)$ over $\fT$
is said to be {\it complete} if $a$ is proper.

\end{enumerate}
 
\end{dfn}

\begin{rem} \label{rem-complete-frame-coup}
The following properties hold.
\begin{enumerate}[(a)]
\item Let $u=(b,a)\colon (Y', X') \to (Y, X)$ be a complete morphism of realizable pairs  over $\fT$.
Then there exists a strictly complete morphism  of frames over $\fT$ of the form
$\theta= (b,a,f) \colon (Y', X', \fP')\to (Y,X,\fP)$.
Indeed, by definition, there exist some frames over $\fT$ of the form
$(Y', X', \fP'')$ and 
$(Y,X,\fP)$. 
There exists an immersion 
$\fP'' \hookrightarrow \fQ ''$ with $\fQ''$ a proper and smooth $\fT$-formal scheme. 
Hence, put $\fP':= \fQ '' \times \fP$ and let $f\colon \fP' \to \fP$ be the projection.
Since $a$ is proper, $X \hookrightarrow \fP$ is proper, 
and $f$ is proper, then the immersion $X' \hookrightarrow \fP'$ is also proper. 

\item Let $u=(b,a)\colon (Y', X') \to (Y, X)$ be a morphism of realizable pairs  over $\fT$.
Similarly, we check that there exists a morphism  of frames over $\fT$ of the form
$\theta= (b,a,f) \colon (Y', X', \fP')\to (Y,X,\fP)$. 

\end{enumerate}

\end{rem}

\begin{ntn}
\label{ntn-6operations}
Let $\mathfrak{C} $ be a data of coefficients over $(\fT,\fT ^\flat)$.
Let  $(Y, X,\fP)$ be a frame over $\fT$. 
We denote by $\mathfrak{C} (Y, \fP/\fT)$
the full subcategory of 
$\mathfrak{C}  (\fP)$ 
of objects $\E$ such that there exists an isomorphism of the form 
$\E \riso \R \underline{\Gamma} ^\dag _{Y} (\E)$.
We remark that $\mathfrak{C} (Y, \fP/\fT)$ only depend on the immersion 
$Y \hookrightarrow \fP$ which explains the notation.
We might choose $X$ equal to the closure of $Y$ in $P$.

\end{ntn}

\begin{ntn}
\label{ntn-t-structureovcoh}
Let $\mathfrak{C}$ be a data of coefficients stable under devissages  and cohomology. 
Let $(Y, X,\fP)$ be a frame over $\fT$.
Choose $\fU$ an open set of $\fP$ such that 
$Y$ is closed in $\fU$.

\begin{enumerate}[(a)]
\item 
We define a canonical t-structure on  $\mathfrak{C} (Y , \fP /\fT)$ as follows.
We denote by 
$\fC  ^{\leq n}    (Y,\fP/\fT)$
(resp. 
$\fC   ^{\geq n}  (Y,\fP/\fT)$)
the full subcategory of 
$\mathfrak{C} (Y , \fP /\fT)$
of complexes 
$\E$ such that
$\E |\fU 
\in 
\fC  ^{\leq n}  (Y,\fU/\fT):=
\fC   (Y,\fU/\fT)
\cap 
\fC  ^{\leq n}  (\fU/\fT)$
(resp. 
$\E |\fU 
\in 
\fC  ^{\geq n}  (Y,\fU/\fT)
:=
\fC   (Y,\fU/\fT)
\cap 
\fC  ^{\geq n}  (\fU/\fT)$),
where the t-structure on 
$\fC  (\fU/\fT)$
is the canonical one (see \ref{t-structure-ovcoh-oc-hol-h}).
The heart of this t-structure
will be denoted by
$\mathfrak{C} ^0   (Y,\fP/\fT)$.
Finally, we denote by 
$\mathcal{H} ^i _{\mathrm{t}}$
the $i$th space of cohomology with respect to this canonical t-structure. 

\item Suppose $Y$ is smooth. 
Then, we denote by
$\fC _{\mathrm{isoc}}   (Y,\fP/\fT)$
(resp. $\fC ^{\geq n}  _{\mathrm{isoc}}   (Y,\fP/\fT)$,
resp. $\fC ^{\leq n}  _{\mathrm{isoc}}   (Y,\fP/\fT)$,
resp. $\fC ^0 _{\mathrm{isoc}}   (Y,\fP/\fT)$)
the full subcategory of 
(resp. $\fC ^{\geq n}  (Y,\fP/\fT)$,
resp. $\fC ^{\leq n}   (Y,\fP/\fT)$,
resp. $\fC ^0  (Y,\fP/\fT)$)
consisting of complexes
$\E ^{(\bullet)}$ such that 
$H ^i  (\E ^{(\bullet)} |\fU ) 
\in 
\mathrm{MIC} ^{(\bullet)} (Y, \fU/K)$.
We refer ``$\mathrm{isoc}$'' as isocrystals. 
The reason is that when we work with varieties over perfect fields, we have the equivalence of categories of \cite[5.4.6.1]{caro-pleine-fidelite} 
In this paper, we avoid trying to check such equivalence of categories (other than the easier case where the partial compactification is smooth).
\end{enumerate}

\end{ntn}

\begin{rem}
\label{rem-tstructure-exact}
Let $\mathfrak{C}$ be a data of coefficients stable under devissages  and cohomology. 
Let $\fP$ be a smooth $\fT$-formal scheme,
$Y$ be a subscheme of $P \times _T T _{\eta}$,
$Z$ be a closed subscheme of $Y$,
and 
$Y':= Y \setminus Z$.

\begin{enumerate}[(a)]
\item We get the t-exact functor 
$(\hdag Z) 
\colon
\mathfrak{C} (Y, \fP/\fT)
\to 
\mathfrak{C} (Y', \fP/\fT)$.
Beware the functor 
$(\hdag Z) 
\colon
\mathfrak{C} (Y, \fP/\fT)
\to 
\mathfrak{C} (Y, \fP/\fT)$
is not always t-exact.

\item 
\label{rem-tstructure-exact2}
We say that $Z$ locally comes from a divisor of $P$ if 
locally in $P$, there exists a divisor $T$ of $P$ such that 
$Z = Y \cap T$ (this is equivalent to saying that locally in $P$, 
the ideal defining $Z \hookrightarrow Y$ is generated by one element).
In that case, 
we get the t-exact functor 
$(\hdag Z) 
\colon
\mathfrak{C} (Y, \fP/\fT)
\to 
\mathfrak{C} (Y, \fP/\fT)$.
Indeed, 
by construction of our t-structures, 
we can suppose $Y$ is closed in $\fP$ (and then we reduce to the case where the t-structure on 
$\mathfrak{C} (Y, \fP/\fT)$ is induced by 
the standard t-structure
of $\smash{\underrightarrow{LD}} ^{\mathrm{b}} _{\Q,\mathrm{coh}} ( \smash{\widehat{\D}} _{\fP} ^{(\bullet)})$).
Since the property is local, 
we can suppose there exists a divisor
 $T$ such that $Z = T \cap Y$. 
Then both functors 
$(\hdag Z) $ and $(\hdag T)$ 
of 
$\mathfrak{C} (Y, \fP/\fT)
\to 
\mathfrak{C} (Y, \fP/\fT)$
are isomorphic. Since 
$(\hdag T)$ is t-exact,
we are done.
\end{enumerate}

\end{rem}

 \subsubsection{Formalism of Grothendieck six operations over realizable pairs}

\begin{thm} [Independence] \label{ind-CYW}
Let $\mathfrak{C}$ be a data of coefficients over $(\fS,\fS ^\flat)$
which contains $\fB _\mathrm{div}$,  which is  stable  under devissages,  pushforwards, extraordinary pullbacks, and  local cohomological functors.
Let $\theta= (b,a,f) \colon (Y', X', \fP')\to (Y,X,\fP)$ be  a morphism of frames over $\fT$ such that $a$ and $b$ are proper.
\begin{enumerate}[(a)]
\item For any 
$\E ^{(\bullet)}  \in  \fC (Y , \fP /\fT)$,  for any $\E ^{\prime (\bullet)} \in \fC (Y ', \fP '/\fT)$  (recall notation \ref{ntn-6operations}), we have 
$$\mathrm{Hom}  _{\fC (Y , \fP /\fT)}
( f ^{(\bullet) } _{+} ( \E ^{\prime (\bullet) }) , \E ^{(\bullet) }) 
\riso 
 \mathrm{Hom}  _{\fC (Y ', \fP '/\fT)}
( \E ^{\prime (\bullet) } ,
\R \underline{\Gamma} ^{\dag} _{Y'}
f ^{!(\bullet )}  ( \E ^{ (\bullet) })).$$

\item Suppose that 
$Y'= Y$ and $b$ is the identity,
and that $\fC$ is stable under cohomology.
 Then, 
for any 
$\E ^{(\bullet)} \in \fC ^0 (Y , \fP /\fT)$,
for any 
$\E ^{\prime(\bullet)} \in 
\fC ^0 (Y , \fP '/\fT)$, 
for any $n \in \Z \setminus \{ 0\}$, we have
$$\mathcal{H} _t ^n \R \underline{\Gamma} ^{\dag} _{Y} f ^{!(\bullet)} ( \E ^{(\bullet)}) =0,
\hspace{1 cm}
\mathcal{H} ^n _t f _{+} ^{(\bullet)} ( \E ^{\prime (\bullet)}) =0.$$

\item Suppose that 
$Y'= Y$ and $b$ is the identity. 
For any 
$\E ^{(\bullet)} 
\in 
\mathfrak{C} (Y , \fP /\fT)$, 
for any $\E ^{\prime (\bullet)} \in \mathfrak{C} (Y , \fP '/\fT)$, 
the adjunction morphisms 
$\R \underline{\Gamma} ^{\dag} _{Y} f ^{!(\bullet)}  f _{+} ^{(\bullet)} (\E ^{\prime (\bullet)})
\to 
\E ^{\prime (\bullet)}$ 
and 
$f _{+} ^{(\bullet)}  \R \underline{\Gamma} ^{\dag} _{Y} f ^{!(\bullet)}  
(\E ^{(\bullet)} )
\to 
\E ^{(\bullet)} $
are isomorphisms.
In particular, the functors
$\R \underline{\Gamma} ^{\dag} _{Y} f ^{!(\bullet)}$
and 
$f _{+} ^{(\bullet)}$ induce quasi-inverse equivalences of categories between 
$\mathfrak{C} (Y , \fP /\fT)$
and 
$\mathfrak{C} (Y , \fP '/\fT)$.

\end{enumerate}

\end{thm}

\begin{proof}
I) Let us check the first statement. 
Replacing $X$ and $X'$ by the closure of $Y$ in $P$ and $Y'$ in $P'$,
we can suppose $Y$ is dense in $X$ and $Y'$ is dense in $X'$.
Let 
$\E ^{(\bullet)} 
\in 
\fC (Y , \fP /\fT)$, 
and $\E ^{\prime (\bullet)} \in \fC (Y ', \fP '/\fT)$.
Since $a$ is proper, using \ref{cor-adj-formulbis-bij2}, 
the stability of $\fC$  under extraodinary pullbacks,
and the equivalence of categories \cite[8.4.5.6]{Car25},
we get the bijection
$$\mathrm{Hom}  _{\smash{\underrightarrow{LD}} ^{\mathrm{b}} _{\Q,\mathrm{coh}} (\smash{\widehat{\D}} _{\fP} ^{(\bullet)})} ( f ^{(\bullet) } _{+} ( \E ^{\prime (\bullet) }) , \E ^{(\bullet) }) 
\riso   \mathrm{Hom}  _{\smash{\underrightarrow{LD}} ^{\mathrm{b}} _{\Q,\mathrm{coh}} (\smash{\widehat{\D}} _{\fP'} ^{(\bullet)})} ( \E ^{\prime (\bullet) } , f ^{!(\bullet )}  ( \E ^{ (\bullet) })).$$
Since $a$ and $b$ are proper, then  the open immersion  $Y' \subset a ^{-1} (Y)$ is proper.  Since $Y'$ is dense in $X'$, this yields  $Y'=a ^{-1} (Y)$.
This implies that the functors  $\R \underline{\Gamma} ^{\dag} _{X'} f ^{!(\bullet )}$ and $\R \underline{\Gamma} ^{\dag} _{Y'} f ^{!(\bullet )}$ 
(resp.  $f ^{(\bullet) } _{+} $  and  $\R \underline{\Gamma} ^{\dag} _{Y} f ^{(\bullet) } _{+} $) are isomorphic on $\fC(Y , \fP /\fT)$ (resp. $\fC(Y ', \fP '/\fT)$). Hence, the functor 
$\R \underline{\Gamma} ^{\dag} _{X'} f ^{!(\bullet )}$ (resp. $f ^{(\bullet) } _{+} $) induces $\R \underline{\Gamma} ^{\dag} _{X'} f ^{!(\bullet )} \colon \fC (Y , \fP /\fT) \to  \fC (Y ', \fP '/\fT)$
(resp.  $f ^{(\bullet) } _{+}  \colon  \fC (Y ', \fP '/\fT) \to  \fC (Y , \fP /\fT)$).
Since $\fC (Y' , \fP' /\fT)$  is a strictly full subcategory of $\smash{\underrightarrow{LD}} ^{\mathrm{b}} _{\Q,\mathrm{coh}}  (\smash{\widehat{\cD}} _{\fP '/\T } ^{(\bullet)}(P ' _s) )$, we conclude using the equality
$$ \mathrm{Hom}  _{\smash{\underrightarrow{LD}} ^{\mathrm{b}} _{\Q,\mathrm{coh}}
(\smash{\widehat{\D}} _{\fP'} ^{(\bullet)})}
( \E ^{\prime (\bullet) } ,
f ^{!(\bullet )}  ( \E ^{ (\bullet) }))
=
 \mathrm{Hom}  _{\fC (Y ', \fP '/\fT)}
( \E ^{\prime (\bullet) } ,
\R \underline{\Gamma} ^{\dag} _{X'}
f ^{!(\bullet )}  ( \E ^{ (\bullet) })).$$

II) Now let us check at the same time the last two statements.
Using the stability properties that $\mathfrak{C}$ satisfies, 
we check that the functors 
$f _+ ^{(\bullet )}\colon 
\mathfrak{C} (Y , \fP' /\fT)
\to 
\mathfrak{C} (Y , \fP /\fT)$
and
$\R \underline{\Gamma} ^{\dag} _{Y} f ^{!(\bullet )}
\colon 
\mathfrak{C} (Y , \fP /\fT)
\to 
\mathfrak{C} (Y , \fP '/\fT)$
are well defined. 
Since $\mathfrak{C}$ is included in $\smash{\underrightarrow{LD}} ^{\mathrm{b}} _{\Q,\mathrm{ovcoh}} $, 
we reduce to check the case where 
$\mathfrak{C}=\smash{\underrightarrow{LD}} ^{\mathrm{b}} _{\Q,\mathrm{ovcoh}} $.
We proceed in the same way as \cite[3.2.6]{caro_surcoherent}:
Choose 
$\fU$ (resp. $\fU'$) an open set of $\fP$ (resp. $\fP '$) such that 
$Y$ is closed in $\fU$ (resp. $Y$ is closed in $\fU'$),
and such that $f (\fU') \subset \fU$.
The functor 
$| \fU \colon 
\smash{\underrightarrow{LD}} ^{\mathrm{b}} _{\Q,\mathrm{ovcoh}} (Y , \fP /\fT)
\to 
\smash{\underrightarrow{LD}} ^{\mathrm{b}} _{\Q,\mathrm{ovcoh}} (Y , \fU /\fT)$
is t-exact, and the same with some primes.
Moreover, 
for any $\E ^{(\bullet)}
\in 
\smash{\underrightarrow{LM}}  _{\Q,\mathrm{ovcoh}} (Y , \fP /\fT)$
(or 
$\E ^{(\bullet)}
\in \smash{\underrightarrow{LD}} ^{\mathrm{b}} _{\Q,\mathrm{ovcoh}} (Y , \fP /\fT)$), 
the property 
$\E ^{(\bullet)}= 0$ is equivalent to 
$\E ^{(\bullet)}|\fU= 0$. Hence, we can suppose 
$\fU= \fP$ and $\fU'= \fP'$, i.e.
$Y \hookrightarrow P$
and 
$Y \hookrightarrow P'$ are closed immersions.
Using the remark \ref{rem-datacoeff-2dfn}, we conclude thanks to \cite[12.2.1]{caro-6operations}.    
\end{proof}

In the same way as \cite[12.2]{caro-6operations} or \cite[19.2.2]{Car25}, we get the following corollaries and lemmas.
\begin{coro}
\label{ind-cat-overcouples}
Let $\mathfrak{C}$ be a data of coefficients over $(\fS,\fS ^\flat)$
which contains $\fB _\mathrm{div}$, 
which is stable  under devissages, 
pushforwards, extraordinary pullbacks,
and local cohomological functors.
Let $\mathbb{Y}:= (Y,X)$ be a realizable pair over $\fT$. 

\begin{enumerate}
\item Choose a frame of the form $(Y, X,\fP)$. 
The category 
$\mathfrak{C} (Y, \fP/\fT)$
 does not depend, up to  canonical equivalence of categories, 
on the choice of the frame $(Y,X, \fP)$ over $\fT$ enclosing $(Y,X)$.
Hence, we can simply write 
$\mathfrak{C} (\mathbb{Y}/\fT)$ 
instead of 
$\mathfrak{C} (Y, \fP/\fT)$
without ambiguity (up to canonical equivalence of categories).

\item If moreover $\fC$ is stable under cohomology, 
then we get a canonical t-structure on 
$\mathfrak{C} (\mathbb{Y}/\fT)$. 
\end{enumerate}

\end{coro}

\begin{ex}
With notation  \ref{ind-cat-overcouples}, if $Y=T$, then $\mathfrak{C} (\mathbb{Y}/\fT)$ is equivalent to the category of $\mathrm{Frac}(\cW [[u]] \{ \frac{1}{u}\})$-vector spaces of finite dimension.
Hence, the objects of $\mathfrak{C} (\mathbb{Y}/\fT)$  have finite de Rham cohomology (the pushforward on the base can be considered as a de Rham cohomology).
\end{ex}

\begin{lem}
\label{ind-dual}
Let $\mathfrak{C}$ be a data of coefficients over $(\fS,\fS ^\flat)$
which contains $\fB _\mathrm{div}$, 
which is stable  under devissages, 
pushforwards, extraordinary pullbacks, 
local cohomological functors, and duals.
Let $\mathbb{Y}:= (Y,X)$ be a realizable pair over $\fT$. 
 Choose a frame of the form $(Y, X,\fP)$. 
The functor 
$\R \underline{\Gamma} ^\dag _{Y} \DD _{\fP} 
\colon 
\mathfrak{C} (Y,\fP/\fT) \to \mathfrak{C} (Y,\fP/\fT)$
does not depend, up to  canonical equivalence, 
on the choice of the frame enclosing $(Y,X)$. 
Hence, 
we will denote by 
$\DD _{\mathbb{Y}}
\colon \mathfrak{C} (\mathbb{Y}/\fT) \to \mathfrak{C} (\mathbb{Y}/\fT)$
the functor 
$\R \underline{\Gamma} ^\dag _{Y} \DD _{\fP}$. 
\end{lem}

\begin{lem}
\label{ind-pushforward-extinv}
Let $\mathfrak{C}$ be a data of coefficients over $(\fS,\fS ^\flat)$
which contains $\fB _\mathrm{div}$, 
which is stable  under devissages, 
pushforwards, extraordinary pullbacks,
and local cohomological functors.
Let  $u=(b,a)\colon (Y', X') \to (Y, X)$ be a morphism of realizable pairs over $\fT$.
Put $\mathbb{Y}:= (Y,X)$ and $\mathbb{Y}':= (Y',X')$.
Let us choose a morphism of frames $\theta= (b,a,f) \colon (Y', X', \fP')\to (Y,X,\fP)$ over $\fT$ enclosing $u$.

\begin{enumerate}[(a)]
\item 
The functor $\theta ^{!(\bullet )} := \R \underline{\Gamma} ^\dag _{Y'} \circ f ^{!(\bullet)}
\colon
\mathfrak{C} (Y,\fP/\fT) \to \mathfrak{C} (Y',\fP'/\fT)$ 
does not depend on the choice of such $\theta$ enclosing $u$
(up to canonical equivalences
of categories).
Hence, it will be denoted by $u ^{!}
\colon \mathfrak{C} (\mathbb{Y}/\fT) \to \mathfrak{C} (\mathbb{Y}'/\fT)$. 

\item Suppose that $u$ is complete, i.e. that $a \colon X' \to X$ is proper. 
The functor $\theta _{+}:= f _+ ^{(\bullet )}
\colon 
\mathfrak{C} (Y',\fP'/\fT) \to \mathfrak{C} (Y,\fP/\fT)$
does not depend on the choice of such $\theta$ enclosing $u$
(up to canonical equivalences
of categories).
Hence, it will be denoted by 
$u _{+} \colon \mathfrak{C} (\mathbb{Y}'/\fT) \to \mathfrak{C} (\mathbb{Y}/\fT)$.
\end{enumerate}

\end{lem}

\begin{lem}
\label{ind-prod-tensor}
Let $\mathfrak{C}$ be a data of coefficients over $(\fS,\fS ^\flat)$
which contains $\fB _\mathrm{div}$, 
which is stable under devissages, 
pushforwards, extraordinary pullbacks, 
and tensor products.
Let $\mathbb{Y}:= (Y,X)$ be a realizable pair over $\fT$. 
 Choose a frame of the form $(Y, X,\fP)$. 
The bifunctor $-\smash{\widehat{\otimes}}^\L
_{\O  _{\fP}} - [-\dim P]
\colon
\mathfrak{C} (Y,\fP/\fT) \times \mathfrak{C} (Y,\fP/\fT) \to \mathfrak{C} (Y,\fP/\fT)$
does not depend, up to the canonical equivalence of categories,
on the choice of the frame enclosing $(Y,X)$.
It will be denoted by 
$\widetilde{\otimes} _{\mathbb{Y}}
\colon 
\mathfrak{C} (\mathbb{Y}/\fT) \times \mathfrak{C} (\mathbb{Y}/\fT) \to \mathfrak{C} (\mathbb{Y}/\fT)$. 
\end{lem}

\begin{empt}
[Grothendieck six operations]
\label{6operations}
Let $\mathfrak{C}$ be a data of coefficients over $(\fS,\fS ^\flat)$
which contains $\fB _\mathrm{div}$, 
which is 
stable  under devissages, 
pushforwards, extraordinary pullbacks, 
duals,
and tensor products.
To sum-up the above Lemmas
we can define Grothendieck six operations on realizable pairs  as follows.
Let  $u=(b,a)\colon (Y', X') \to (Y, X)$ be a morphism of realizable pairs over $\fT$.
Put $\mathbb{Y}:= (Y,X)$ and $\mathbb{Y}':= (Y',X')$.
\begin{enumerate}[(a)]
\item We have the dual functor $\DD _{\mathbb{Y}}
\colon \mathfrak{C} (\mathbb{Y}/\fT) \to \mathfrak{C} (\mathbb{Y}/\fT)$ (see \ref{ind-dual}).

\item 
We have the extraordinary pullback $u ^{!}
\colon
\mathfrak{C} (\mathbb{Y}/\fT) \to \mathfrak{C} (\mathbb{Y}'/\fT)$ (see \ref{ind-pushforward-extinv}).
We get the pullbacks $u ^{+}:= \DD _{\mathbb{Y}'} \circ u ^{!} \circ \DD _{\mathbb{Y}}$.

\item Suppose that $u$ is complete.
Then, we have the functor $u _+
\colon \mathfrak{C} (\mathbb{Y}'/\fT) \to \mathfrak{C} (\mathbb{Y}/\fT)$ (see \ref{ind-pushforward-extinv}). 
We denote by  $u _{!}:= \DD _{\mathbb{Y}} \circ u _{+} \circ \DD _{\mathbb{Y}'}$, 
the extraordinary pushforward by $u$.

\item 
We have the tensor product 
$-\widetilde{\otimes} _{\mathbb{Y}}-
\colon 
\mathfrak{C} (\mathbb{Y}/\fT) \times \mathfrak{C} (\mathbb{Y}/\fT) \to \mathfrak{C} (\mathbb{Y}/\fT)$
(see \ref{ind-prod-tensor})
\end{enumerate}

\end{empt}

\begin{coro}\label{nota-h-ovhol} 
For any realizable pair $\bbY$ over $\fT$, via the data of coefficients $\smash{\underrightarrow{LD}} ^{\mathrm{b}} _{\Q,\mathrm{ovhol}}$
and $\smash{\underrightarrow{LD}} ^{\mathrm{b}} _{\Q,\mathrm{h}}$
(defined respectively at \ref{ex-cst-surcoh}.\ref{hstab} and \ref{ex-cst-surcoh}.\ref{ovholstab})
we get the well defined categories 
$\smash{\underrightarrow{LD}} ^{\mathrm{b}} _{\Q,\mathrm{h}}(\mathbb{Y}/\fT)$
or $\smash{\underrightarrow{LD}} ^{\mathrm{b}} _{\Q,\mathrm{ovhol}}(\mathbb{Y}/\fT)$
which are endowed with five of Grothendieck cohomological operations (the tensor product is a priori missing).
\end{coro}

\begin{coro}
For any realizable pair $\bbY$ over $\fT$, via the data of coefficients $T$ (defined at \ref{dfnquprebis}) we get the well defined categories $T(\mathbb{Y}/\fT)$ which are endowed with a t-structure and Grothendieck six operations.
\end{coro}

\subsubsection{Grothendieck six operations over realizable varieties}

\begin{dfn} [Proper compactification]
\begin{enumerate}[(a)]
\item A frame $(Y,X,\fP)$ over $\fT$ is said to be {\it proper} if $\fP$ is proper. 
The category of proper frames over $\fT$ is the full subcategory of the category
of frames over $\fT$ whose objects are proper frames over $\fT$.

\item The category of {\it proper realizable pairs} over $\fT$ is the full subcategory of 
the category of realizable pairs over $\fT$ whose objects $(Y, X)$ are such that 
$X$ is proper. 
We remark that if $(Y, X)$ is a proper realizable pair over $\fT$  then 
there exists a proper frame over $\fT$ of the form $(Y,X, \fP)$.

\item A {\it realizable variety} over $T _\eta$ is a $T _\eta$-scheme $Y$ such that there exists a proper frame over $\fT$ of the form 
$(Y,X,\fP)$. For such frame $(Y,X,\fP)$, we say that 
the proper frame $(Y,X,\fP)$ encloses $Y$ or that the proper realizable pair $(Y,X)$ encloses $Y$. 
\end{enumerate}

\end{dfn}

\begin{empt} [Grothendieck six operations] \label{6operations-variety}
Let $\mathfrak{C}$ be a data of coefficients over $(\fS,\fS ^\flat)$ which contains $\fB _\mathrm{div}$,  which is  stable  under devissages,  pushforwards, extraordinary pullbacks,  duals, and tensor products.
In the same way as Lemma \ref{ind-cat-overcouples}, we check using Theorem \ref{ind-CYW} that the category  $\mathfrak{C} (Y, \fP/\fT)$ (resp. $\mathfrak{C} (Y,X/\fT)$)
 does not depend, up to  canonical equivalence of categories,  on the choice of the proper frame $(Y,X, \fP)$ (resp. the proper realizable pair $(Y,X)$) over $\fT$ enclosing $Y$.
Hence, we simply denote it by $\mathfrak{C} (Y/\fT)$. As for \ref{6operations}, we can define Grothendieck six operations on realizable varieties  as follows. Let  $u\colon Y'\to Y$ be a morphism of realizable varieties over $\fT$.
\begin{enumerate}[(a)]
\item We have the dual functor $\DD _{Y} \colon \mathfrak{C} (Y/\fT) \to \mathfrak{C} (Y/\fT)$ (see \ref{ind-dual}).

\item  We have the extraordinary pullback $u ^{!} \colon \mathfrak{C} (Y/\fT) \to \mathfrak{C} (Y'/\fT)$ (see \ref{ind-pushforward-extinv}). We get the pullbacks $u ^{+}:= \DD _{Y'} \circ u ^{!} \circ \DD _{Y}$.

\item  We have the functor $u _+ \colon \mathfrak{C} (Y'/\fT) \to \mathfrak{C} (Y/\fT)$ (see \ref{ind-pushforward-extinv}).  We denote by  $u _{!}:= \DD _{Y} \circ u _{+} \circ \DD _{Y'}$,  the extraordinary pushforward by $u$.

\item  We have the tensor product  $-\widetilde{\otimes} _{Y}- \colon  \mathfrak{C} (Y/\fT) \times \mathfrak{C} (Y/\fT) \to \mathfrak{C} (Y/\fT)$ (see \ref{ind-prod-tensor})
\end{enumerate}
\end{empt}

\bibliographystyle{alpha}

\begin{thebibliography}{Bou61b}

\bibitem[Abe18]{Abe-Langlands}
Tomoyuki Abe.
\newblock Langlands correspondence for isocrystals and the existence of
  crystalline companions for curves.
\newblock {\em J. Amer. Math. Soc.}, 31(4):921--1057, 2018.

\bibitem[BBD82]{BBD}
A.~A. Be{\u \i}linson, J.~Bernstein, and P.~Deligne.
\newblock Faisceaux pervers.
\newblock In {\em Analysis and topology on singular spaces, {I} ({L}uminy,
  1981)}, volume 100 of {\em Ast{\'e}risque}, pages 5--171. Soc. Math. France,
  Paris, 1982.

\bibitem[Ber90]{Be0}
Pierre Berthelot.
\newblock {Cohomologie rigide et th\'eorie des $\mathcal{D}$-modules}.
\newblock In {\em $p$-adic analysis (Trento, 1989)}, pages 80--124. Springer,
  Berlin, 1990.

\bibitem[Ber96]{Be1}
Pierre Berthelot.
\newblock ${\mathcal{d}}$-modules arithm\'etiques. {I}. {O}p\'erateurs
  diff\'erentiels de niveau fini.
\newblock {\em Ann. Sci. \'Ecole Norm. Sup. (4)}, 29(2):185--272, 1996.

\bibitem[Ber00]{Be2}
Pierre Berthelot.
\newblock {$\mathcal{D}$}-modules arithm\'etiques. {I}{I}. {D}escente par
  {F}robenius.
\newblock {\em M\'em. Soc. Math. Fr. (N.S.)}, (81):vi+136, 2000.

\bibitem[Ber02]{Beintro2}
Pierre Berthelot.
\newblock {Introduction \`a la th\'eorie arithm\'etique des
  {$\mathcal{D}$}-modules}.
\newblock {\em Ast\'erisque}, (279):1--80, 2002.
\newblock Cohomologies {$p$}-adiques et applications arithm\'etiques, {II}.

\bibitem[Bou61a]{bourbaki3-4}
N.~Bourbaki.
\newblock {\em \'{E}l\'ements de math\'ematique. {F}ascicule
  {X}{X}{V}{I}{I}{I}. {A}lg\`ebre commutative. {C}hapitre 3: {G}raduations,
  filtra- tions et topologies. {C}hapitre 4: {I}d\'eaux premiers associ\'es et
  d\'ecomposition primaire}.
\newblock Hermann, Paris, 1961.

\bibitem[Bou61b]{bourbaki}
N.~Bourbaki.
\newblock {\em {\'{E}l\'ements de math\'ematique. {F}ascicule {X}{X}{V}{I}{I}.
  {A}lg\`ebre commutative. {C}hapitre 1: {M}odules plats.{C}hapitre 2:
  {L}ocalisation}}.
\newblock Herman, Paris, 1961.

\bibitem[Car]{caro-semistable}
Daniel Caro.
\newblock Arithmetic $\mathcal{D}$-modules over laurent series fields: absolute
  case.
\newblock {\em To appear in Memoirs of the AMS}.

\bibitem[Car04]{caro_surcoherent}
Daniel Caro.
\newblock {$\mathcal{D}$}-modules arithm{\'e}tiques surcoh{\'e}rents.
  {A}pplication aux fonctions {L}.
\newblock {\em Ann. Inst. Fourier, Grenoble}, 54(6):1943--1996, 2004.

\bibitem[Car05]{caro_comparaison}
Daniel Caro.
\newblock Comparaison des foncteurs duaux des isocristaux surconvergents.
\newblock {\em Rend. Sem. Mat. Univ. Padova}, 114:131--211, 2005.

\bibitem[Car06]{caro_courbe-nouveau}
Daniel Caro.
\newblock Fonctions {L} associ{\'e}es aux {$\mathcal{D}$}-modules
  arithm{\'e}tiques. {C}as des courbes.
\newblock {\em Compositio Mathematica}, 142(01):169--206, 2006.

\bibitem[Car09a]{caro-construction}
Daniel Caro.
\newblock {Arithmetic $\mathcal D$-modules associated with overconvergent
  isocrystals. Smooth case. ($\mathcal D$-modules arithm{\'e}tiques
  associ{\'e}s aux isocristaux surconvergents. Cas lisse.)}.
\newblock {\em Bull. Soc. Math. Fr.}, 137(4):453--543, 2009.

\bibitem[Car09b]{caro_log-iso-hol}
Daniel Caro.
\newblock {Overconvergent log-isocrystals and holonomy. (Log-isocristaux
  surconvergents et holonomie.)}.
\newblock {\em Compos. Math.}, 145(6):1465--1503, 2009.

\bibitem[Car11a]{caro-holo-sansFrob}
Daniel Caro.
\newblock Holonomie sans structure de {F}robenius et crit\`eres d'holonomie.
\newblock {\em Ann. Inst. Fourier (Grenoble)}, 61(4):1437--1454 (2012), 2011.

\bibitem[Car11b]{Caro:2011fk}
Daniel Caro.
\newblock Overcoherence implies holonomicity.
\newblock 03 2011.

\bibitem[Car11c]{caro-pleine-fidelite}
Daniel Caro.
\newblock {Pleine fid{\'e}lit{\'e} sans structure de Frobenius et isocristaux
  partiellement surconvergents}.
\newblock {\em Math. Ann.}, 349:747--805, 2011.

\bibitem[Car11d]{caro-stab-holo}
Daniel Caro.
\newblock Stabilit\'e de l'holonomie sur les vari\'et\'es quasi-projectives.
\newblock {\em Compos. Math.}, 147(6):1772--1792, 2011.

\bibitem[Car16a]{surcoh-hol}
Daniel Caro.
\newblock La surcoh\'{e}rence entra\^{\i}ne l'holonomie.
\newblock {\em Bull. Soc. Math. France}, 144(3):429--475, 2016.

\bibitem[Car16b]{caro-stab-sys-ind-surcoh}
Daniel Caro.
\newblock Syst\`emes inductifs coh\'{e}rents de {$\mathcal D$}-modules
  arithm\'{e}tiques logarithmiques, stabilit\'{e} par op\'{e}rations
  cohomologiques.
\newblock {\em Doc. Math.}, 21:1515--1606, 2016.

\bibitem[Car18]{caro-unip}
Daniel Caro.
\newblock Unipotent monodromy and arithmetic {$\mathcal {D}$}-modules.
\newblock {\em Manuscripta Math.}, 156(1-2):81--115, 2018.

\bibitem[Car19]{caro-6operations}
Daniel Caro.
\newblock Arithmetic $\mathcal{D}$-modules over algebraic varieties of
  characteristic $p>0$.
\newblock {\em ArXiv Mathematics e-prints}, 2019.

\bibitem[Car25]{Car25}
Daniel Caro.
\newblock \textit{Arithmetic $\mathcal{D}$-modules}, available at
  \url{https://carod.users.lmno.cnrs.fr/Arithmetic D-modules.pdf}.
\newblock 2025.

\bibitem[CT12]{caro-Tsuzuki}
Daniel Caro and Nobuo Tsuzuki.
\newblock Overholonomicity of overconvergent {$F$}-isocrystals over smooth
  varieties.
\newblock {\em Ann. of Math. (2)}, 176(2):747--813, 2012.

\bibitem[CV22]{Caro-Vauclair}
Daniel Caro and David Vauclair.
\newblock Logarithmic {$p$}-bases and arithmetical differential modules.
\newblock {\em Ann. Inst. Fourier (Grenoble)}, 72(3):1011--1096, 2022.

\bibitem[dJ96]{dejong}
A.~J. de~Jong.
\newblock Smoothness, semi-stability and alterations.
\newblock {\em Inst. Hautes \'Etudes Sci. Publ. Math.}, (83):51--93, 1996.

\bibitem[FK18]{FujiwaraKatoBookI}
Kazuhiro Fujiwara and Fumiharu Kato.
\newblock {\em Foundations of rigid geometry. {I}}.
\newblock EMS Monographs in Mathematics. European Mathematical Society (EMS),
  Z\"{u}rich, 2018.

\bibitem[Ful69]{Fulton-noetherian}
William Fulton.
\newblock A note on weakly complete algebras.
\newblock {\em Bull. Amer. Math. Soc.}, 75:591--593, 1969.

\bibitem[Gar98]{Garnier-descente-frov-expl}
L.~Garnier.
\newblock {Descente par Frobenius explicite pour les
  ${\mathcal{D}}^\dag$-modules}.
\newblock {\em J. Algebra}, 205(2):542--577, 1998.

\bibitem[Gro60]{EGAI}
A.~Grothendieck.
\newblock \'{E}l\'ements de g\'eom\'etrie alg\'ebrique. {I}. {L}e langage des
  sch\'emas.
\newblock {\em Inst. Hautes \'Etudes Sci. Publ. Math.}, (4):228, 1960.

\bibitem[Gro61]{EGAIII1}
A.~Grothendieck.
\newblock \'{E}l\'ements de g\'eom\'etrie alg\'ebrique. {III}. \'{E}tude
  cohomologique des faisceaux coh\'erents. {I}.
\newblock {\em Inst. Hautes \'Etudes Sci. Publ. Math.}, (11):167, 1961.

\bibitem[Gro64]{EGAIV1}
A.~Grothendieck.
\newblock \'{E}l\'ements de g\'eom\'etrie alg\'ebrique. {IV}. \'{E}tude locale
  des sch\'emas et des morphismes de sch\'emas. {I}.
\newblock {\em Inst. Hautes \'Etudes Sci. Publ. Math.}, (20):259, 1964.

\bibitem[Gro66]{EGAIV3}
A.~Grothendieck.
\newblock \'{E}l\'ements de g\'eom\'etrie alg\'ebrique. {IV}. \'{E}tude locale
  des sch\'emas et des morphismes de sch\'emas. {III}.
\newblock {\em Inst. Hautes \'Etudes Sci. Publ. Math.}, (28):255, 1966.

\bibitem[Gro67]{EGAIV4}
A.~Grothendieck.
\newblock \'{E}l\'ements de g\'eom\'etrie alg\'ebrique. {IV}. \'{E}tude locale
  des sch\'emas et des morphismes de sch\'emas {IV}.
\newblock {\em Inst. Hautes \'Etudes Sci. Publ. Math.}, (32):361, 1967.

\bibitem[Har66]{HaRD}
Robin Hartshorne.
\newblock {\em Residues and duality}.
\newblock Springer-Verlag, Berlin, 1966.

\bibitem[Isa09]{Isaacs}
I.~Martin Isaacs.
\newblock {\em Algebra: a graduate course}, volume 100 of {\em Graduate Studies
  in Mathematics}.
\newblock American Mathematical Society, Providence, RI, 2009.
\newblock Reprint of the 1994 original.

\bibitem[Kat89]{Kato-logFontaine-Illusie}
Kazuya Kato.
\newblock Logarithmic structures of {F}ontaine-{I}llusie.
\newblock In {\em Algebraic analysis, geometry, and number theory (Baltimore,
  MD, 1988)}, pages 191--224. Johns Hopkins Univ. Press, Baltimore, MD, 1989.

\bibitem[Liu02]{Liu-livre-02}
Qing Liu.
\newblock {\em Algebraic geometry and arithmetic curves}, volume~6 of {\em
  Oxford Graduate Texts in Mathematics}.
\newblock Oxford University Press, Oxford, 2002.
\newblock Translated from the French by Reinie Ern{{\'e}}, Oxford Science
  Publications.

\bibitem[LP16]{Lazda-Pal-Book}
Christopher {Lazda} and Ambrus {P\'al}.
\newblock {\em {Rigid cohomology over Laurent series fields.}}
\newblock Cham: Springer, 2016.

\bibitem[LS07]{LeStum-livreRigCoh}
Bernard Le~Stum.
\newblock {\em Rigid cohomology}, volume 172 of {\em Cambridge Tracts in
  Mathematics}.
\newblock Cambridge University Press, Cambridge, 2007.

\bibitem[LvO96]{ZarFilt}
Huishi Li and Freddy van Oystaeyen.
\newblock {\em Zariskian filtrations}, volume~2 of {\em $K$-Monographs in
  Mathematics}.
\newblock Kluwer Academic Publishers, Dordrecht, 1996.

\bibitem[Mat89]{matsumura}
Hideyuki Matsumura.
\newblock {\em Commutative ring theory}, volume~8 of {\em Cambridge Studies in
  Advanced Mathematics}.
\newblock Cambridge University Press, Cambridge, second edition, 1989.
\newblock Translated from the Japanese by M. Reid.

\bibitem[Mer72]{meredith-weakformalschemes}
David Meredith.
\newblock Weak formal schemes.
\newblock {\em Nagoya Math. J.}, 45:1--38, 1972.

\bibitem[Mon02]{these_montagnon}
Claude Montagnon.
\newblock {\em {G{\'e}n{\'e}ralisation de la th{\'e}orie arithm{\'e}tique des
  $\mathcal{D}$-modules {\`a} la g{\'e}om{\'e}trie logarithmique}}.
\newblock PhD thesis, Universit{\'e} de {R}ennes {I}, 2002.

\bibitem[MW68]{MonskyWashnitzer}
P.~Monsky and G.~Washnitzer.
\newblock Formal cohomology. {I}.
\newblock {\em Ann. of Math. (2)}, 88:181--217, 1968.

\bibitem[NvO82]{GradedRingTheory}
C.~N{\u{a}}st{\u{a}}sescu and F.~van Oystaeyen.
\newblock {\em Graded ring theory}, volume~28 of {\em North-Holland
  Mathematical Library}.
\newblock North-Holland Publishing Co., Amsterdam-New York, 1982.

\bibitem[Ogu18]{Ogus-Logbook}
Arthur Ogus.
\newblock {\em Lectures on logarithmic algebraic geometry}, volume 178 of {\em
  Cambridge Studies in Advanced Mathematics}.
\newblock Cambridge University Press, Cambridge, 2018.

\bibitem[Rot09]{Rotman-HomAlg}
Joseph~J. Rotman.
\newblock {\em An introduction to homological algebra}.
\newblock Springer-Verlag, New York, 2009.

\bibitem[Ser68]{Serre-corpslocaux}
Jean-Pierre Serre.
\newblock {\em Corps locaux}.
\newblock Hermann, Paris, 1968.
\newblock Deuxi\`eme \'edition, Publications de l'Universit\'e de Nancago, No.
  VIII.

\bibitem[sga03]{sga1}
{\em Rev\^etements \'etales et groupe fondamental ({SGA} 1)}.
\newblock Documents Math\'ematiques (Paris) [Mathematical Documents (Paris)],
  3. Soci\'et\'e Math\'ematique de France, Paris, 2003.
\newblock S\'eminaire de g\'eom\'etrie alg\'ebrique du Bois Marie 1960--61.
  [Geometric Algebra Seminar of Bois Marie 1960-61], Directed by A.
  Grothendieck, With two papers by M. Raynaud, Updated and annotated reprint of
  the 1971 original [Lecture Notes in Math., 224, Springer, Berlin; ].

\bibitem[sta]{stackproject}
Stack project.

\bibitem[Tsu19]{Tsuji-saturated}
Takeshi Tsuji.
\newblock Saturated morphisms of logarithmic schemes.
\newblock {\em Tunisian Journal of Mathematics}, 1(2):185--220, 2019.

\bibitem[Vir04]{Vir04}
Anne Virrion.
\newblock Trace et dualit\'e relative pour les {$\mathcal{D}$}-modules
  arithm\'etiques.
\newblock In {\em Geometric aspects of Dwork theory. Vol. I, II}, pages
  1039--1112. Walter de Gruyter GmbH \& Co. KG, Berlin, 2004.

\bibitem[Wei94]{Weibel-HomologicalAlg}
Charles~A. Weibel.
\newblock {\em An introduction to homological algebra}, volume~38 of {\em
  Cambridge Studies in Advanced Mathematics}.
\newblock Cambridge University Press, Cambridge, 1994.

\end{thebibliography}

\def\cprime{$'$}

\bigskip
\noindent Daniel Caro\\
Laboratoire de Mathématiques Nicolas Oresme\\
Université de Caen Normandie
Campus 2\\
14032 Caen Cedex\\
France.\\
email: daniel.caro@unicaen.fr

\end{document}